\documentclass[12pt]{article}
\usepackage{geometry}                % See geometry.pdf to learn the layout options. There are lots.
\usepackage{graphicx}
\usepackage{amssymb}
\usepackage{psfrag}

\usepackage{color}

\usepackage{amsthm}

\def\Bbb{\mathbb}
\def\R{{\Bbb R}}
\def\N{{\Bbb N}}

\def\mbf0{\mathbf{0}}

\newcommand{\bmat}{\begin{pmatrix}}
\newcommand{\emat}{\end{pmatrix}}

\theoremstyle{theorem}
\theoremstyle{theorem}
\theoremstyle{theorem}
\theoremstyle{definition}
\theoremstyle{remark}

\title{\bf The Geometry of Oriented Cubes\\
{ \small   Macquarie University Research Report No: 86-0082}}

\author{\bf \Large  Iain R. Aitchison\\
}
\date{}                                           % Activate to display a given date or no date

\def\two#1#2{$\, \begin{array}{rr} 
\hbox to 0pt{$#1$}&  
\hbox to 0pt{$#2$}   
 \end{array}\ $}

\def\ftwo#1#2{$\hskip-6pt\begin{array}{rr} 
\hbox to -2pt{$#1$}&  
\hbox to -2pt{$#2$}   
 \end{array} $}

\def\three#1#2#3{$\begin{array}{rrr} 
\hbox to -1pt{$#1$}&  
\hbox to -1pt{$#2$} & 
\hbox to -1pt{$#3$} \end{array}\ $\ }

\def\shree#1#2#3{\small$\begin{array}{rrr} 
\hbox to -1pt{\small${\small #1}$}&  
\hbox to -1pt{\small${\small #2}$} & 
\hbox to -1pt{\small${\small #3}$} \end{array}\ $\ }

\def\free#1#2#3{\small$\begin{array}{rrr} 
\hbox to -3pt{\small${\small #1}$}&  
\hbox to -3pt{\small${\small #2}$} & 
\hbox to -3pt{\small${\small #3}$} \end{array}\ $\ }

\def\free#1#2#3{$\hskip-6pt\begin{array}{rrr} 
\hbox to -2pt{$#1$}&  
\hbox to -2pt{$#2$} & 
\hbox to -2pt{$#3$} \end{array}$ }

 \def\four#1#2#3#4{$\hskip -0.5pt  \begin{array}{rrrr} 
{\centering\hbox to  -0.6pt{$#1$}}&  
\hbox to -0.6pt{${\centering{#2}}$} & 
\hbox to -0.6pt{${\centering{#3}}$} &
\hbox to -0.6pt{$#4$} 
\end{array}\  $\ }

 \def\ffour#1#2#3#4{$\hskip -8pt  \begin{array}{rrrr} 
{\centering\hbox to  -2pt{$#1$}}&  
\hbox to -2pt{${\centering{#2}}$} & 
\hbox to -2pt{${\centering{#3}}$} &
\hbox to -2pt{$#4$} 
\end{array}  $}

\def\five#1#2#3#4#5{$\, \begin{array}{rrtrr} 
\hbox to 1pt{$#1$}&  
\hbox to 1pt{$#2$} & 
\hbox to 1pt{$#3$} &
\hbox to 1pt{$#4$} &
\hbox to 1pt{$#5$} 
\end{array}\  $}

\begin{document}
\maketitle
 
 \begin{abstract}

 This reports on the fundamental objects revealed by Ross Street,
which he called `orientals'. Street's work was in part inspired by Robert's attempts 
\cite{Ro} to use $n$-category ideas to  construct nets of $C^*$-algebras in Minkowski space for applications to relativistic quantum field theory: Roberts'  additional  challenge was that ``no
amount of staring at the
low dimensional cocycle conditions would reveal
the pattern for higher dimensions''.  

This report takes up this challenge,  presents a natural inductive construction of explicit cubical cocyle conditions,  and gives three ways in which the simplicial ones can be derived from these. 
%The dual `string diagram' version of this work, giving rise to a Pascal's triangle of diagrams for cocycle conditions, will be exposited elsewhere: part of this work has been described in \cite{St1}. 
A consequence of this work is that the Yang-Baxter equation, the `pentagon of pentagons', and higher simplex equations, are in essence different manifestations of the same underlying abstract structure.

There has been recent interest in higher-categories, by  computer scientists investigating concurrency theory,
as well as by physicists, among others. The dual `string' version of this paper, originally presented in \cite{Ai1}, 
and which motivated much of the role of cubes in
\cite{St1}, makes clear the relationship with higher-dimensional simplex equations in physics. 

Much work in this area has been done since these notes were written: no attempt has been made to update the original report. However, all diagrams have been redrawn by computer, replacing all original hand-drawn pictures.

\end{abstract}

\pagebreak

We obtain a geometric interpretation of Street's simplicial orientals, by describing
an oriental structure on the $n$-cube. This corresponds to describing nested
sequences of embedded disks in the boundary of the $n$-cube, obtained
consecutively as ordered deformations keeping the boundaries fixed. 
Canonical forms for the order are obtained, leading to definitions 
of
cocycle conditions for application to non-abelian cohomology in 
a category.

\section{Introduction }

Street \cite{St} has defined an $n$-category structure on the
$k$-simplex, leading to his notion of $orientals$.
This was motivated by an approach to non-abelian cohomology, arising from
quantum field theory \cite{Ro}, but also seems
to have intriguing connections with homotopy coherence, loop space
desuspension and the cobar construction, as well as the realisation of
homotopy types. Intimately related to Street's structure on simplices is an
analogous one on cubes, from which the simplicial results can be
derived. This 
paper takes up Roberts' challenge \cite{Ro} that ``no
amount of staring at the
low dimensional cocycle conditions would reveal
the pattern for higher dimensions''. The result is surprisingly simple,
natural and beautiful. We define and describe this finer
structure of oriented cubes, from a geometric point of view.
This makes obvious the relevance to homotopy theory.

This 
paper is self-contained,
presupposing minimal familiarity with either category theory
or homotopy theory.
Categorical aspects and applications will be discussed in a subsequent
paper. 

Consider the following geometric reasoning: 
In obstruction theory applied to simplicial or cubical complexes, a
standard problem is the continuous extension of data defined on
the boundary of a cell (usually a simplex or cube). 
Although the boundary of such a cell is a geometrically decomposed sphere, the
only structure used involves assigning alternating signs appropriately to
the sub-cells of the decomposition. 
If we wish to keep the spherical structure in mind, we can imagine the
boundary sphere of the cell decomposed as two hemispheres, with an
extension over the cell of the data assigned to the hemispheres now viewed
as a continuous deformation, keeping the boundary data fixed, between the assignments of data to each of the hemispheres. Since the boundary of each
hemisphere is also  a sphere, it is natural to seek splittings of
spheres in every dimension.

Suppose we have some ``data'' defined over a sphere, and two possibly distinct
extensions over a cell. By viewing two cells as upper and lower hemispheres, 
we may instead consider this as an extension, over a sphere of one
dimension higher, of data assigned to the ``equatorial'' sphere. We can then
think of the two extensions over the hemispheres as being ``equivalent'' 
in some sense if the data assigned to this higher sphere extends
over a cell. But there may be
another such extension, tantamount to an alternative equivalence, which
enables us to bootstrap the procedure. At any stage, we see data on a
sphere, split as two hemispheres, with data on the equatorial sphere, which
in turn splits as two hemispheres, whose boundary sphere
splits as two hemispheres, and so on.

\begin{figure}[htbp] %  figure placement: here, top, bottom, or page
   \centering
   \includegraphics[width=4in]{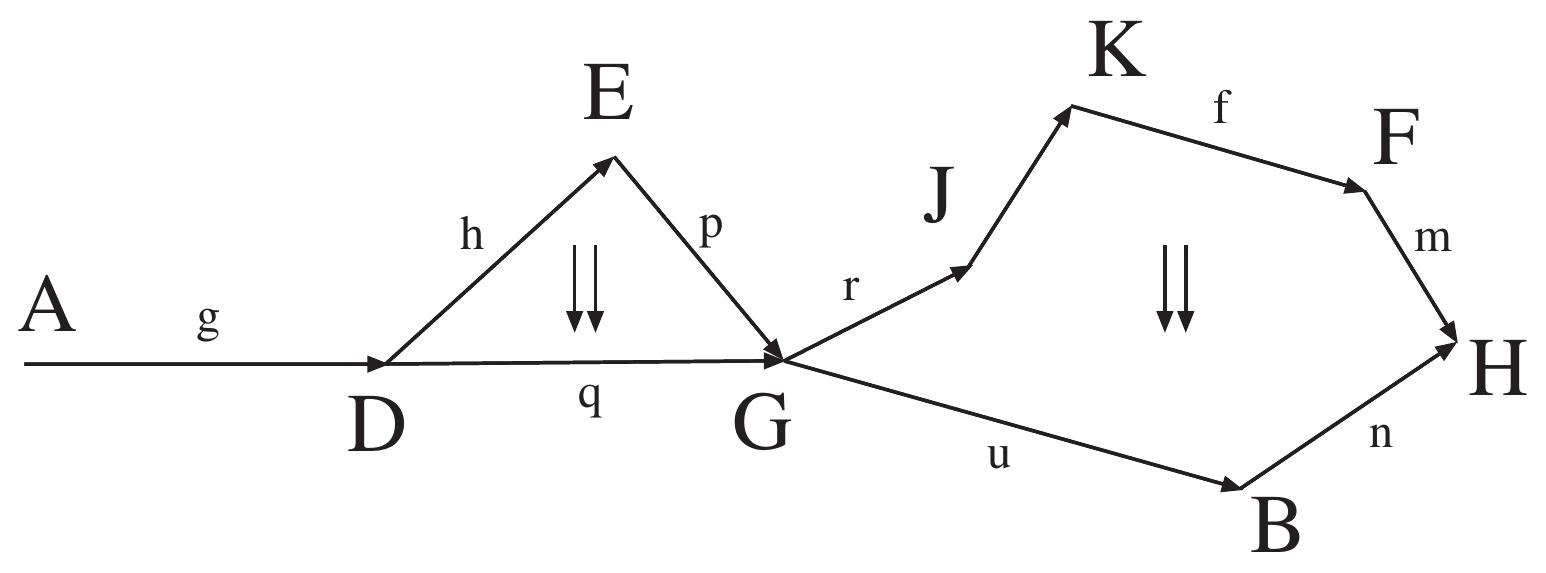} 
   \caption{Typical 2-category pasting scheme: boundaries of 2-cells naturally split as source- and target-hemispheres}
   \label{fig:1}
\end{figure}

Consider now category theory, viewing arrows (or morphisms) as
one dimensional cells whose boundary 0-spheres split as the union of a
source and a target - which we view as 0-dimensional points. If we
consider a 2-category (see \cite{K-S} for example),
%\gapp2
%end page 2
%\noindent
we see that 2-cells have a 1-dimensional source and
target, each of 
which is an ordinary arrow in a category, and both of which
have the same 0-dimensional source and target.
Thus we see the boundary of a 2-cell split into two
hemispheres, and the boundary of each hemisphere similarly split. 
(Figure \ref{fig:1})

Street's
$oriental$
structure on simplices \cite{St} allows an interpretation as a coherent version,
across all dimensions, of a structure as described above.
Our purpose is to exhibit such a structure on the 
$n$-dimensional cube. 
This leads to the notion of a ``$k$-dimensional source'' and
``$k$-dimensional target'' of the $n$-cube, for each $k$,
just as a 2-cell in a 2-category has sources in dimensions 0 and 1.
Amusingly, the answer suggests a curious connection with $p$-forms on $\R^n$.
From a categorical viewpoint this decomposition 
leads to the notion of an
$n$-$category$. 
%\gapp3
%end page 3

%\pagebreak

\section {\bf  The structure of the $n$-cube}

Denote by $\underline n $ the (ordered) set $\{1,\ 2,\dots ,\ n\}$,
and by $ {\cal I } = \{-,\ 0,\ +\}\  \cong \ \underline 3$. 
Let $ {I} $ denote the interval 
$[-1,\ 1] $ in the real number axis, and
$ {I}^n $ the 
(geometric) $n$-dimensional 
cube. We think of
a point as the 0-dimensional cube or simplex.

\begin{figure}[htbp] %  figure placement: here, top, bottom, or page
   \centering
   \includegraphics[width=5in]{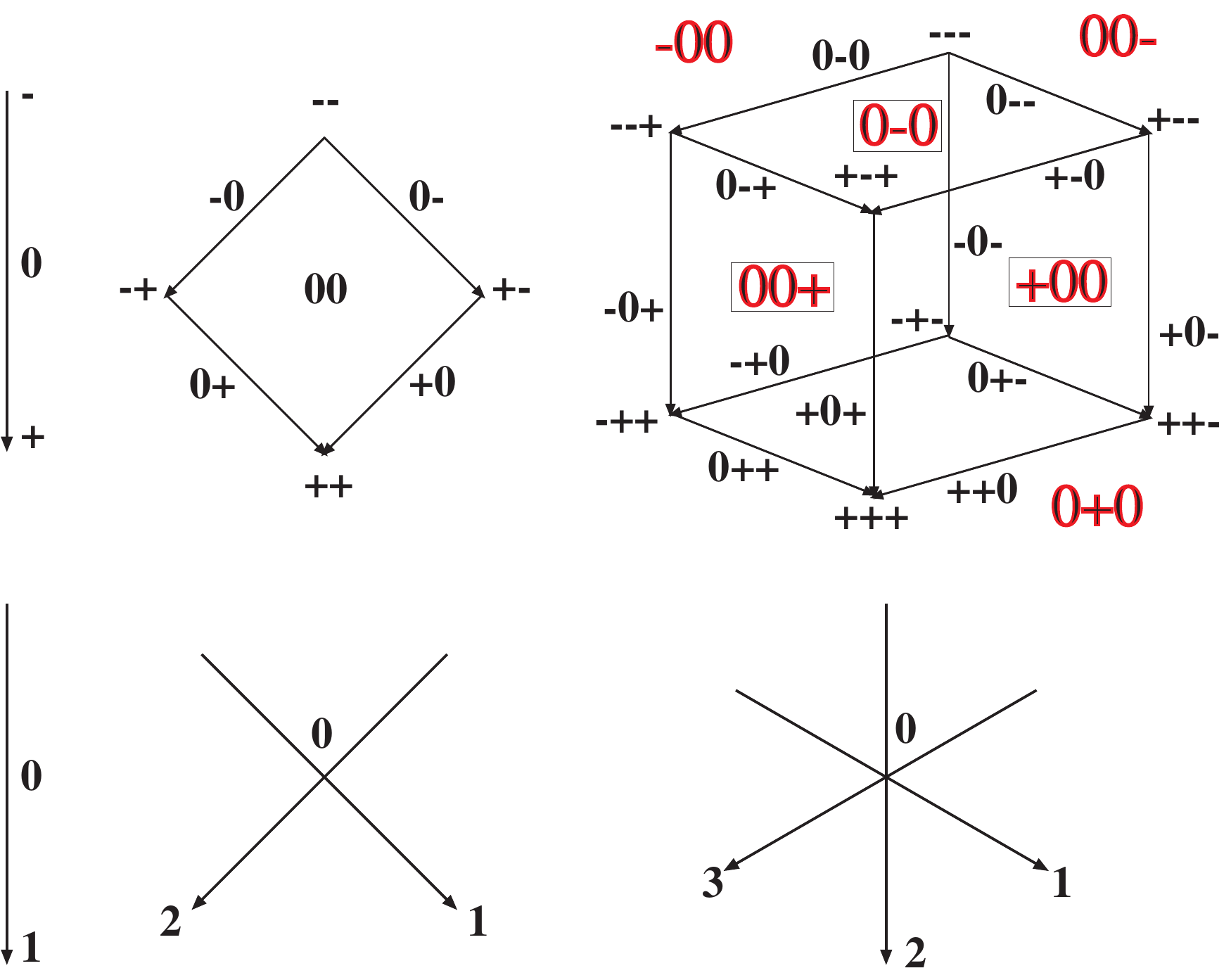} 
   \caption{Labeling of cubes in dimensions 1,2,3, and their dual representations}
   \label{fig:2}
\end{figure}

\medskip
\medskip
\noindent
{\bf Dimension 1}: 
In Figure \ref{fig:2} we have a single $ {I}^1 $, with a notion of 
source (``0-source'') and target. We may think of this as an arrow 
assigned to the interval 
$I  $. As a convention, 
we choose 
``$-1$'' 
to be considered the source and ``$+1$'' the target. Hence the dichotomy 
``source/target'' arises from the disconnectedness of the 
0-sphere 
$ S^0 = \partial I$. It is convenient for ``$-$'' (respectively ``$+$'') to be 
considered synonymous with ``$-1$'' (respectively $+1 = 1$). Moreover, we
use ``$0$" to denote the whole interval:

Observe that $ {I}^1 $ has three sub-cubes, $+,\ 0$ 
and $-$.
When we take the product with the unit interval we obtain 
${I }^2 $, 
to which each sub-cube
of ${I}^1 $ contributes three sub-cubes. 
These can be thought of as a new
copy, an old copy and a thickened copy. Proceeding inductively we 
see that
${I}^n $ has $ 3^n $ sub-cubes. More precisely, 
let ${\cal I}^{ n } = \{ x:\ \underline n \ \longrightarrow \  
{\cal I}  \} $. 

\medskip
\medskip
\noindent
{\bf Proposition 2.1}:
 There is a 1:1 correspondence between (geometric) sub-cubes of 
${I }^n $ and elements of ${\cal I}^n$. 
\medskip
\medskip

Thus we may
{\it  define}
the
n-cube
to be ${\cal I}^n $, and 
call an element $x$ of ${\cal I}^n  $ a 
$sub$-$cube$
of $ {\cal I}^n $. We will henceforth use the notation  
${\cal I}^n  $, or occasionally [$n$],
to denote the $n$-cube, and no longer use $I^n$.

\medskip
\medskip
\noindent
{\bf Notation}:
Since each $ \underline n$ is ordered, an element 
$x$ of ${\cal I}^n $ can be specified as
a word of length $n$ in the symbols of $ \cal I $.

%\gapp4
%end page 4
\medskip
\medskip

\noindent
{\bf Example 2.2}:
In Figure \ref{fig:2} we draw the 1, 2, and 3 dimensional cubes, 
labeled according to
the usual convention arising from the unit interval in $ R^1 $.

\medskip
\medskip
\noindent
Define maps  $ \lambda , \ \mu ,\ \nu  : {\cal I}^n 
\ \longrightarrow \ {\cal I}^{n+1}$ 
by 
$$(\lambda (x))(i) = x(i)\ \ {\rm for}\ \  i \leq n,\ \ \ \ 
(\lambda (x))(n+1)\ =\ -$$
$$(\mu(x))(i) = x(i)\ \ {\rm for}\ \  i \leq n,\ \ \ \ 
(\mu(x))(n+1)\ = \ 0$$
$$(\nu(x))(i) = x(i)\ \ {\rm for}\ \  i \leq n,\ \ \ \ 
(\nu(x))(n+1)\ = \ +$$
\medskip

These add to the right of a word the symbols  $-,\ 0$ and $+$  
respectively, and
can be thought of as the old, thickened and new copies of $x$.

\medskip
\medskip
\noindent
{\bf Definition 2.3}:
The 
$dimension$ of $x$, denoted
by $|x|$, is the cardinality of $x^{-1}(0)$.

\medskip
\medskip
This is the number of 0's occurring in the word-form of $x$,
and induces a grading  
$ {\cal I}^n\  = \ {\Sigma}_{j=o}^{j=n}\ 
{\cal I}^n_j $ 
by setting $x  \in  {\cal I}^n_j \ 
\Longleftrightarrow \ |x| = j$.

%\gapp5
% end page 5

\bigskip

\noindent 
There are three distinguished subcubes of ${\cal I}^n$, which in word form
are $0\cdots 0$,
$-\cdots -$ and
$+\cdots +$.
We shall respectively refer to the latter two as the
0-source and 0-target of ${\cal I}^n$, denoted by
${\sigma}_0[n]$ and ${\tau}_0[n]$. More generally we will
define 
$ {\sigma}_k[x] $ and $ {\tau}_k[x] $ for each $x\in {\cal I}^n $.
These are to be thought of as ``$k$-dimensional source'' and 
``$k$-dimensional 
target'' of $x$ respectively. We shall suppress some symbols 
when the context is
clear. 

For each $x \in  {\cal I}^n$ of dimension $p$ and 
$ y \in {\cal I}^p $, 
we can define a new element $ x*y $.
Set ($ x*y )(k) = x(k)$ if $x(k)\ \not=  0$. If 
$ x^{-1} (0) = \{ k_1,\ \dots \ ,\ k_p \}$, 
listed in increasing order, define $ x*y $ by performing $x$ on 
$ x^{-1}(-,\ +) $ and $y$ otherwise, using the natural
order: ($x*y)( k_j ) = y(j)$. 
This gives a pairing
$$ {\cal I}^n_p\  \times \ {\cal I}^p_q\ \longrightarrow \ 
{\cal I}^n_q. $$
Note that this defines an associative operation as with matrix
multiplication. 

We may think of any such $x$ as a $p$-cube placed in 
${\cal I}^n $.
More generally we can use $*$ to place any set of sub-cubes of the 
$p$-cube
into $ {\cal I}^n$.
Once we define $ {\sigma }_k $ and $ {\tau}_k $ on 
$ {\cal I}^p$, we can extend to $x$ by 
$ {\sigma}_k [x] = x* ( {\sigma}_k [ p ]) $, 
and similarly for $ {\tau}_k $.

%\gapp6
% end page 6

%\pagebreak

  \section{\bf Motivational examples}

\begin{figure}[htbp] %  figure placement: here, top, bottom, or page
   \centering
   \includegraphics[width=4in]{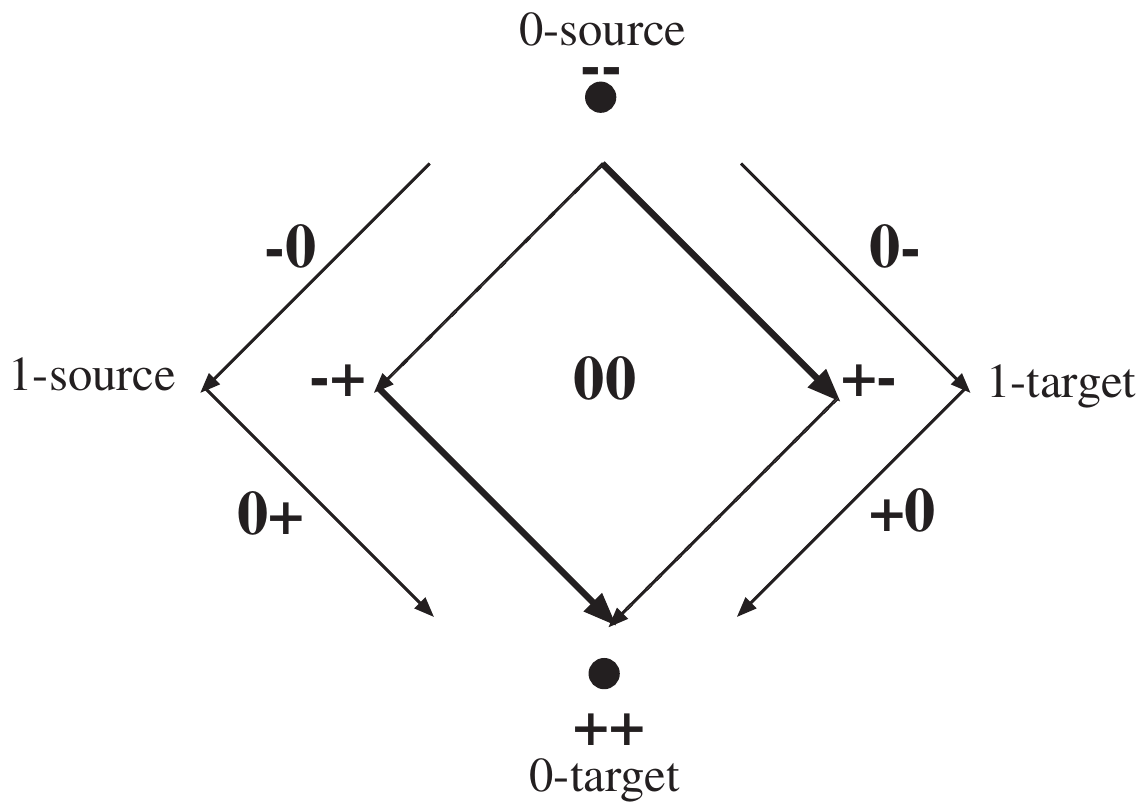} 
   \caption{The 2-cube with all sources and targets, which are $00$ for dimension 2 or more.}
   \label{fig:3}
\end{figure}

We examine the lowest dimensional cubes to elucidate their 
structure:

\medskip
\medskip
\noindent
{\bf Dimension 2}: 
The boundary of $ {\cal I 
}^2 \ \cong \ I^2 $ consists of 
four edges and
four vertices. To each of the edges is assigned naturally, 
given our convention, an orientation. This canonically determines 
the source and target in
dimension 0, as shown in Figure \ref{fig:3}. Since the edge $\{ -1 \} \times \ I $
(respectively $\{ 1\} \times \ I $) arises as the thickened 0-source
(respectively 0-target) of $ {\cal I}^1 $, a natural criterion 
for choosing
the one dimensional
source and target for ${\cal I}^2 $ is determined. Hence we 
may assign a
``double arrow'' to $ {\cal I}^2 $, pointing from the 1-source 
to the 1-target.
Note that geometrically the boundary circle (1-sphere) of 
${\cal I}^2$ is 
decomposed into two hemispheres, each having as boundary the 
0-source and 
0-target.
 %\begin{minipage}{.5\columnwidth}
%\end{minipage}
\medskip
\medskip

\begin{figure}[htbp] %  figure placement: here, top, bottom, or page
   \centering
   \includegraphics[width=5in]{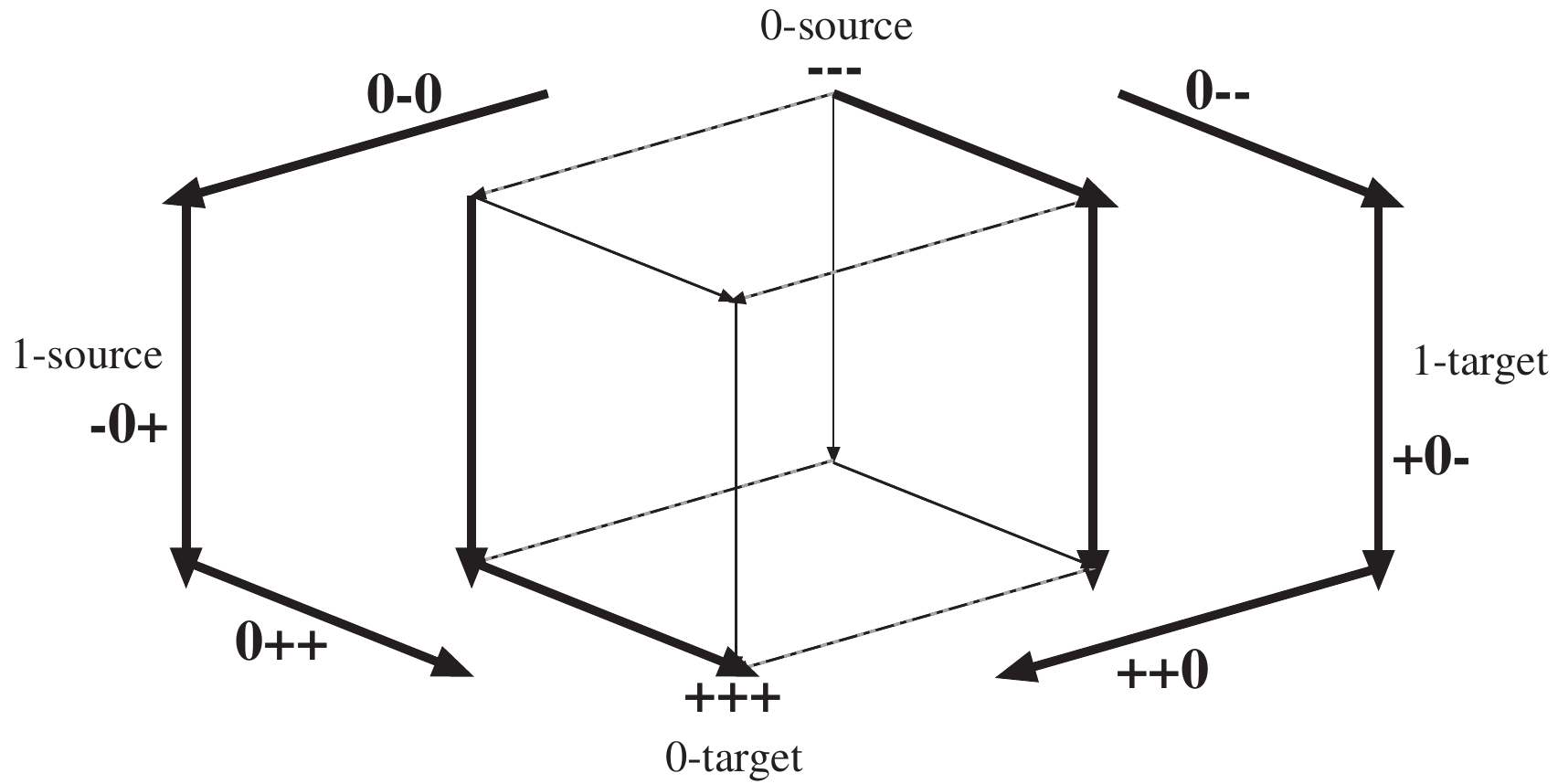} 
   \caption{0- and 1-dimensional source and target of the 3-cube}
   \label{fig:4}
\end{figure}

\begin{figure}[htbp] %  figure placement: here, top, bottom, or page
   \centering
   \includegraphics[width=4.5in]{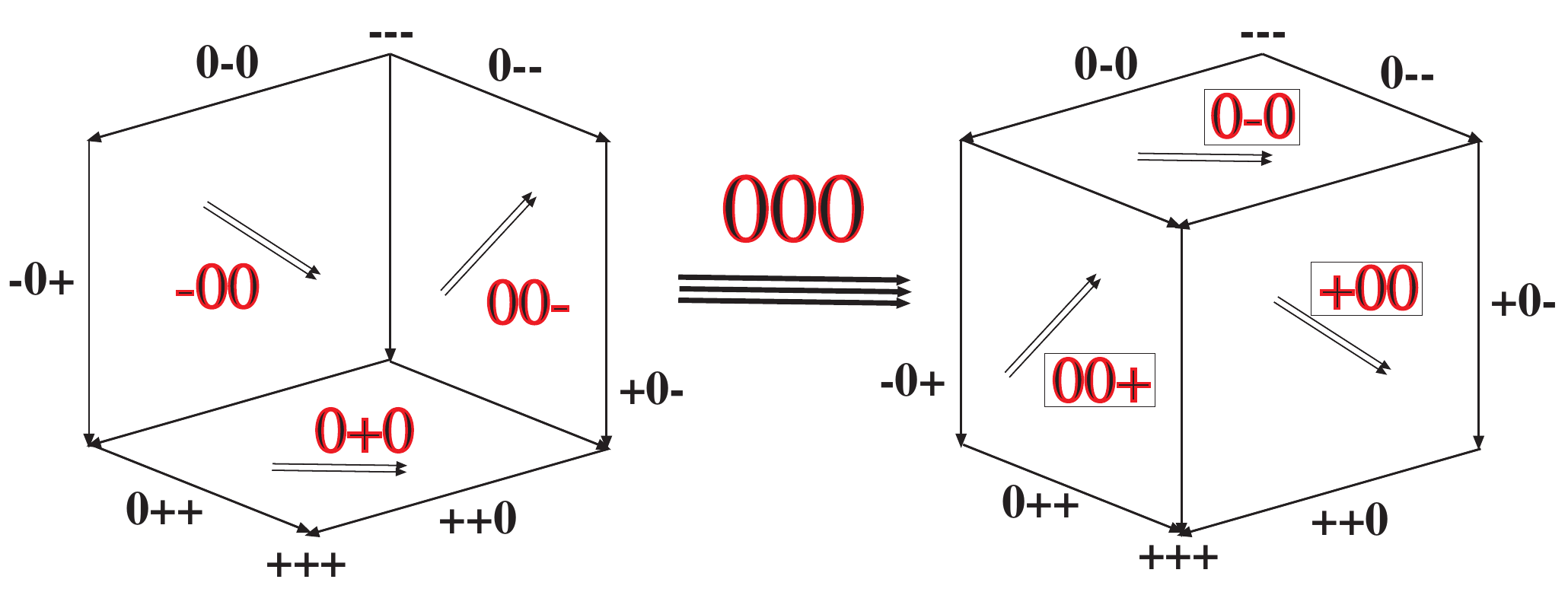} 
   \caption{The 2-source and 2-target of the 3-cube}
   \label{fig:5}
\end{figure}

\noindent
{\bf Dimension 3}: 
In Figure \ref{fig:4}, we thicken the various subcubes of 
$ {\cal I}^2$ to obtain the usual
description of ${\cal I}^3$. Again the 0-dimensional sources 
and targets are
determined, as are the 1-dimensional ones. Observe that the 1-source 
and
1-target are symmetrically defined, and that each consists of a 
copy of
the 2-dimensional ones, suitably shifted by the third dimension, 
together
with a thickened 0-dimensional source or target. Applying the same
reasoning to the definition of 2-source and 2-target as we did for 
dimension 1,
we see the 2-source and
2-target of $ {\cal I}^3 $ as depicted in Figure \ref{fig:5}. The boundary 
of $ {\cal I}^3 $ 
is a 
%\gapp7
% end page 7
2-sphere splitting as two hemispheres, the 2-dimensional source
(at the back) and target (at front) of ${\cal I}^3 $, which 
in turn have 
boundary 1-sphere splitting as two hemispheres, the 1-dimensional 
source and target of $ {\cal I}^3 $, which in
turn have boundary 0-sphere splitting as two hemispheres, the 
0-dimensional
source and target of ${\cal I}^3 $.

\begin{figure}[htbp] %  figure placement: here, top, bottom, or page
   \centering
   \includegraphics[width=4.5in]{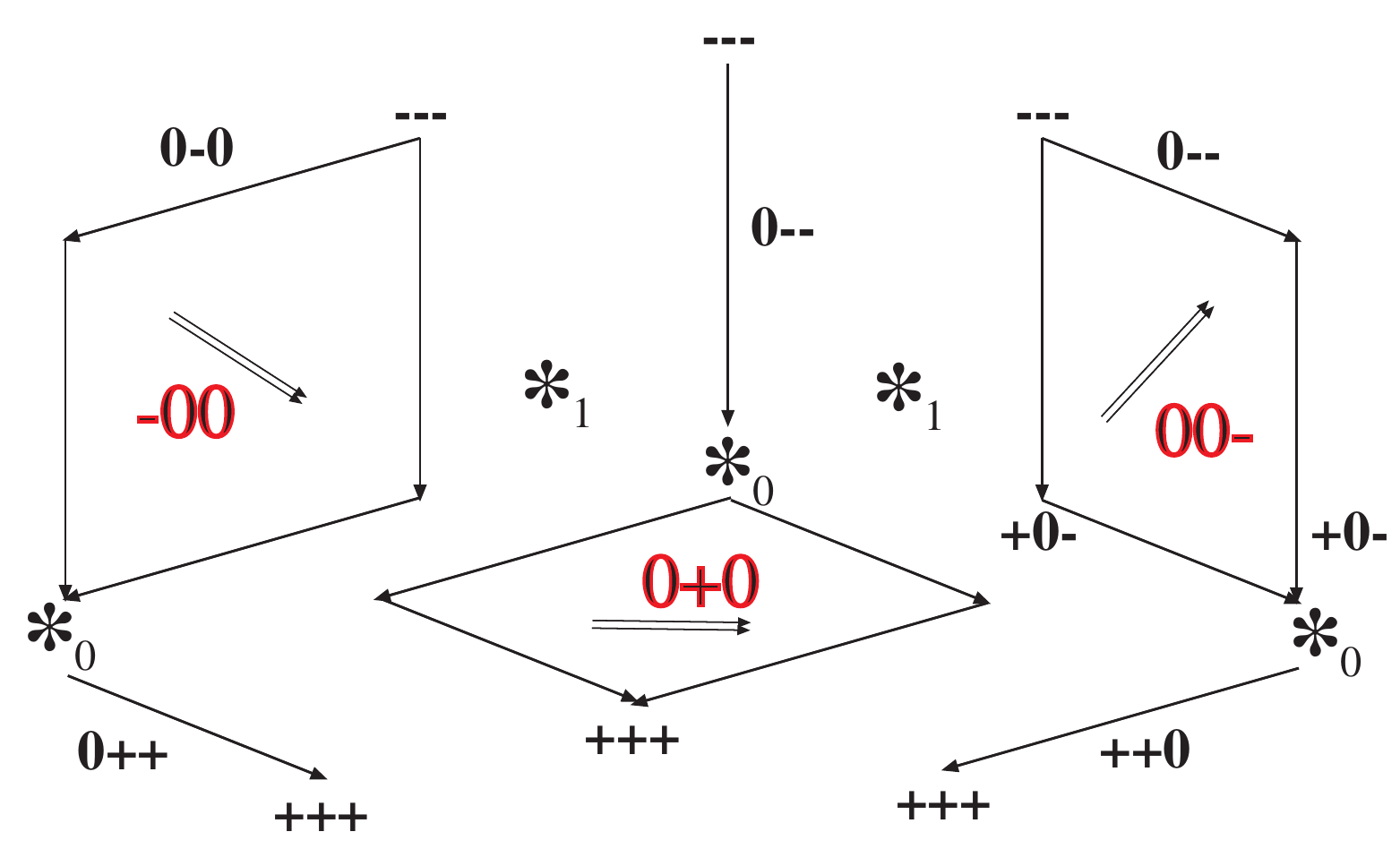} 
   \caption{The 2-source of the 3-cube: we add 1-cells to the 2-cells, enabling matching}
   \label{fig:6}
\end{figure}

\begin{figure}[htbp] %  figure placement: here, top, bottom, or page
   \centering
   \includegraphics[width=6in]{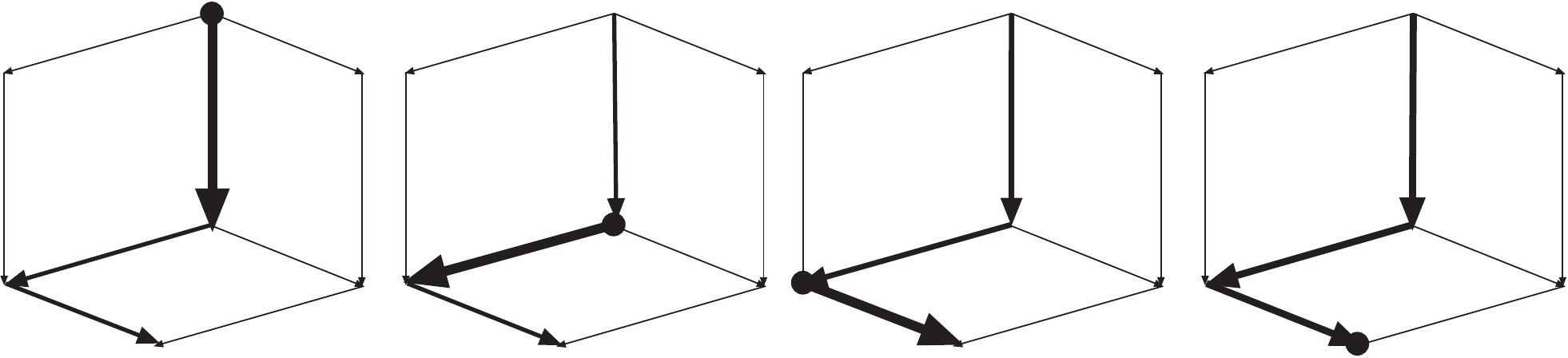} 
   \caption{Each 1-path provides a route from the 0-source to the 0-target}
   \label{fig:7}
\end{figure}

\medskip
\medskip
\noindent
{\bf Remark}:
Although there is obvious appeal in the simplicity of the preceding
description of the 2-source and 2-target of 
$ {\cal I}^3 $ as the union of
three $ {\cal I}^2$s each, observe that the 
${\cal I}^2$s concerned do not 
have the same 0-source and 0-target, nor the same 1-source 
%\gapp8
% end page 8
and 
1-target. Hence a little care is required to get the definitions 
correct. At this level,
the answer both algebraically and geometrically is provided by the 
notions
of 2-category theory and the middle-interchange law 
(\cite{ K-S, St}). 
This is described in Figure \ref{fig:6}. Each of the 
$ {\cal I}^2$s comprising the 2-source 
of $ {\cal I}^3 $ must be stretched out by an appropriate 
$ {\cal I}^1 $ to the 
0-source or 0-target of $ {\cal I}^3 $.
An interpretation of this goes as follows: each 1-dimensional path 
from the
0-source to the 0-target can be thought of as the trace of a 
deformation or
homotopy of the 0-source to the 0-target, consecutively across the 
arrows
forming the path. See Figure \ref{fig:7}. In ${\cal I}^2$, there are only 
two possible paths, the 1-source and the 1-target, and 
$ {\cal I}^2$ itself 
can be thought of as providing a deformation from the 1-source to 
the 1-target. For $ {\cal I}^3 $, there are four possible 
paths for each of the 2-source and
2-target. 

\begin{figure}[htbp] %  figure placement: here, top, bottom, or page
   \centering
   \includegraphics[width=6in]{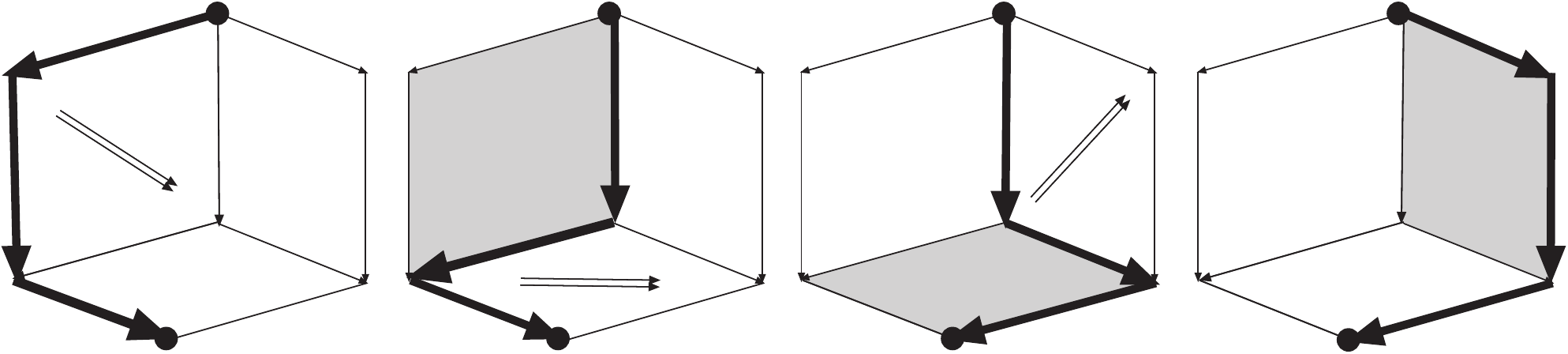} 
   \caption{2-cubes provide deformations of the 1-source to the 1-target; each such deformation first requires the presence of all  of its 1-source cubes. White dots represent 0-sources of 2-cubes; grey dots the target 0-cube of the previous 2-cube, which subsequently are moved by 1-cubes.}
   \label{fig:8}
\end{figure}

In the same way, each of the 2-source and 2-target can be thought 
of as 
deformations, in three stages each, of the 1-source to the 1-target. 
The 
source deformations must be carried out in the order shown in Figure \ref{fig:8}. 
This is 
because we cannot deform a 1-path across a 
$ {\cal I}^2$ unless all of the 
1-source of the ${\cal I}^2$ in question is part of the 
1-path.  
Each such deformation leaves unchanged the rest of the 1-path, 
and the number 
of segments in the path remains the same. 
Each component of the 2-source consists of a 
$ {\cal I}^2 $ across which part of the 1-path is being 
deformed, and the 
remainder 
%\gapp9 
of the 1-path. From another viewpoint, the 
${\cal I}^1 $s string out 
the ${\cal I}^2 $s, providing the path by which the 0-source 
of $ {\cal I}^3 $ 
can be deformed to the 0-source of the relevant 
$ {\cal I}^2$. This gives rise 
to the $ *_0 $ composition of 2-category theory.

\noindent
Note that in this case we are composing cubes stuck end-to-end, 
source-to-target, along the 0-dimensional (geometric) sources and 
targets.
The ${\cal I}^2$s, on the other hand, are glued together along 
1-dimensional
source/target paths. This corresponds to $ *_1 $ in the language of 
2-categories.The geometry is exactly that of 2-categorical pasting 
diagrams. 
(See for example \cite{K-S}.)

Thus not only are the constituent parts of the 2-source of 
${\cal I}^3 $
derived from the structure of 
$ {\cal I}^2 $, but a natural order emerges for
their 2-categorical order of composition. 
%\gapp{10}
This order in turn derives 
entirely 
from our initial convention for source and target of 
$ {\cal I}^1 $. We also
note that the 2-cubes of the source or target fit together in such a 
way
that the deformations are possible, and compatibly with the 
structure
defined on ${\cal I}^2 $.
%\gapp{11}

%\pagebreak

\section{\bf Source and target faces - Piles}

We require $ {\cal I}^n$ to have sources and targets in all dimensions 
$ \leq $ $n$. 
These allow us to deform the ($k-1$)-source, across either the
$k$-dimensional source or target, onto the ($k$-1)-target. 
The ($k$-1)-source is made up of 
${\cal I}^{k-1} $s as well as lower dimensional cubes. 
The lower dimensional cubes stretch out each 
$ {\cal I}^{k-1} $
to the sources and targets of ${\cal I}^n$ in all dimensions 
$ \leq$ ($k-3$). 

In deforming the $(k-1)$-source across to the ($k-1)$-target, such 
lower 
dimensional cubes play a role; we will not be able to push across a 
given 
$ {\cal I}^{j} $ unless all of
$its$ 
lower dimensional sources have been reached. The order in which we 
proceed is 
determined inductively, necessitated by the need to expose all lower 
dimensional sources as required (compare Figure \ref{fig:8}). The notion of 
``pile'' enables us to describe
the procedure to any desired degree of refinement. A pile is a list 
of all 
${\cal I}^k $s we have reached at any stage of the 
deformation. Thus there are 
source and target piles, and intermediate piles. Each time we deform 
across a 
${\cal I}^k $ we alter the pile in all dimensions less than 
$k$, since we deform from the source of that 
${\cal I}^k $ to its target in
$all$
lower dimensions. In this way sub-cubes act as ``operators'' on 
piles. The
order of deformation is then determined by groupings of sub-cubes 
as ``blocks''. 
These blocks are characterised in that the initial and final
piles on which they act respectively contain the $j$-dimensional 
source and target faces
of the
${\cal I}^n$ in question. 
The $j$-dimensional sources and targets of this
$ {\cal I}^n $ are themselves blocks, acting on the piles 
corresponding to the
still lower dimensional sources and targets. 

Our immediate concern is to develop an appropriate formalism and 
notation
to make these ideas concise for the $n$-cube. We take an 
inductive approach 
starting with the raw material above.

For ${\cal I}^3 $ we saw that there are three 
$ {\cal I}^2 $s making up the 
2-dimensional source. If we ignore the order in which these cubes 
are written,
as well as their position in $ {\cal I}^3 $, we have a set of 
$ {\cal I}^2 $s
we shall call the
$source\ 2$-$face$
of ${\cal I}^3 $. More generally we will 
%\gapp{12}
%\noindent
define the source and 
target $k$-$faces$ of ${\cal I}^n $, denoted $ {\psi}_k (n)$ and 
$ {\omega}_k (n) $ respectively. For given $n$, we obtain two 
sequences of sets 
$ {\psi}_0(n) ,\dots , {\psi}_{n-1} (n) ,\ {\psi}_n (n) $, 
and
${\omega}_0(n),\dots , {\omega}_{n-1} (n), \  {\omega}_n (n) $. 
The
subscript indicates the dimension of elements of the set.

\medskip
\medskip
\noindent
{\bf Definition 4.1}:
Let $ {\psi}_0[n]  = {\sigma}_0 [n] $, 
$ {\omega}_0[n] = {\tau}_0 [n] $ and 
${\psi}_j [n] = [n] = {\omega}_j [n] $, for each $j\ \geq \ n$. 
Inductively define
\hfil
\break
\noindent 
a)   For $k$ even:  
$$
{\psi }_k [ n ] = {\mu} {\psi}_{k-1} [n-1]  \ \cup \  \lambda 
{\psi}_k [n-1] 
$$
$$
{\omega}_k [ n ]  = \nu {\omega}_k [n-1] \ \cup \ \mu 
{\omega}_{ k-1}  [n-1]  
$$
\hfil
\break
\noindent
b)   For $k$ odd: 
$$
{\psi}_k [ n ]  = \mu {\psi}_{k-1} [n-1]  \ \cup \ 
\nu {\psi}_k [n-1]
$$
$$
{\omega}_k [ n ] = \lambda {\omega}_k [n-1] \ \cup \ 
\mu {\omega}_{ k-1}  [n-1]  
$$

\medskip
\medskip
We have deliberately written these as (inductively) ordered, 
since the
order reflects their sequence of operation on piles.
A simple induction from the motivational examples proves 

\medskip
\medskip
\noindent
{\bf Theorem 4.2}:
(i)   Identifying elements of either $ {\psi}_k [n]  $ or
$ {\omega}_k [n] $ along common boundaries produces a 
$k$-dimensional disk
$D_k $.
\hfil
\break
\noindent
\ \ \ \ \ \ \ \ \ \ \ \ \ \         (ii)   $\partial D_k [n] = 
{\psi}_{k-1} [n] \ \cup \ 
{\omega}_{k-1} [n]  = \partial {\omega}_k [n] .$
\medskip
\medskip

\noindent
{\bf Examples 4.3}:
In the following table we list the source and target faces for 
$n = 1,\dots ,
4$.
%\gapp{13}

\centerline{
\vbox{\offinterlineskip
\hrule
\halign{
&\vrule#&\strut\ \ \hfil#\ \ 
&\vrule#&\strut\ \ \hfil#\ \ 
&\vrule#&\strut\ \ \hfil#\ \ 
&\vrule#&\strut\ \ \hfil#\ \ 
&\vrule#&\strut\ \ \hfil#\ \ 
&\vrule#&\strut\ \ \hfil#\ \ 
&\vrule#&\strut\ \ \hfil#\ \ \cr
height2pt
&\omit&
&\omit&&\omit&
&\omit&&\omit&
&\omit&&\omit&
&\omit&&\omit&
&\omit&&\omit&\cr
&[n]&
&${\psi}_0$ \hfil &&${\omega}_0$\hfil &
&${\psi}_1$ \hfil &&${\omega}_1$\hfil &
& ${\psi}_2$\hfil && ${\omega}_2$\hfil &
& ${\psi}_3$\hfil &&${\omega}_3$\hfil &
&${\psi}_4$\hfil  &&${\omega}_4$\hfil&\cr
height2pt 
&\omit&
&\omit&&\omit& 
&\omit&&\omit&
&\omit&&\omit&
&\omit&&\omit&
&\omit&&\omit&\cr
\noalign{\hrule}
height2pt
&\omit&
&\omit&&\omit&
&\omit&&\omit&
&\omit&&\omit&
&\omit&&\omit&
&\omit&&\omit&\cr
\noalign{\hrule}
height2pt
&\omit&
&\omit&&\omit&
&\omit&&\omit&
&\omit&&\omit&
&\omit&&\omit&
&\omit&&\omit&\cr
&    [1]  &&   $-$  &&   $+$   && 0 && 0 && 0 && 0 && 0 && 0 && 0 && 0 &\cr
&      &&      &&      &&  &&  &&  &&  &&  &&  &&  &&  &\cr
&  [2]     &&  \two--    && \two++     && \two-0 &&   \two0-  &&  \two00 &&  \two00 &&  \two00 &&  \two00 &&  \two00 &&  \two00 &\cr
&      &&      &&      &&\two0+  &&  \two+0  &&    &&    &&    &&    &&    &&    &\cr
&      &&      &&      &&  &&  &&  &&  &&  &&  &&  &&  &\cr
&  [3]    &&    \free---   &&  \free+++&&  \free--0  &&  \free0--  &&  \free-00  &&  \free00+  && \free000 && \free000 && \free000 && \free000 &\cr
&      &&      &&      &&  \free-0+  &&  \free+0-  &&  \free0+0  &&  \free0-0  &&     &&     &&     &&     &\cr
&      &&      &&      &&  \free0++  &&  \free++0  &&  \free00-  &&  \free+00  &&     &&     &&     &&     &\cr
&      &&      &&      &&  &&  &&  &&  &&  &&  &&  &&  &\cr
&   [4]   &&   \ffour----   && \ffour++++    &&  \ffour---0  &&  \ffour0---  &&  \ffour--00  &&  \ffour00++  &&  \ffour-000  &&  \ffour000-  && \ffour0000 && \ffour0000 &\cr
&      &&      &&      &&  \ffour--0+  &&  \ffour+0--  &&  \ffour-0+0  &&  \ffour0-0+  &&  \ffour0+00  &&  \ffour00+0  &&      &&      &\cr
&      &&      &&      &&  \ffour-0++  &&  \ffour++0-  &&  \ffour0++0  &&  \ffour+00+  &&  \ffour00-0  &&  \ffour0-00  &&      &&      &\cr
&      &&      &&      &&  \ffour0+++  &&  \ffour+++0  &&  \ffour-00-  &&  \ffour0--0  &&  \ffour000+  &&  \ffour+000  &&      &&      &\cr
&      &&      &&      &&      &&      &&  \ffour0+0-  &&  \ffour+0-0  &&      &&      &&      &&      &\cr
&      &&      &&      &&      &&      &&  \ffour00--  &&  \ffour++00  &&      &&      &&      &&      &\cr
height2pt
&\omit&
&\omit&&\omit&
&\omit&&\omit&
&\omit&&\omit&
&\omit&&\omit&
&\omit&&\omit&\cr}
\hrule}
}

Denote by $ C^n_k $ the binomial coefficients arising in 
Pascal's triangle. 

\medskip
\medskip
\noindent
{\bf Definition 4.4}:
An
$(m,\ n)$-$pile$
is a sequence of ($m+1$) sets $ v_j $, $j= 0,\dots ,m$ such that
\hfil
\break
\noindent
(a)  $ v_j \ \subset \ {\cal I}^n $
\hfil
\break
\noindent
(b)  each $ v_j $ has cardinality $ C^m_j $
\hfil
\break
\noindent
(c)  each $x \in v_j $ has dimension $j$.

Fix $n$, and denote the collections of source and target faces by
$$
 {\Psi }_n = \cup \  {\psi}_k [n]
$$
$$
{\Omega}_n  = \cup \ {\omega}_k [n] 
$$

\medskip
\medskip

\noindent
{\bf Proposition 4.5}:
$ {\Psi}_n $ and $ {\Omega}_n $ are $(n,\ n)$-piles.

\medskip
\medskip
\noindent
{\bf Remark 4.6}:
Each sub-$k$-cube of the $n$-cube is parallel to a subspace of 
$ R^n $. The
subspaces arising in this fashion correspond to choices of $k$ 
elements of
$\{ 1,\dots ,\ n\} $ as subscripts for coordinate axes. 
Thus there are $ C^n_k $ 
possible such subspaces, and for each, 
$ 2^{n-k} $ parallel copies
in the boundary of ${\cal I}^n $. 
We shall call the collection of all such 
%\gapp{14}
%\noindent
parallel copies a
$parallel$ $set$. Observe that the inductive definition of
source and target faces canonically chooses a representative from 
each
parallel set. 

\medskip
\medskip
\noindent
{\bf Definition 4.7}:
We will call any choice of a unique representative from each
parallel set of $k$-dimensional subcubes a 
$k$-$section$.

\medskip
\medskip

Hence Theorem 4.3 tells us that there are at least two distinct
$k$-sections whose topological union is a disc. Moreover, these 
unions are
coherent across dimensions.

Note also that differential $p$-forms in $n$-space admit a basis analogous to a $p$-section. This suggests some relationship to de Rham cohomology.

For each $x \in {\cal I}^n$, $|x| = p$, we define $ {\Psi}_x $ and 
$ {\Omega}_x$ as ($p,\ n)$-piles by using $x$ to embed elements of 
$ {\Psi}_p $
and $ {\Omega}_p $ into ${\cal I}^n$.

The definitions given above lead to the structure referred to
in the introduction. We describe the relationship between the source- and 
target-$n$-piles, and how these relate to the cocycle conditions. 

\medskip
\medskip
\noindent
{\bf Proposition 4.8}:
For each $k$, $ {\psi}_k [n]  = - {\omega}_k [n] $, where the 
antipodal
map \ $-:{\cal I}^n \longrightarrow {\cal I}^n$ 
is the involution sending a word $x$ to $-x$, obtained from $x$ by 
interchanging
each $-$ and $+$. Hence $ {\psi}_k [n] $ and $ {\omega}_k [n] $ 
are
interchanged by the antipodal map on [$n$]. 

\medskip
\medskip
\noindent
{\bf Remark 4.9}:
Since the parallel set to which a sub-cube belongs is determined by 
the
position of 0's in its word form, the antipodal map preserves 
parallel
sets.

If now we have an $(n,\ n)$-pile $ \Pi $ into which there exists an 
embedding 
$ {\Psi}_x $ for some $x \in {\cal I}^n$, we obtain a new list 
by 
replacing $ {\Psi}_x $ by $ {\Omega}_x $. Note that if $|x| = j$, 
$x \in {\psi}_j [n] $. We have a map $ {\pi}_x $, whose
source and target are ($n,\ n$)-piles. The top dimensional 
%\gapp{15}
%\noindent
cube 
of
$x$ (merely $x$ itself) remains in the pile, only its lower 
dimensional
faces changing. (Here we abuse notation somewhat as $ {\Psi}_x $ 
may occur in 
many different ($n,\ n$)-piles.) 

Denote by $ *_k $ the operation which applies the involution 
``$-$'' to every
element of $ v_j $, $j \leq k-1 $. The operation $ *_0 $ will 
denote
the identity. If now we begin with a pile $ {\Pi}_0 $ and an 
embedding $ {\Psi}_{x_0} \hookrightarrow {\Pi}_0$, it may be 
possible to find sequences 
$$ x_0, \ *_{i_0}  , \ x_1 , \ *_{i_1} ,\dots ,\ x_j ,\ *_{i_j} $$
such that a sequence of piles is generated from $ {\Pi}_0 $, 
each obtained from the previous by applying consecutively the 
$ {\pi}_{x_i} $ 
and $ *_{i_k} $. Observing that $ {\pi}_x $ is parallel set 
preserving,
%we have
it is not difficult to prove

\medskip
\medskip
\noindent
{\bf Proposition 4.10}:
If $ {\Pi}_0$ is a pile for which each $ v_k $ is a $k$-section 
whose topological
union, identifying along common boundaries, is a disk with 
boundary invariant under the antipodal
map, the same is true for $ *_j ({\Pi}_0 $)
or $ {\pi}_x ({\Pi}_0 $) whenever these are defined.
Furthermore, the boundary of each such disc remains the same.
\medskip
\medskip

The
$(n-1)$-$source$ of [n], 
denoted $ {\sigma}_{n-1} [n]$, will be such a sequence containing 
$ {\psi}_{n-1} [n] $, such that the composite of the sequence has 
source
$ {\Psi}_n $ and target $ {\cal T}_n  =  {\Omega}_n $ but with 
$ {\psi}_{n-1} $
replacing $ {\omega}_{n-1} $. By applying the involution ``$-$'' 
to the sequence
written in reverse order we will find another such sequence, which 
we call
the
$(n-1)$-$target$ $ {\tau}_{n-1} [n] $
of $[n]$. This has source pile $ {\Psi}_n $ with $ {\psi}_{n-1} $
replaced by $ {\omega}_{n-1}$. The sequence uses each element of 
$ {\omega}_{n-1} [n]$, and has target $ {\Omega}_n $.

\medskip
\medskip
\noindent
{\bf Definition 4.11}:
The formal equality $ {\sigma}_{n-1} [n]  = {\tau}_{n-1} [n] $ is 
called
the
($n-1)$-$cocycle$ $condition$.

%\gapp{16}

%\pagebreak

\section {Modifying piles - Blocks} 

Intuitively a pile describes a nested set of disks $ V_k $ of 
descending dimension, 
geometrically embedded in the $n$-cube. 
$ V_{n-2}$ is the ($n-2$)-source in 
some stage of (cubical) isotopy or deformation across the 
interior of 
$ V_{n-1} $, keeping the boundary fixed. The disc $ V_{n-1}$ in
practice is either the ($n-1$)-source or target of $ {\cal I}^n $. 
In turn each
$ V_{k-1} $ is the ($k-1$)-source of $ {\cal I}^n $
after some deformation across $ V_k $, keeping the boundary fixed.
Each disk $ V_k $ decomposes as a set $ v_k $, the list of 
``currently encountered'' $k$-dimensional faces of 
$ {\cal I}^n $ as we deform 
across the cubes of $ v_{k+1} $. (See Figure \ref{fig:8}.) 

\begin{figure}[htbp] %  figure placement: here, top, bottom, or page
   \centering
   \includegraphics[width=2in]{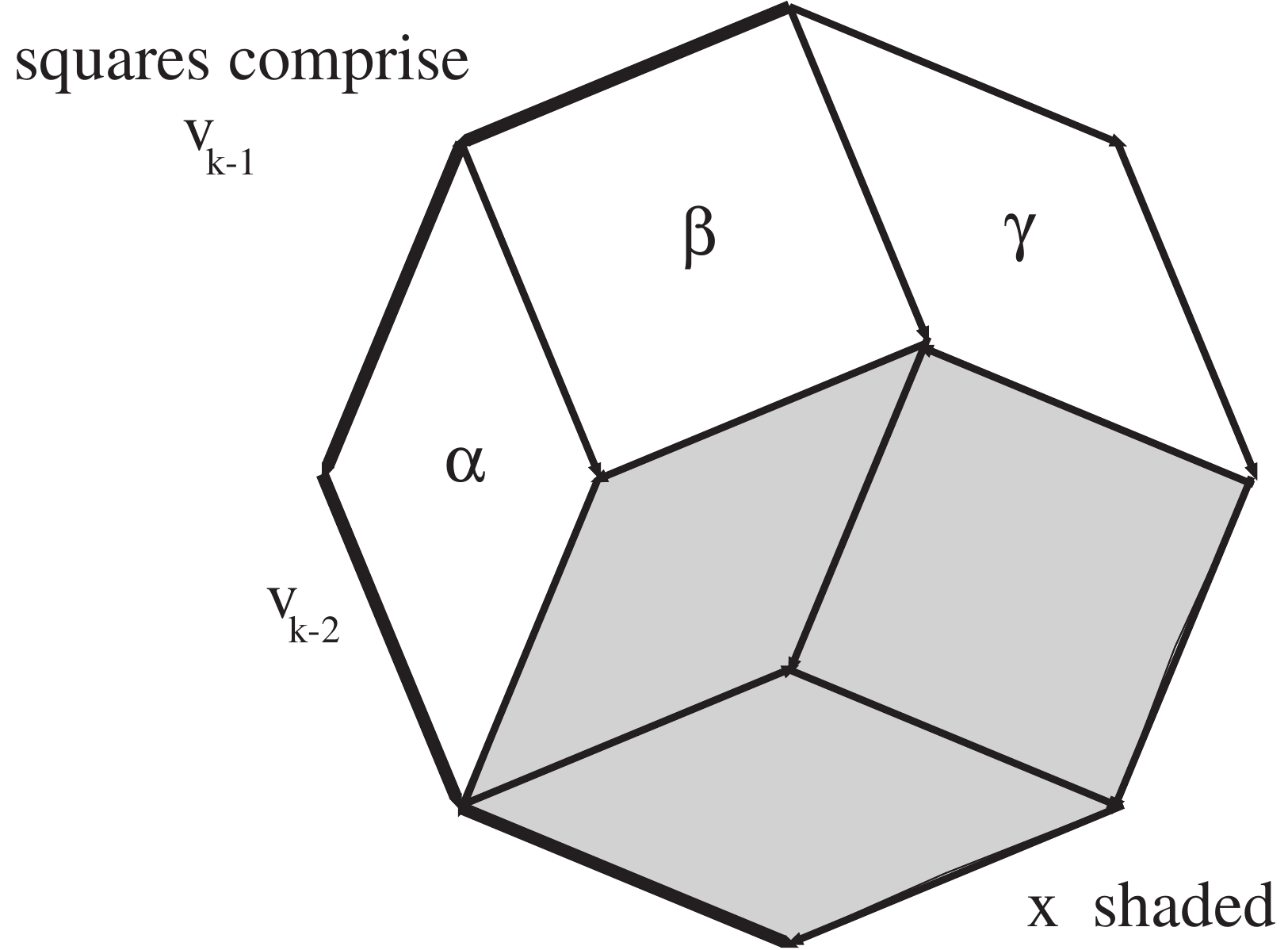} 
   \caption{The 2-cubes $\alpha,\ \beta$  can be used to move the 1-source of $[4]$ to include the  1-source of $x$.}
   \label{fig:9}
\end{figure}

Suppose at some stage we
have a pile $ {\Pi} $ such that $ v_j  =  {\psi}_j [n] $ for 
$j \leq (k-2)$,
and wish to apply $ {\pi}_x $ for $ x \in v_k $. This requires that 
the ($k-1$)-, ($k-2$)-,.., 0-faces of $x$ appear in the 
appropriate $ v_j $. 
In general this will not be the case. By looking at $ v_{k-1} $ 
we can decide
if there is any chance of applying $ {\pi}_x $. If $x$ has 
$(k-1)$ source faces appearing in $ v_{k-1} $,
it should be possible to
use the remaining $(k-1)$-cubes in $ v_{k-1}$ to
change the lower $ v_j $'s until the rest of the source faces of 
$x$
appear. We illustrate this schematically in Figure \ref{fig:9}. 
The $(k-1)$-source
cubes of $x$ are those shaded among the cubes of $ v_{k-1} $.

\begin{figure}[htbp] %  figure placement: here, top, bottom, or page
   \centering
   \includegraphics[width=5.5in]{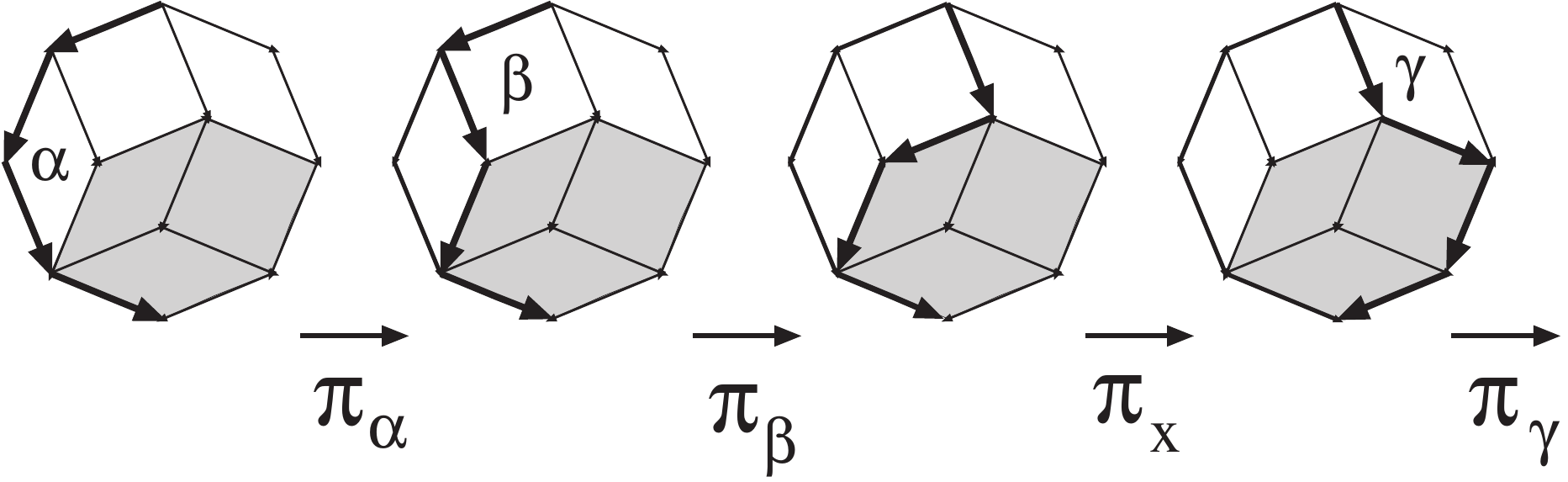} 
   \caption{Lower-dimensional completion of a cell deformation: using 2-cubes to move the 1-source to include the 1-source of a higher $3$-cube $x$, so that $\pi_x$ can be applied.}
   \label{fig:10}
\end{figure}

The heavy path
depicts $ v_{k-2} $, from which it is clear the ($k-2)$-source 
faces of
$x$ do not all appear. Hence we apply consecutively 
$ {\pi}_{\alpha} $, 
$ {\pi}_{\beta} $, then $ {\pi}_x $, and finally 
$ {\pi}_{\gamma} $ to obtain a pile
with 
$ v_{ k-2 }  =   t{\psi}_{k-2} $ replaced by $ {\omega}_{k-2} $. If
we now apply $ *_{ k-1} $, having used all of the faces of 
$ v_{k-1} $ 
in the sequence of applications, we can replace ${\omega}_{k-2} $ 
by 
$ {\psi}_{k-2 } $ and try to find the next cube 
$ x' $ higher up the
list to apply. We illustrate in Figure \ref{fig:10}.

%\gapp{17}

\medskip
\medskip
\noindent
{\bf Remark 5.1}:
There are two things to note here. The first is that to apply 
$ {\pi}_{\alpha} $
all of $ {\alpha} $'s lower dimensional source faces must appear in 
the pile. We 
may first have to apply the same process as above, at a lower 
dimension.
Secondly we have the ``moral'' of the heuristic description: 
to apply 
$ {\pi}_ x $, for $ x \in v_k $, we must work in 
conjunction 
with the remaining $ {\pi}_{\alpha} $s of $ v_{k-1 } $. Hence 
$x$ occurs as 
part of a ``block'' of cubes, namely $x$ and the rest of the
cubes in $ v_{k-1} $ which are not in $x$. This ``block'' is on 
the one hand 
derived from the pile, and on the other, composed of cubes which 
serve to modify the pile so that the process can continue.
The conclusion is that every time we apply $ {\pi}_x $, $x
\in  v_k $,
there is an associated ``sweep'' of the ($k-2$)-source of
$ {\cal I}^n $ across to the ($k-2$)-target,
at some stage sweeping across $x$.
All of these ``sweeps'' occur at different dimensions simultaneously.
%\gapp{18}

We can try to find an order in which to proceed directly from the 
pile 
$ {\psi}_n $, in order to use all of the elements of 
$ {\psi}_{ n-1} $ to change each lower $ {\psi}_k $ into the 
corresponding $ {\omega}_k $. This order will 
not be unique. However,
knowing an order from the previous dimension will enable
us to make inductively a natural choice.

\begin{figure}[htbp] %  figure placement: here, top, bottom, or page
   \centering
   \includegraphics[width=5in]{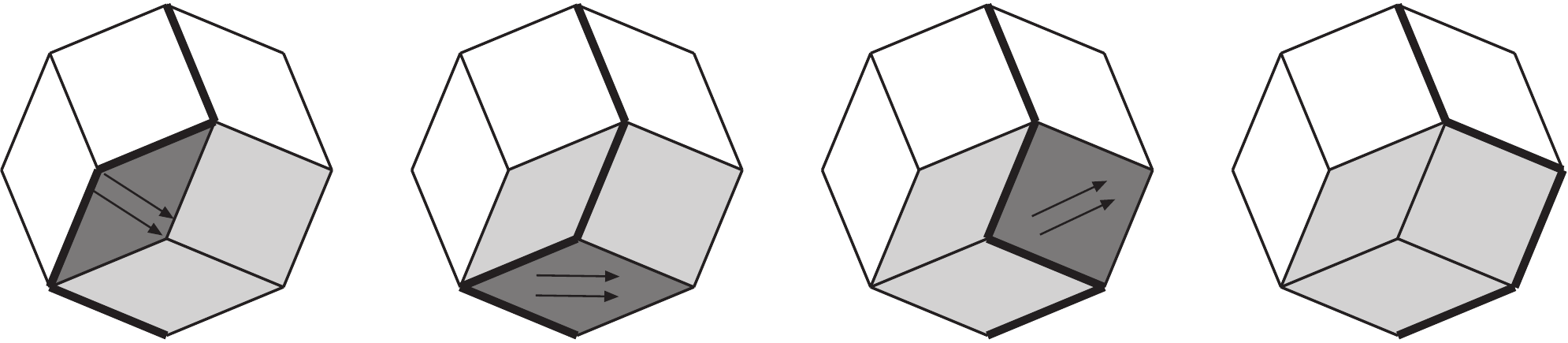} 
   \caption{Decomposition of the action of a $k$-cell into the action of ($k-1$)-cells}
   \label{fig:11}
\end{figure}

Note that the effect of $ {\pi}_x $, $|x|=k$, on 
$ v_j ,\ j \leq k $) can be decomposed : $ {\pi}_x $ alters 
$ C_k^j $ elements, an effect which can be obtained in $ C_k^{k-1} $
stages corresponding to the effect of each ($k-1$)-source-face of 
$x$ applied
consecutively (and similarly target faces). Each of these can be 
further 
decomposed, all the way down to
dimension $j+1$. In Figure \ref{fig:11} we break up the effect of 
$ {\pi}_x $ of the 
previous figure.

To present the procedure by which a canonical order 
of pile modifications is determined by the cocycle conditions,
we seek an order for the application of elements of a ``block''. 

\medskip
\medskip
\noindent
{\bf Definition 5.2}:
For $k$ odd, an $ (n,\  k)$-$block $ is a column of 
($n$, $k-1$)-blocks. For 
$k$ even, an $(n,\ k)$-$block$ is a row of ($n$, $k-1$)-blocks.  
An $(n,\  0)$-$block$
is an element of ${\cal I}^n$. 

\medskip
\medskip
\noindent
{\bf Note 5.3}:
We do not require that the number of blocks occurring at any stage 
has any
%\gapp{19}
%\noindent
specific value, nor that they have equal size.

\medskip
\medskip
\noindent
{\bf Example 5.4}:
Below we give examples of a $k$-block for $k = 0,\ 1,\ 2,\ 3$. 
Note that vertical
and horizontal bars serve well to divide up the blocks. 
A block is ``read''
from left to right, top to bottom, and from smaller to larger 
sub-blocks. These give the $(k-1)$-sources for the $k=1$-, 2- ,3- and 4-cubes respectively:

\bigskip

%\centering{
\begin{tabular}{ccccccccc }
 $-$&&
 \begin{tabular}{cc}
 \two-0\\
 \two0+\\
 \end{tabular}
 & 
 {\begin{tabular}{c|c|c}
 \three-00\  & \three-0-\  & \three00-\\
  \three0++\  & \three0+0\  & \three++0\\
 \end{tabular}
 }
  & 
 {\begin{tabular}{c|c|c|c}
  \four-000   & \four-0--       & \four-0--  &\four00--\\
  \four0+++    & \four-+0-    & \four0+0-  &\four++0- \\
                      & \four0++0  & \four+++0  &\four+++0 \\
  \hline\hline
   \four---0  & \four-0-0   & \four-0--&\four00--\\
    \four-00+    & \four-+0+  & \four0+00 &\four++0- \\
  \four0+++  & \four0+++&  &\four+++0\\
  \hline\hline
   \four---0  & \four---0 & \four00-0&\four0---\\
    \four-00+    & \four-0-+  & \four++0+ &\four+0-- \\
  \four0+++ & \four0+0+  & &\four++00\\
  \hline\hline
   \four---0  & \four0--0 & \four0---&\four0---\\
    \four000+    & \four+0-+  & \four+0-0 &\four+0--\\
                    & \four++0+  & \four++0+&\four++00\\
 \end{tabular}
 }
 \\
 \\
$k-1=0$  &&  $k-1=1$  & $k-1=2$ &$k-1=3$
\end{tabular}
\bigskip
\bigskip

%\gapp{20}

 %\pagebreak
 
Blocks can be written in a unique linear form, where when two
$k$-sub-blocks are separated by a bar, we write $*_k $.
We obtain, for each of the `sources' just given in the previous table,

\begin{enumerate}
\item %$k=1$
%\subitem 
$-$
\item %$k=2$
%\subitem  
\two-0 $*_0 $ \two0+
\item %$k=3$
%\subitem 
\three-00 $*_0$   \three0++  $*_1$   \three-0-  $*_0$   \three0+0 
 $*_1$   \three00-  $*_0 $  \three++0
\item % $k=4$:
%\subitem 
 \hskip 10pt \ffour-000  $*_0$    \ffour0+++    $*_1$    \ffour-0--    $*_0$   
 \ffour-+0-    $*_0$    \ffour0++0    $*_1$    \ffour-0--    $*_0$    
 \ffour0+0-   $*_0$    \ffour+++0    $*_1$    \ffour00--    $*_0$    
 \ffour++0-    $*_0$    \ffour+++0   $*_2$    \ffour---0    $*_0$   
  \ffour-00+   $*_0$   \ffour0+++    $*_1$    \ffour-0-0    $*_0$    
  \ffour-+0+    $*_0$    \ffour0+++    $*_1$    \ffour-0--    $*_0$   
   \ffour0+00    $*_1$    \ffour00--   $*_0$    \ffour++0-    $*_0$   
    \ffour+++0   $*_2$    \ffour---0     $*_0$    \ffour-00+    $*_0$    
    \ffour0+++    $*_1$    \ffour---0    $*_0$    \ffour-0-+    $*_0$    
    \ffour0+0+    $*_1$    \ffour00-0    $*_0$    \ffour++0+    $*_1$   
     \ffour0---    $*_0$    \ffour+0--    $*_0$    \ffour++00   $*_2$    
     \ffour---0    $*_0$    \ffour000+    $*_1$    \ffour0--0    $*_0$   
  \ffour+0-+    $*_0$    \ffour++0+    $*_1$    \ffour0---   $*_0$  
 \ffour+0-0    $*_0$    \ffour++0+    $*_1$    \ffour0---    $*_0$   
  \ffour+0--    $*_0$    \ffour++00 \qquad\qquad\ \ \ \  \hfill
\end{enumerate}

\bigskip

\medskip
\medskip
\noindent
{\bf Remark 5.5:}
An equivalent definition of a $k$-block on a set $S$ can be given. 
A $k$-block is a pair $(f, g)$ where $f: \underline n\to S$  and 
$g: \underline{n-1} \to \underline k  $, for some $n$. Note that for $S$ an
ordered set, we obtain something very much like a span.

\medskip
\medskip
\noindent
{\bf Definition 5.6:}
By a $j$-sub-block of a $k$-block we will mean any sub-block characterized
by
elements between a pair $ *_r , \  *_s  $, for  $r, s\geq j-1$, such that no $ *_l $
occurs between them for $l\geq j$. A full $j$-sub-block will be a
$j$-sub-block maximal with respect to this condition.
 
For each word $x$ above, replace $x$ in the string by $\pi_x $. Interpret 
each $ *_i $ as ``rewrite $ \omega_j [n] $ as $ {\psi}_j [n] $ for each $j\leq (i-1)$'', where $n$ is the word length. 
 
For (c) we obtain the sequence of piles in the following Table 2.
The 12 columns correspond to the geometric figures in Figure \ref{fig:12}. Those on the
left and right are 
%\gapp{21}
%\noindent
%\end{document}
those of $ {\Psi}_3 $ and $ \Omega_3 $. The elements
of each $ v_k $ are shown in bold, and the hatched element is the
appropriate $x$ for which $ \pi_x $ is about to be applied. This
illustrates the ``nestings'' of disks referred to earlier.

\medskip
\medskip

\noindent
{\bf
Proposition 5.7:}
The 0, 1-, 2-, 3-cocycle conditions are given by (1)--(4) above.

%\pagebreak 
\medskip
\medskip

{\small
\hskip 6pt \begin{tabular}{cccccccccccccccccccccccc}
\hline
\hskip 2pt& {\three-00}& & \hskip 7pt $*_0$& &\hskip 0pt \three0++& &\hskip -2pt  $*_1$& &\hskip  -4pt \three-0-& & \hskip -3pt $*_0$& &\hskip -5pt \three0+0& &\hskip -3pt $*_1$& &\hskip -3pt \three00-& &\hskip -3pt  $*_0$& &\hskip -3pt \three++0&&\\
\hline
\hline
 \end{tabular}
 }
 
% \vskip -6pt
%\leftline
{
  \begin{tabular}{cccccccccccc}
%\hline
%\hline
\shree---&\shree-++&\shree-++&\shree+++&\shree---&\shree-+-&\shree-+-&\shree+++&\shree---&\shree++-&\shree++-&\shree+++
\\
\\
\shree--0&\shree-+0&\shree-+0&\shree-+0&\shree-+0&\shree-+0&\shree-+0&\shree++0&\shree++0&\shree++0&\shree++0&\shree++0
\\
\shree-0+&\shree-0-&\shree-0-&\shree-0-&\shree-0-&\shree-0-&\shree-0-&\shree-0-&\shree-0-&\shree+0-&\shree+0-&\shree+0-
\\
\shree0++&\shree0++&\shree0++&\shree0++&\shree0++&\shree0++&\shree0++&\shree0+-&\shree0+-&\shree0--&\shree0--&\shree0--
\\
\\
\shree-00&\shree-00&\shree-00&\shree-00&\shree-00&\shree-00&\shree-00&\shree-00&\shree-00&\shree-00&\shree-00&\shree-00
\\
\shree0+0&\shree0+0&\shree0+0&\shree0+0&\shree0+0&\shree0+0&\shree0+0&\shree0+0&\shree0+0&\shree0+0&\shree0+0&\shree0+0
\\
\shree00-&\shree00-&\shree00-&\shree00-&\shree00-&\shree00-&\shree00-&\shree00-&\shree00-&\shree00-&\shree00-&\shree00-
\\
\\
\shree000&\shree000&\shree000&\shree000&\shree000&\shree000&\shree000&\shree000&\shree000&\shree000&\shree000&\shree000
\\
\hline
 \end{tabular}
  
 \medskip
\centerline{\bf Table 2 - the 2-cocycle source of the 3-cube}

 \medskip
  \medskip
  
We can tabulate part of the effect of the block (d) on $ {\psi}_4 $ by 
recording the changes at the 2-dimensional level. (We can always ignore the 
changes lower down in the pile to save space in description.) 
Since the highest dimension of $x$ occurring 
for any application of $ \pi_x $ is 3, the set $ {\psi}_3 $ remains unaltered.
The reader may wish to refer forward to the ``octagon of octagons'' depicted in
Figure \ref{fig:17}, comparing the lists of 2-cubes above with the left hand side 
configurations in the figure.
%\gapp{22}

$$
\begin{array}{cccccccc}
&
 \includegraphics[width=0.7in]{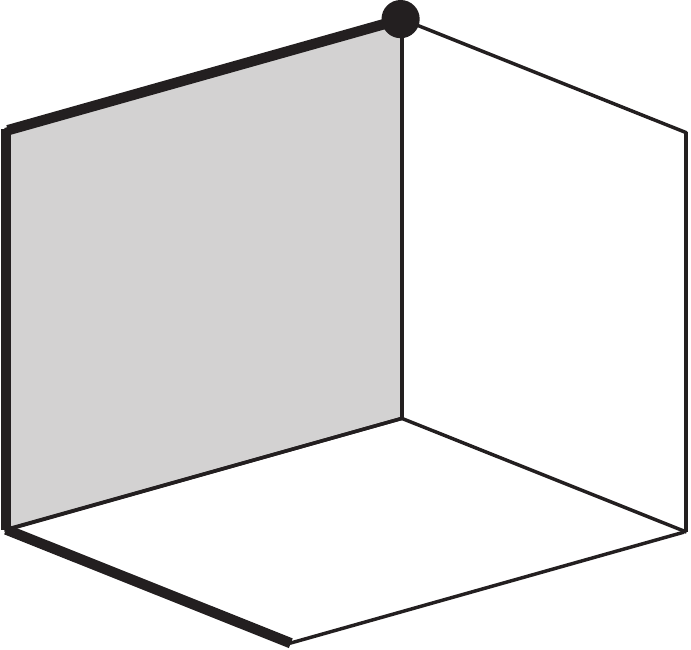} 
 & \begin{array}{c}\longrightarrow \\ \pi_{-00} \end{array}
&
 \includegraphics[width=0.7in]{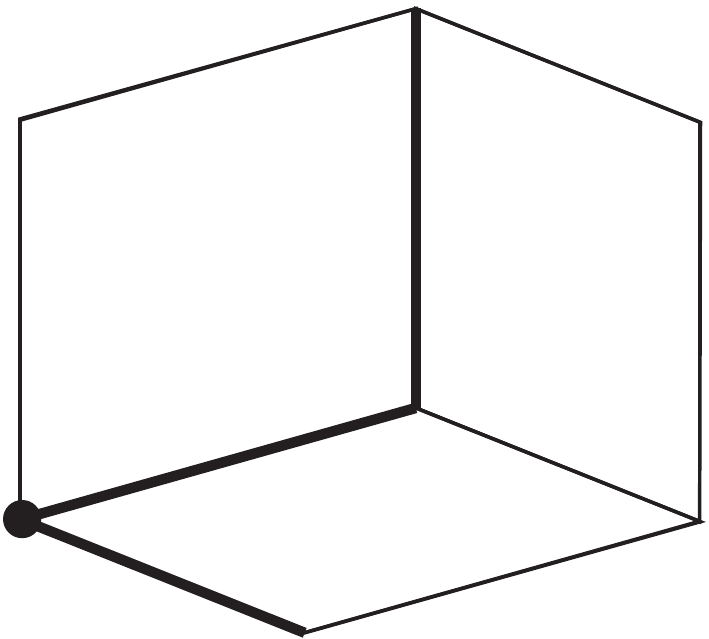} 
 & \begin{array}{c}\longrightarrow \\ {*_0} \end{array}
&
 \includegraphics[width=0.7in]{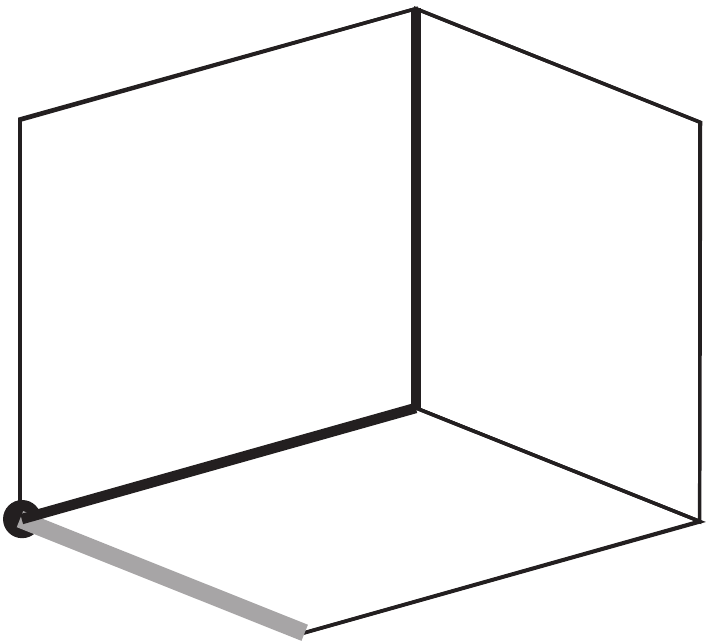} 
 & \begin{array}{c}\longrightarrow \\ \pi_{0++} \end{array}
&
 \includegraphics[width=0.7in]{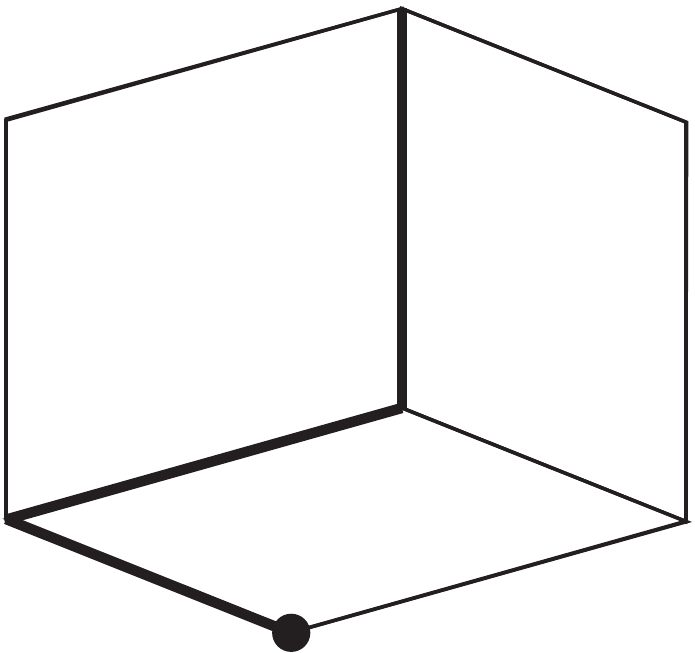}
  \\
  %$$\vskip -0.5in$$
 %&&\pi_{-00}&&*_0&&\pi_{0++}\\
  \begin{array}{c}\longrightarrow \\ {{*_1}}\end{array}&
 \includegraphics[width=0.7in]{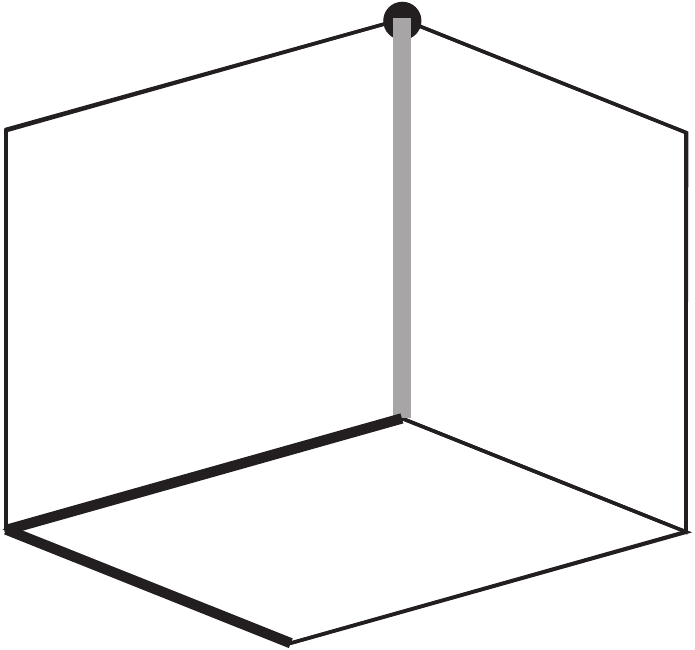} 
 & \begin{array}{c}\longrightarrow \\ \pi_{-0-}\end{array}&
 \includegraphics[width=0.7in]{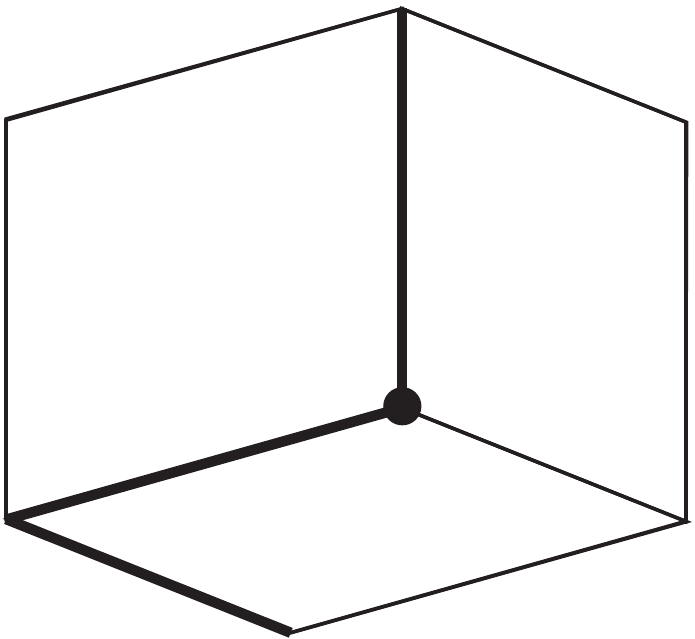} 
 & \begin{array}{c}\longrightarrow \\ {*_0}\end{array}&
 \includegraphics[width=0.7in]{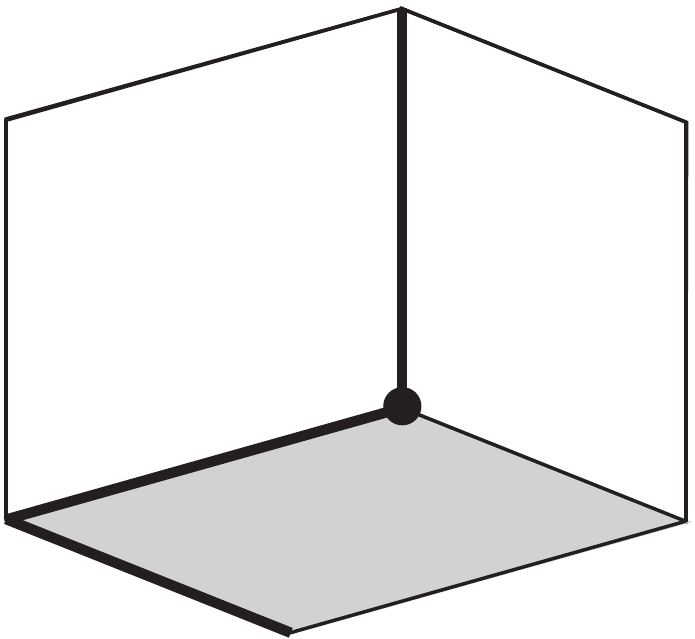} 
 & \begin{array}{c}\longrightarrow \\ \pi_{0+0} \end{array}&
 \includegraphics[width=0.7in]{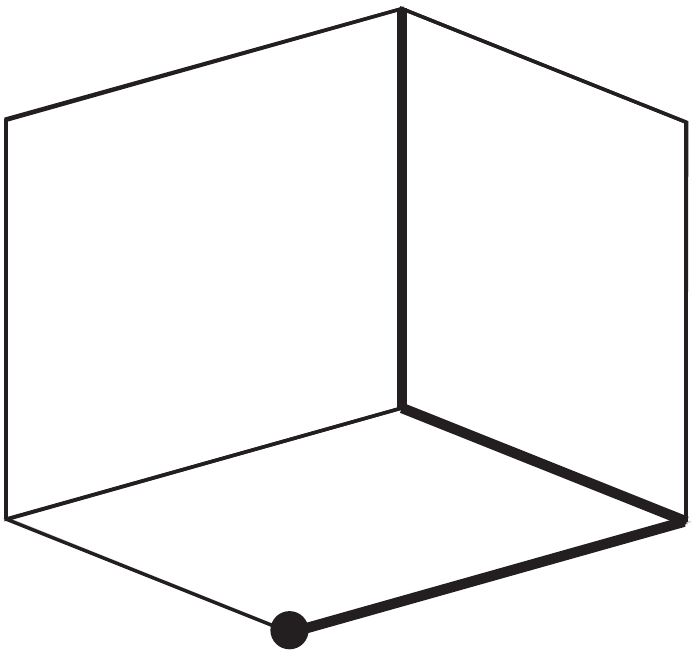}
  \\
%{*_1}&&\pi_{-0-}&&*_0&&\pi_{0+0}\\
  \begin{array}{c}\longrightarrow \\ {*_1}\end{array}&
 \includegraphics[width=0.7in]{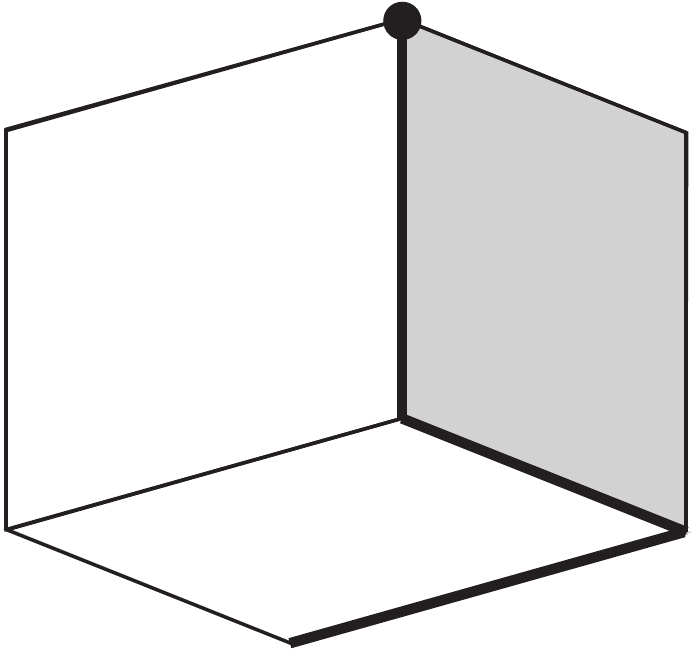} 
 & \begin{array}{c}\longrightarrow \\ \pi_{00-}\end{array}&
 \includegraphics[width=0.7in]{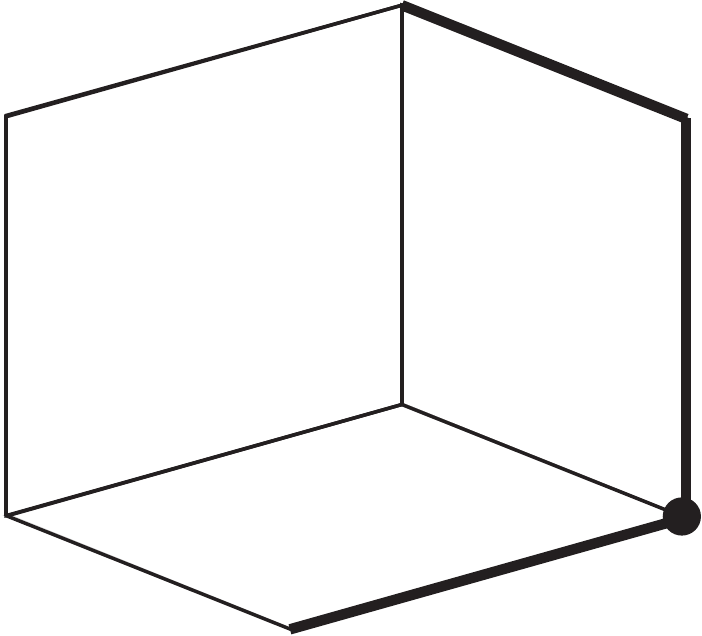} 
 & \begin{array}{c}\longrightarrow \\ {*_0} \end{array}&
 \includegraphics[width=0.7in]{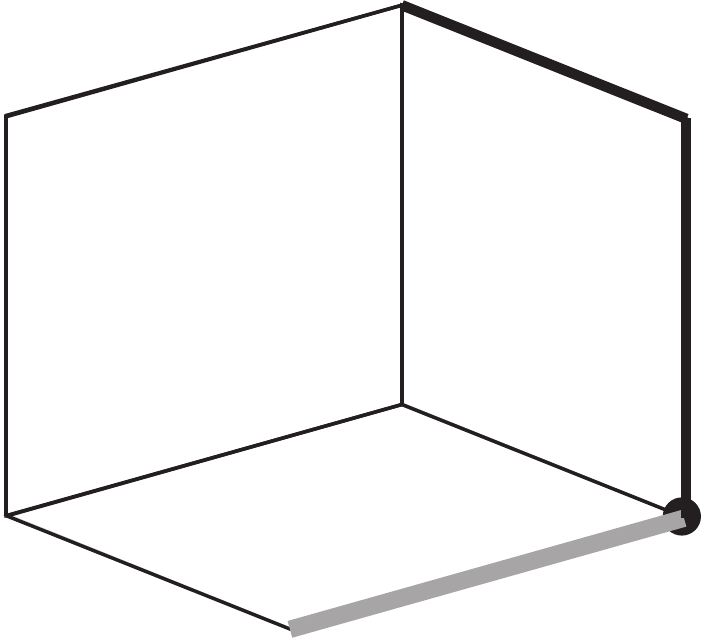} 
 & \begin{array}{c}\longrightarrow \\ \pi_{++0} \end{array}&
 \includegraphics[width=0.7in]{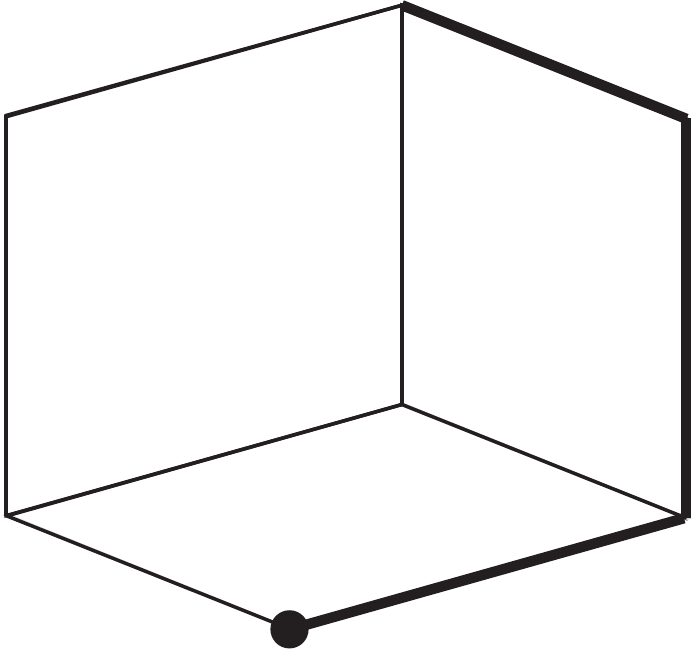}
 \\
 % \vskip -0.5in
%{*_1}&&\pi_{00-}&&*_0&&\pi_{++0}\\
\end{array}
$$
\begin{figure}[htbp] %  figure placement: here, top, bottom, or page
   \centering
   \caption{Rewriting: the action of a 2-block on unions of lower-dimensional subdiscs}
   \label{fig:12}
\end{figure}
\pagebreak
\medskip
\medskip

  \begin{tabular}{ccccccccc}
\hline
&\four-000&&\four0+00&&\four00-0&&\four000+&\\
\hline
\hline
\four--00&&\four-+00&&\four++00&&\four++00&&\four++00\\
\four-0+0&&\four-0-0&&\four-0-0&&\four+0-0&&\four+0-0\\
\four0++0&&\four0++0&&\four0+-0&&\four0--0&&\four0--0\\
\four-00-&&\four-00+&&\four-00+&&\four-00+&&\four+00+\\
\four0+0-&&\four0+0-&&\four0+0+&&\four0+0+&&\four0-0+\\
\four00--&&\four00--&&\four00--&&\four00-+&&\four00++\\
\hline
\end{tabular}

  \medskip
\centerline{\bf Table 3 -- the source  3-cube action on source  2-cubes of the 4-cube}

 \medskip
\medskip
\medskip
%\gapp{23}

We denote by $ B_n^k $ the set of $k$-blocks of 
elements of ${\cal I}^n$. Note that there are natural inclusions  
$ B_n^k  \to B_n^{ k+1} $,
and an isomorphism  ${\cal I}^n \cong  B_n^0 $. The $ B_n^k $ can be used to define a structure analogous to 
the $ N_\omega $ occurring in Street \cite{St}. Note
that there is a ``forgetful'' map which sends a $k$-block $ {\alpha} $ to its
underlying set of elements of ${\cal I}^n$. We can thus define the 
{\it 
dimension
}
of $ {\alpha} $ to be the maximum of the dimensions of its elements.

\medskip
\medskip
\noindent
{\bf 
Definition 5.8:}
A
{\it 
distinguished
}
$k$-block is a $k$-block for which each sub-$j$-block has all but one
sub-($j$-1)-blocks of dimension $j$, the other ($j$-1)-block of dimension $ \geq j$. We will denote by $ D_n^k \subset   B_n^k $ the subset of
distinguished blocks. A 0-block is always considered distinguished.

\medskip
\medskip
\noindent
{\bf 
Examples 5.9:}
All of the examples given above are distinguished. Our inductive definition
of cocycles will demonstrate that all cocycles are distinguished. The
notion makes precise how an element $x$ of a pile can only be applied in
conjunction with those elements of a pile of the next lower dimension which
do not occur as source faces of $x$. Thus $x$ must be ``stretched out''
appropriately to the source and target of the whole pile.

Writing out the terms using the $ *_i $ provides not only the order in
which to modify the ($n,\ n)$-piles, but also tells us at which stage we
must rewrite the pile from below the $i$-dimensional level. Moreover, a 
$ *_k $ separating two blocks $ {\alpha}_i $ and $ {\alpha}_{ i+1}$ indicates
that geometrically the respective collections of cubes are ``glued together'' 
along a $k$-dimensional disk $D$.
$D$ itself is a collection of subcubes, and has a pile/block description as
the the $k$-target of $ {\alpha}_i $ and the $k$-source of $ {\alpha}_{ i+1}$. 

%\gapp{24}

\section{Sources and targets - inductive cocycles}
 
We use the order defining the ($n-2$)-source of $ {\cal I}^{n-1} $ to obtain
the ($n-1$)-source of $ {\cal I}^n $. Since $ {\cal I}^n $ is obtained
geometrically by thickening the $(n-1$)-cube, all of the cubes in the
expression of the ($n-2$)-source of $ {\cal I}^{n-1} $ will play a role in
the form of either old, new or thickened copies, via the sources and
targets they themselves have in lower dimensions.
All $k$-blocks subsequently considered will be distinguished, unless
otherwise stated. 
Our intention is to describe how to complete the table for all $n$, $ {\sigma}_k $:

\medskip
\medskip

%\leftline
{\small %\footnotesize
%\centering{

\begin{tabular}{c||c|c|c|c }
& $\sigma_0[n]$ &   $\sigma_1[n]$ & $\sigma_2[n]$ & $\sigma_3[n]$   \\
\hline
\hline
 [1]&$-$ &0&0&0 \\
 \hline
 [2]&\ftwo-- &
  \begin{tabular}{cc}
 \ftwo-0\\
\ftwo0+\\
  \end{tabular} 
  &\ftwo00
  &\ftwo00\\
  \hline
 [3]&\free--- &
  \begin{tabular}{cc}
  \free--0\\
  \free-0+\\
  \free0++
  \end{tabular} 
  &
 {\begin{tabular}{c|c|c}
\free-00\  & \free-0-\  & \free00-\ \\
  \free0++\  & \free0+0\  & \free++0\ \\
 \end{tabular}
 }
  & 
  \free000
  \\
  \hline
 [4]&\ffour----&
 \begin{tabular}{cc}
 \ffour---0\\
 \ffour--0+\\
  \ffour-0++\\
   \ffour0+++\\
 \end{tabular}
 & 
 {\begin{tabular}{c|c|c|c|c|c}
\ffour--00   & \ffour--0-   & \ffour--0- &\ffour-00-  &\ffour-0--  & \ffour00--\\
\ffour-0++  & \ffour-0+0 &\ffour-0+- &\ffour0++- &\ffour0+0- & \ffour++0-\\
\ffour0+++ & \ffour0+++ &\ffour0++0 &\ffour+++0 &\ffour+++0 & \ffour+++0\\
 \end{tabular}
 }
  & 
 {\begin{tabular}{c|c|c|cc}
  \ffour-000   & \ffour-0--       & \ffour-0--  &\ffour00--&\\
  \ffour0+++    & \ffour-+0-    & \ffour0+0-  &\ffour++0- \\
                      & \ffour0++0  & \ffour+++0  &\ffour+++0 \\
  \hline
  \hline
   \ffour---0  & \ffour-0-0   & \ffour-0--&\ffour00--\\
    \ffour-00+    & \ffour-+0+  & \ffour0+00 &\ffour++0- \\
  \ffour0+++  & \ffour0+++&  &\ffour+++0\\
  \hline
  \hline
  \ffour---0  & \ffour---0 & \ffour00-0&\ffour0---\\
    \ffour-00+    & \ffour-0-+  & \ffour++0+ &\ffour+0-- \\
  \ffour0+++ & \ffour0+0+  & &\ffour++00\\
  \hline
  \hline
   \ffour---0  & \ffour0--0 & \ffour0---&\ffour0---\\
    \ffour000+    & \ffour+0-+  & \ffour+0-0 &\ffour+0--\\
                    & \ffour++0+  & \ffour++0+&\ffour++00\\
 \end{tabular}
 }
 \\
 \hline
\end{tabular}
}
 }
 \medskip
\centerline{\bf Table 3 -- $\sigma_k[n]$ for small $k,n$}

 \medskip
  \medskip

We shall only need to consider ${\sigma}_k[n] $, $ k \  < \  n $. It is easy to
list the sub-$k$-cubes of $ {\cal I}^n $ which occur, but the tricky part is
%\gapp{25}
%noindent
to pin down and order those lower dimensional subcubes which stretch out
the $ {\cal I}^k $s to $ \sigma_j [n] $ and $ \tau_j [n] $ for $ j < k $.

We shall inductively define functions

$$
{\sigma}_k, \ {\tau}_k \ :  D_n^m  \longrightarrow  D_n^k  
$$
$$\ \ \  
{\lambda}_k  ,   {\nu}_k  :   D_n^k    \longrightarrow     D_{n+1}^k  
$$
$$\qquad \ \
 {\mu}_k : D_n^k \longrightarrow  D_{n+1}^{k+1}
$$

\begin{description}
\item[(a)]  
 For the $n$-cube [$n$], define $ {\sigma}_n [n] $ = [$n$] = 
$ {\tau}_n [n] $. Hence these are 0-blocks (and $k$-blocks by default.)

\item[(b)]   For 0-blocks, define maps $ {\mu}_0 ,\  {\nu}_0 \  $ and $ {\lambda} _0  :
  D_n^0 \  \to \  D_{n+1}^0 $  by $x \mapsto {\mu }x ,\ {\nu}x ,\  \lambda x $ respectively. Extend the definition to an arbitrary $k$-block by
operating on each sub-0-block. (In practice symbolically this will partly involve adding   $-,0,+$.)

\item[(c)]   Define $ {\sigma} _0 [ n ] $ = $ {\lambda} _0 ( 
{\sigma} _0 [n-1] ) $, and $ {\tau}_0 [ n ] \  = \  
{\nu} _0 ( {\tau}_0[n-1]) $. This is the geometrically obvious
choice. Extend the definition to  $ {\cal I}^n = D_n^0 $, as described using
$*$. For an arbitrary $k$-block, define $\sigma_0 $ and $ \tau_0 $ by restricting to
the first 0-block. 
 
\item[(d)]  Define the 1-blocks $ {\sigma}_1 [n ] $ and $ {\tau}_1 
[ n ] $ inductively by  

$$
\begin{array}{ccc}
  {\sigma}_1 [n]   &  =  &   {\mu}_0 {\sigma}_0 [n-1] \\
 &&  {\nu}_0{\sigma}_1 [n-1] \\
\end{array}
\qquad
\qquad
\begin{array}{ccc}
 {\tau}_1 [n]     &  =  & {\lambda}_0 {\tau}_1 [n-1]    \\
 &&  {\mu}_0 {\tau}_0 [n-1] \\
\end{array}
$$

These coincide with $ {\psi}_1 [n] $ and $ \omega_1 [n] $ respectively.  
Extend the definition of $ {\sigma}_1 $ and $ {\tau}_1 $ to 
$ D_n^0 \to D_{n+1}^1 $. 
For $ {\alpha}\in D_n^1 $, set 
$$
\begin{array}{ccccc}
\sigma_1 ( {\alpha} )   &  =  & \sigma_1
\left|\begin{array}{c}
\alpha_1\\
\vdots\\
\alpha_t\\
\end{array} 
\right|&:=& 
\left|\begin{array}{c}
\sigma_1\alpha_1\\
\vdots\\
\sigma_1\alpha_t\\
\end{array} 
\right| \\
% &&  {\mu} _0 {\tau}_0{[n-1] }}\\
\end{array}
\qquad
\qquad
\begin{array}{ccccc}
\tau_1 ( {\alpha} )   &  =  &\tau_1
\left|\begin{array}{c}
\alpha_1\\
\vdots\\
\alpha_t\\
\end{array} 
\right|&:=& 
\left|\begin{array}{c}
\tau_1\alpha_1\\
\vdots\\
\tau_1\alpha_t\\
\end{array} 
\right| \\
\end{array}
$$
%\gapp{26}

This means that all of the columns are to be strung together, one on top 
of the next. For $ {\alpha} \in  D_n^k $, we define the 1-source
and 1-target by restricting to the first column. Geometrically think of the
column as a string of subcubes joined together along their 0-sources and
0-targets.
Taking the 1-source then strings together the $ {\cal I}^1 $s making up
the 1-sources of the elements of the initial column. 
Note that for an arbitrary 1-block, 0-sources and 0-targets will not match, 
but those with which we are inductively concerned do continue to enjoy this property. 

\item[(e)]   Define $ {\lambda}_1 \  and \  {\nu}_1 $: 
$ D_n^1 \  -> \  D_{n+1}^1 $ by

$$
\begin{array}{ccccc}
 \lambda_1 ( {\alpha} )  &  =  &  \lambda_1
\left|\begin{array}{c}
\alpha_1\\
\vdots\\
\alpha_t\\
\end{array} 
\right|&:=& 
\left|\begin{array}{c}
 \lambda_0\alpha_1\\
\vdots\\
 \lambda_0\alpha_t\\
\mu_0\tau_0[n]\\
\end{array} 
\right| \\
% &&  {\mu} _0 {\tau}_0{[n-1] }}\\
\end{array}
\qquad
\qquad
\begin{array}{ccccc}
\nu_1 ( {\alpha} )   &  =  &\nu_1
\left|\begin{array}{c}
\alpha_1\\
\vdots\\
\alpha_t\\
\end{array} 
\right|&:=& 
\left|\begin{array}{c}
\mu_0\sigma_0[n] \\
\nu_0\alpha_1\\
\vdots\\
\nu_0\alpha_t\\
\end{array} 
\right| \\
\end{array}
$$
 
Hence $ \sigma_1 [n]\  =\  {\nu}_1 \sigma_1 [n-1] $ and 
$ \tau_1 [n] \  = \  {\lambda}_1 \tau_1 [n-1]$. Extend the domain of 
$ {\lambda}_1 $ and $ {\nu}_1 $ to $ D_n^k $ by application to each 
sub-1-block. Note that $ {\mu} _0 \tau_0 $ and $ {\mu} _0 \sigma_0 $ ``stretch out''  the 
ends of a string of $ {\cal I}^1 $s arising from $ {\cal I}^{n-1} $,
which stretched between the 0-source and 0-target of $ {\cal I}^{n-1} $, to
the 0-source and 0-target of $ {\cal I}^n $.
 
\item[(f)]  Now define $ {\mu}_1  :  D_n^1\to D_{n+1}^2 $ by
$$
\begin{array}{ccc}
\mu_1  ( {\alpha} )  &  =  & \mu_1 
\left|\begin{array}{c}
\alpha_1\\
\vdots\\
\alpha_t\\
\end{array} 
\right|
\end{array}:= 
 \begin{array}{c|c|c|c}
 \begin{array}{c}
 \mu_0\alpha_1\\
\nu_0\sigma_1\alpha_2\\
\vdots\\
\nu_0\sigma_1\alpha_t\\
\end{array} 
&
 \begin{array}{c}
\lambda_0\tau_1\alpha_1\\
\mu_0\alpha_2\\
\vdots\\
\nu_0\sigma_1\alpha_t\\
\end{array} 
&
\cdots
& \begin{array}{c}
\lambda_0\tau_1\alpha_1\\
\lambda_0\tau_1\alpha_2\\
\vdots\\
\mu_0 \alpha_t\\
\end{array} 
% &&  {\mu} _0 {\tau}_0{[n-1] }}\\
\end{array}
$$
This is illustrated in Figure \ref{fig:13}.

%\gapp{27}
 \begin{figure}[htbp] %  figure placement: here, top, bottom, or page
    \centering
  \includegraphics[width=1.1in]{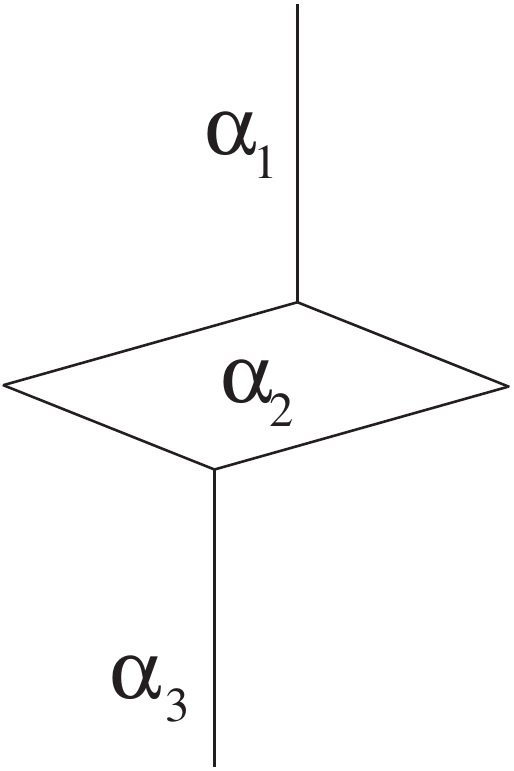} \ \ \ 
%${\longrightarrow}$
 \ \ \ \qquad \qquad \qquad \qquad
    \includegraphics[width=3in]{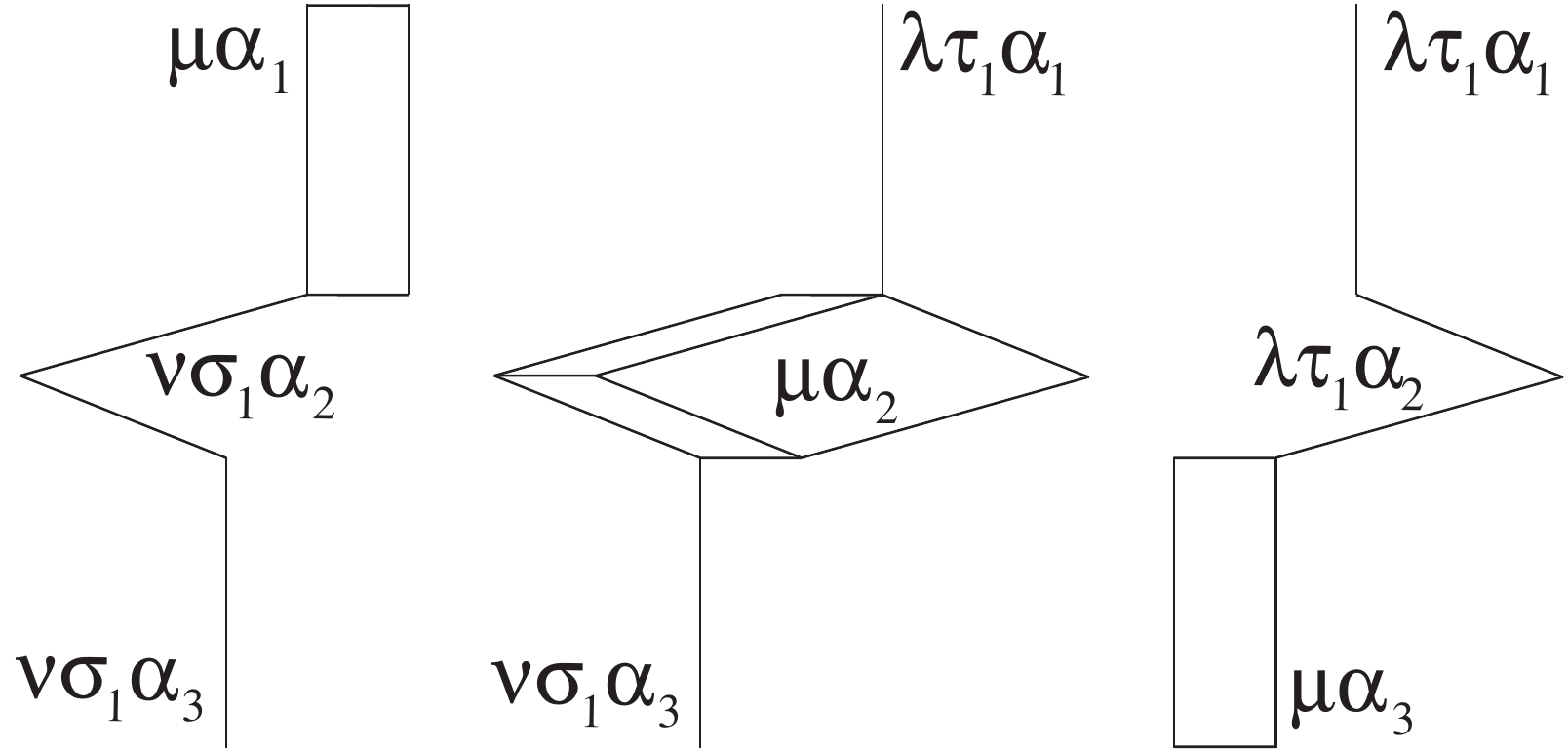}  
 %   $$\longrightarrow$$
  %\label{fig:13}
    \caption{ The source thickens}
    \label{fig:13}
 \end{figure}

\item[(g)]   To define $ \sigma_j $ and $ \tau_j $ on $ D_n^k $ we use triple
induction, first on $n$, then $k$, then $j$.
For $j$ = $n$, $k$ = 0, we have already given the
definition for ``pure'' cubes. Extension to $ D_n^0 $ is carried out
as above. 
 
All that is required is to avoid loops in our definitions. For $ \sigma_2 $ we
extend the domain to $ D_n^0 $ by consecutively defining $ \sigma_2 $[2], 
$ \tau_2 $[2], and 
$$
     {\sigma}_2 [ n ] \ \ = \ \  {\mu}_1 \sigma_1 [n-1]  \, | \, 
{\lambda}_1 ( \sigma_2 [n-1]) 
;
\qquad
\qquad
{\tau}_2 [ n ] \ \ = \ \  {\nu}_1 ( {\tau}_1
[n-1])   \, | \,   {\mu}_1 \tau_1  [n-1]  .
$$
This makes precise the way in which the 1-dimensional sources and targets
of $ {\cal I}^{n-1} $ contribute to the 2-dimensional sources of $ {\cal I}^n $.

We now wish to extend the domain to $ D_n^1 $.
It is here that the notion of distinguished k-block is crucial.
In each 1-block, every sub-0-block except possibly one, $ {\alpha}_i $, say,
has dimension 1. These sub-blocks have trivial 2-dimensional sources, and
serve to connect the 0-sources of the 1-block to the 0-source of 
$ {\alpha}_i $, and similarly for the 0-targets of $ {\alpha} $ and $ {\alpha}_i $.
By induction, $ {\alpha}_i $ has sources and targets in lower dimensions. For
%\gapp{28}
%\noindent 
$ \sigma_2 $($ {\alpha}_i )$ we obtain a 2-block, each 1-block of which stretches 
between the 0-source and 0-target of $ {\alpha}_i $. In order to define $ \sigma_2 $ 
of the 1-block $ {\alpha} $, we must stretch out each of the 1-blocks of $ \sigma_2
( {\alpha}_i $) by adding to each column a copy of the preceding 
$ {\alpha}_j $s of $ {\alpha} $ at the top, and those succeeding at the bottom. We illustrate this in Figure \ref{fig:14}. 
\begin{figure}[htbp] %  figure placement: here, top, bottom, or page
   \centering
   \includegraphics[width=1.5in]{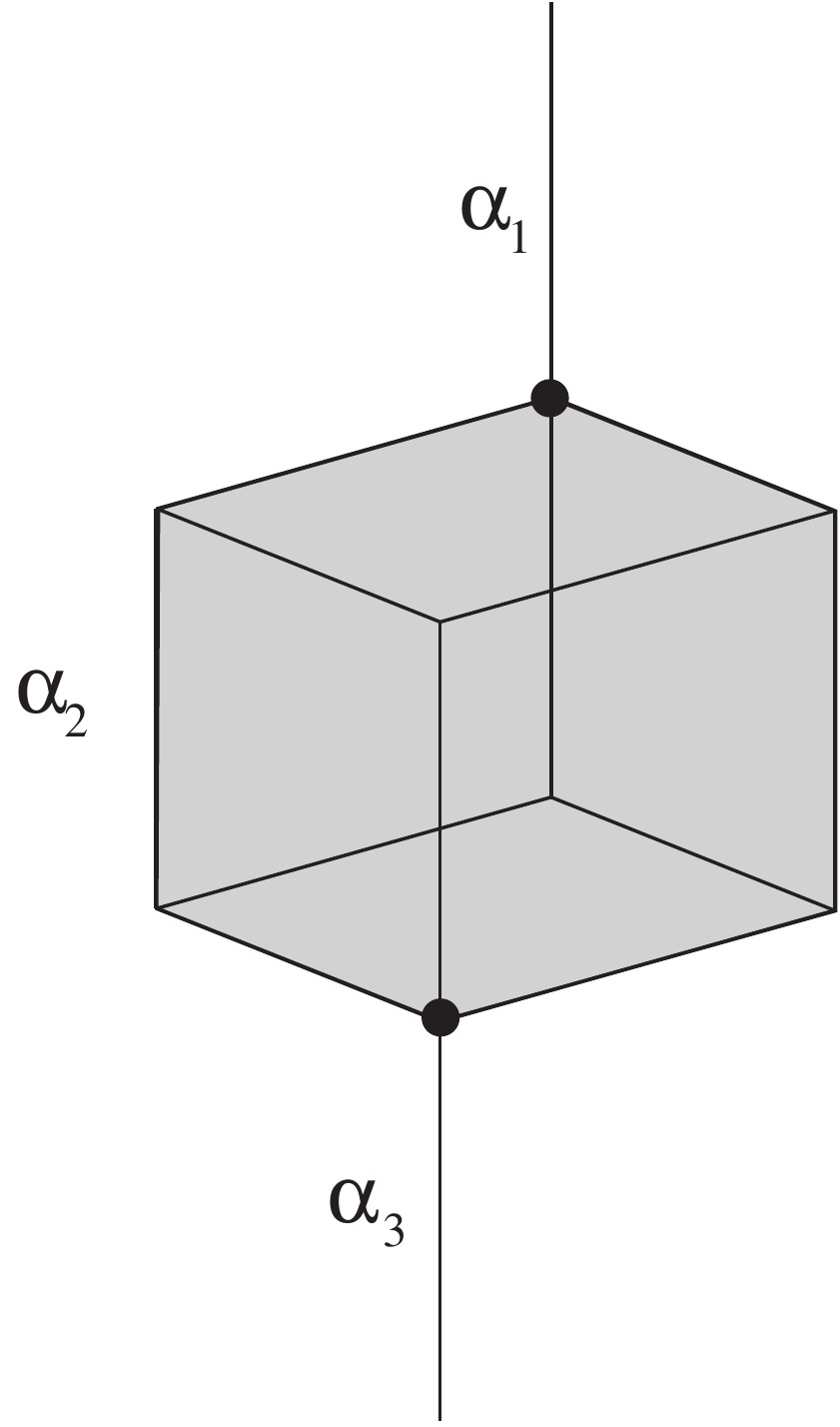} 
   \qquad\qquad 
      \includegraphics[width=3in]{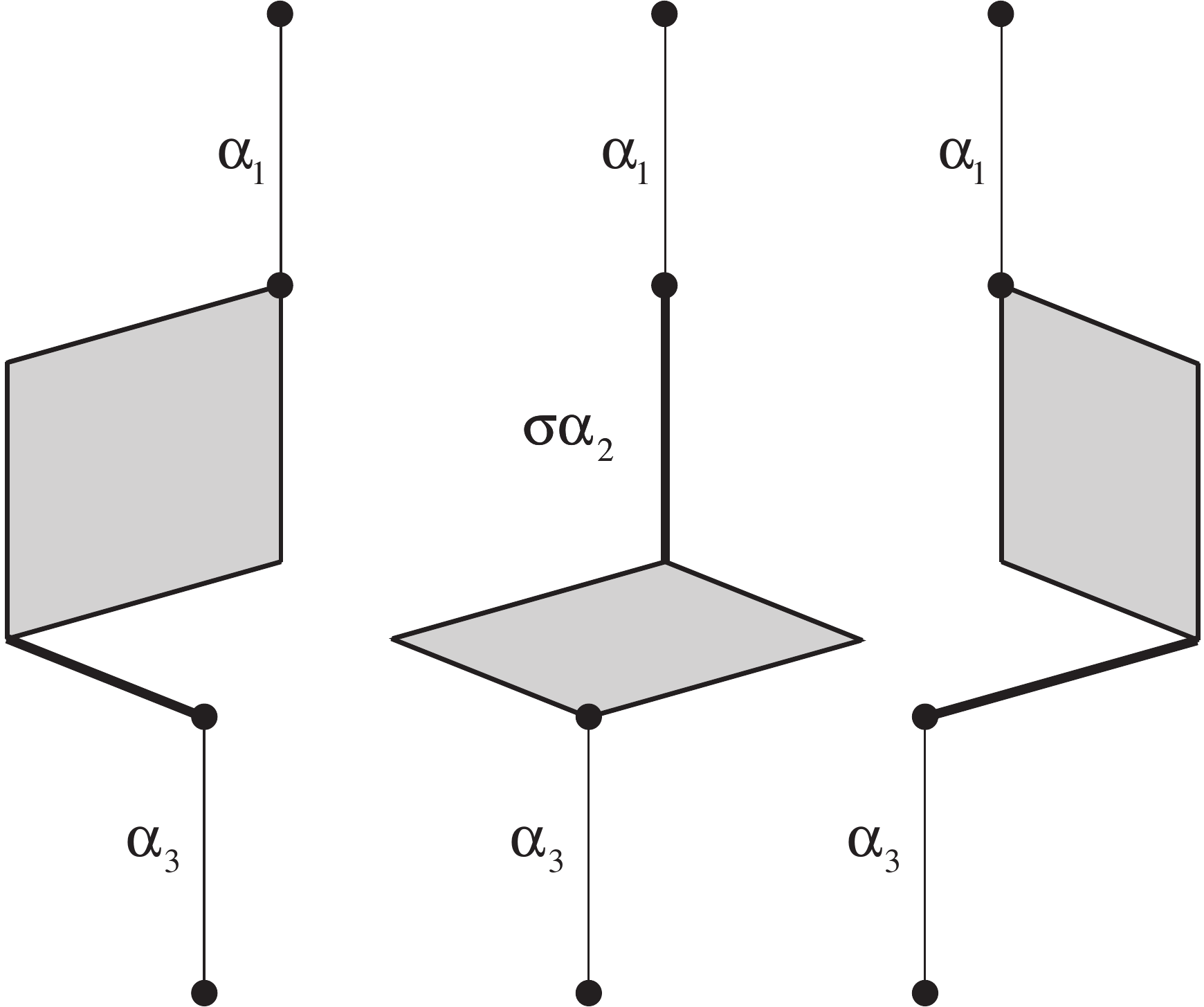} 
   \caption{ Stretching out lower dimensional sources and targets of distinguished blocks : $\alpha_2\to\sigma_2\alpha_2$. Left is the block with distinguished 3-cube; on the right, three 2-cubes are stretched out.}
   \label{fig:14}
\end{figure}

Note that the resulting 2-block is again distinguished. The same
procedure will be repeated in higher dimensions. Thus we extend the domain of 
$ \sigma_2 $ to $ D_n^1 $ by setting
$$
\begin{array}{ccccc}
 \sigma_2 ( {\alpha} )  &  =  &  \sigma_2 
\left|\begin{array}{c}
\alpha_1\\
\vdots\\
\alpha_i\\
\vdots\\
\alpha_t\\
\end{array} 
\right|&:=& 
\left|\begin{array}{c}
 \alpha_1^{(m)}\\
\vdots\\
 \sigma_2 \alpha_i\\
 \vdots\\
 \alpha_t^{(m)}\\
\end{array} 
\right| \\
% &&  {\mu} _0 {\tau}_0{[n-1] }}\\
\end{array}
$$
Here we have written the superscript ``$m$'' to indicate multiple copies will
be taken, since $ \sigma_2 ( {\alpha}_i )$ is a 2-block.
This is 2-category pasting again.

To define $\sigma_2$ on 2-blocks, apply it as defined above consecutively to each 
sub-1-block, juxtaposing the 2-blocks so obtained from left to right in the 
original
order 
%\gapp{29}
%\noindent
of the 1-blocks. For $ {\alpha} \in D_n^k , \  k\geq 2 $,
define $ \sigma_2 ( {\alpha} ) $ to be $ \sigma_2 $ of the first sub-2-block. 

\item[(h)]   We now presume all maps defined for subscripts $j  <  k < n$, on
all $ D_n^i $.  
To define ${\sigma}_k, \ {\tau}_k :   D_n^k \longrightarrow
  D_{n+1}^k $: 
  
  \noindent
For $k$ even, define
$$
     {\sigma}_k [ n ] \ \ = \ \  {\mu}_{k-1} \sigma_{k-1} [n-1]  \, | \, 
{\lambda}_{k-1}  \sigma_k [n-1] 
;
\qquad
\qquad
{\tau}_k [ n ] \ \ = \ \  {\nu}_{k-1}  {\tau}_{k}
[n-1]  \, | \,   {\mu}_{k-1} \tau_{k-1}  [n-1]  
$$
For $k$ odd we set
$$
     {\sigma}_k [ n ] \ \ = \ \left|  
     \begin{array}{c}
     {\mu}_{k-1} \sigma_{k-1} [n-1]  \\ 
{\lambda}_{k-1}  \sigma_k [n-1] \\
\end{array}
\right|
;
\qquad
\qquad
{\tau}_k [ n ] \ \ = \ \left|  
     \begin{array}{c}
      {\nu}_{k-1}  {\tau}_{k}[n-1]   \\ 
      {\mu}_{k-1} \tau_{k-1}  [n-1]  \\
\end{array}
\right|
$$

This allows us to define ${\sigma}_k[n] $ for all blocks after we set 
$ {\sigma}_n [n] =  \tau_n [n] := [n]$ for $ k \  \geq \  n $.
To extend the definition to higher blocks, we use the intuition we have of
the underlying geometry.
A $k$-block corresponds to a possibly higher dimensional cube $x$, together 
with a number of others of dimension $k$, each ``filled out'' by lower dimensional
cubes. Then ${\sigma}_k$ and ${\tau}_k$ should correspond to the $k$-cubes already 
present, together with ${\sigma}_k$ or ${\tau}_k$ of $x$, with everything still appropriately stretched out by lower dimensional cubes. 
We merely break up the action of $x$ on the pile sets $v_j$, $j \ <\  k$, as 
the composition of the action of its $k$-dimensional faces.

Specifically, to define ${\sigma}_k( {\alpha} ) $, where $ {\alpha} $ is a $j$-block,
 
\subitem ($h_1$)   If $j \geq k$, define ${\sigma}_k( {\alpha} ) $ to be ${\sigma}_k$ applied to  the
first $k$-block of $ {\alpha} $.

\subitem ($h_2$)   For $j\ =\ k$, define ${\sigma}_k$ to be the juxtaposition of ${\sigma}_k$
 applied to 
the consecutive ($k-1$)-blocks of $ {\alpha} $. Hence we need only 
be specific for $j \  \leq \  k-1$. 
 
\subitem ($h_3$)   If  $ {\alpha} $ has dimension $k$, we set ${\sigma}_k( {\alpha} ) := {\alpha} $. 
Otherwise there is a 0-block 
%\gapp{30}
%\noindent
contained in a unique
sub-($k$-1)-block $ \gamma $ of dimension $m\ >\ k$. Since ${\sigma}_k( \gamma ) $ is a
$k$-block,  we stretch it out by the same sub-blocks which stretch out $ \gamma  $ inside $ {\alpha} $. Now $ \gamma $ is contained in a 1-block, with 0-blocks
above and below. We add these 0-blocks above and below every 1-block of 
${\sigma}_k( {\alpha} ) $. 

\subitem ($h_4$)   Similarly, $ \gamma $ occurs in a unique 2-block, with 1-blocks occurring 
before and after the 1-block in which $ \gamma $ occurs.
We thus add each of these 1-blocks appropriately before and after every
1-block of ${\sigma}_k( \gamma ) $.

\subitem ($h_5$)   Proceeding in this way, up to the highest dimensional block, we fill out 
${\sigma}_k( \gamma ) $ to obtain the $k$-block we define to be ${\sigma}_k( {\alpha} ) $.

The procedure/definition merely tells us that if an element $x$ has been
stretched out by lower dimensional cubes, we use these to stretch
out the components of $x$.
 $$
\begin{array}{ccccc}
 \sigma_k ( {\alpha} )  &  =  &  \sigma_k 
\left|\begin{array}{c}
\alpha_1\\
\vdots\\
\alpha_i\\
\vdots\\
\alpha_t\\
\end{array} 
\right|&:=& 
\left|\begin{array}{c}
 \alpha_1^{(m)}\\
\vdots\\
 \sigma_k \alpha_i\\
 \vdots\\
 \alpha_t^{(m)}\\
\end{array} 
\right| \\
% &&  {\mu} _0 {\tau}_0{[n-1] }}\\
\end{array}
$$
where again $ {\alpha}_i $ is the distinguished sub-block, and repetitions
are required for the other blocks. We proceed similarly for ${\tau}_k$ and for the 
case $k$ even. 

\item[(i)]   Now define, for $k$ odd, $ {\lambda}_k$   and   ${\nu}_k :
 D_n^k \  \to \  D_{n+1}^k $ by
$$
\begin{array}{ccccc}
\lambda_k ( {\alpha} )  &  =  &  \lambda_k 
\left|\begin{array}{c}
\alpha_1\\
\vdots\\
\alpha_t\\
\end{array} 
\right|&:=& 
\left|\begin{array}{c}
\lambda_{k-1} \alpha_1\\
\vdots\\
\lambda_{k-1} \alpha_t\\
\mu_{k-1}\tau_{k-1}[n]\\
\end{array} 
\right| \\
% &&  {\mu} _0 {\tau}_0{[n-1] }}\\
\end{array}
;
\qquad\qquad
\begin{array}{ccccc}
 \nu_k ( {\alpha} )  &  =  &  \nu_k 
\left|\begin{array}{c}
\alpha_1\\
\vdots\\
\alpha_t\\
\end{array} 
\right|&:=& 
\left|\begin{array}{c}
\mu_{k-1}\sigma_{k-1}[n]\\
\nu_{k-1} \alpha_1\\
\vdots\\
\nu_{k-1} \alpha_t\\
\end{array} 
\right| \\
\end{array}
$$
% \gapp{31}

For k even, set
$$
\nu_k ( {\alpha} )    =    \nu_k 
\left|\begin{array}{ccc}
\alpha_1&
\cdots&
\alpha_t\\
\end{array} 
\right|
:=
\left|\begin{array}{c|c|c|c}
\nu_{k-1} \alpha_1&
\cdots&
\nu_{k-1} \alpha_t&
\mu_{k-1}\tau_{k-1}[n]\\
\end{array} 
\right| 
$$
$$
\lambda_k ( {\alpha} )    =    \lambda_k 
\left|\begin{array}{ccc}
\alpha_1&
\cdots&
\alpha_t\\
\end{array} 
\right|
:=
\left|\begin{array}{c|c|c|c}
\mu_{k-1}\sigma_{k-1}[n]&
\lambda_{k-1} \alpha_1&
\cdots&
\lambda_{k-1} \alpha_t\\
\end{array} 
\right| 
$$

\item[(j)]   Now define for k odd 
$$
\begin{array}{ccc}
\mu_k  ( {\alpha} )  &  =  & \mu_k 
\left|\begin{array}{c}
\alpha_1\\
\hline
\vdots\\
\hline
\alpha_t\\
\end{array} 
\right|
\end{array}:= 
\left|
 \begin{array}{c|c|c|c}
 \begin{array}{c}
 \mu_{k-1}\alpha_1\\
 \hline
\nu_{k-1}\sigma_k\alpha_2\\
\hline
\vdots\\
\hline
\nu_{k-1}\sigma_k1\alpha_t\\
\end{array} 
&
 \begin{array}{c}
\lambda_{k-1}\tau_k\alpha_1\\
\hline
\mu_{k-1}\alpha_2\\
\hline
\vdots\\
\hline
\nu_{k-1}\sigma_k\alpha_t\\
\end{array} 
&
\cdots
& \begin{array}{c}
\lambda_{k-1}\tau_k\alpha_1\\
\hline
\lambda_{k-1}\tau_k\alpha_2\\
\hline
\vdots\\
\hline
\mu_{k-1} \alpha_t\\
\end{array} 
% &&  {\mu} _0 {\tau}_0{[n-1] }}\\
\end{array}
\right|
$$
 
For k even, 
 $$
\mu_k ( {\alpha} )    =    \mu_k 
\left|\begin{array}{c|c|c}
\alpha_1&
\cdots&
\alpha_t\\
\end{array} 
\right|
:=
\begin{array}{c}
\left|\begin{array}{c|c|c|c}
\mu_{k-1}\alpha_1&
\lambda_{k-1} \sigma_{k}\alpha_2&
\cdots&
\lambda_{k-1}\sigma_{k} \alpha_t\\
\end{array} 
\right| \\
\hline
\hline
\left|\begin{array}{c|c|c|c}
\nu_{k-1}\tau_{k}\alpha_1&
\mu_{k-1} \alpha_2&
\cdots&
\lambda_{k-1}\sigma_{k} \alpha_t\\
\end{array} 
\right| \\
\hline
\hline
\vdots\\
\hline
\hline
\left|\begin{array}{c|c|c|c}
\nu_{k-1}\tau_{k}\alpha_1&
\nu_{k-1}\tau_k \alpha_2&
\cdots&
\mu_{k-1} \alpha_t\\
\end{array} 
\right| \\
\end{array}
$$

\end{description}

Having defined ${\sigma}_k$ and ${\tau}_k$ for $ {\cal I}^n $, it is possible to define a structure on the infinite dimensional cube, characterised as
$$ {\cal I}^\infty     :=  \{ f : \N \to {\cal I} \ | \  f^{-1} (0)\ \ {\rm  is\ finite }\}.
$$
 
Notions of well-formedness of subsets of $  {\cal I}^\infty$ can be given as in
\cite{St}, using the cocycle conditions, and it is straightforward to prove

 \medskip
 \medskip
 \noindent
 {\bf
Proposition 6.1:}
There is an $n$-category structure defined on $ {\cal I}^\infty $.

%\gapp{32}

  \medskip
 \medskip
 \noindent
 {\bf
Proposition 6.2:}
(i)  If  $ {\alpha} $ is a full $j$-sub-block of $ {\sigma}_k [n] $, then $ {\sigma}_i ({\alpha} ) $ = $ {\sigma}_i [n] $  for $i  < j < k$, and similarly for
targets. \hfill\break
          (ii)   If $ {\alpha}   *_j  \beta $ is an
occurence of consecutive $j$-blocks of $ {\sigma}_k [n] $, 
then $ {\sigma}_j ( \beta ) =  {\tau}_j ( {\alpha} )$. Moreover, topologically the union
of sub-cubes of $ {\sigma}_j ( \beta ) $ is a disk of dimension $j$. \hfill\break
          (iii)   The piles generated by applications of elements of
the cocycle blocks are nested embeddings of $k$-disks $ V_k $ with geometric boundary $ \partial V_k = \sigma_{k-1} \  \cup \  \tau_{k-1} $.

\section{Low dimensional examples}
 
We gave diagrams indicating how the first three cocycles arise,
repeating their description in block form accompanied by the geometric
decomposition corresponding to the block. Note that the 2-cocycle condition
can be thought of as saying ``replace the interior of the hexagon on the
left by the configuration on the right''. In other words, flip from the front 
faces to the back ones in the standard Necker cube. 
 
\begin{figure}[htbp] %  figure placement: here, top, bottom, or page
   \centering
   \includegraphics[width=4in]{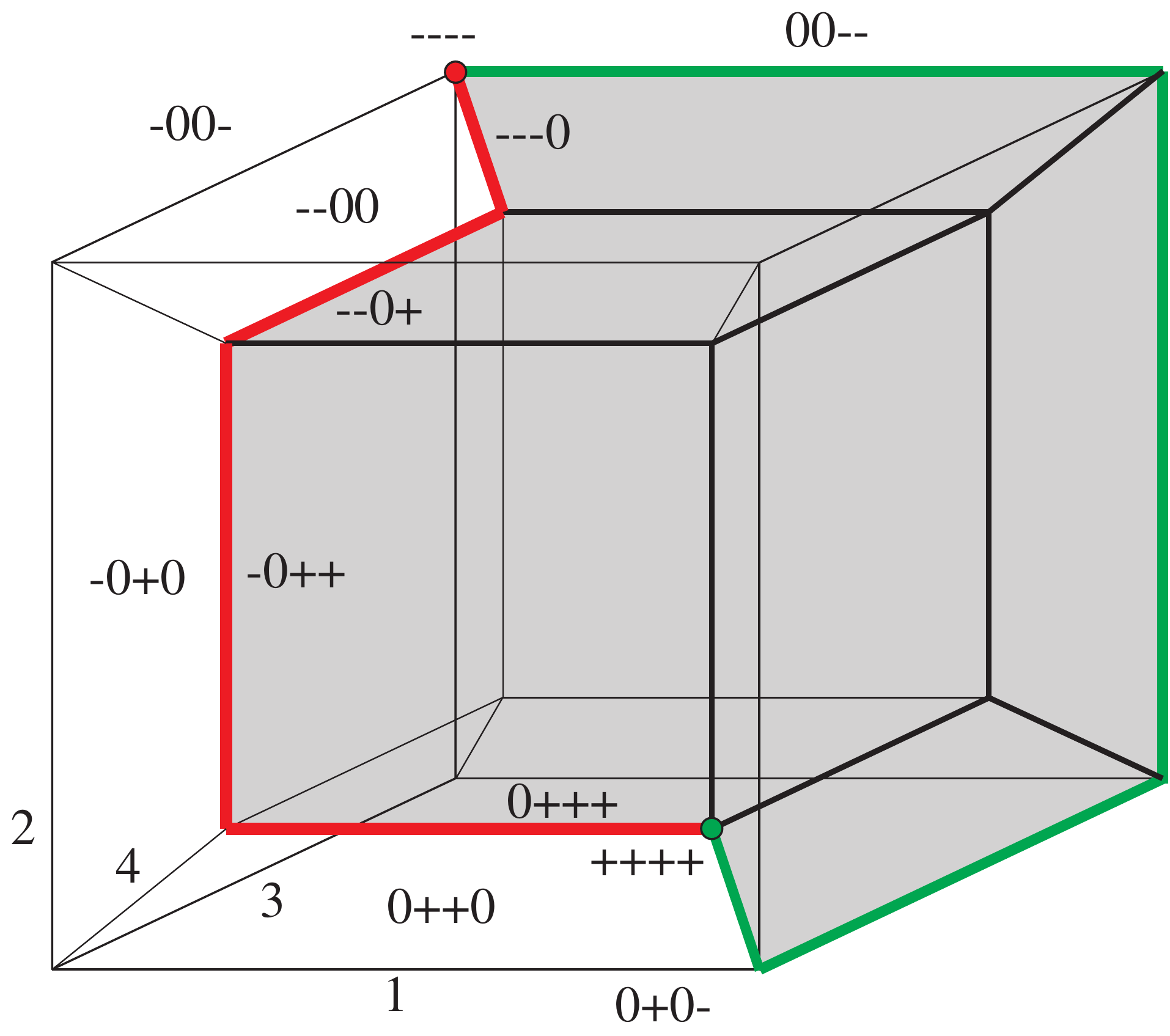} 
   \caption{The shaded squares give $\tau_2[4]$: Left edges of the shaded squares give $\sigma_1[4]$; the right edges give
   $\tau_1[4]$. Labels   reflect the structure arising from $[3]$.}
   \label{fig:15}
\end{figure}

\begin{figure}[htbp] %  figure placement: here, top, bottom, or page
   \centering
   \includegraphics[width=1.45in]{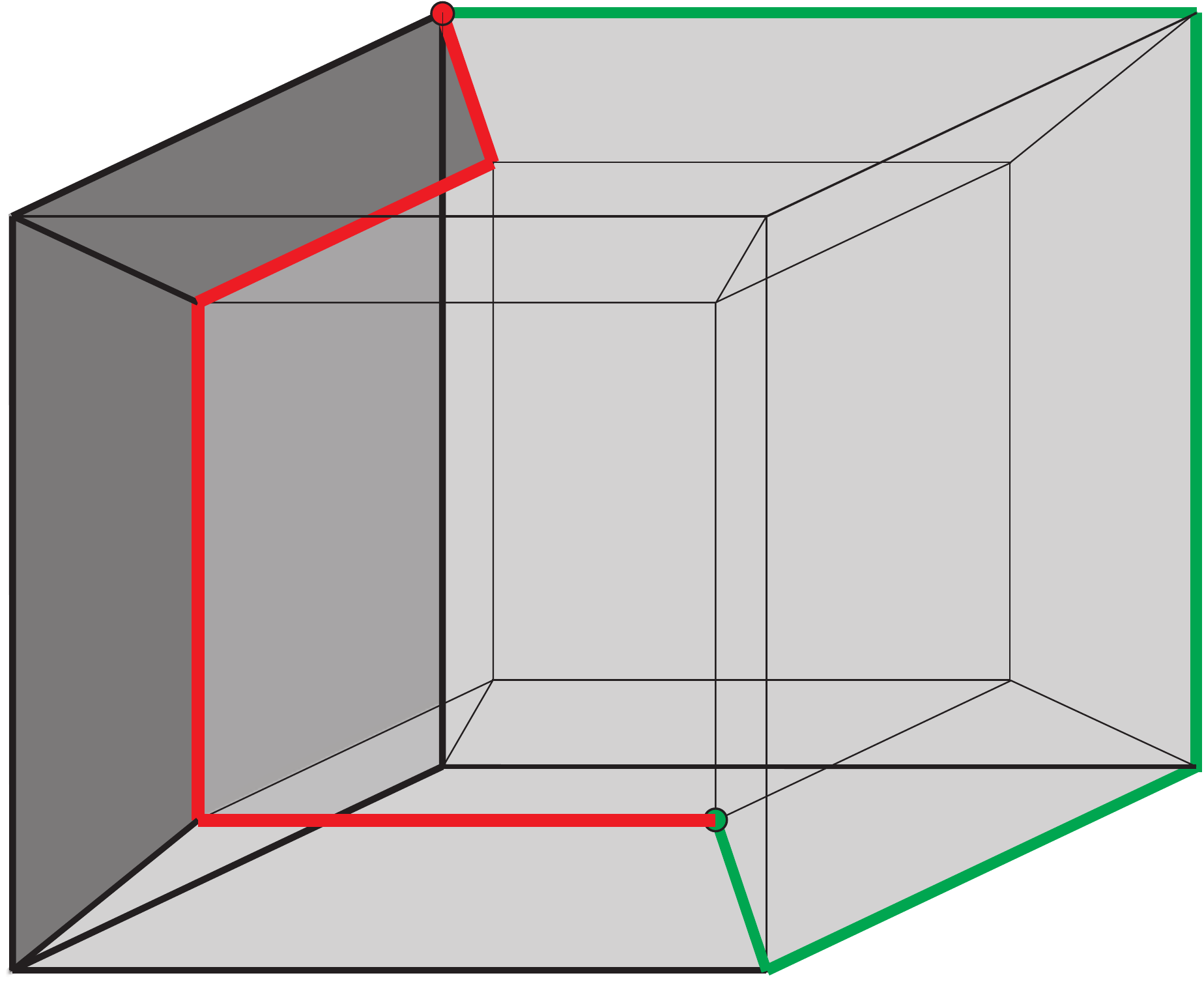} 
      \includegraphics[width=1.45in]{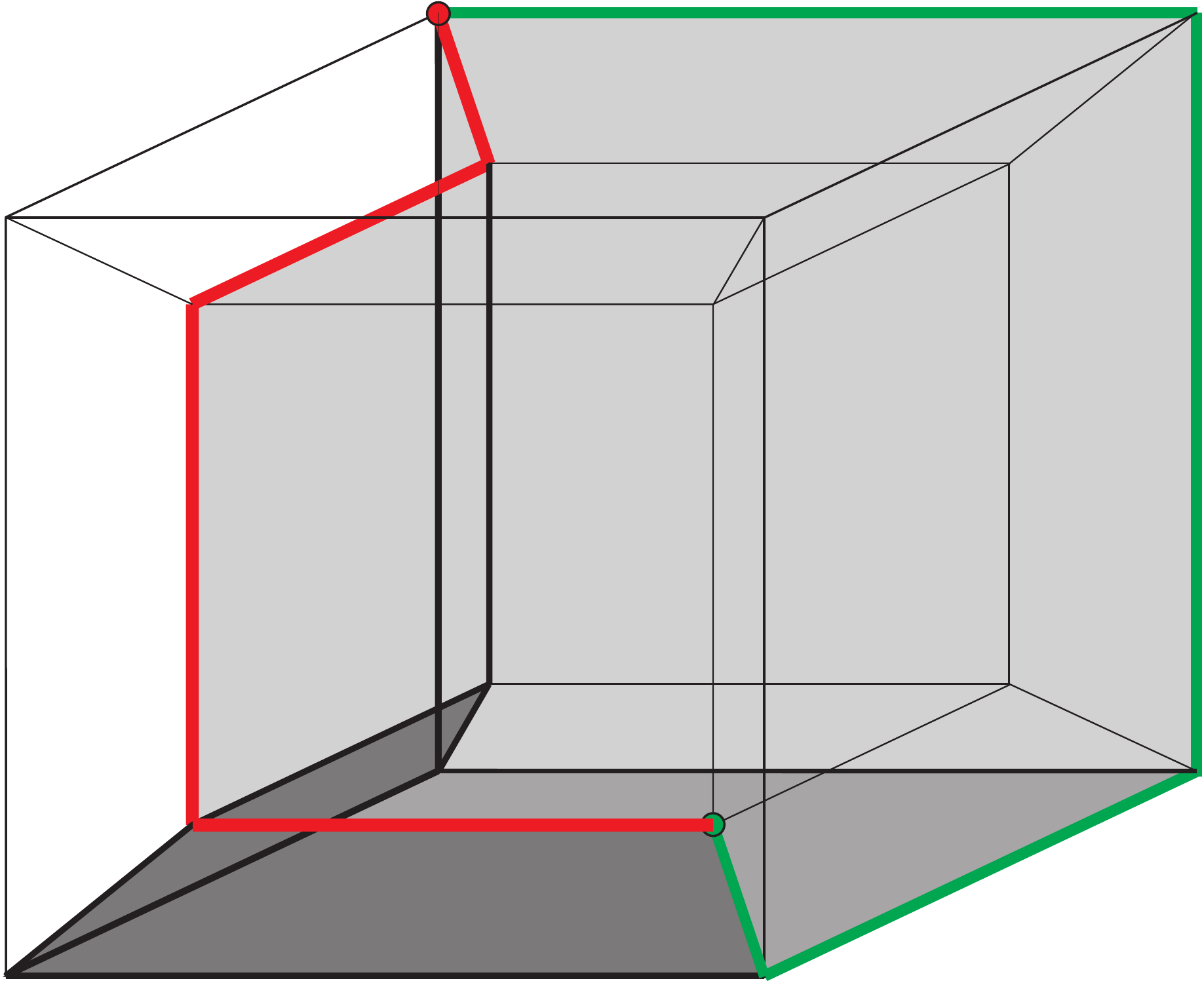} 
         \includegraphics[width=1.45in]{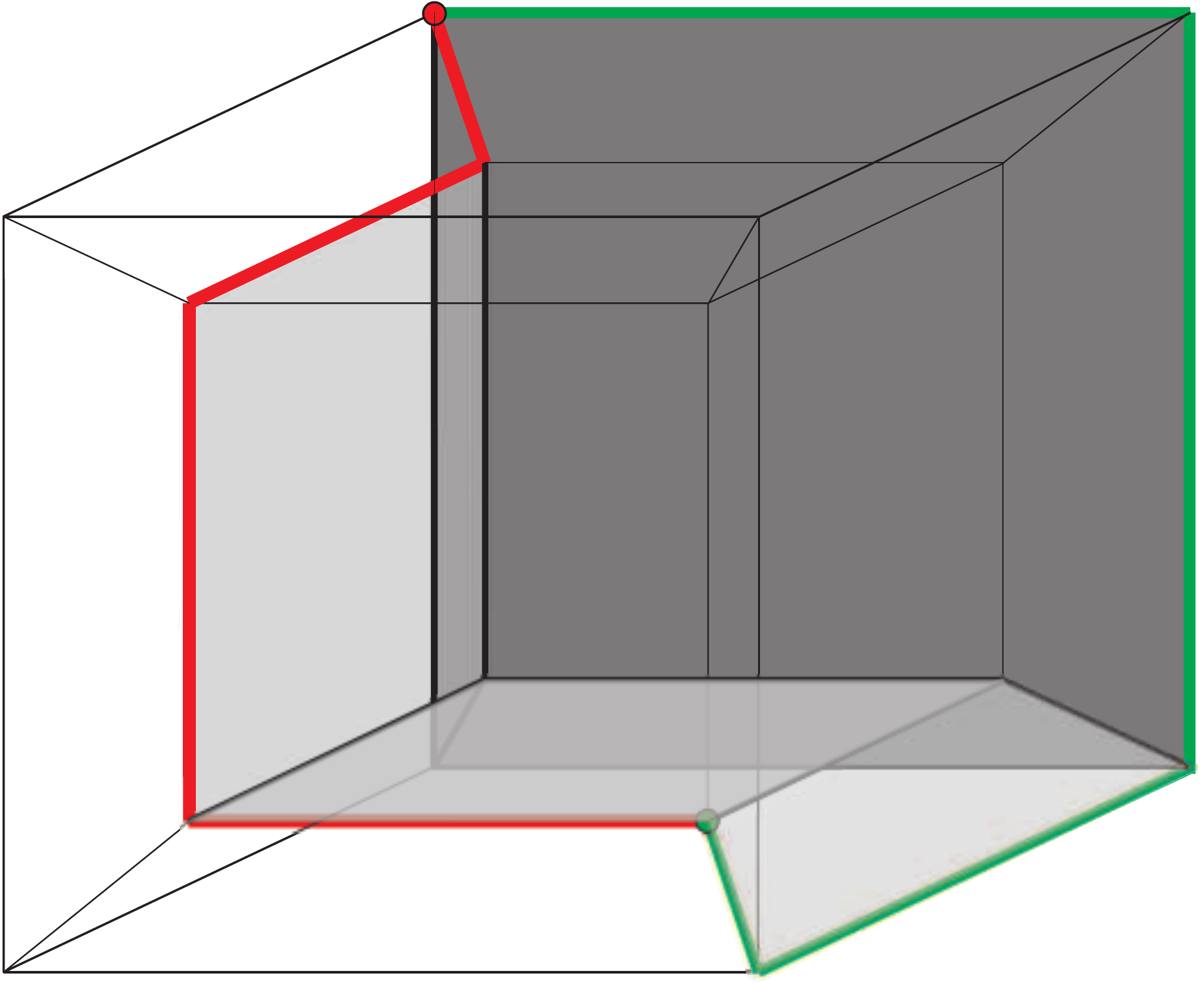} 
            \includegraphics[width=1.45in]{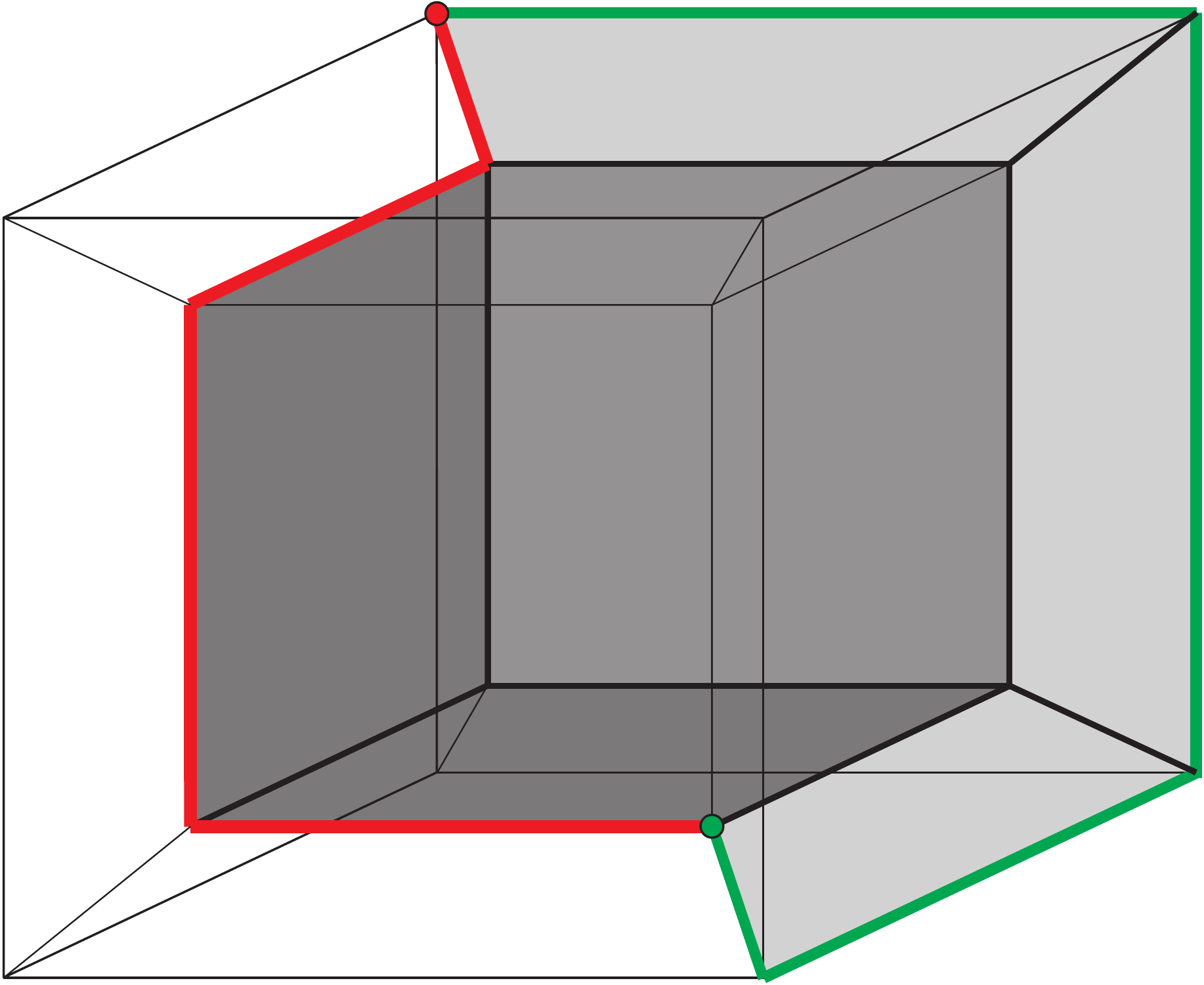} 
               \includegraphics[width=1.45in]{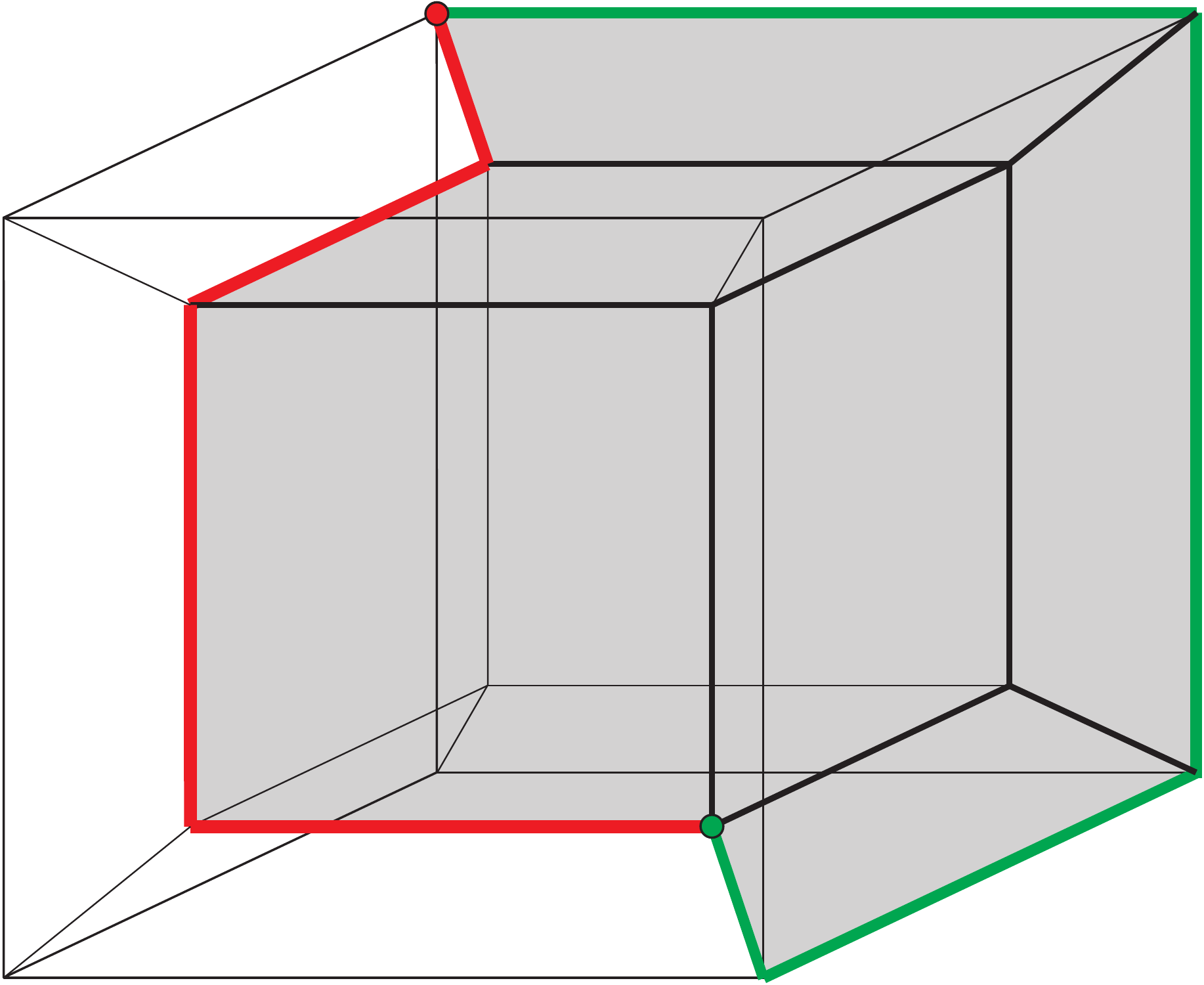} 
   \caption{The action of $\sigma_3[4]$ on $\sigma_2[4]$, producing $\tau_2[4]$.}
   \label{fig:16}
\end{figure}

Figure \ref{fig:15} shows the boundary of the 4-dimensional cube, labelled by our
convention. The 2-dimensional target faces are shaded, the source faces
labelled, and the heavy line 
depicts the 1-source. To depict the 3-cocycle condition, we present 
the five stages of deformation of the 2-face cubes across the four 3-dimensional cubes of the 3-dimensional source faces. This corresponds to looking at the 
set of $ {\cal I}^2 $s arising in the respective piles. These pictures decompose as indicated in Figure \ref{fig:16}, where the order is that determined by the corresponding block.

\begin{figure}[htbp] %  figure placement: here, top, bottom, or page
   \centering
   \includegraphics[width=4.5in]{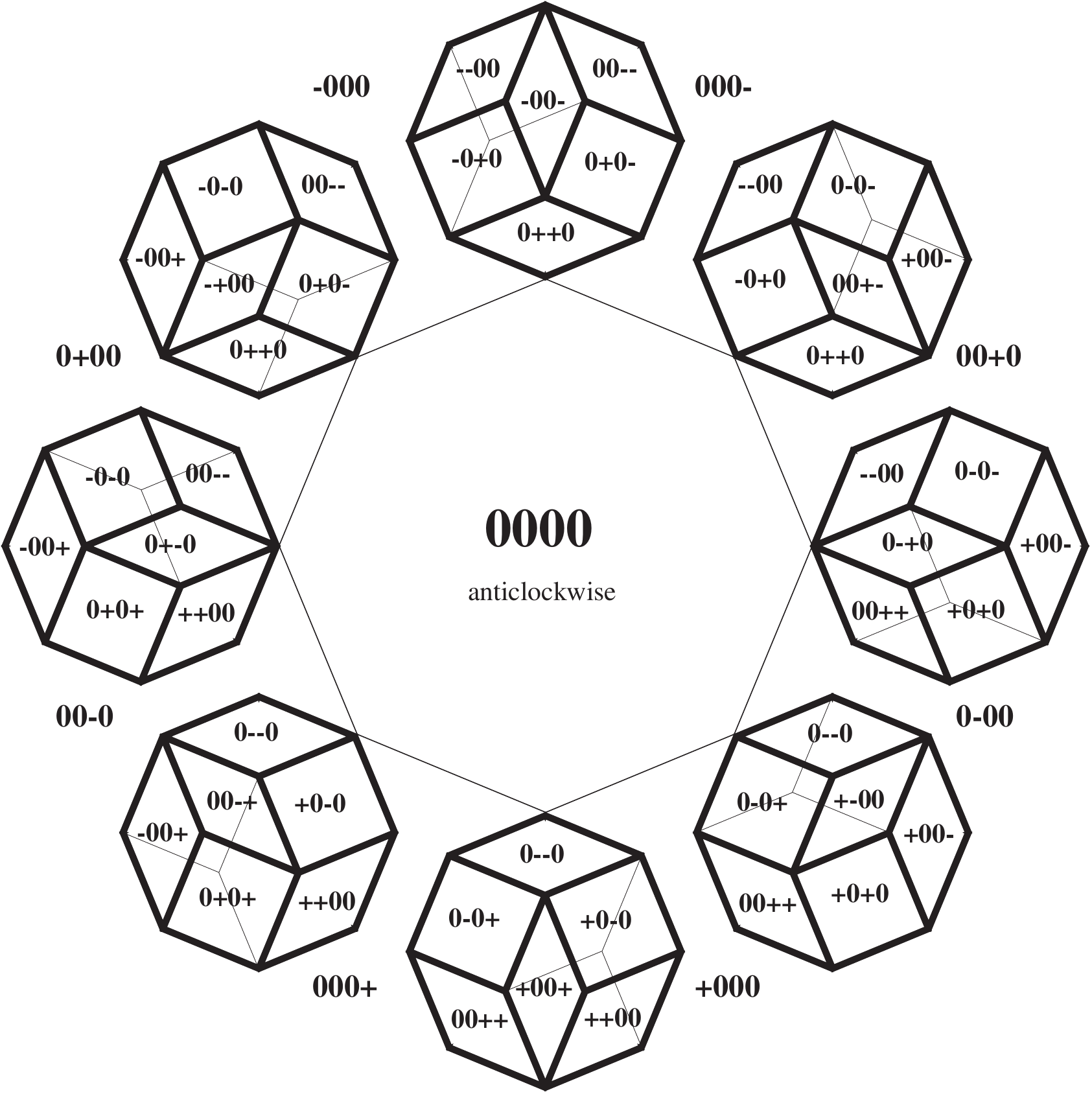} 
  \caption{The Octagon of Octagons: 3-cubes in $[4]$ acting on oriented 2-cells in the 4-cube. Observe the $\pi$-rotational symmetry and interchange of $\pm$. Faintly shaded 3-cubes are labeled correctly when read anticlockwise. Each octagon is rotated by $3\pi/4$ at each stage, with a consequent (not described here) action on the labelings of squares.}
%   \label{fig:example}
\end{figure}
\begin{figure}[htbp] %  figure placement: here, top, bottom, or page
   \centering
   \includegraphics[width=4.5in]{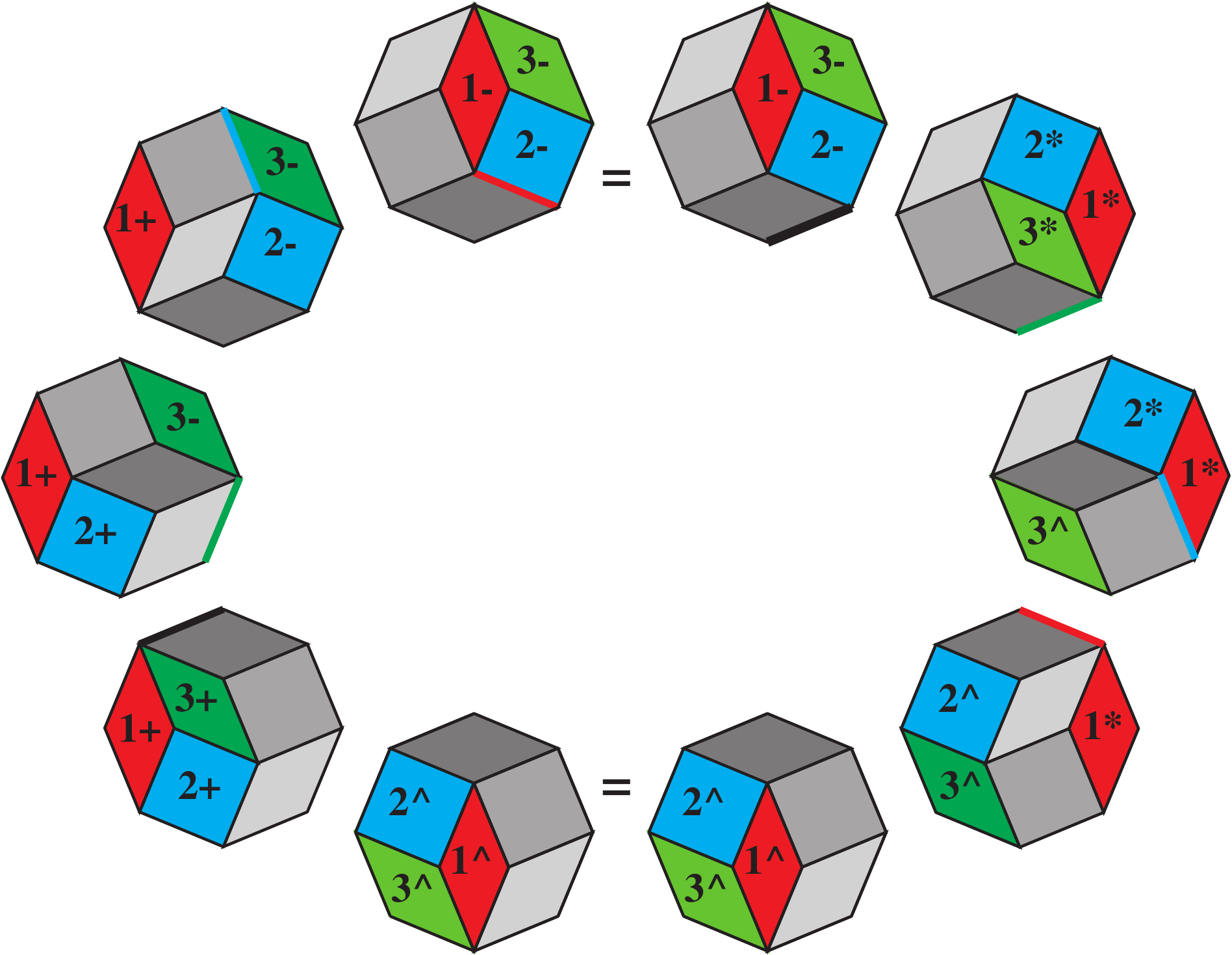} 
 %  \caption{example caption}
   \label{fig:17}
\end{figure}
\vfill
%\pagebreak

\begin{figure}[htbp] %  figure placement: here, top, bottom, or page
   \centering
   \includegraphics[width=6in]{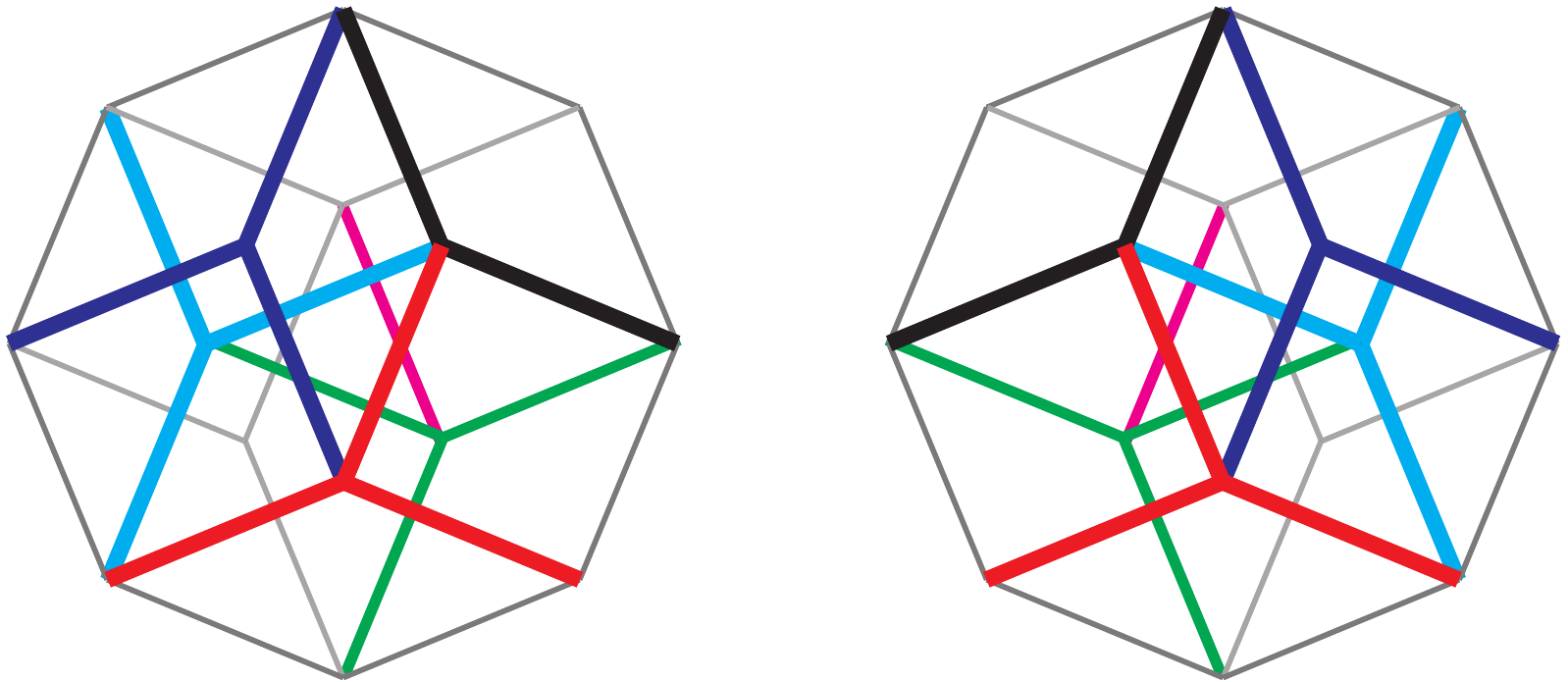} 
   \includegraphics[width=7in]{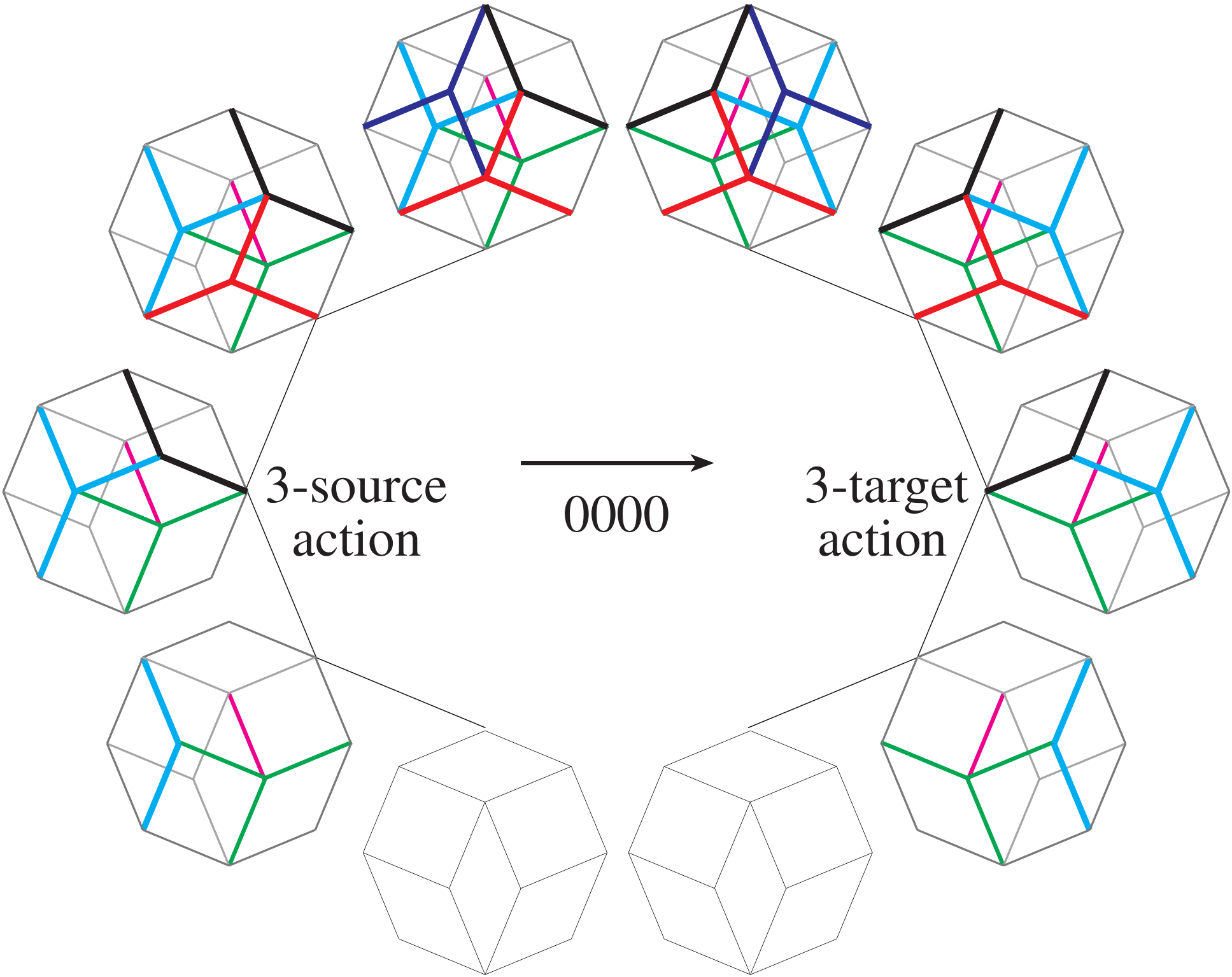} 
   \caption{3-cubes in the source and target of $[4]$}
   \label{fig:18}
\end{figure}

  \vskip -1.5in
 
 {\footnotesize
  \[
%\left(
\begin{array}{ccccc}  
      \includegraphics[width=1in]{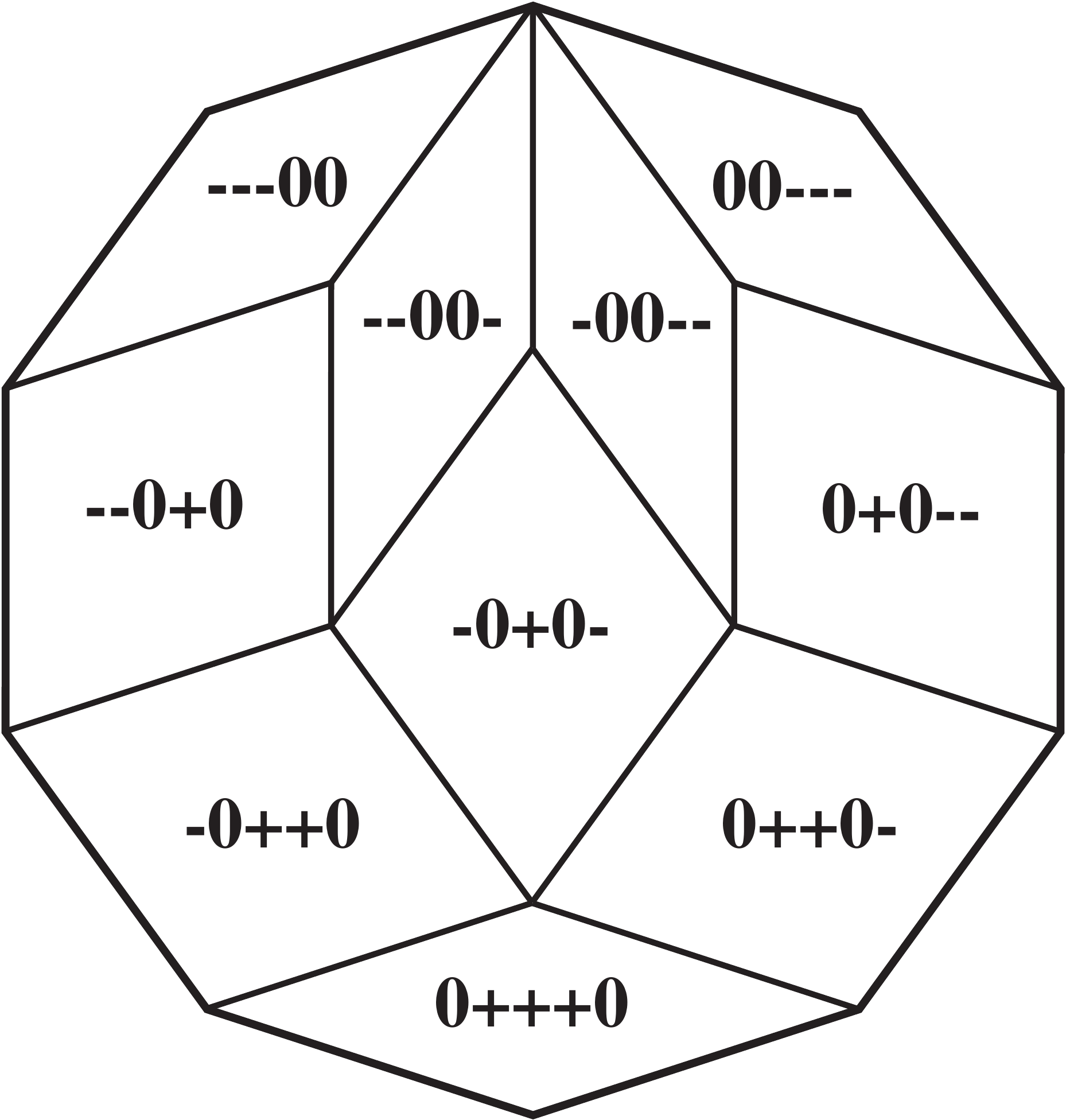}     
  &      \includegraphics[width=1in]{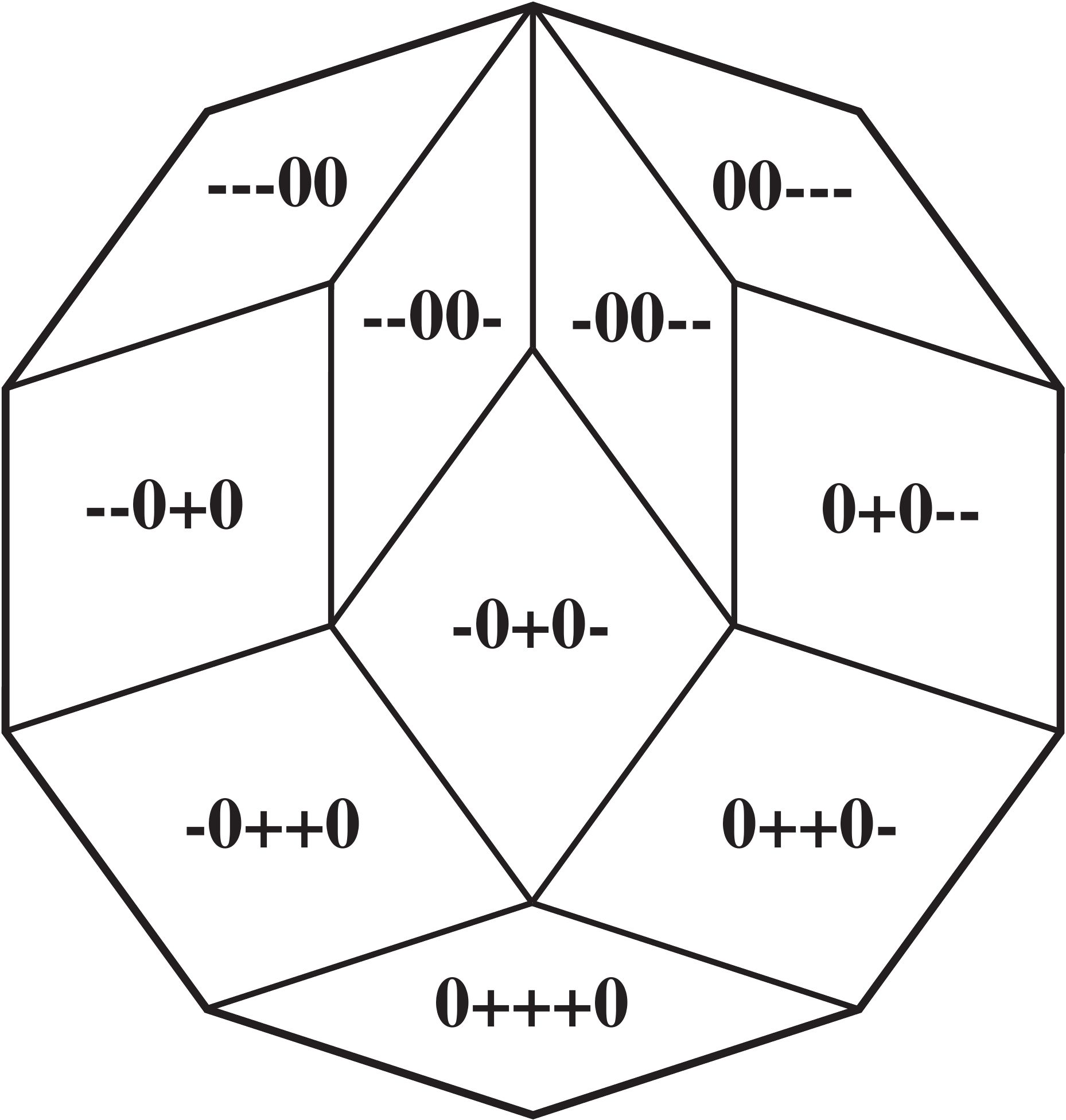}     
 &    \includegraphics[width=1in]{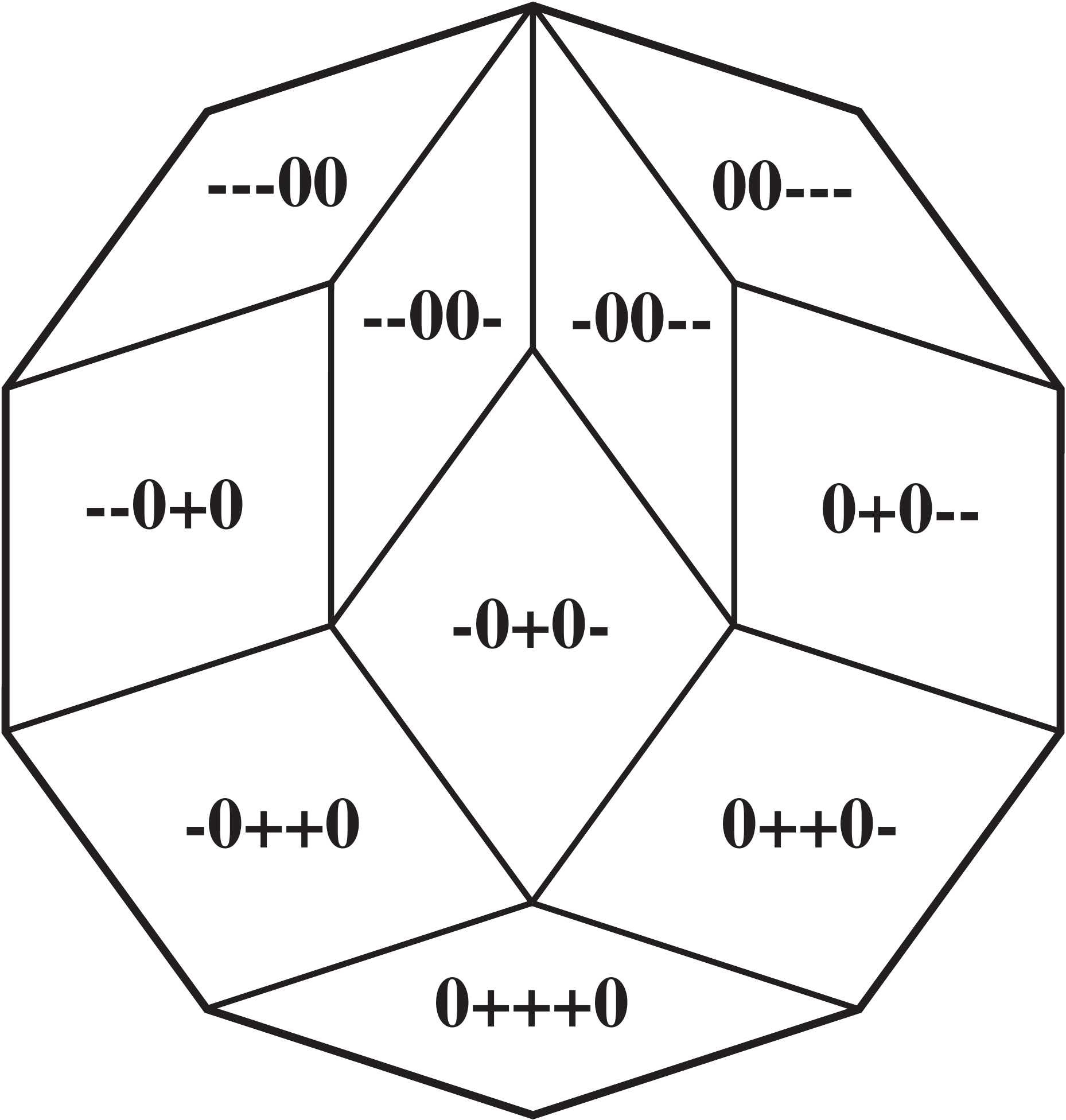}     
&     \includegraphics[width=1in]{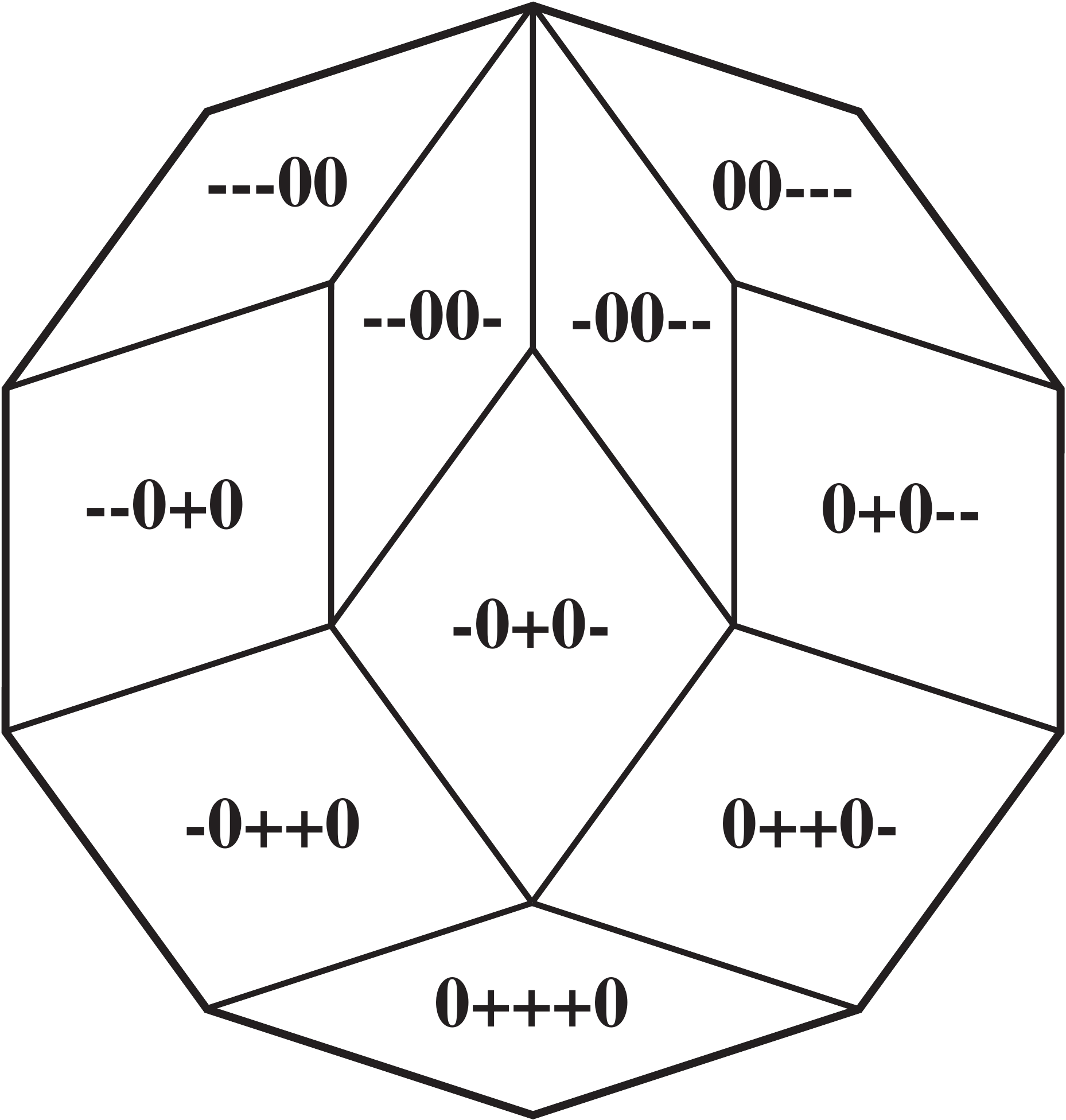}     
&     \includegraphics[width=1in]{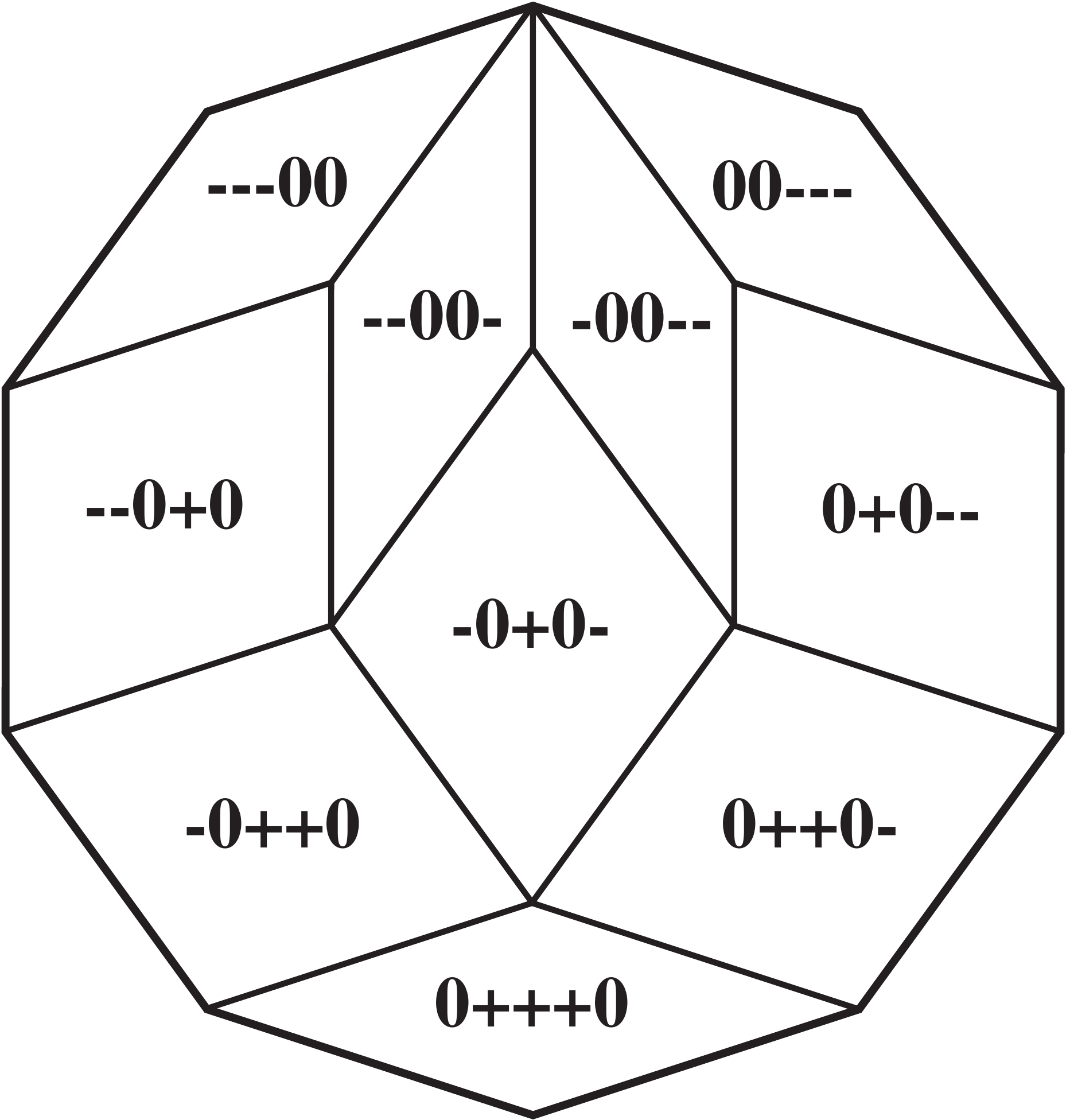} \\

  -0000 &  -000- &  -000- &  -000- & 0000-   \\

     \includegraphics[width=1in]{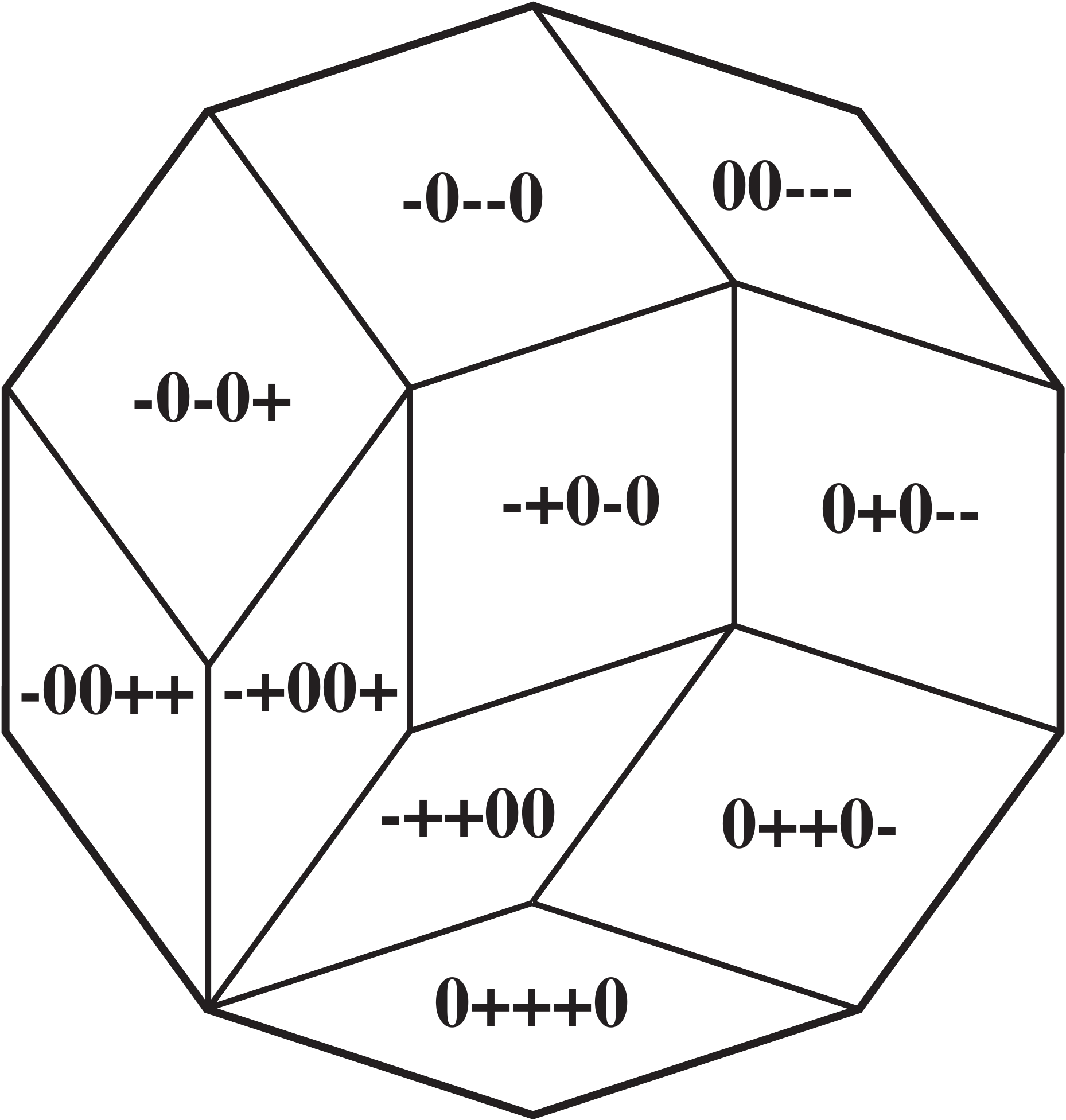}     
  &      \includegraphics[width=1in]{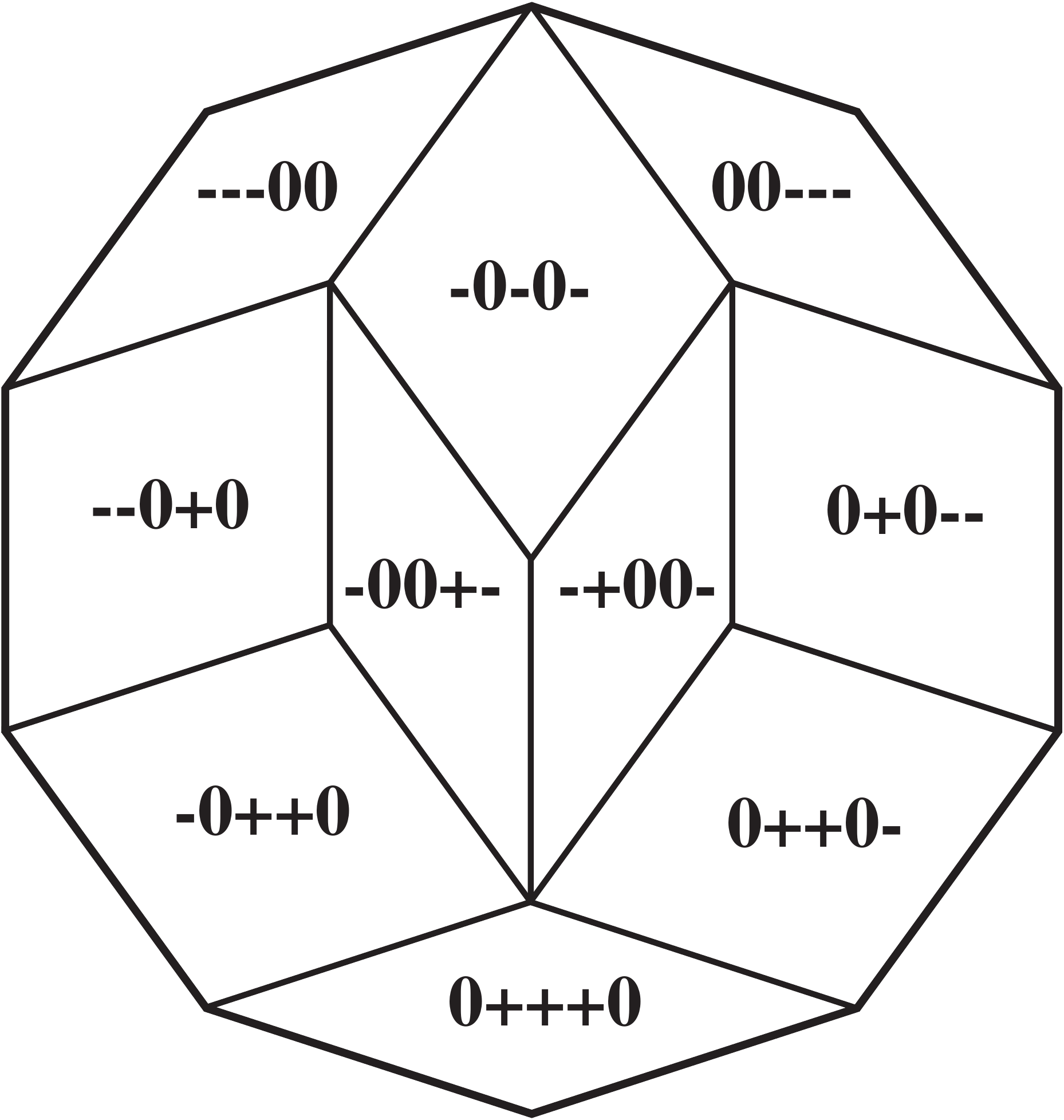}     
 &    \includegraphics[width=1in]{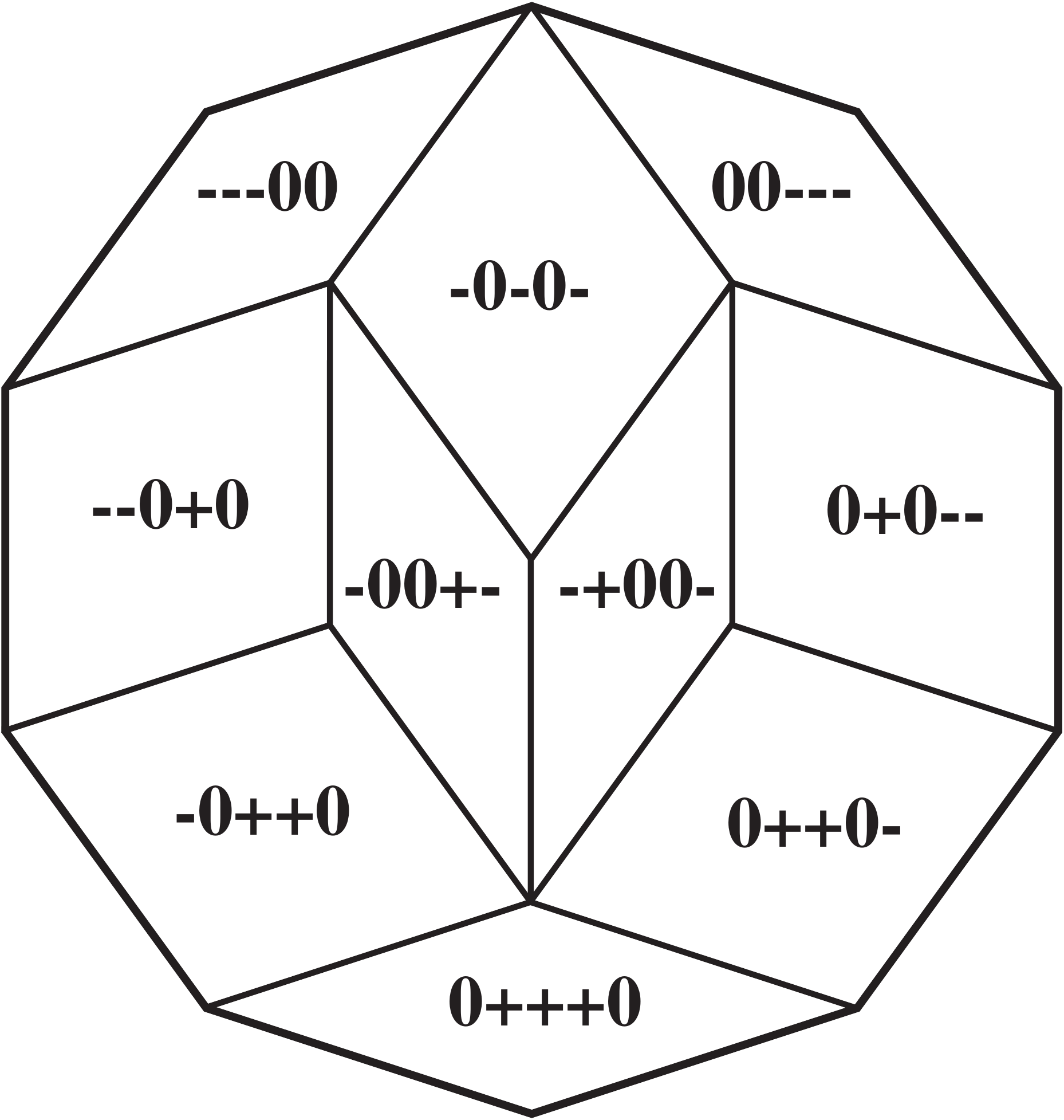}     
&     \includegraphics[width=1in]{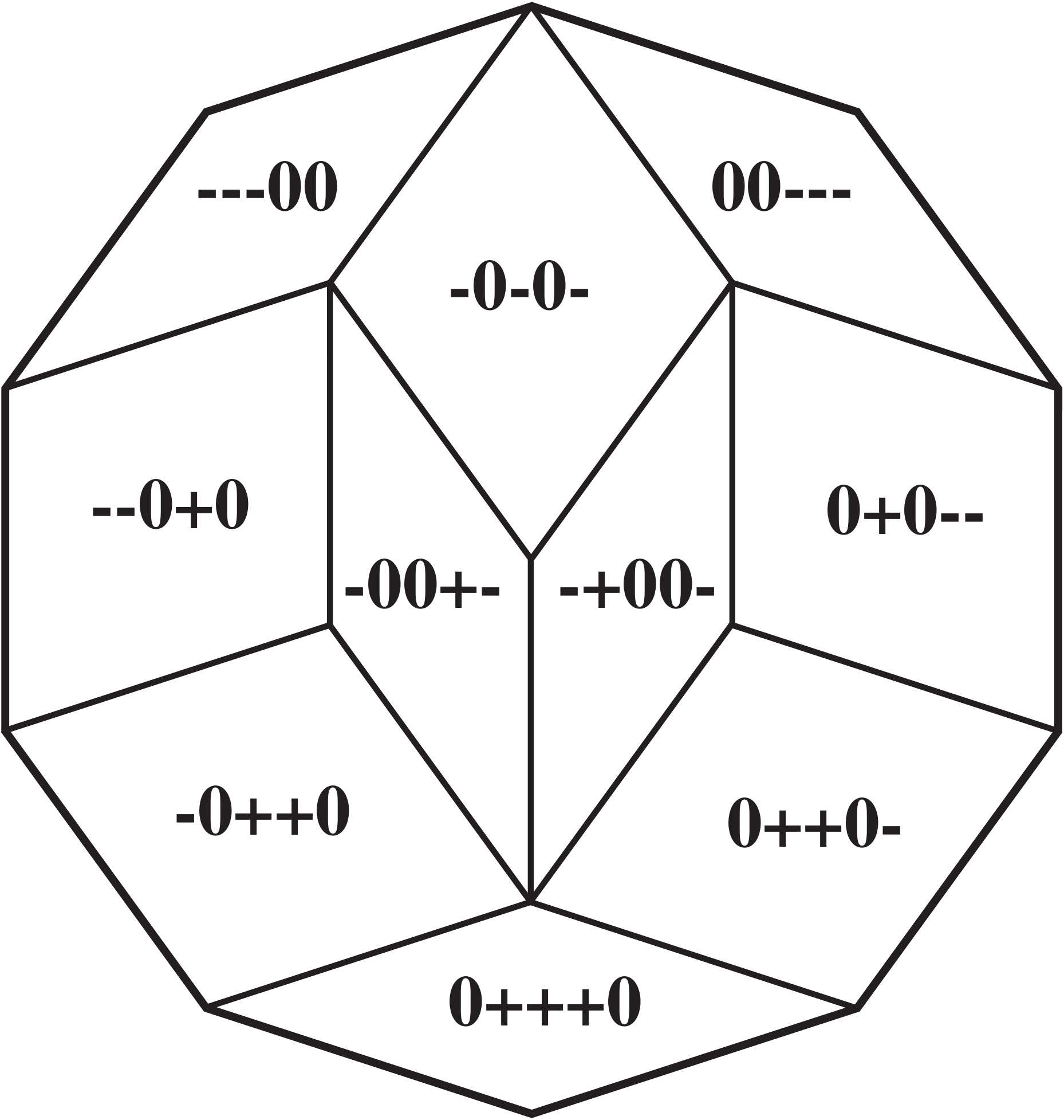}     
&     \includegraphics[width=1in]{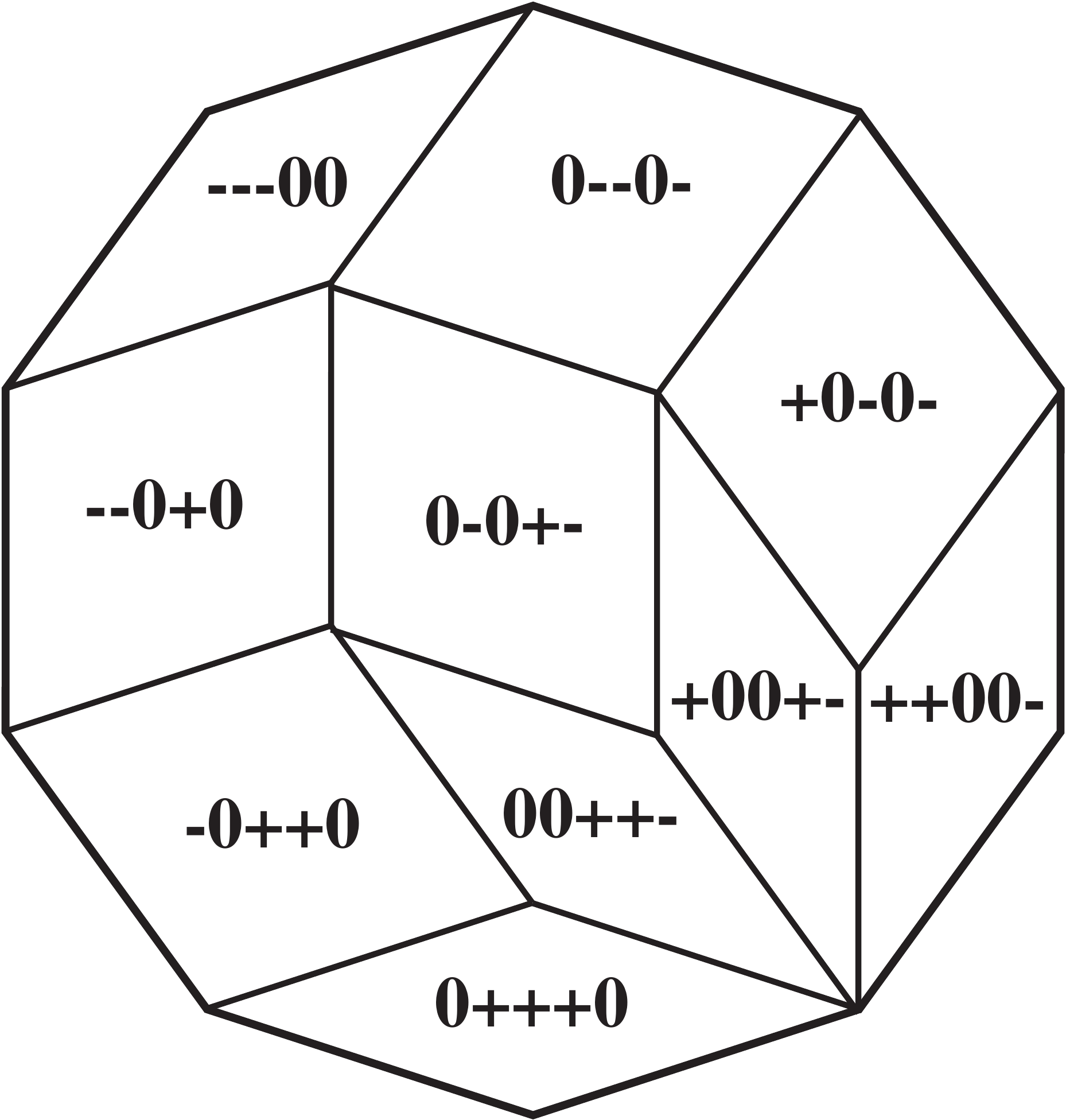} \\

 0++00 &  -00+0 &  0+00- &  0+00-   &   00++0\\

      \includegraphics[width=1in]{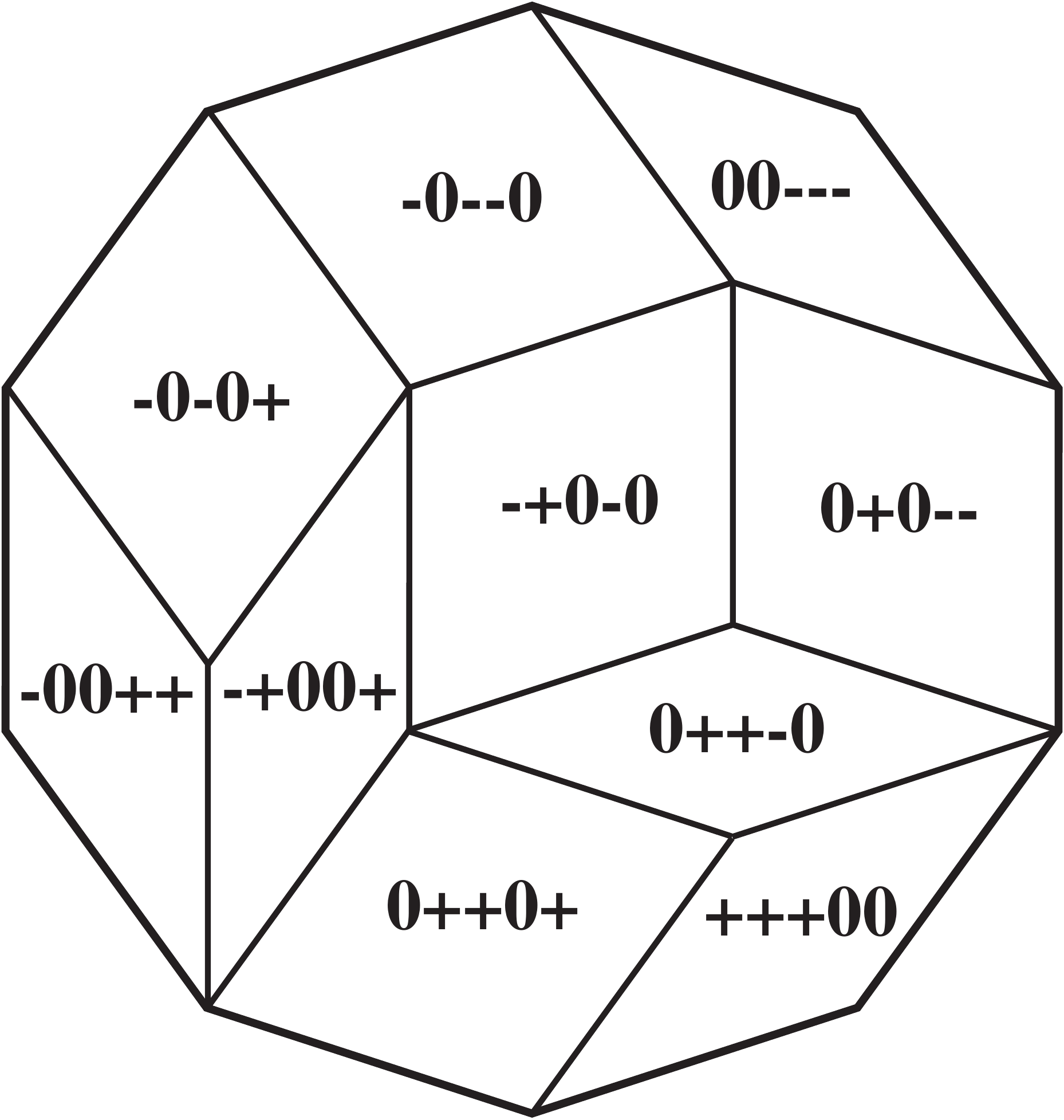}     
  &      \includegraphics[width=1in]{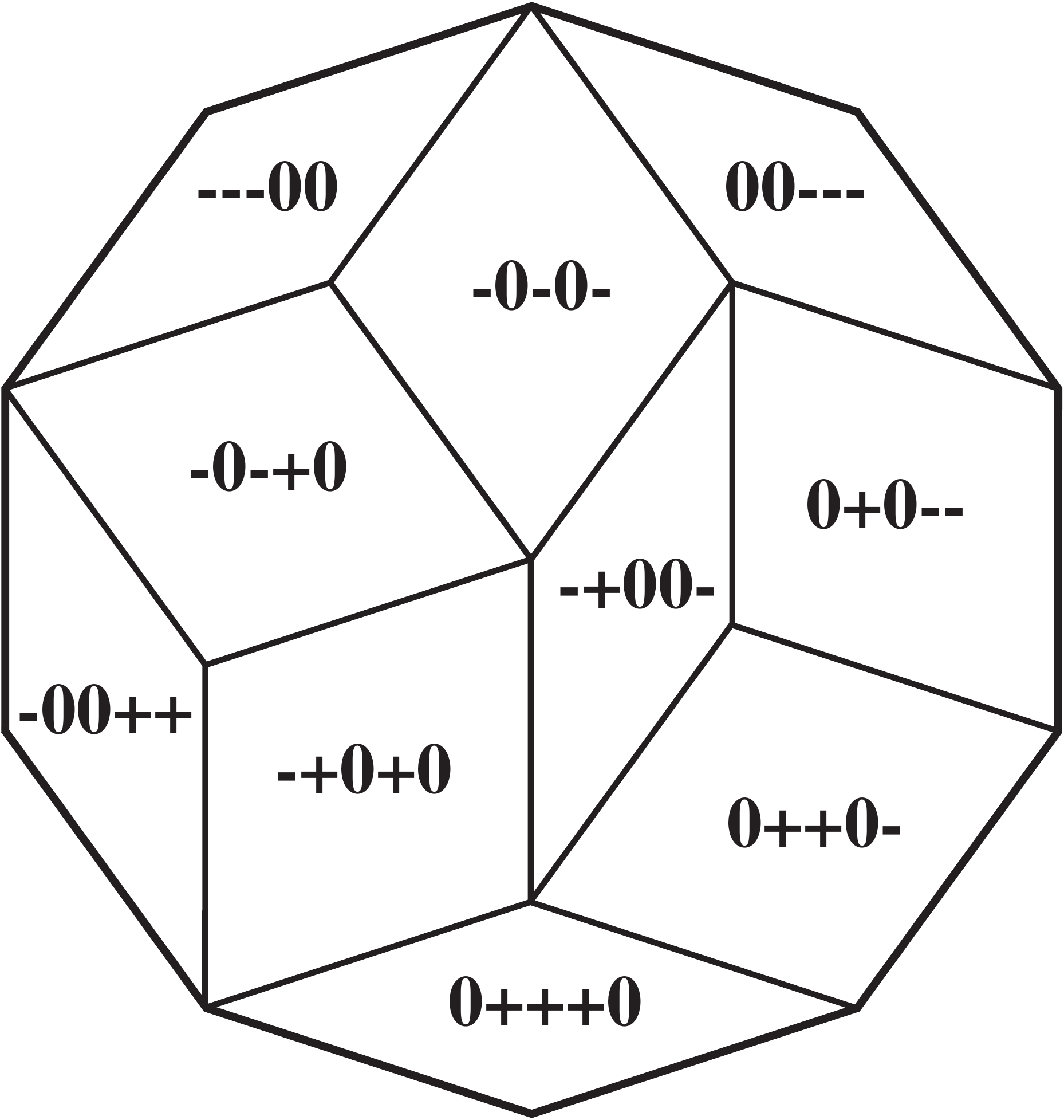}     
 &    \includegraphics[width=1in]{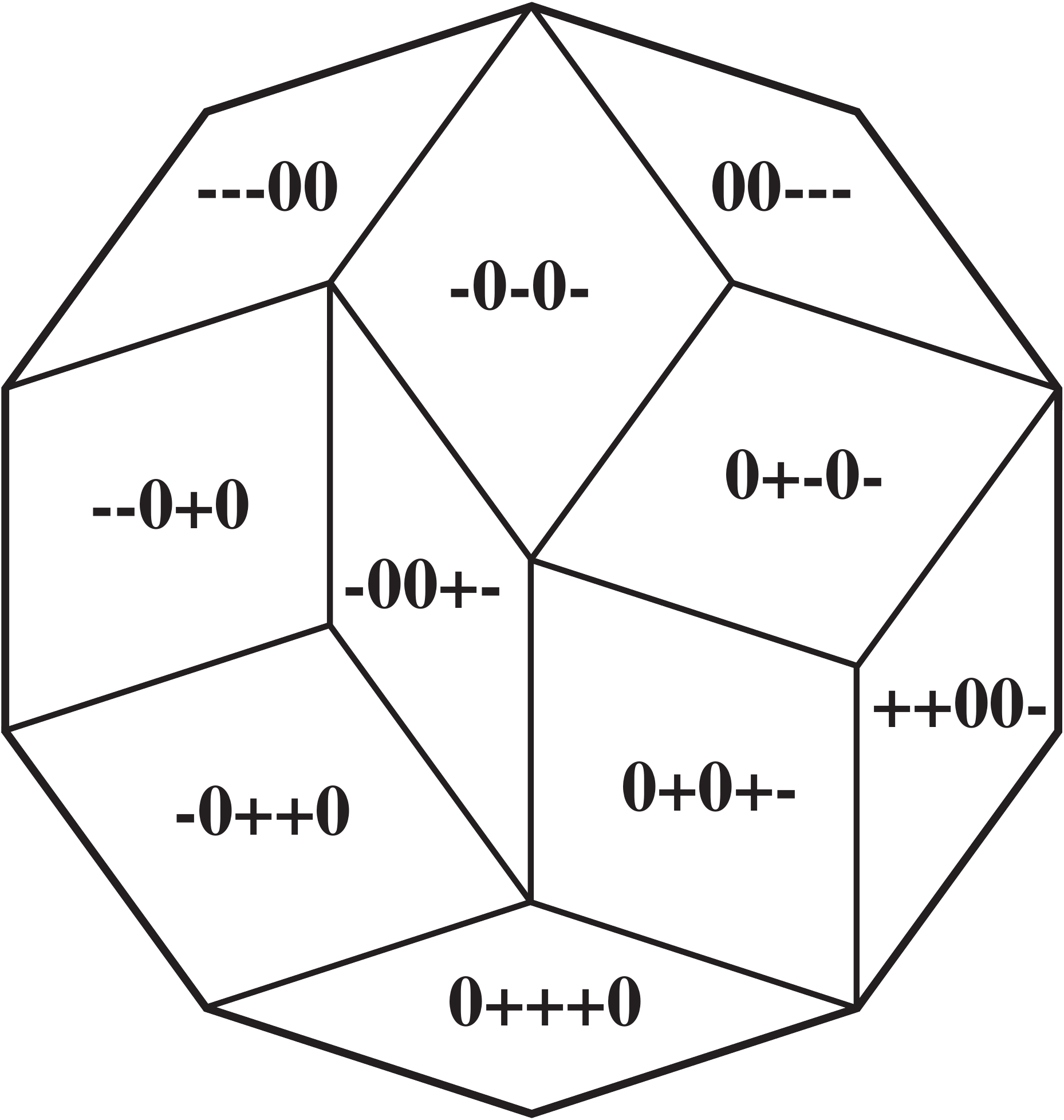}     
&     \includegraphics[width=1in]{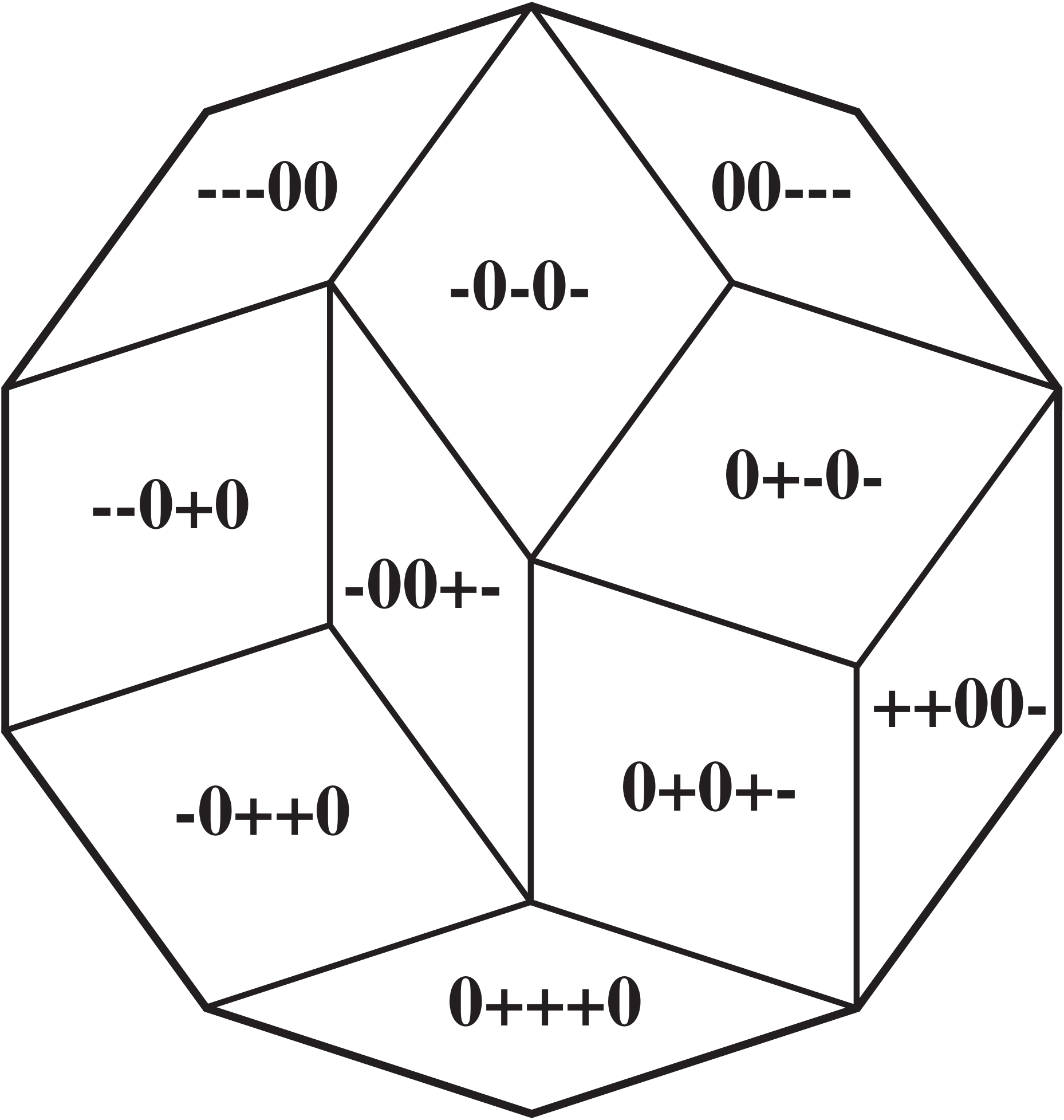}     
&     \includegraphics[width=1in]{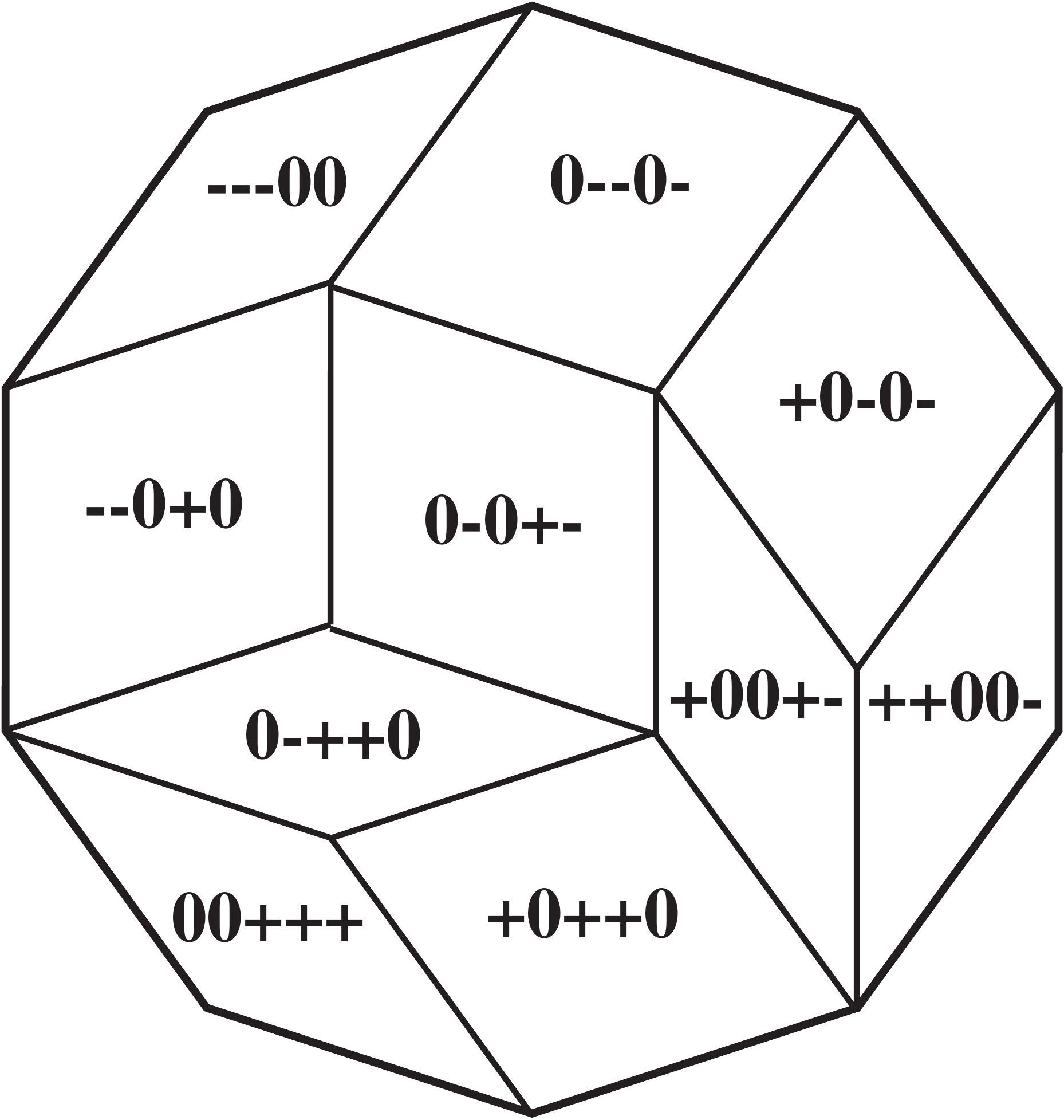} \\

 0+0-0 &  -0-00 &  -00+0 &  00-0- &   0-0+0 \\

      \includegraphics[width=1in]{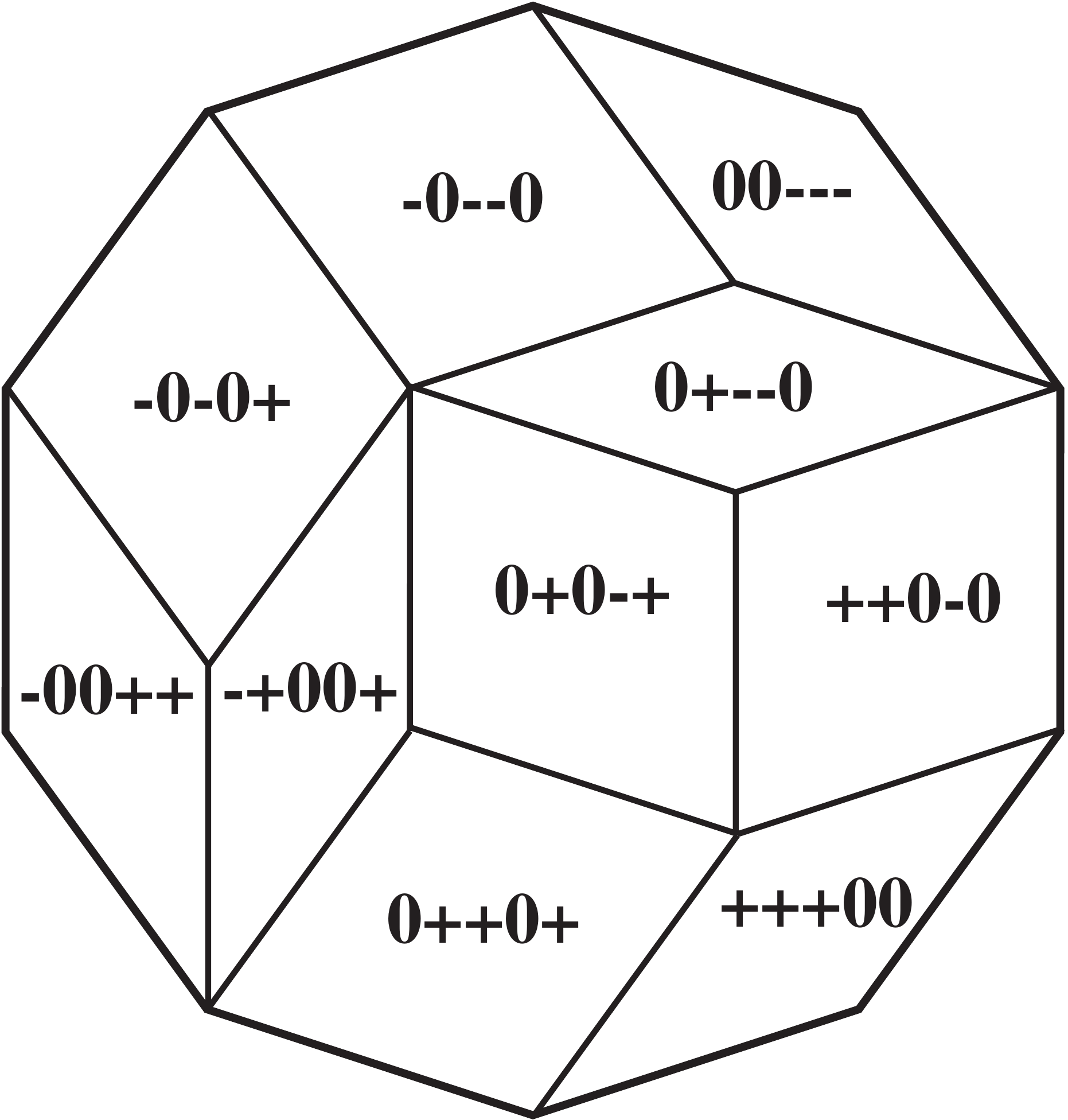}     
  &      \includegraphics[width=1in]{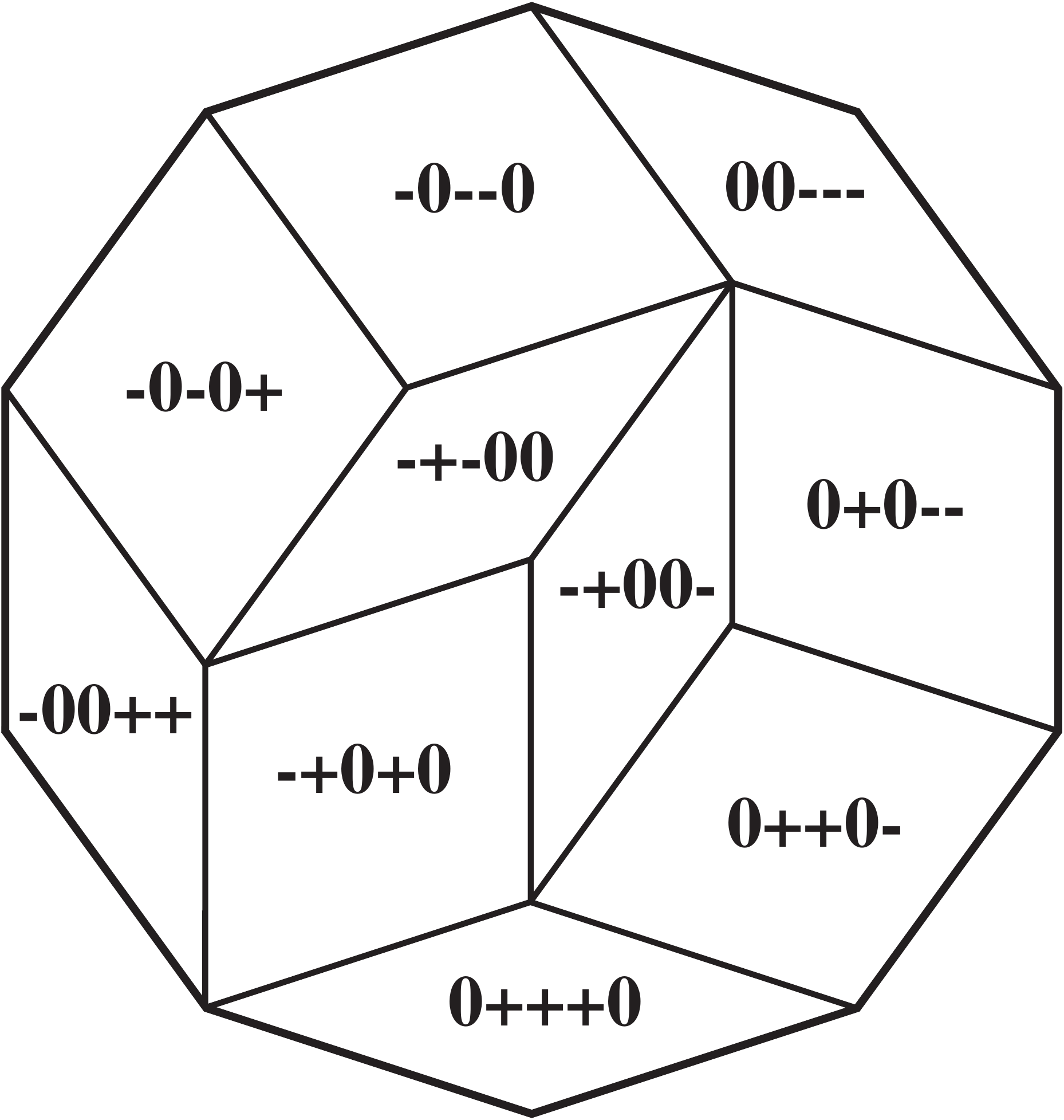}     
 &    \includegraphics[width=1in]{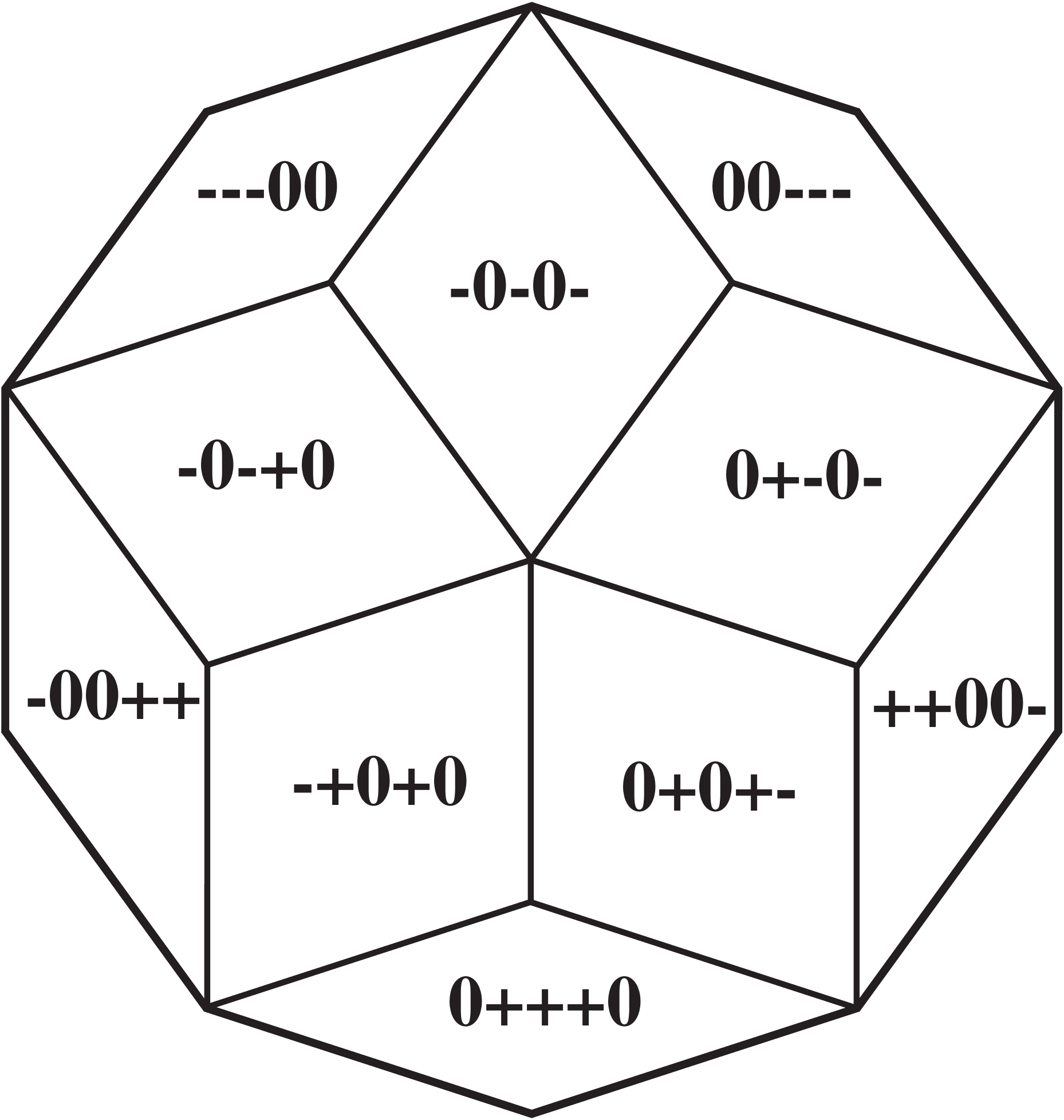}     
&     \includegraphics[width=1in]{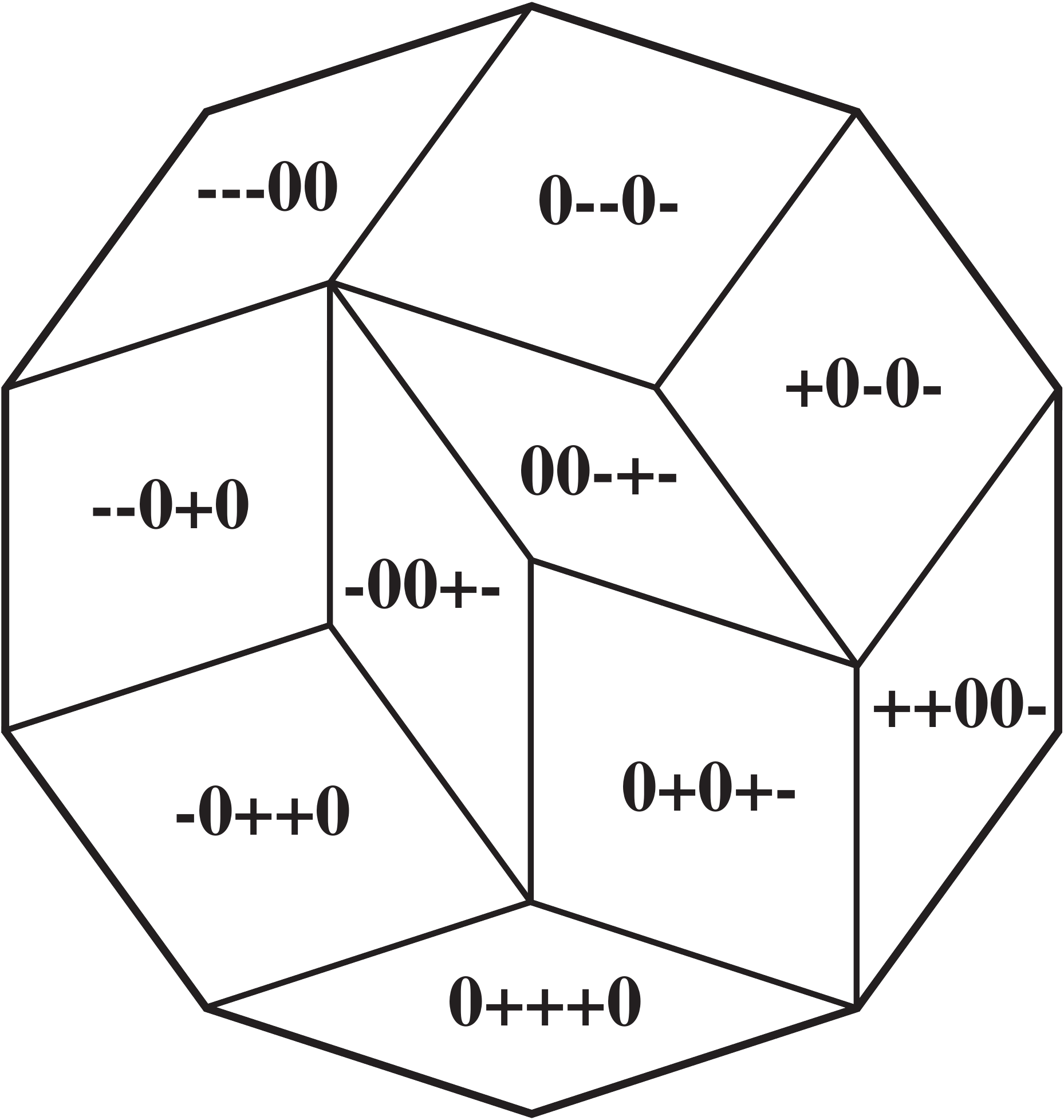}     
&     \includegraphics[width=1in]{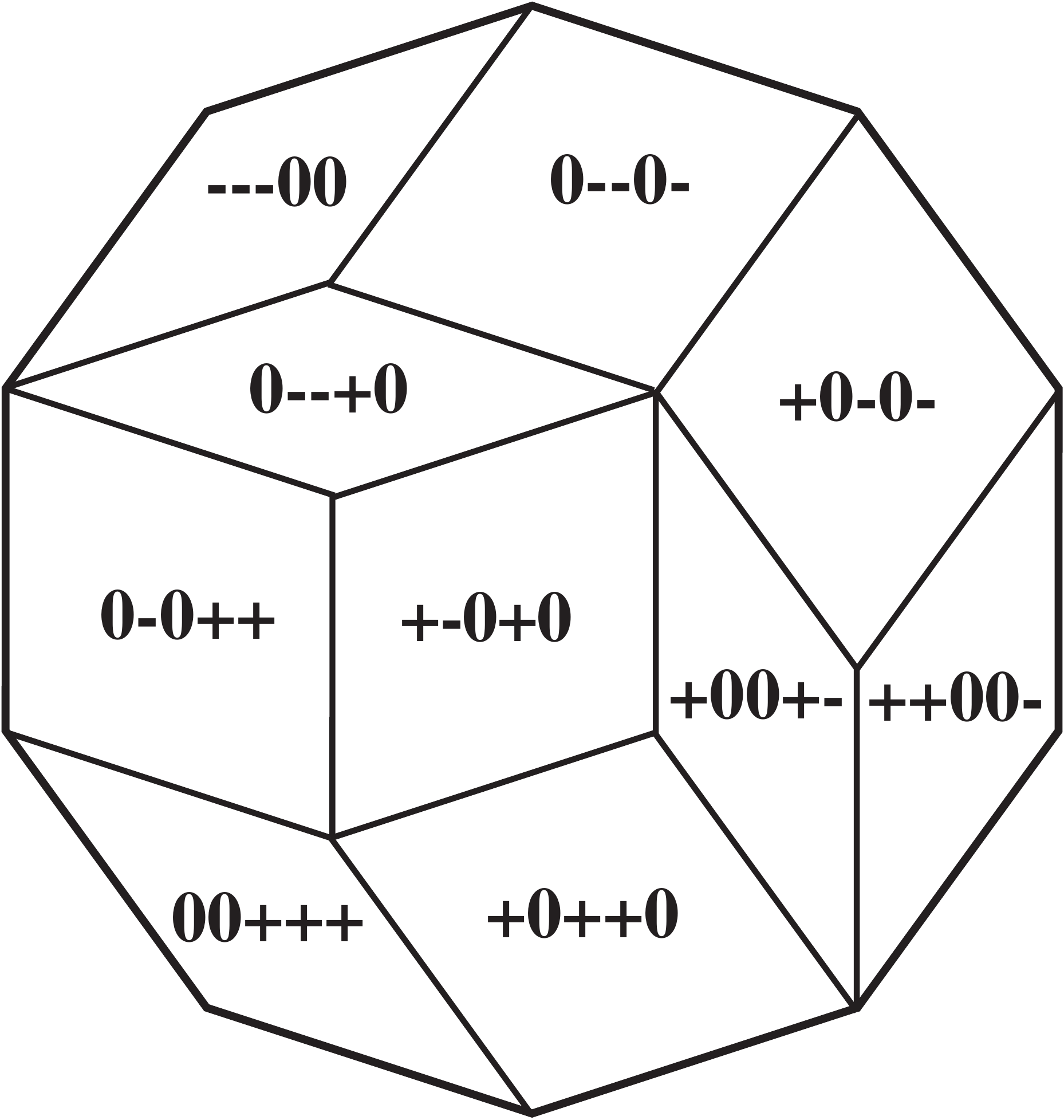} \\

 00--0 &  0+000 &  0+0+0 &  000+0 &  0--00\\

     \includegraphics[width=1in]{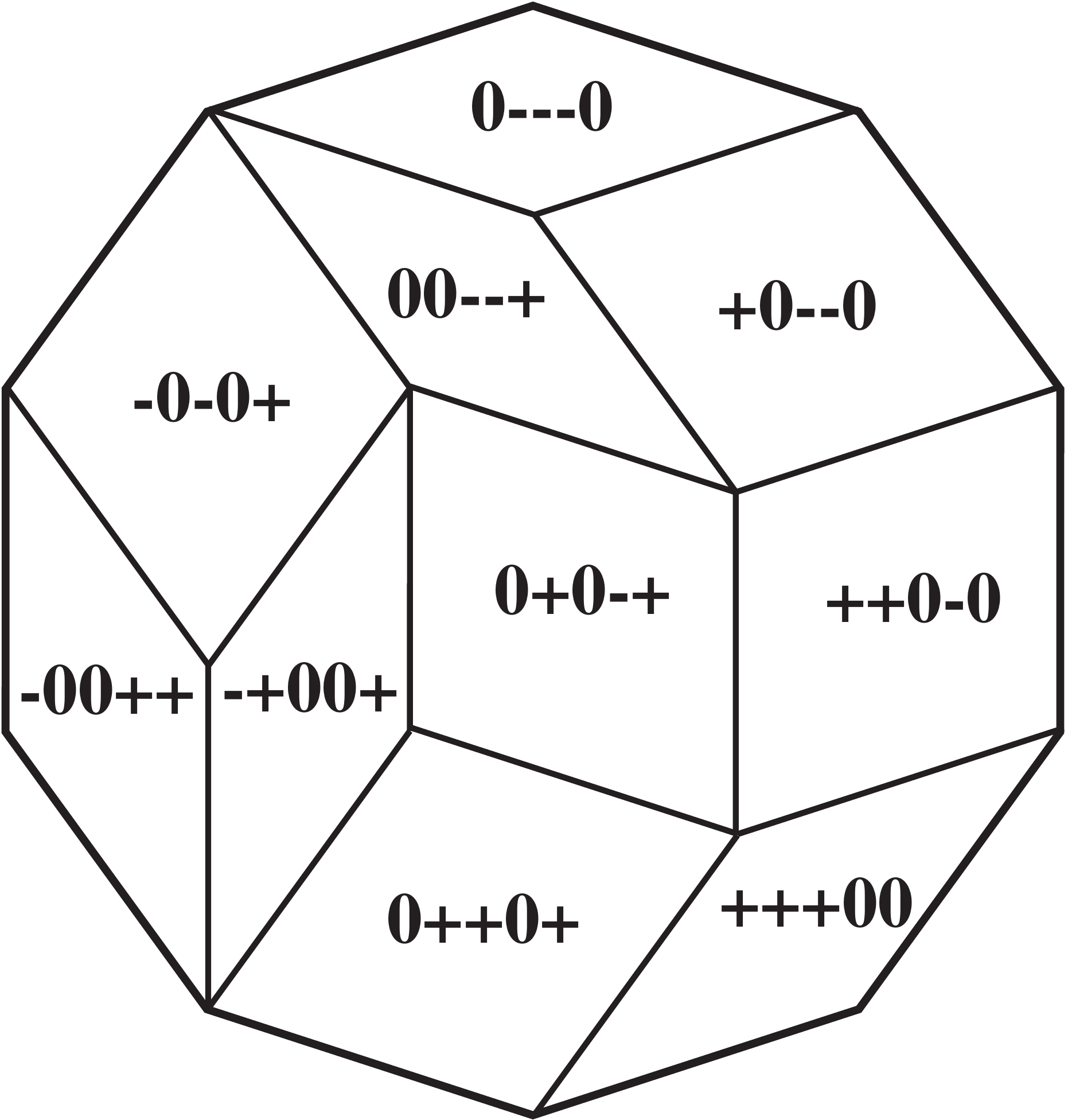}     
  &      \includegraphics[width=1in]{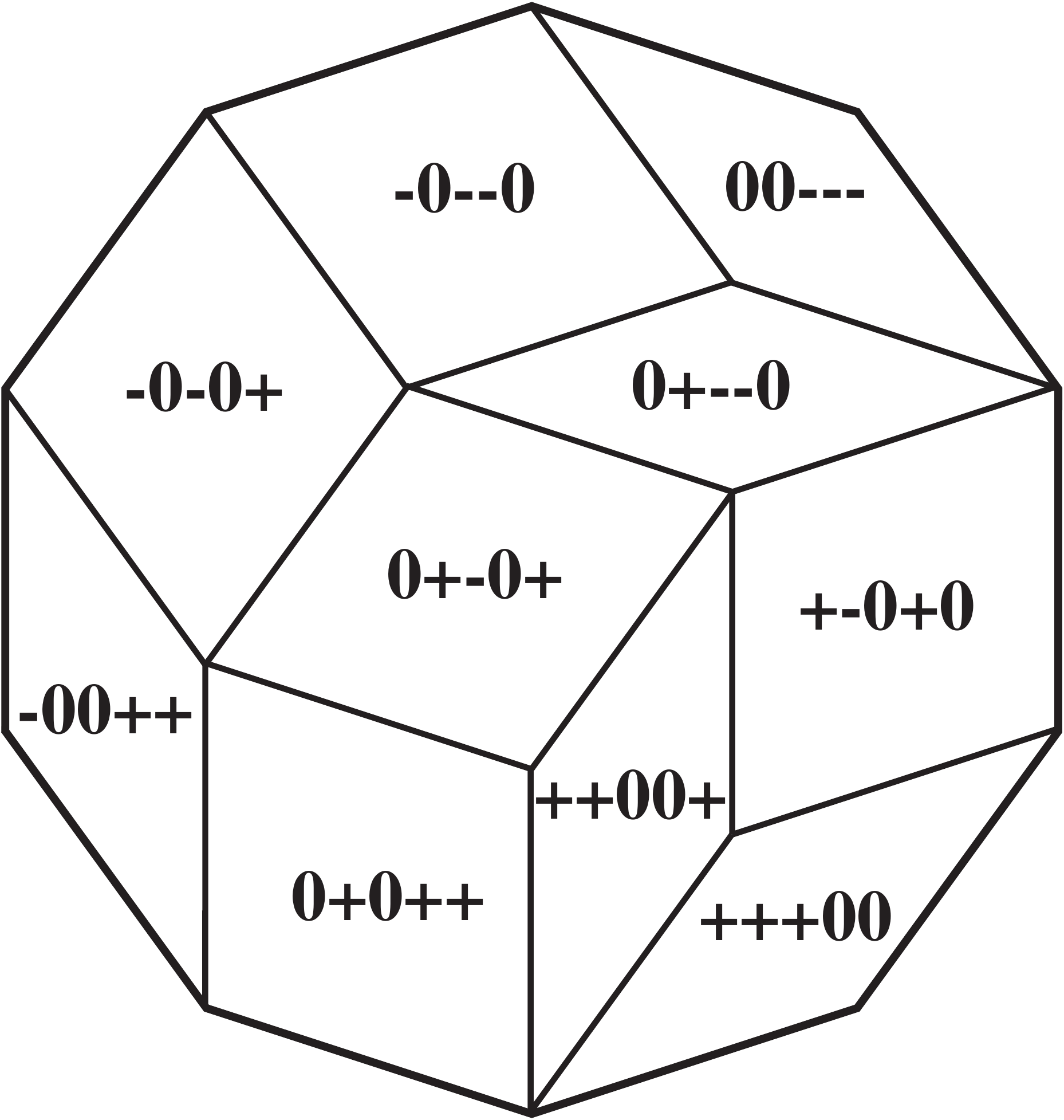}     
 &    \includegraphics[width=1in]{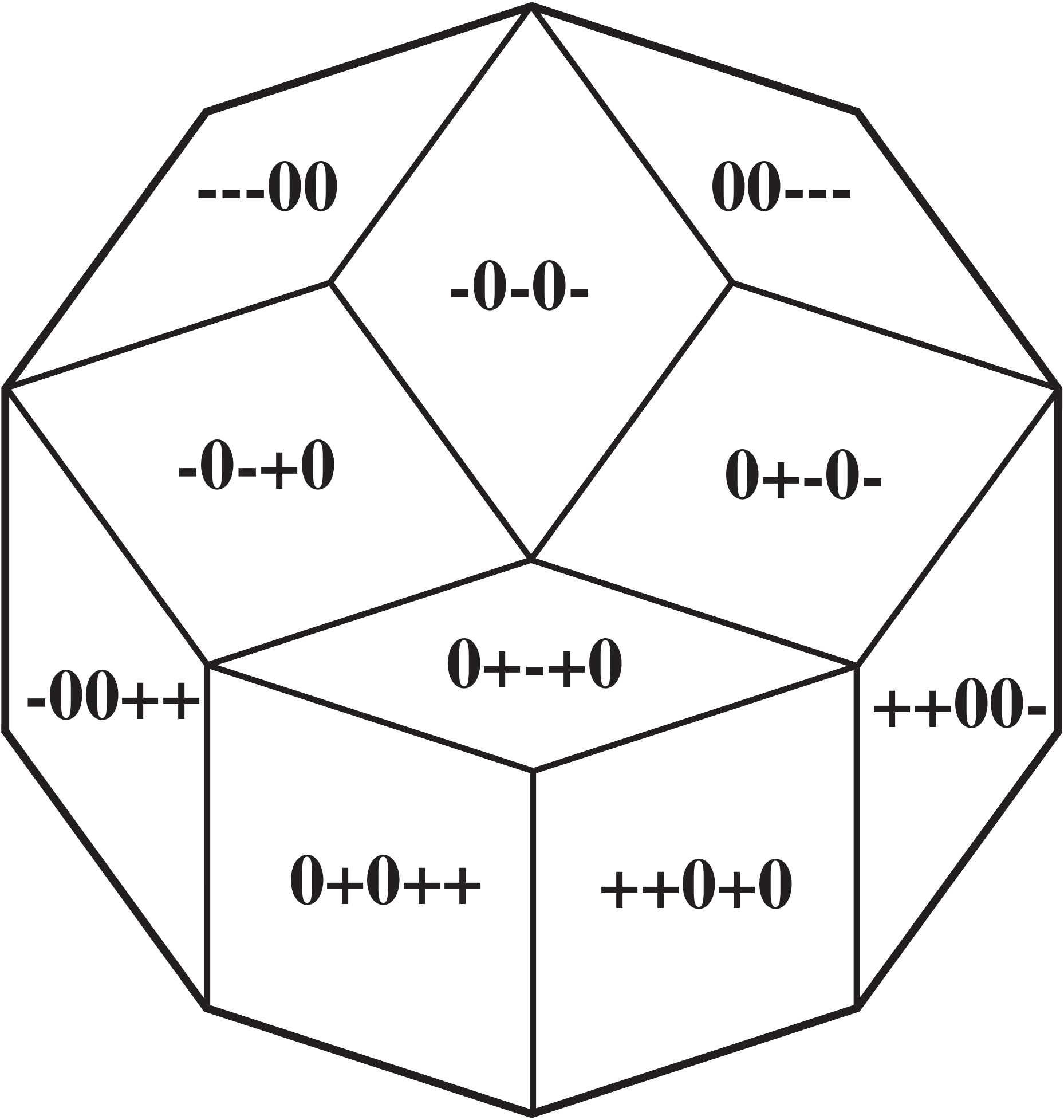}     
&     \includegraphics[width=1in]{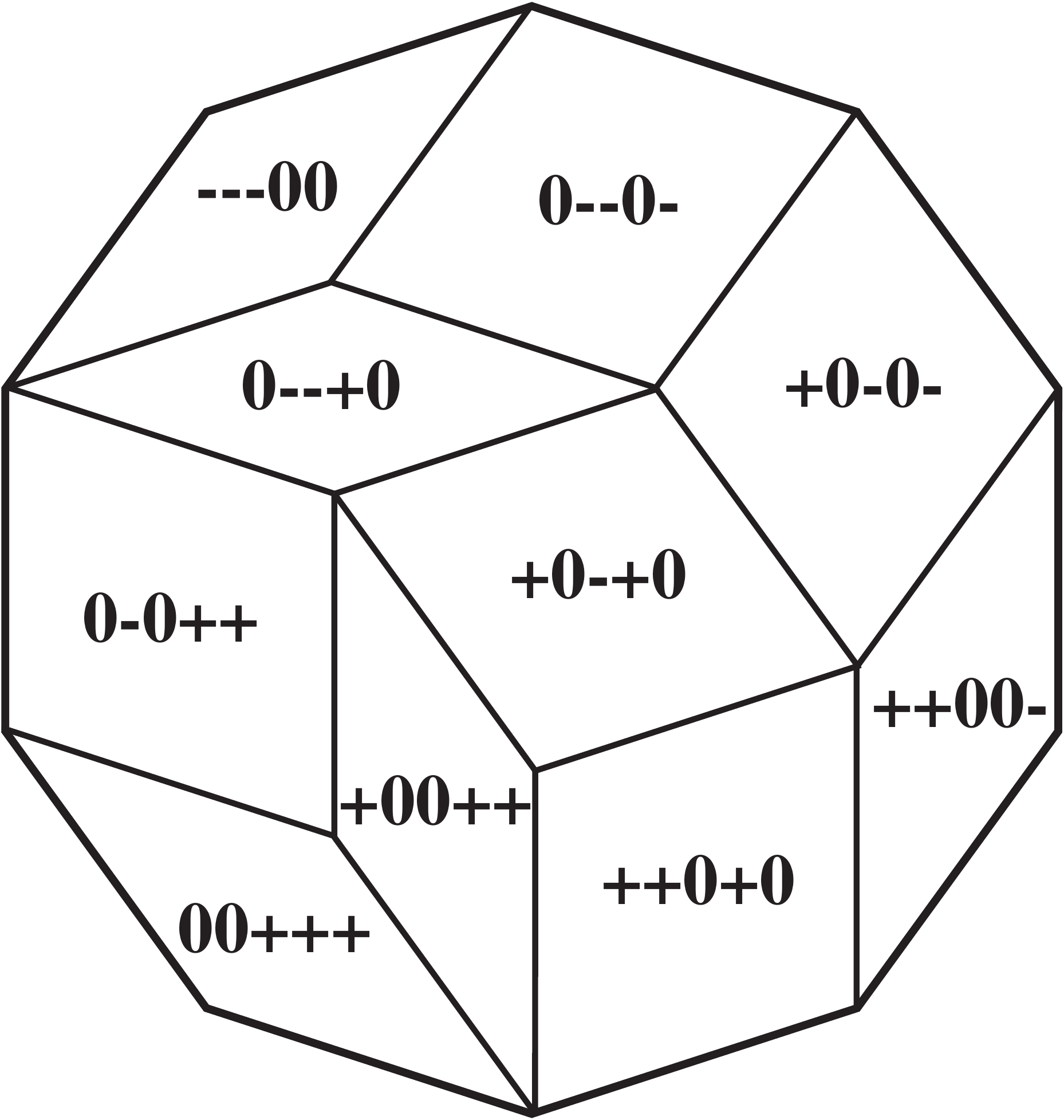}     
&     \includegraphics[width=1in]{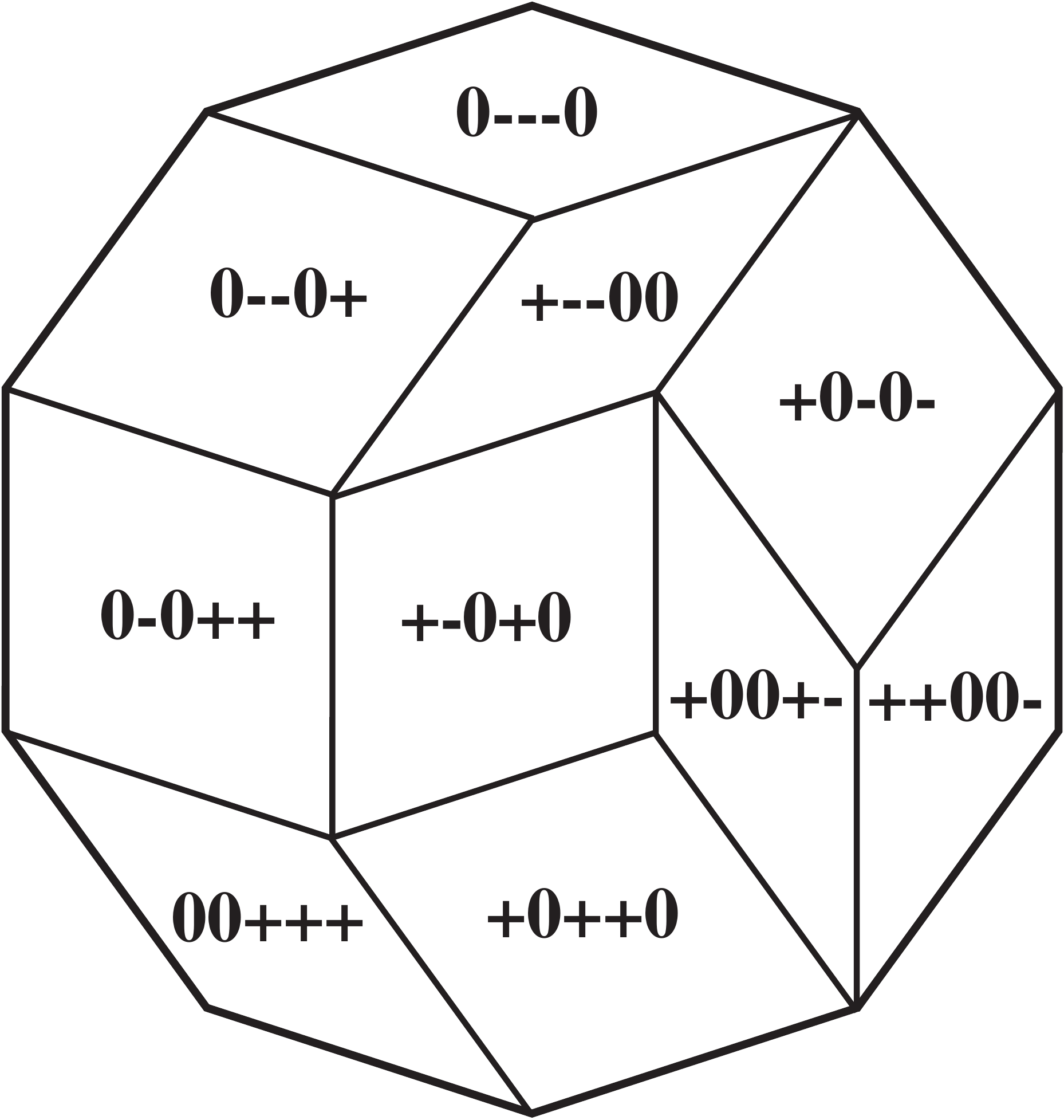} \\

 0+00+ &  00--0 &  00-00 &  0--00 &  +00+0\\

     \includegraphics[width=1in]{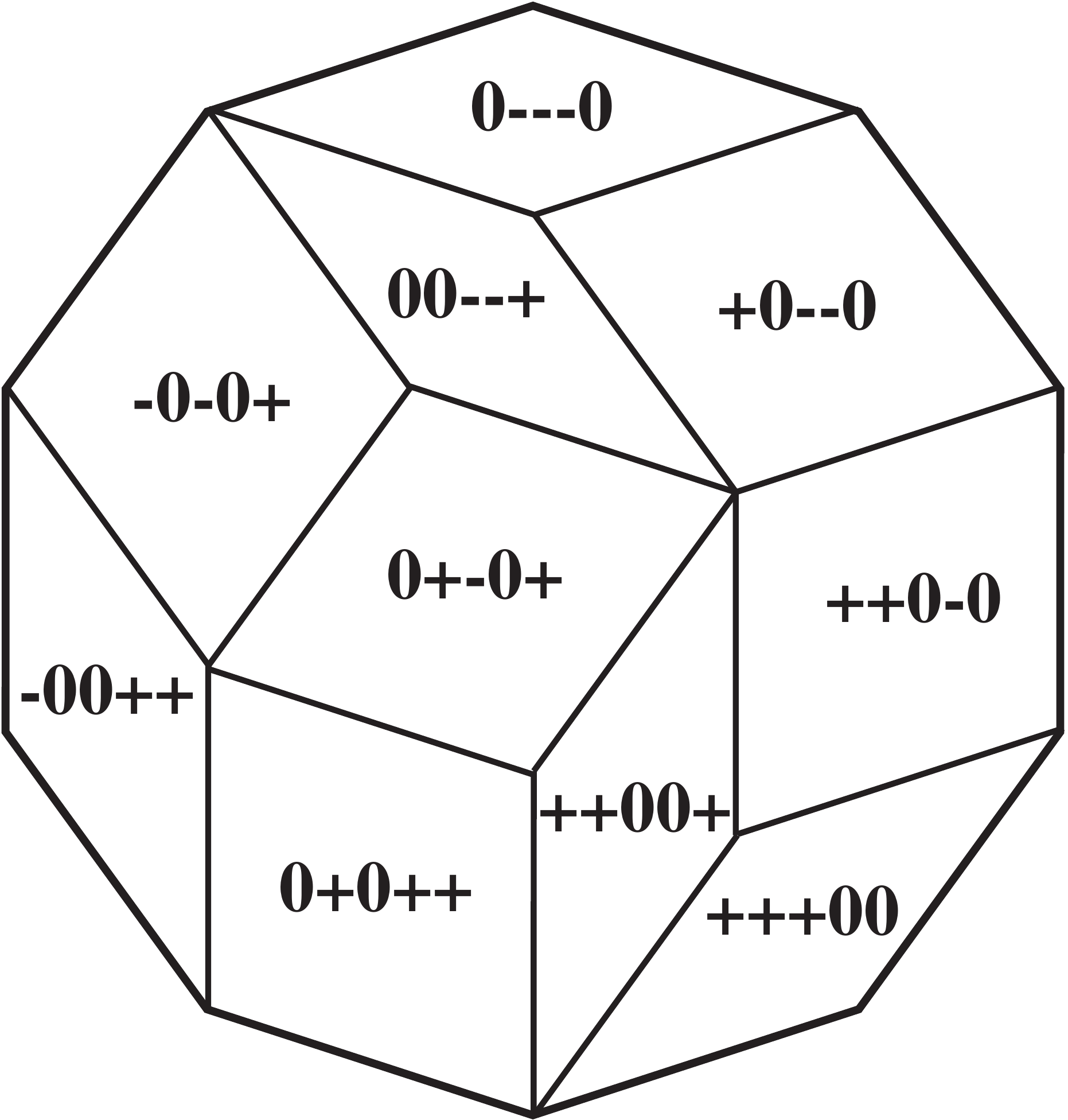}     
  &      \includegraphics[width=1in]{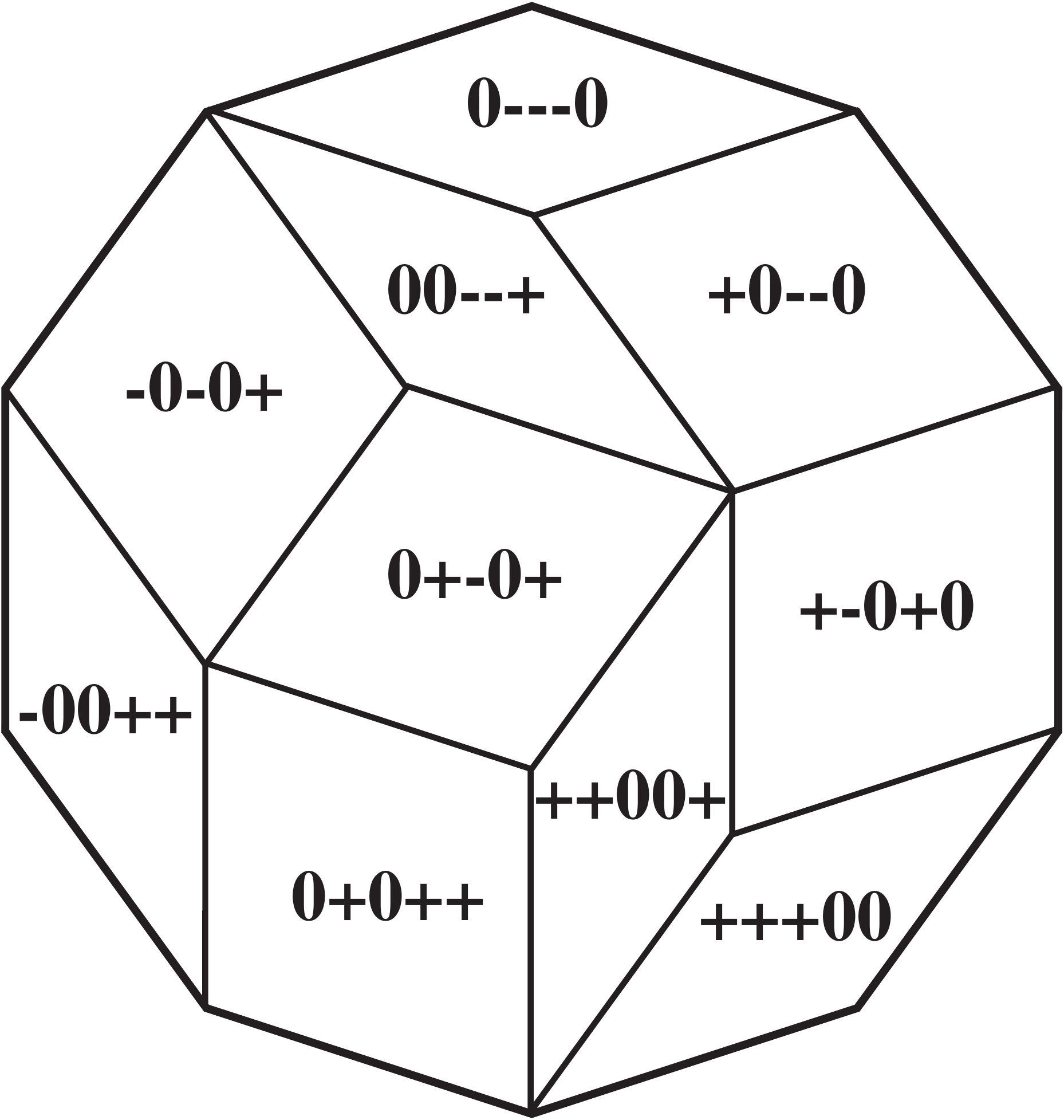}     
 &    \includegraphics[width=1in]{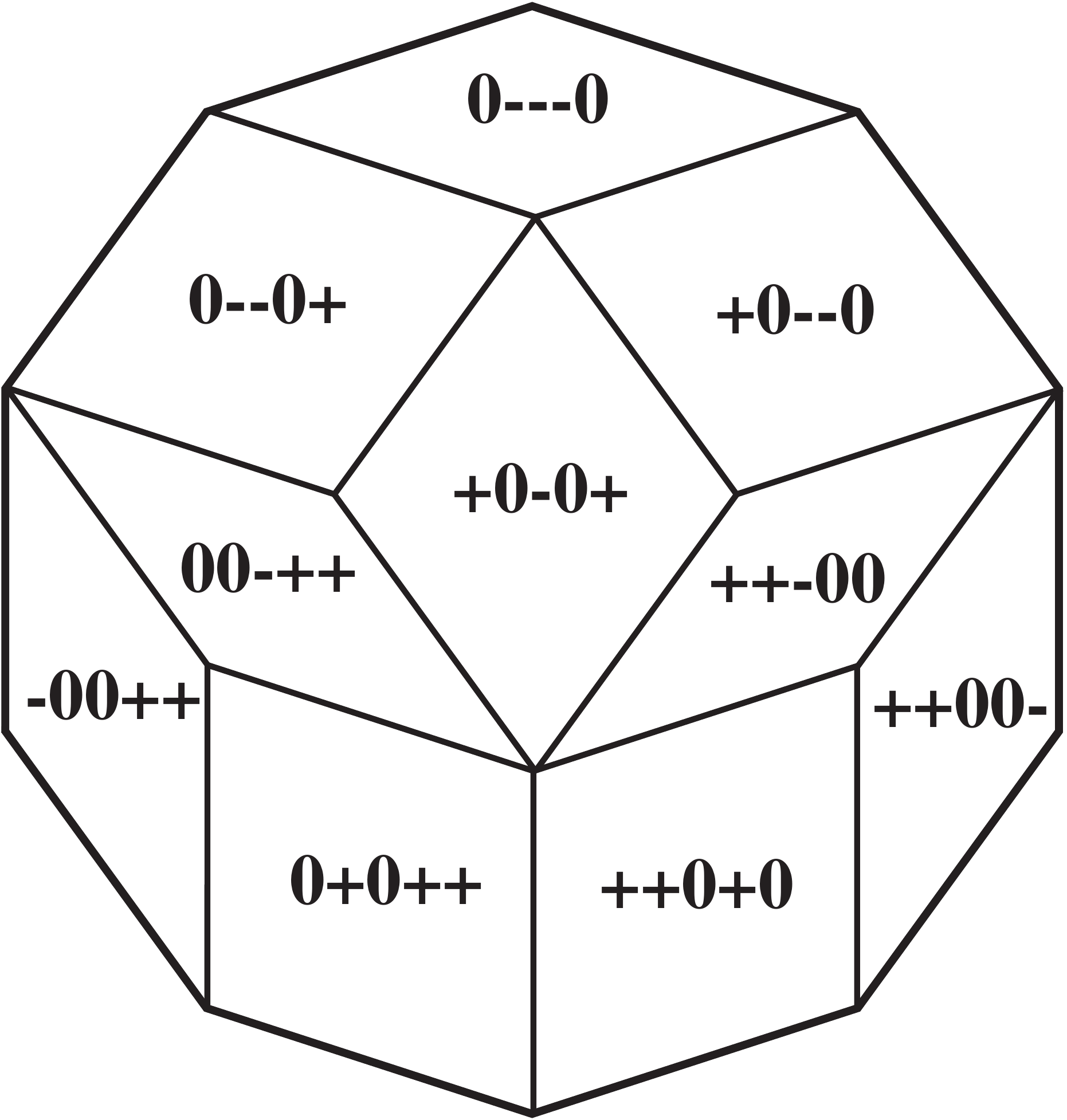}     
&     \includegraphics[width=1in]{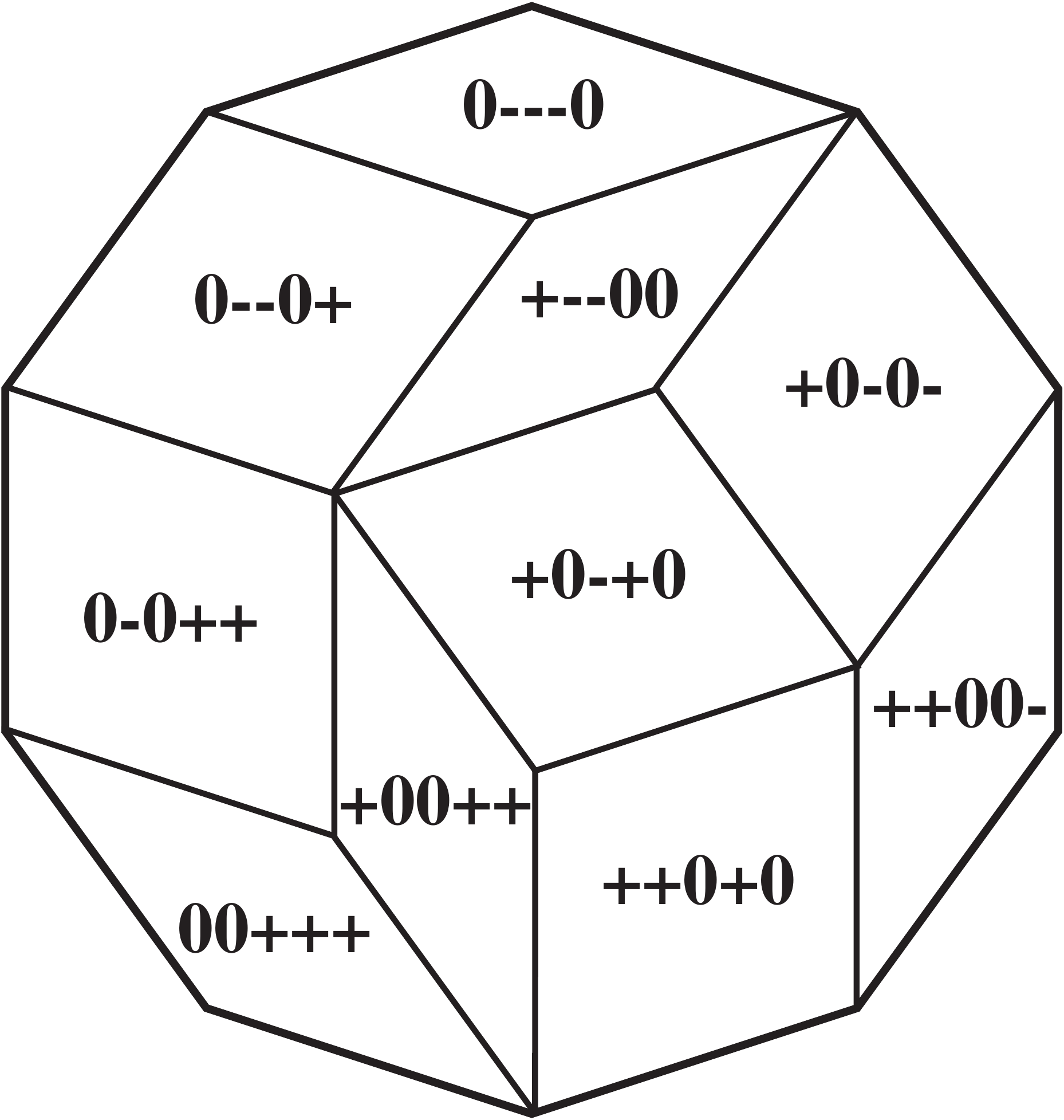}     
&     \includegraphics[width=1in]{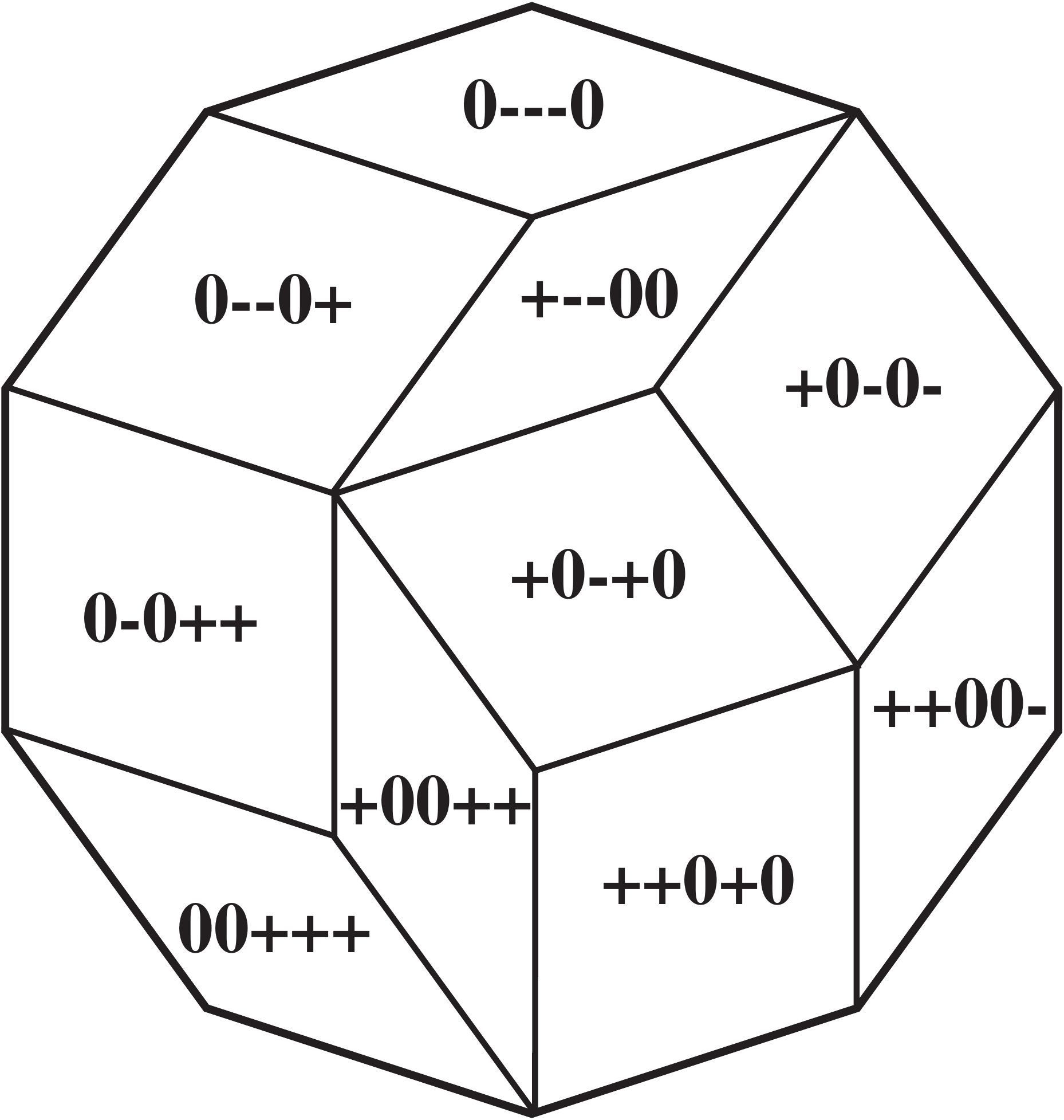} \\

 00-0+ &  00-0+ &  ++000 &  +0-00 &  +0-00\\

      \includegraphics[width=1in]{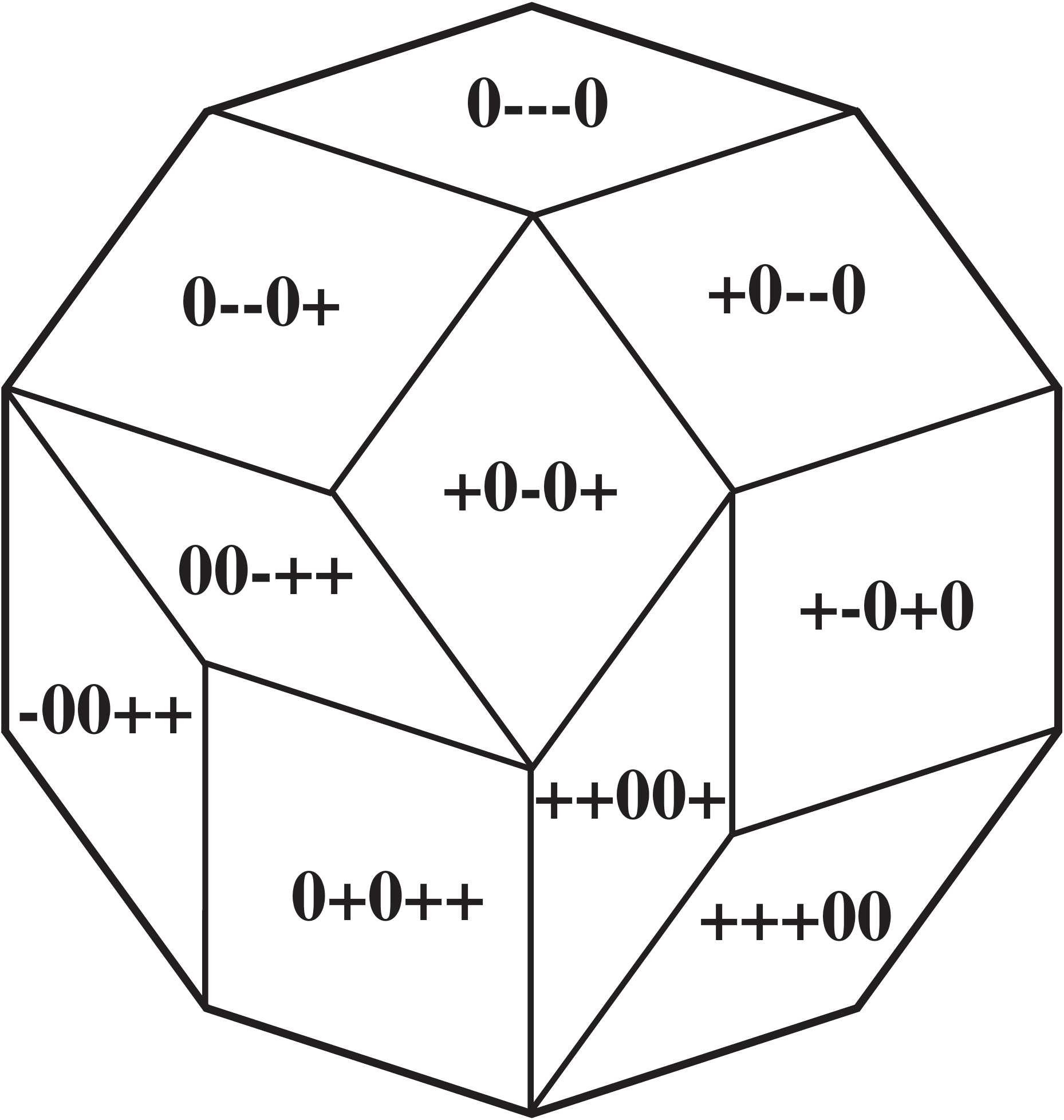}     
  &      \includegraphics[width=1in]{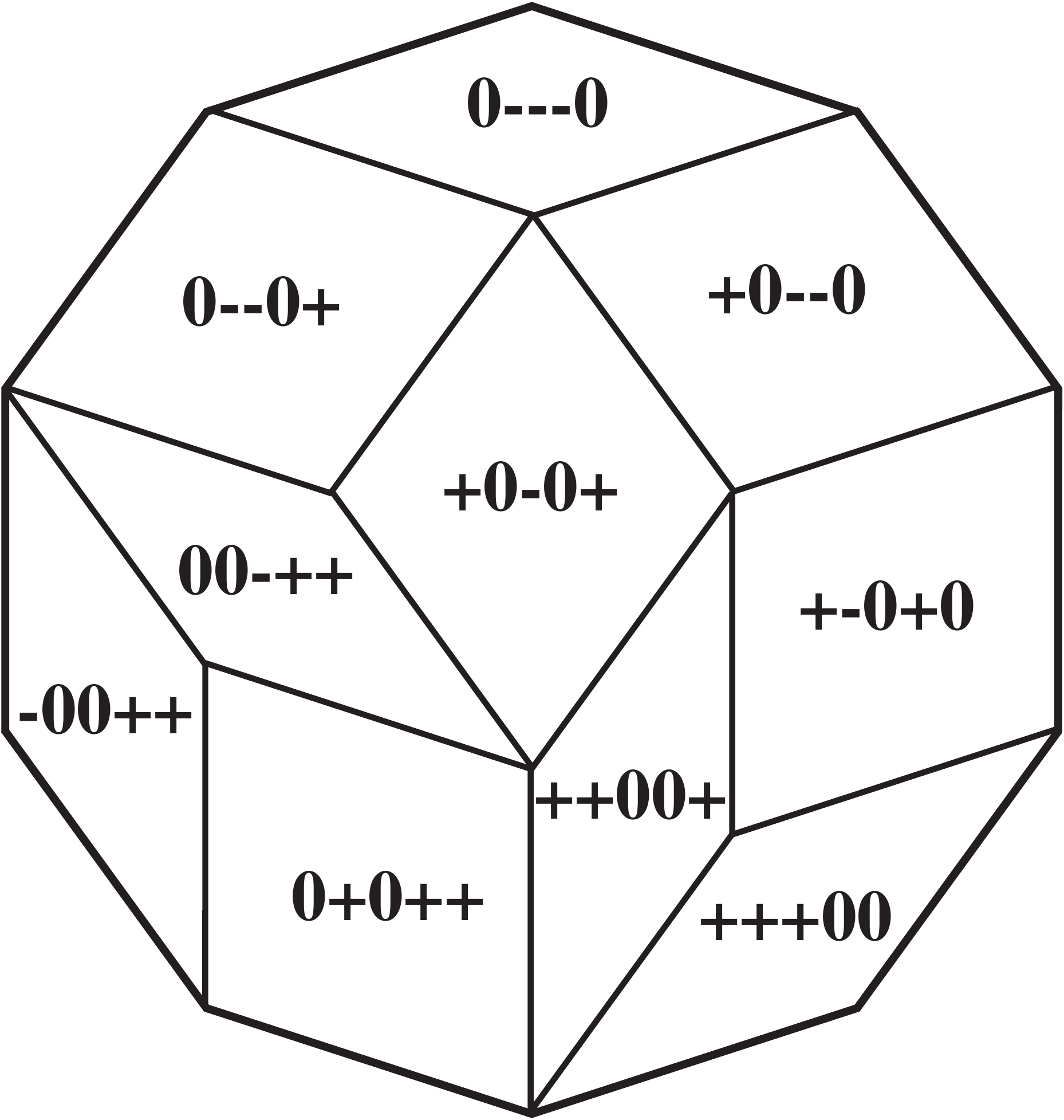}     
 &    \includegraphics[width=1in]{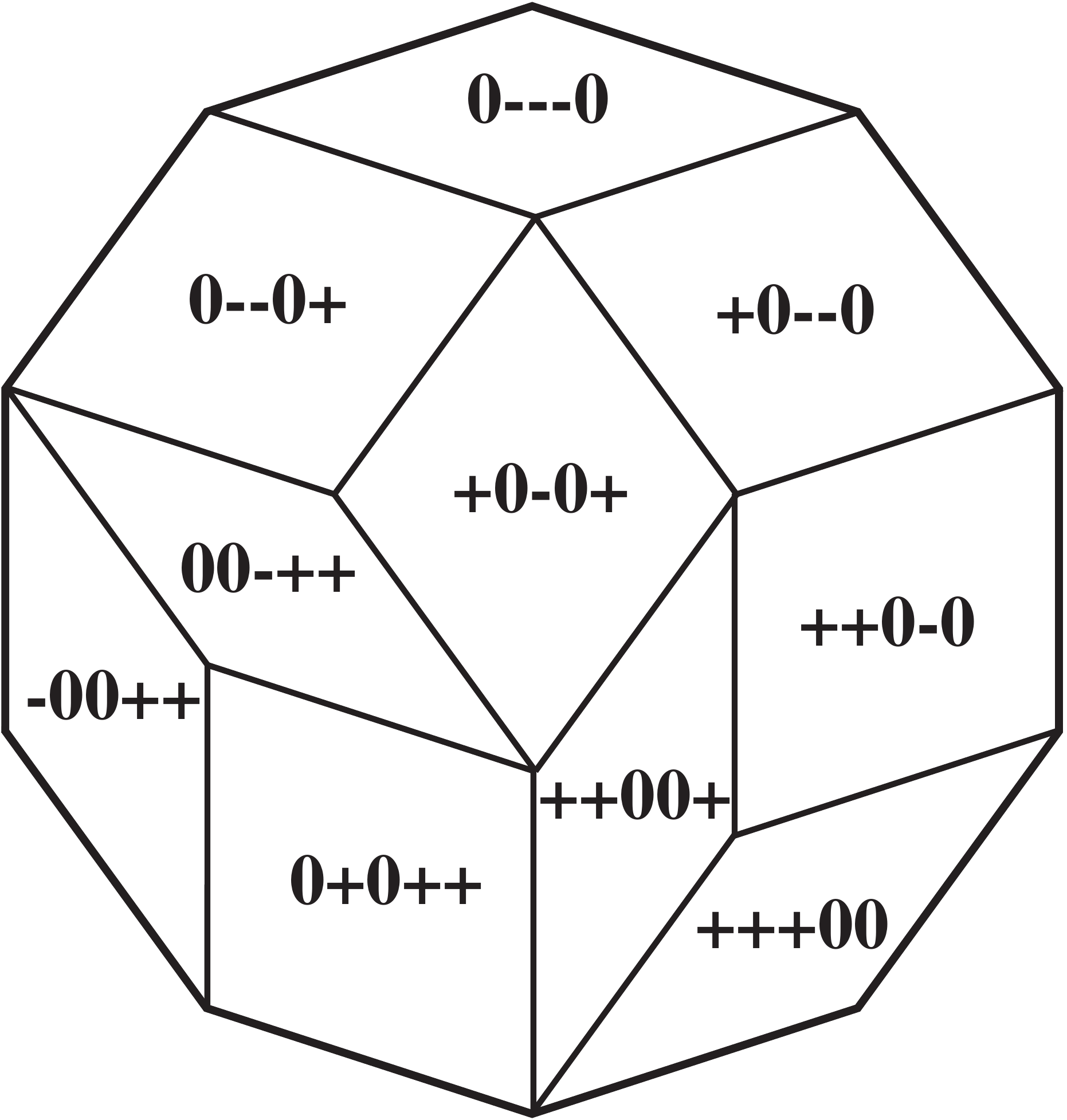}     
&     \includegraphics[width=1in]{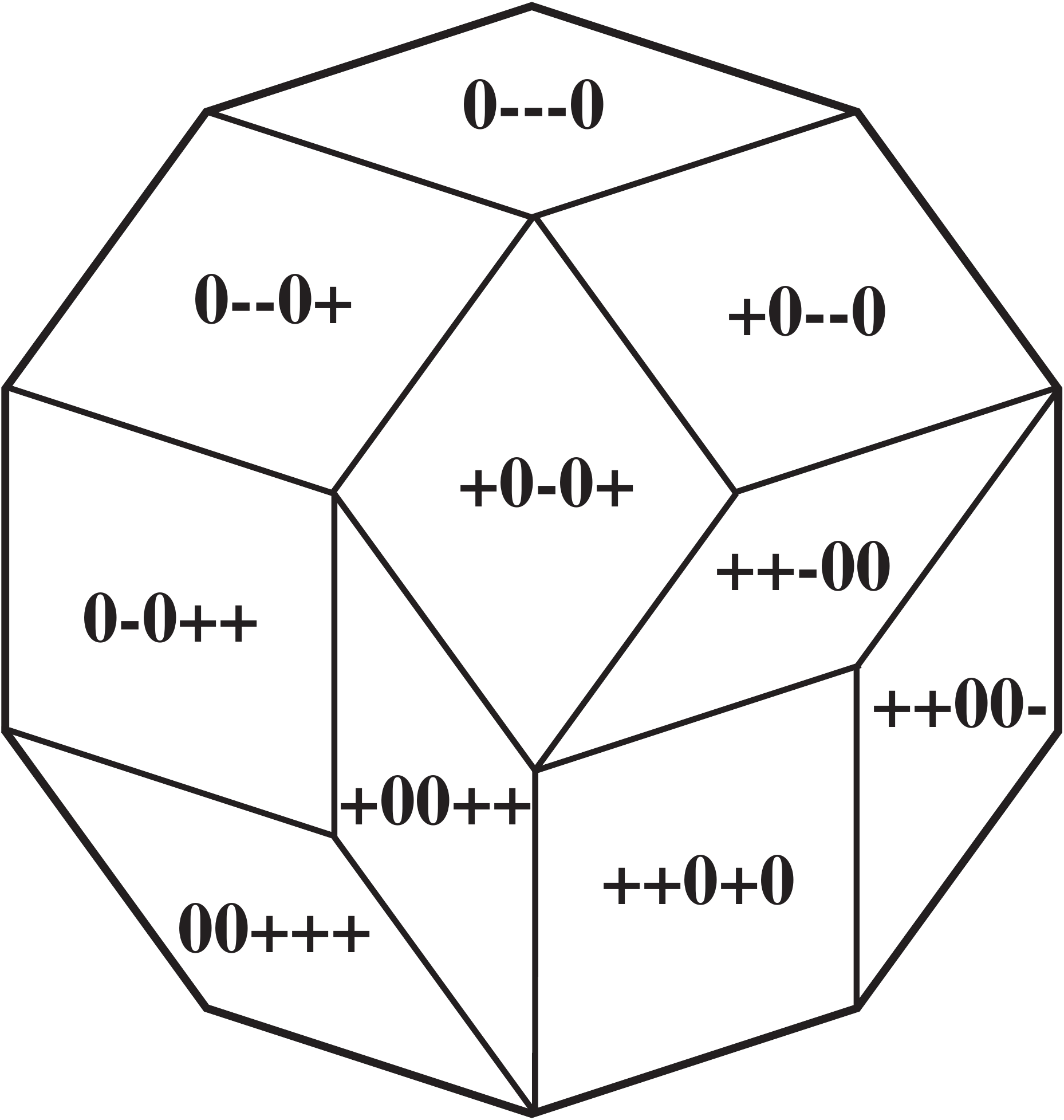}     
&     \includegraphics[width=1in]{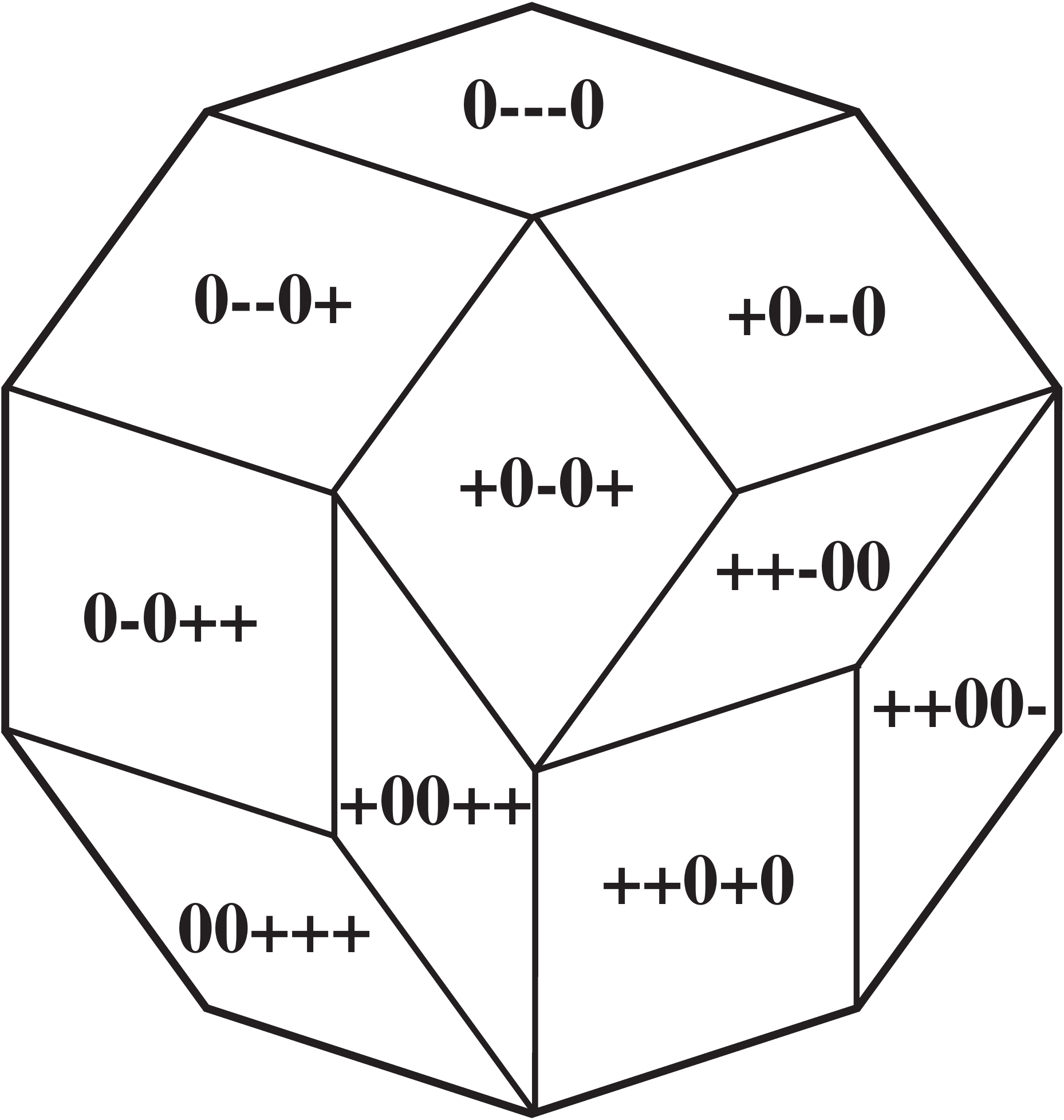} \\

 000++ &  000++ &  000++ &  ++000 &  ++000 \\
 %&& {Figure 19: 5-cube source}&&
 %6-0  &  6-1  &  6-2  &  6-3  &  6-4
 \end{array}
%\right)
\]
\vskip -0.5cm
\begin{figure} \includegraphics[width=0in]{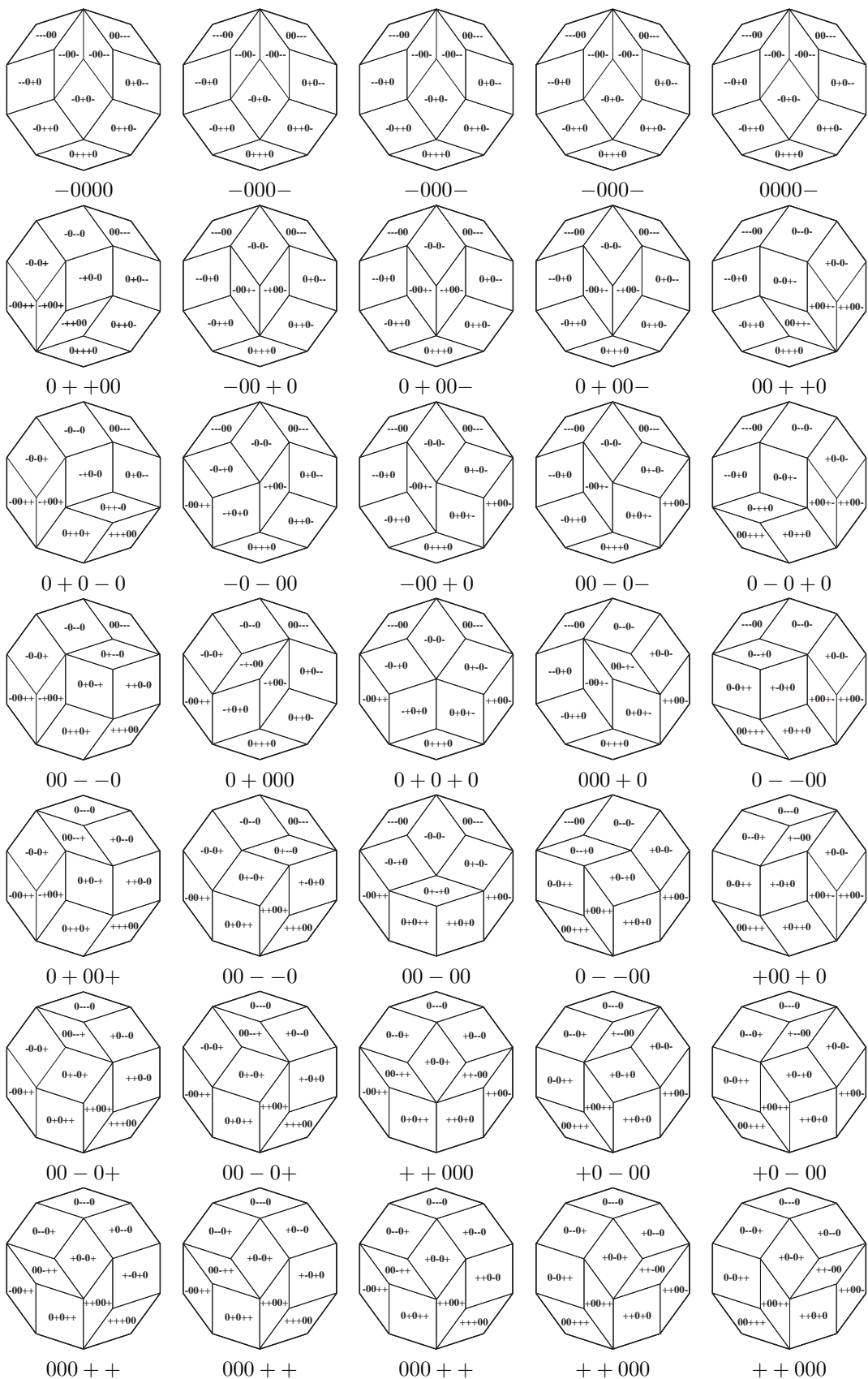}
\vskip -1cm\caption{5-cube source: and target next page}
\label{fig:19}
\end{figure}
 }
%\pagebreak 
 
 \vskip 0.25in
 \centerline{\bf 5-cube target}
  \[
%\left(
\begin{array}{ccccc}  
      \includegraphics[width=1in]{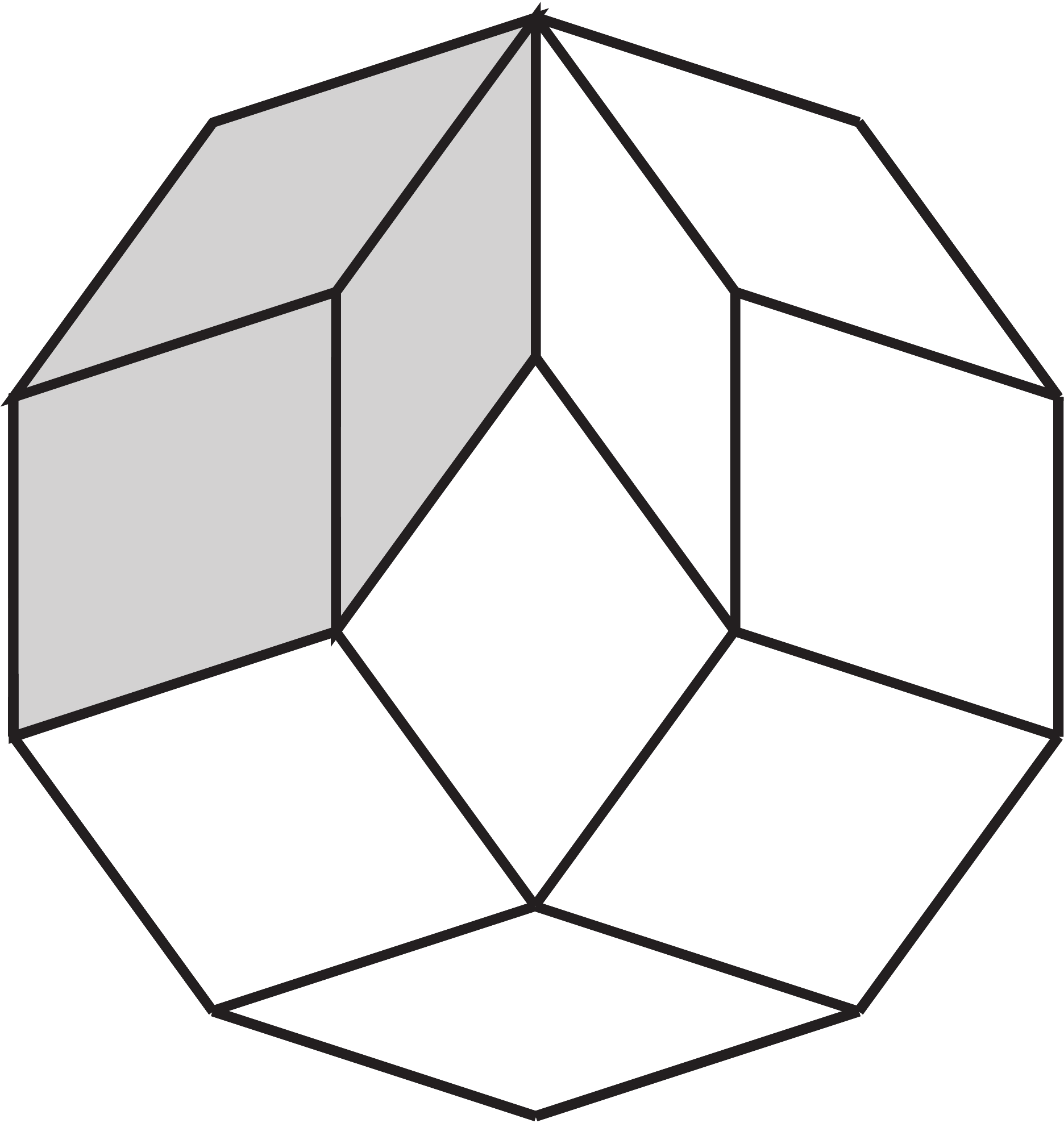}      
  &      \includegraphics[width=1in]{r11}     
 &    \includegraphics[width=1in]{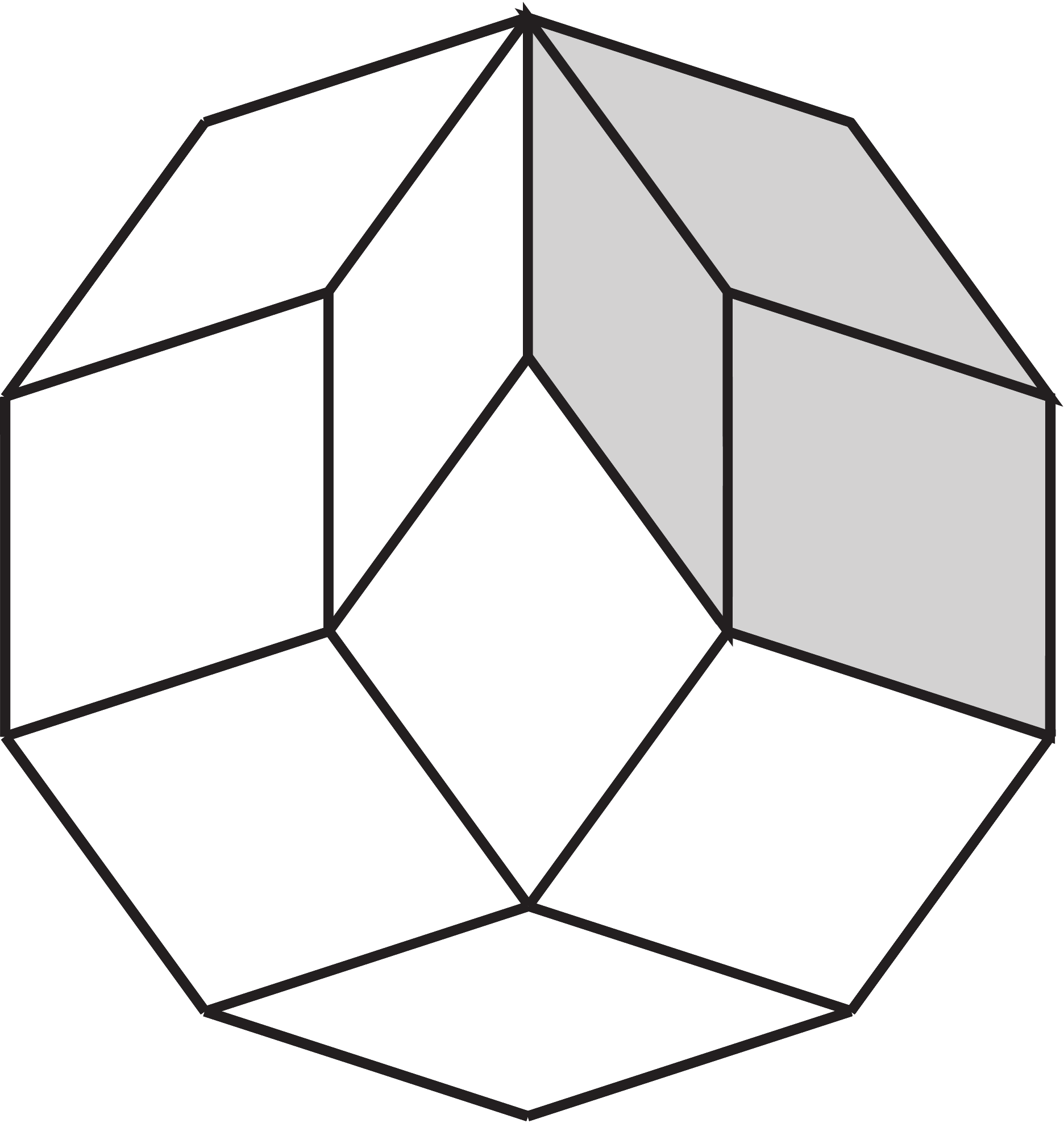}     
&     \includegraphics[width=1in]{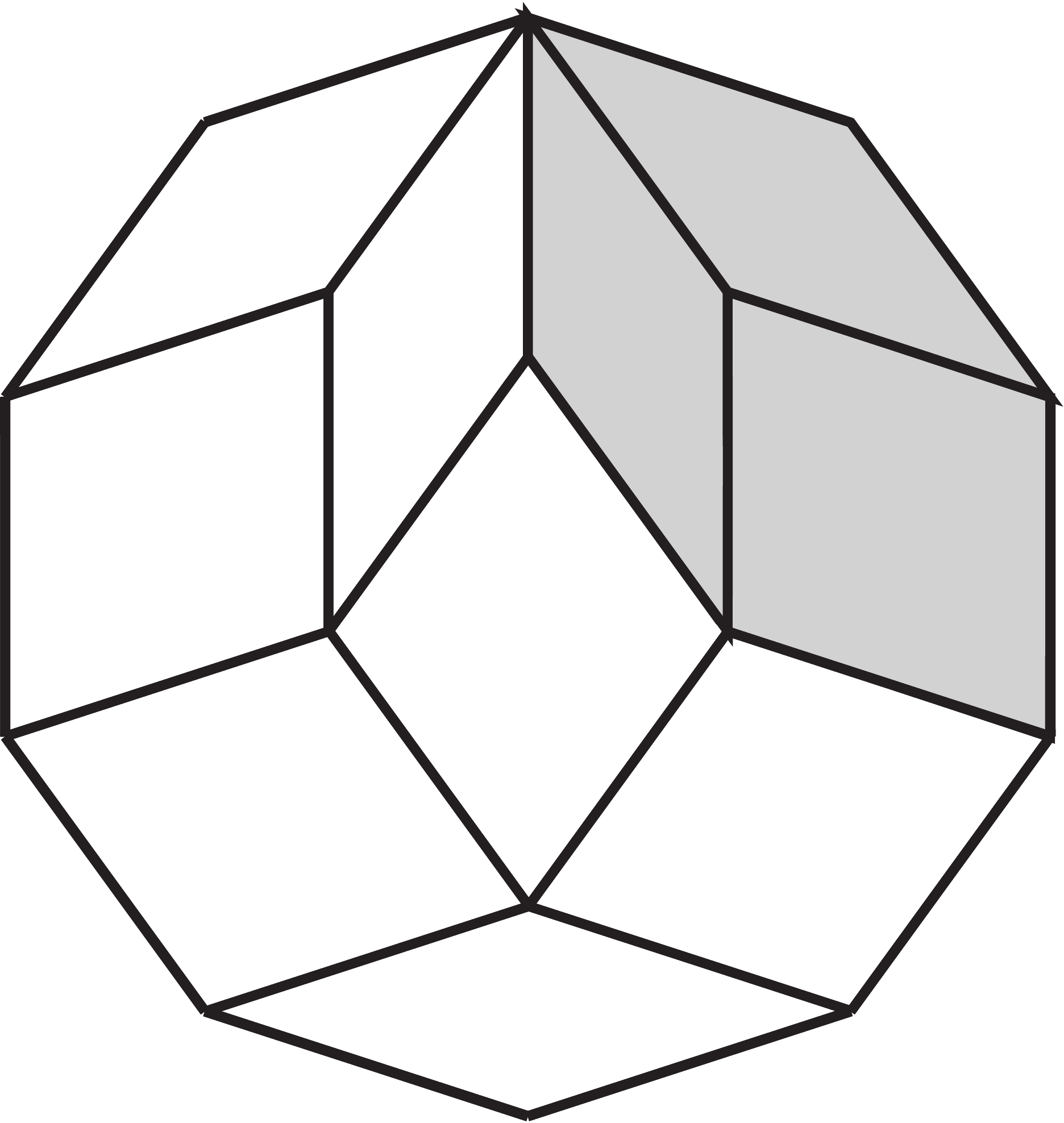}     
&     \includegraphics[width=1in]{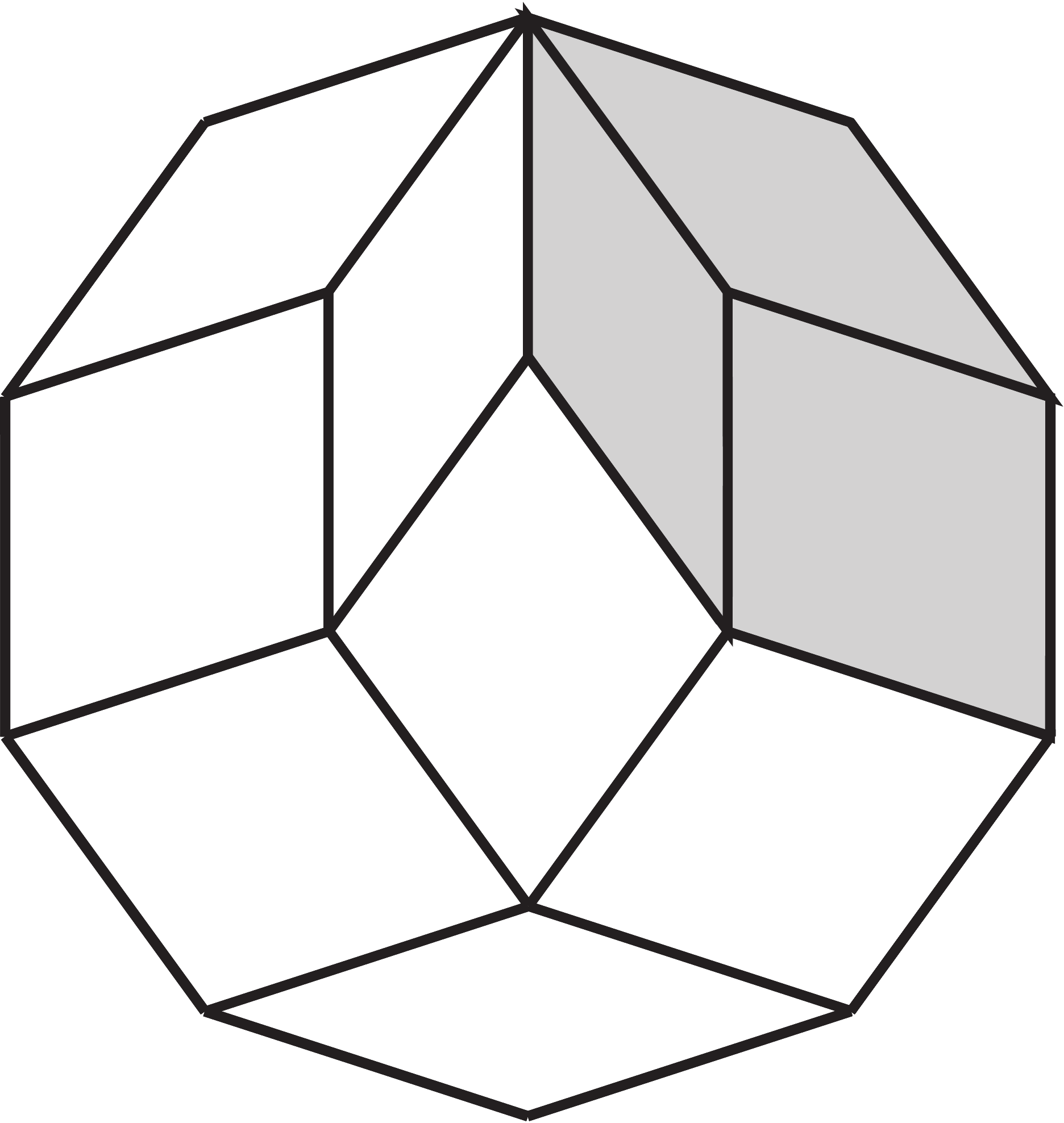} \\

 --000 & --00-  & 000-- & 000-- &000--\\

     \includegraphics[width=1in]{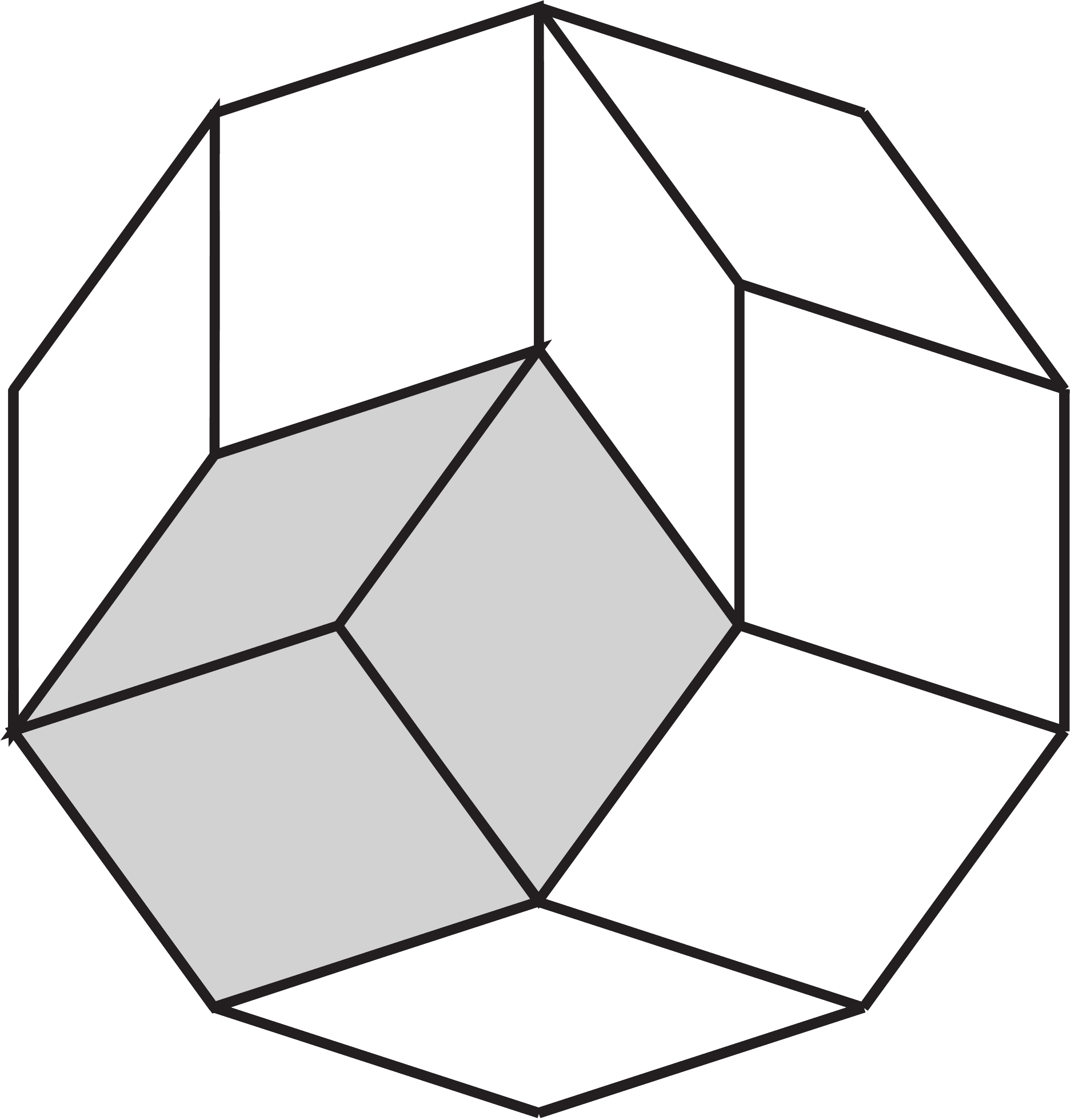}     
  &      \includegraphics[width=1in]{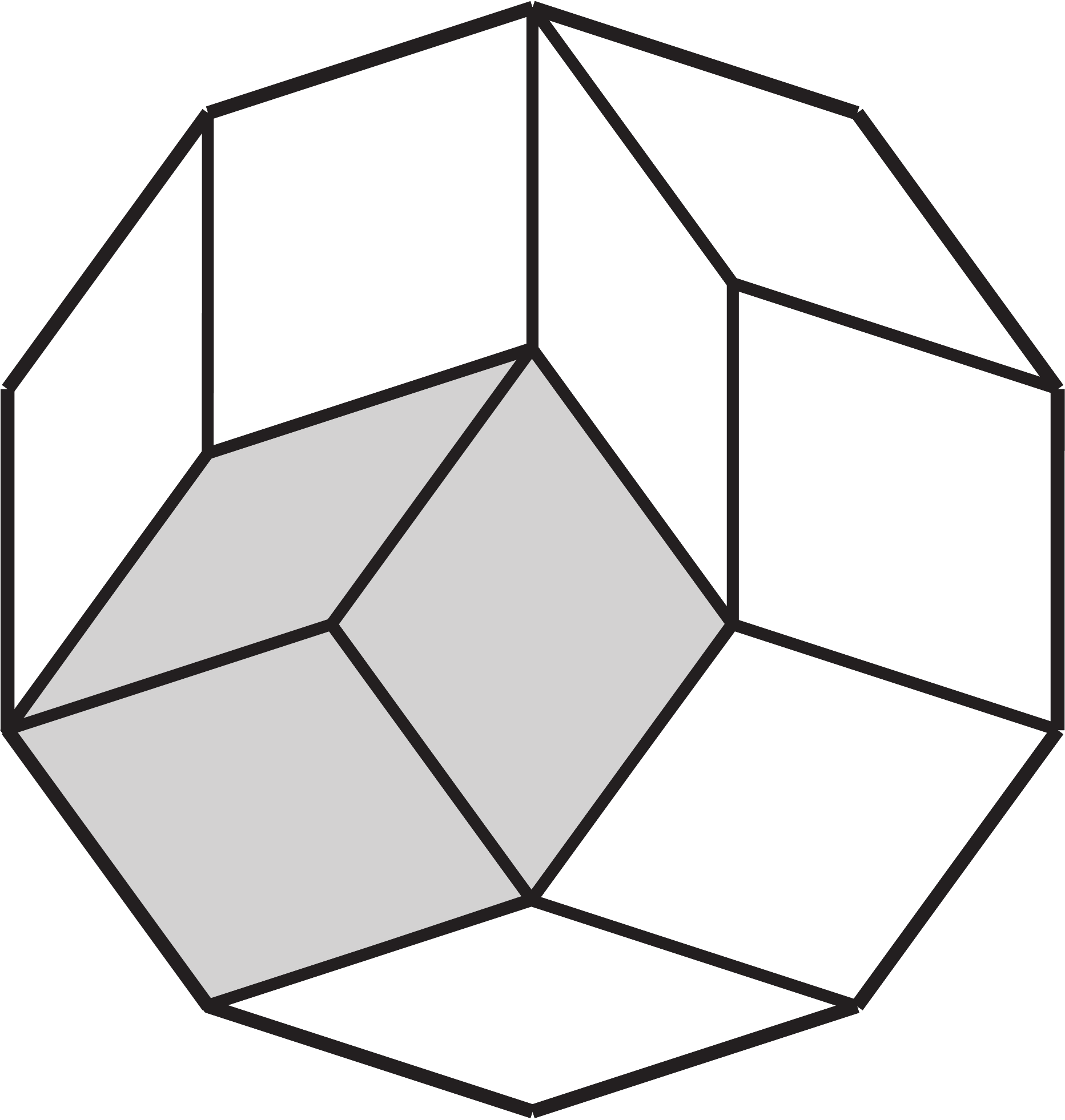}     
 &    \includegraphics[width=1in]{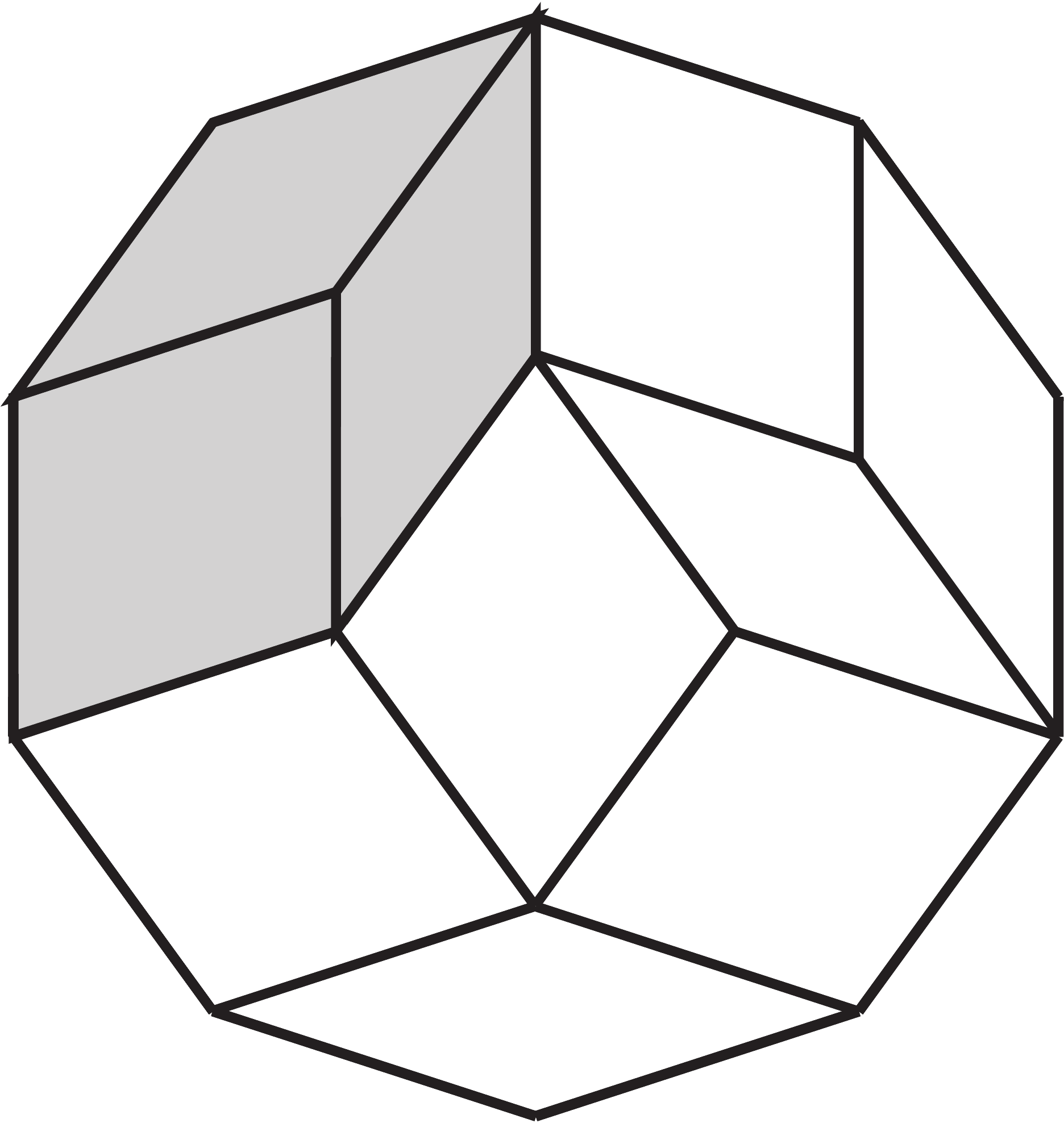}     
&     \includegraphics[width=1in]{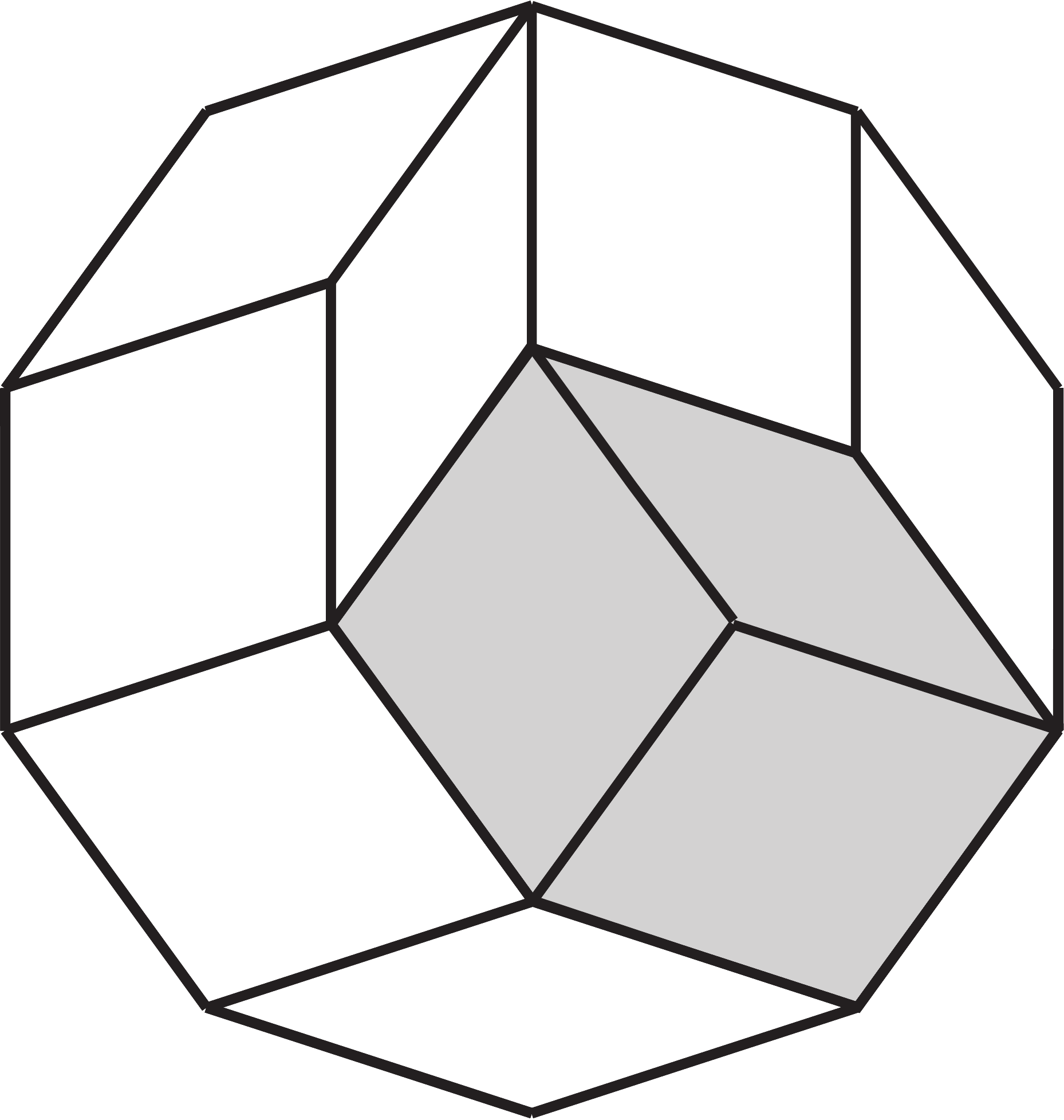}     
&     \includegraphics[width=1in]{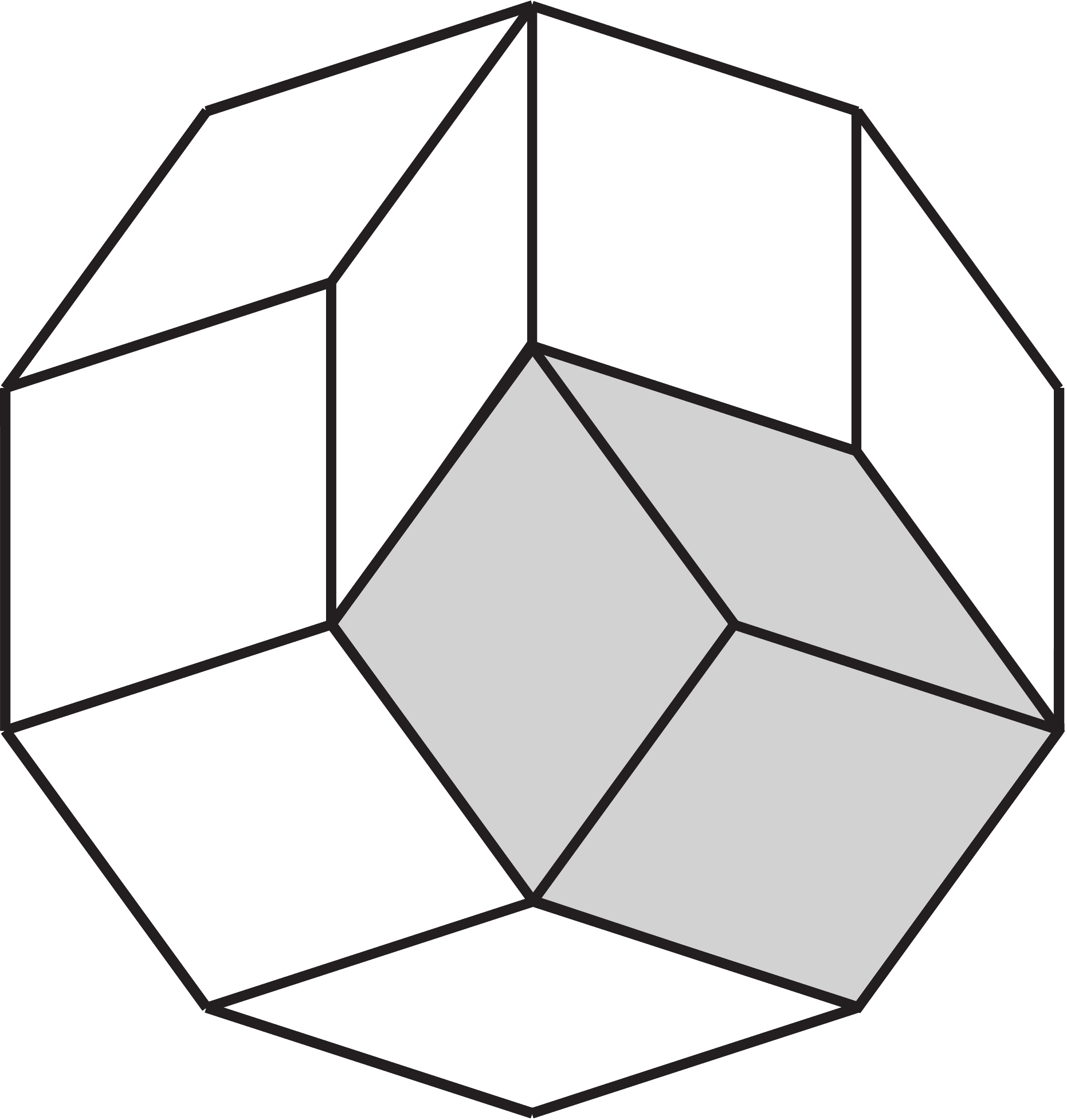} \\

-0+00 &-0+00 & --000 & 00+0-   &  00+0-\\

      \includegraphics[width=1in]{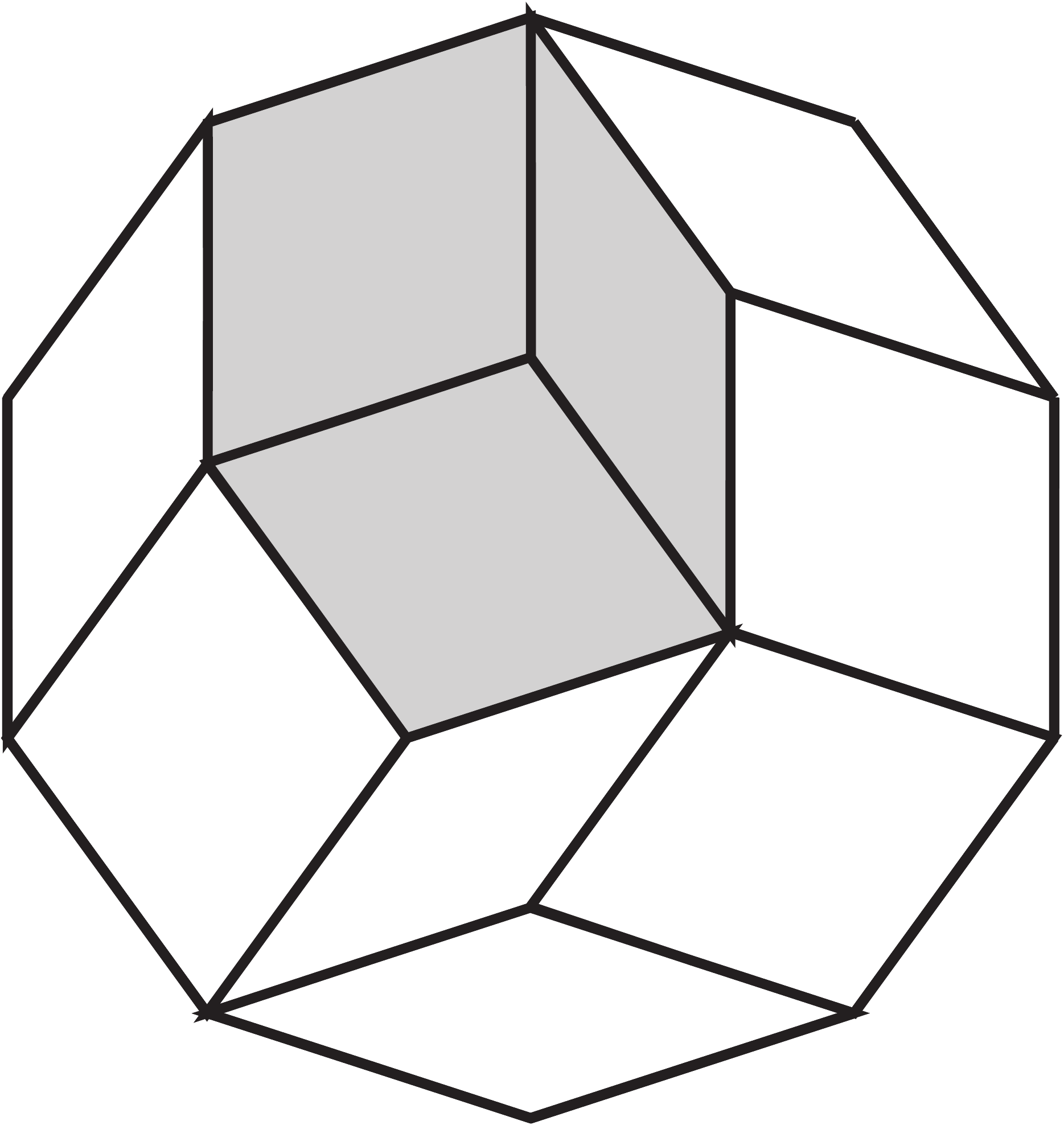}     
  &      \includegraphics[width=1in]{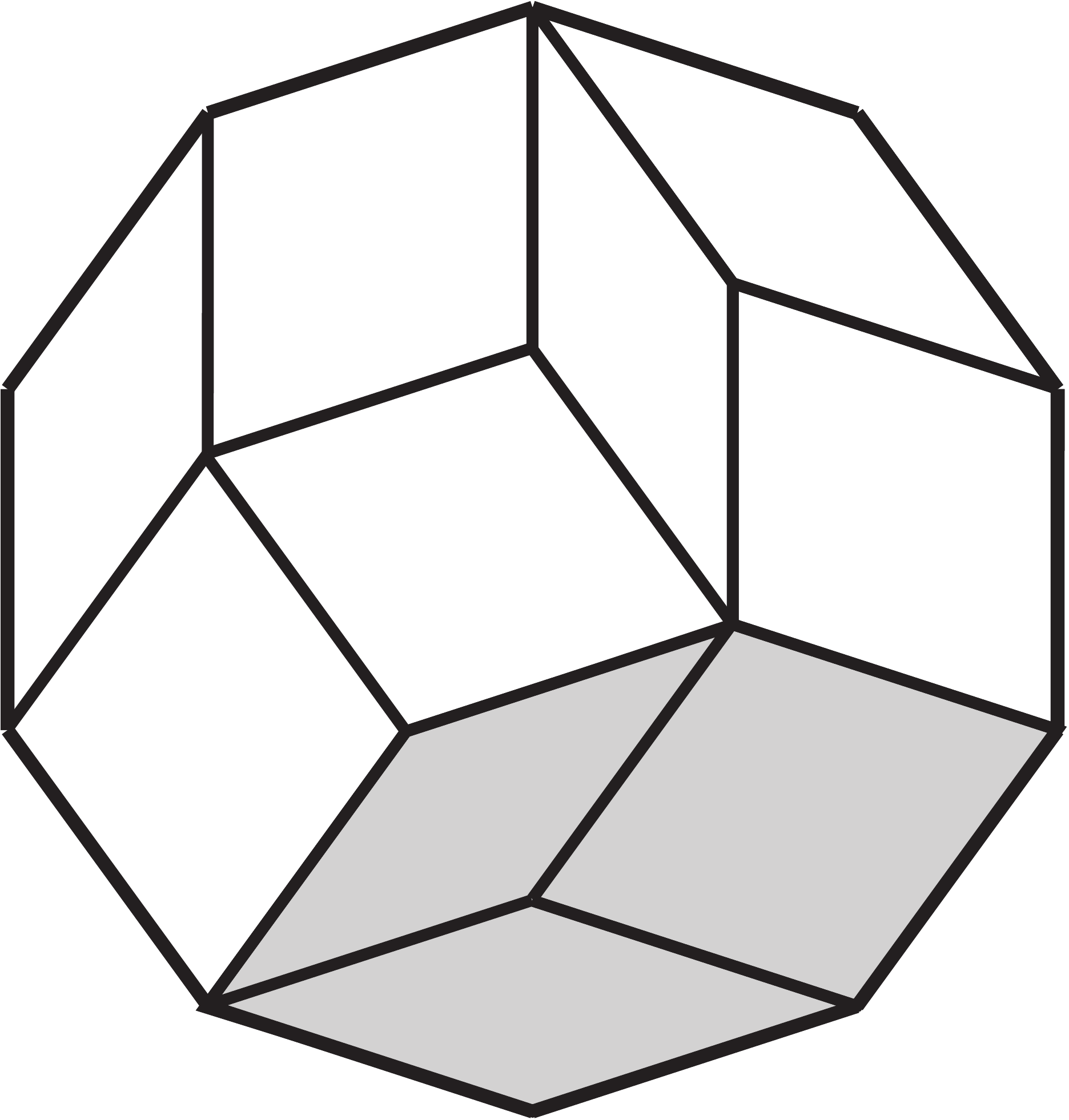}     
&    \includegraphics[width=1in]{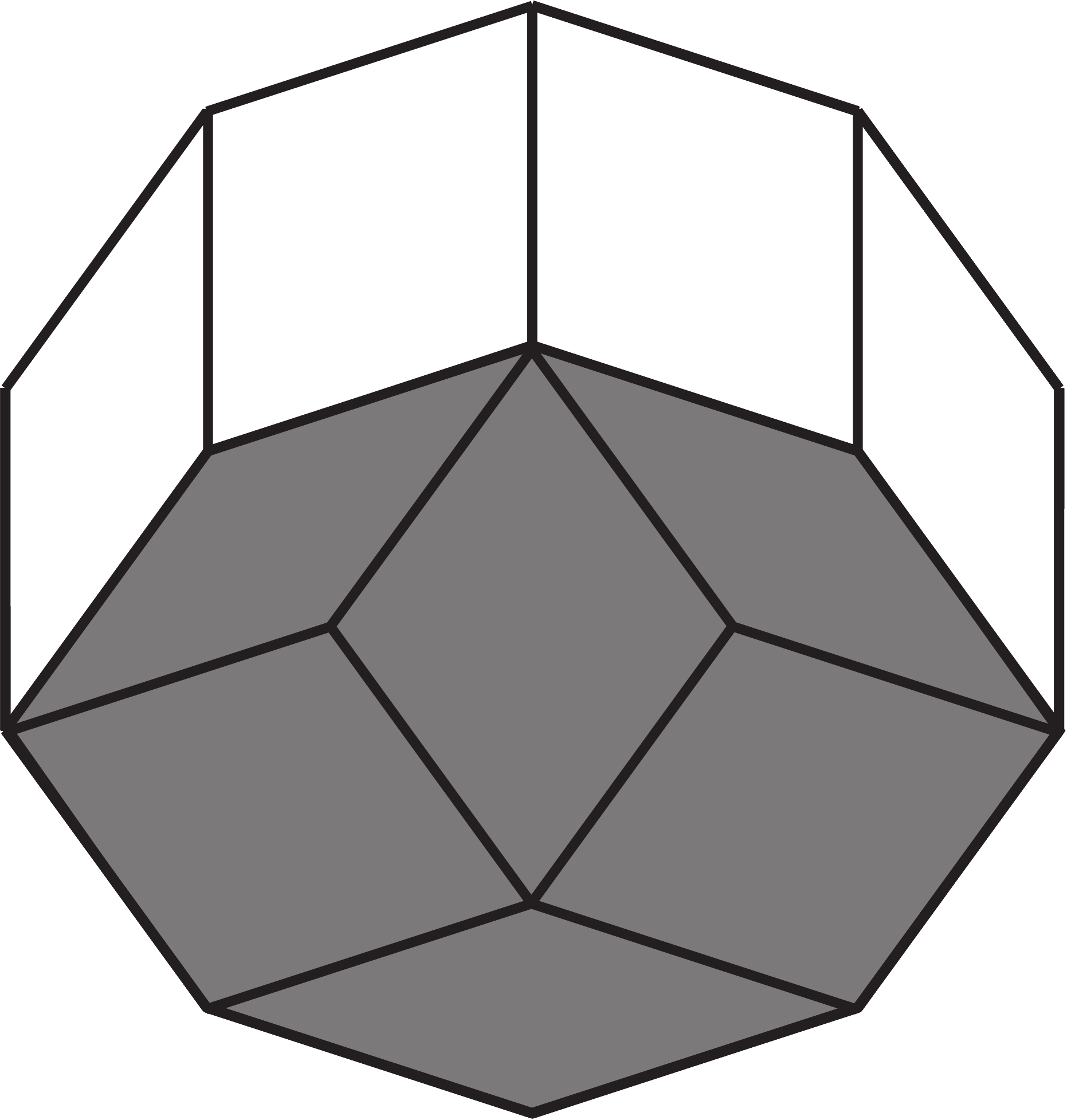}     
&     \includegraphics[width=1in]{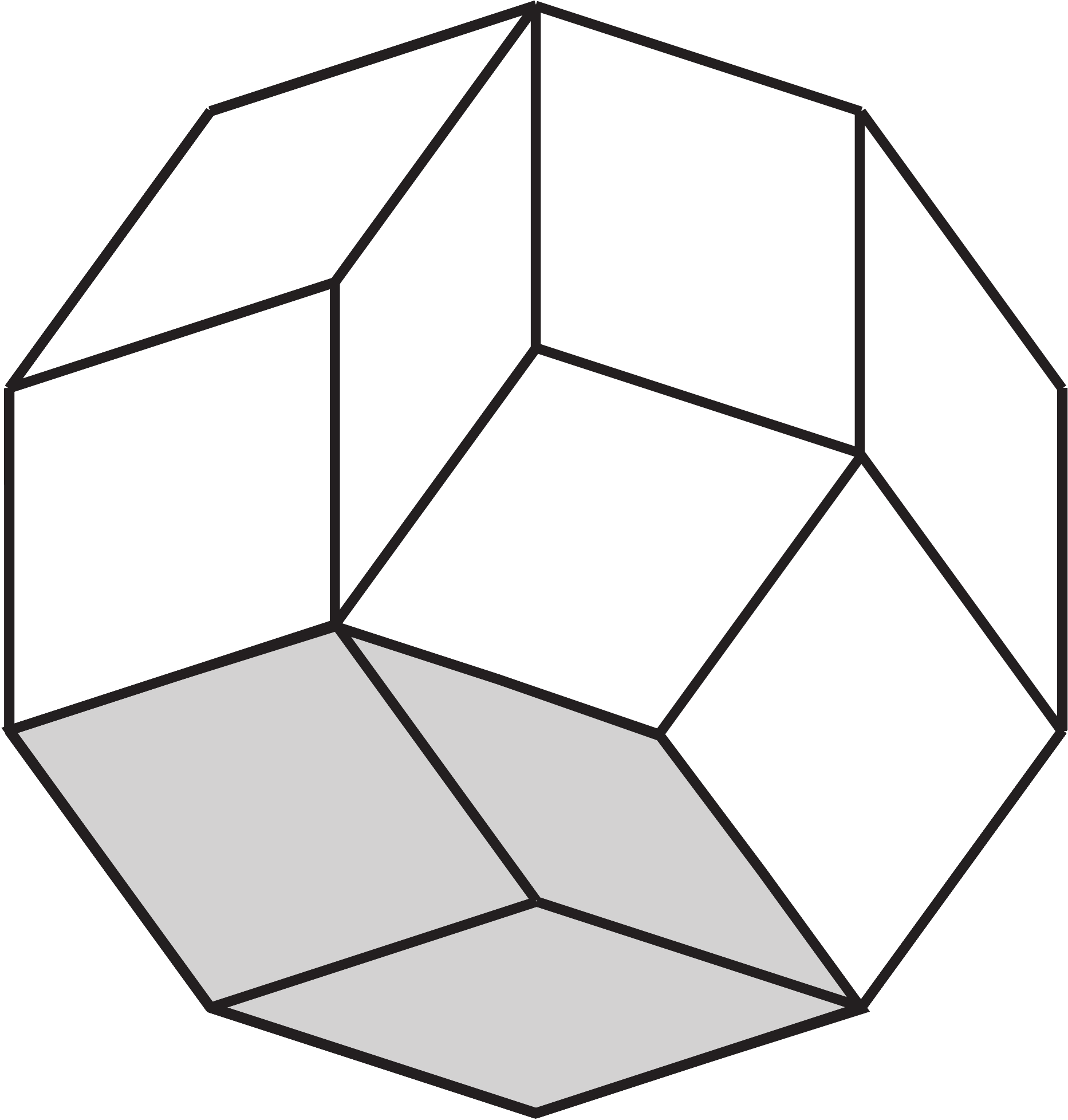}     
&     \includegraphics[width=1in]{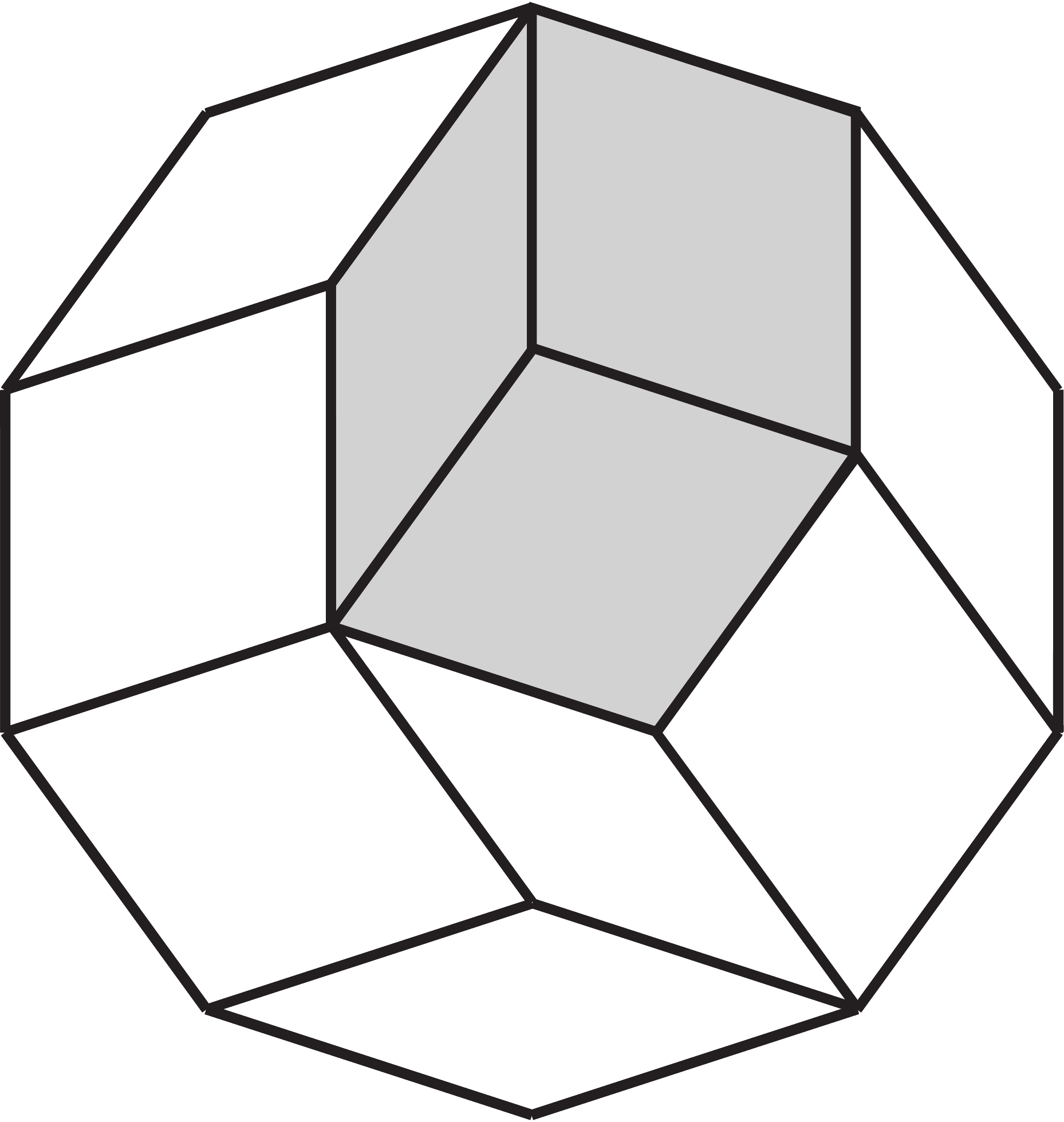} \\

-00-0 & 0++00  & 00+00 & 00++0 &  0-00- \\

      \includegraphics[width=1in]{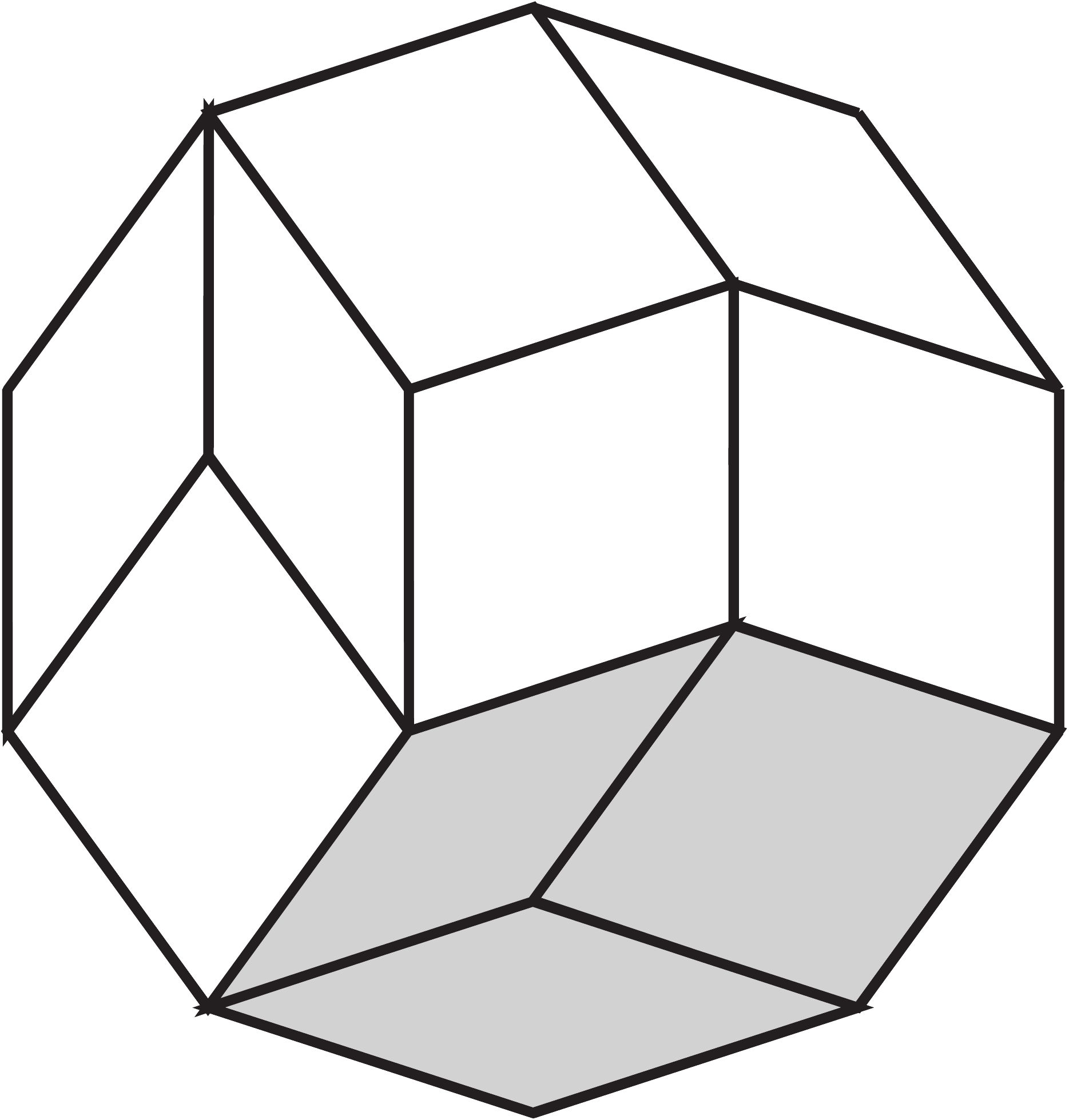}     
  &      \includegraphics[width=1in]{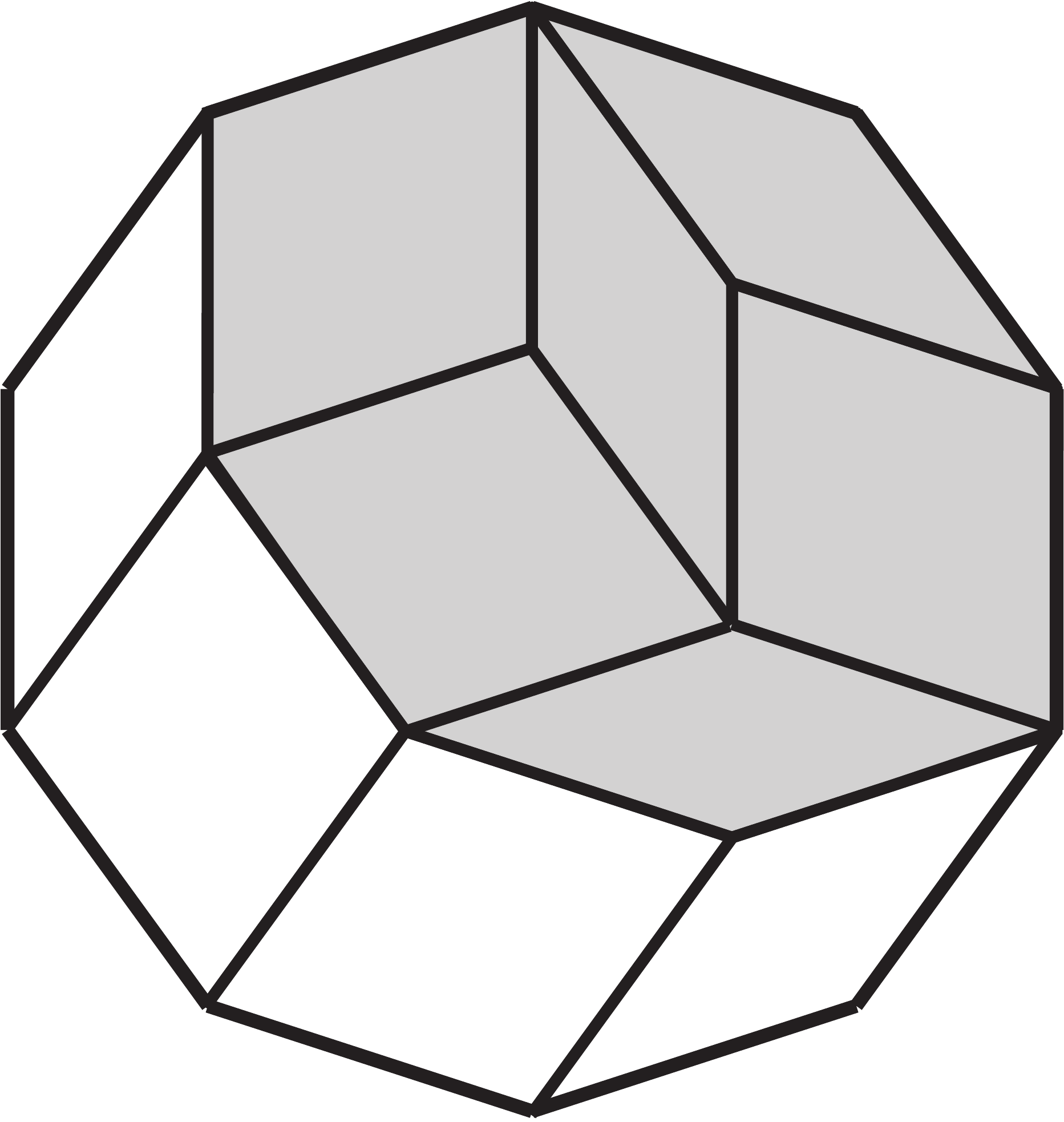}    
 &    \includegraphics[width=1in]{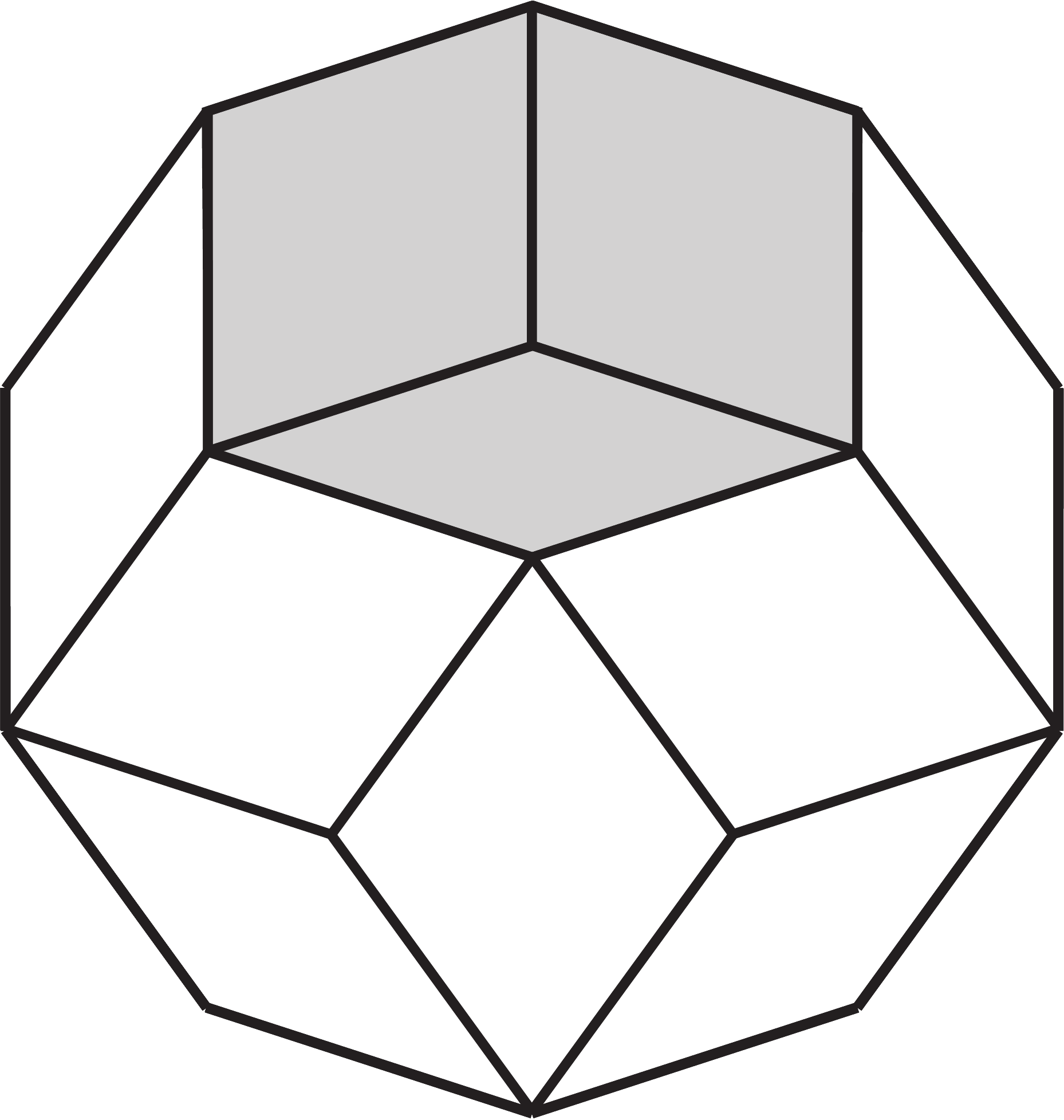}     
&     \includegraphics[width=1in]{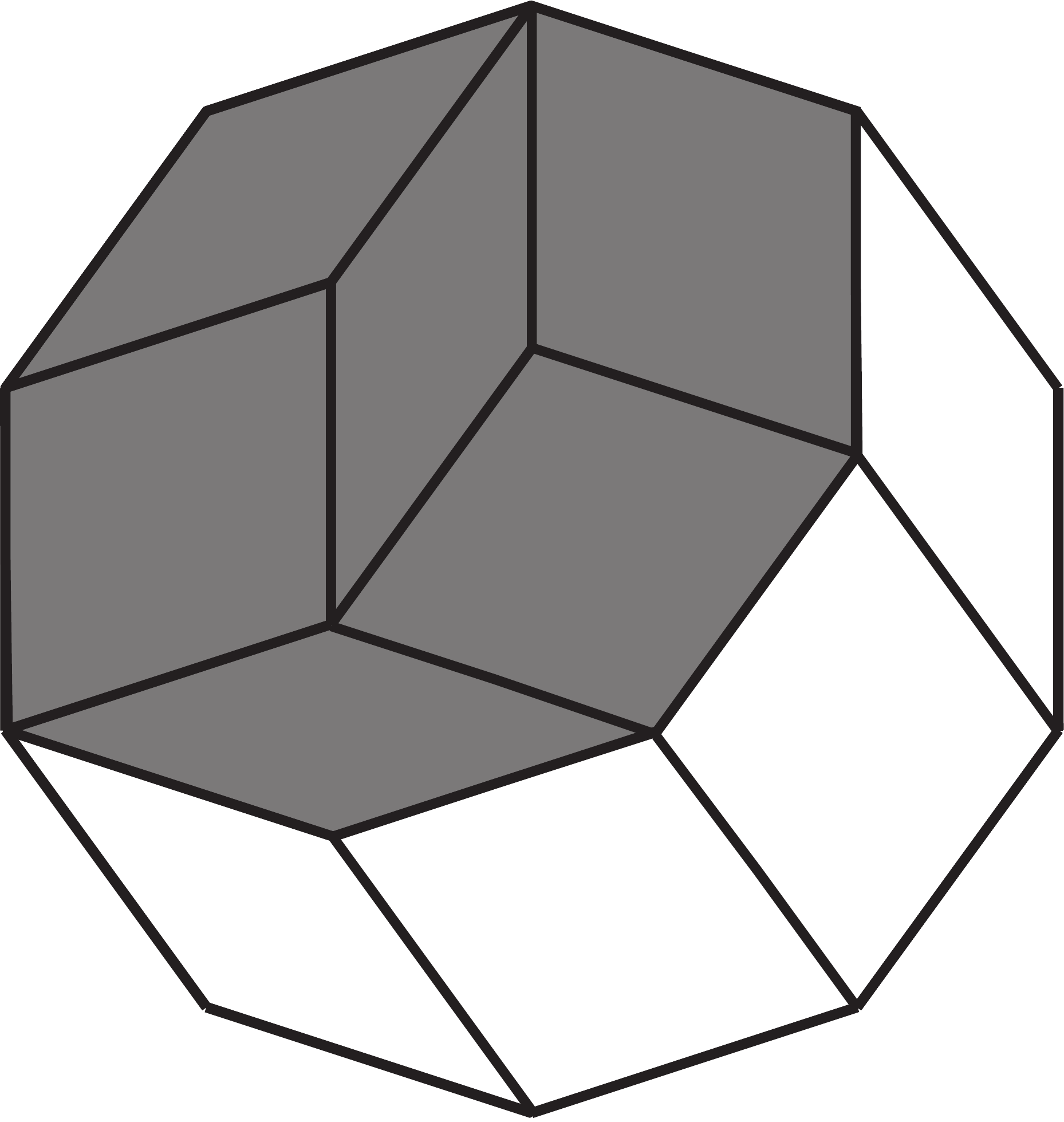}     
&     \includegraphics[width=1in]{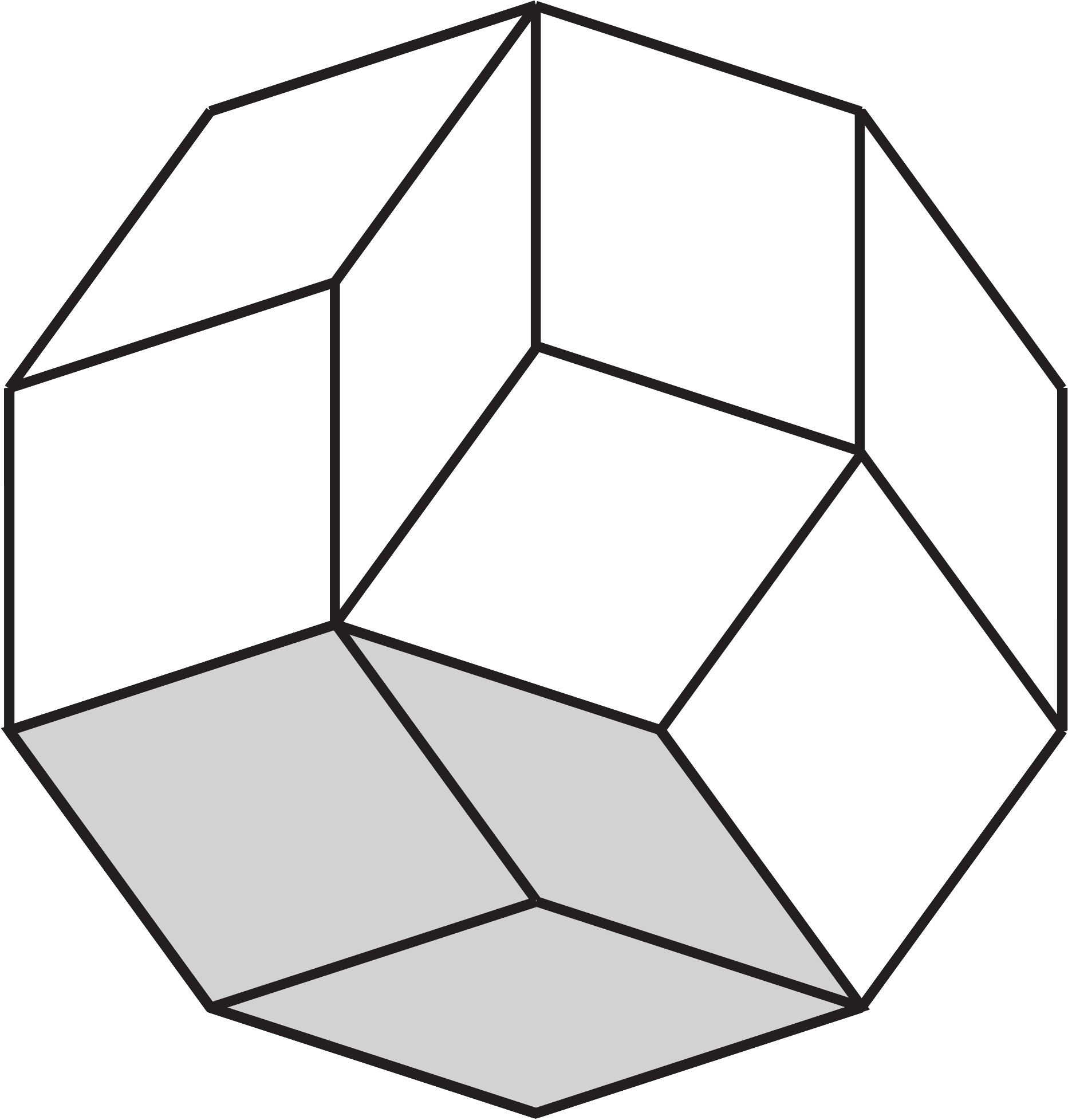} \\

0++00 & 000-0  & 0-0-0 & 0-000 & 0++00\\

     \includegraphics[width=1in]{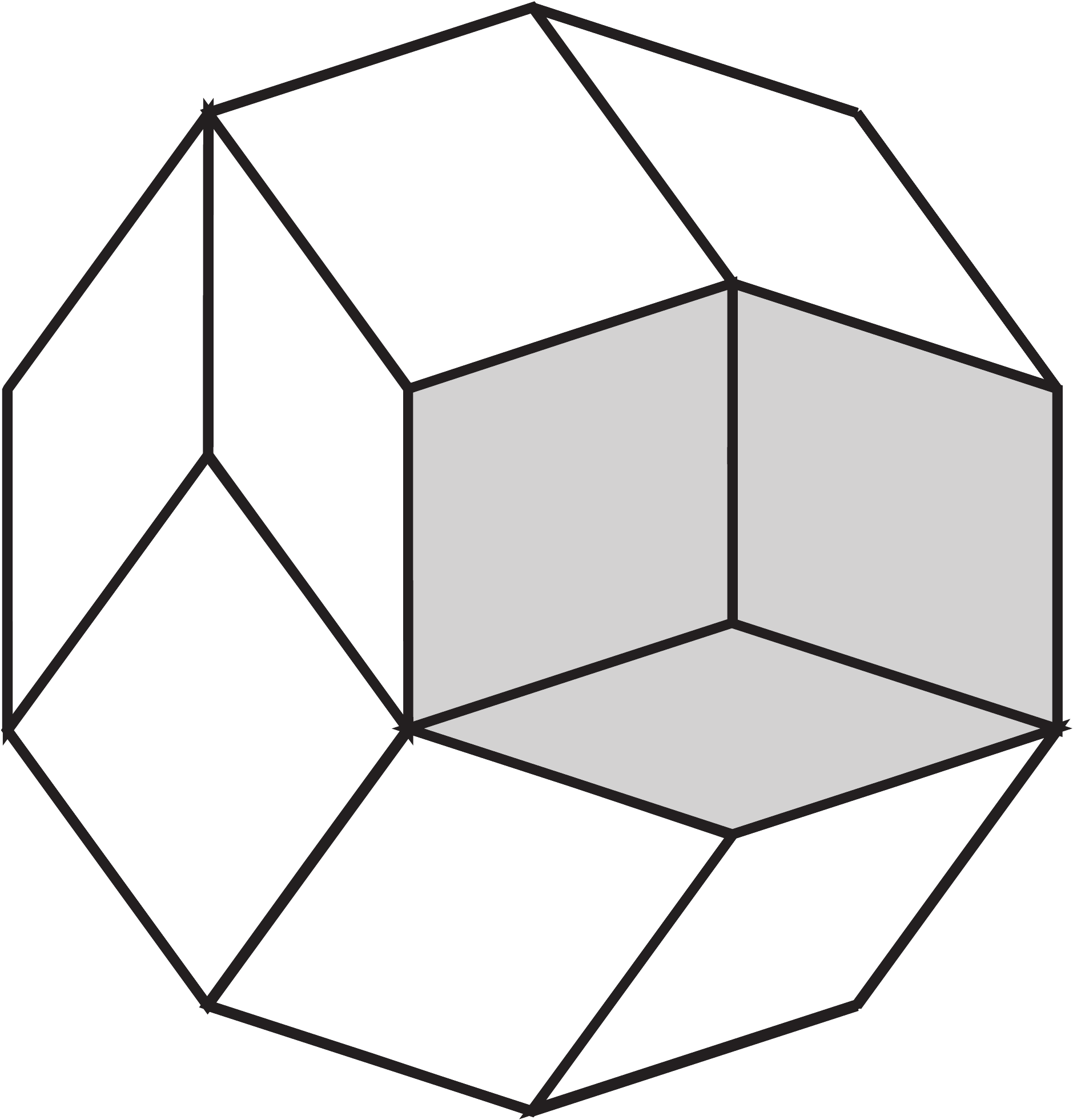}     
  &      \includegraphics[width=1in]{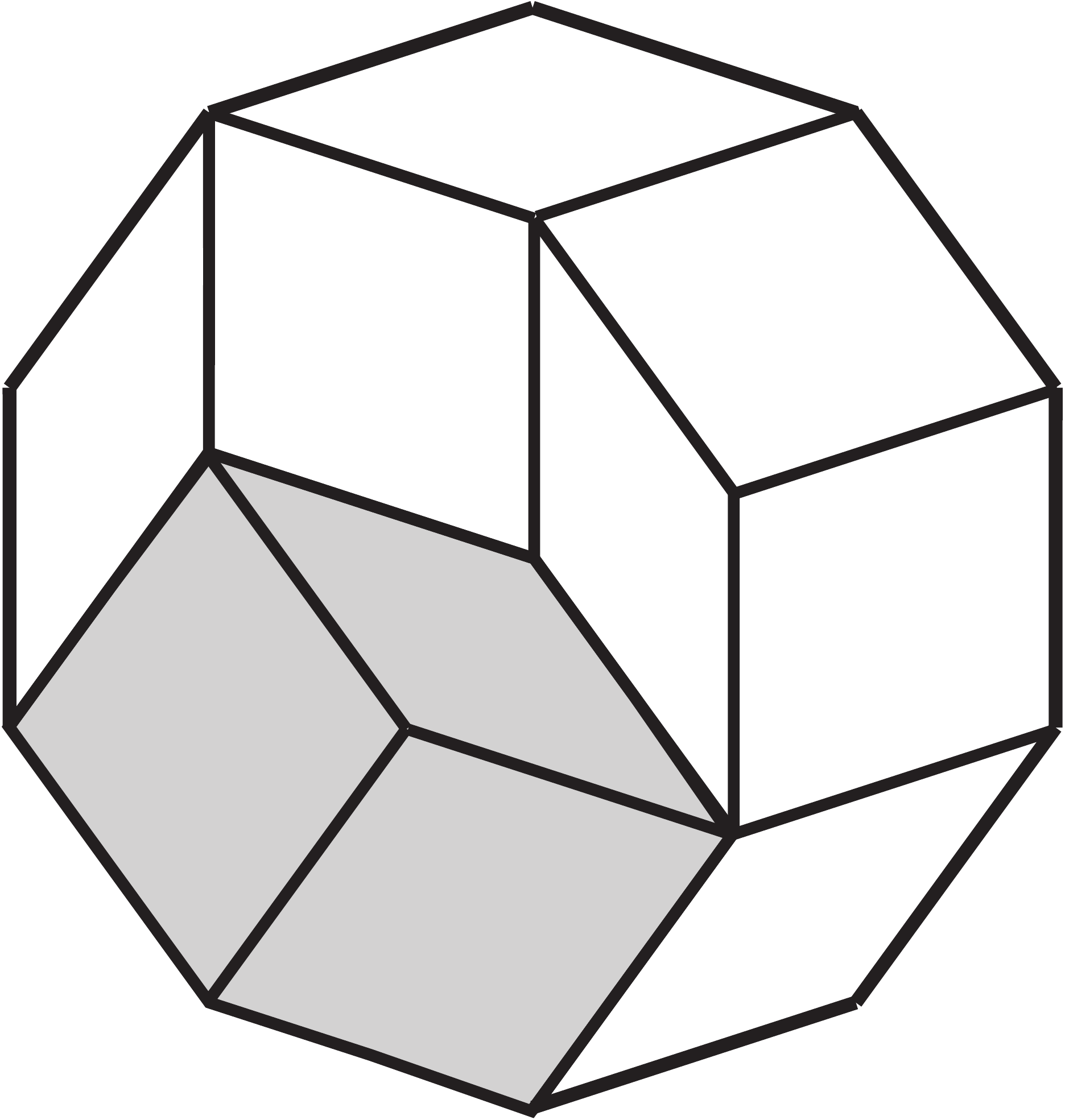}   
 &    \includegraphics[width=1in]{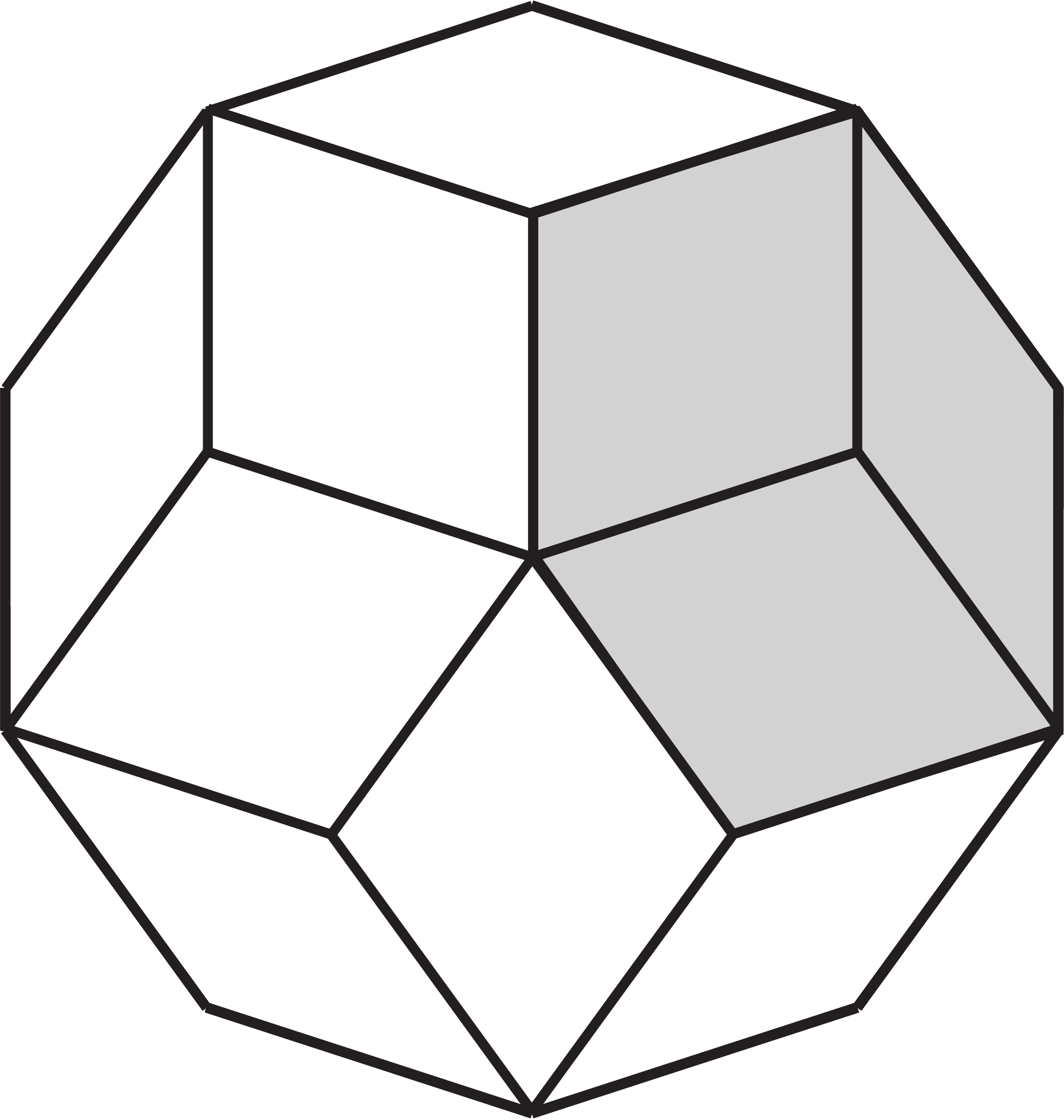}     
&     \includegraphics[width=1in]{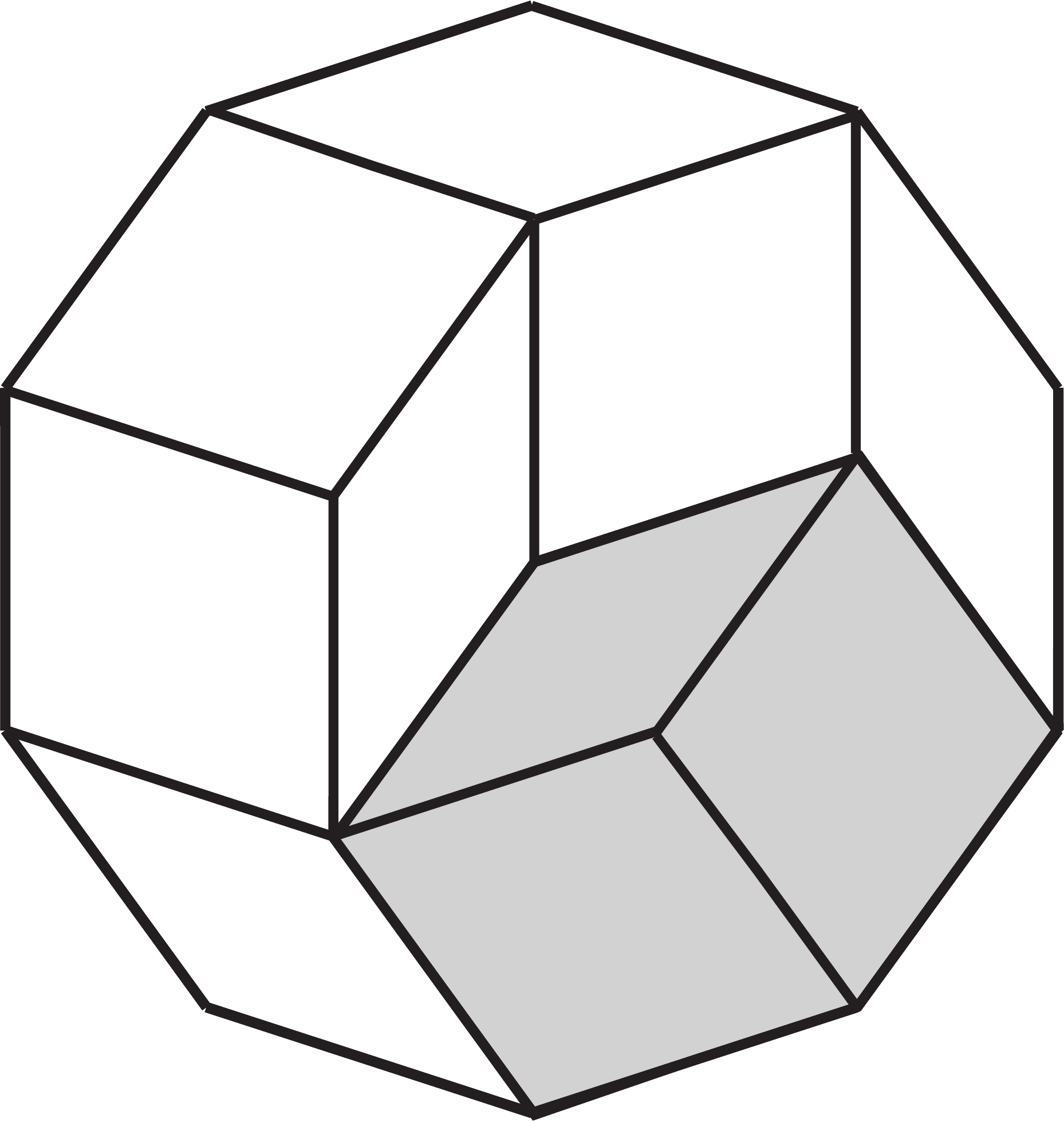}     
&     \includegraphics[width=1in]{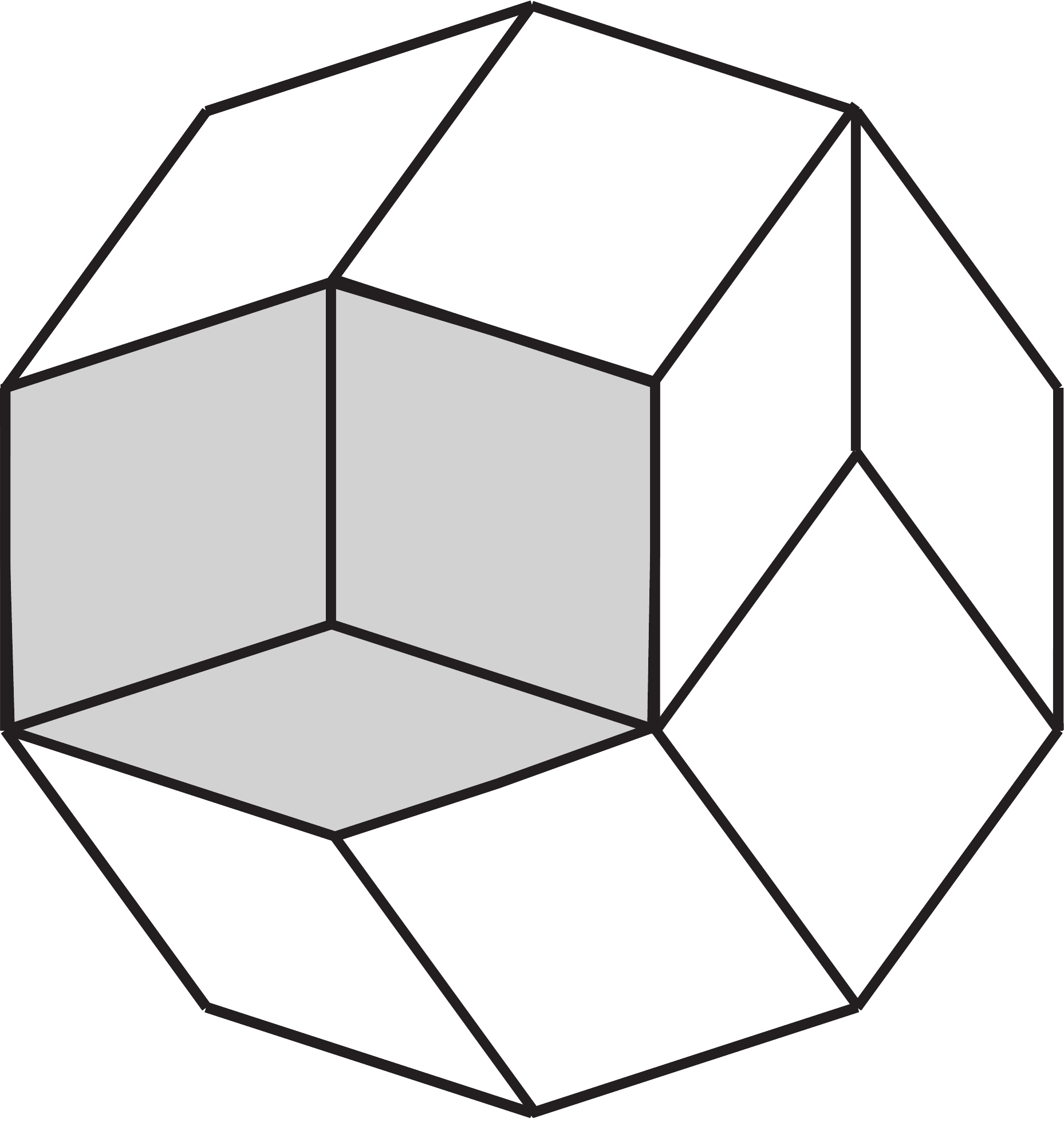} \\

0+0-0 & 00+0+  & +00-0 & +0+00 & 0-0+0\\

     \includegraphics[width=1in]{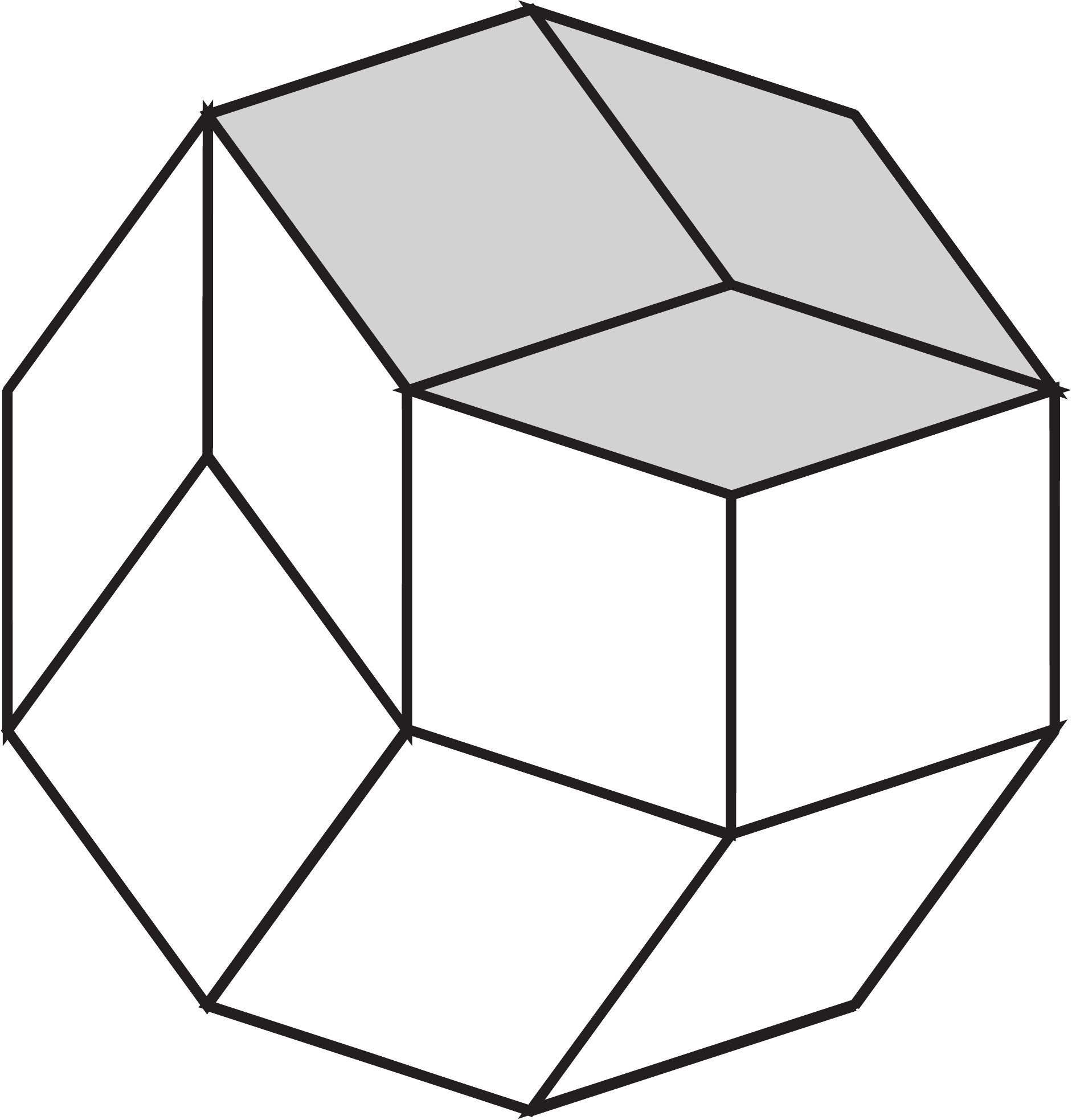}     
  &      \includegraphics[width=1in]{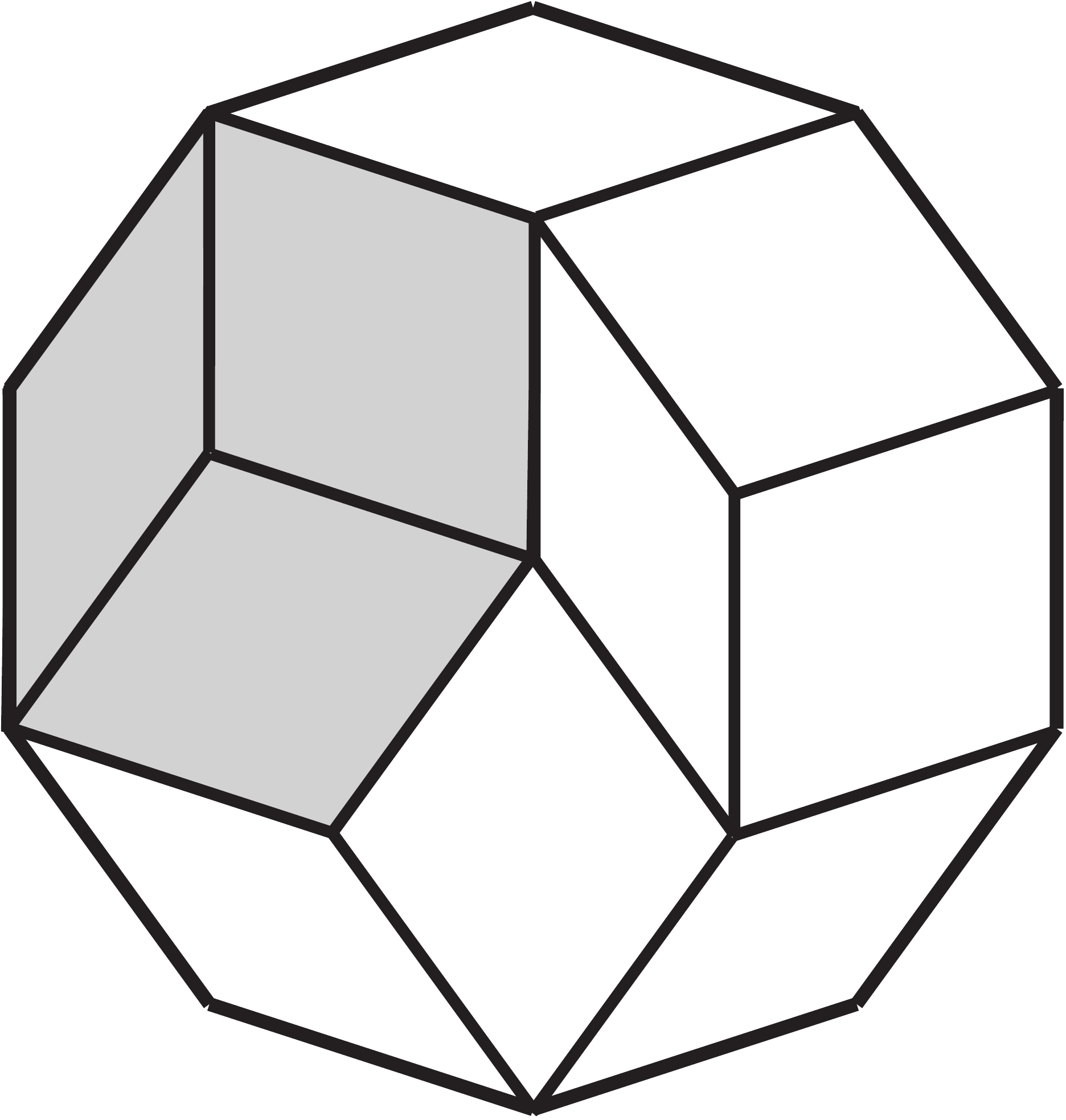}    
 &    \includegraphics[width=1in]{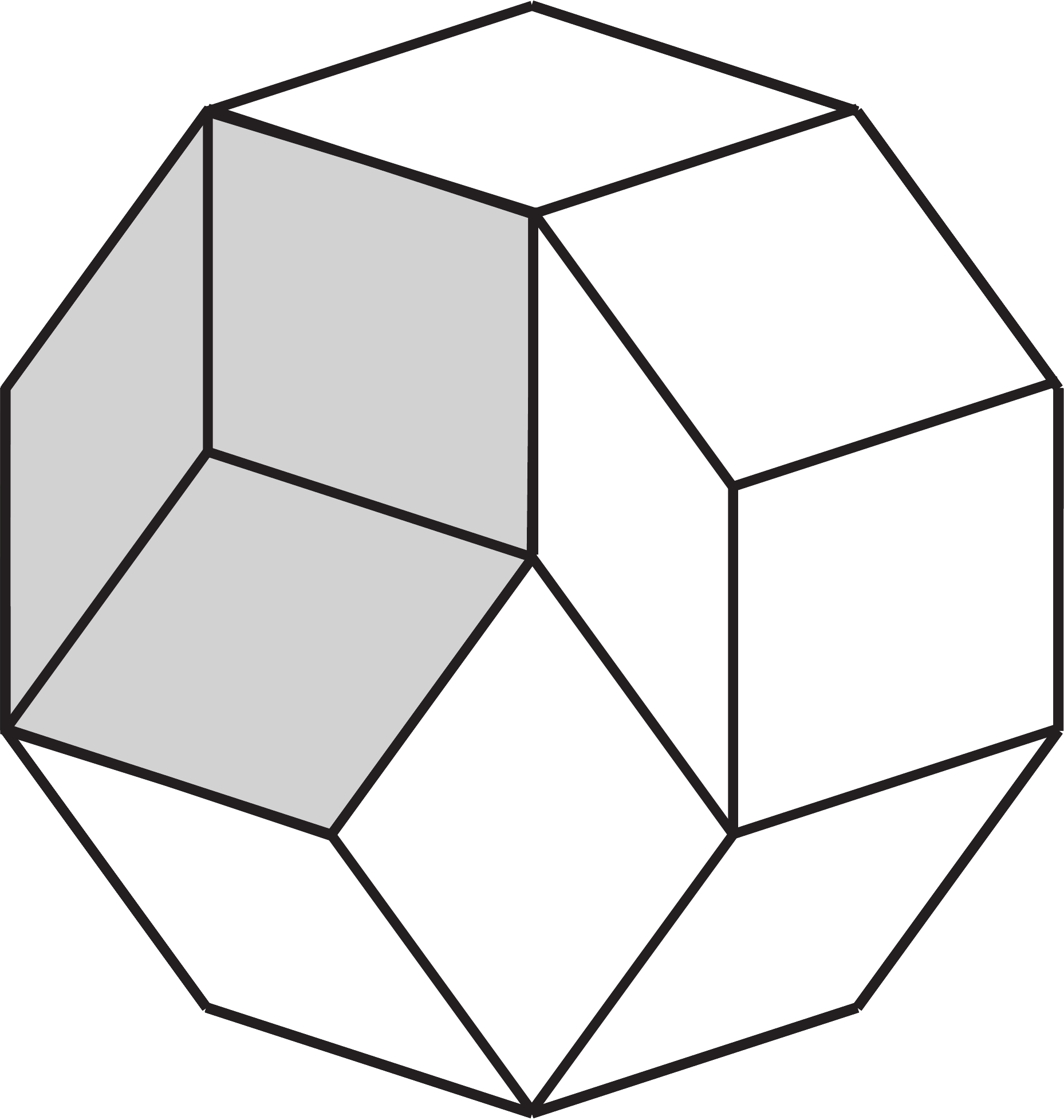}     
&     \includegraphics[width=1in]{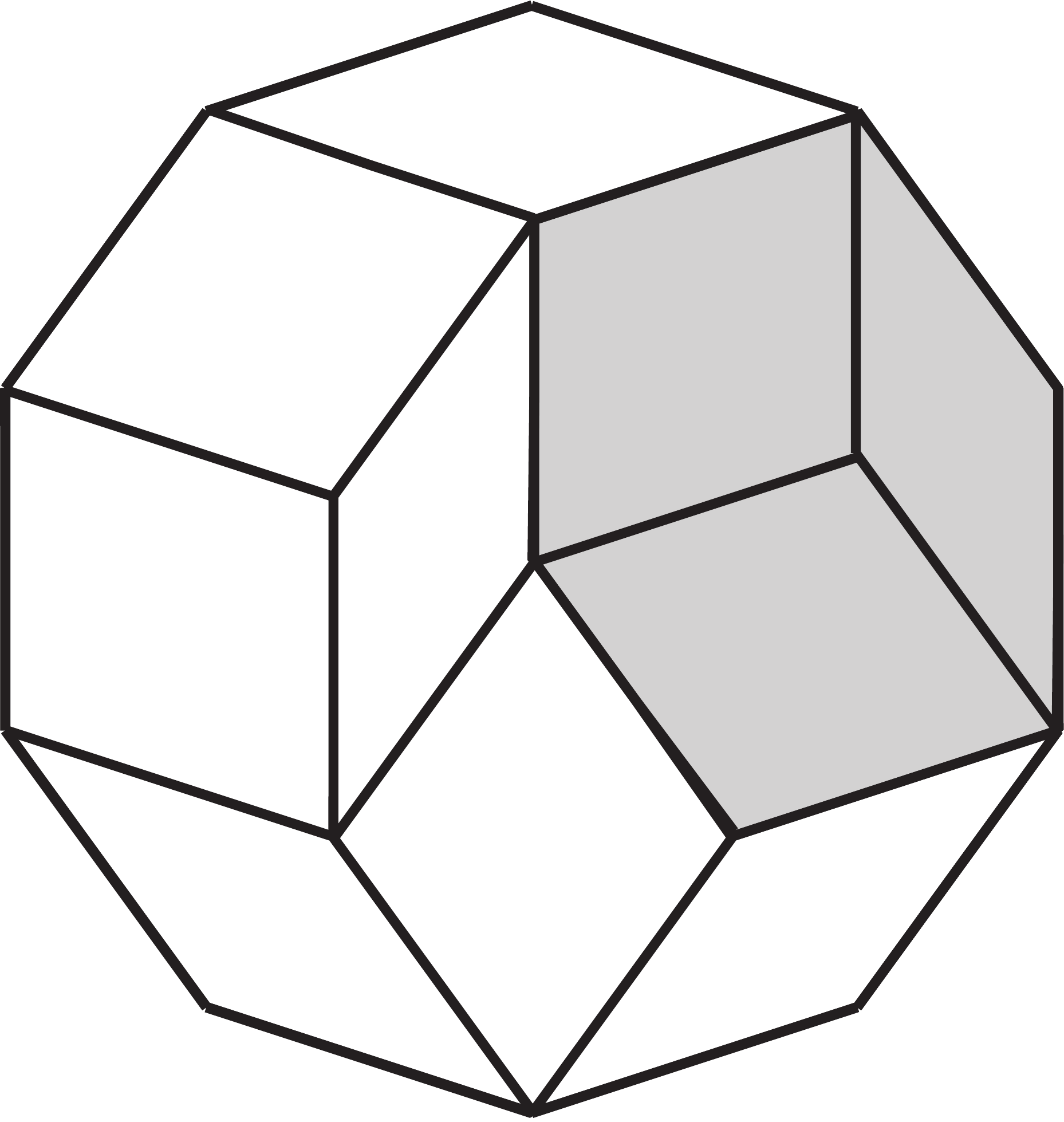}     
&     \includegraphics[width=1in]{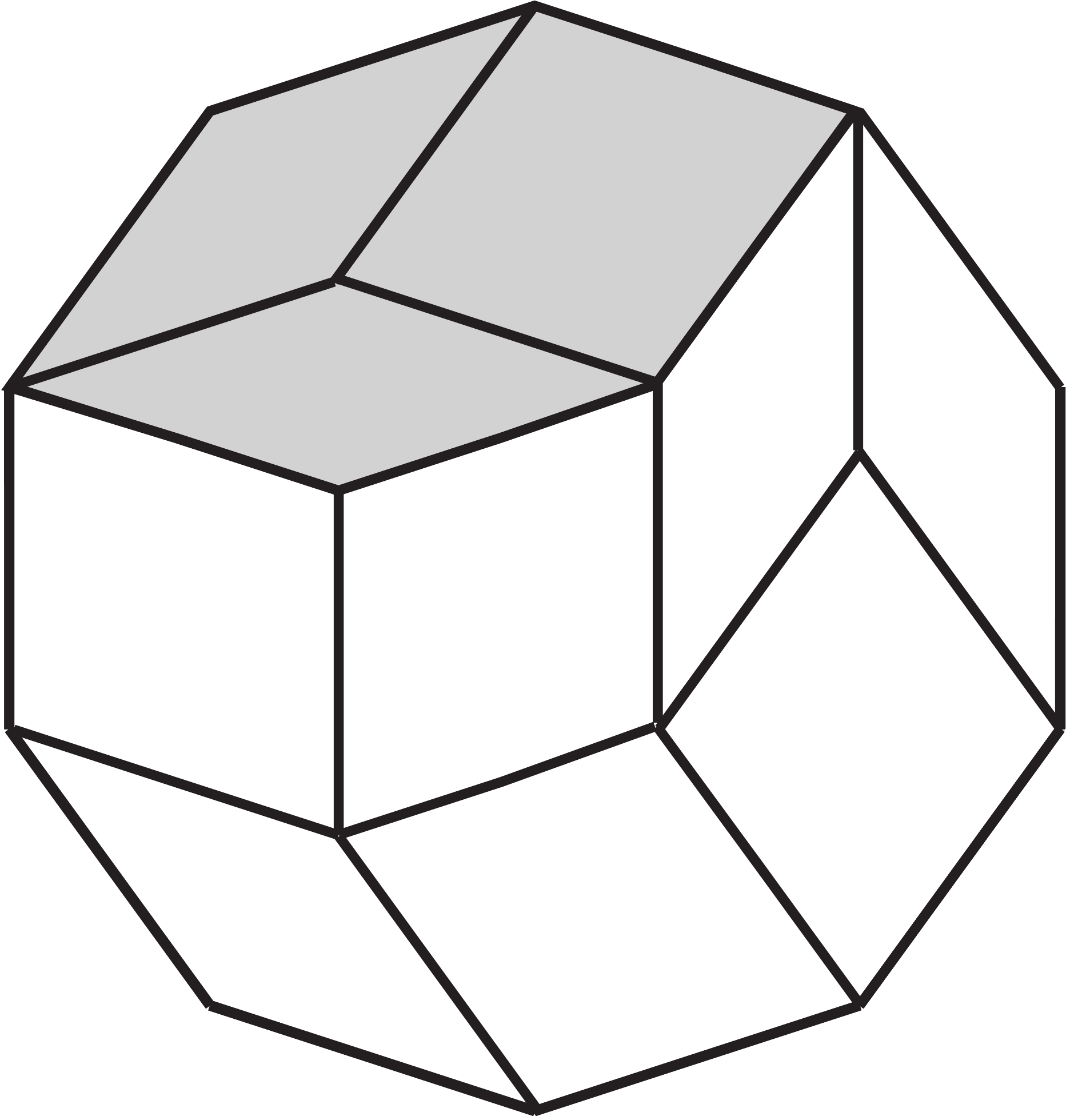} \\

00--0 & 0-00+  & 0-00+ & +00-0 & 0--00\\

      \includegraphics[width=1in]{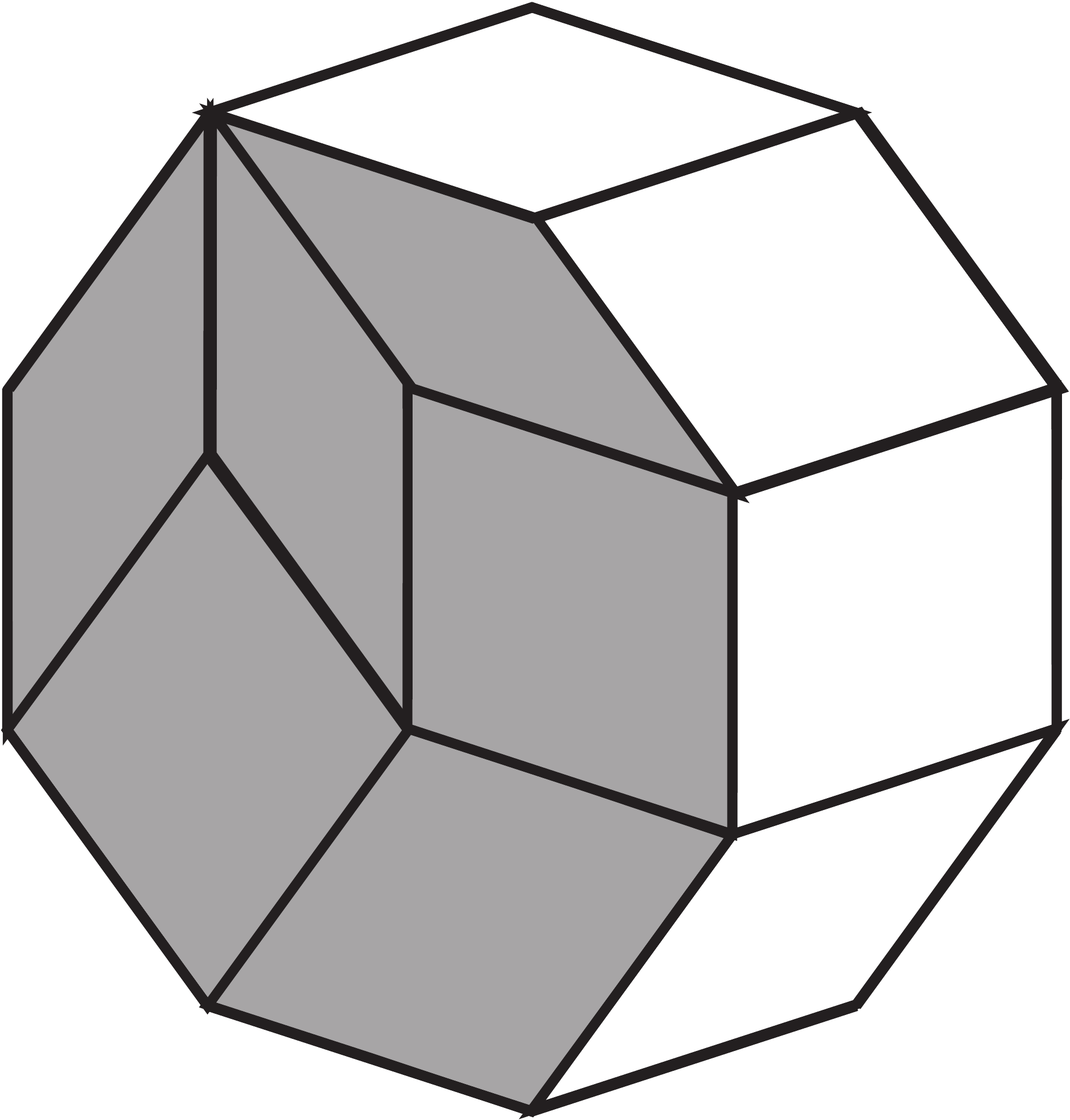}     
  &      \includegraphics[width=1in]{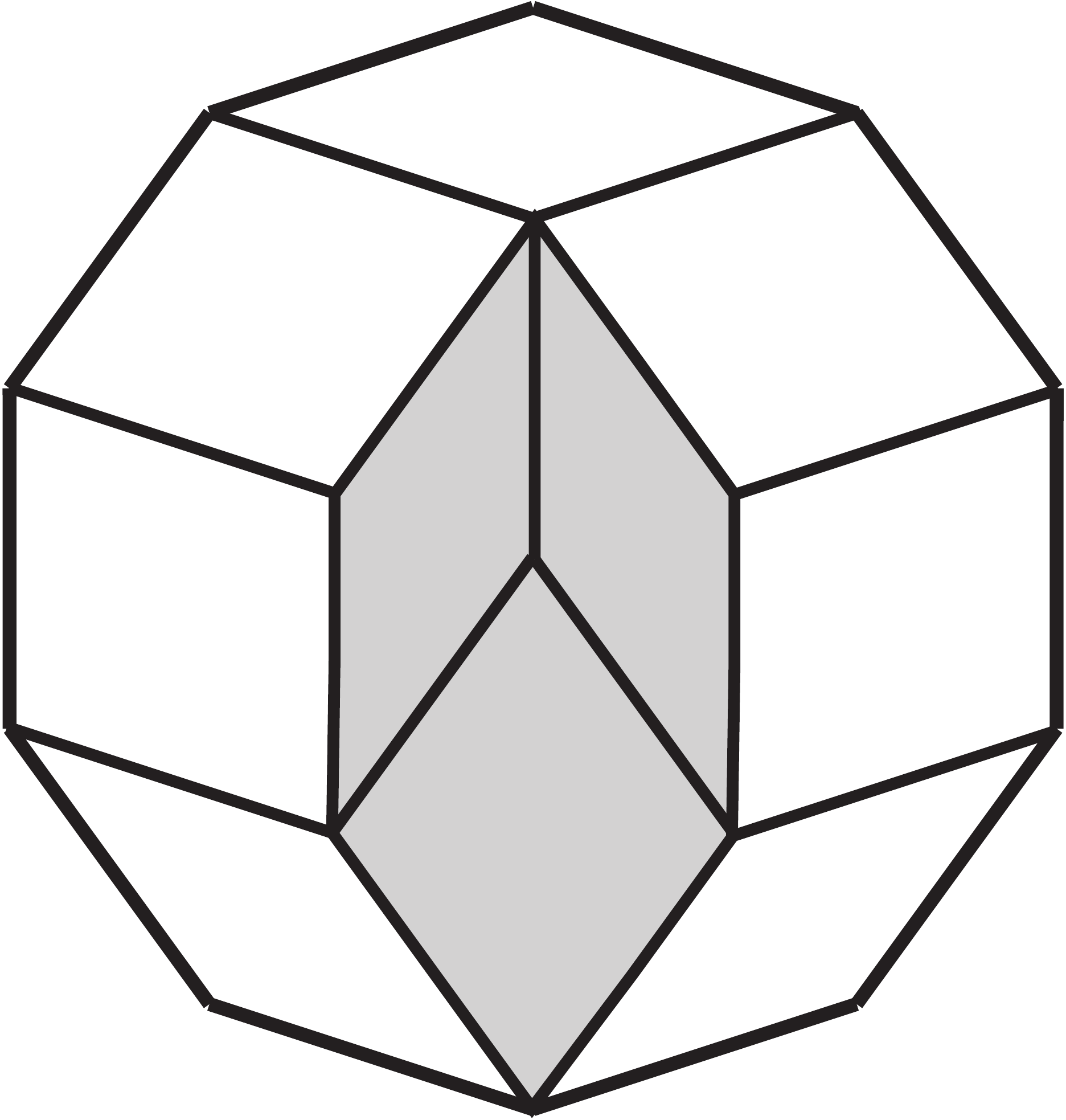}     
 &    \includegraphics[width=1in]{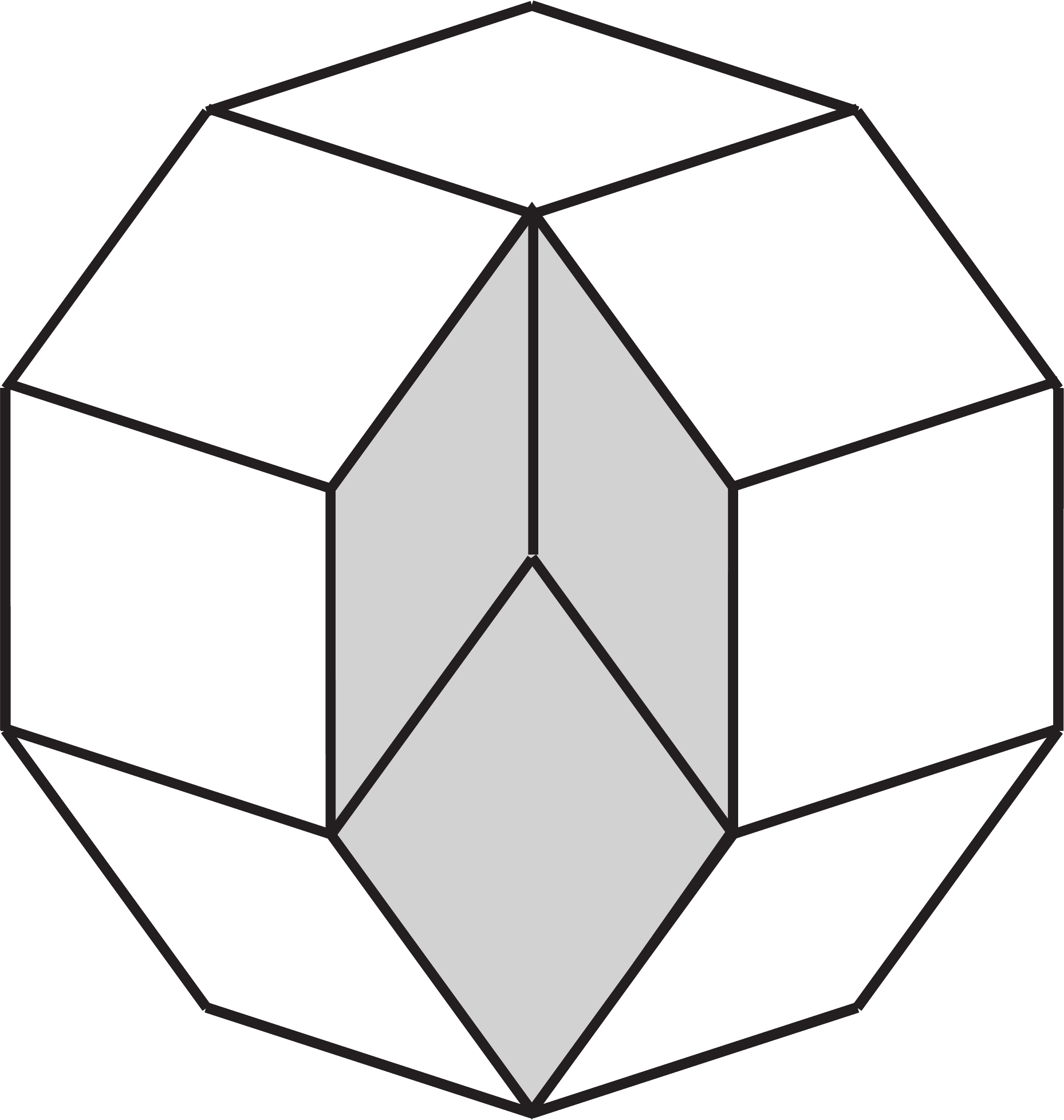}     
&     \includegraphics[width=1in]{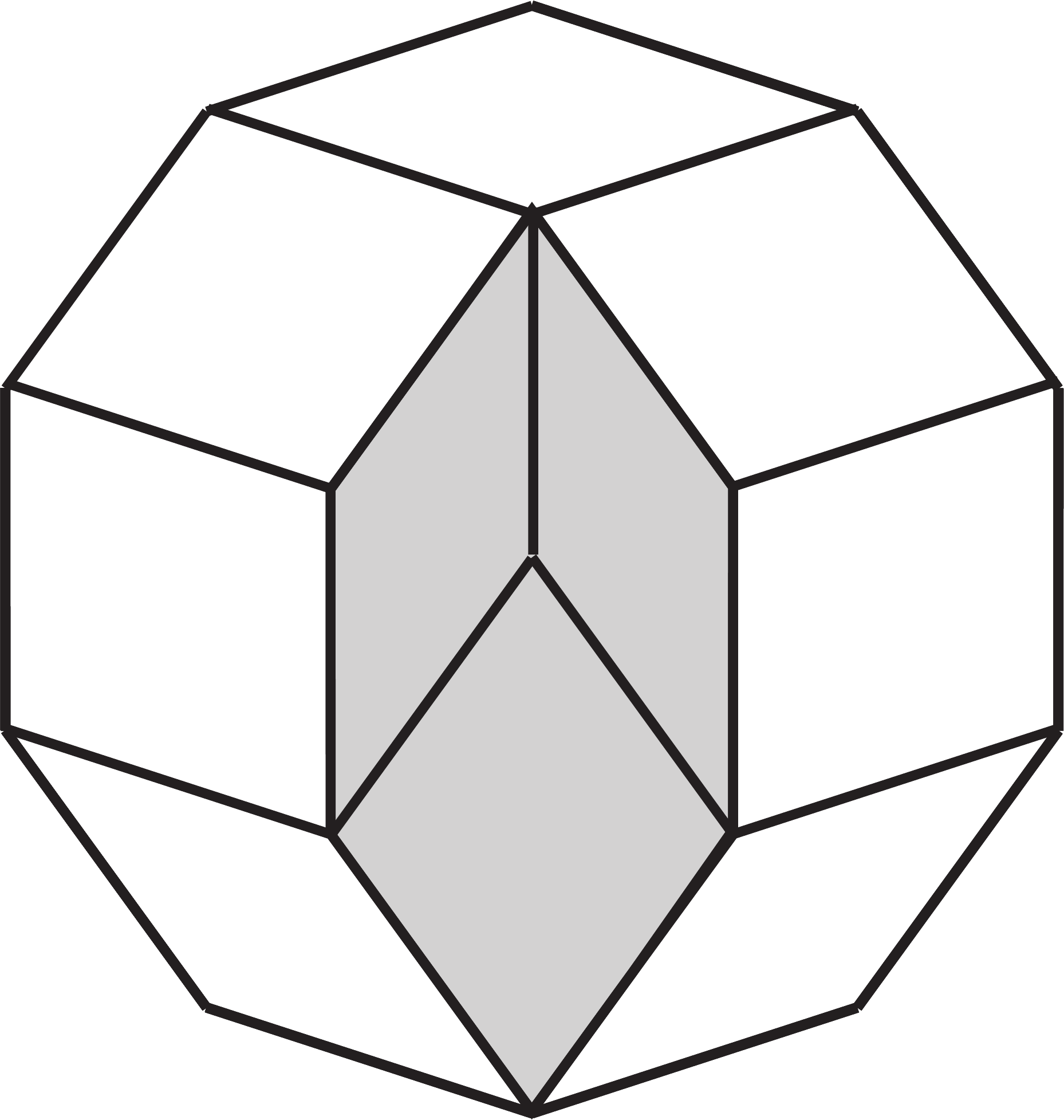}     
&     \includegraphics[width=1in]{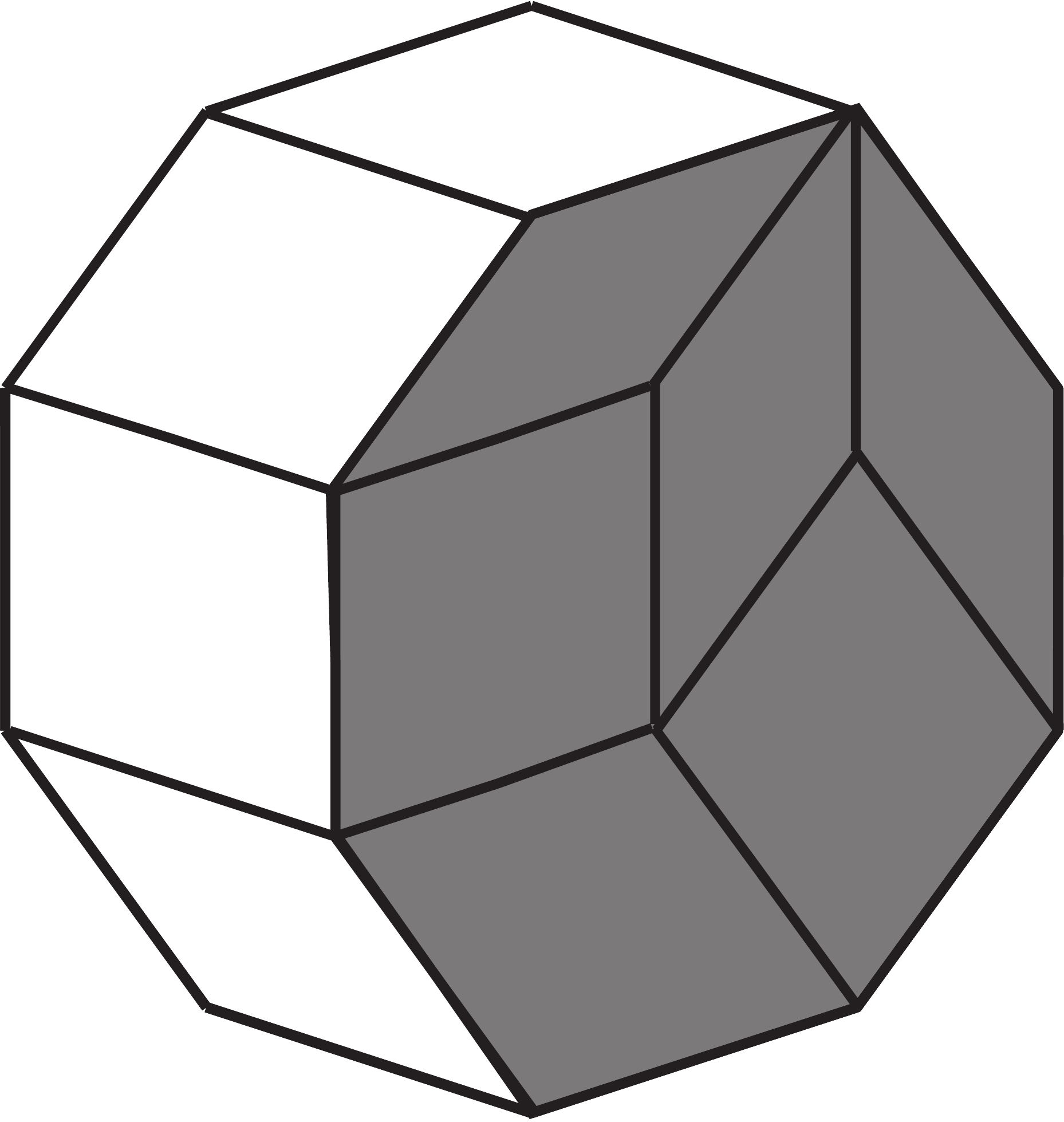} \\

0000+ & +000+   & +000+ & +000+ & +0000\\
% IRA   LIGHTEN OR DARKEN 4-cell  IRA
 %6-0  &  6-1  &  6-2  &  6-3  &  6-4
% \label{fig:19-2}
%\end{figure}
 \end{array}
%\right)
\]

  \[
%\left(
\begin{array}{ccccc}  
      \includegraphics[width=0.9in]{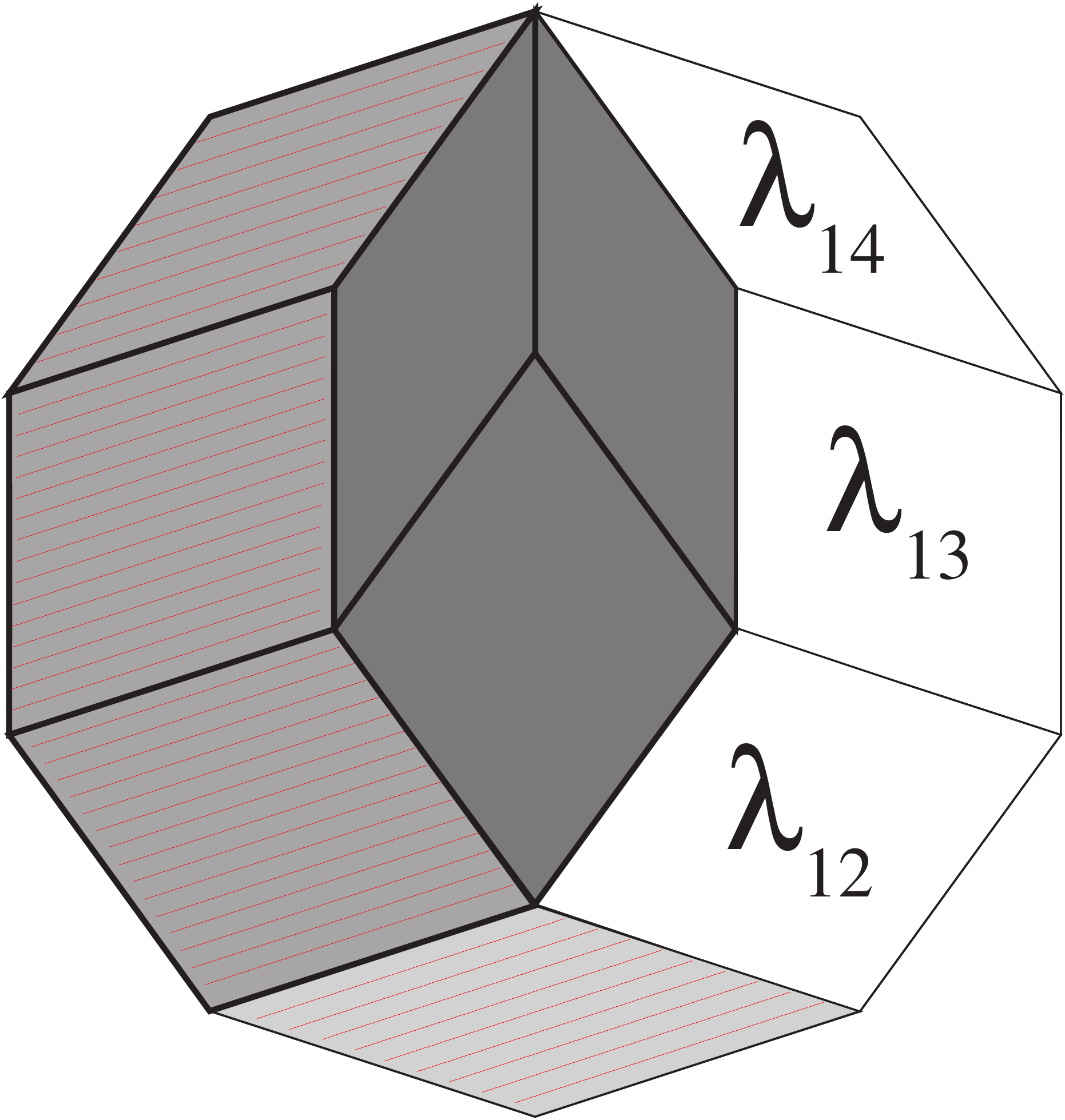}      
  &      \includegraphics[width=0.9in]{m11}     
 &    \includegraphics[width=0.9in]{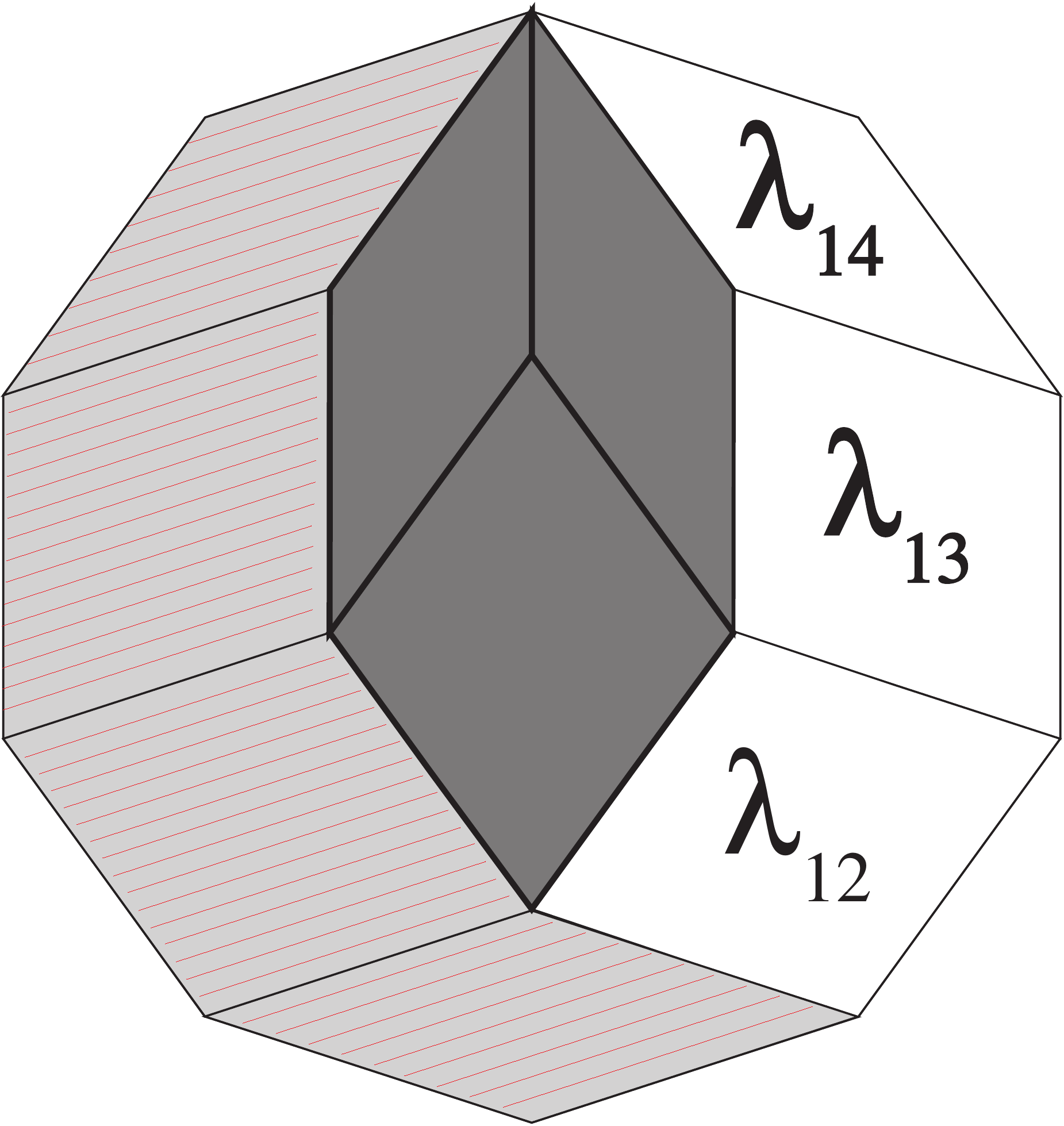}     
&     \includegraphics[width=0.9in]{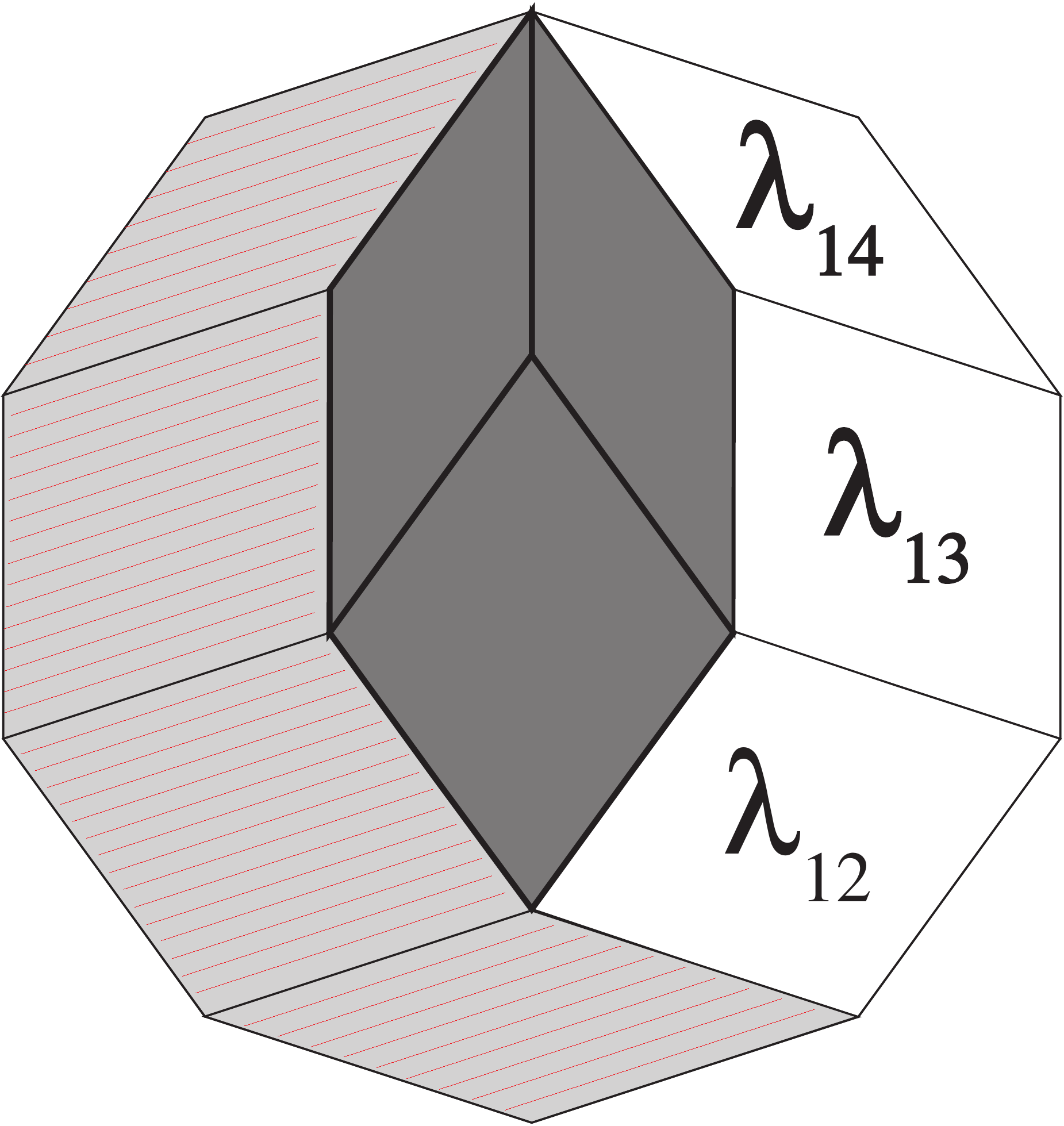}     
&     \includegraphics[width=0.9in]{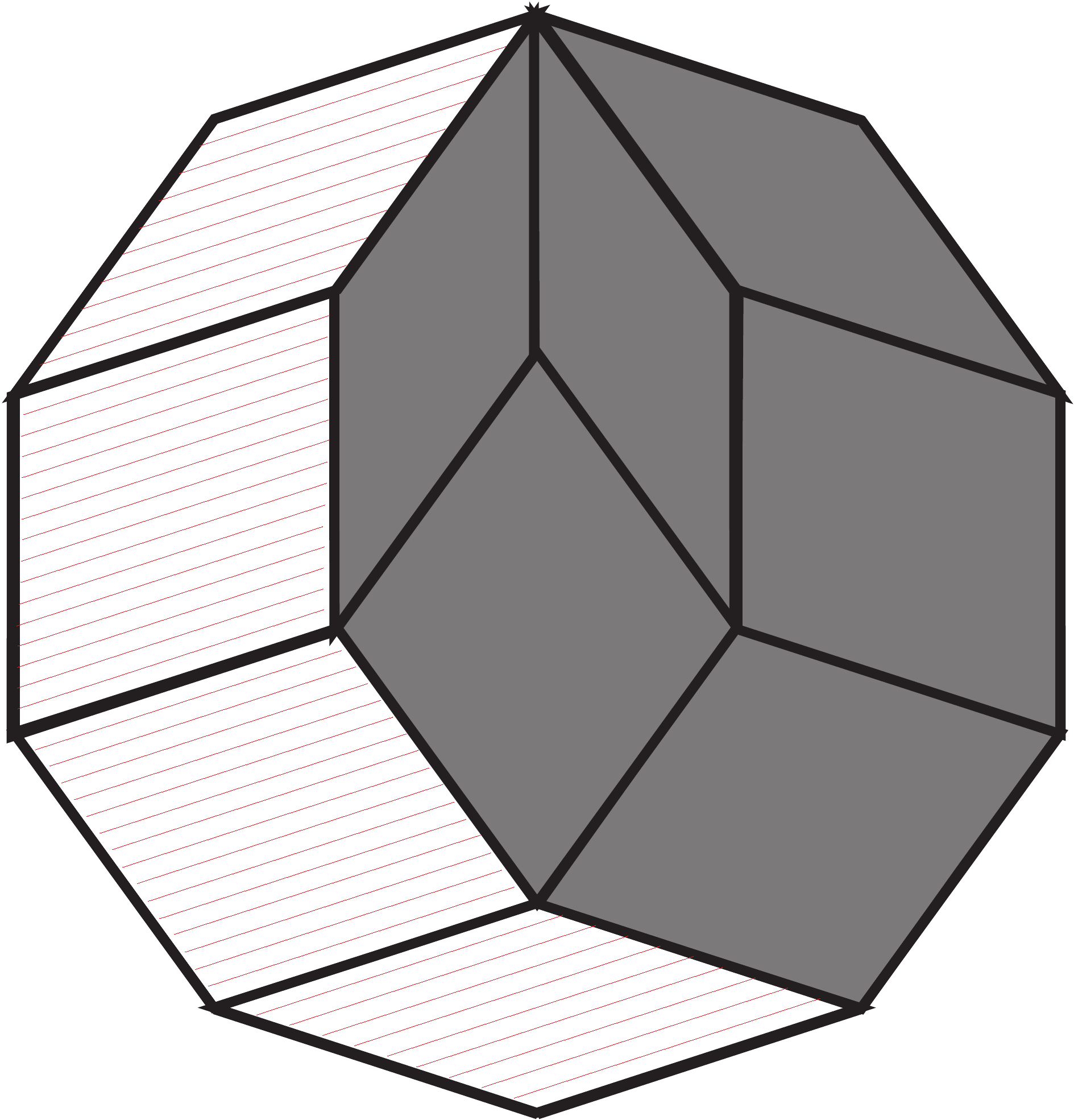} 
\\

\mu_{11} & \lambda_{11} &  \lambda_{11}  &  \lambda_{11}  & \lambda    \\

     \includegraphics[width=0.9in]{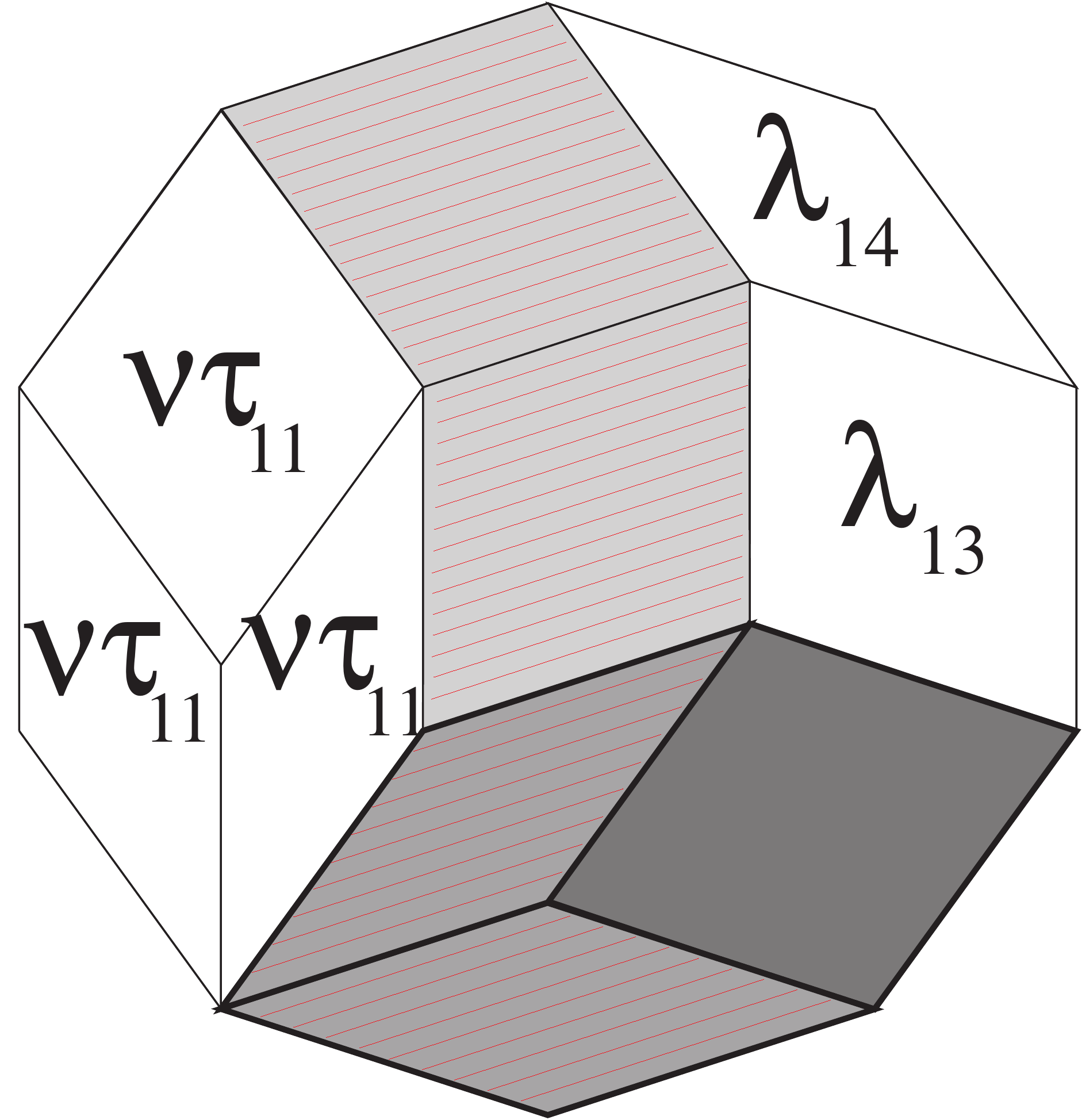}     
  &      \includegraphics[width=0.9in]{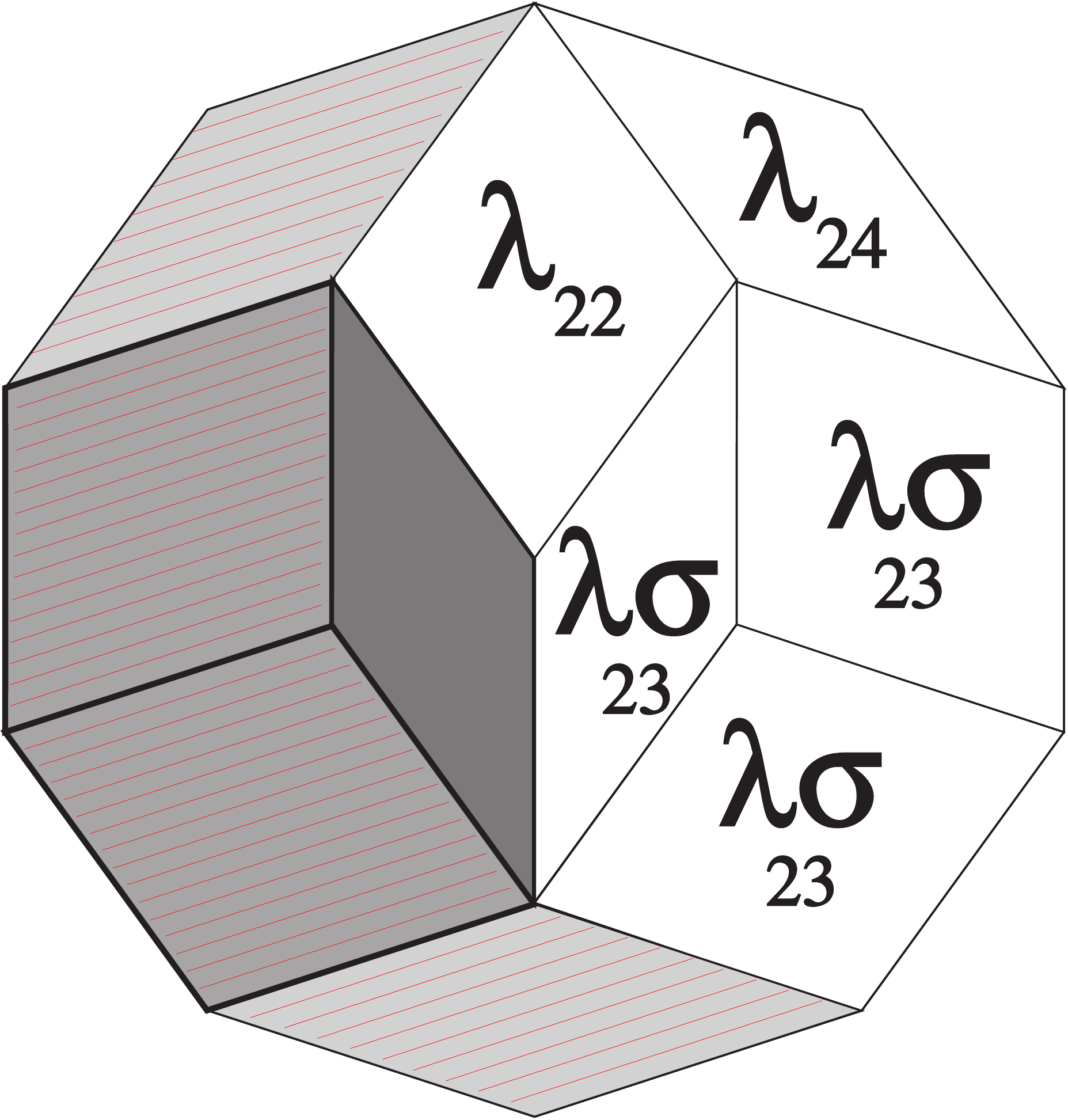}     
 &    \includegraphics[width=0.9in]{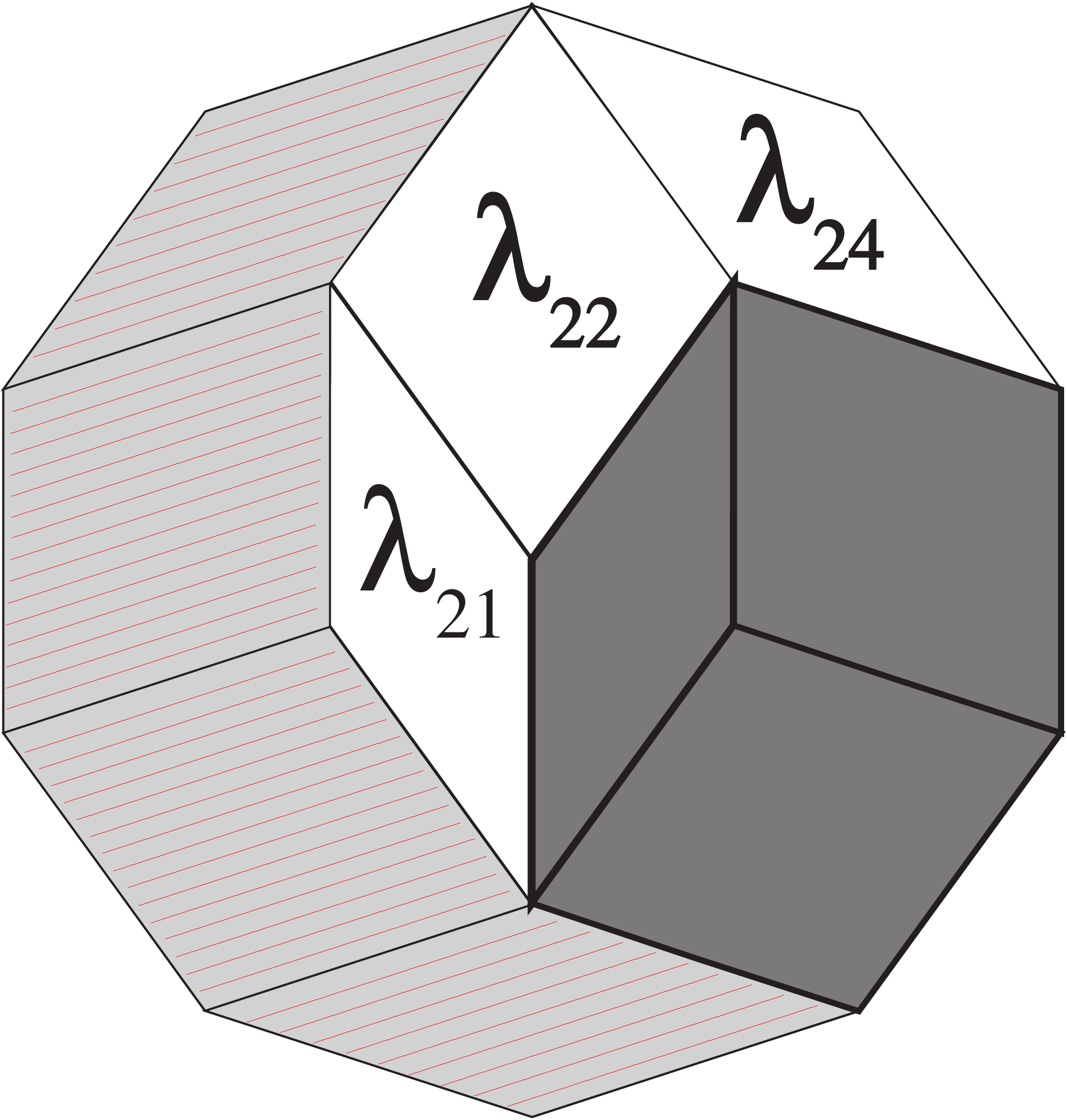}     
&     \includegraphics[width=0.9in]{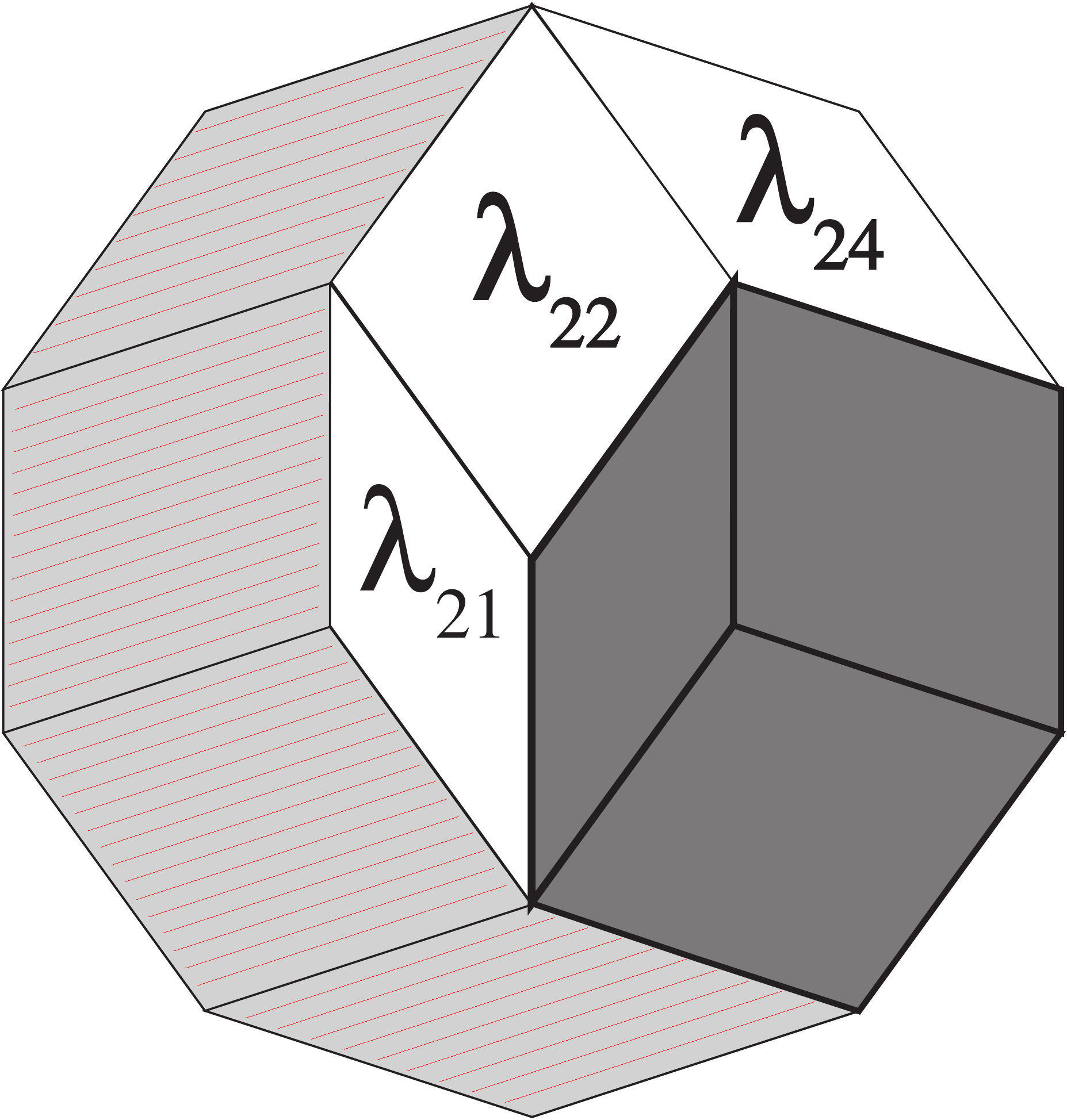}     
&     \includegraphics[width=0.9in]{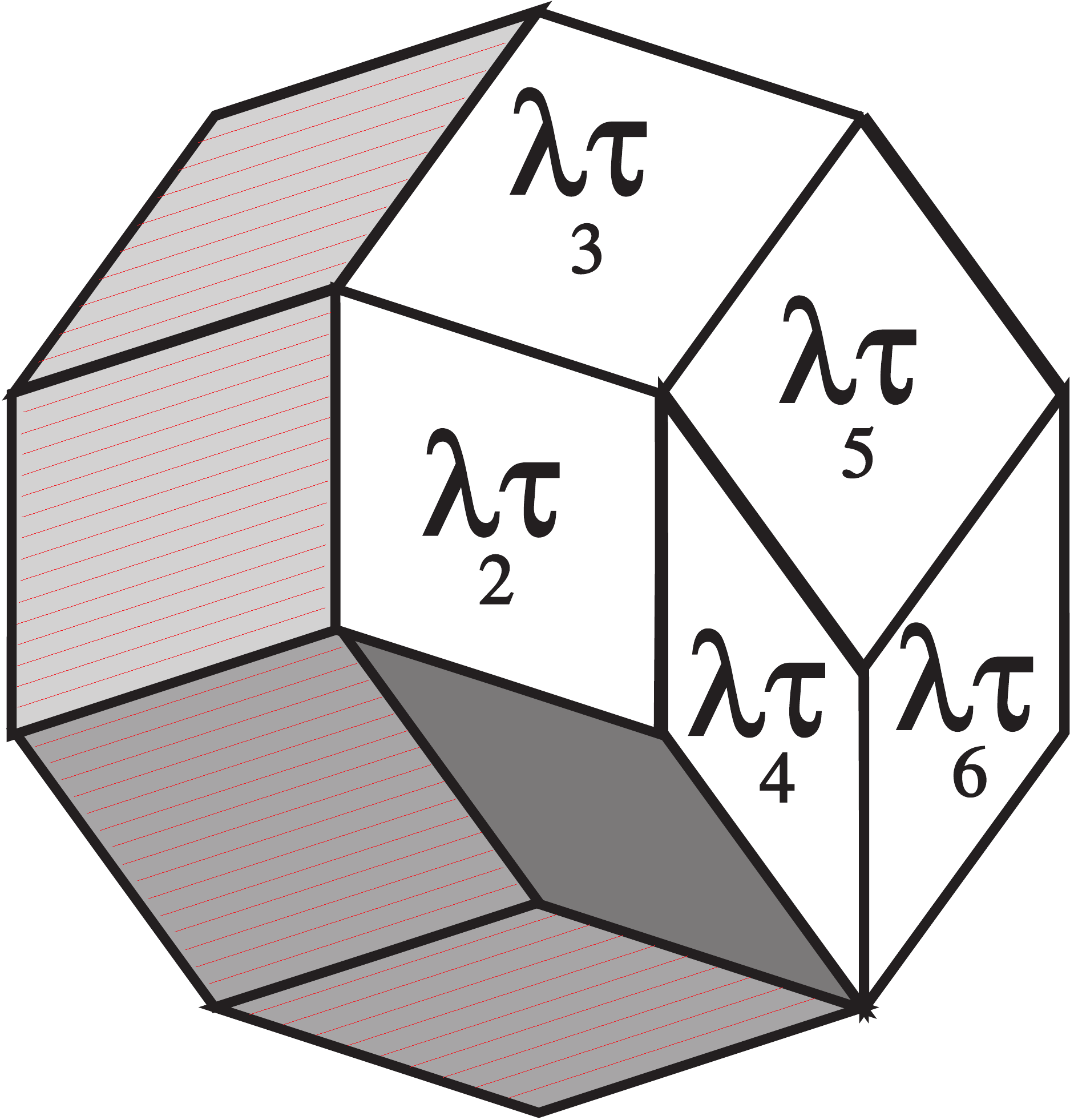} 
\\

\mu_{12} & \mu_{21} &  \lambda_{23}  &  \lambda_{23}  & \mu_{1}    \\

      \includegraphics[width=0.9in]{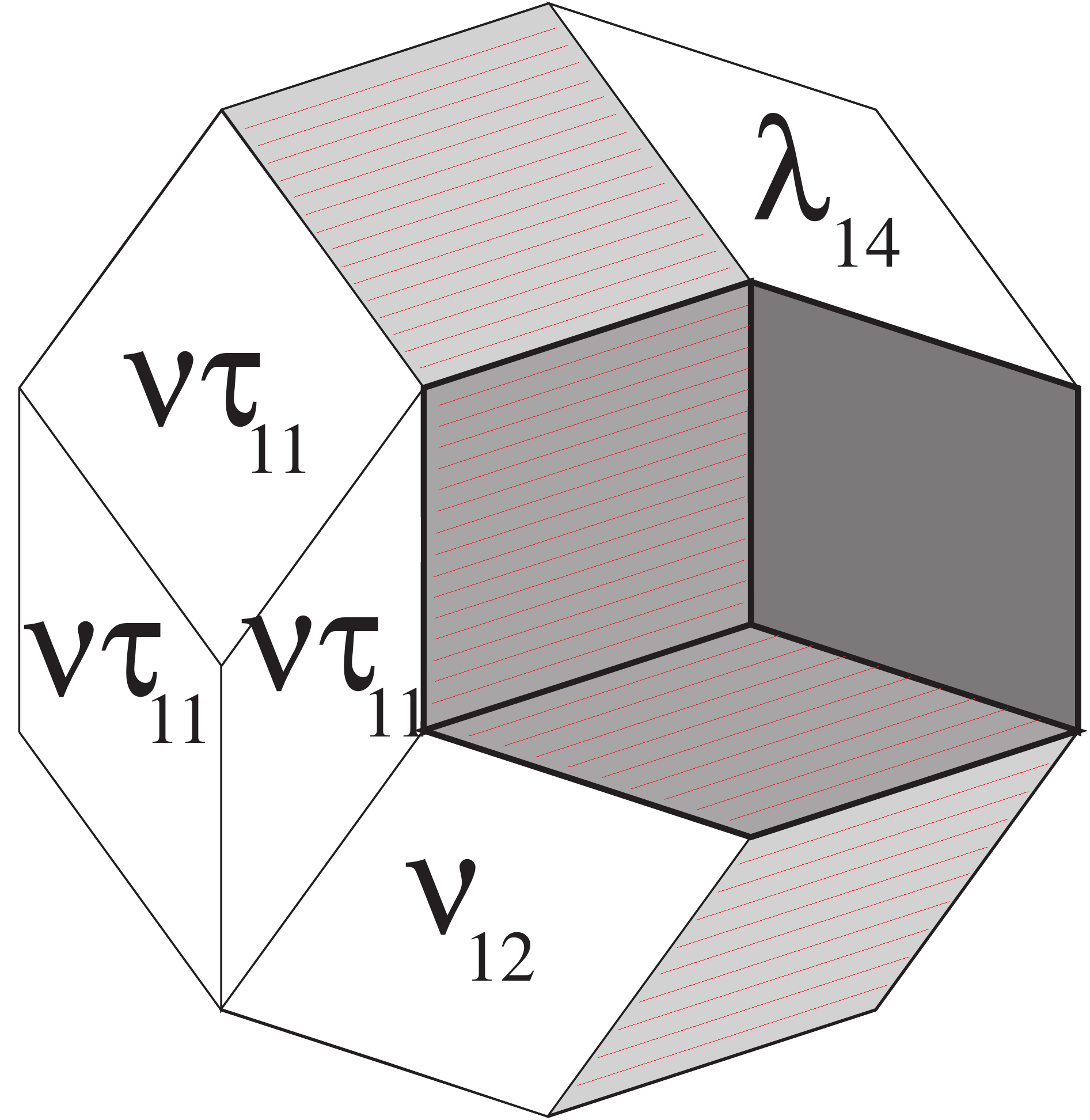}     
  &      \includegraphics[width=0.9in]{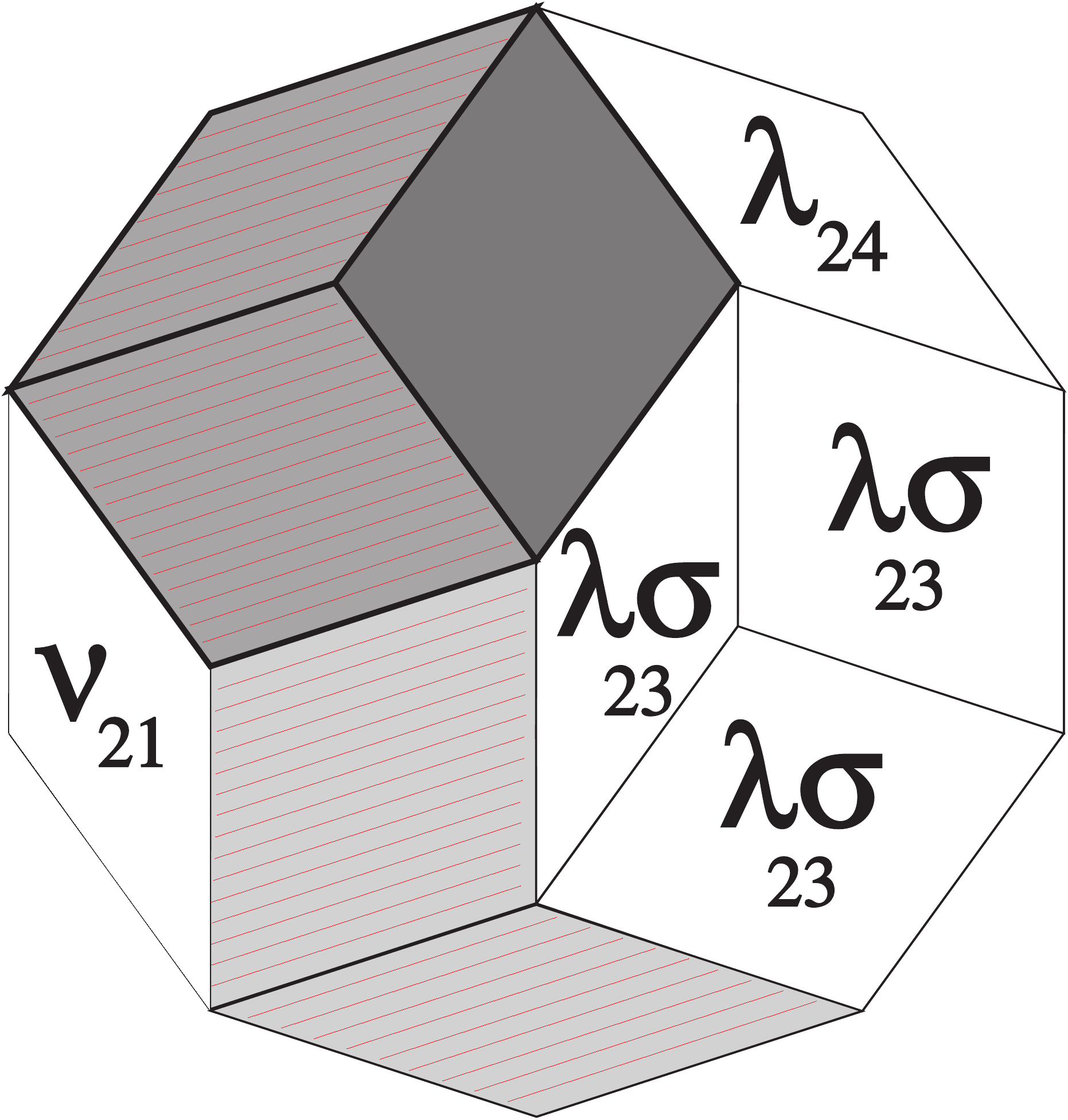}     
&    \includegraphics[width=0.9in]{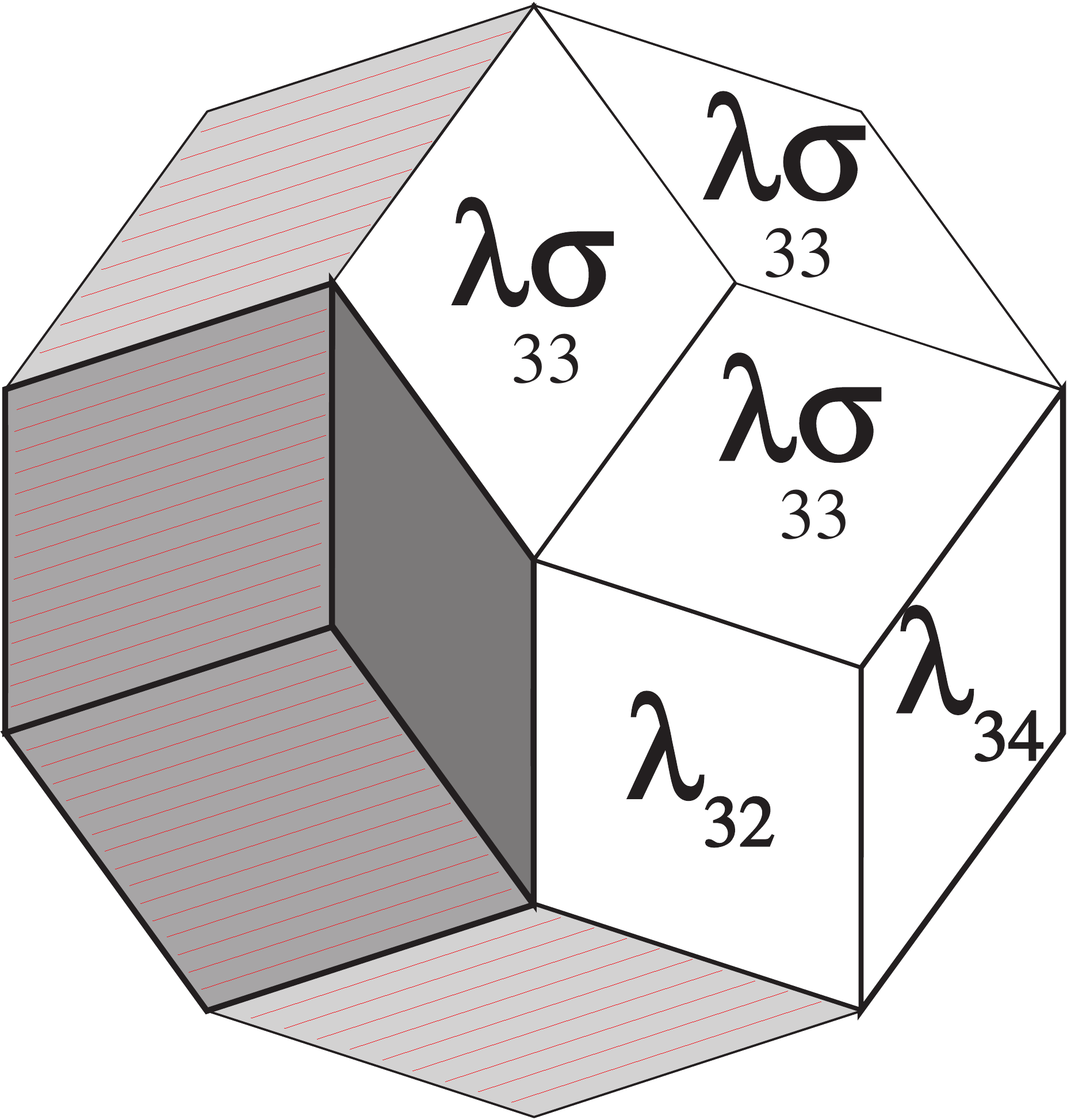}     
&     \includegraphics[width=0.9in]{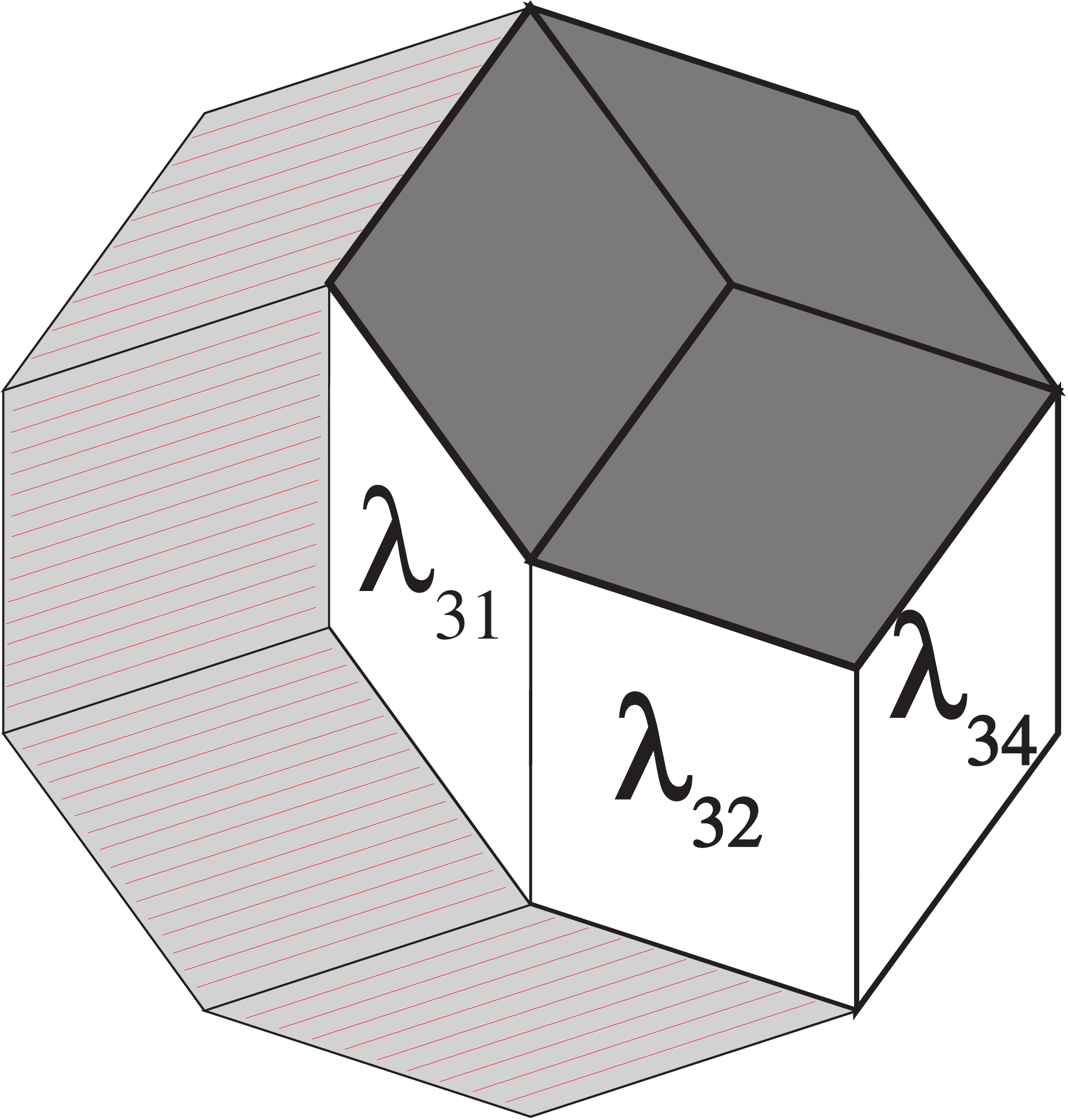}     
&     \includegraphics[width=0.9in]{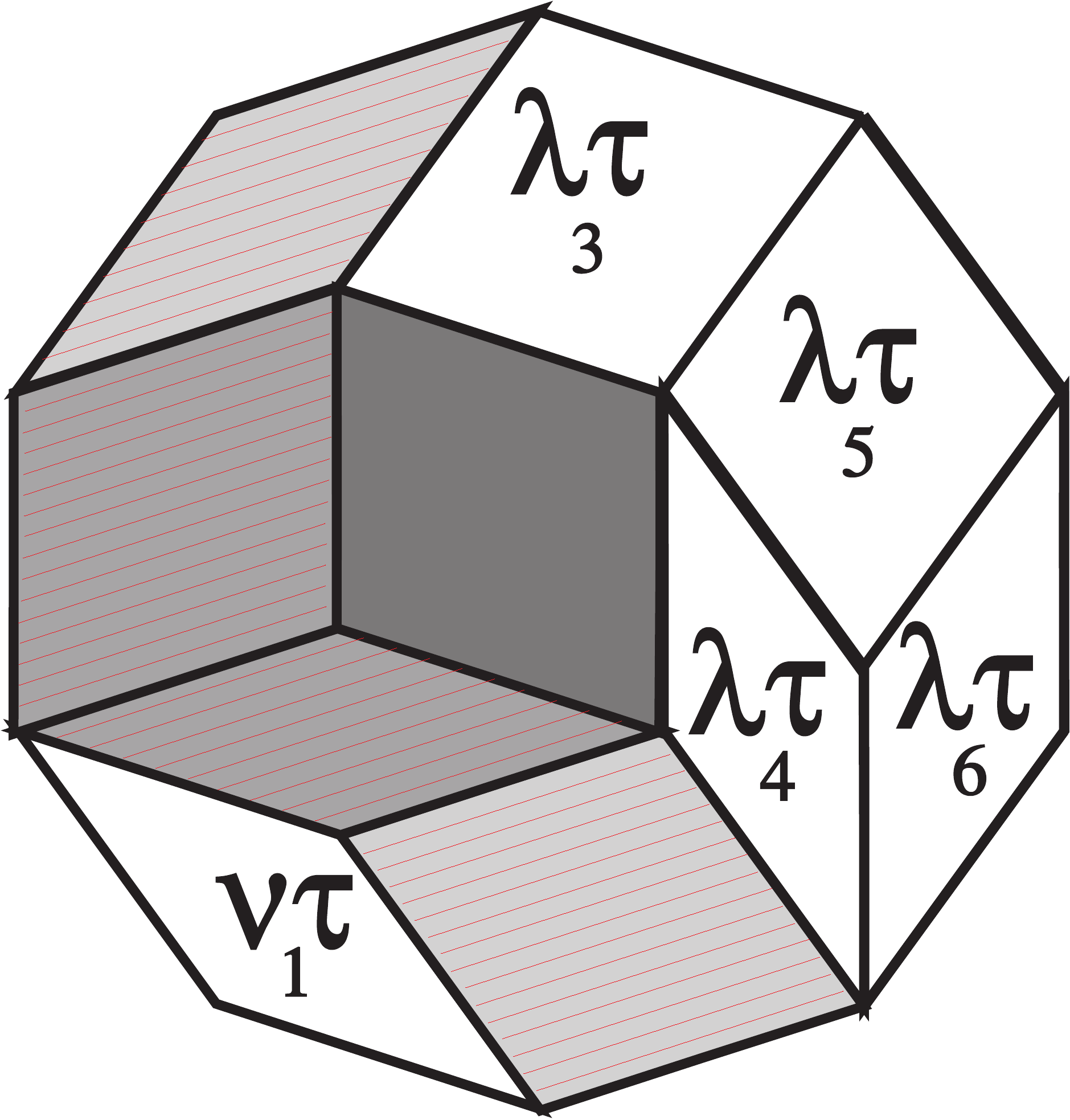} 
\\

\mu_{13} & \mu_{22} & \mu_{31}  &  \lambda_{33}  & \mu_{2}  \\

      \includegraphics[width=0.9in]{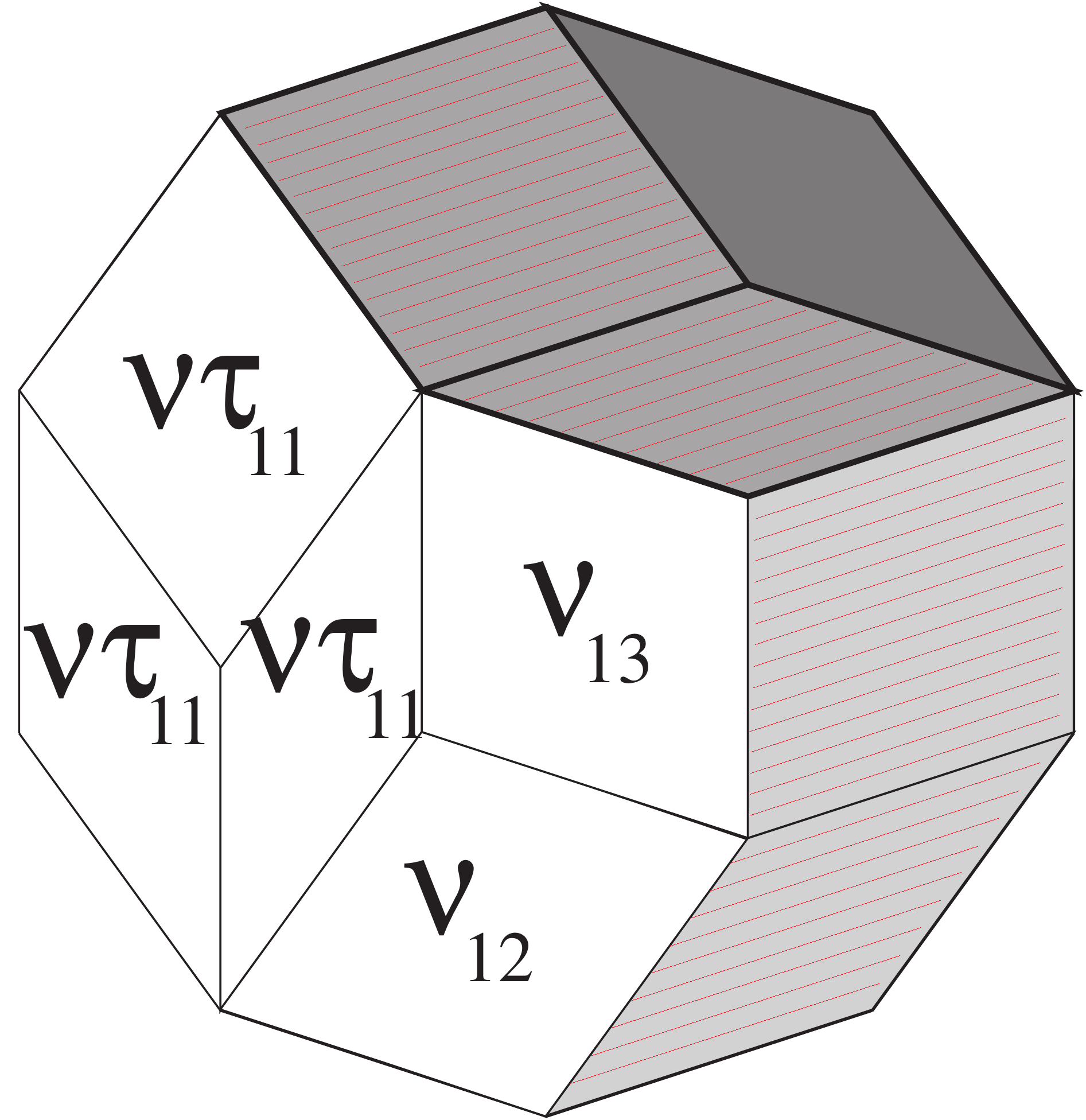}     
  &      \includegraphics[width=0.9in]{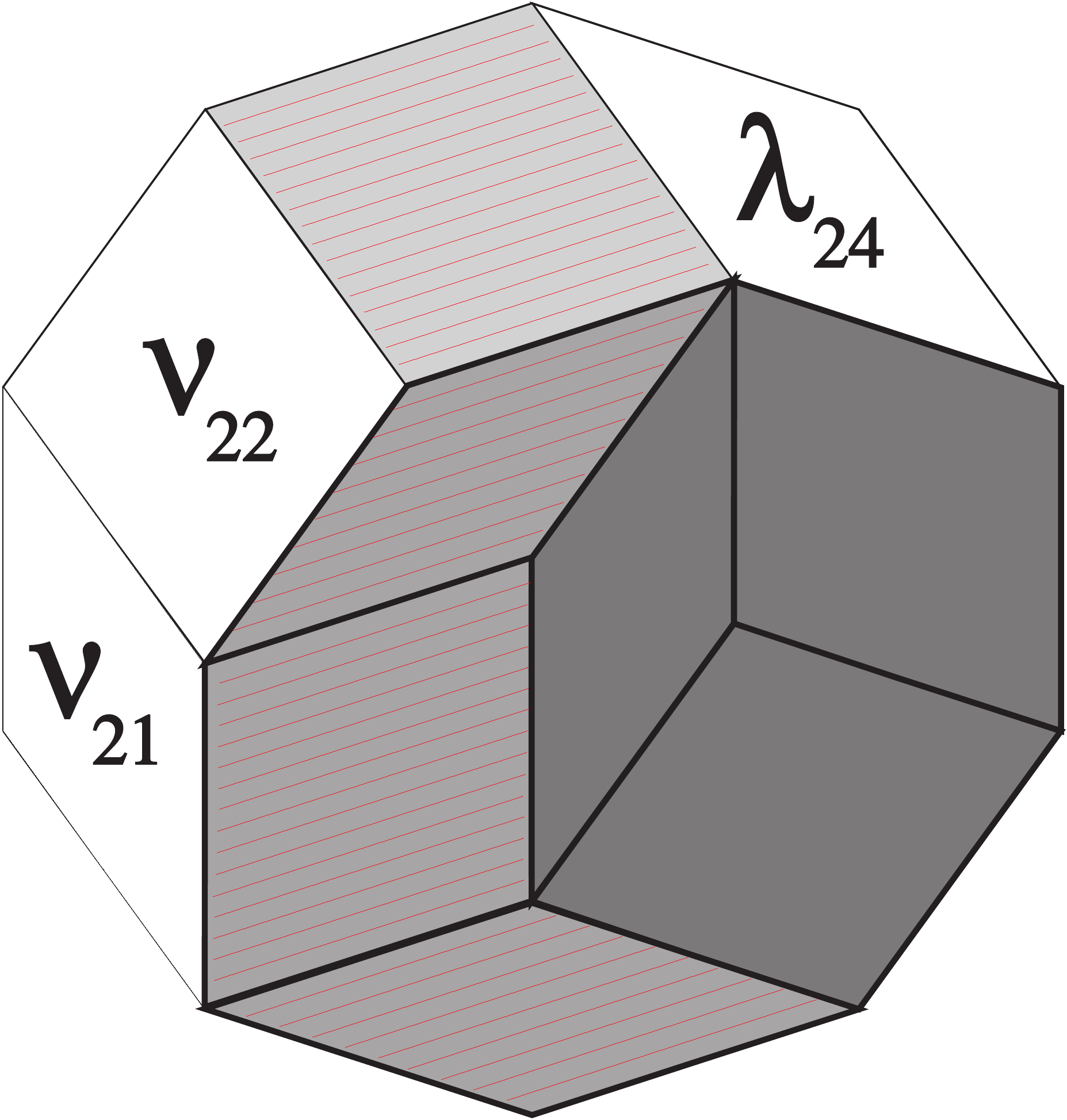}    
 &    \includegraphics[width=0.9in]{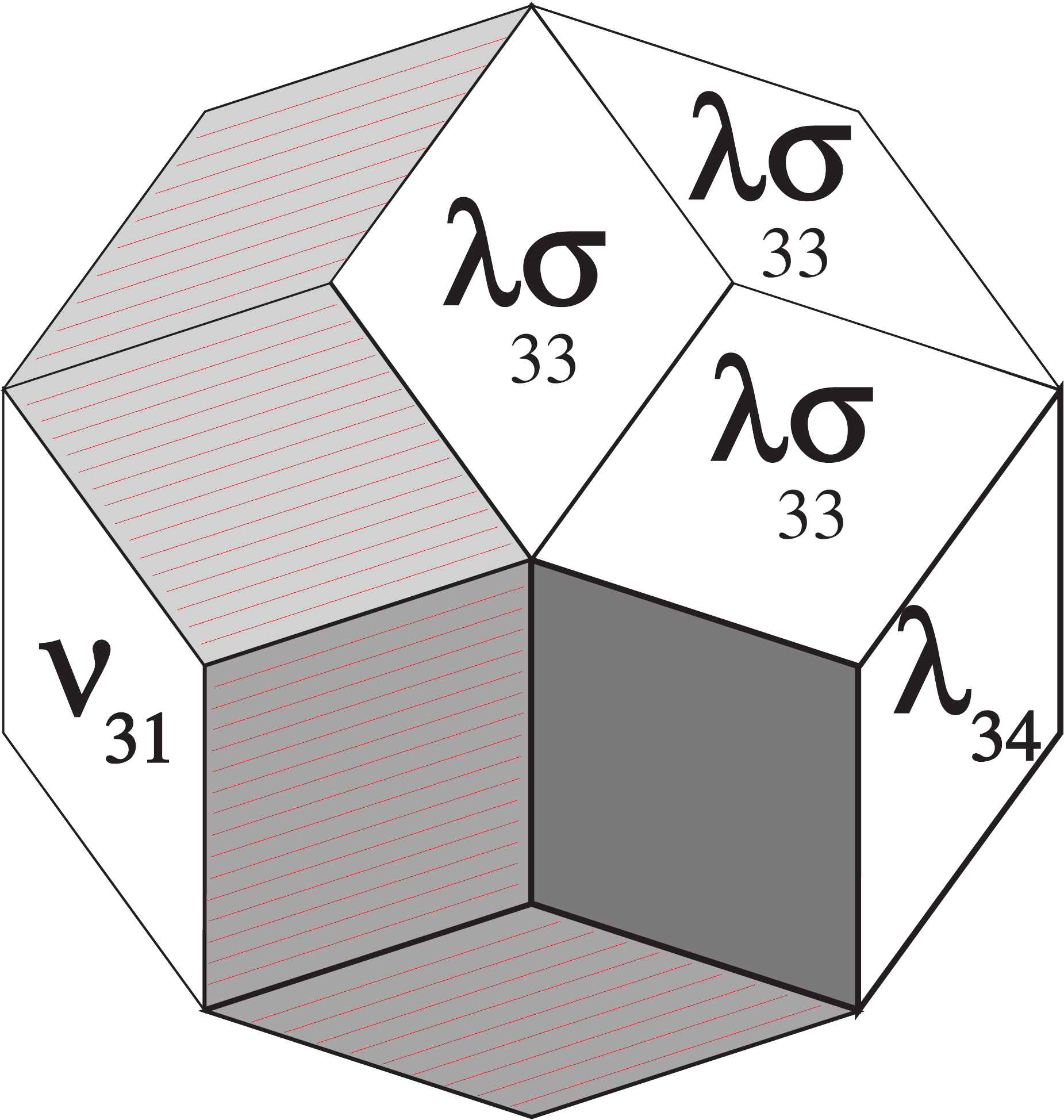}     
&     \includegraphics[width=0.9in]{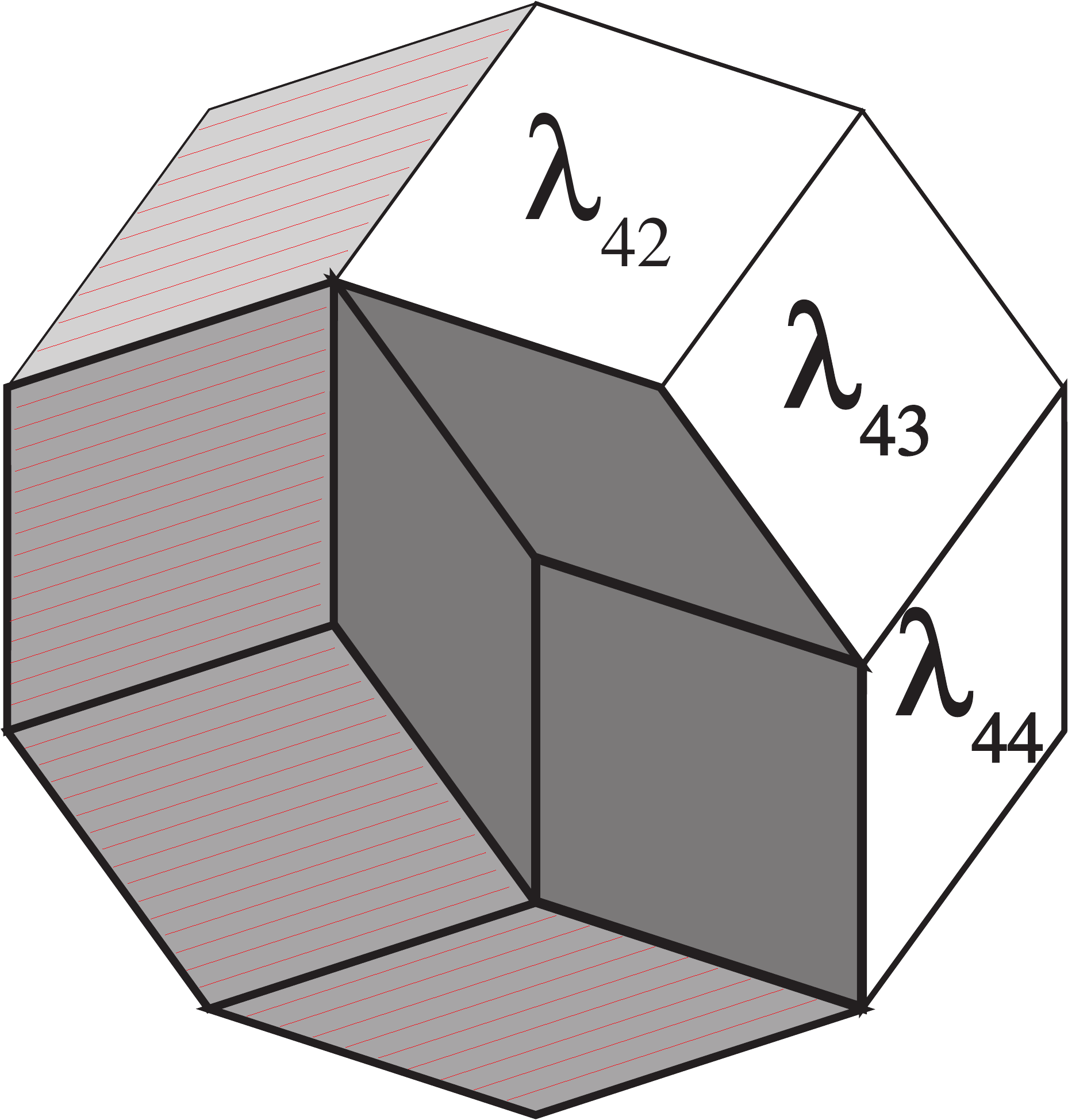}     
&     \includegraphics[width=0.9in]{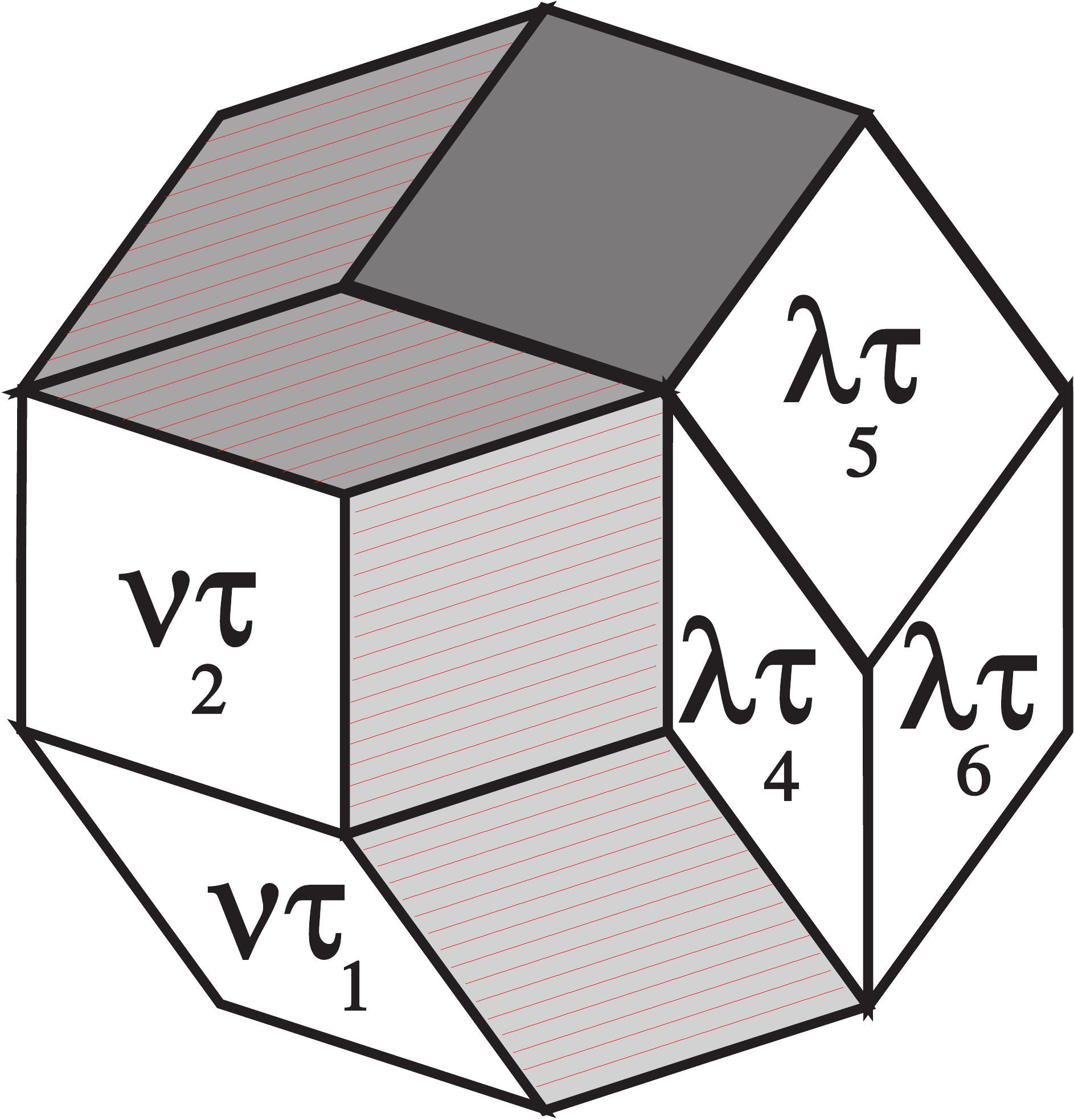} 
\\

\mu_{14} & \mu_{23} & \mu_{32}  & \mu_{41} & \mu_{3}    \\

     \includegraphics[width=0.9in]{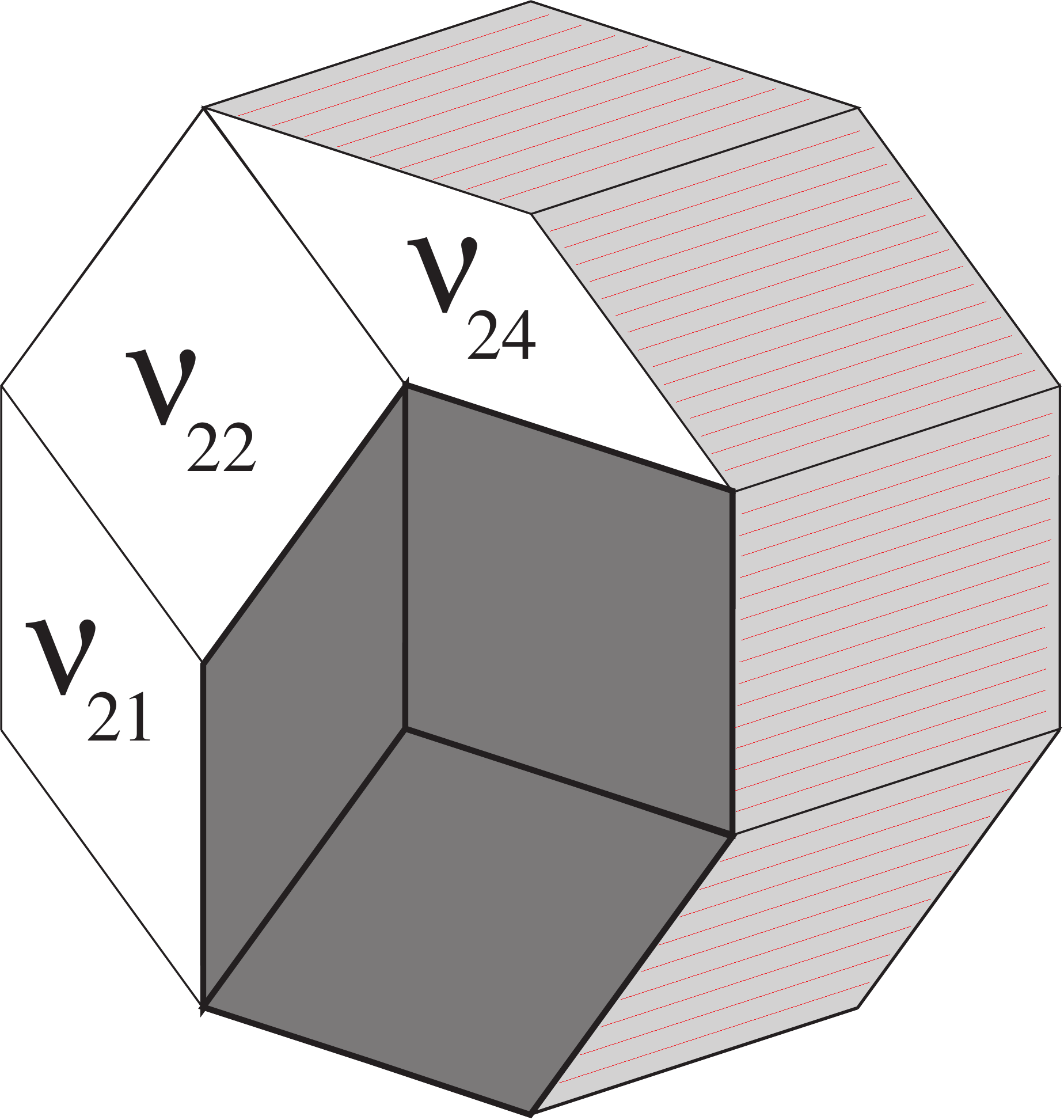}     
  &      \includegraphics[width=0.9in]{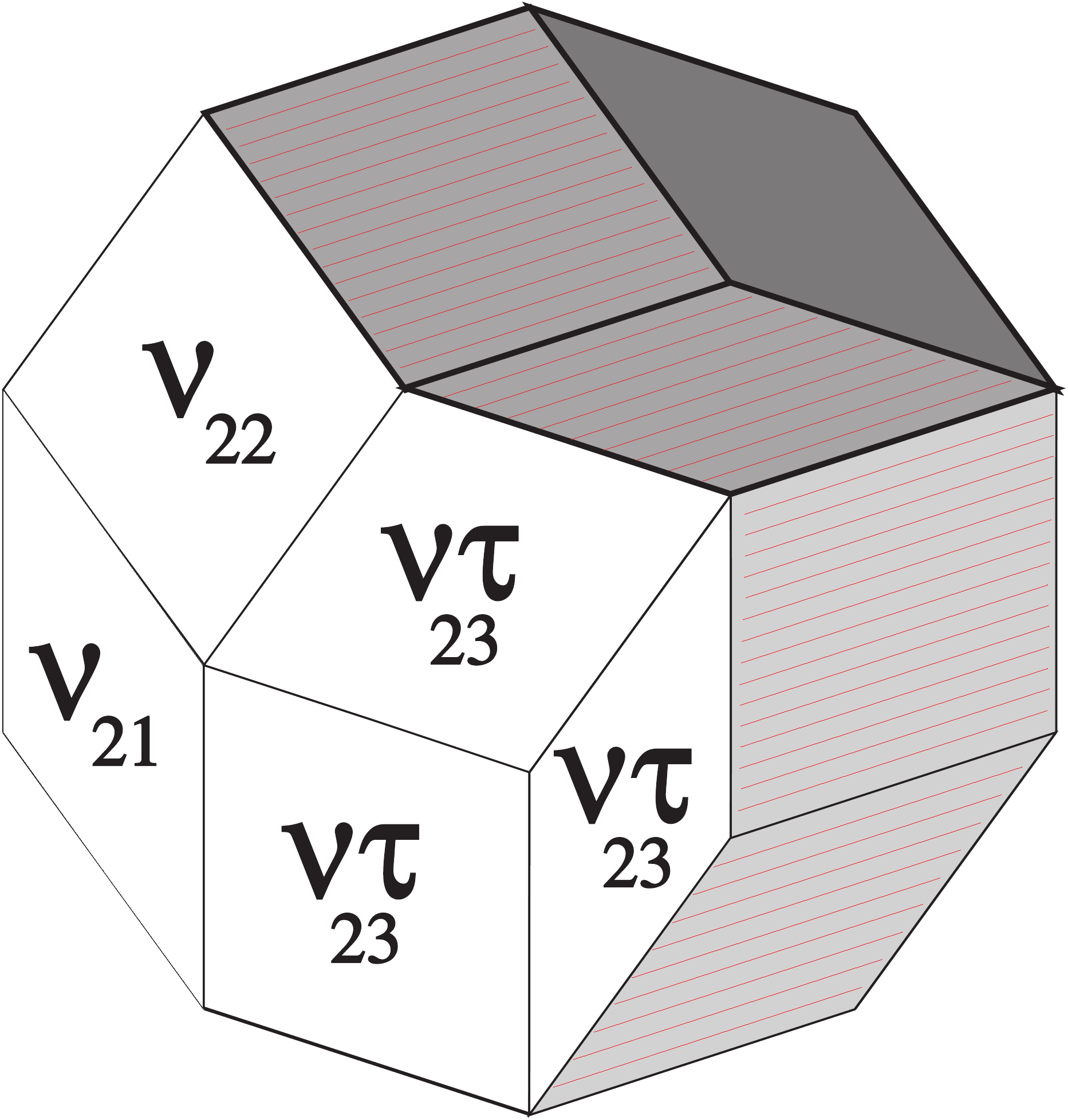}   
 &    \includegraphics[width=0.9in]{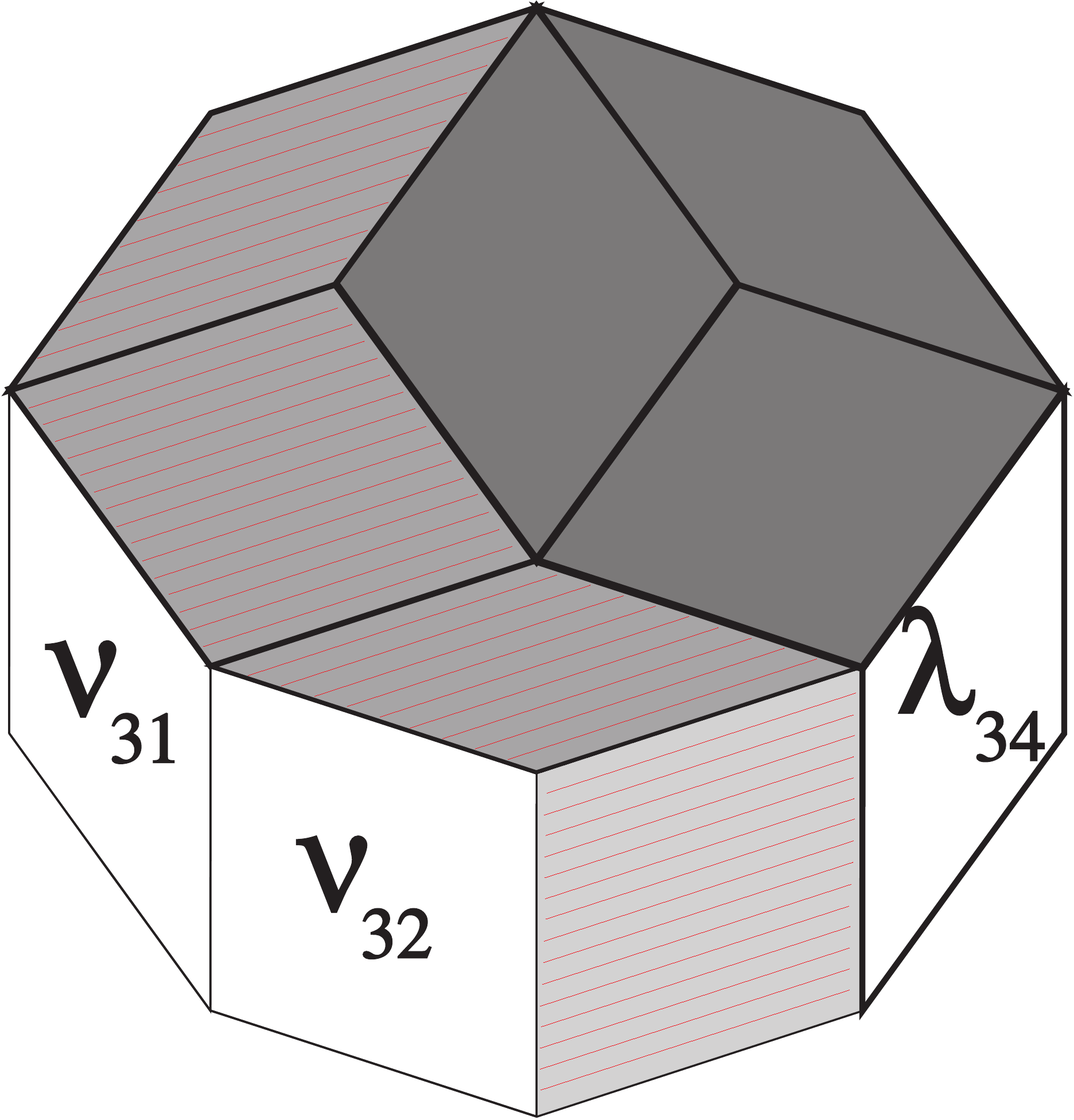}     
&     \includegraphics[width=0.9in]{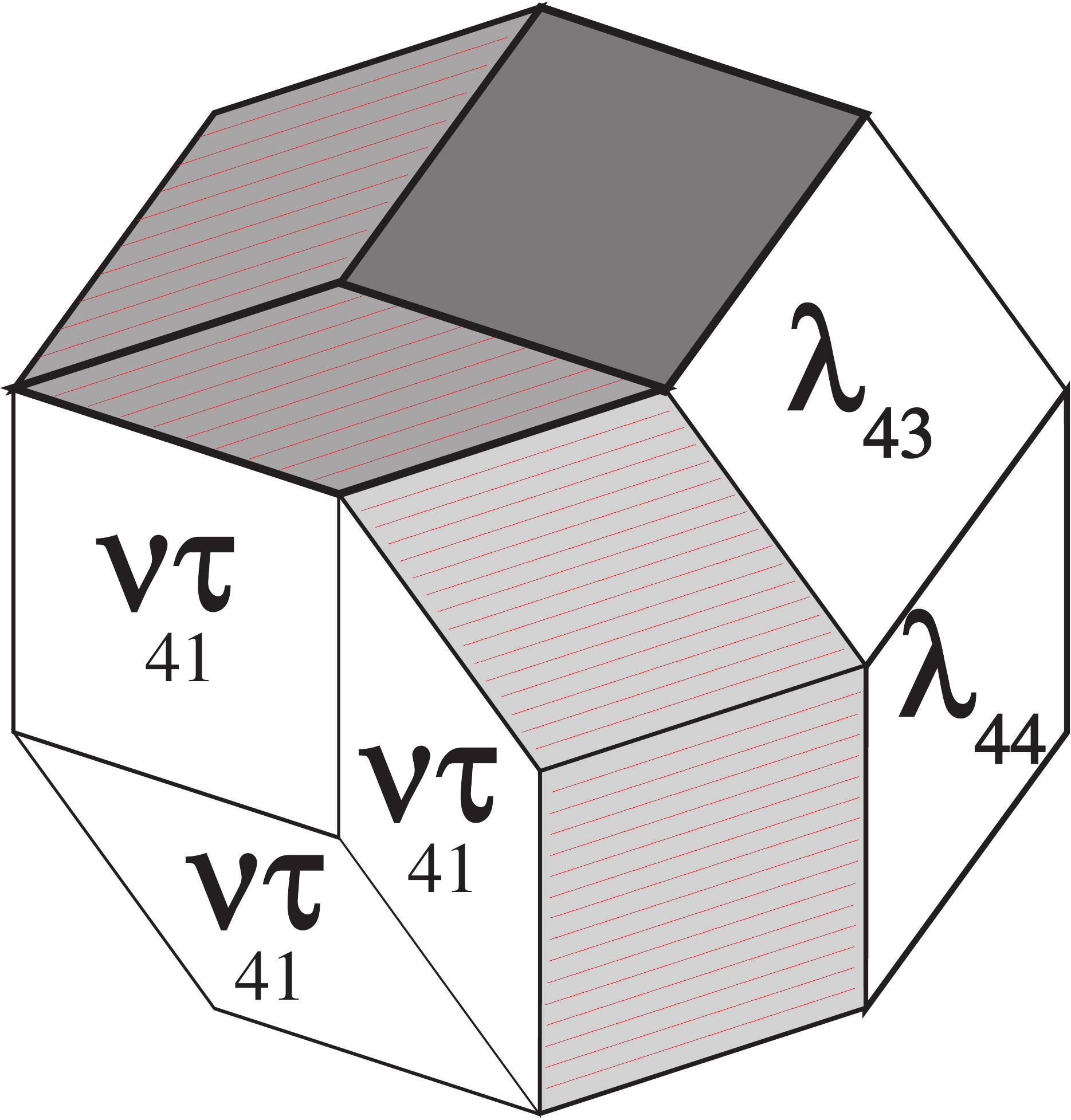}     
&     \includegraphics[width=0.9in]{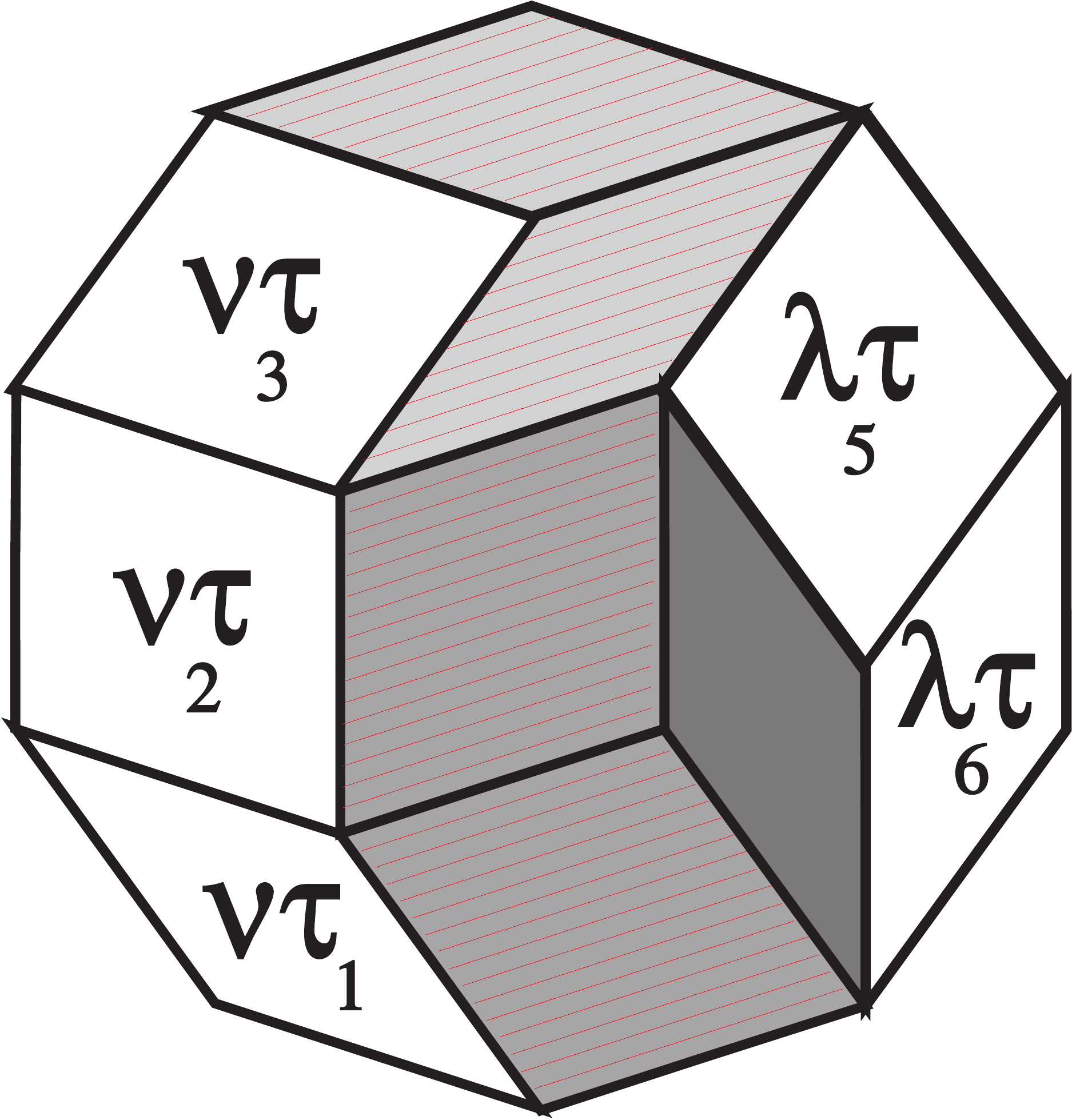} 
\\

\nu_{23} & \mu_{24} & \mu_{33}  & \mu_{42}  & \mu_{4}     \\

     \includegraphics[width=0.9in]{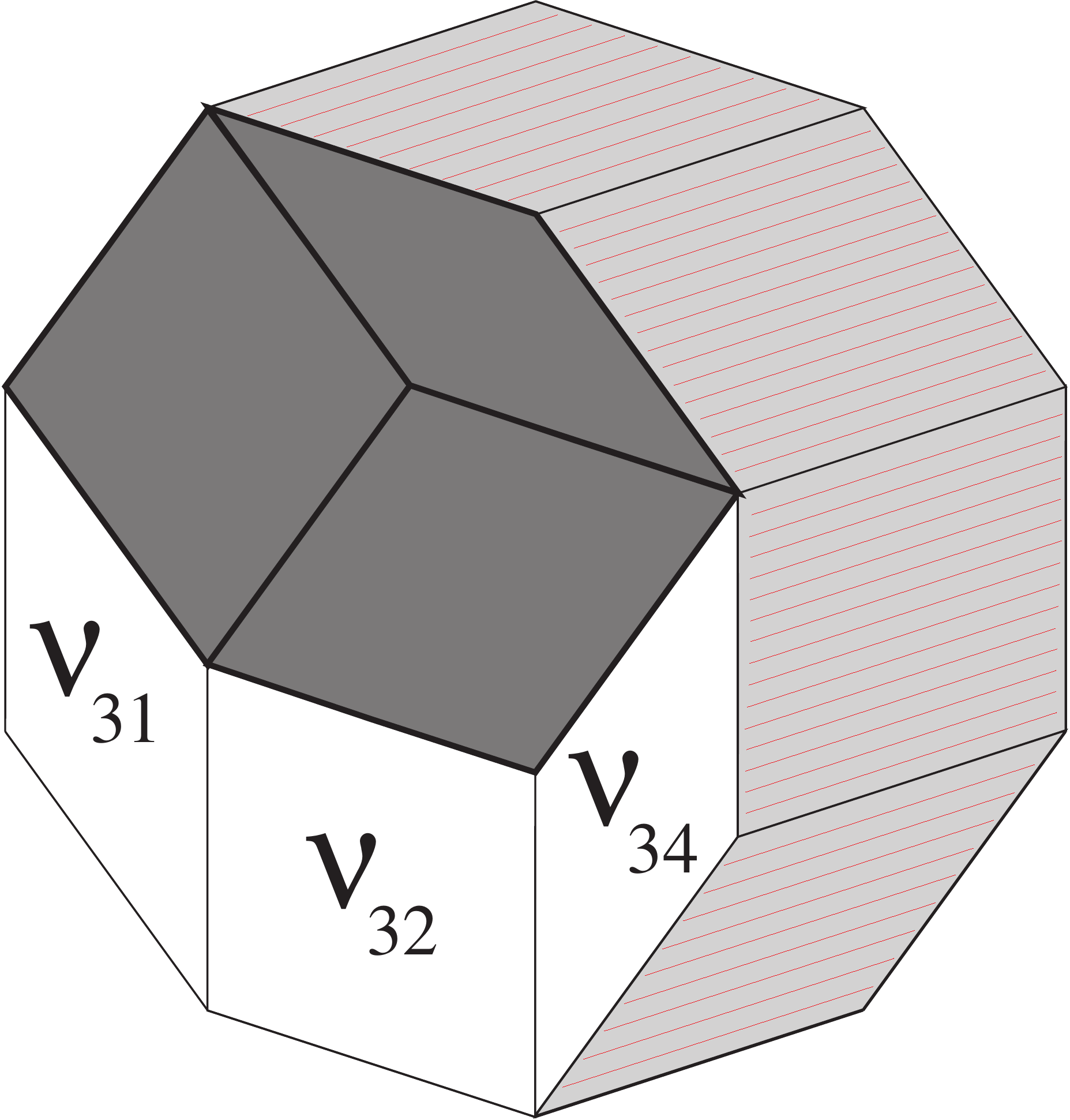}     
  &      \includegraphics[width=0.9in]{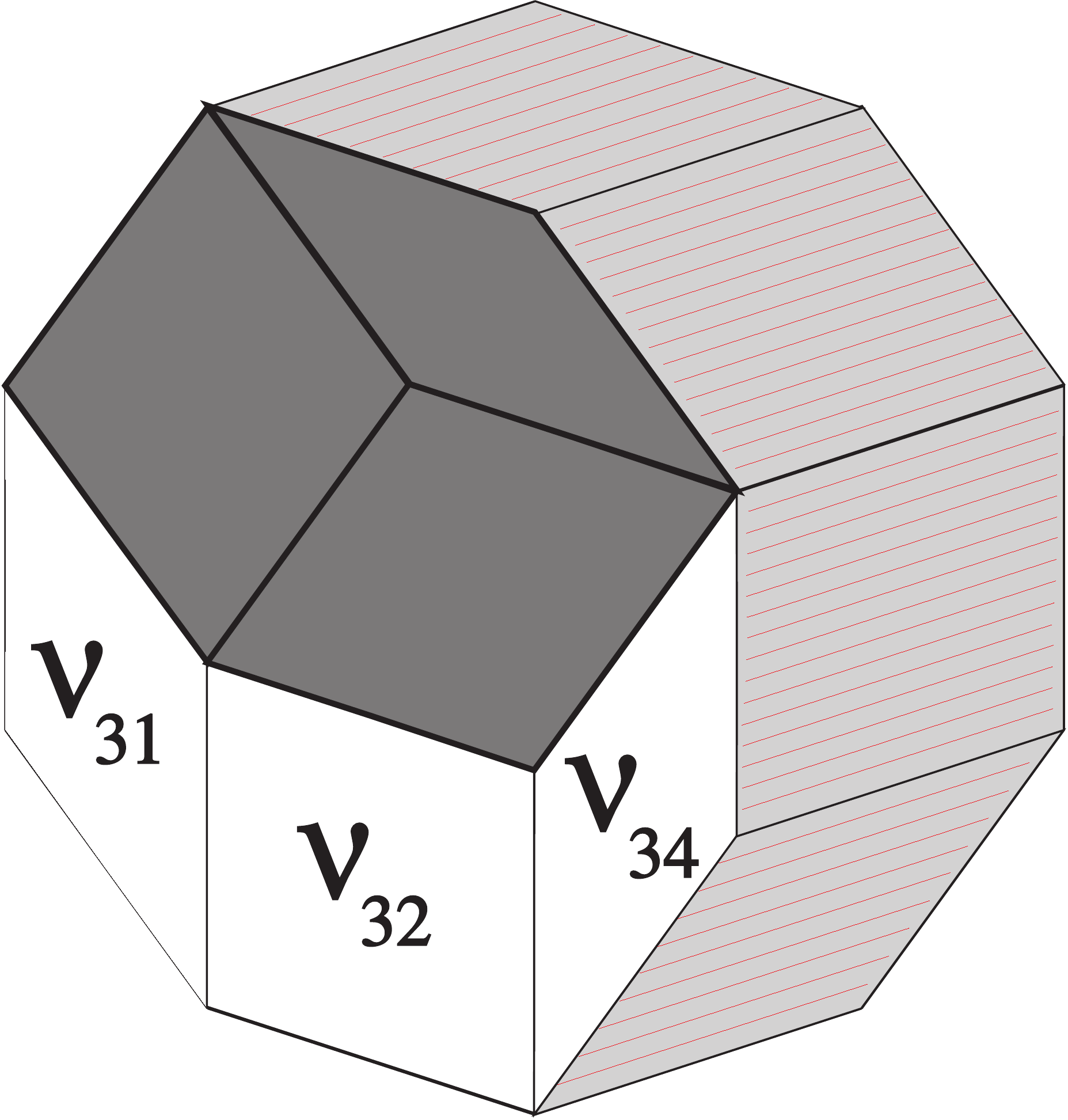}    
 &    \includegraphics[width=0.9in]{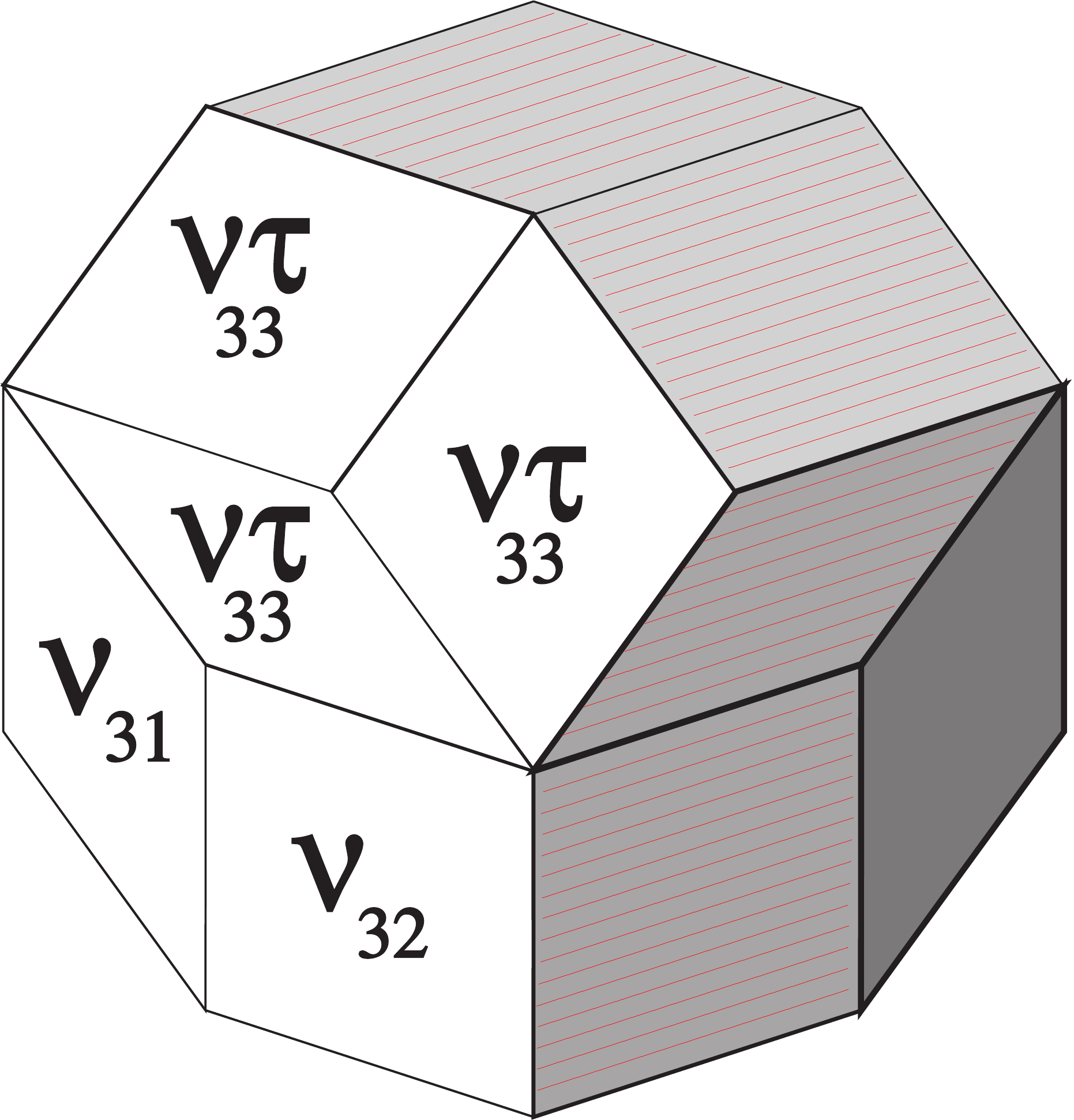}     
&     \includegraphics[width=0.9in]{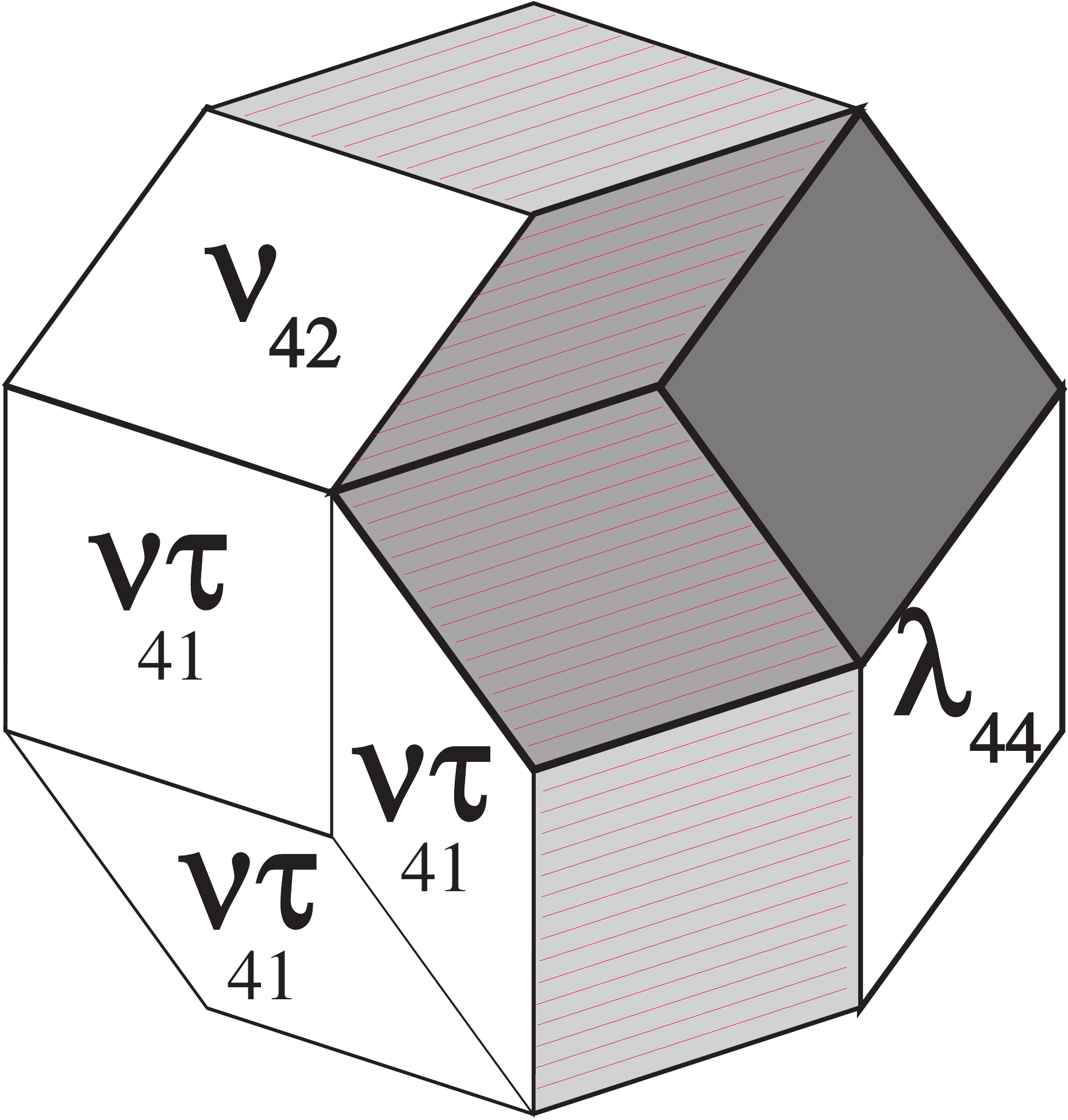}     
&     \includegraphics[width=0.9in]{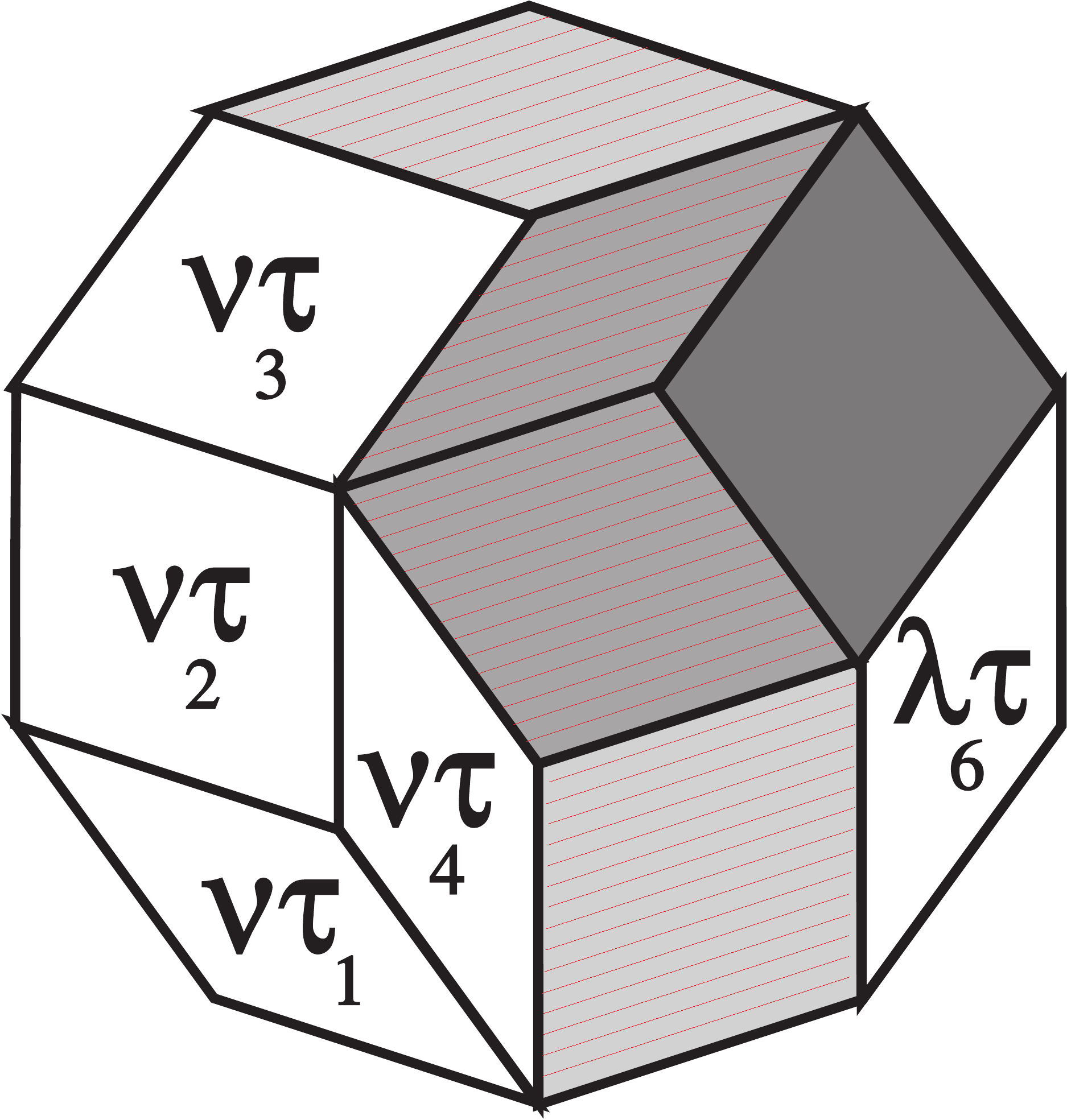} 
\\

\nu_{33} & \nu_{33} & \mu_{34}  & \mu_{43}& \mu_{5}      \\

      \includegraphics[width=0.9in]{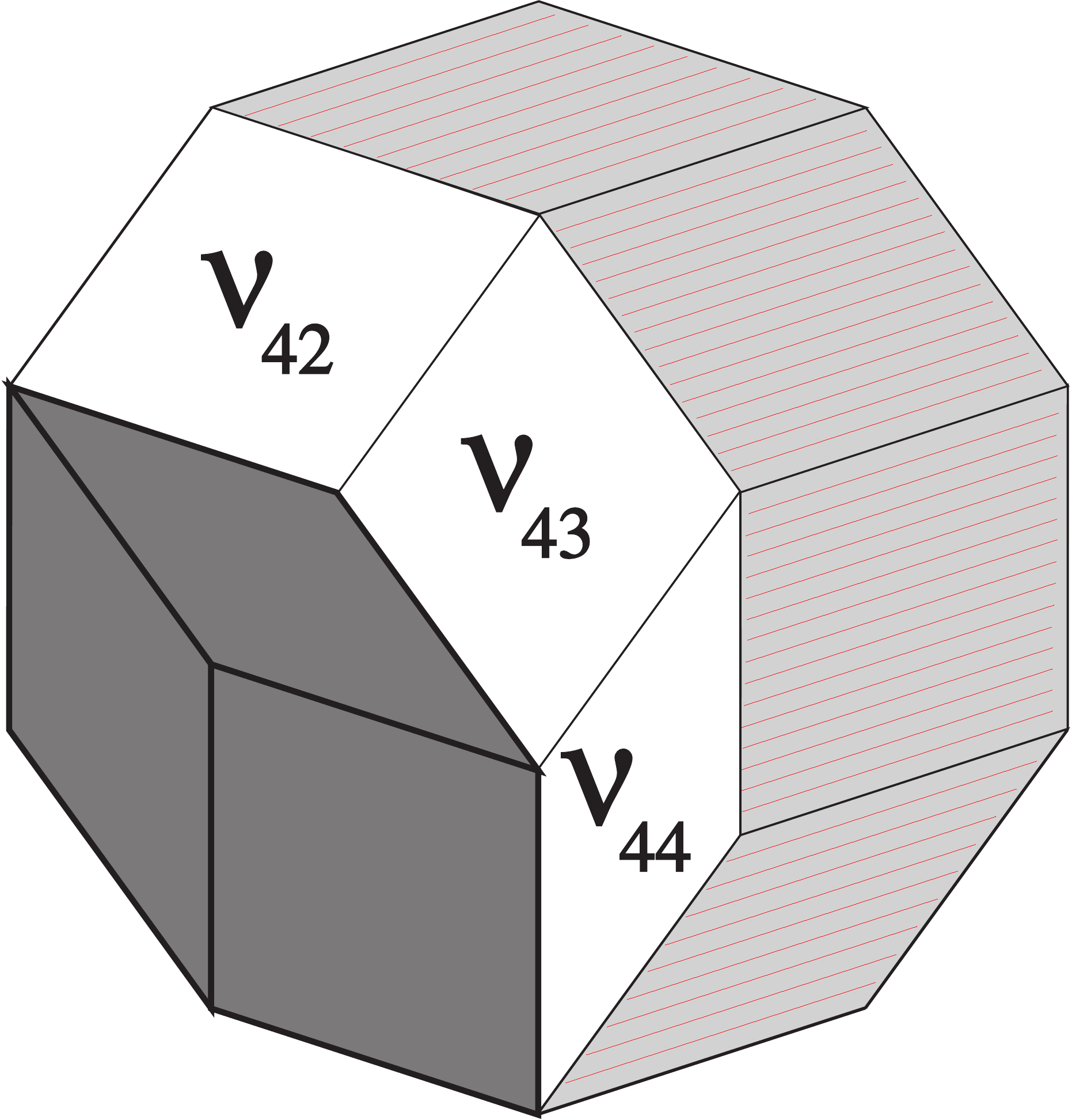}     
  &      \includegraphics[width=0.9in]{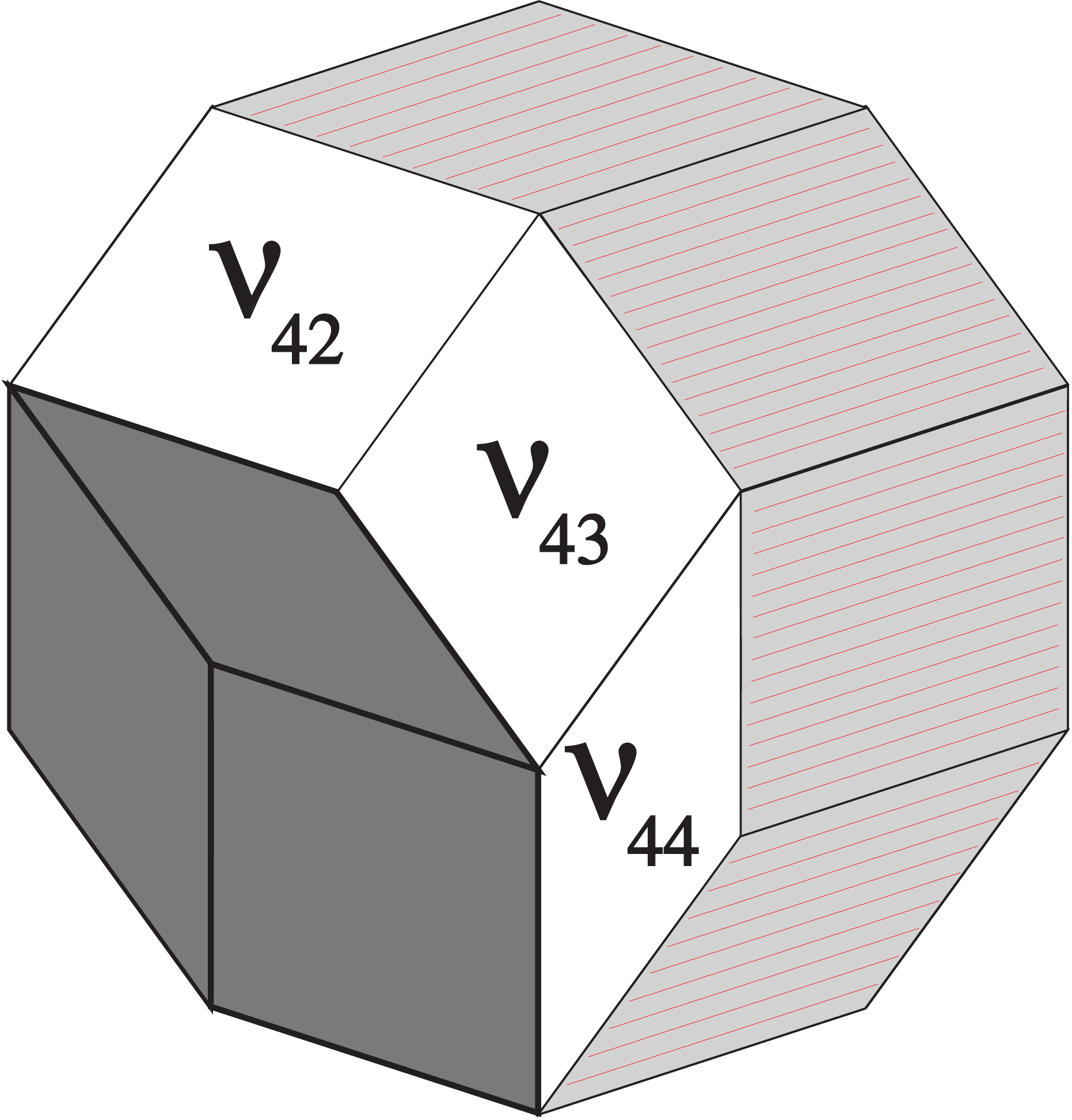}     
 &    \includegraphics[width=0.9in]{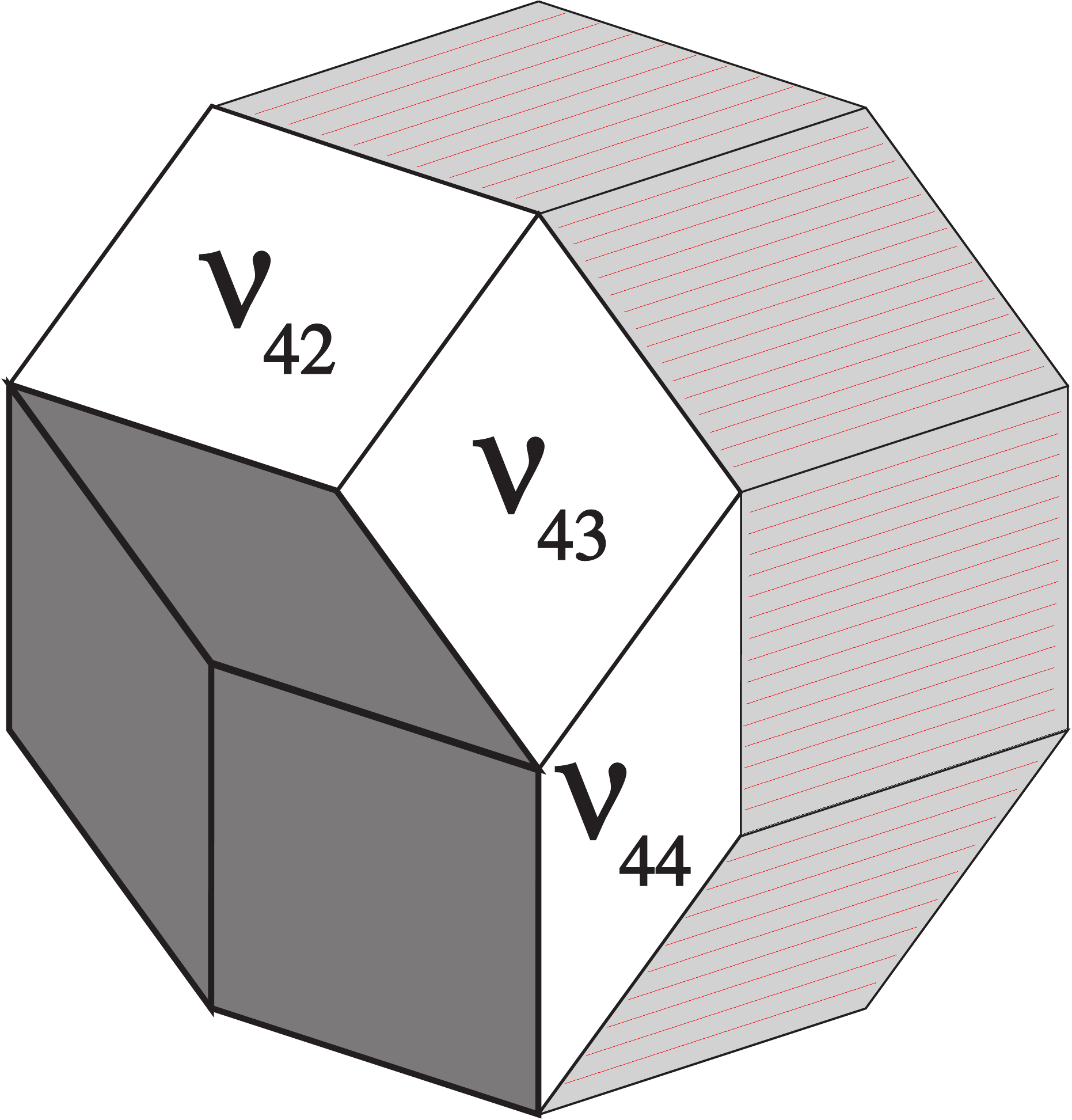}     
&     \includegraphics[width=0.9in]{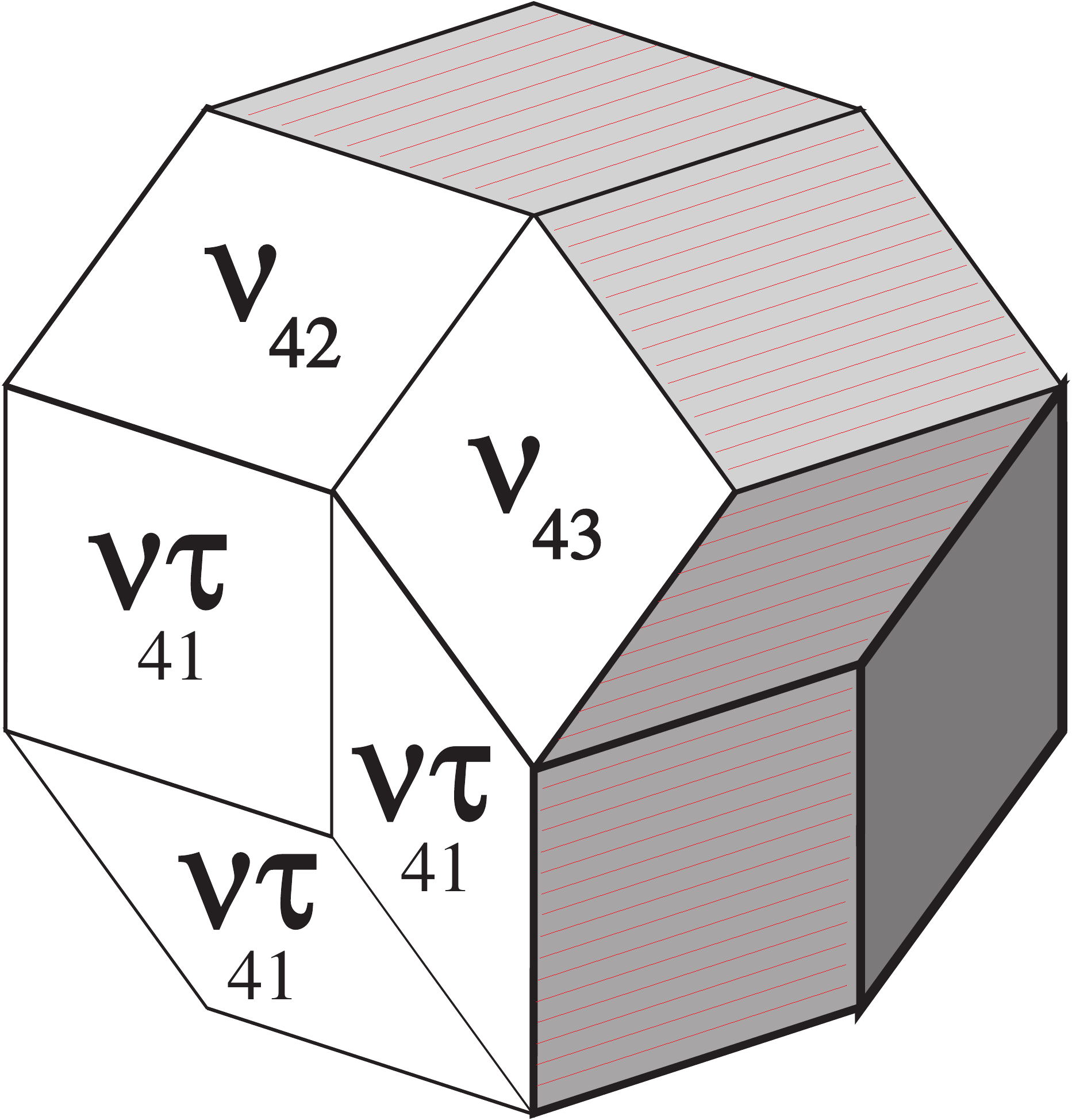}     
&     \includegraphics[width=0.9in]{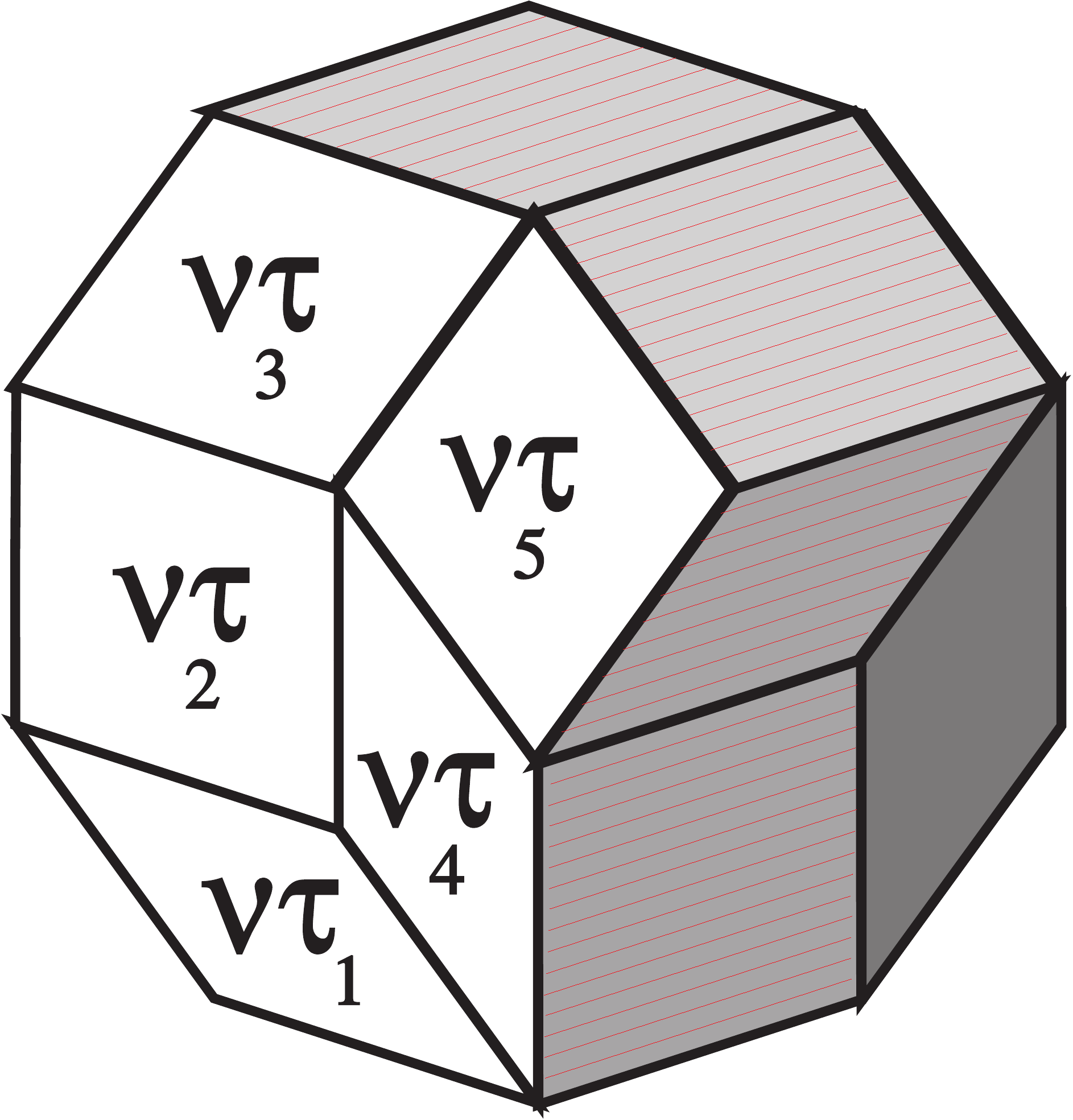} 
\\

\nu_{41} & \nu_{41} & \nu_{41} & \mu_{44}   & \mu_{6}    \\
 \end{array}
%\right)
\]
%\vskip -0.5cm
\begin{figure} \includegraphics[width=0in]{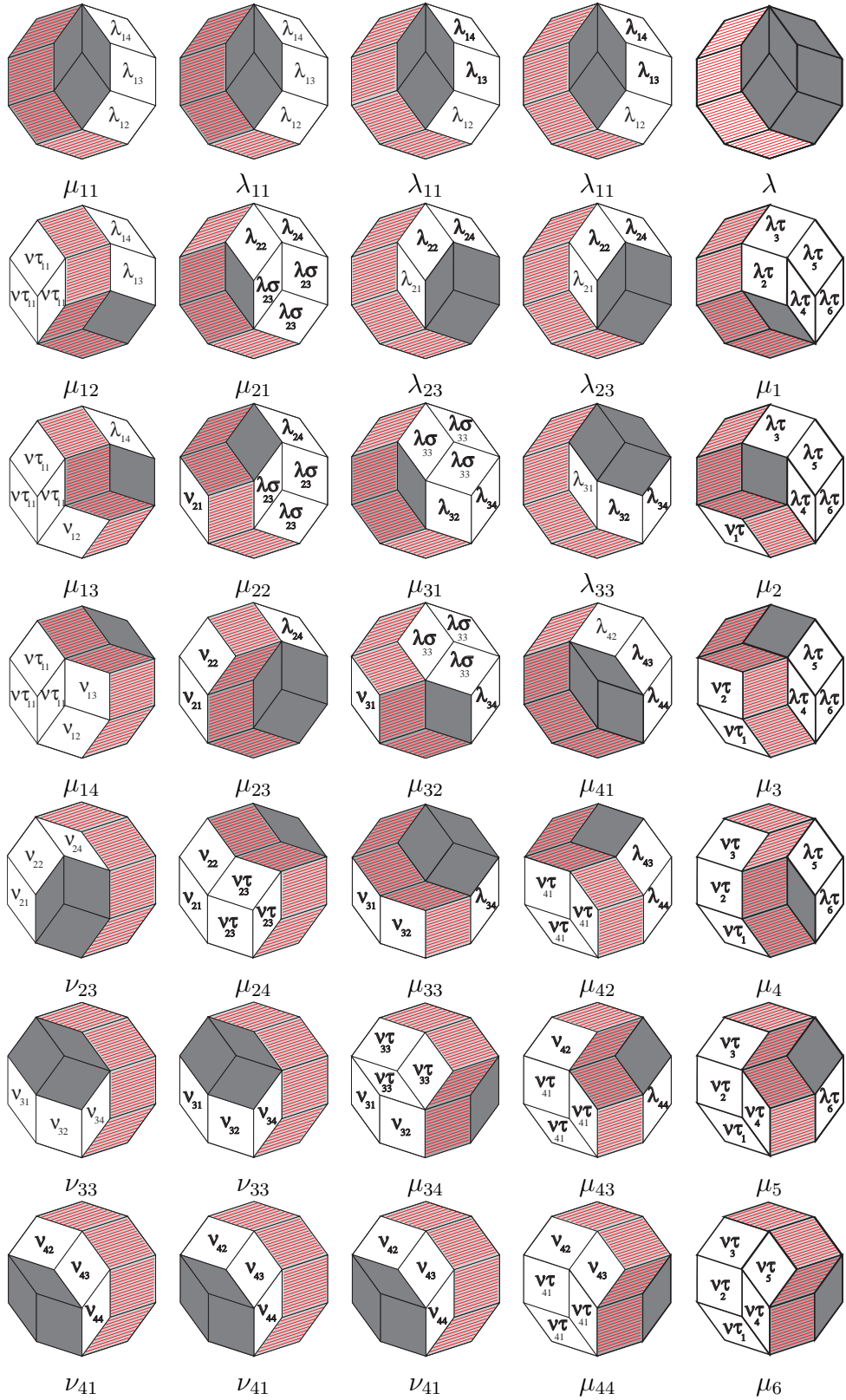}
\vskip -1cm\caption{Inductive structure of the 5-cube source}
\label{fig:20}
\end{figure} 
 
 \begin{figure}[htbp] %  figure placement: here, top, bottom, or page
    \centering
    \includegraphics[width=6in]{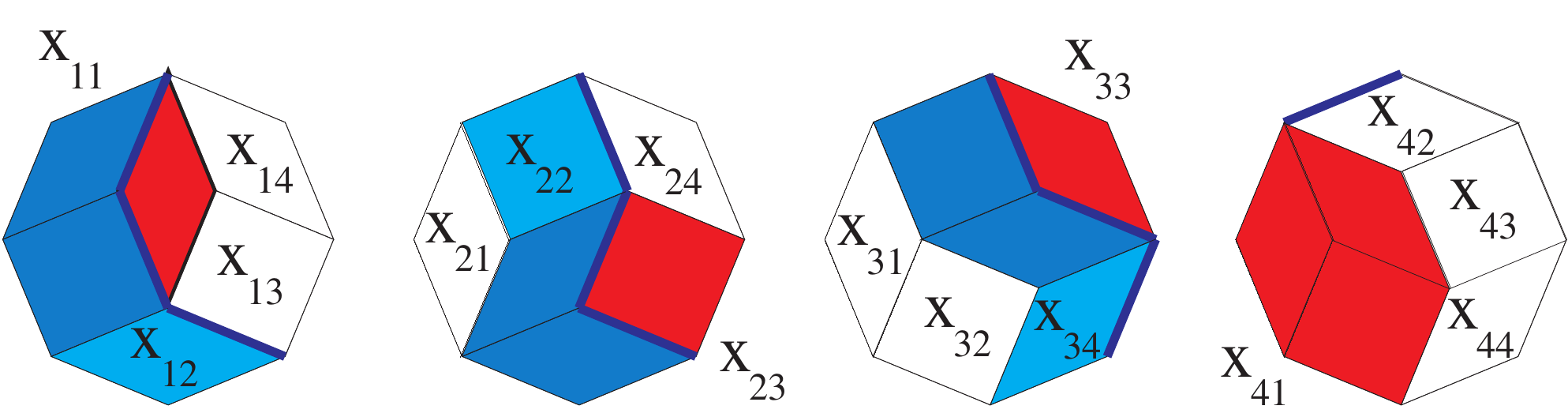} 
    \caption{Labeling of the 4-cube, used for the 5-cube, arising from the 3-cube. The first of the pair of subscript labels refers to the diagram number; the second  to the order in which deformations must be applied to 2-cubes/3-cubes. The subscripts for the 5-cube have been written in reverse order.}
    \label{fig:21}
 \end{figure}

Hence if we draw the 
{\it  two}
possible paths from the 2-source to the 2-target, across the 3-source and
3-target, we obtain what we call the 
{\it
octagon of octagons.
}
Note the bilateral symmetry of the figure, as well as its antipodal
symmetry. If all symbols were to be deleted, the octagon of octagons shows that each figure, going clockwise around the octagon, differs from the preceding by 
a $ { 3 \pi } / 4 $ rotation. 
 
 In Figure 17 %\ref{fig:17} 
 we show how the 3-cocycle arises from the 2-cocycle by ordering the 
thickenings and shiftings of cubes that occur. The left hand side of the octagon arises by considering the order of 2-cell composition for the 2-cocycle source,
indicated by the numbers, and consecutively moving these across the octagon
via their thickenings into 3-cubes as indicated by the bold-face hexagons.
The target side of the octagon of octagons can be obtained either by
symmetry considerations or by the same process. 
 
It is clear that the unlabelled configurations of the octagon of octagons
can be derived entirely without labelling 
{\it 
any
}
of the faces! Thus the only purpose of the labellings of the cube is to
allow us to write down the cocycle conditions in a linear symbolic form,
without the need for these diagrams.

We may replace the octagon by Figure \ref{fig:18}, which depicts the splitting of the eight 
$ {\cal I}^3 $s in the boundary of $ {\cal I}^4 $ into the source and target sets.
The octagonal ``outline'' we see is the union of the 1-source and 1-target, 
splitting the boundary 2-sphere, which itself is the union of the 2-source 
and 2-target, into these two hemispheres. The 3-source of $ {\cal I}^4 $ is
on the left, the 3-target on the right. The front faces are the 2-source 
$ \sigma_2 [4] $, the back face $ \tau_2 [4] $. To find the upper 4 octagons of the 
octagon of octagons, consecutively ``puncture'' the $ {\cal I}^3 $s having three 
2-faces exposed, thinking of the configuration as cubical soap bubbles. The 
order to proceed is uniquely determined at this dimensional level, and the 
collection of $ {\cal I}^2 $s we see at each stage gives the configurations of
the octagon. We illustrate the puncturing procedure by redrawing the
``octagon of octagons'' as a sequence of 3-dimensional configurations of cubes.
 
Although it is possible to draw the ``source'' of the 4-cell in a single 
configuration -- that of four cubes glued together as half of a tessaract -- by
doing so it is not clear how to write down a  ``canonical'' ordering of the
2- and 3-cubes at each stage of the octagon, as a path of deformations of
the 1-source to the 1-target.

\begin{figure}[htbp] %  figure placement: here, top, bottom, or page
   \centering
   \includegraphics[width=3.5in]{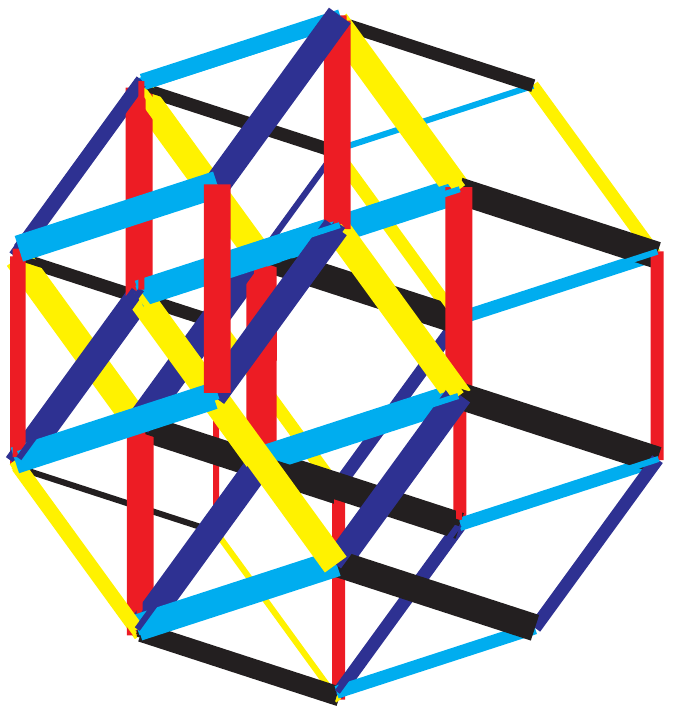} \ = \   \includegraphics[width=3.5in]{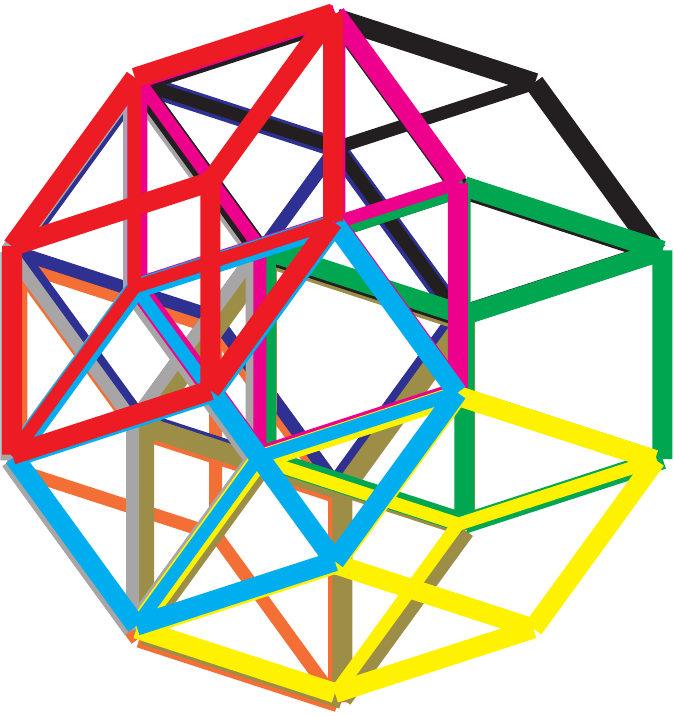}
       \caption{The 5-cube, $\psi_3[5]$. Left: Five edge-types, coloured.  Right:  ten 3-cubes, coloured. We present the interior structure of the 5-cube 3-source 3-cell in both ways for the benefit of the reader: the 10 colours refer to the number of ways three distinct edge colours can be chosen from five. Parallel edges are labeled from 1 to 5, clockwise from the top of the figure. Thus `145' identifies the 3-cube at the bottom right (coloured yellow). The 0-source and target are respectively top and bottom; the 1-source and target are left and right in each decagon.}
      \label{fig:22}
\end{figure}
\begin{figure}[htbp] %  figure placement: here, top, bottom, or page
   \centering
\includegraphics[width=3.5in]{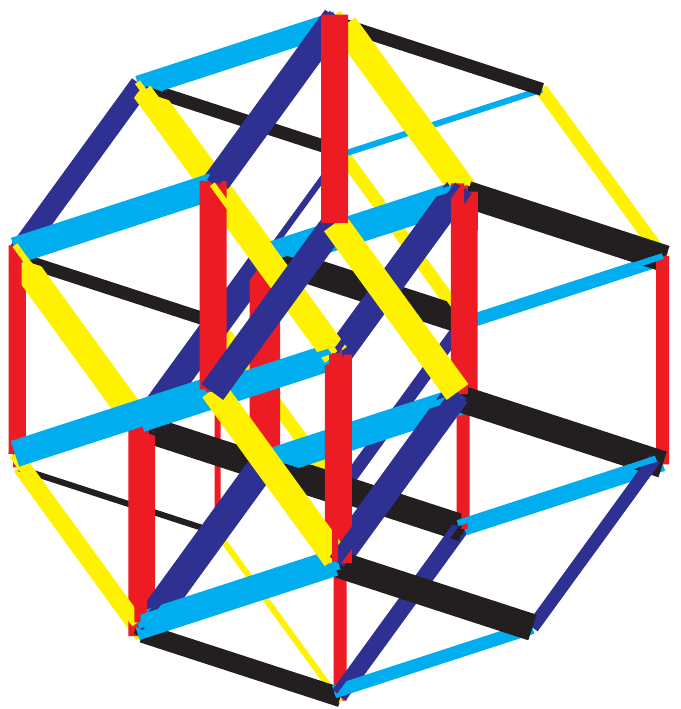} \ = \    \includegraphics[width=3.5in]{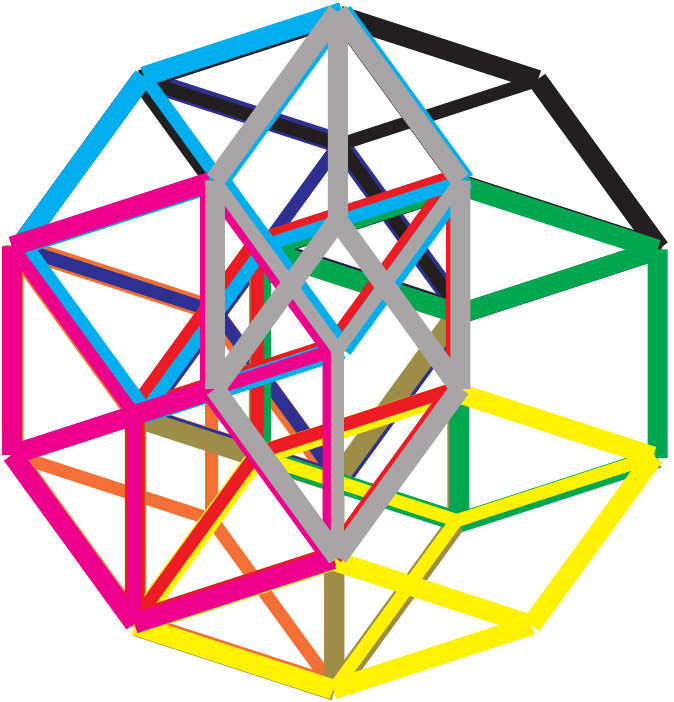} 
    \caption{   --0000($ \psi_3[5]$): front cells $ \psi_2[5]$, rear cells $ \omega_2[5]$}
   \label{fig:23}
\end{figure}
\begin{figure}[htbp] %  figure placement: here, top, bottom,, rear cells  or page
   \centering  
\includegraphics[width=3.5in]{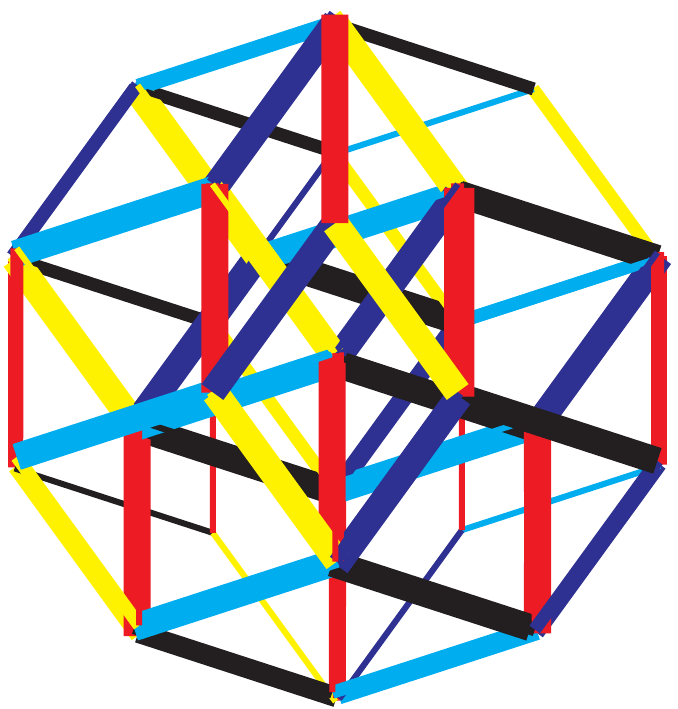} \ = \      \includegraphics[width=3.5in]{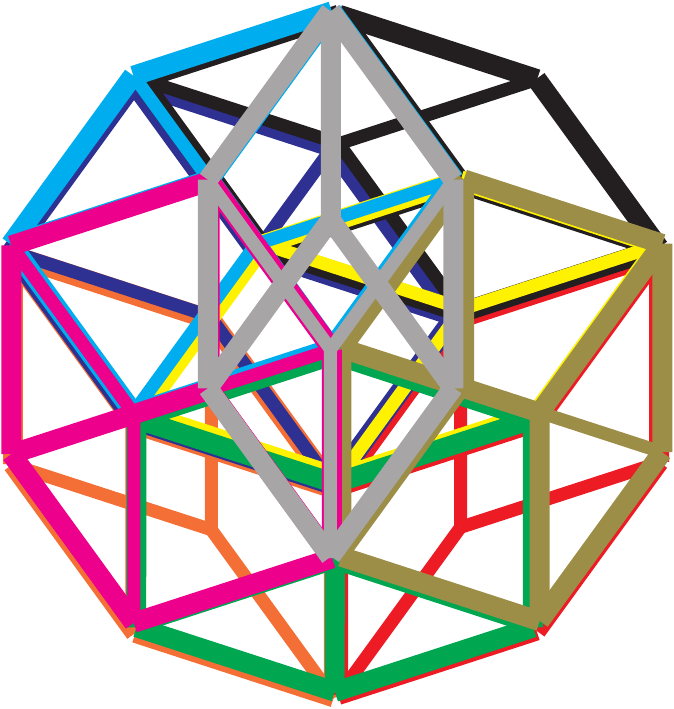} 
  \caption{ 0+000[--0000($ \psi_3[5]$)]}  
   \label{fig:24}
\end{figure}
%\pagebreak

\begin{figure}[htbp] %  figure placement: here, top, bottom, or page
   \centering
   \includegraphics[width=3.5in]{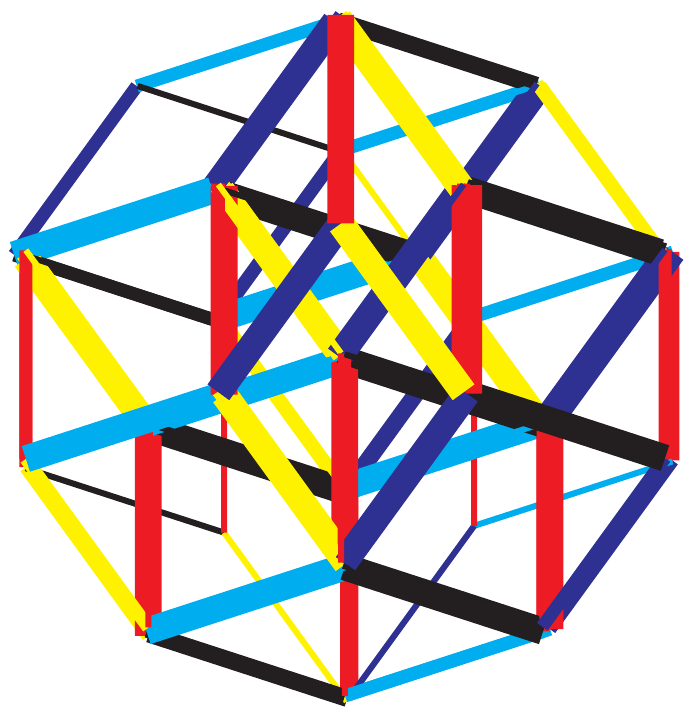} \ = \    \includegraphics[width=3.5in]{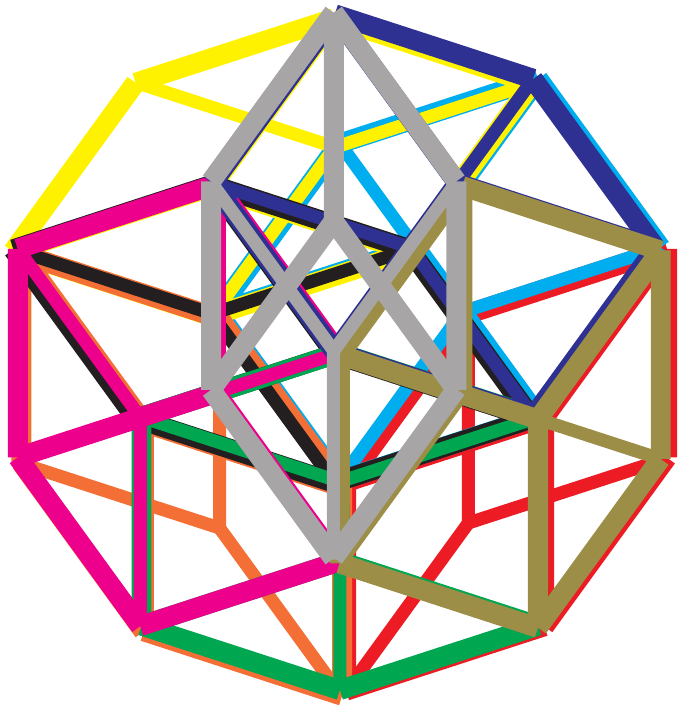}
  \caption{ 00--00(0+000[--0000($ \psi_3[5]$)])}  
      \label{fig:25}
\end{figure}
\begin{figure}[htbp] %  figure placement: here, top, bottom, or page
   \centering 
\includegraphics[width=3.5in]{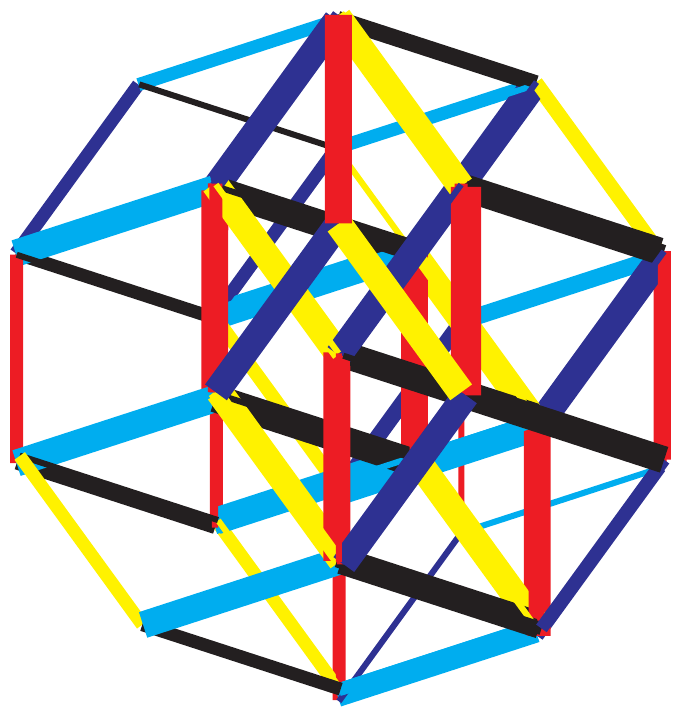} \ = \   \includegraphics[width=3.5in]{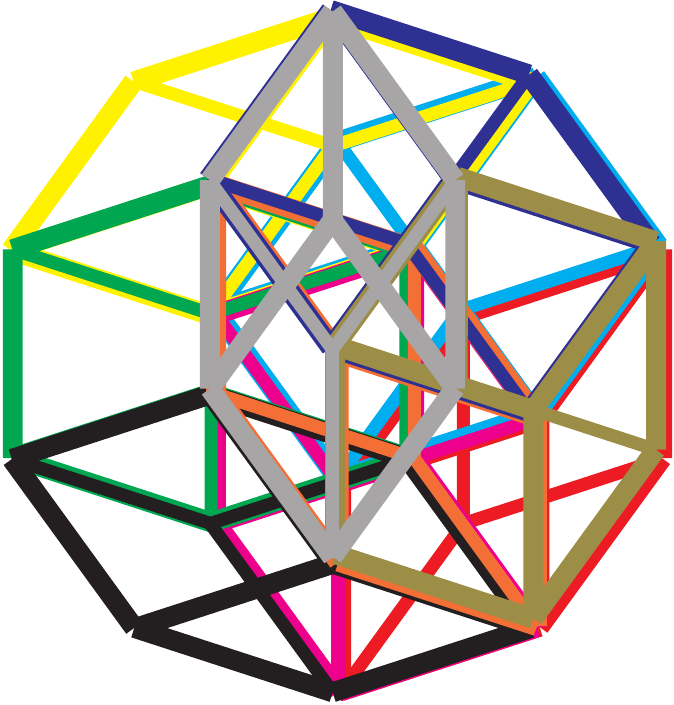} 
  \caption{ 000+0[00--00(0+000[--0000($ \psi_3[5]$)])]}  
   \label{fig:26}
\end{figure}
\begin{figure}[htbp] %  figure placement: here, top, bottom, or page
   \centering
\includegraphics[width=3.5in]{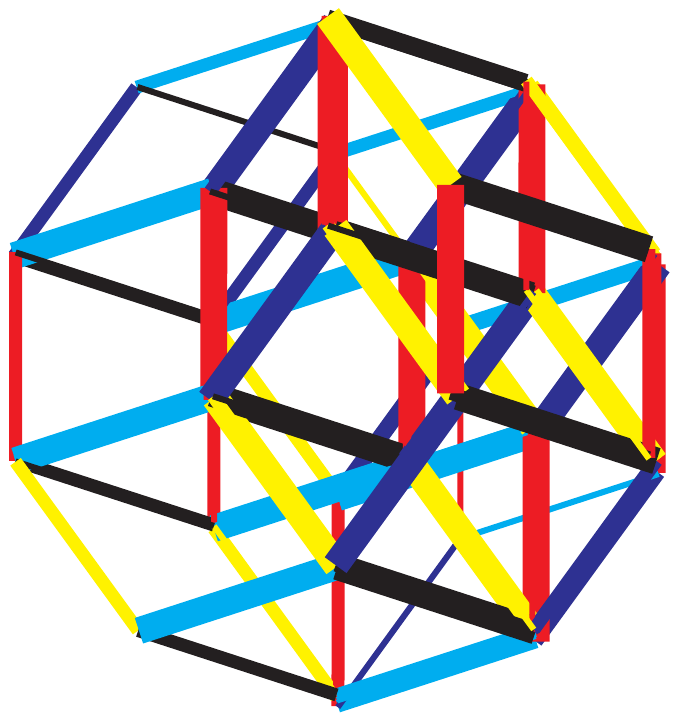} \ = \       \includegraphics[width=3.5in]{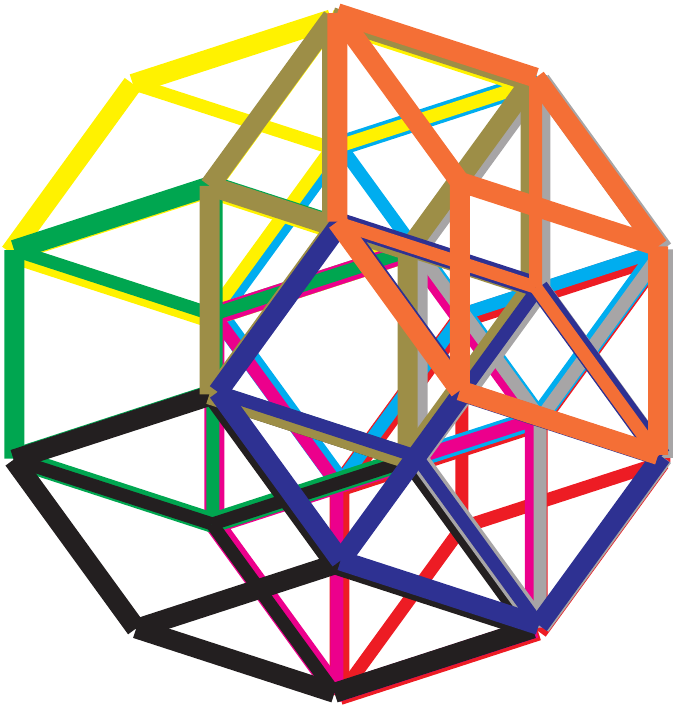}
  \caption{ 0000--(000+0[00--00(0+000[--0000($ \psi_3[5]$)])]) = $\omega_3[5]$.}  
   \label{fig:27}
\end{figure}

The 5-dimensional cube can be treated in the same manner. It is easiest to
present the data as pictures of consecutive $ v_2 $-collections arising
in the sequence of piles. At the top and bottom of each column is $\sigma_2 [5]$
and $ \tau_2 [5]$, which are again mirror images. We obtain 5 columns, corresponding to the five $ {\cal I}^4 $s of the source (Figure 19),
% \ref{fig:19}), 
and another 5 columns 
correponding to the target. Note the bilateral symmetry of the decagonal 
configurations for each of $ {\sigma}_4 [5] $ and $ {\tau}_4 [5]$, and that 
the configurations of $ {\tau}_4 [5] $ are those of the source read backwards 
with $+$ and $-$ interchanged. Each column of decagons provides a path from 
$ \sigma_2 [5] $ to $ \tau_2 [5]$. 
 
Between the decagons of any column we have written the word $x$ whose 
application $ \pi_x $ transforms the decagons one to the next. Hence the 
shaded subsets of Figure 19 %\ref{fig:19} 
indicate $ \sigma_2 (x) $ for each subsequent $x$. We 
merely replace the interior configuration by its mirror image, with the symbols behaving accordingly. 
 
Each decagon itself gives a path from $ \sigma_1 [5] $ to $\tau_1 [5] $. 
Together the five columns of either $ {\sigma}_4 [5]$ or $ {\tau}_4 [5]$ 
give a path from $ {\sigma}_3 [5] $ to $ {\tau}_3 [5]$. 
In each column occurs one $ {\cal I}^4 $ and six $ {\cal I}^3 $s. These 
latter $ {\cal I}^3 $s, together with the four of the source of $ {\cal I}^4 $, make up the ten $ {\cal I}^3 $s of the pile at that stage. The shaded 
$ {\cal I}^2 $s together with those remaining make up the ten of $ v_2 $ at 
that pile stage, accordingly giving the block structure at the 2-dimensional 
level, if we ignore considerations of order. From these figures we can write 
down a sequence of piles with their 1-dimensional sets either present or 
deleted.

In Figure \ref{fig:21}, we indicate how the first four 
columns of $ {\sigma}_4 [5]$
arise from $ {\sigma}_3 [4]$. We use the symbols ${\mu }, \ {\nu}, \ {\lambda}$
and hatching to bring out the sub-cubes of relevance. The dotted squares
correspond to the sixteen $ x_{ij} $ terms from the 4-cube source, depicted 
at the bottom. These exist in the 5-cube source as either old, new or
thickened versions. The lined squares indicate thickenings. Note that one
of the 4-dimensional configurations arises if we collapse all lined squares
of a given decagon along the direction of the lines; each decagon arises by 
splitting an octagon along a path from the 0-source to the 0-target and 
inserting a square for each segment of the path (the lined squares).
Observe that each
column corresponds to doing some of the 4-dimensional operations in the
``old'' ($\lambda $) copy, then moving   each of the four parts of an octagon
across to the ``new''  copy, then completing the octagon operations in the
``new'' (${\nu}$) copy. The subscripts refer to the corresponding objects in
4-dimensions, with 2-dimensional source and target faces indicated
accordingly.

In the sequence of figures (Figure \ref{fig:22} --  \ref{fig:27}) we show how the 3-faces $ {\psi}_3 [5]$ 
are glued together, and how the configuration is transformed into 
$ \omega_3 [5]$ by a sequence of modifications of the interiors of the
3-dimensional sources of the five $ {\cal I}^4 $s making up $ {\psi}_4 [5]$.
These collections of sub-$ {\cal I}^3 $s are delineated in bold.
Each copy of $ \sigma_3 [4] $ is replaced by a copy of $ \tau_3 [4] $, with the same
2-dimensional boundary.
As with the 2-dimensional configurations, the interiors are replaced by
their mirror images, the labels changing accordingly. This enables us to
write down the sequence of pile sets $ v_2 $ and $ v_3 $, but neither
the order nor the intermediate stages of $ v_1 $'s.

%\gapp{37}

\section{The cubical/simplicial dichotomy}

We make some remarks on the relationships between the simplicial 
and cubical approaches. 
We describe how the structure on $ {\cal I}^n $ gives 
rise to structure on 
$ {\Delta}^k $, the $k$-simplex, for $ k \  = \ n+1, \  n , \  n-1 $. 
We shall be
more detailed for the last case, deferring the others to a
subsequent paper as the questions are deeper.

\medskip
\medskip
%\ni
\subsection{$ {\bf   {\Delta}^{n-1} \
as\ a\ sliced\ corner\ of\ 
 {\cal I}^n }$}

The ($n-1$)-simplex has vertices $0,\ 1, \dots,\ (n-1)$. 
Each sub-simplex can
be described as a word in symbols 0 and +, the presence of 0 in the 
$ i^{\rm th} $ position indicating that the vertex $i-1$ is present
in the sub-simplex. Since all such words of length $n$ 
occur in our description of
sub-cubes of the $n$-cube, we have the standard 
description of $ {\Delta}^{n-1}$ as the slice off the corner 
$+\cdots +$ of 
$ {\cal I}^n $ by a hyperplane.

Consider the sequence of piles occurring in the 
cubical cocycle condition.
This gives a sequence of nested embedding of $k$-disks into  
the cube, at each dimension having fixed boundary. Each pile is a
collection of sets of words from ${\cal I}^n$, 
some of which do not contain `$-$'. If we restrict to the subsets
within the pile consisting of those sub-cubes in whose word
description no `$-$ occurs, we see that the action of $ \pi_x $ for
$x$ containing a `$-$' symbol leaves the sub-pile unaltered. Only for
those $x$ consisting entirely of 0's and +'s do we change the
sub-set. Since geometrically every such $x$ intersects 
the slicing hyperplane in one dimension less, 
we see that cubical pile
modifications induce a structure on the simplex also 
describable by piles
and blocks, but where now a word in 0's and +'s corresponds to a
sub-simplex of dimension $|x|-1$. 
In particular, if we look at the source
and target ($n-1$)-face cubes of $ {\cal I}^n $, 
we pick out all
($n-2$)-faces of $ {\Delta}^{ n-1}$. 
These are accordingly divided into two
sets, the source and target faces. Note that the
former correspond to the
consecutive deletion of 
{\it odd}
vertices, whereas to obtain the target faces we delete 
{\it even}
vertices from $ {\Delta}^{n-1}$. 
There is a simple procedure to derive the
induced cocycle condition for $ {\Delta}^{n-1}$ from that of the
cube $ {\cal I}^n $: Rewrite according to the rules
\begin{enumerate}

\item 
Starting with $ *_{n-2}$ and proceeding down to $ * _0 $, consider
the highest dimensional cube between successive 
$ *_k $'s.

\item 
For each such pair, if the highest dimensional term occurring between
them contains a `$-$', delete all terms between the pair,
and one of the $ *_k $ terms.

\item 
Delete all vertices and occurrences of $ *_0 $.

\item 
For the words remaining, write $i-1$ for each occurence of 0 in the 
$ i${th}  position.

\item 
Replace each $ *_k $ by $ *_{k-1}$.
 
\end{enumerate}

\begin{figure}[htbp] %  figure placement: here, top, bottom, or page
   \centering
   \includegraphics[width=6in]{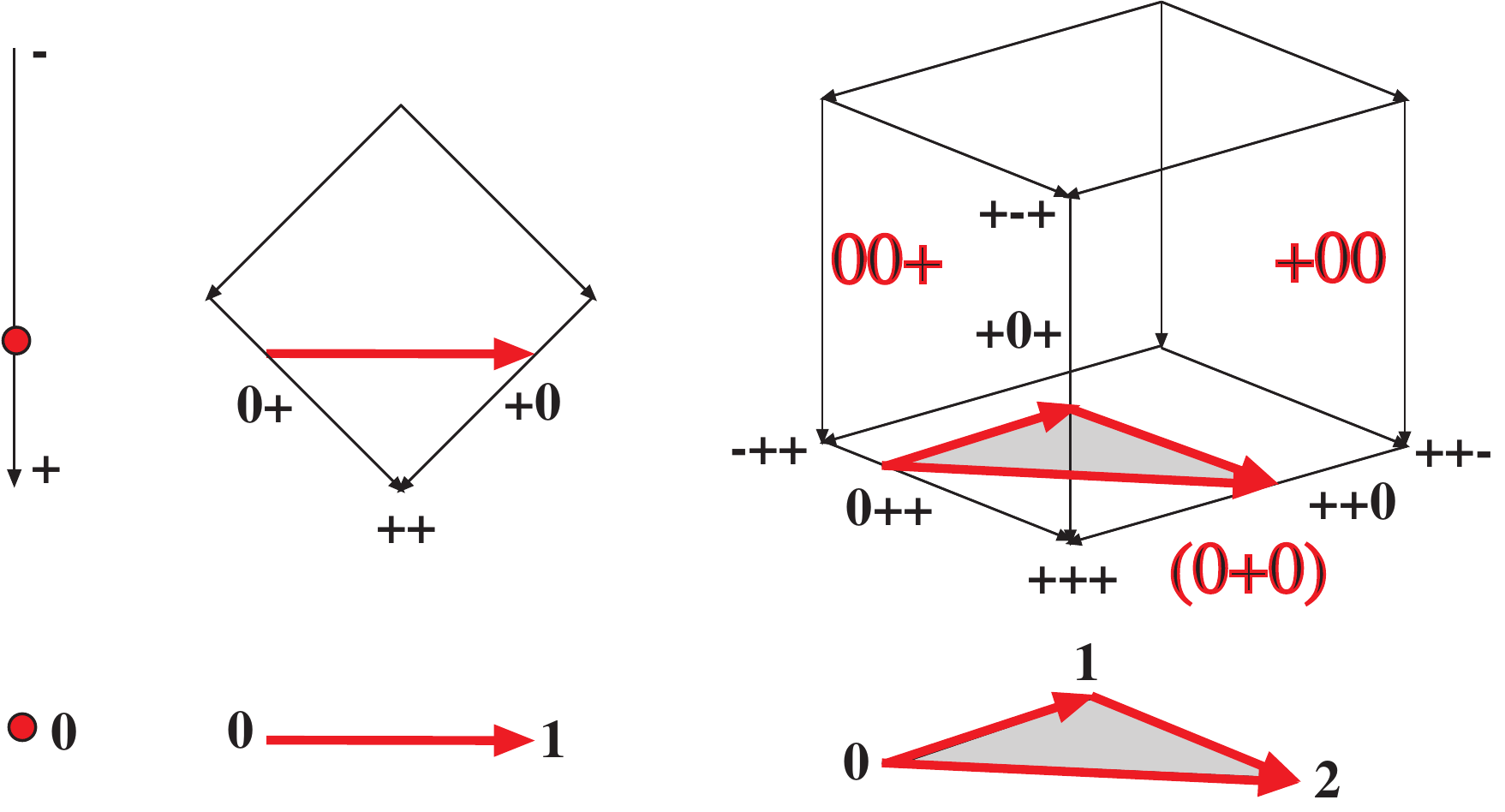} 
   \caption{Slicing off a corner of  cube with labeled sub-cubes gives a simplex with labeled sub-simplices}
   \label{fig:28}
\end{figure}

\begin{figure}[htbp] %  figure placement: here, top, bottom, or page
   \centering
   \includegraphics[width=5in]{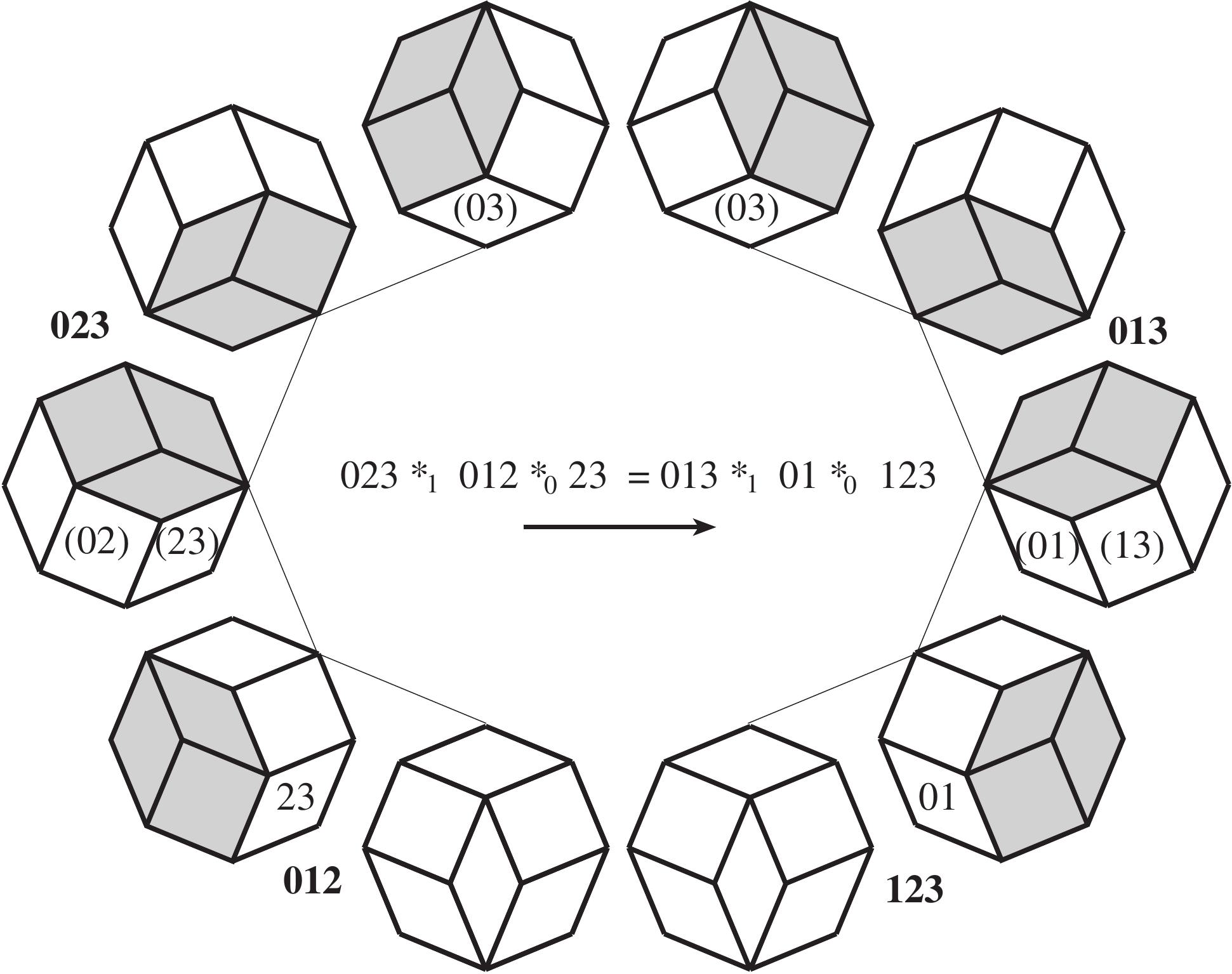} 
%   \caption{Fig29}
%   \label{fig:29}
%\end{figure}

\bigskip

%\begin{figure}[htbp] %  figure placement: here, top, bottom, or page
%   \centering
   \includegraphics[width=5in]{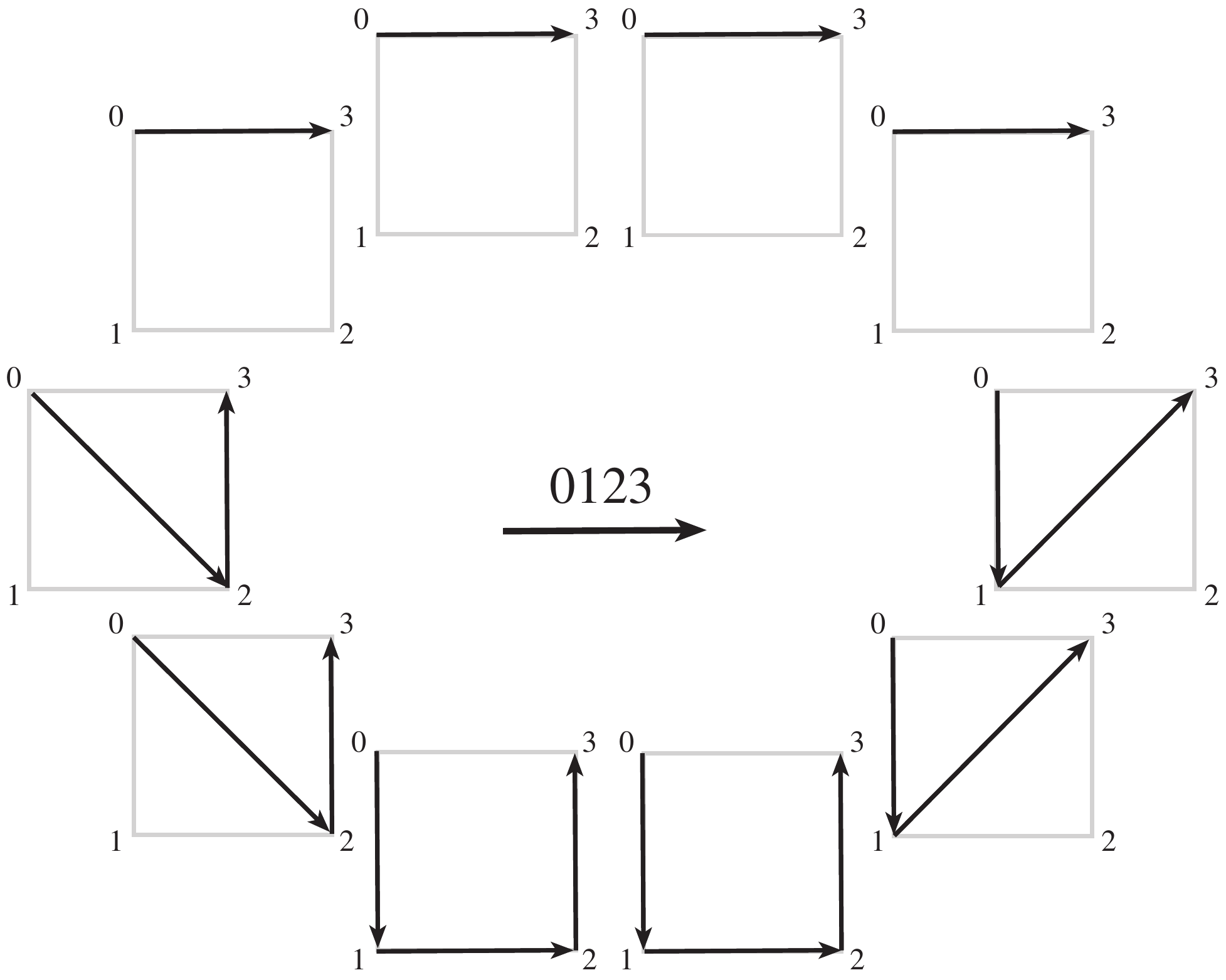} 
   \caption{Deriving simplicial cocycles form   cubical   cocycle}
   \label{fig:29}
\end{figure}

\medskip
\medskip
\noindent
{\bf 
Examples {8.1}:}
For the cocycle conditions previously described, 
we obtain exactly the 
description of simplicial cocycles obtained by Street \cite{St}. 
For the 1-,
2-, 3-cubes, we obtain Figure \ref{fig:28}. In Figure \ref{fig:29} we do the same for
the 4-cube, the upper diagram giving the data in terms of the
`octagon of octagons', the lower the interpretation on the
tetrahedron slice $ {\Delta}^3 $ acting on paths.

\begin{figure}[htbp] %  figure placement: here, top, bottom, or page
   \centering
   \includegraphics[width=6in]{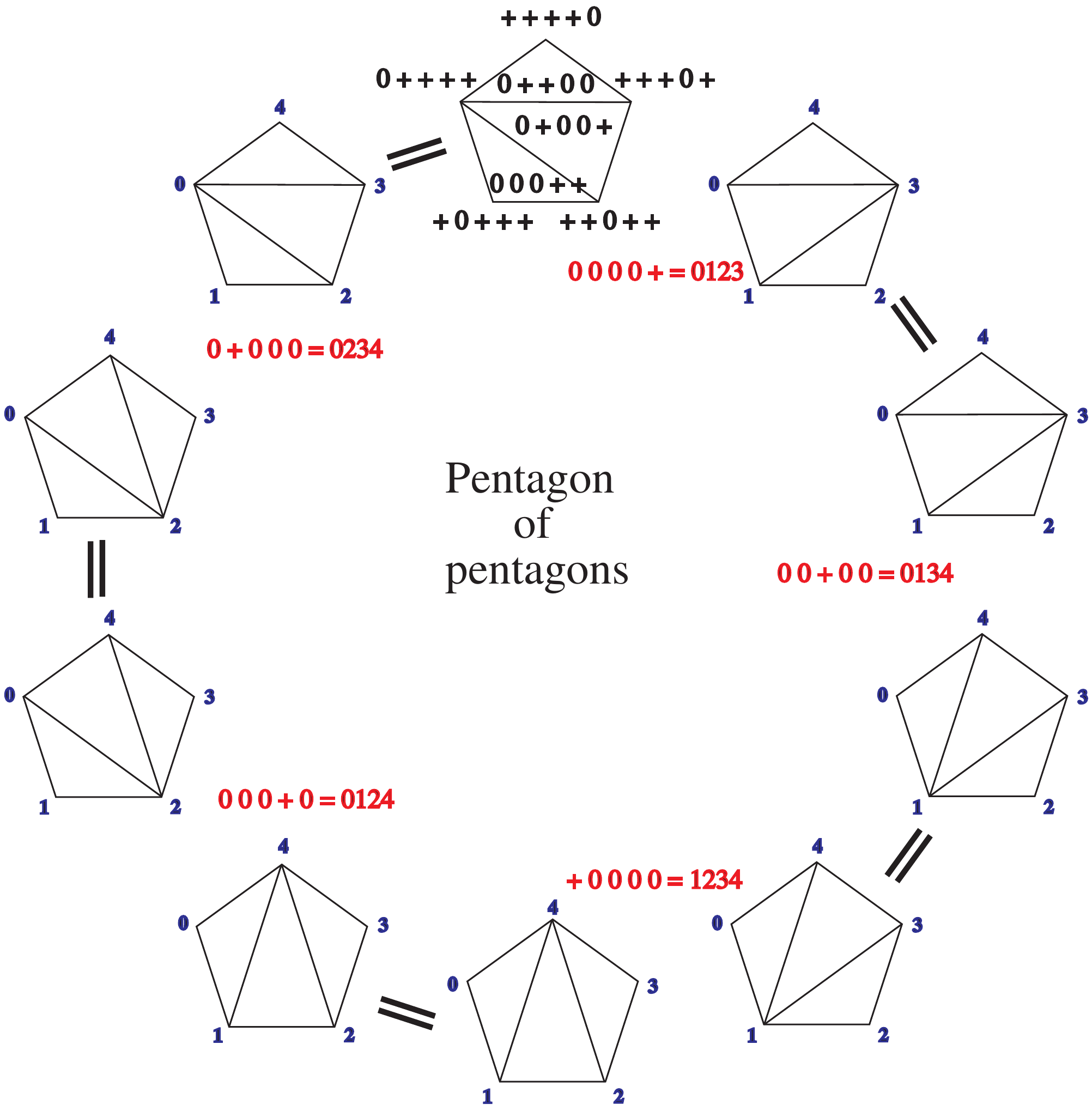} 
   \caption{Slicing the 5-cube, giving the 
   pentagon of pentagons}   \label{fig:30}
\end{figure}

In Figure \ref{fig:30} we show that the 5-dimensional cube gives rise to the
famous `pentagon of pentagons'
arising in homotopy coherence theory as well as in the definition of
bicategories. We refer to \cite{St}  for further discussion. The
pentagon of pentagons can be seen directly from Figures \ref{fig:22} -- \ref{fig:27}, by
slicing the bottom of each figure with a horizontal plane. The
vertices of the pentagon are determined as the intersection between
the plane and the edges
emanating from the 0-target of $ {\cal I}^5 $, 
the internal configuration of
the triangles arising from intersections with $ {\cal I}^2 $s. 
Alternatively we can interpret the decagons as providing 
deformations of
the 0-source of $ {\Delta}^4 $ to the 0-target, 
each column corresponding
to a configuration of triangles in $ {\Delta}^4 $.
We shall shortly obtain two other derivations of the pentagon of 
pentagons.

In \cite{St}   the structure of an {\it 
$n$-category}
is shown to exist freely on the $n$-simplex. 
This derives from Street's
notion of sources and targets in terms of even and odd faces. 
From his results on 
`excision of extremals', it is possible to derive the cocycle
conditions for this structure, non-canonically since order is not
determined. Alternatively, the
$n$-category structure can  be defined 
if we have the cocycle conditions. 
Since at the level of ($n-1$)-faces of 
the $n$-simplex we have derived the same structure as in [St], 
we conclude 

\medskip
\medskip
\noindent
{\bf 
Theorem 8.2 (\cite{St}):}
{\quad {\bf (a)}  } $ {\Delta}^n $ admits the structure of
an $n$-category.
\begin{description}
 
\item[(b)] Street's structure of orientals can be derived 
from the 
structure of oriented cubes, by slicing.
\item[(c)]  There is a canonical form inductively defined
for the  expression of simplicial cocycle conditions. 
\item[(d)]  The simplicial cocycle conditions can be written
in a  form involving neither brackets nor the 2-categorical middle
interchange law.
\end{description}

\medskip
\medskip

Roberts had conjectured
that there should be a natural way of expressing the cocycle
conditions without brackets.

\medskip
\medskip

The two approaches that follow, 
relating cubes and simplices, also lead to
induced coycle conditions, 
indicating just how `rigid' and canonical is the
structure.

\subsection{  { \bf Stretching\ a\ simplex\ into\ a\ cube}}

We show how the cocycle condition on $ {\cal I}^n $ 
gives rise to one on $ {\Delta}^n $. 
To do this we briefly describe how to `stretch' the
$n$-simplex into a cube.
\begin{enumerate}
\item For dimension one, $ {\Delta}^1 $ and $
{\cal I}^1 $ coincide. 
\item  Higher dimensional simplices arise by
`coning', whereas cubes arise by taking the product with $I$, the
unit interval. We use labellings of vertices of simplices by the
symbols $0,\ 1,\dots ,\ n$. For the 1-simplex, vertices $0$ and $1$
(Figure \ref{fig:31}), the 1-cell is labelled $01$.

\begin{figure}[htbp] %  figure placement: here, top, bottom, or page
   \centering
   \includegraphics[width=1.5in]{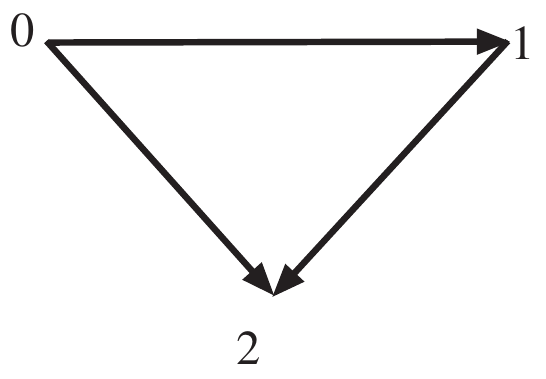} \qquad\qquad
      \includegraphics[width=1in]{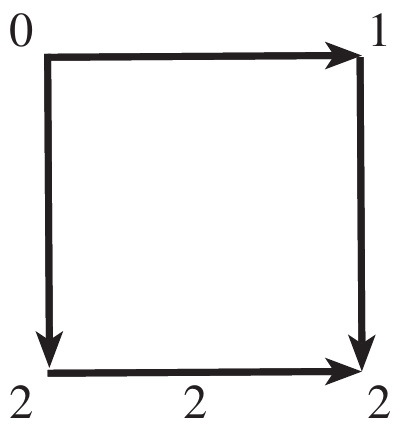} 
        \caption{Expanding a labeled 2-simplex into a labeled 2-cube}
   \label{fig:31}
\end{figure}
\begin{figure}[htbp] %  figure placement: here, top, bottom, or page
   \centering
      \includegraphics[width=2in]{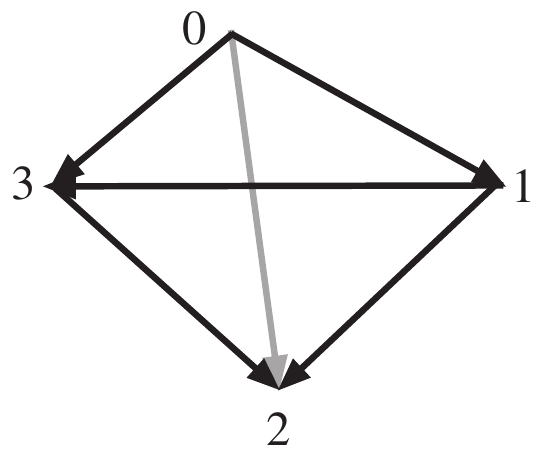}  \qquad\qquad
         \includegraphics[width=1.5in]{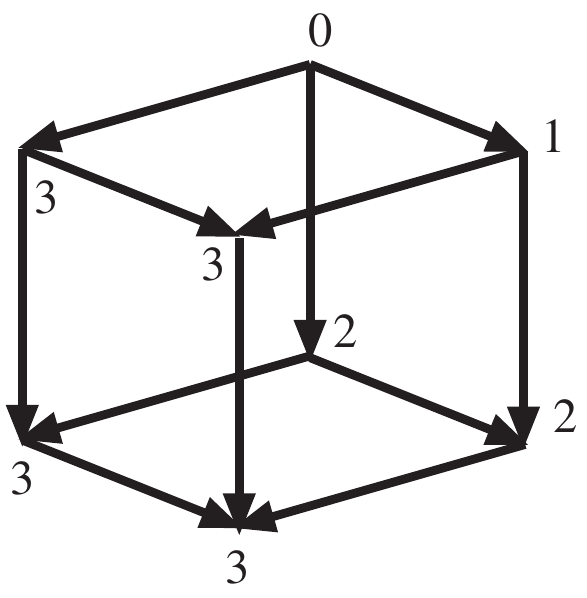} 
   \caption{Expanding a labeled 3-simplex into a labeled 3-cube}
   \label{fig:32}
\end{figure}

Take the 2-cube whose other
two vertices are labelled `2'. The edges are now $01$, $02$, $12$, and
`$2$'. The latter edge is really a stretched out vertex,
corresponding to an `identity arrow' in the categorical setting. 
Note that if we identify everything labelled `2' together, 
we obtain a 
triangle (2-simplex).
\item  Taking the product with the interval, 
and labelling all new vertices
by `3' leads to the 3-cube labelled as in Figure \ref{fig:32}. Again,
collapsing together like-labelled objects yields the tetrahedron.
The source and target configurations for the 3-cube yield Figure \ref{fig:33},
in which each 2-cell is labelled by its vertices. We obtain the
simplicial 2-cocycle  condition from the cubical one by direct
substitution, deleting  degenerate ($\equiv$ identity) cells. Note
that the only change from the source to the  target configuration
is the internal shape and re-labelling of the vertex. 

\begin{figure}[htbp] %  figure placement: here, top, bottom, or page
   \centering
$   \begin{array}{ccccc}
   \includegraphics[width=1.5in]{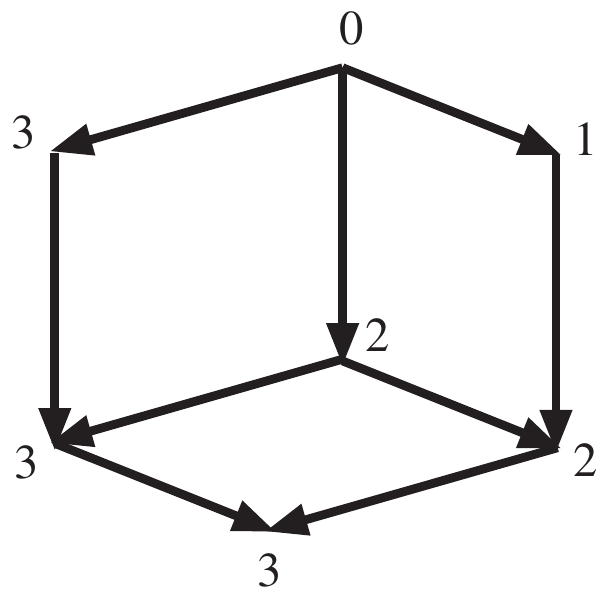} 
     \begin{array}{c}
     \longrightarrow 
      \\
     \\
     \\
          \\
     \\
     \\
          \\
     \end{array}
   &
      \includegraphics[width=1.5in]{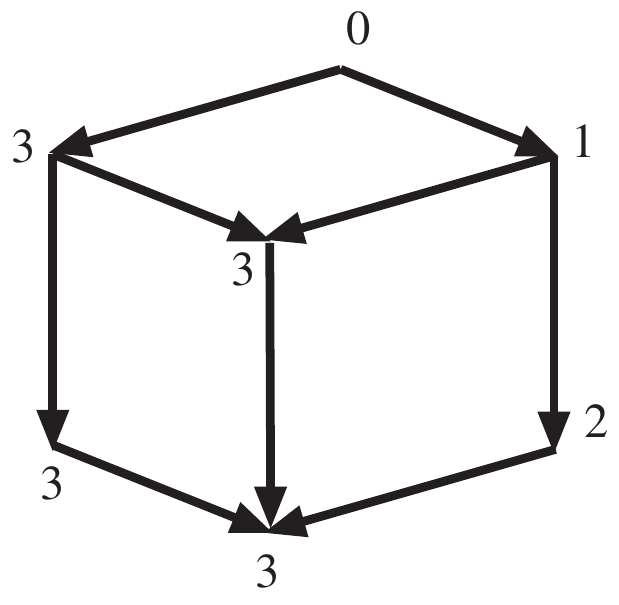} &\qquad
         \includegraphics[width=1in]{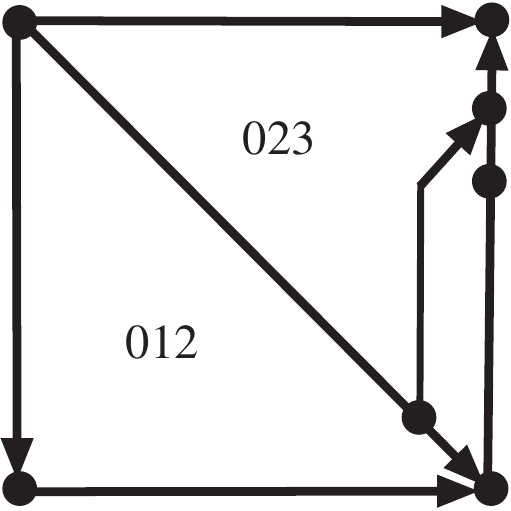} &
   \begin{array}{c}
       \\
       \\
     \longrightarrow  \\
       \includegraphics[width=0.4in]{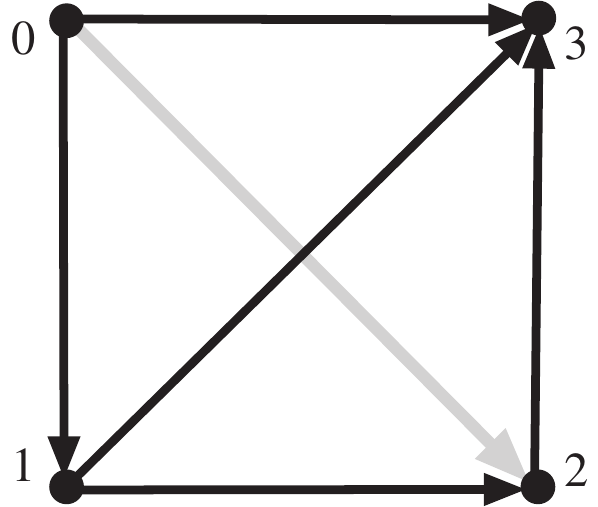} \\
    0123  \\
     \\
     \\
     \\
     \\
     \\
     \\
     \end{array}
   &
      \includegraphics[width=1in]{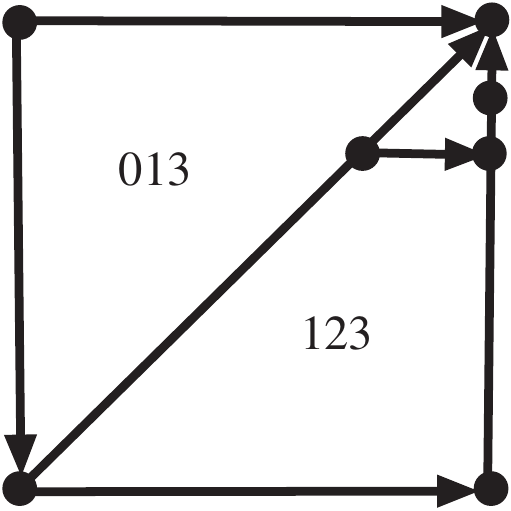} \\
   \end{array}$
   \vskip -3cm
   \caption{The 3-cubical 2-cocycle becomes the 3-simplicial 2-cocycle}
   \label{fig:33}
\end{figure}
%\vskip -1cm

\begin{figure}[htbp] %  figure placement: here, top, bottom, or page
   \centering
   \includegraphics[width=4.8in]{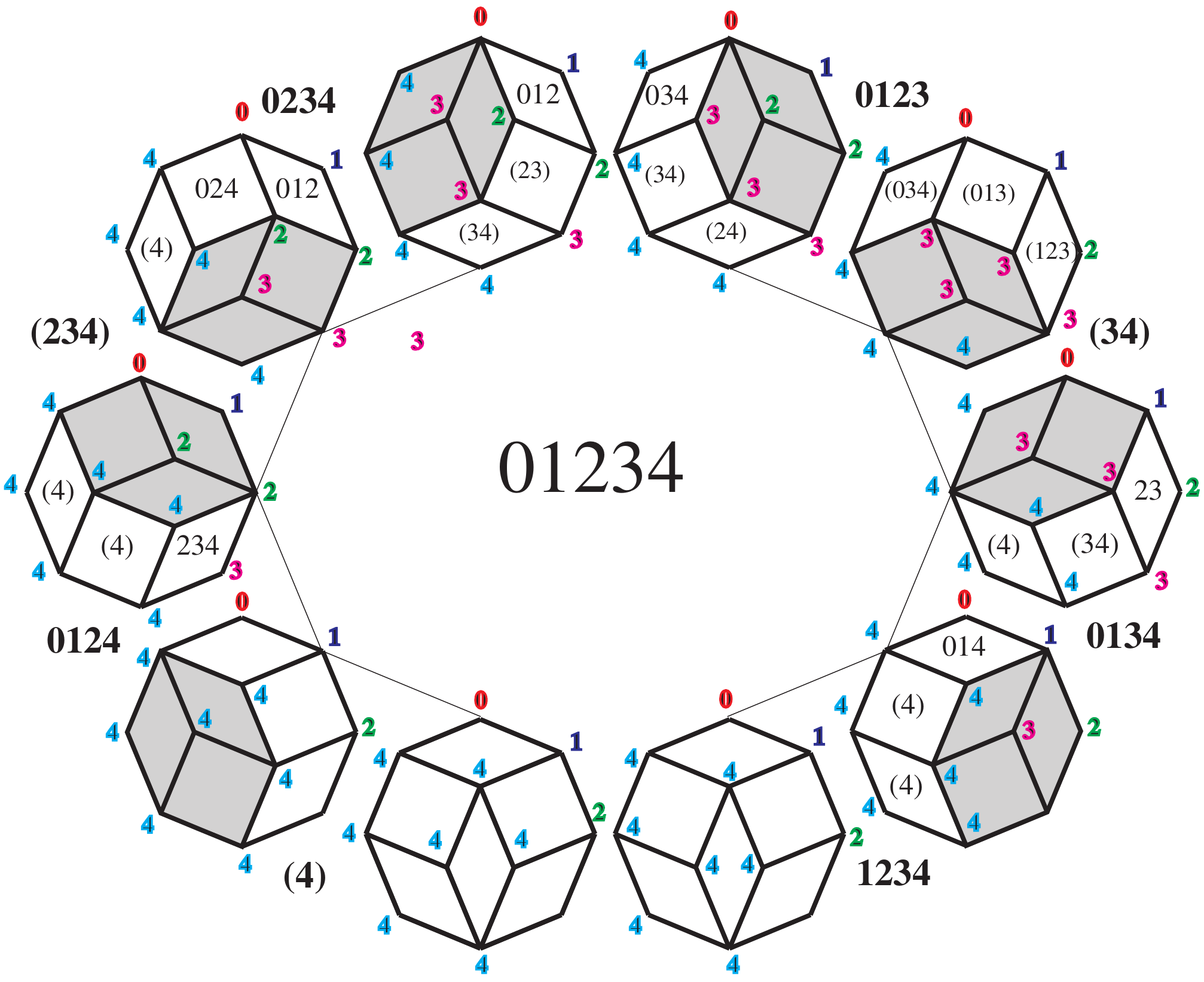} 
   \includegraphics[width=5in]{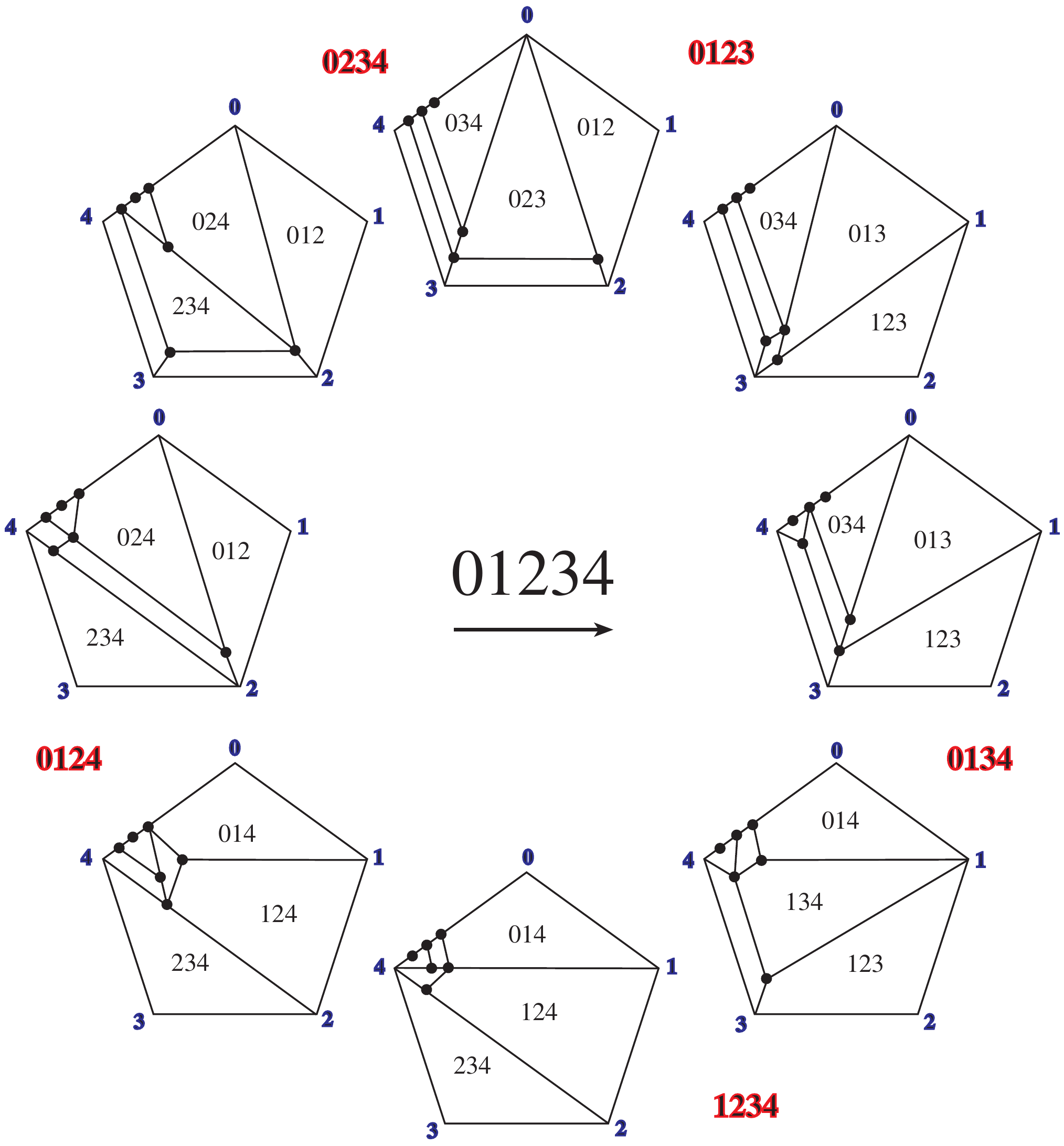} 
   \caption{Deriving the pentagon of pentagons from the octagon of octagons}
   \label{fig:34}
\end{figure}

\item Proceed in the same way to obtain $
{\Delta}^n $ `exploded' into  $ {\cal I}^n $. 
\item The rule to reverse the process is to
replace a vertex of the cube, as a  word in  $-$ and $+$, by the
position of the highest occurring $+$, or $0$ if none occur. For an
arbitrary word in $-$, $0$ or $+$, we also record all positions of $0$
beyond the highest $+$. 
This enables us to relabel each sub-cube as an exploded sub-simplex.
\item  Now rewrite cubical cocycles by replacing
words in this fashion, forgetting those (as in $1$ of the slicing
approach) correponding to `identities'. These identities are
those simplices of dimension less than the cube from which they
arise in this correspondence. 
\end{enumerate}

\medskip
\medskip
In Figure \ref{fig:34}, we show how the octagon of octagons 
gives rise directly
to the pentagon of pentagons. The upper figure gives direct
substitution, the lower a deformation of the figures into a more
recognisable shape.

The cubes arising in this fashion have a number of faces which are
`identities'. The same happens in the following 
procedure:

\subsection{\bf Strings\ on\ simplices}

Any (ordered) path from $0$ to $n$ on $ {\Delta}^n $ is 
uniquely determined
by listing the vertices through which it passes. 
Since the list for any
such path contains $0$ and $n$, 
we need not include these in a specification
of the path. If we insist on listing vertices in increasing order of
magnitude, and $i$ and $j$ are two consecutive integers 
in the list with
$j\not= i+1 $, then for each triangle $ikj$ with $ i < k <j $ 
there is a triangle in the $n$-simplex and a new path
obtained by including $k$. In this way triangles act on paths, and
thus label morphisms in the (discrete) space of paths.

The resulting complex, with paths as objects and
triangles as morphisms 
(acting as deformations of one path to another), is
well known to form $ {\cal I}^{n-1} $. 
We refer the reader to \cite{Ba} or \cite{Le}.
Not all faces of the cube correspond to faces of the simplex, 
some as in
the previous procedure corresponding to `identities'. 
There is fortunately
a well defined procedure to see which sub-cubes correspond to
sub-simplices. Some repetitions occur, as 
can be understood by considering
the effect of a triangle on a path: 
different paths may abut the same
triangle, giving rise to distinct vertices of the cube, 
but with edges emanating having the same triangular correspondence.
The application of this approach to the understanding of simplicial
oriental cocycles was initiated by
Duskin \cite{Du}. % [Du].

\medskip

Here is the prescription for converting the cubical 
($n-1$)-cocycle condition into the simplicial ($n-1$)-cocycle condition:
\begin{enumerate}
\item  Consider terms of dimension $n-1$. 
A unique
such occurs between consecutive $ *_{n-2} $'s.  
\item If any $0$'s in the word are separated by a
$+$, delete all terms between the prior and subsequent $ *_{n-2} $.
Delete one of the $ *_{n-2}$'s. 
\item  If no $+$ occurs before a $0$, write $0$.
\item  If a $+$ occurs before a $0$, record the
highest position in which a $+$ occurs before a $0$ is encountered.
\item  Record all positions in which a $0$ occurs.
\item  If no $+$ is then reached before the end of
the word, write $n+1$. 
\item  If a $+$ does occur to the right of all $0$'s,
record the position of the first such encountered. 
\item  Now look at $(n-2)$-dimensional words between
the $ *_{n-3} $, and follow the analogues of steps (2) and (3).
\item Repeat this process for every dimensional
word. 
\item  Replace each $ *_i $ by $ *_{ i+1} $.
\item  Insert the terms 
$\ (0.1) *_0 (1.2) *_0 \dots *_0 (i-1.i) *_0 $
before the word, if the first number of a word is $i$. If the last
letter is $j$, add after the word the symbols $ *_0 (j.j+1) *_0
\dots *_0 (n.n+1).$
\end{enumerate}

\begin{figure}[htbp] %  figure placement: here, top, bottom, or page
   \centering
 $
   \begin{array}{ccc}
   \includegraphics[width=3in]{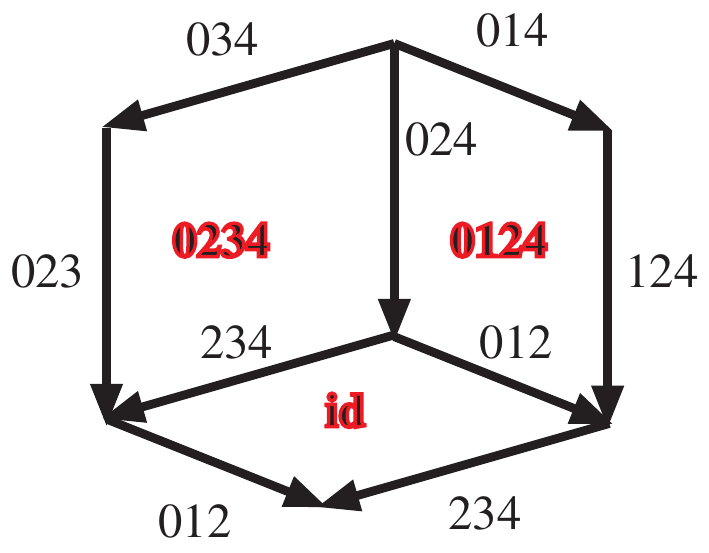} &
   \begin{array}{c}
    \longrightarrow\\
   \\
   \\
   \\
   \\
   \\
   \\
   \\
   \\
   \end{array}
      \includegraphics[width=3in]{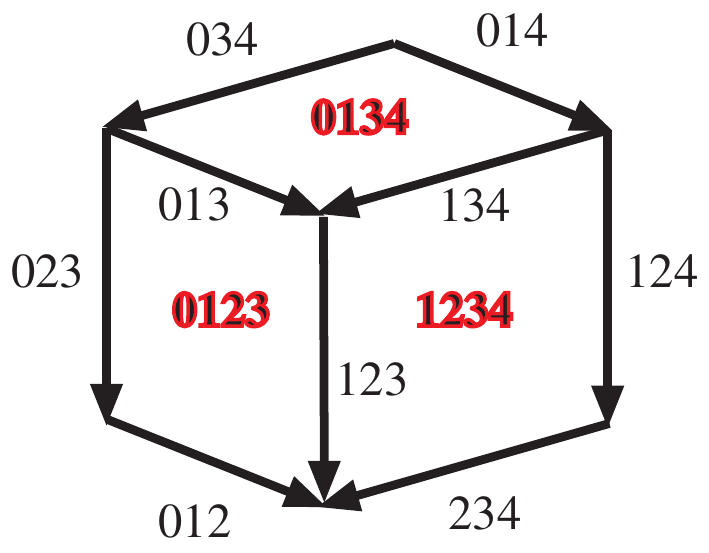} &
   \end{array}
      $
      \includegraphics[width=6in]{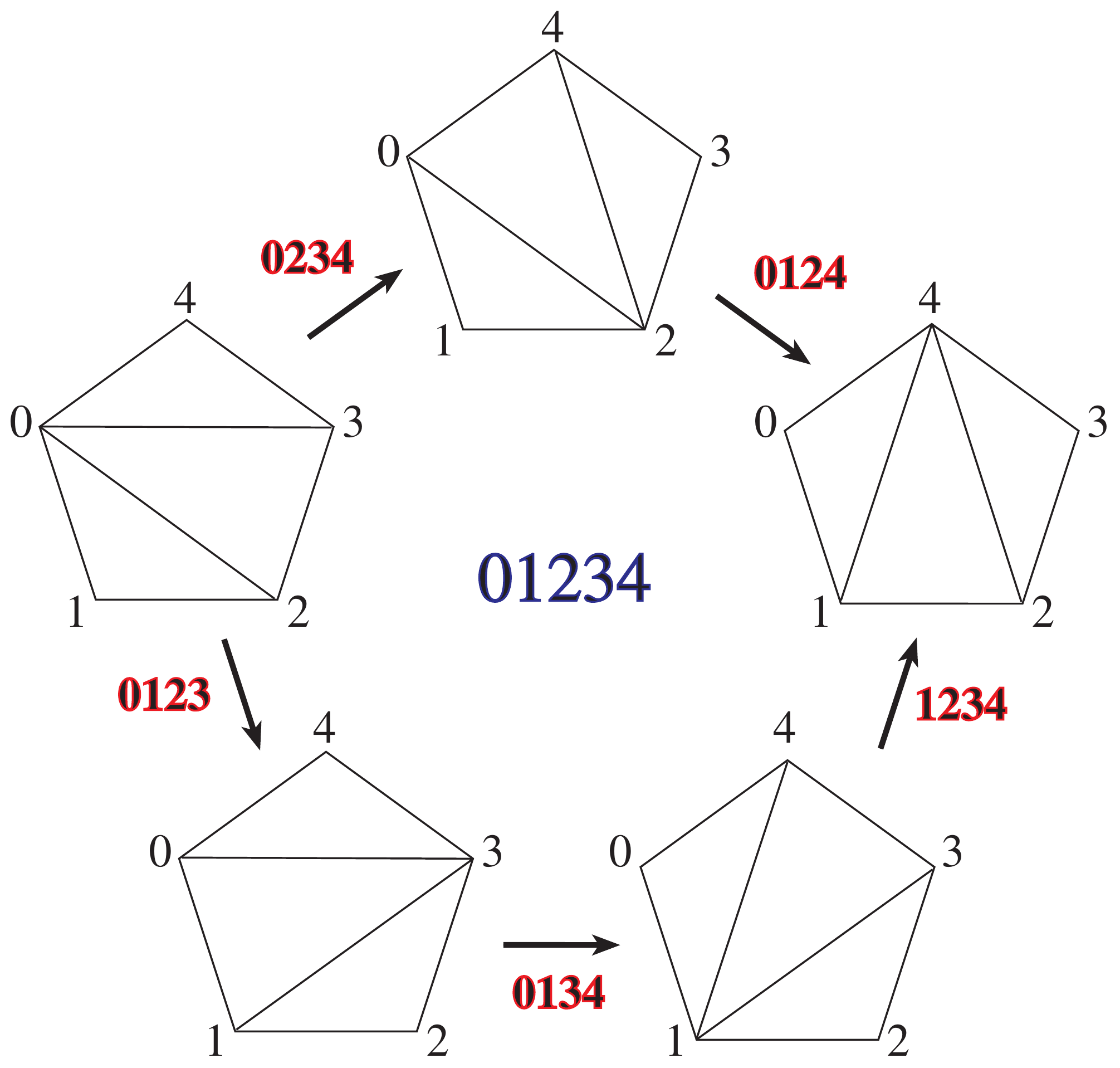} 
   \caption{The pentagon of pentagons as the 3-cube  `Yang-Baxter' equation}
   \label{fig:35}
\end{figure}

\begin{figure}[htbp] %  figure placement: here, top, bottom, or page
   \centering
      \includegraphics[width=6in]{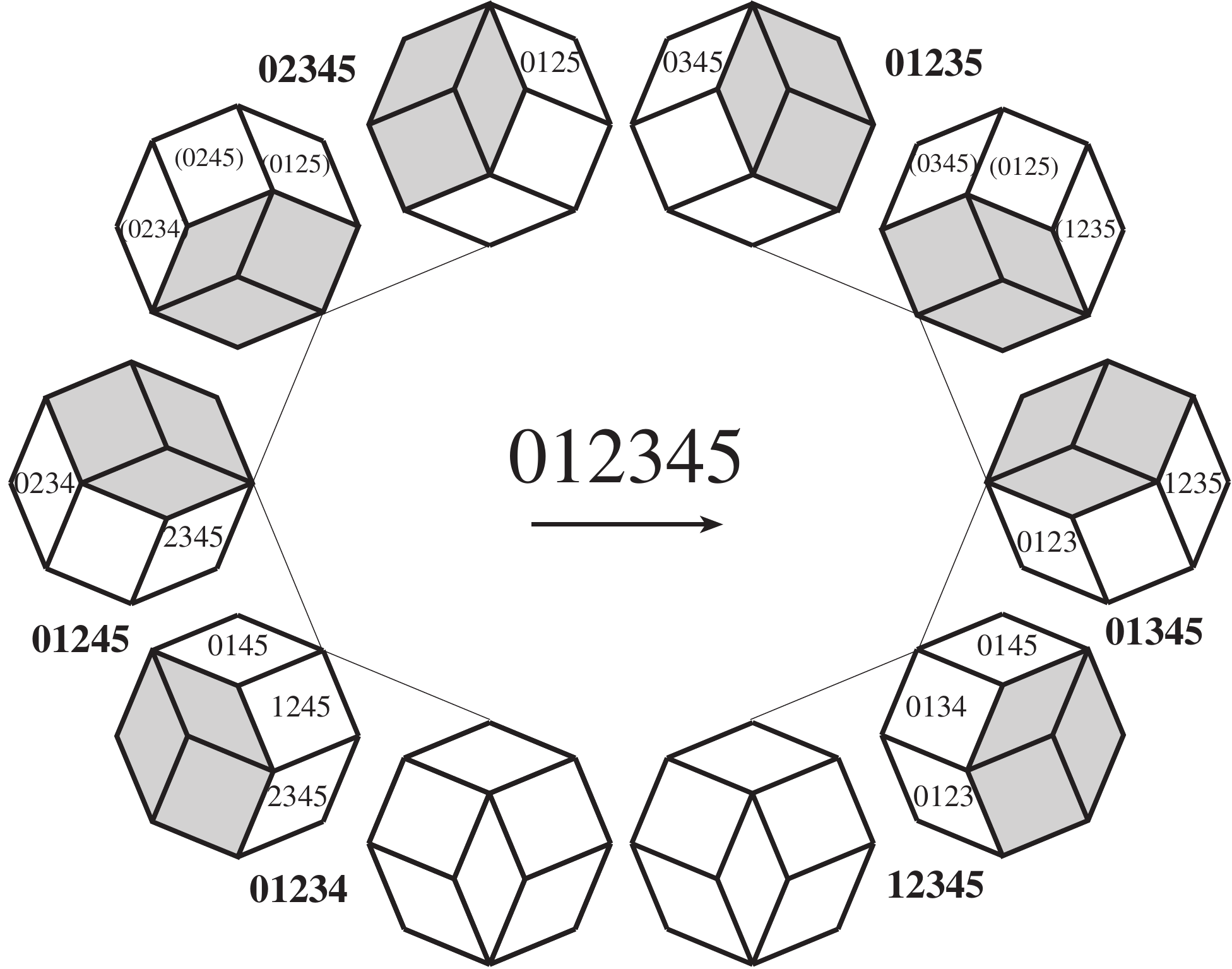} 
   \caption{The 4-cube  octagon of octagons giving rise to the 5-simplicial 4-cocyle }
   \label{fig:36}
\end{figure}

\medskip
\medskip
\noindent
{\bf
Examples {8.3}:}
In Figures \ref{fig:35} and \ref{fig:36} we show again how the usual simplicial
cocycle conditions arise. From the
3-cube we obtain the pentagon of pentagons, 
from which it is also possible
to see that the `identity' face of the 3-cube correponds to a 
{\it
re-ordering
}
of the triangles across which we are deforming the 1-source of 
$ {\Delta}^4 $ 
onto the 1-target. Hence the cube contains some additional information
about structure manifest at the simplicial level. 

 It is possible to iterate this procedure, and consider strings on
the cube. Doing so leads to the {\it parallelahedra}, which in
dimensions 2- and 3- are the hexagon and truncated octahedron. The
parallelahedra are known to tesselate Euclidean space, and like the
cube and simplex have an intimate relationship with the symmetric
groups. We refer the reader to Baues \cite{Ba}  for applications to the
geometric cobar construction. Another discussion occurs in Leitch \cite{Le}.
 We also remark that to carry over  to the parallelahedra the
$n$-category structures discussed above should be possible in a
geometric fashion, there being an inductive definition for the
construction of parallelahedra. This is more difficult to describe
directly as the parallelahedra are no longer regular figures. A
categorical approach should exist. 

These descriptions for obtaining simplicial cocycles 
from cubical ones
show that the cubical and simplicial cocycles/structures are deeply
interlocked. The procedures `take strings', `slice' and `explode'
provide a sequence of correspondences of structures: three modifications from each of  the $n$-cubical cocycles, and three ways to obtain each of the $(n-1)$-simplicial  cocycles:
$$\matrix{
\dots &\longrightarrow  & {\cal I}^{n+1} & \longrightarrow
& {\cal I}^{n} & \longrightarrow & {\cal I}^{n-1} &
\longrightarrow
&  {\cal I}^{n-2} &  \longrightarrow  & \dots \cr
&&&&&&&&&&\cr
  &  &  &  \buildrel {string} \over \swarrow & ex\downarrow  &
 \buildrel {slice} \over \searrow & ex \downarrow 
& \buildrel {string} \over \swarrow &  &  &  \cr
&&&&&&&&&&\cr
\dots &\longrightarrow &  {\Delta}^{n+1} & \longrightarrow
&  {\Delta}^n &  \longrightarrow &  {\Delta}^{n-1} & 
\longrightarrow  &  {\Delta}^{n-2} & \longrightarrow  & \dots \cr}
$$

There is an analogue of the `$k$-blocks' for cubes in the simplicial
context, and a manifestation of the symmetries of the cubical 
cocycles.
On the other hand, the three relationships described 
above show how the
structure on the cubes gives rise to a structure
on simplices, 
but it is not clear how the procedure can be reversed for 
{\it any
}
of the approaches! 

It would be of interest to investigate 
further the commutativity of the above sequences, 
with respect to face and degeneracy maps for cubes or
simplices. In fact, to put the whole structure in a more explicitly
categorical framework would be desirable. 
The symmetries of many of the
objects also suggests possible extensions of Connes' cyclic
homology, for which we refer to Goodwillie \cite{Go}  for a categorical
approach.

\bigskip

\bigskip

\noindent
{\bf Added 2009: 3-cell recombination for the 5-cube [5], figures after references}

The 3-source of the 5-cube [5] contains $10 =
\left(\begin{array}{c}5\\3\\ \end{array}\right)$ 3-cubes. Each cube is a parallel class determined by the location of three 0's in a string of length 5, there being $2^2=4$ possibilities in total corresponding to the remaining $\pm $ entries. An application of a 4-cube in the 5-cube 4-source may change the representative: each 3-cube class is involved twice in changing the two $\pm$ entries to $\mp$ corresponding to the 3-cubes in the 3-target, which are the first and last colummns. We only record the position of 0's in intermediate columns. 

Each column bold indicates source 3-cubes of the operational 4-cube:
Observe that in each case reordering can be achieved before application, so that these source 3-cubes occur consecutively. Partial order rearrangement is possible when two 3-cubes do not share a common face (pair of numbers): two cubes intersect on at most 1 vertex. (This is clearer in the dual Pascal's triangle version). The columns before the double vertical lines, read from bottom to top, enable the filling of the 3-ball by 10 distinct 3-cubes.

 {\footnotesize
 $$
 \begin{array}{cc|c||c|c||c|c||c|c||c|c||c|cc}
  \begin{array}{c}
--000\\  
-0+00\\ 
0++00\\
-00-0\\ 
0+0-0 \\ 
00--0 \\ 
-000+\\ 
0+00+\\ 
00-0+\\ 
000++\\ 
\end{array}
&
 \begin{array}{c}
   {\bf  345}\\
  {\bf   245}\\
 145\\
  {\bf  235}\\
135\\
125\\
  {\bf   234}\\
134\\
124\\
  123
\end{array}
 &
  \begin{array}{c}
   {\bf  345}\\
  {\bf   245}\\
  {\bf  235}\\
    {\bf   234}\\
   145\\
135\\
125\\
134\\
124\\
  123
\end{array}
 &
 \begin{array}{c}
 234\\
 235\\
 245\\
  {\bf   345}\\
   {\bf   145}\\
    {\bf   135}\\
  125\\
    {\bf   134}\\
  124\\
 123
\end{array}
&
 \begin{array}{c}
 234\\
235\\
 245\\
 {\bf   345}\\
   {\bf   145}\\
    {\bf   135}\\
    {\bf   134}\\
      125\\
  124\\
 123
\end{array}
&
 \begin{array}{c}
  234\\
  235\\
   {\bf   245}\\
 * { 134}\\
*  {135}\\
*  {\bf 145}\\
 * {345}\\
    {\bf   125}\\
   {\bf    124}\\
  123
\end{array}
&
 \begin{array}{c}
  234\\
  235\\
{ 134}\\
 {135}\\
   {\bf   245}\\
 {\bf 145}\\
    {\bf   125}\\
   {\bf    124}\\
 {345}\\
  123
\end{array}
&
 \begin{array}{c}
  234\\
 {\bf  235}\\
 {134}\\
{\bf  135}\\
 * { 124}\\
 * {\bf 125}\\
 * { 145}\\
 *   { 245}\\
  {345}\\
{\bf   123}
\end{array}
&
 \begin{array}{c}
  234\\
 {134}\\
 { 124}\\
  {\bf  235}\\
{\bf  135}\\
 {\bf 125}\\
 {\bf   123}\\
{ 145}\\
 { 245}\\
  {345}\\
\end{array}
&
 \begin{array}{c}
 {\bf   234}\\
{\bf  134}\\
{\bf   124}\\
 *  {\bf    123}\\
 *   { 125}\\
*{  135}\\
*  {235}\\
{ 145}\\
 { 245}\\
  {345}\\
\end{array}
&
 \begin{array}{c}
 {\bf   234}\\
{\bf  134}\\
{\bf   124}\\
 {\bf    123}\\
 { 125}\\
{  135}\\
  {235}\\
{ 145}\\
 { 245}\\
  {345}\\
\end{array}
&
 \begin{array}{c}
 {   123}\\
*{ 124}\\
*{ 134}\\
* {  234}\\
 { 125}\\
{  135}\\
  {235}\\
{ 145}\\
 { 245}\\
  {345}\\
\end{array}
&
 \begin{array}{c}
 {   123}\\
{ 124}\\
{ 134}\\
 {  234}\\
 { 125}\\
{  135}\\
  {235}\\
{ 145}\\
 { 245}\\
  {345}\\
\end{array}
&
  \begin{array}{c}
  000--\\ 
00+0-\\ 
0-00-\\ 
+000-\\ 
00++0 \\ 
0-0+0 \\ 
+00+0\\ 
0--00\\
+0-00\\ 
++000\\  \end{array}
 \end{array}
 $$
}

\bigskip

The following figures, from the original, indicate how the simplicial cocycle ${\cal O}_6$ from the 6-simplex derives from the  
5-cube: three nontrivial columns occur in the 5-cube source, and four in the target.  Observe from the first 5-cube source column the four- and five-simplices, or considering only source four-simplices:
{\footnotesize
$$
<023456>*
<02356>* 
<01245> *  
<01234> 
 \Longrightarrow  
<03456> * 
<02356> *
<02345> *
<01256> *
<01245> * 
<01234>
$$
}
which are the labeled 2-cells in $F_{16} $ of Street \cite{St}. Adding in the
cubically-ordered listed three-simplices gives 
{\footnotesize
$$\begin{array}{l}
<023456>*<0126> \\
\qquad **
{<0234>*}<0245>*<2345>*<01256>*<2356>*<3456> \\
\qquad **
{<0234>*}<01245>*<2345>*<0156>*<1256>*<2356>*<3456> \\
\qquad ** 
<01234> *<0145>*<1245>*<2345>*<0156>*<1256>*<2356>*<3456>
\end{array}
$$
}
{\footnotesize
$$\begin{array}{l}
= \  (<03456>*<0236> 
** 
<0345>
*
<02356>*<3456> 
** <02345>
*
<0256>*<2356>*<3456> 
)*<0126>  \\
\qquad ** 
<0245>*<2345>*<02356>*<01256>*<2356>*<3456> \\
\qquad ** 
<0234>*<01245>*<2345>*<0156>*<1256>*<2356>*<3456>\\
\qquad **
<01234> *<0145>*<1245>*<2345>*<0156>*<1256>*<2356>*<3456>
\end{array}
$$
}
  \[
%\left(
\begin{array}{ccccc}  
      \includegraphics[width=1in]{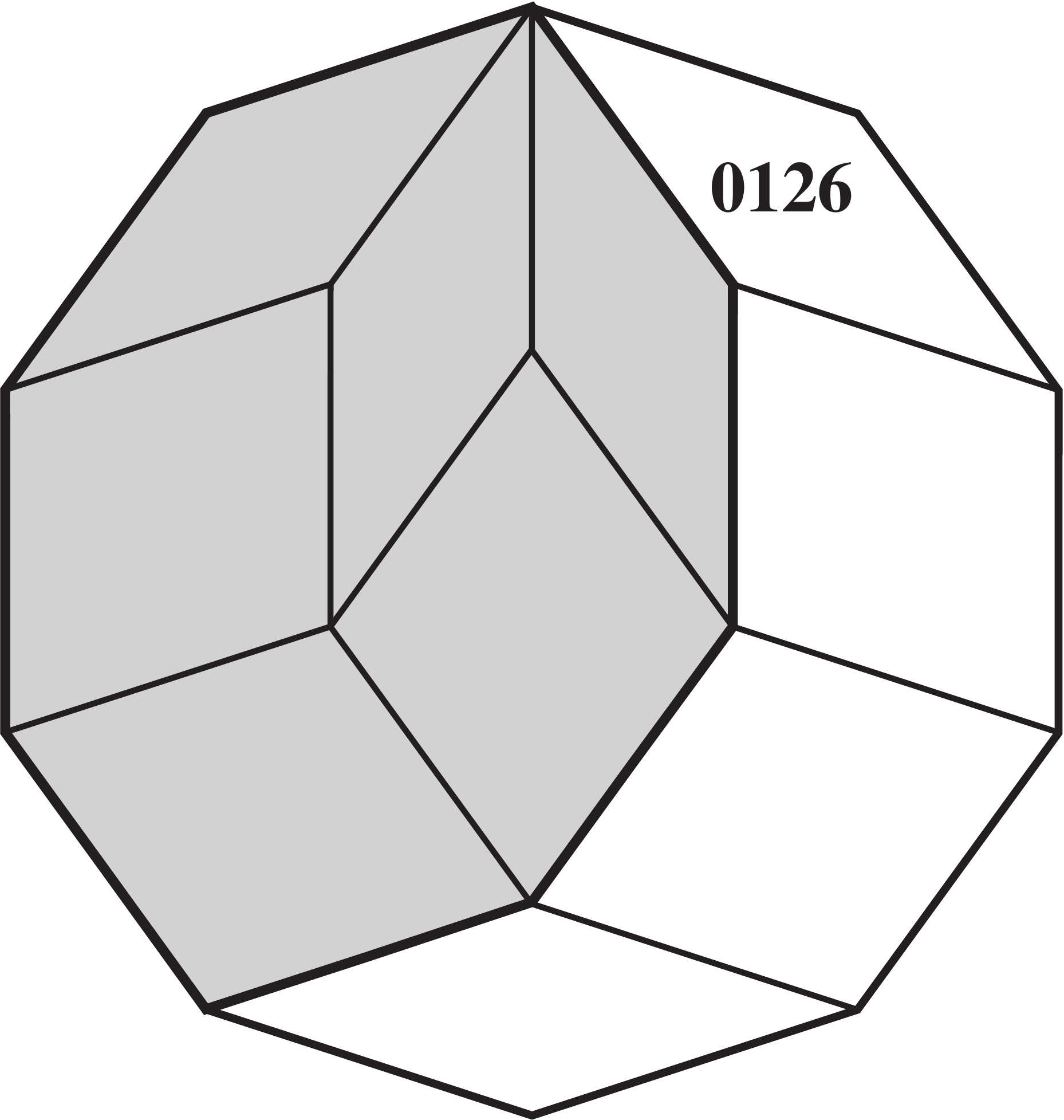}     
  &      \includegraphics[width=1in]{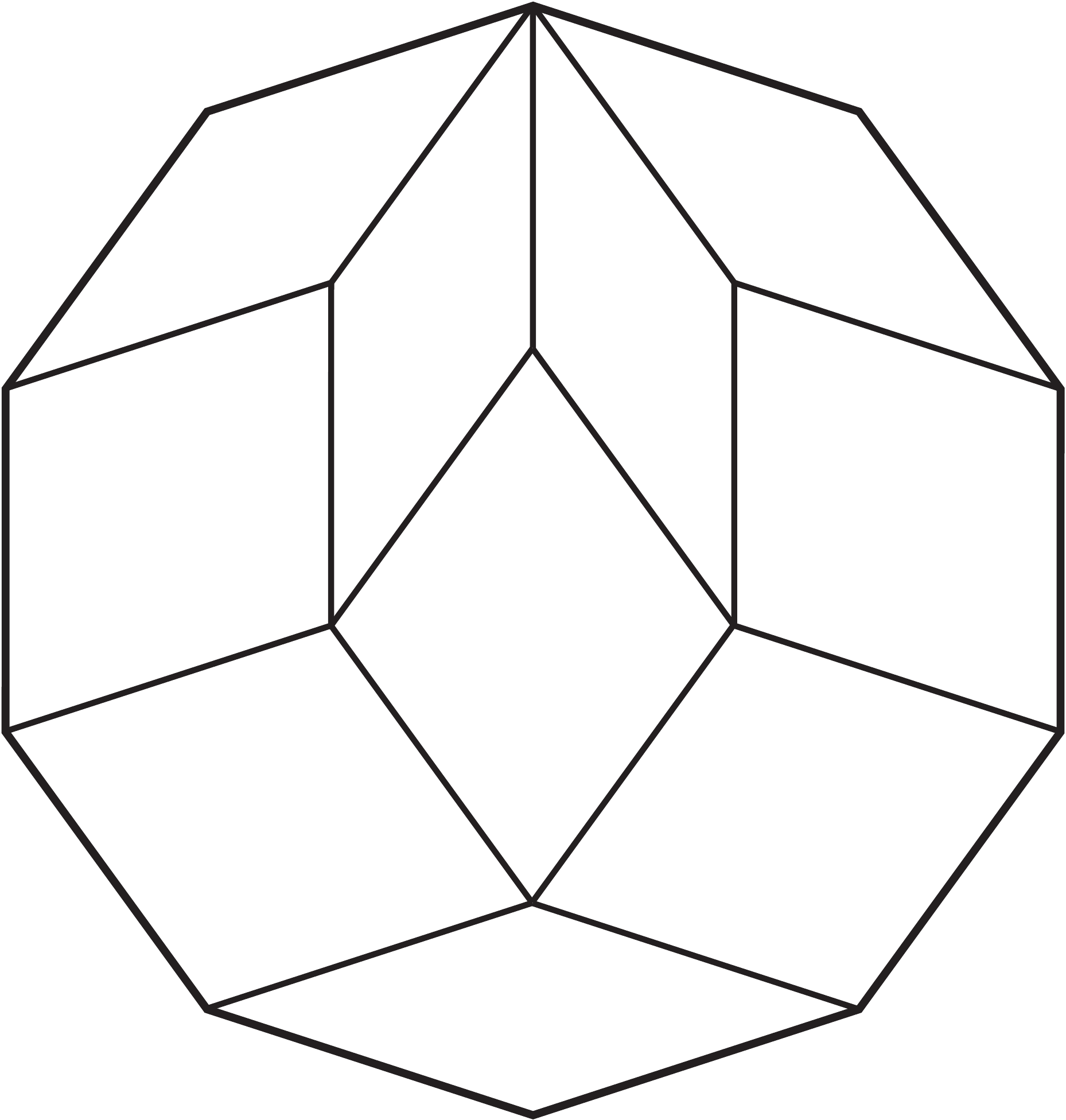}     
 &    \includegraphics[width=1in]{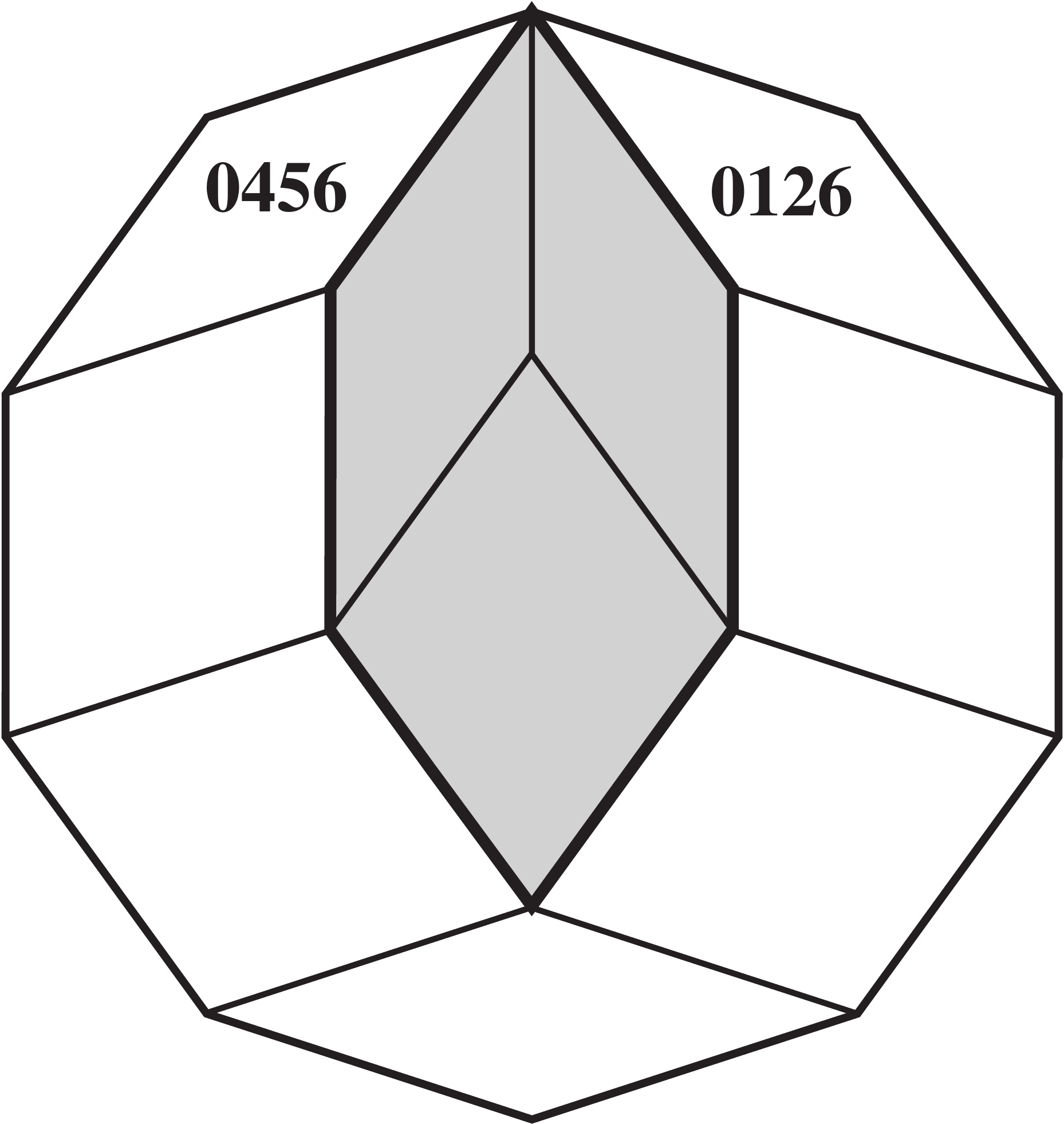}     
&     \includegraphics[width=1in]{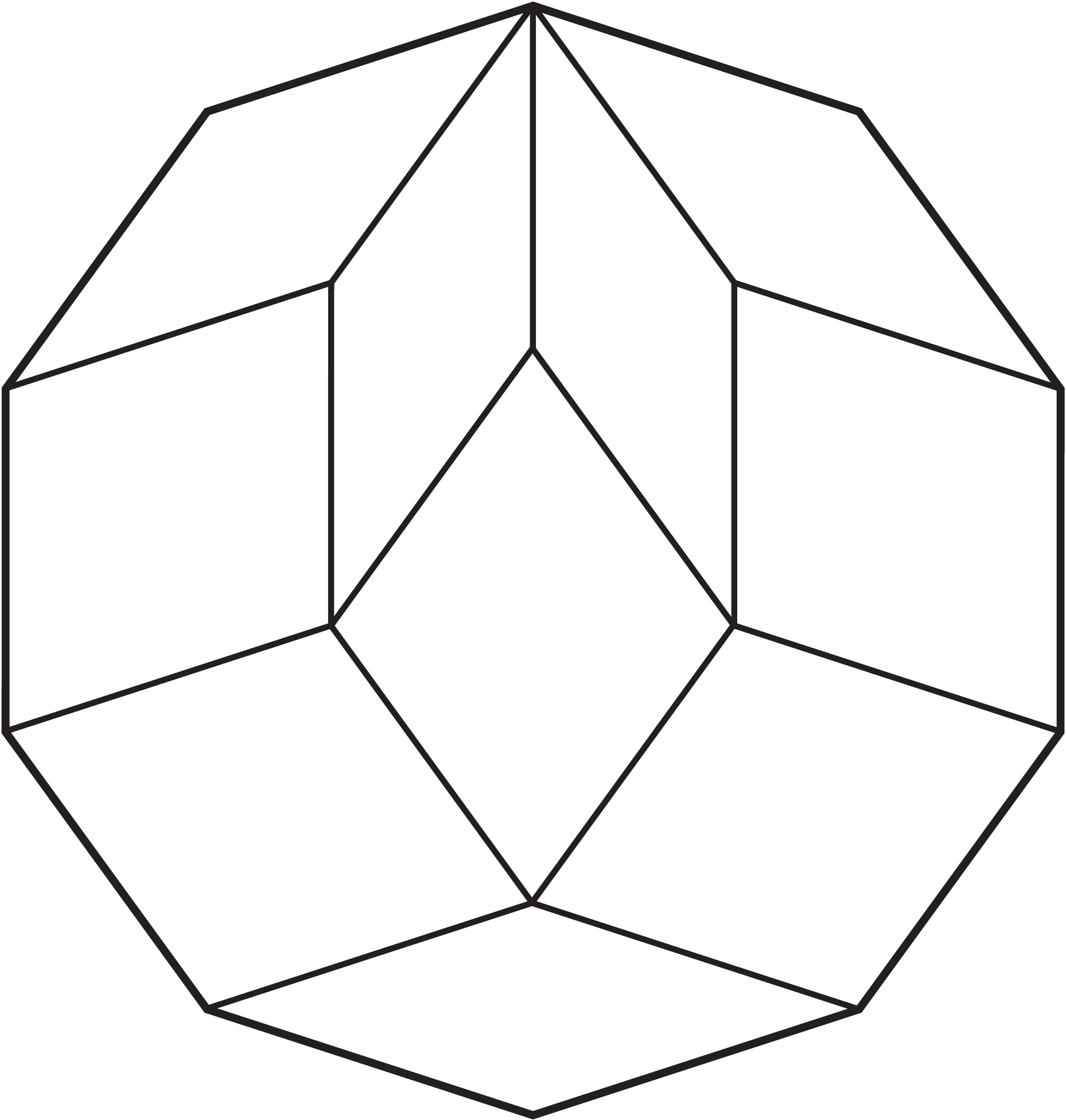}     
&     \includegraphics[width=1in]{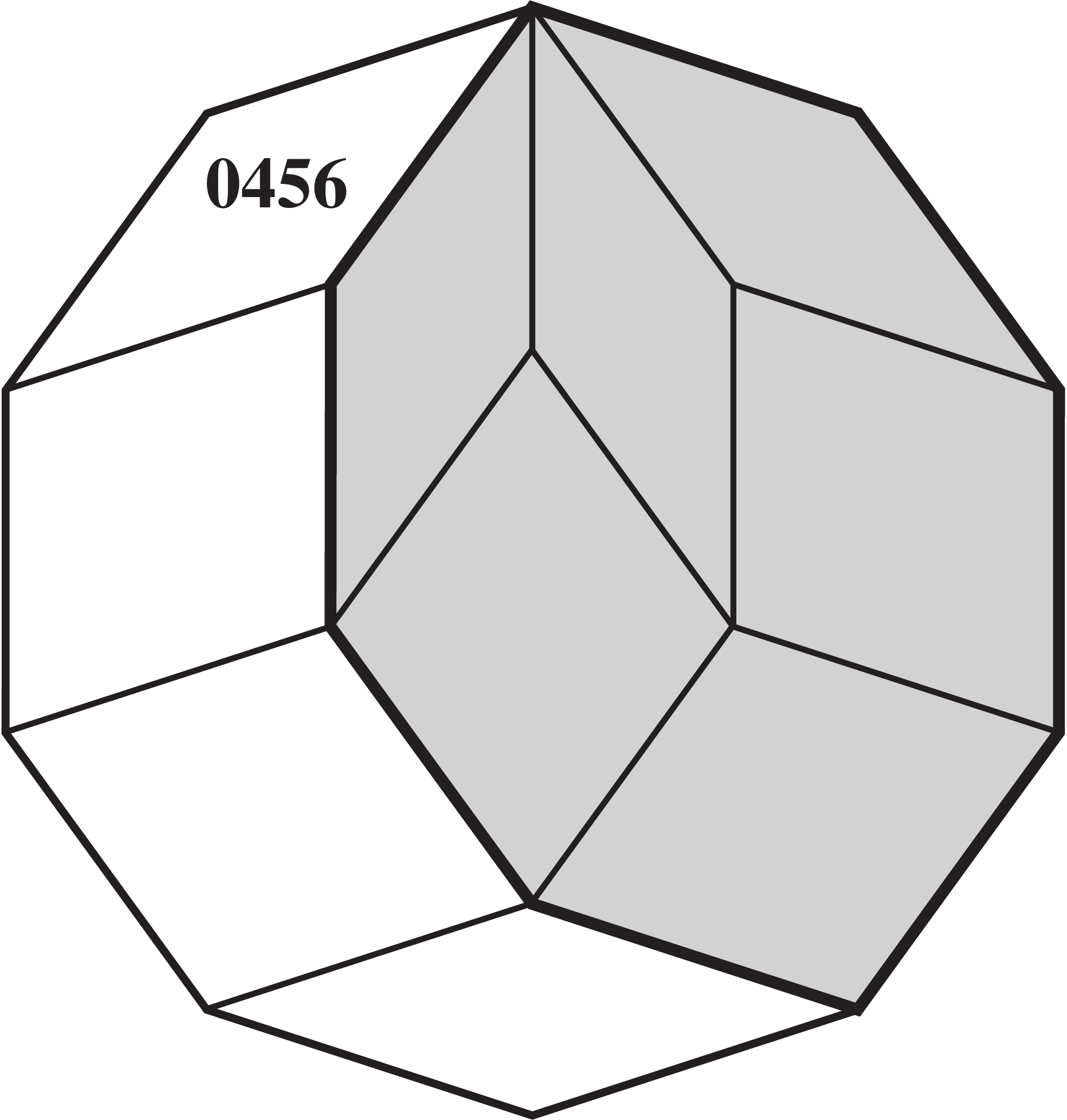} \\

023456 &    &  02346 &    & 012346  \\

     \includegraphics[width=1in]{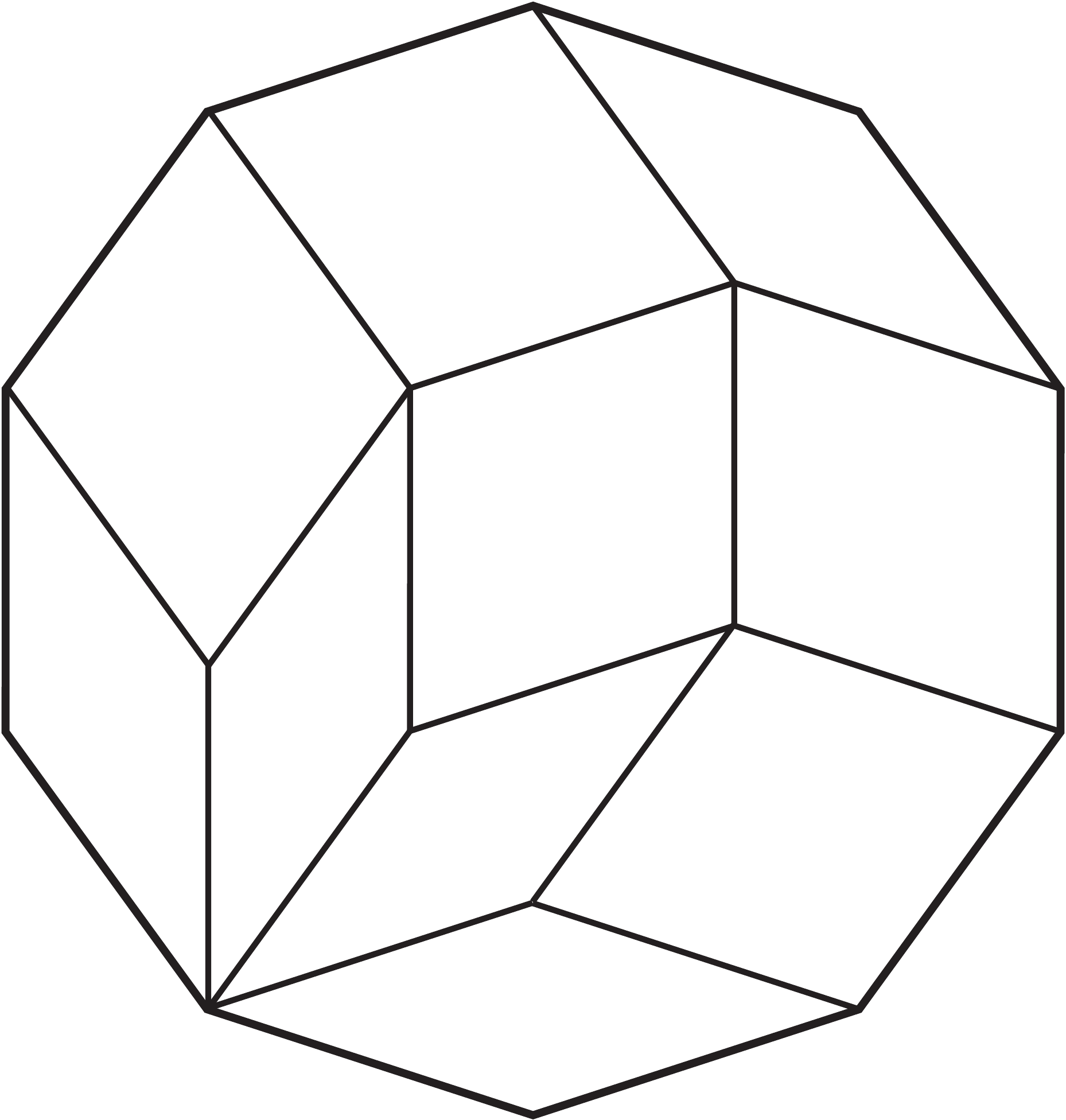}     
  &      \includegraphics[width=1in]{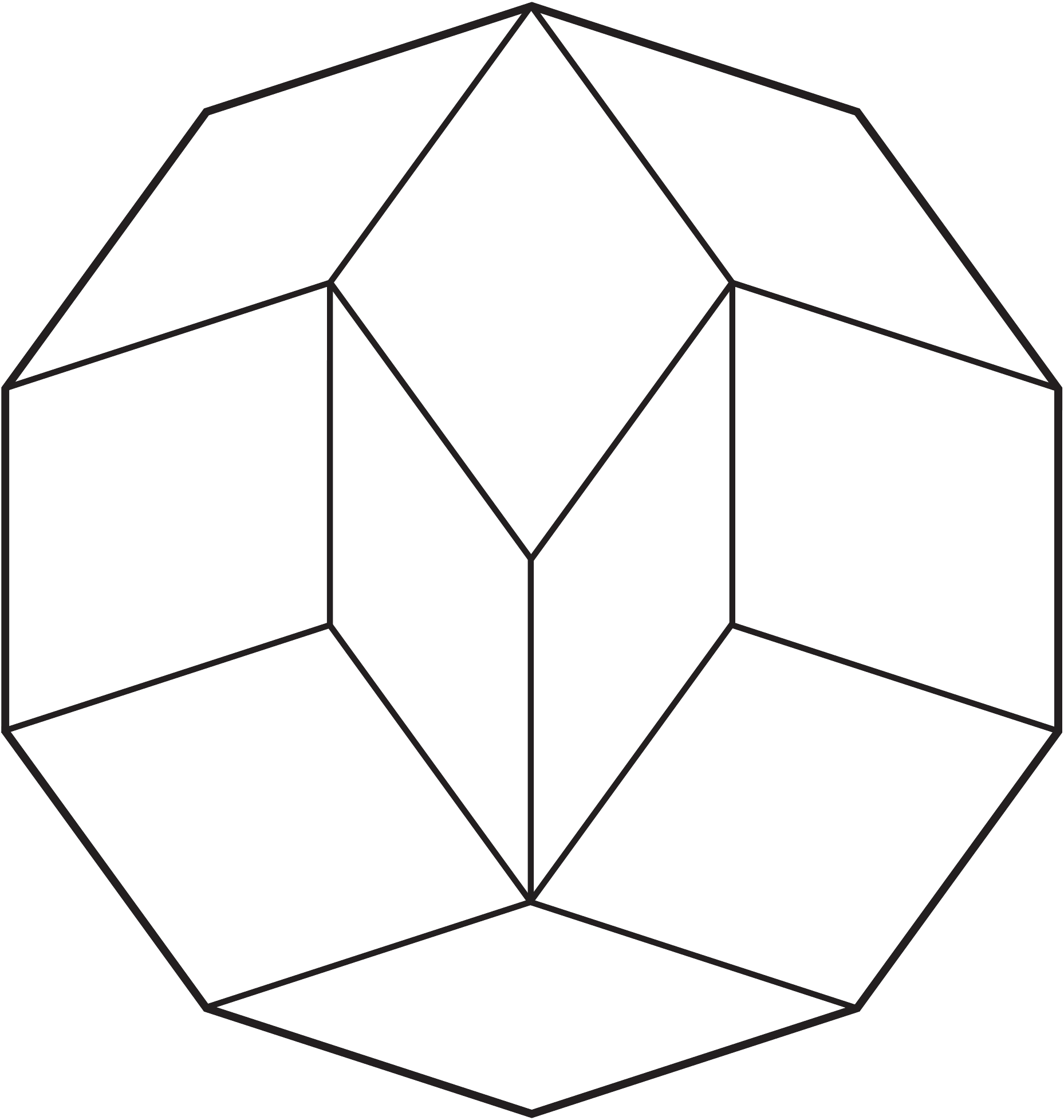}     
 &    \includegraphics[width=1in]{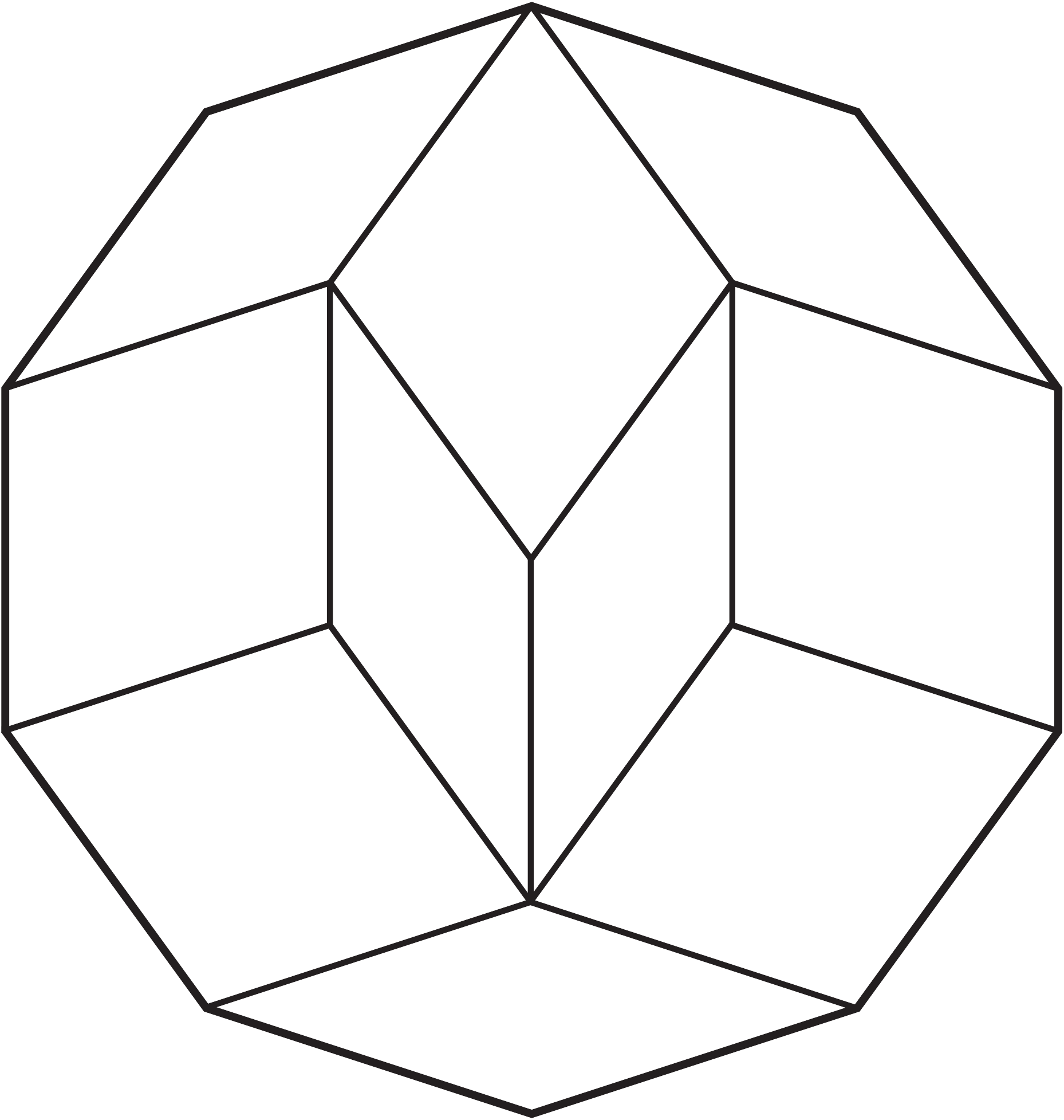}     
&     \includegraphics[width=1in]{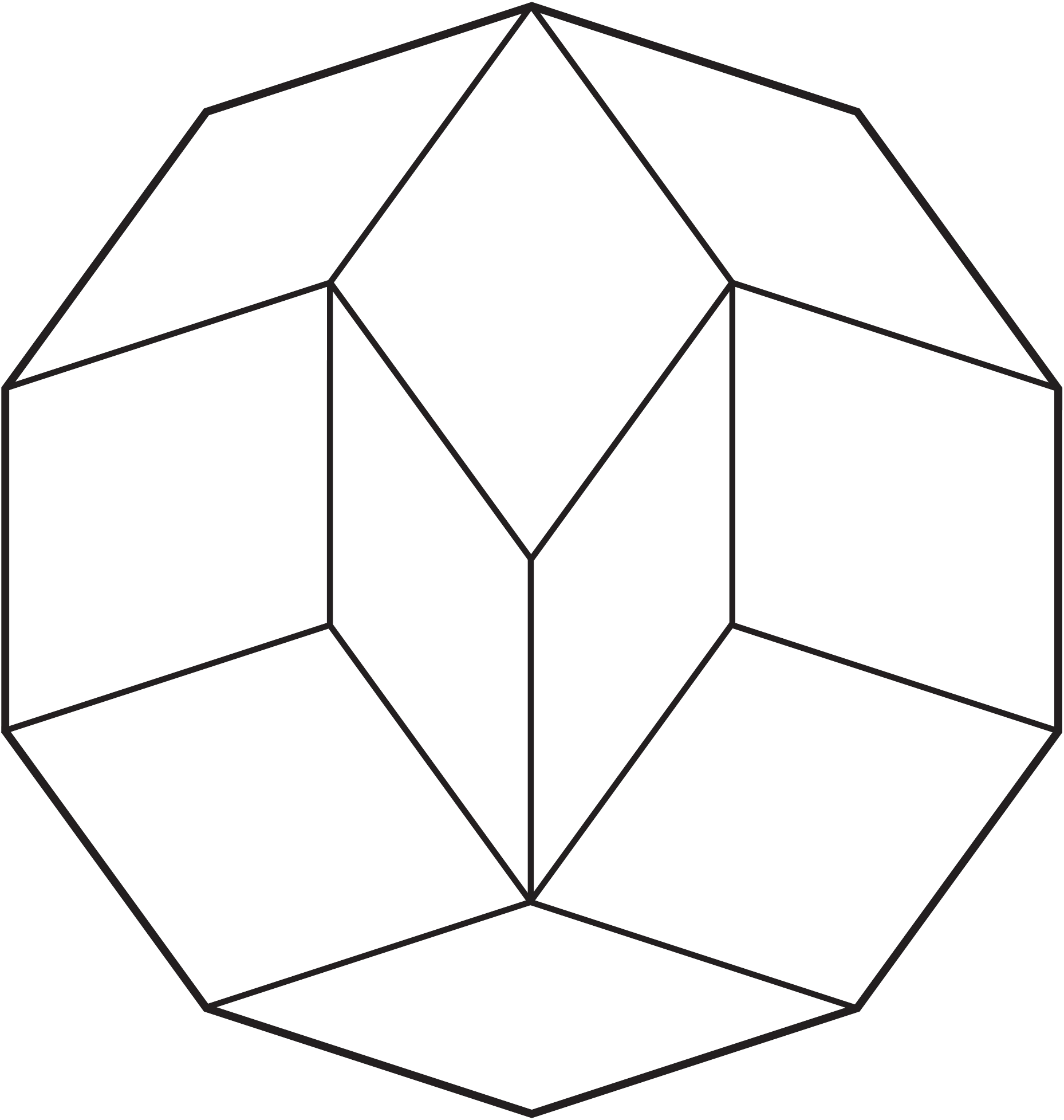}     
&     \includegraphics[width=1in]{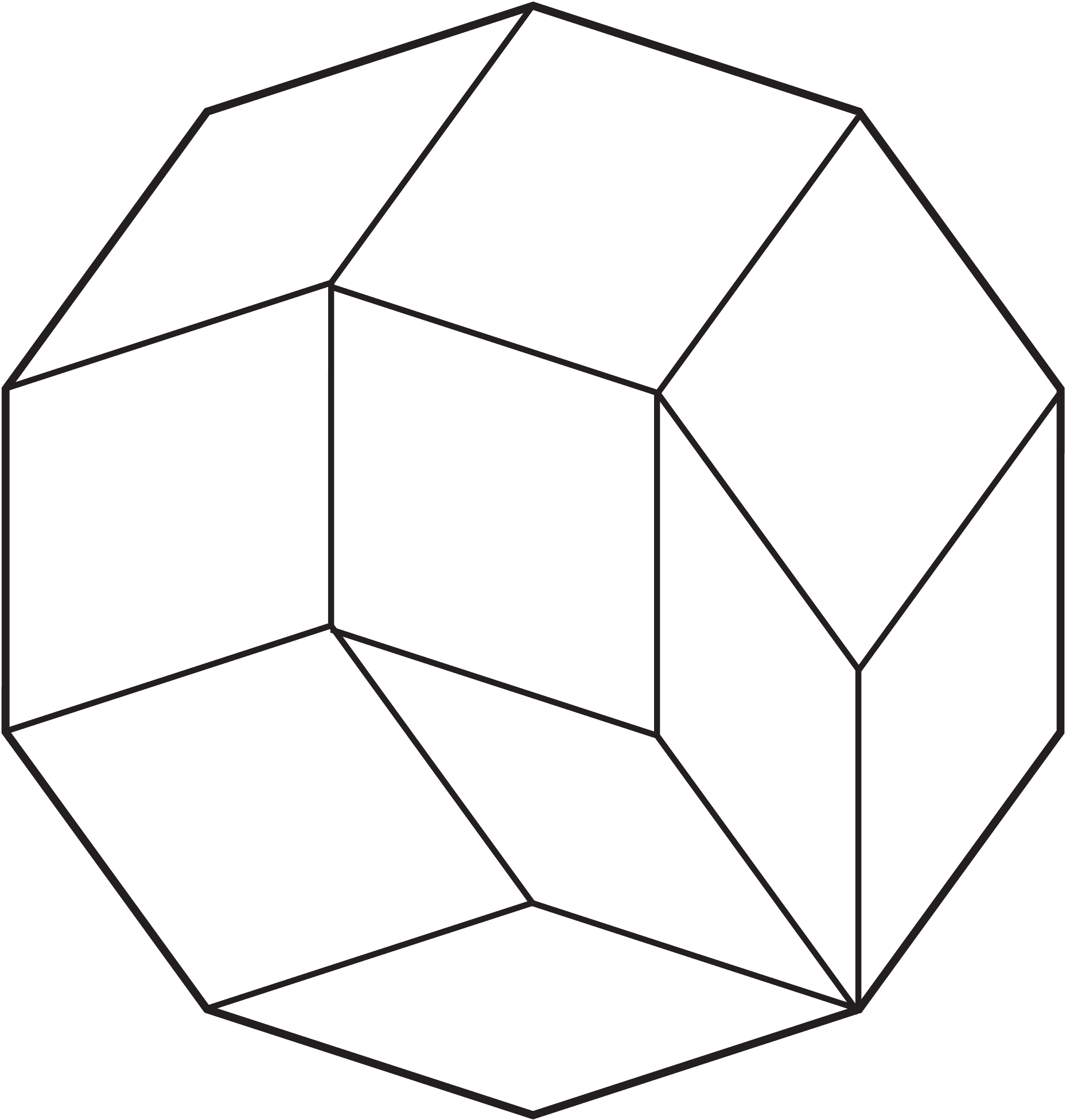} \\

  &    &   &     &    \\

      \includegraphics[width=1in]{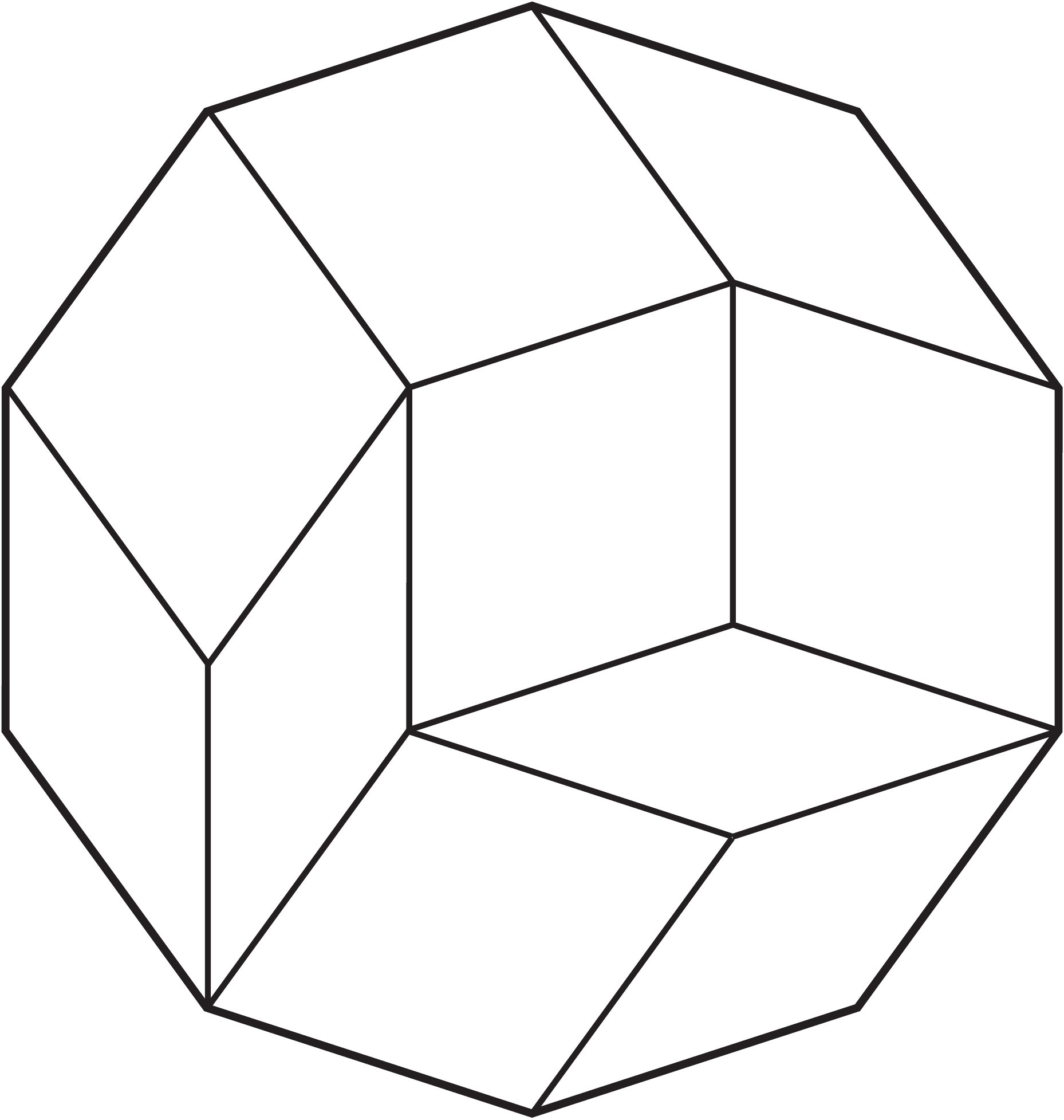}     
  &      \includegraphics[width=1in]{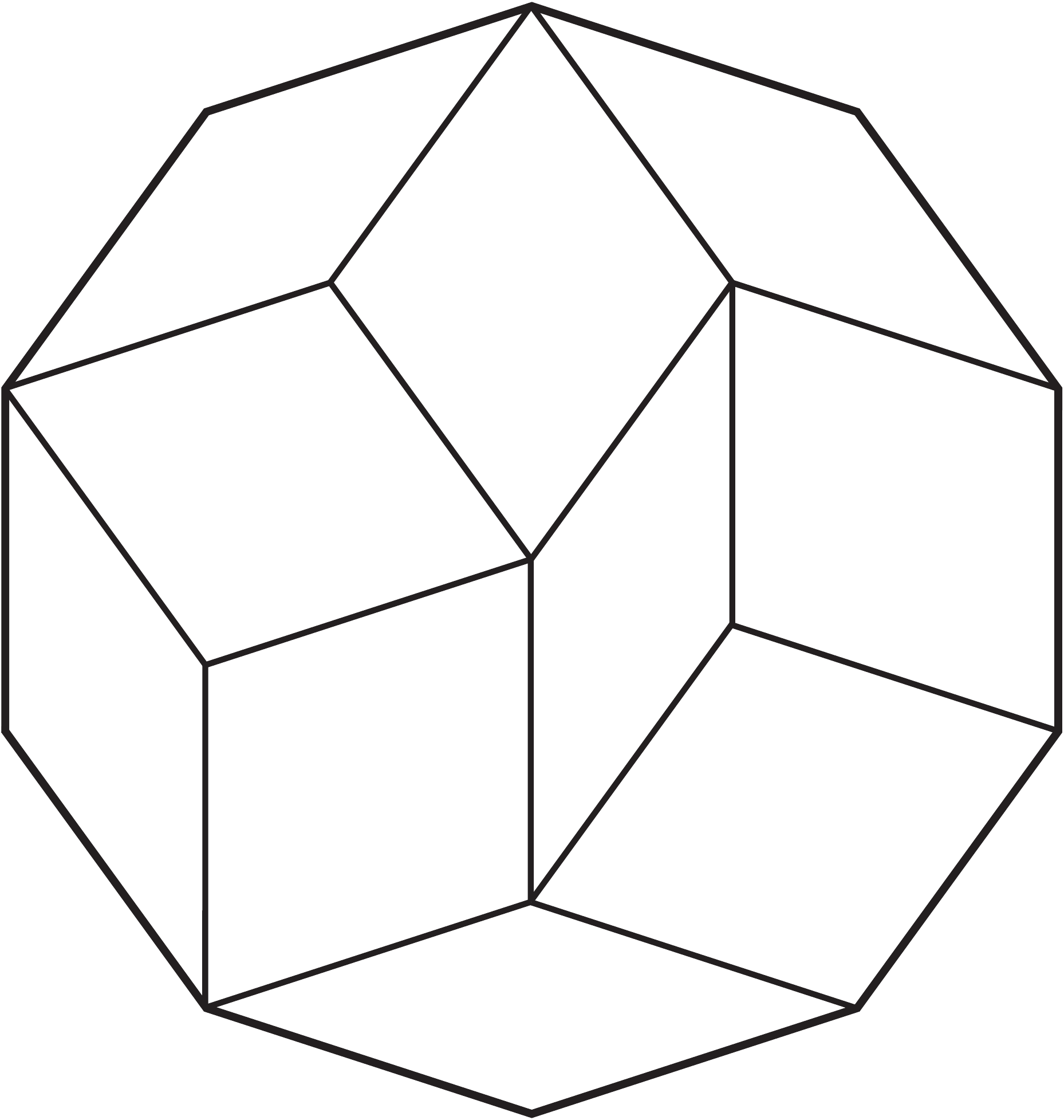}     
 &    \includegraphics[width=1in]{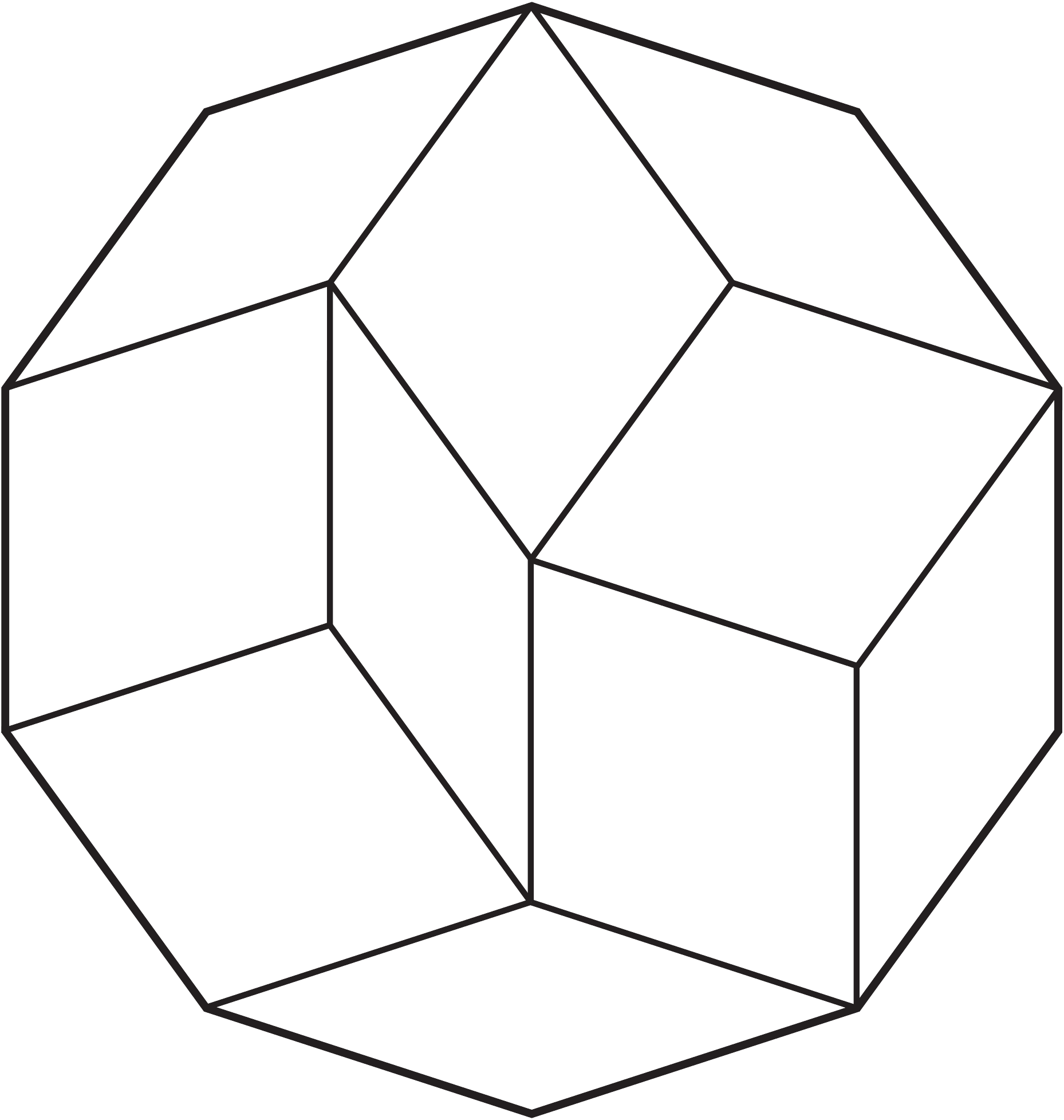}     
&     \includegraphics[width=1in]{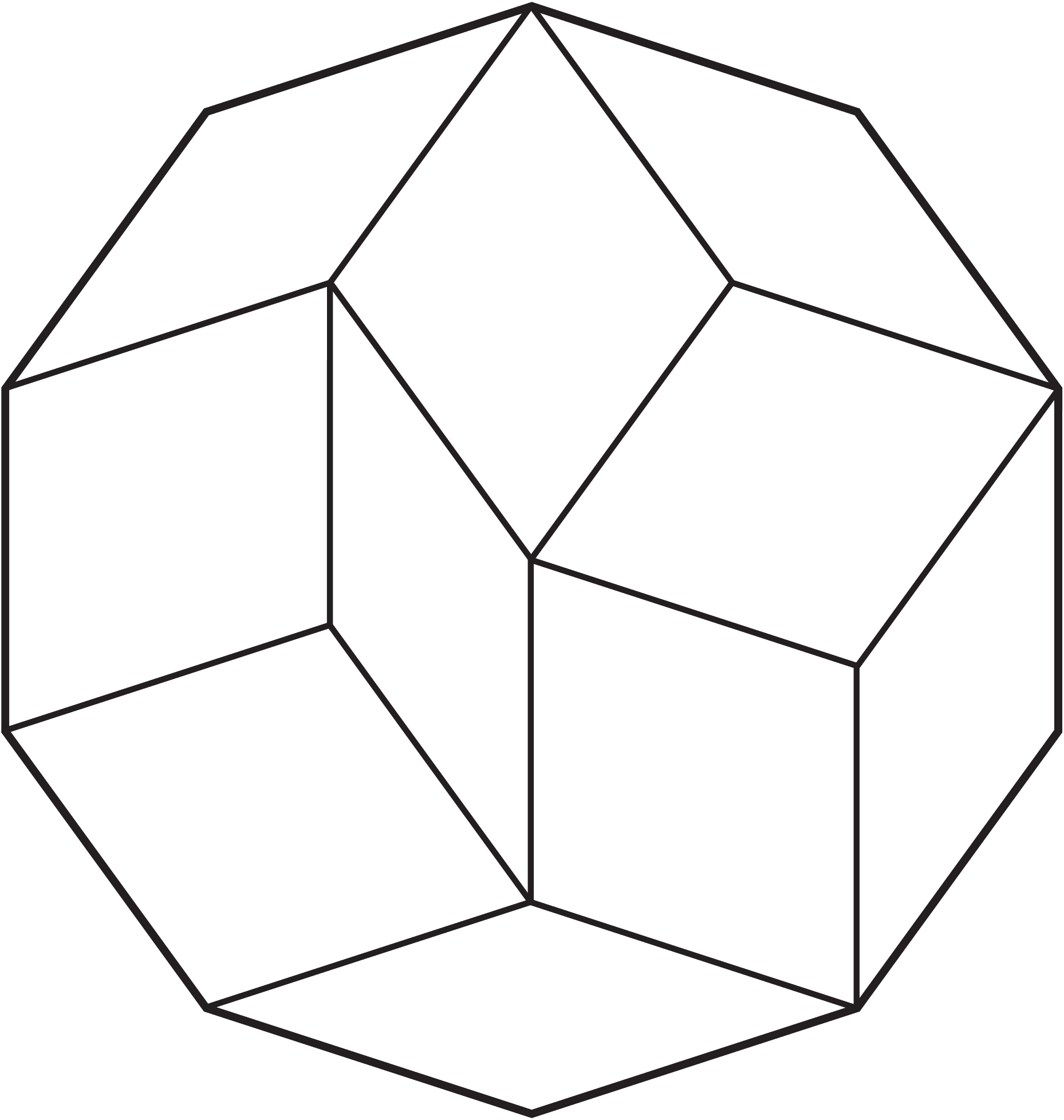}     
&     \includegraphics[width=1in]{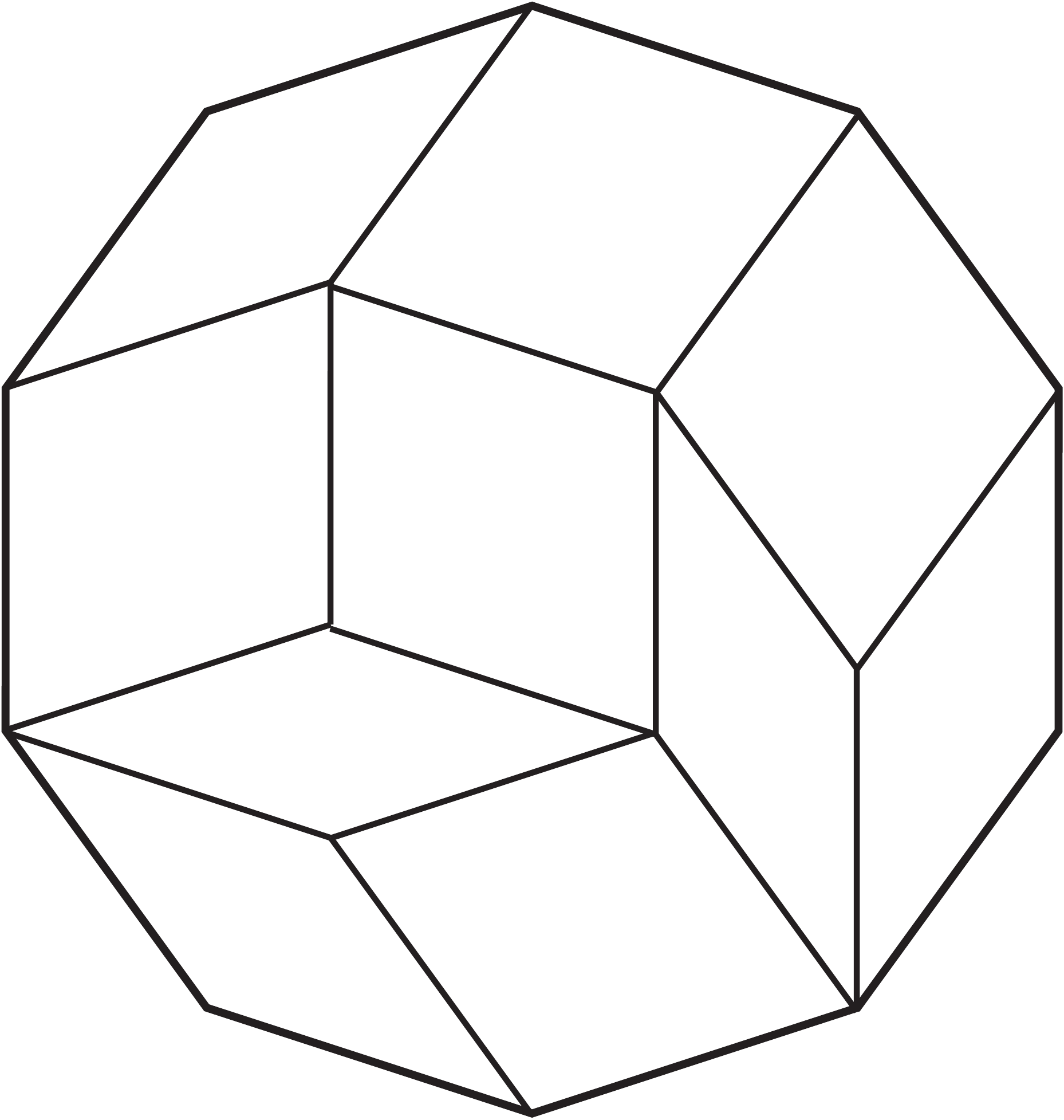} \\

  &    &   &     &    \\
  
      \includegraphics[width=1in]{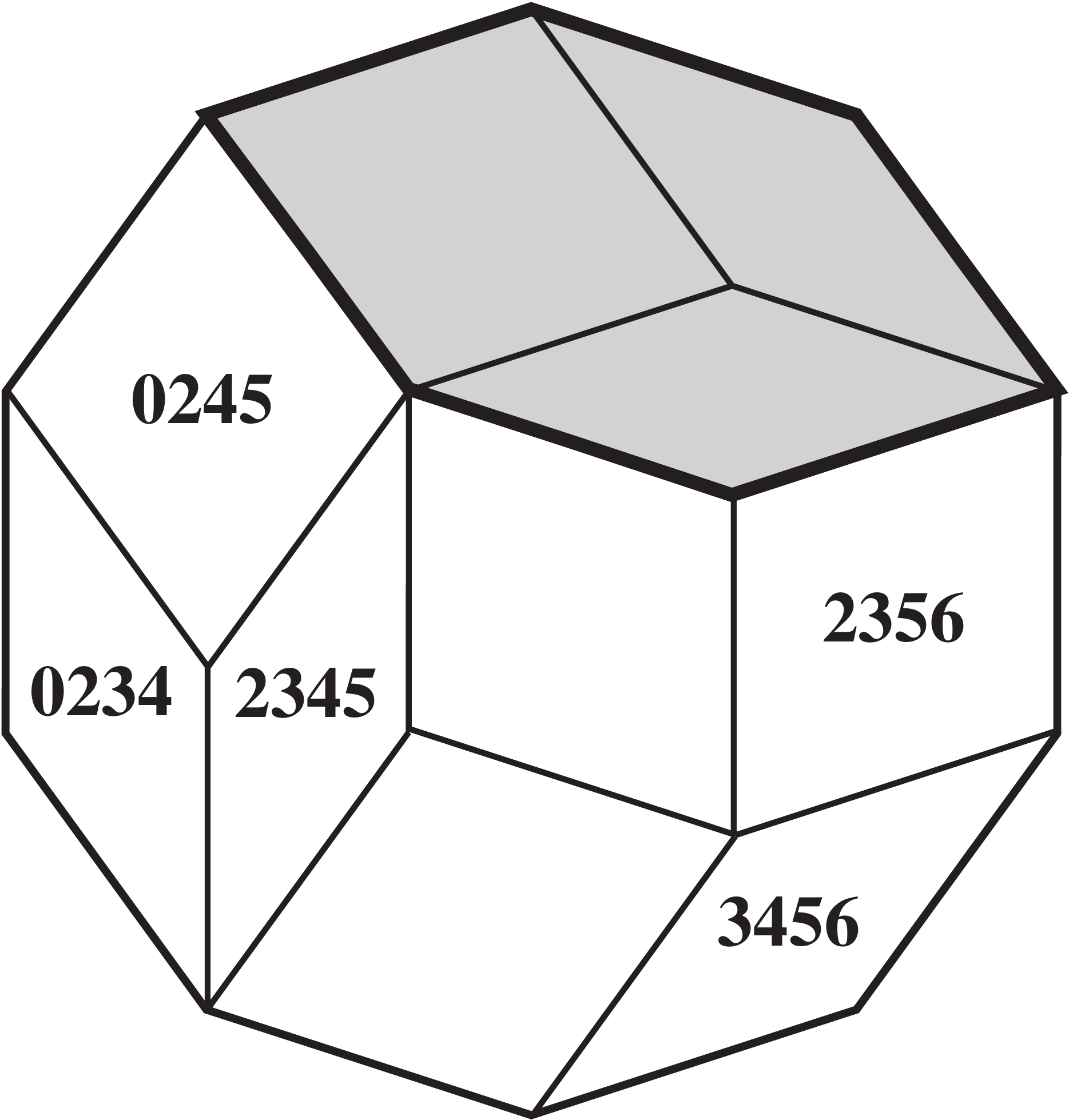}     
  &      \includegraphics[width=1in]{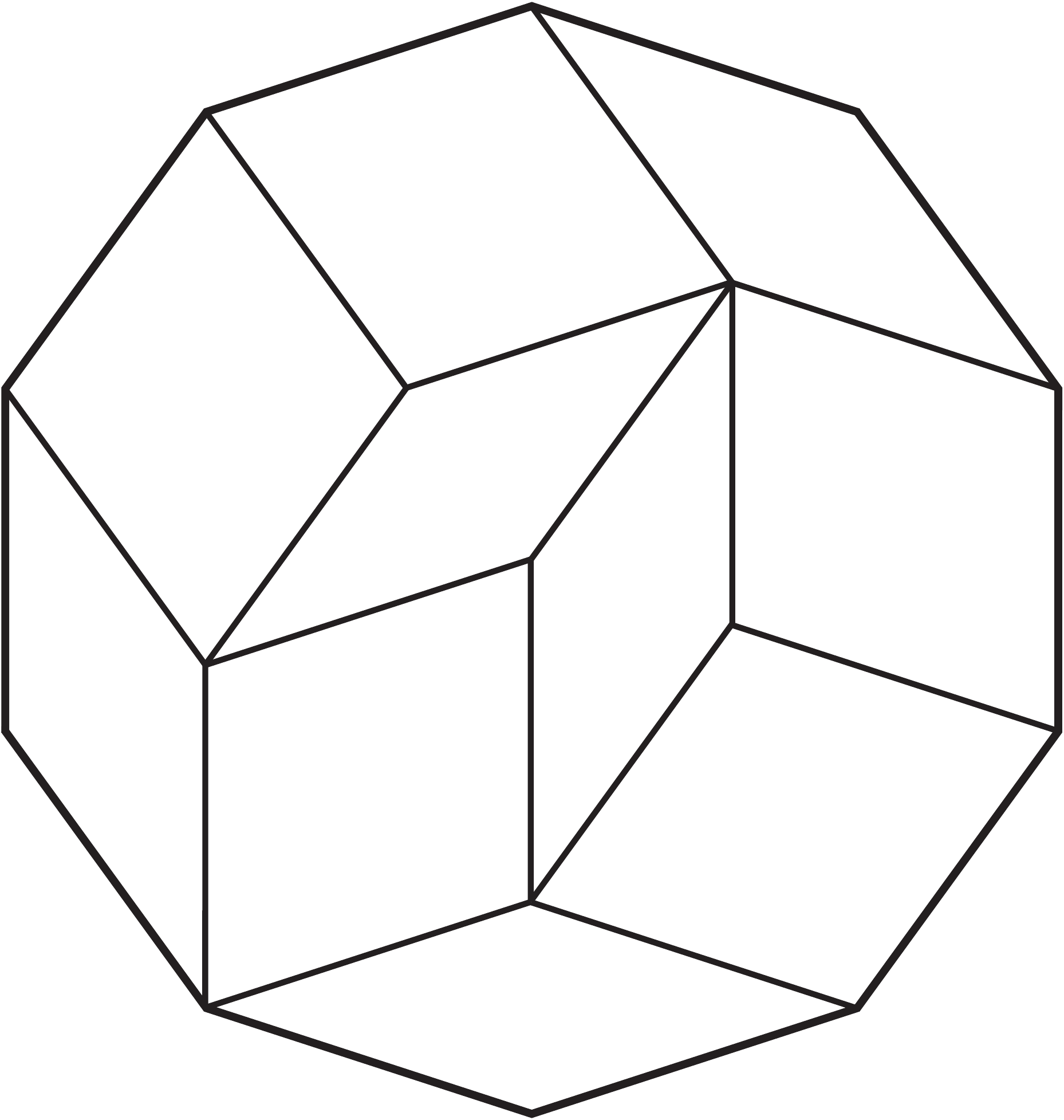}     
 &    \includegraphics[width=1in]{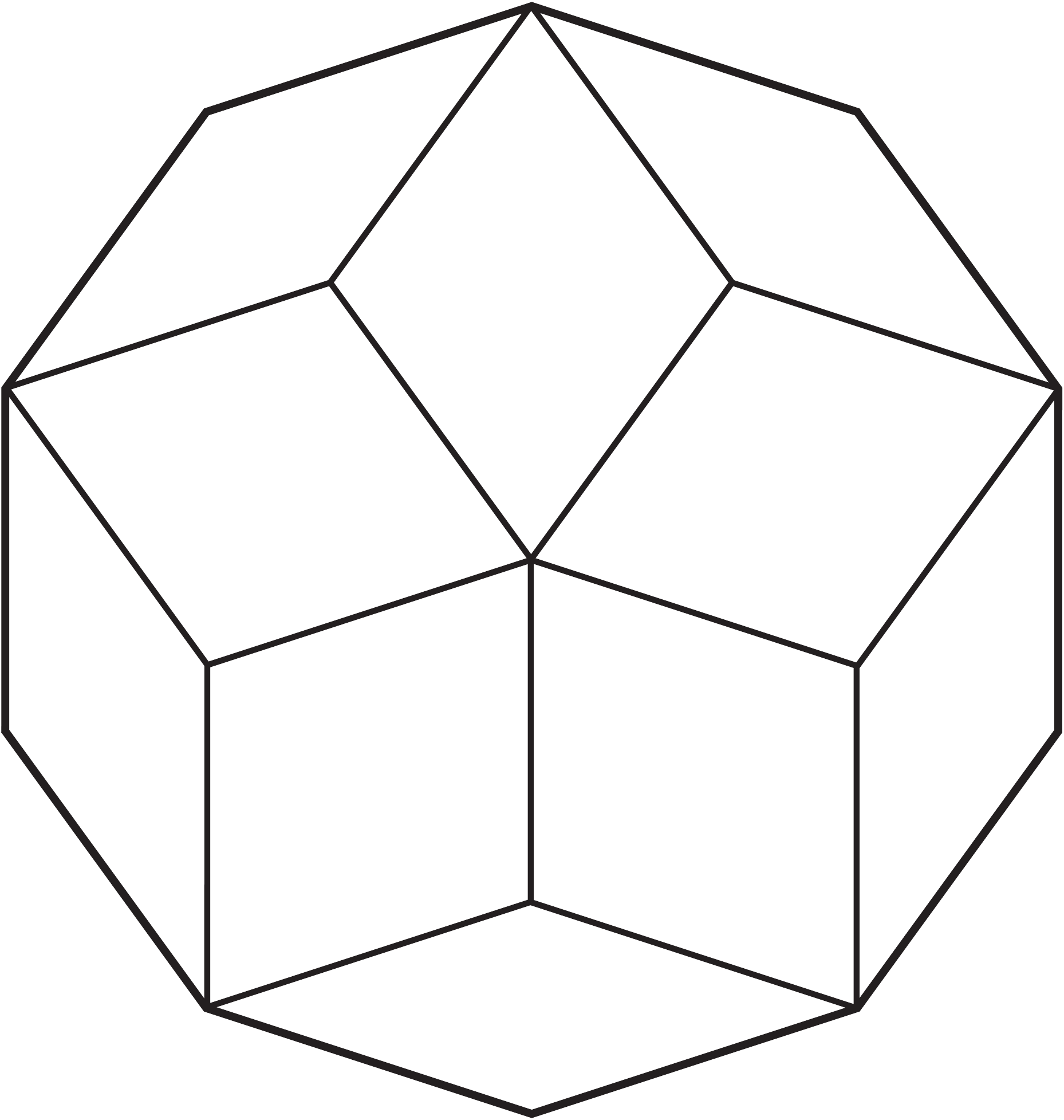}     
&     \includegraphics[width=1in]{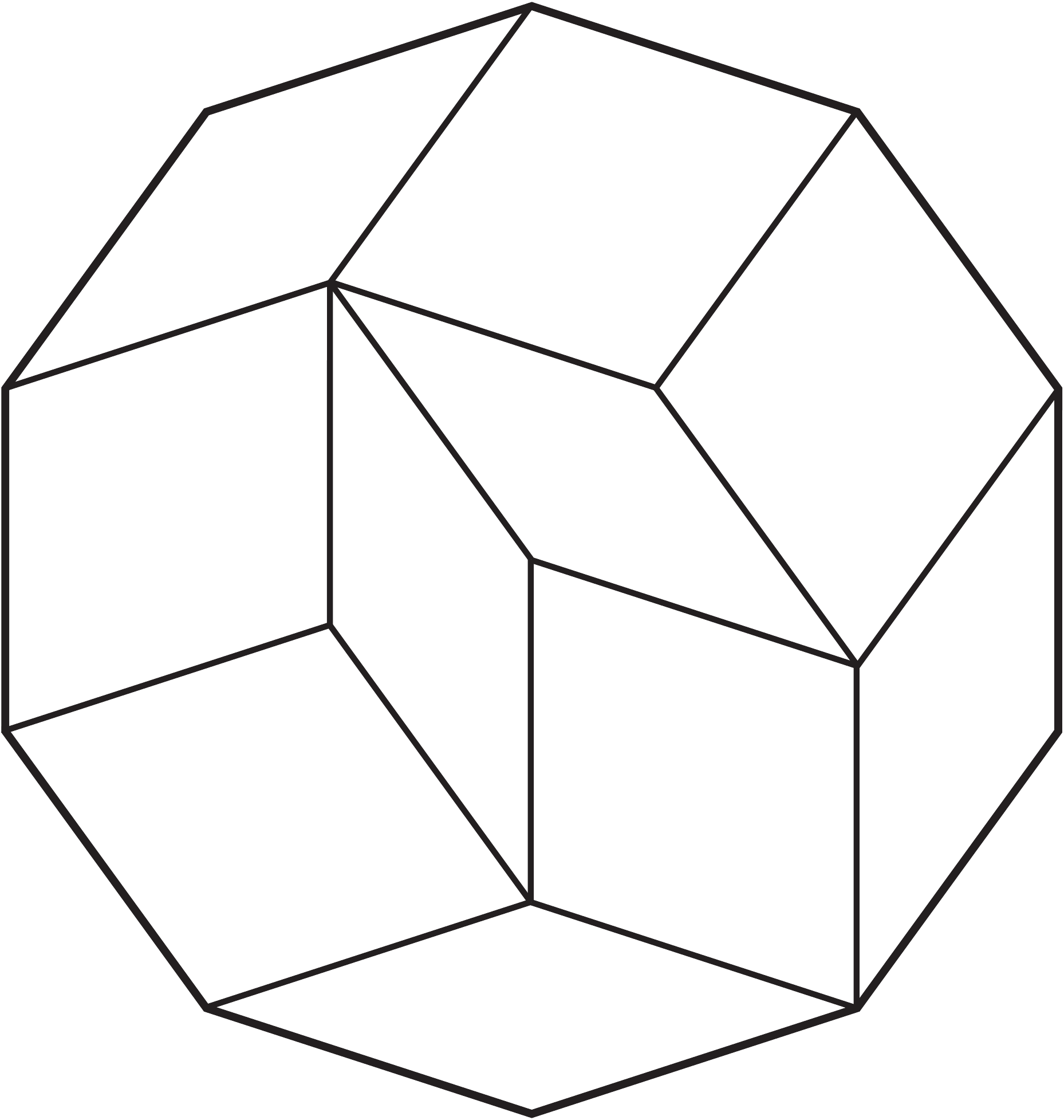}     
&     \includegraphics[width=1in]{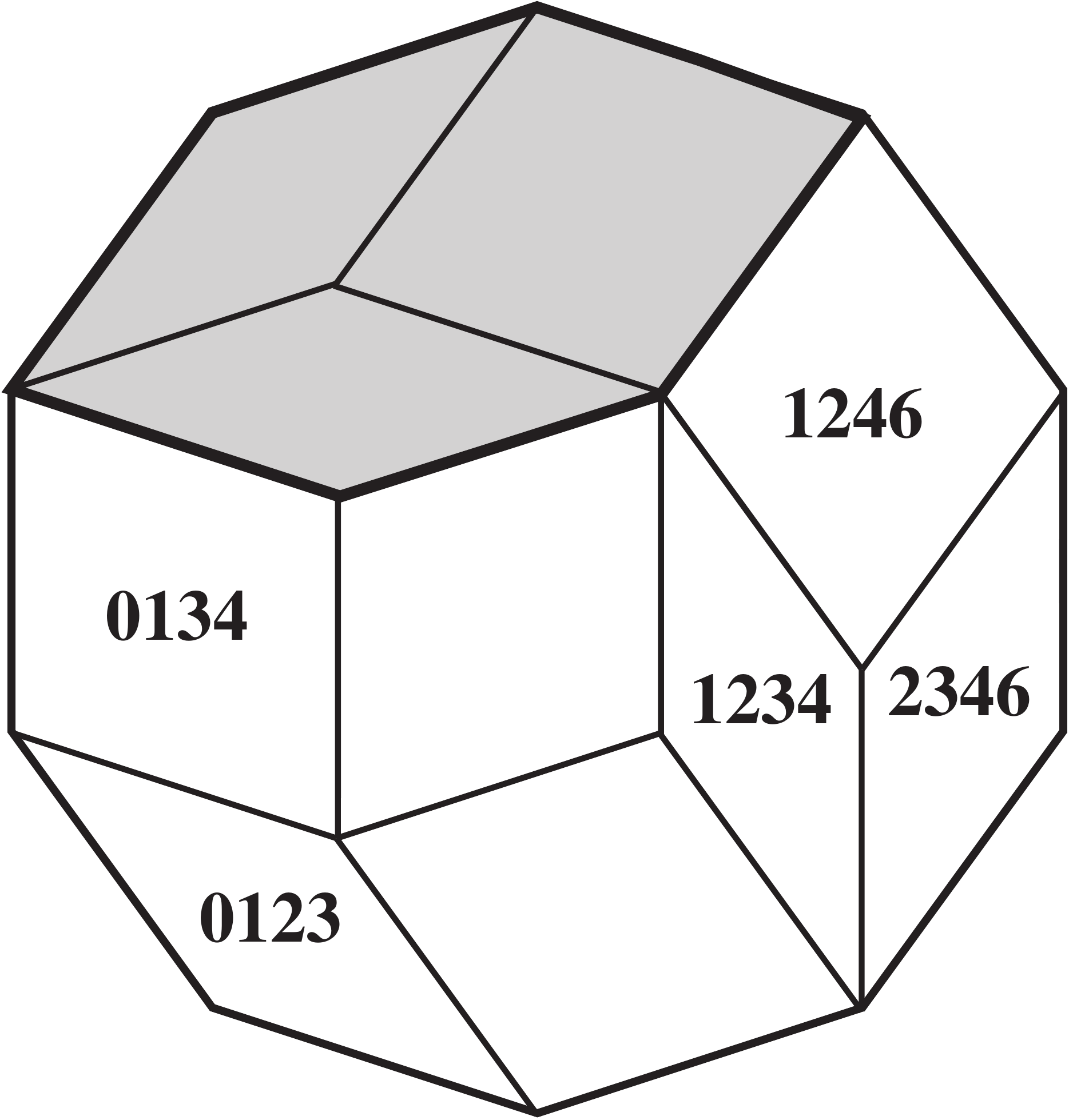} \\

01256  &    &   &     & 01456   \\

     \includegraphics[width=1in]{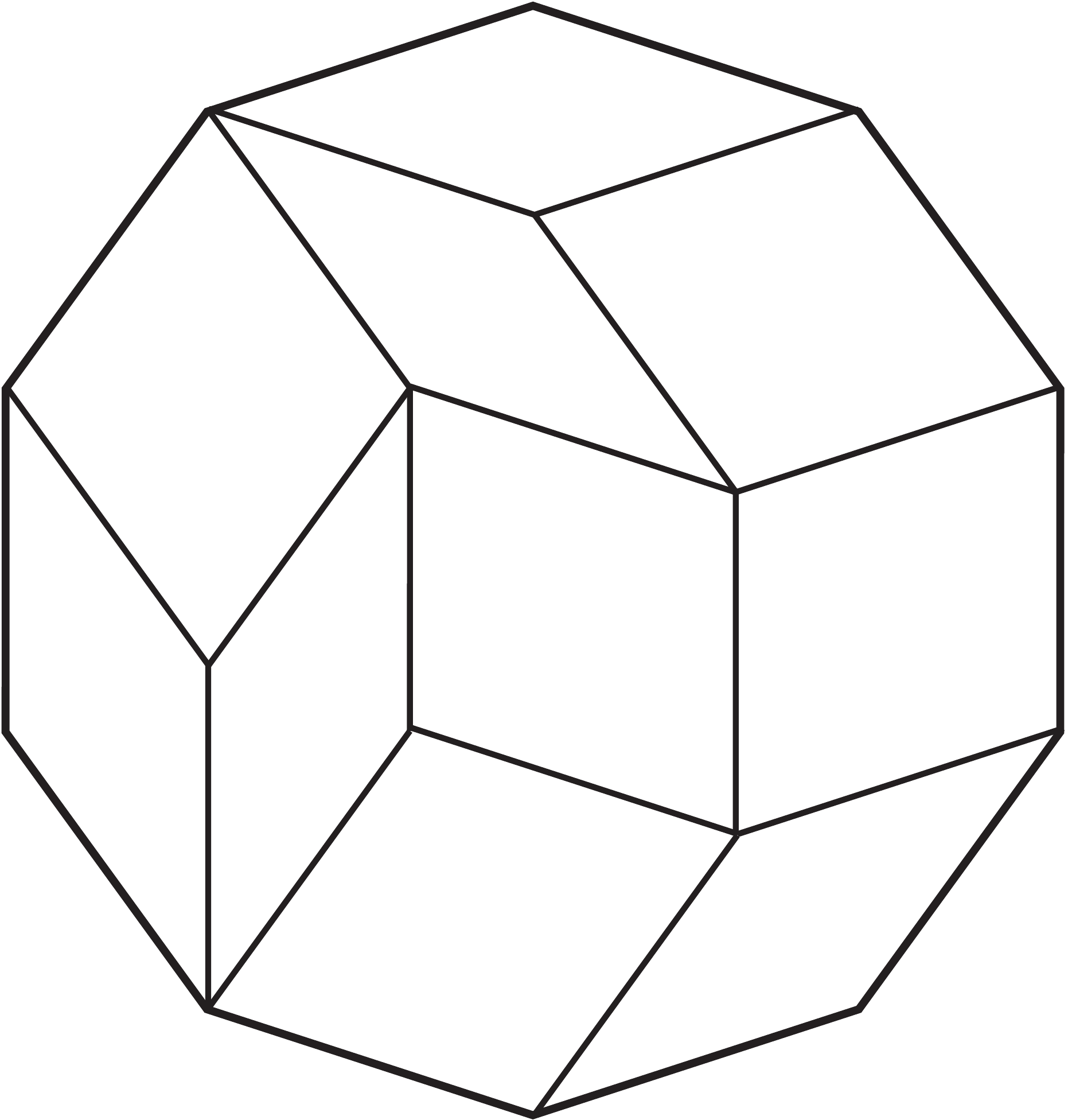}     
  &      \includegraphics[width=1in]{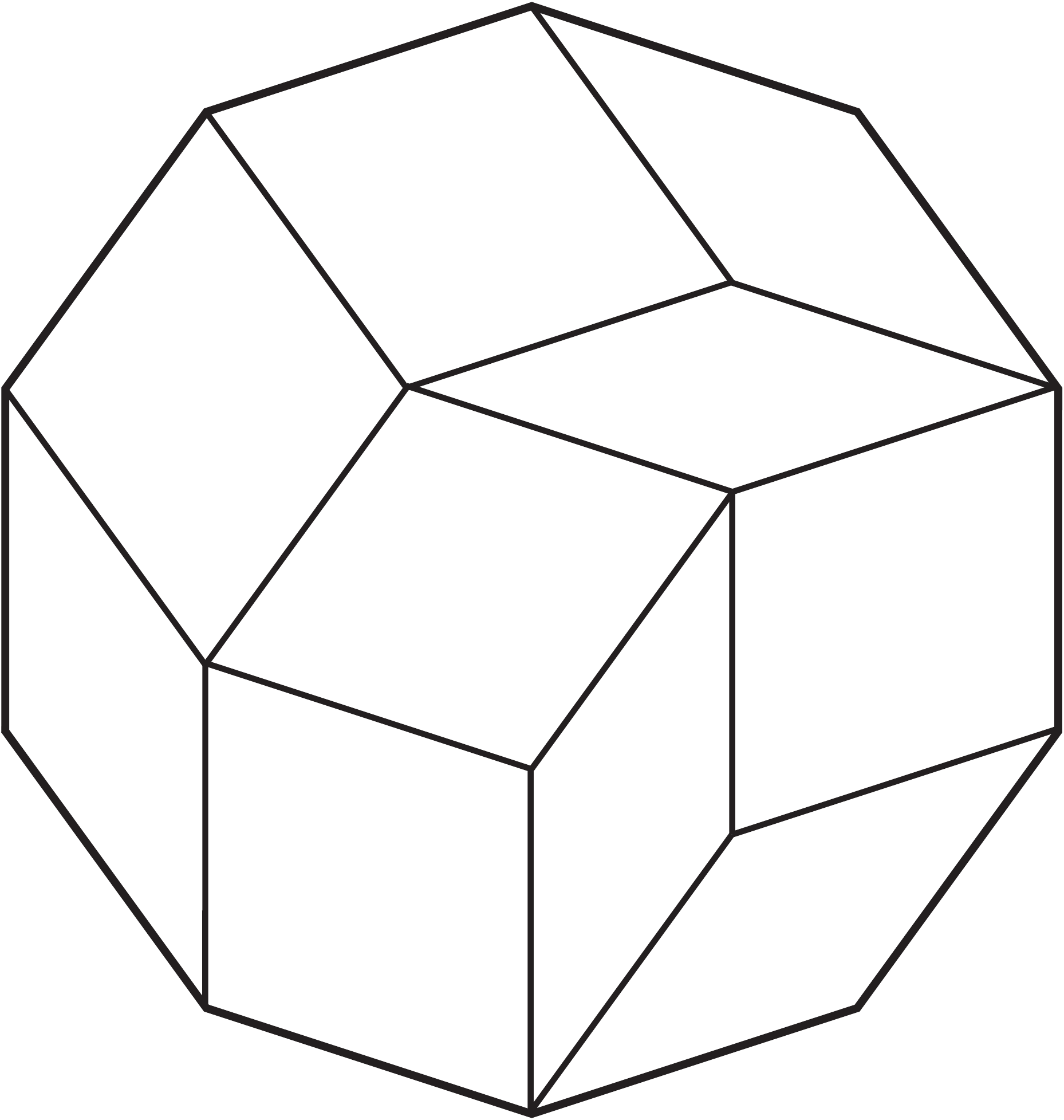}     
 &    \includegraphics[width=1in]{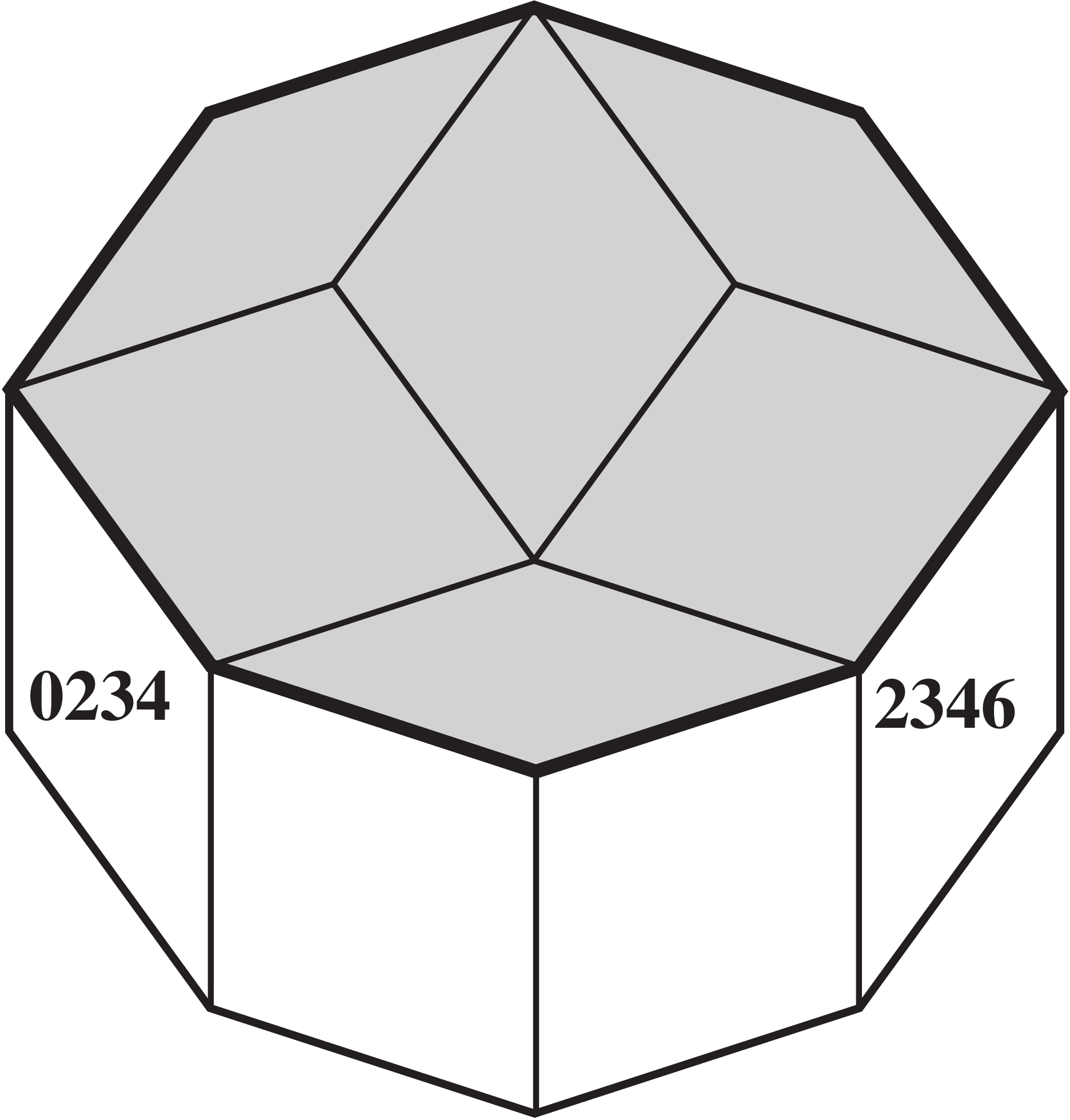}     
&     \includegraphics[width=1in]{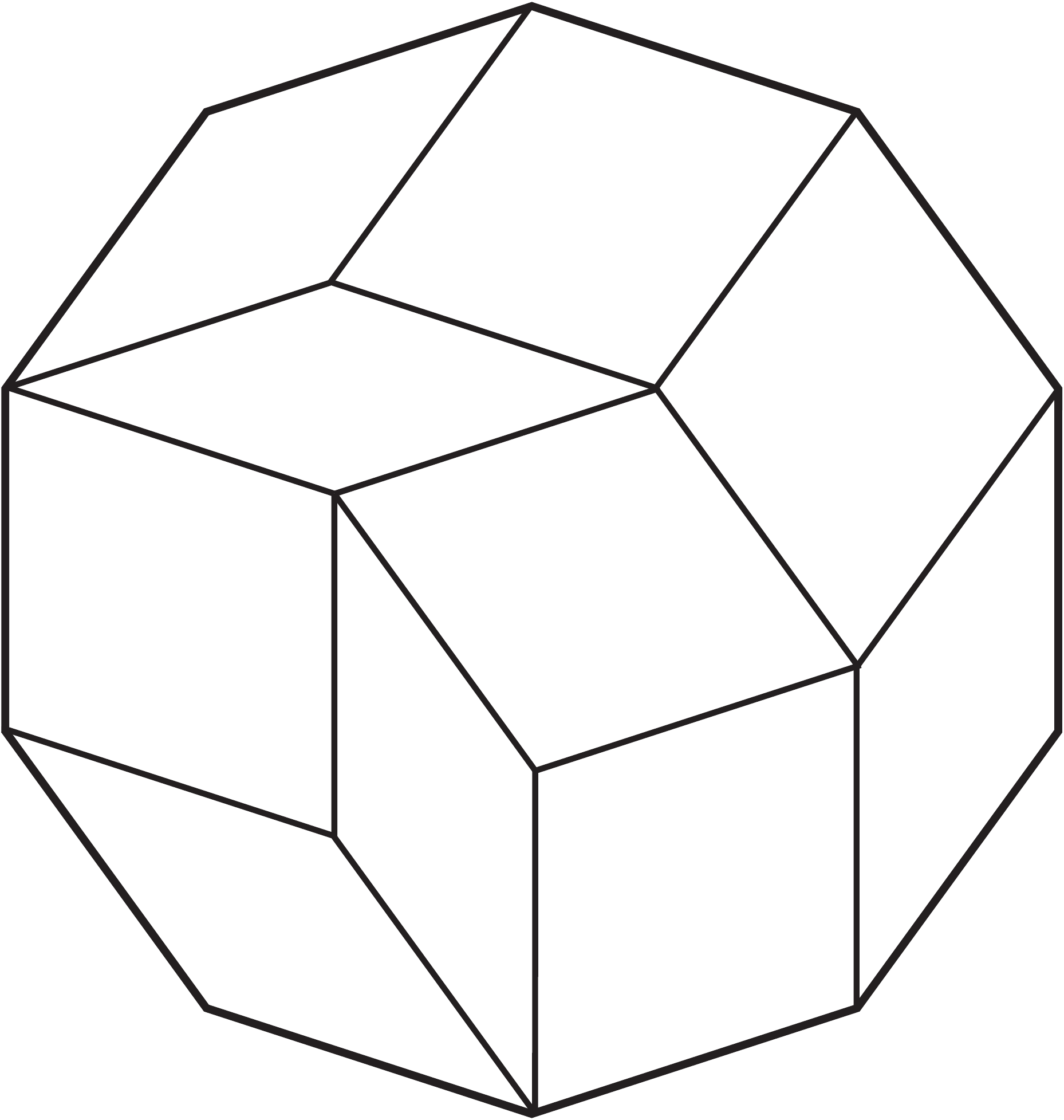}     
&     \includegraphics[width=1in]{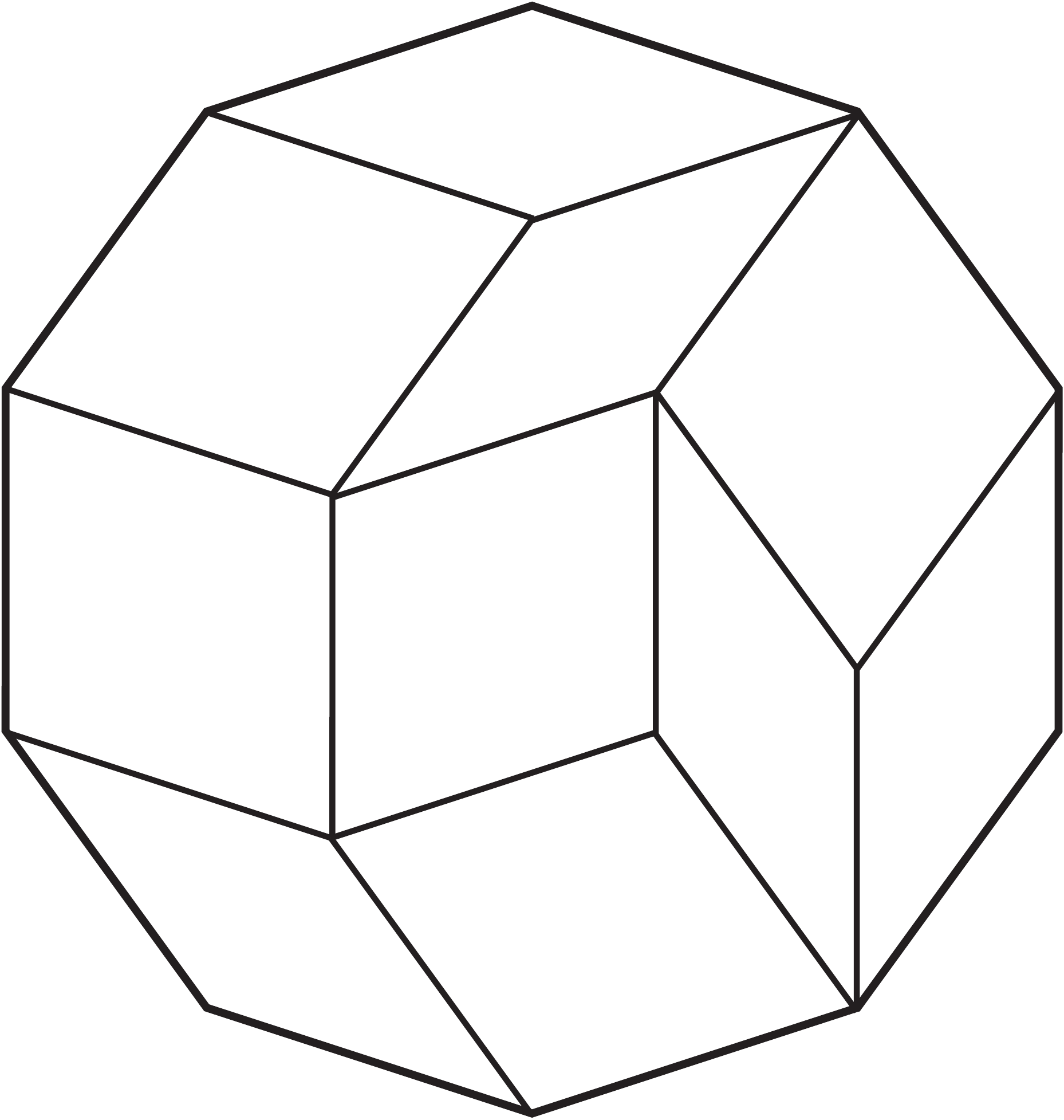} \\

   &    & 012456  &     &    \\

     \includegraphics[width=1in]{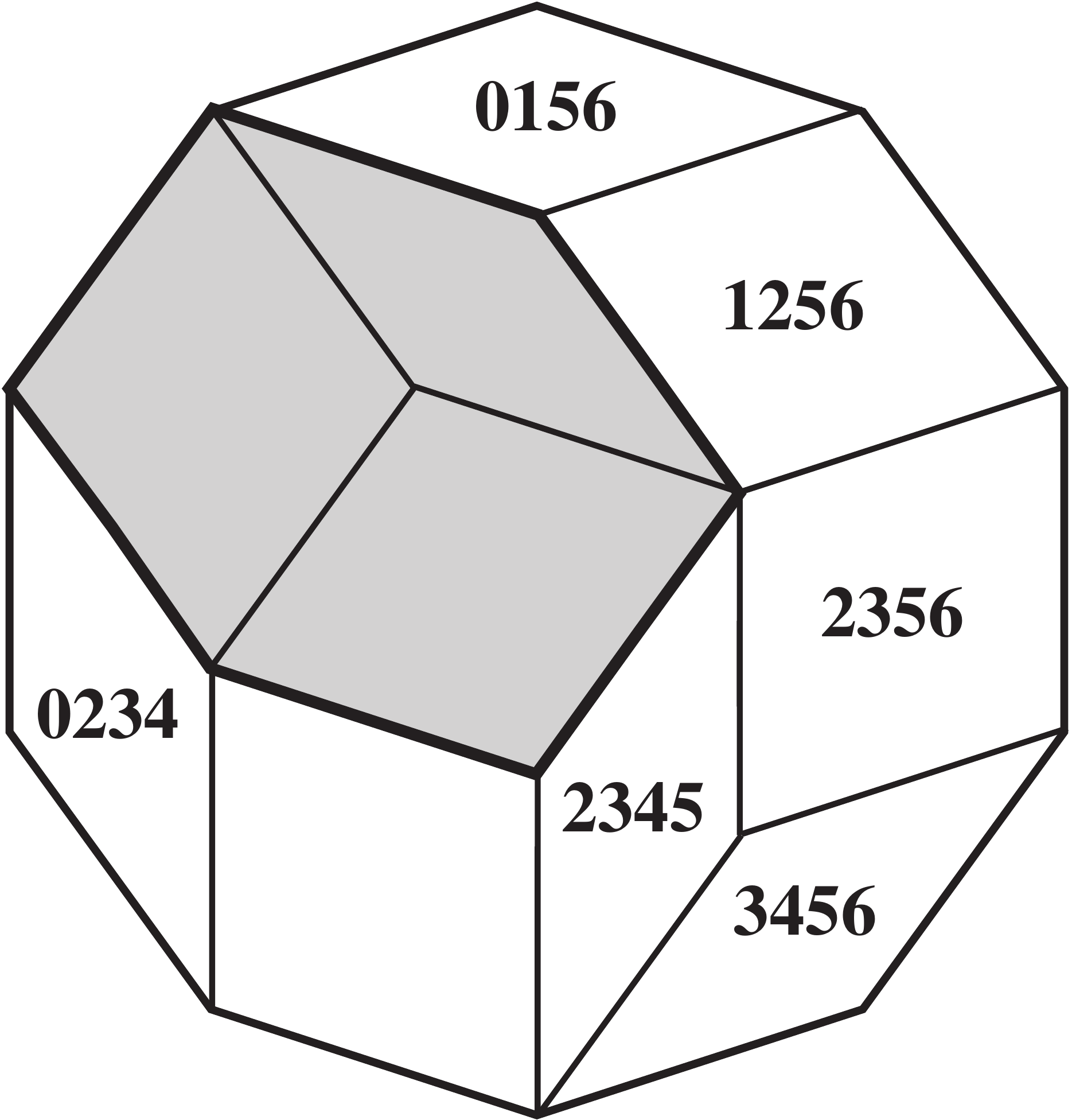}     
  &      \includegraphics[width=1in]{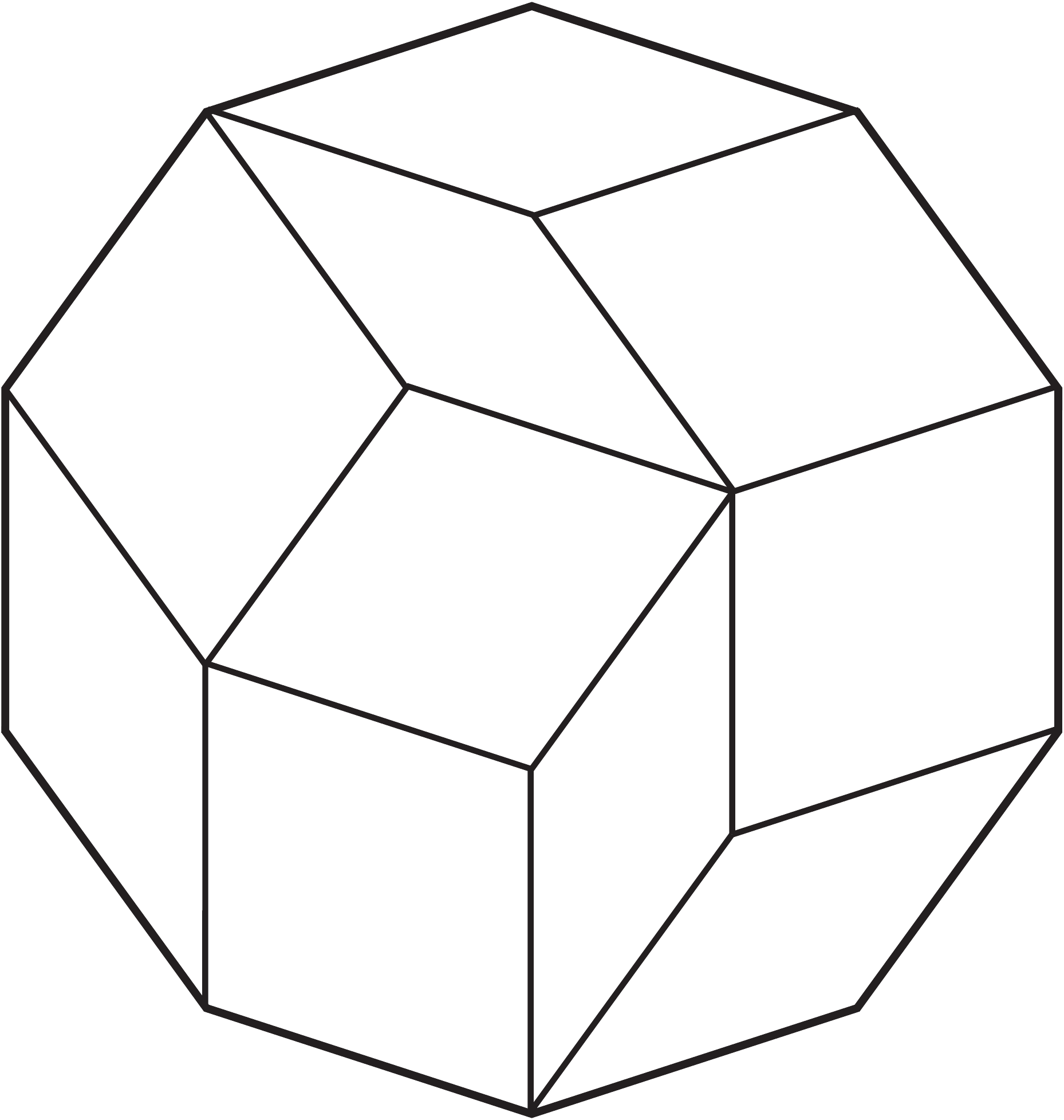}     
 &    \includegraphics[width=1in]{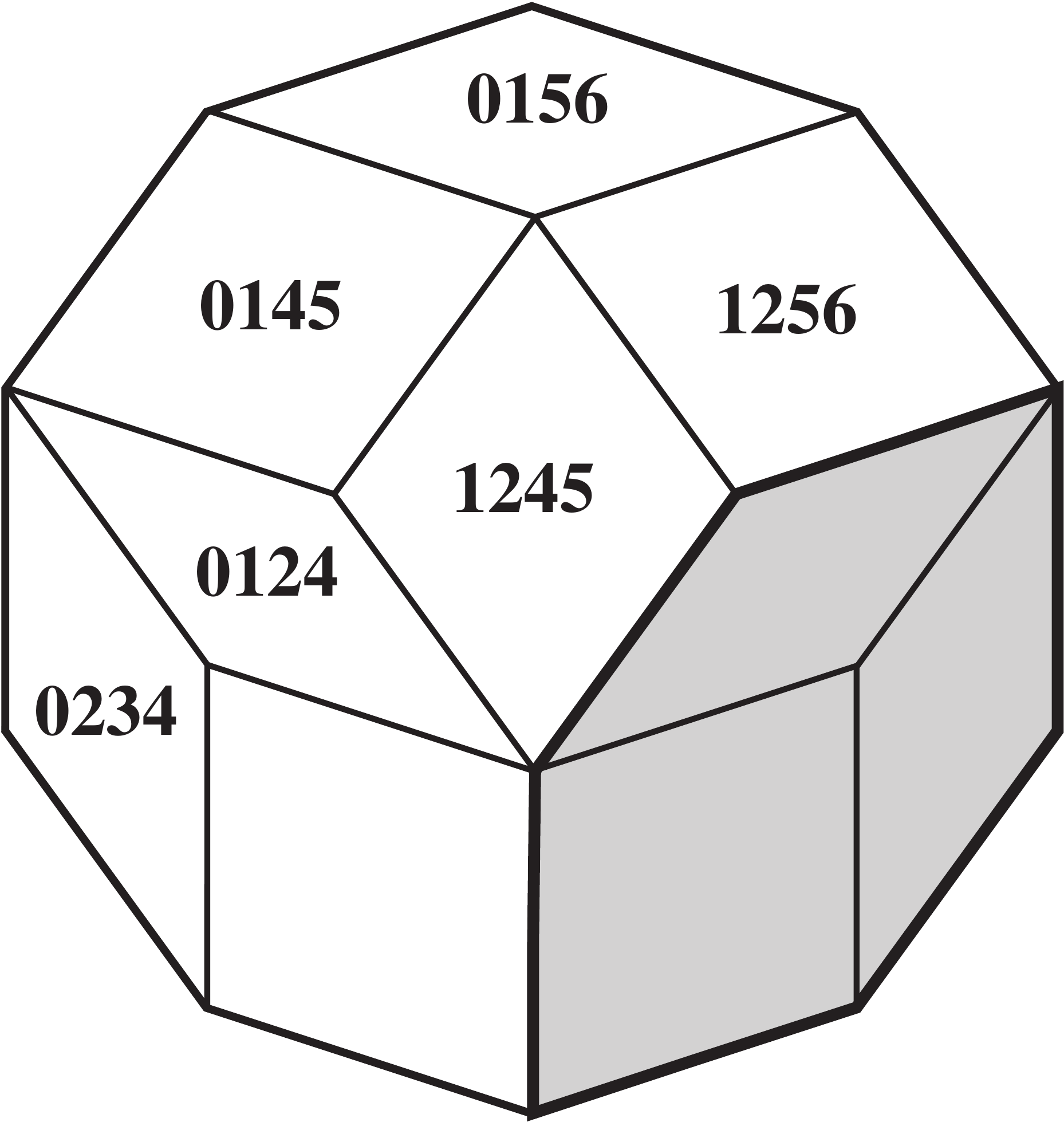}     
&     \includegraphics[width=1in]{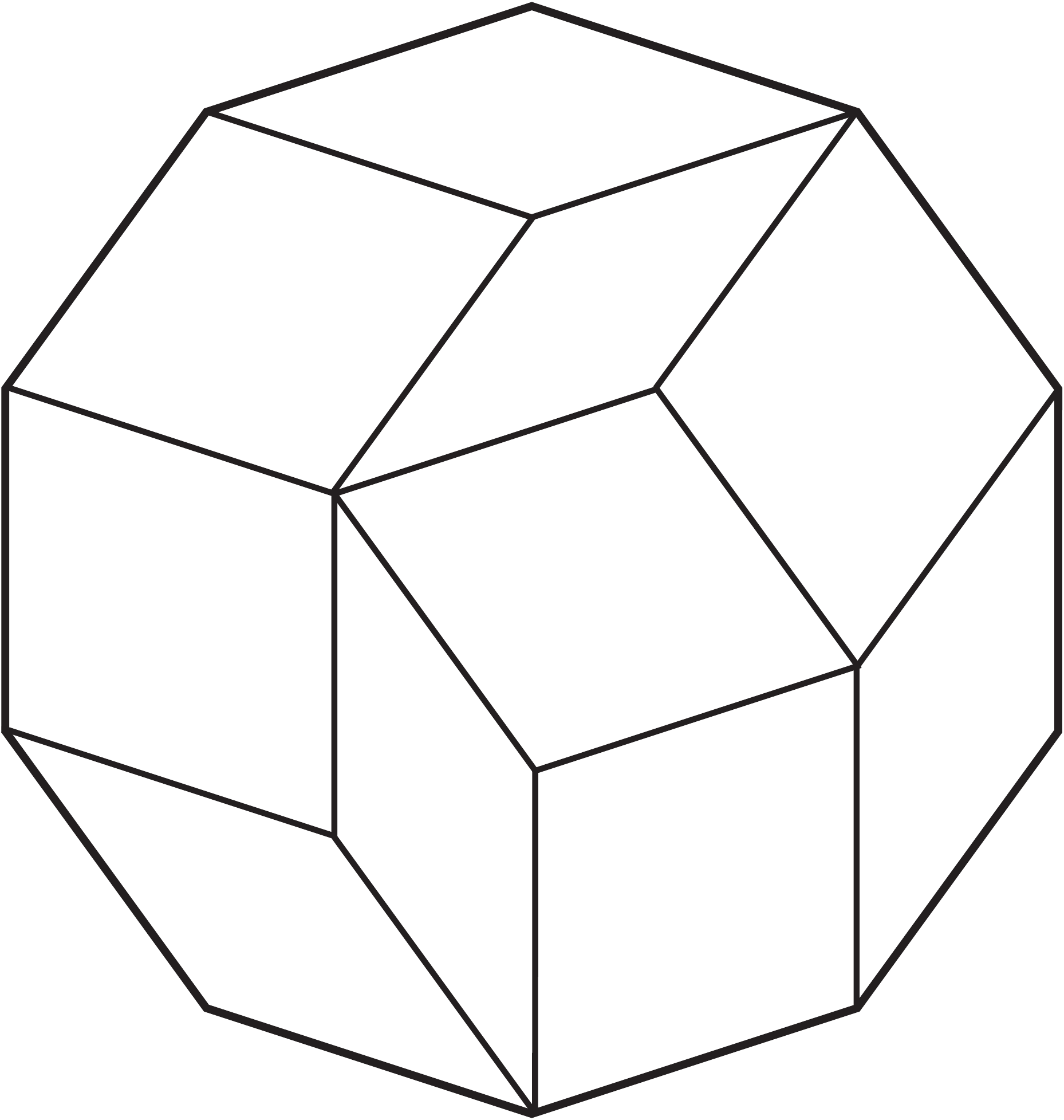}     
&     \includegraphics[width=1in]{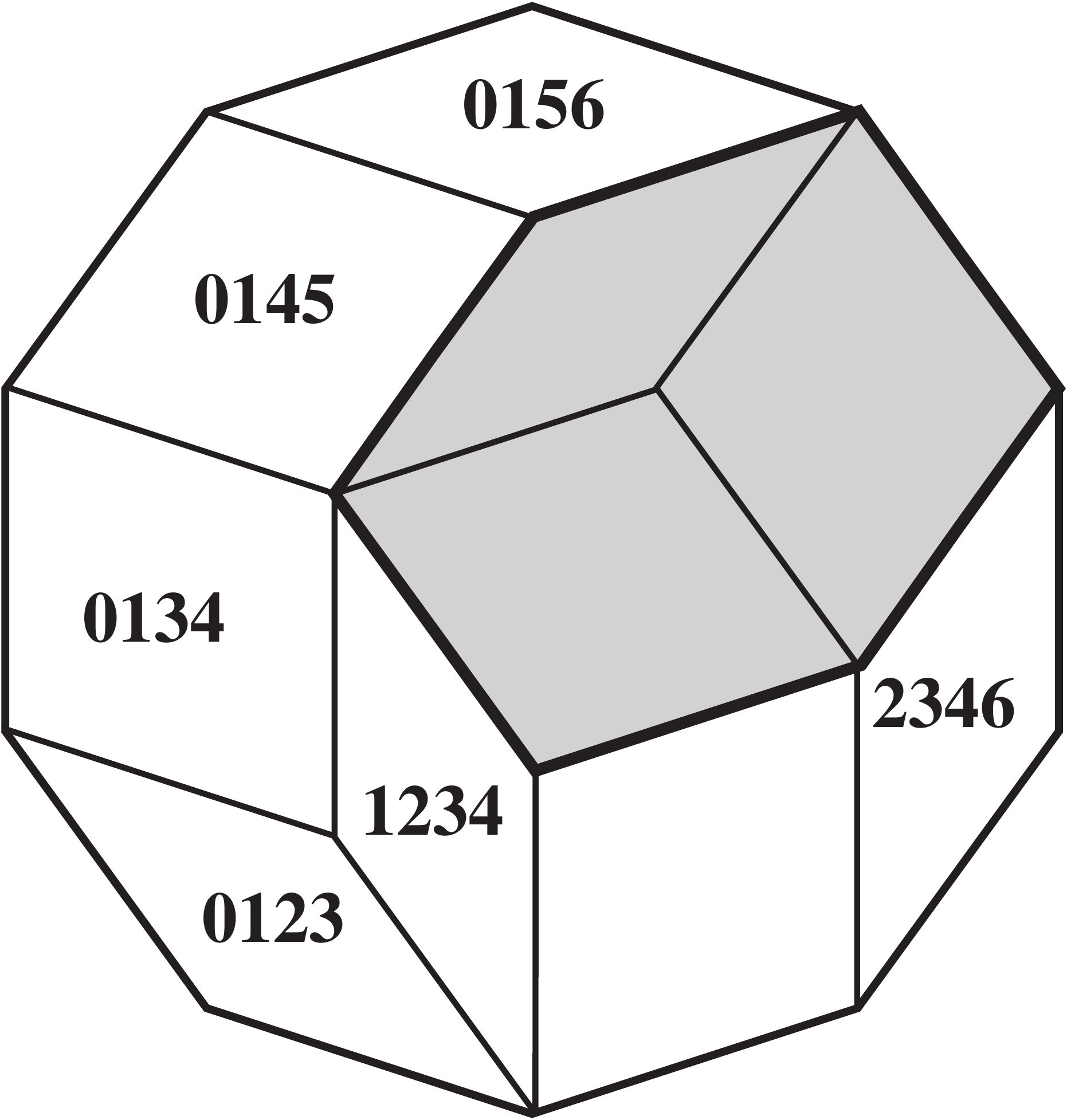} \\

  01245 &    & 23456  &     &  12456  \\

      \includegraphics[width=1in]{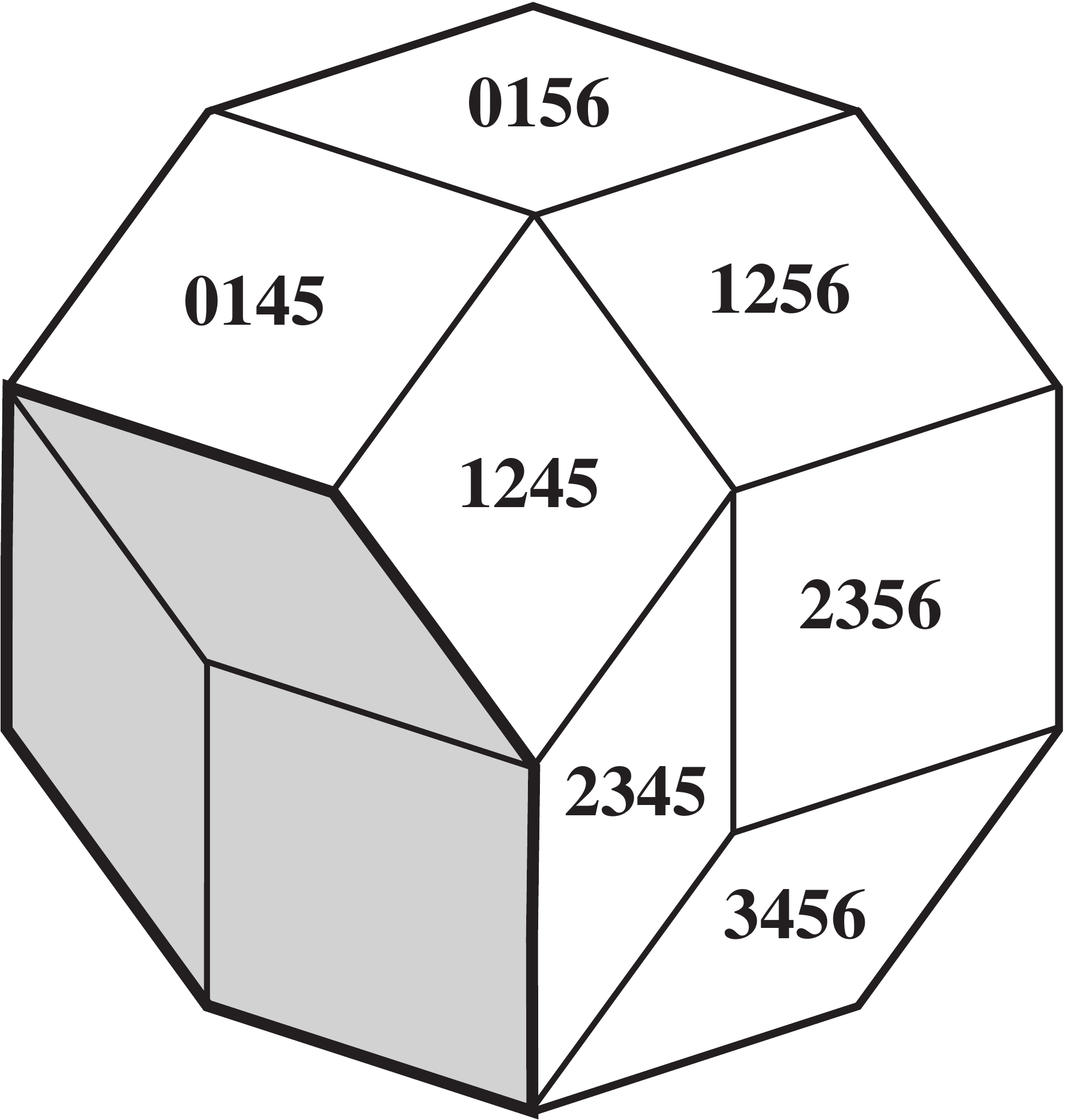}     
  &      \includegraphics[width=1in]{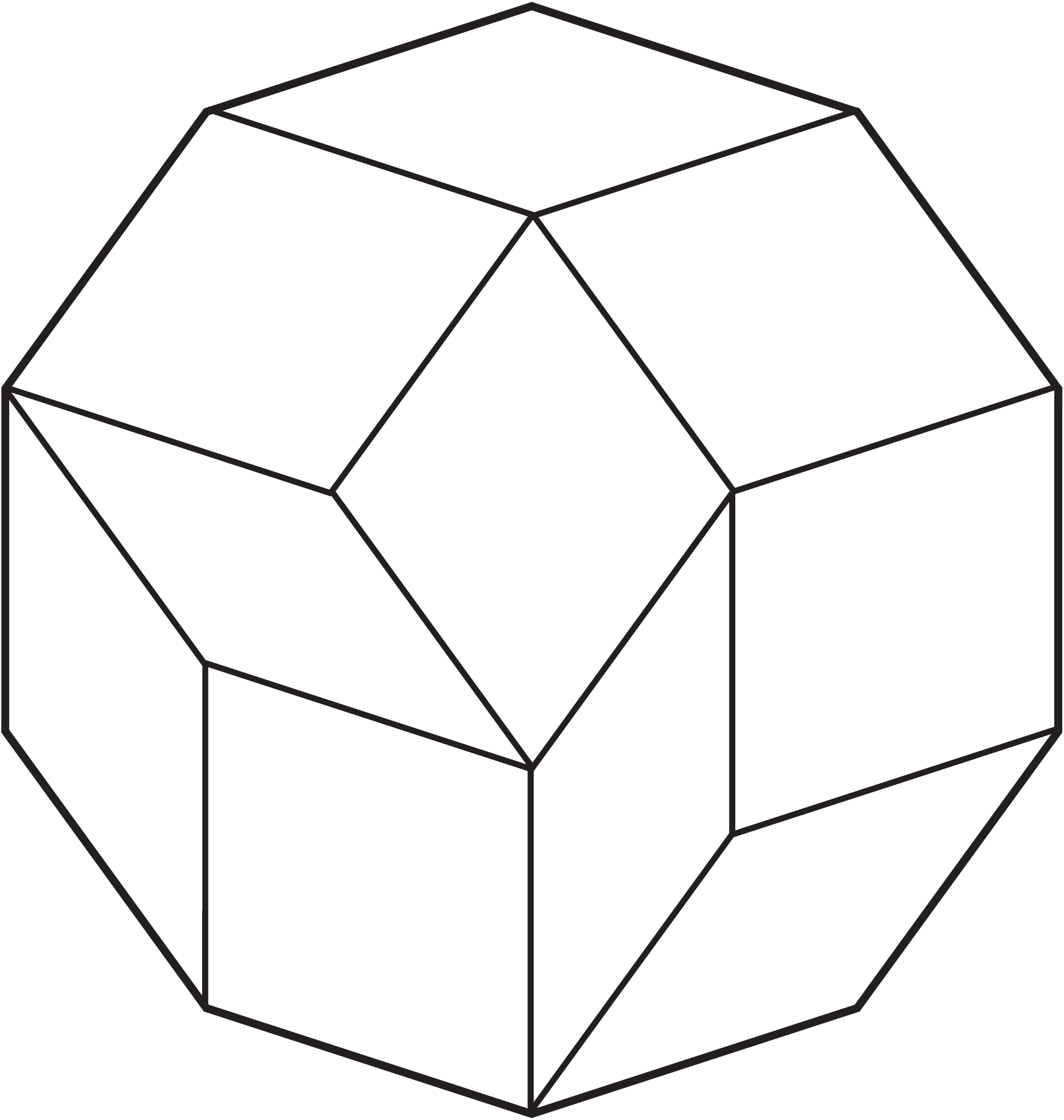}     
 &    \includegraphics[width=1in]{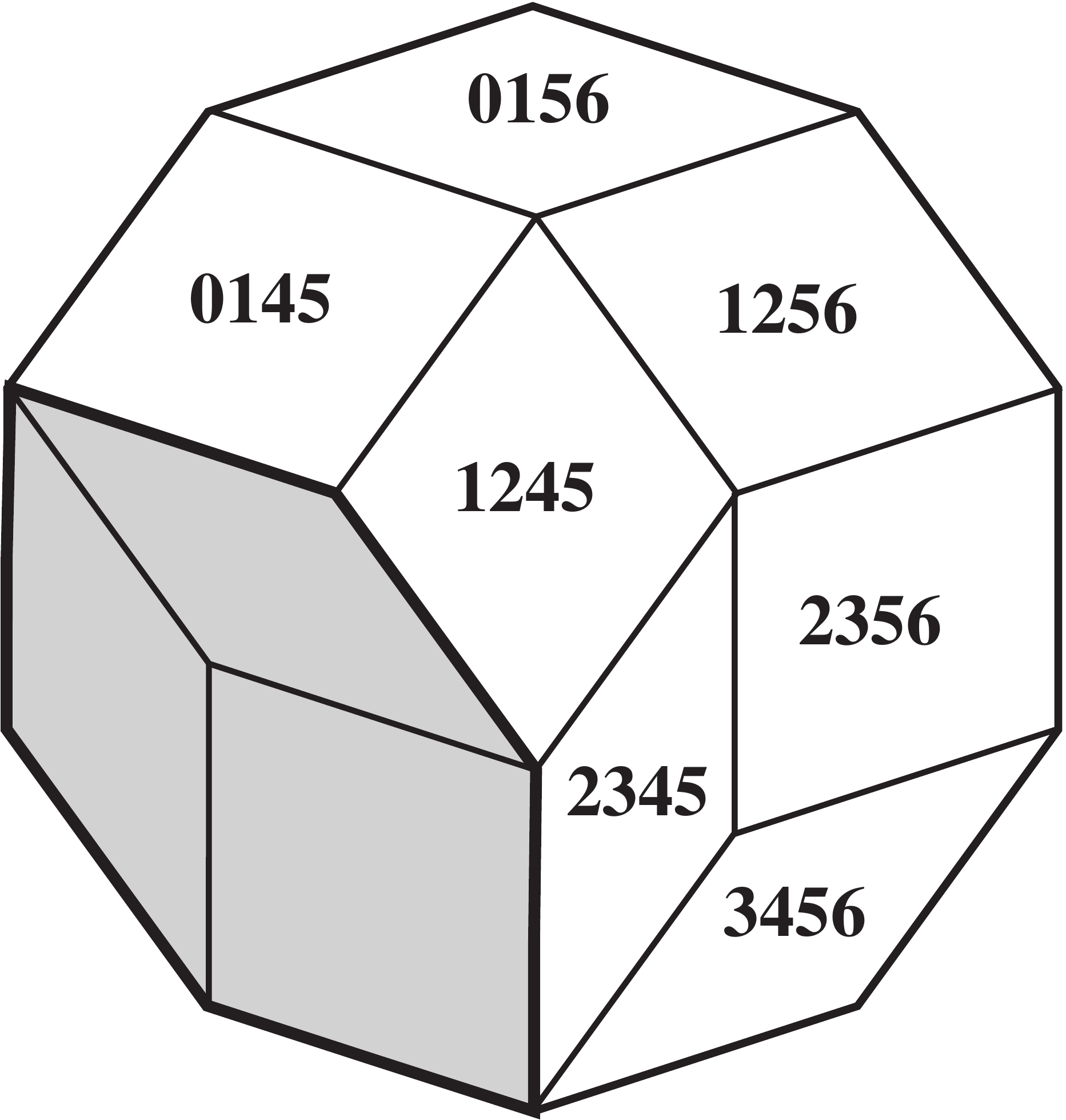}     
&     \includegraphics[width=1in]{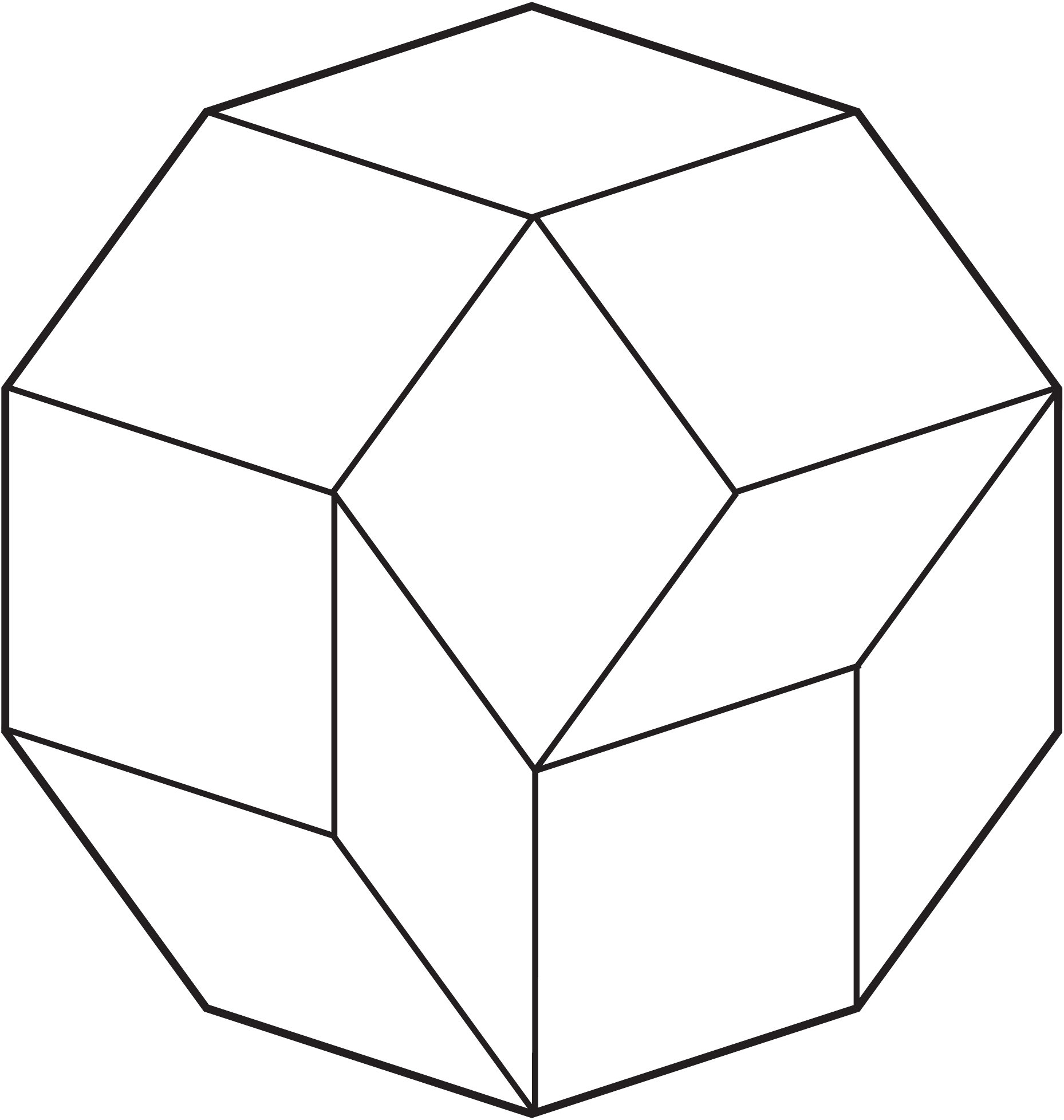}     
&     \includegraphics[width=1in]{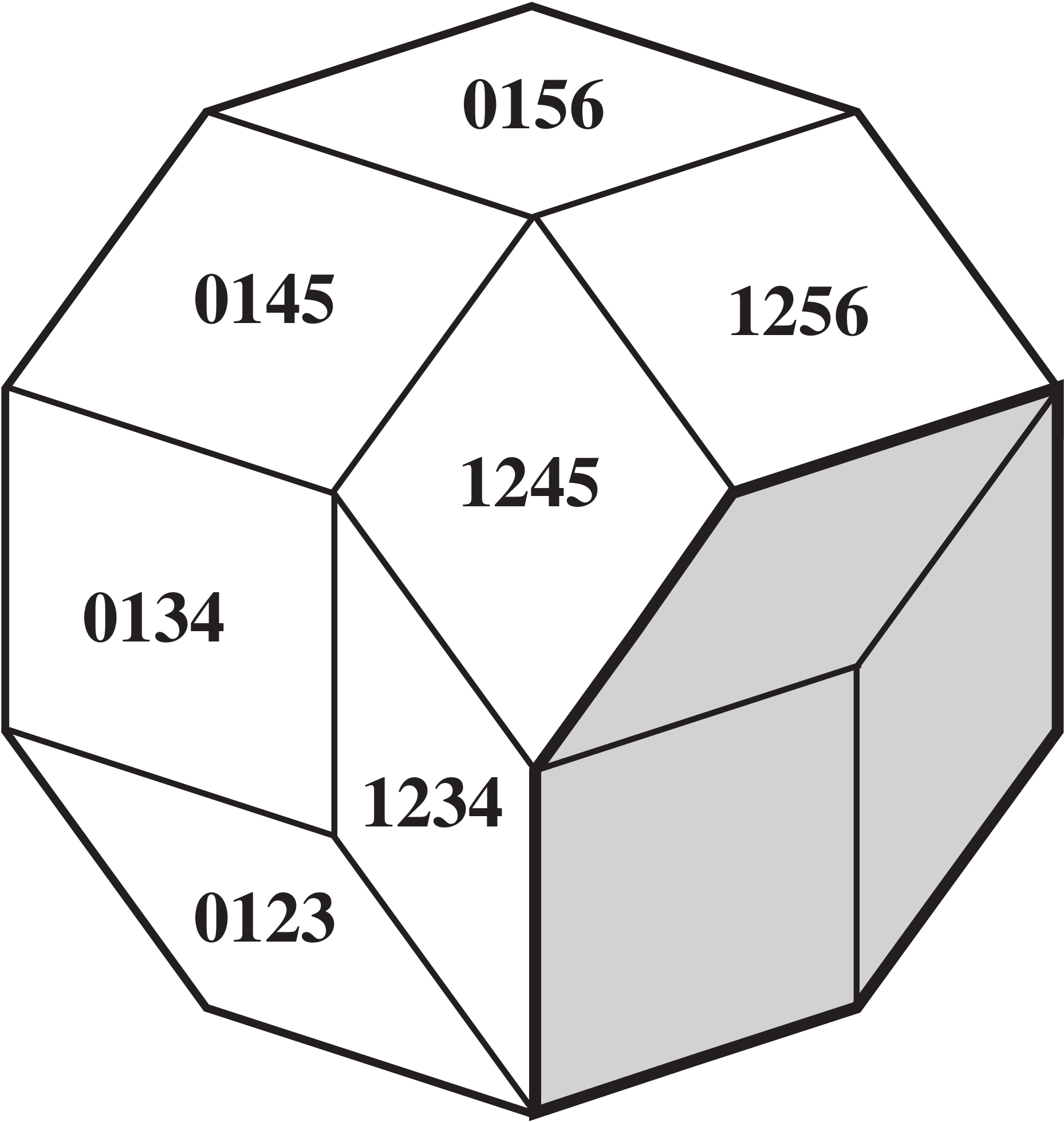} \\

   01234&    & 01234  &     &  23456  \\
 %6-0  &  6-1  &  6-2  &  6-3  &  6-4
  \label{fig:19}
%\end{figure}
 \end{array}
%\right)
\]

  \[
%\left(
\begin{array}{ccccc}  
      \includegraphics[width=1in]{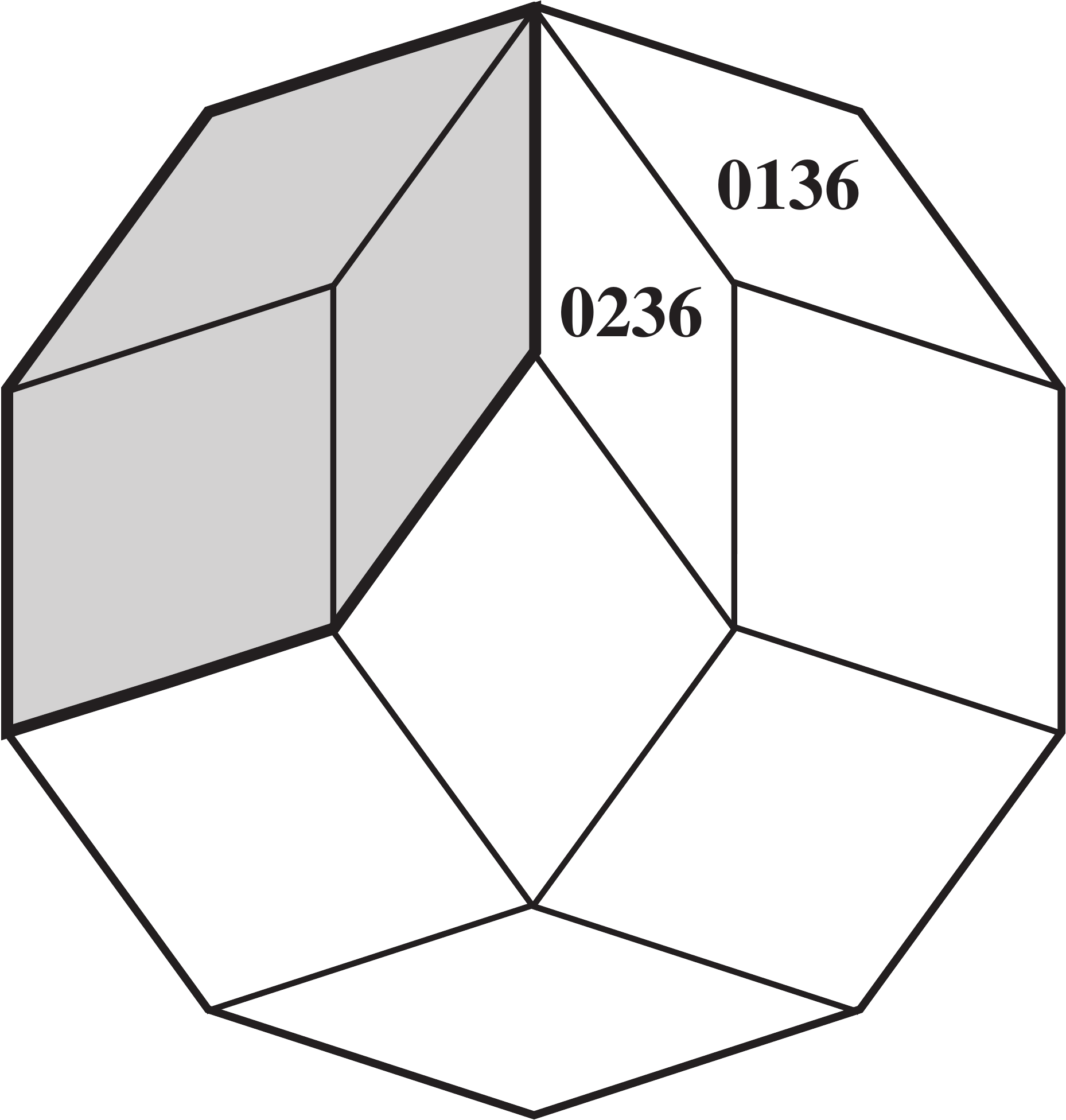}     
  &      \includegraphics[width=1in]{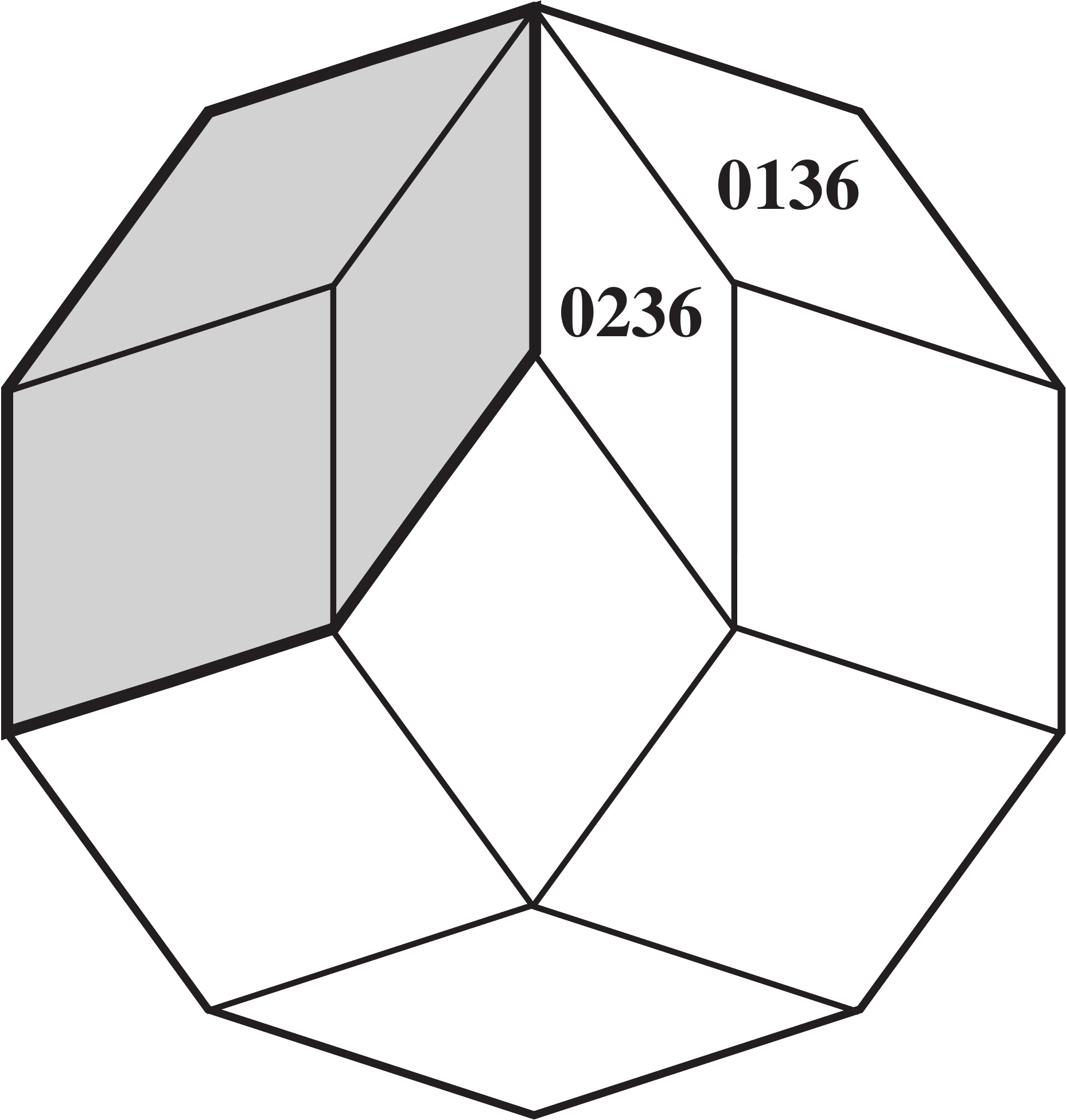}     
 &    \includegraphics[width=1in]{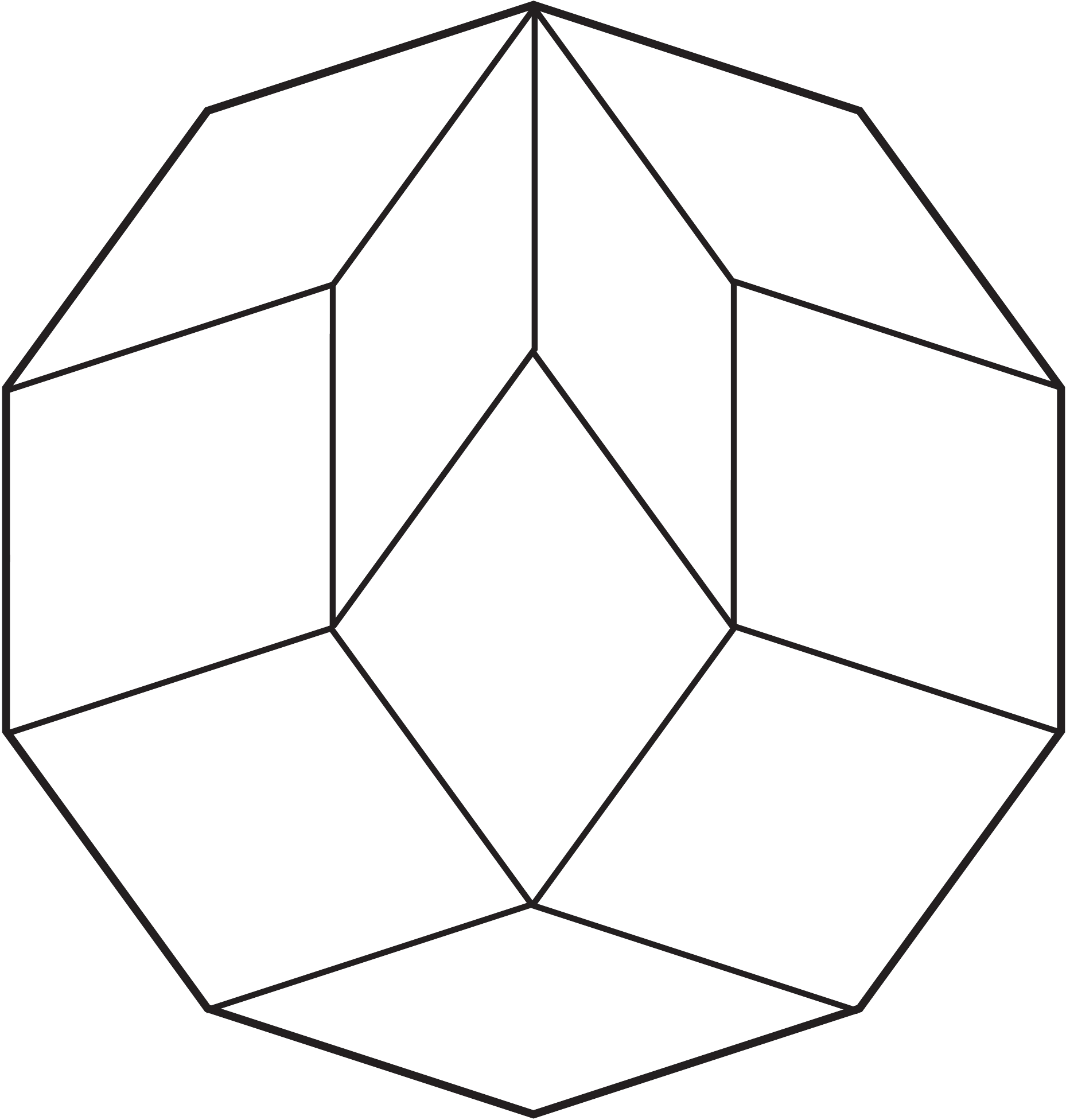}     
&     \includegraphics[width=1in]{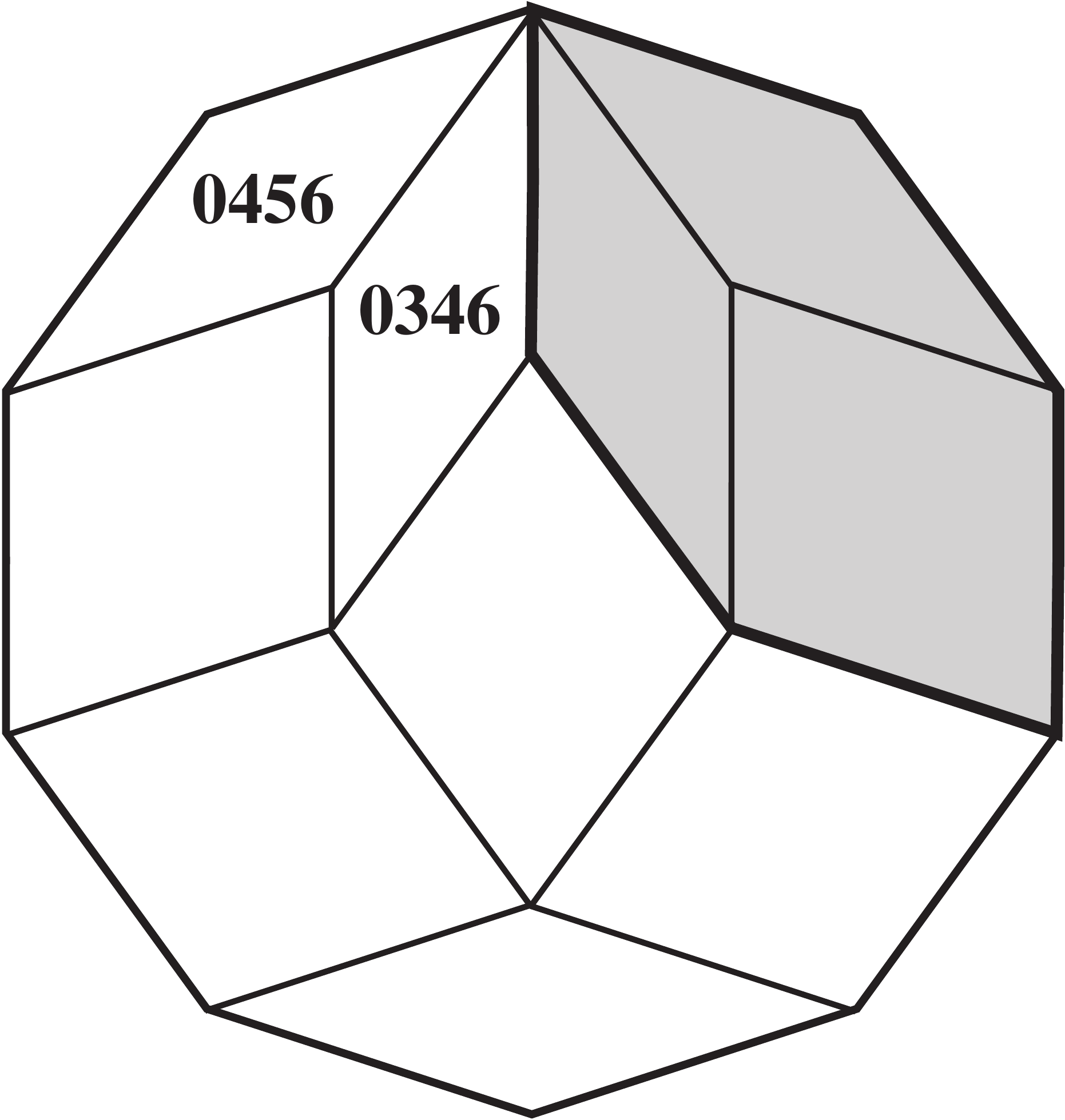}     
&     \includegraphics[width=1in]{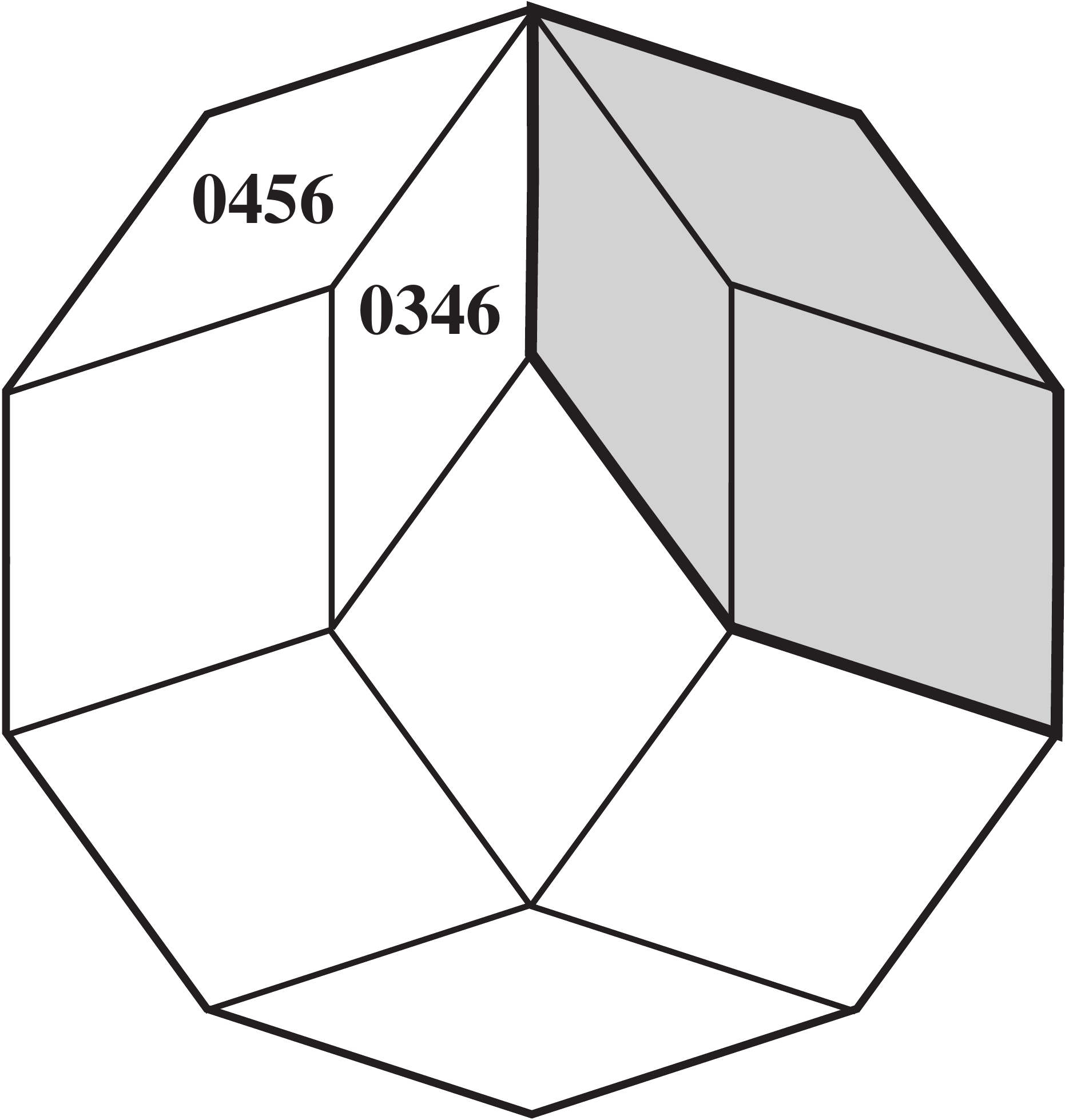} \\

03456  & 03456   &   & 01236    &  01236  \\

     \includegraphics[width=1in]{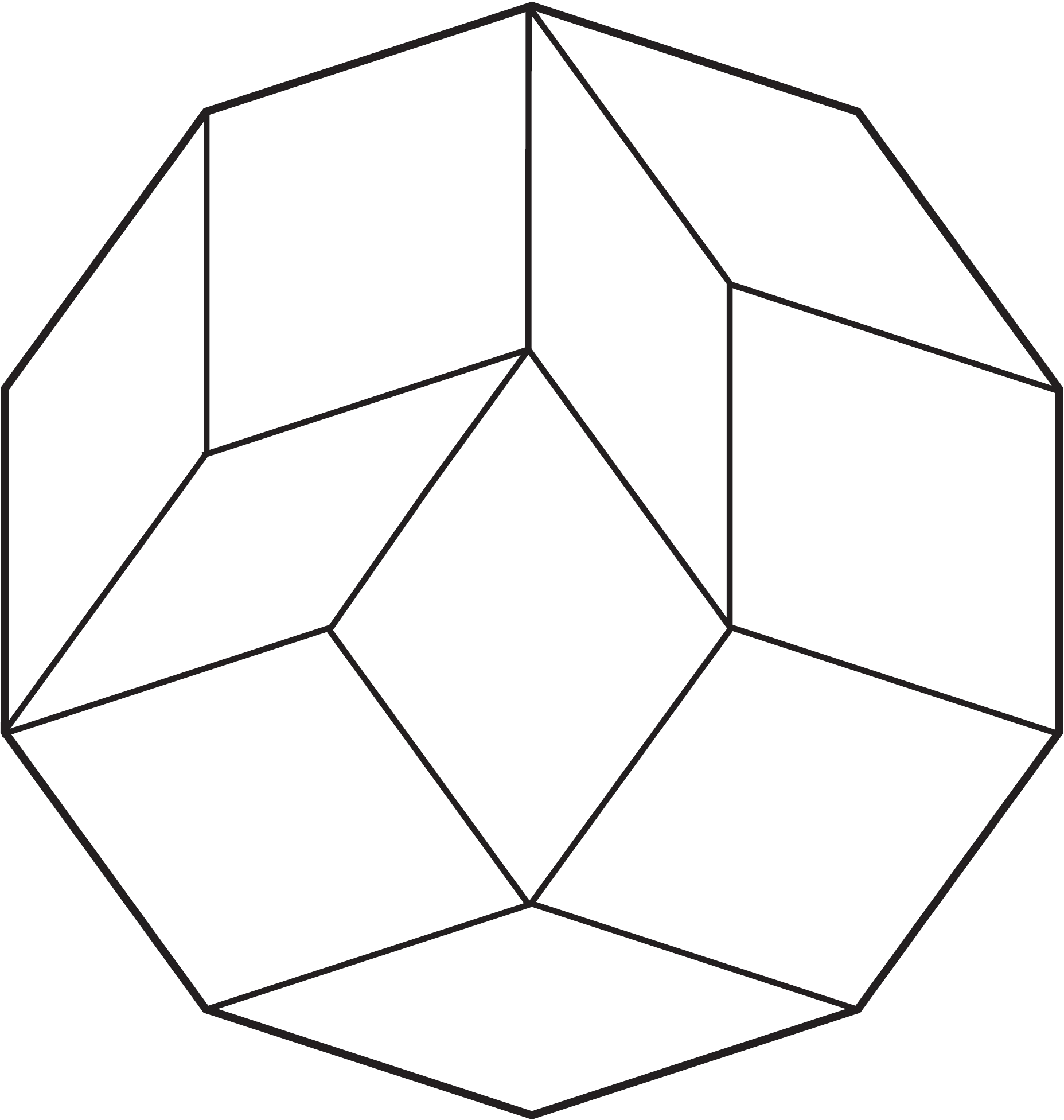}     
  &      \includegraphics[width=1in]{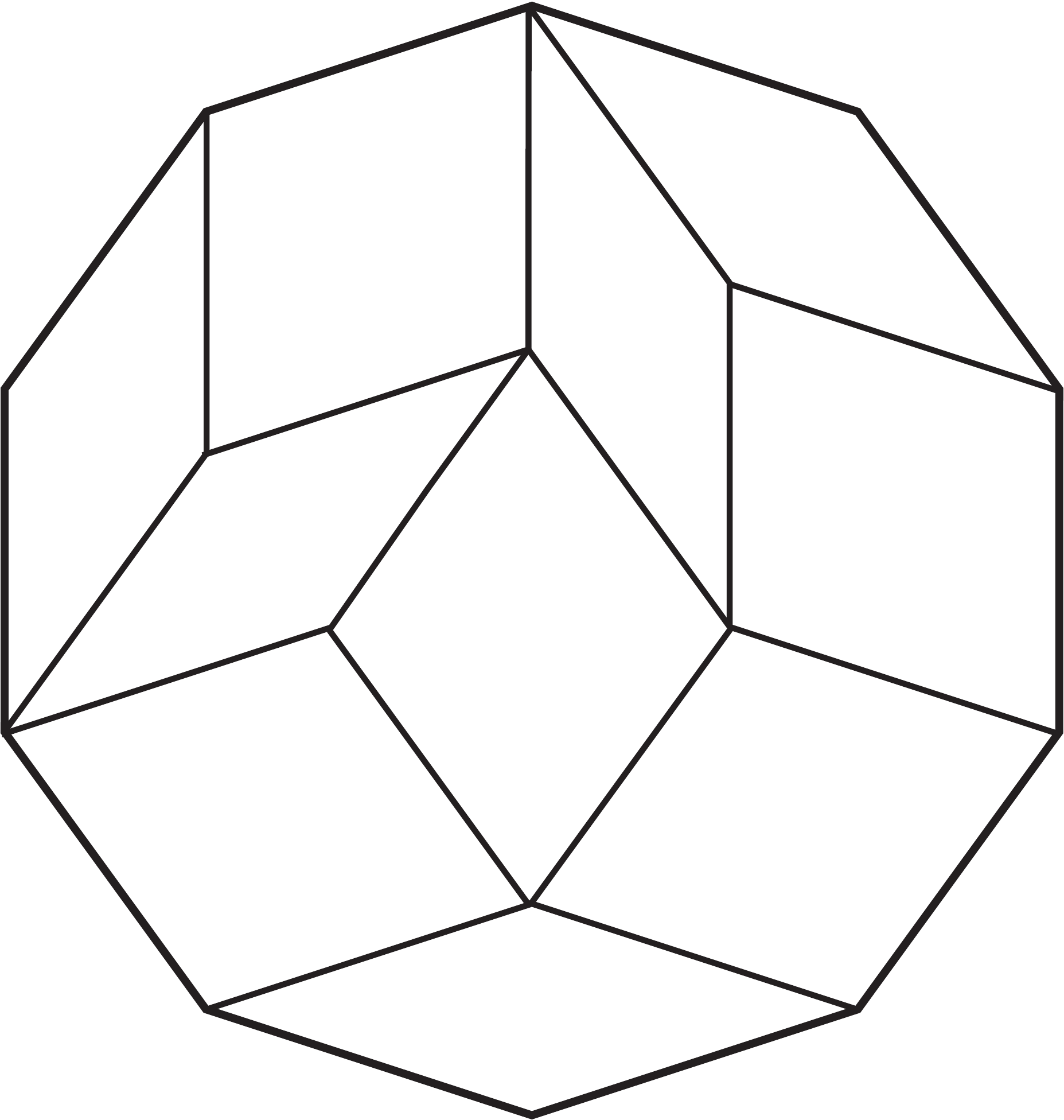}     
 &    \includegraphics[width=1in]{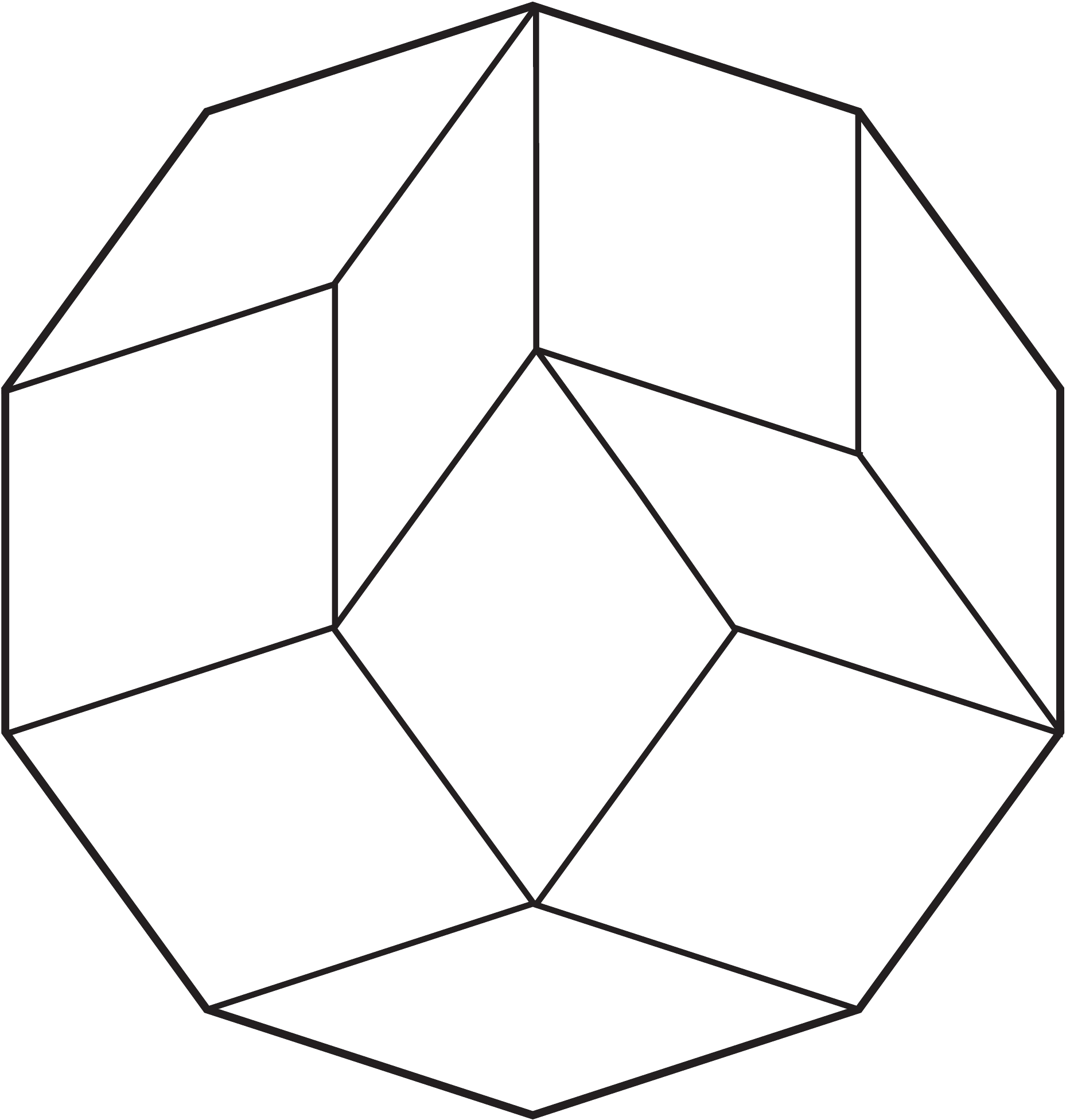}     
&     \includegraphics[width=1in]{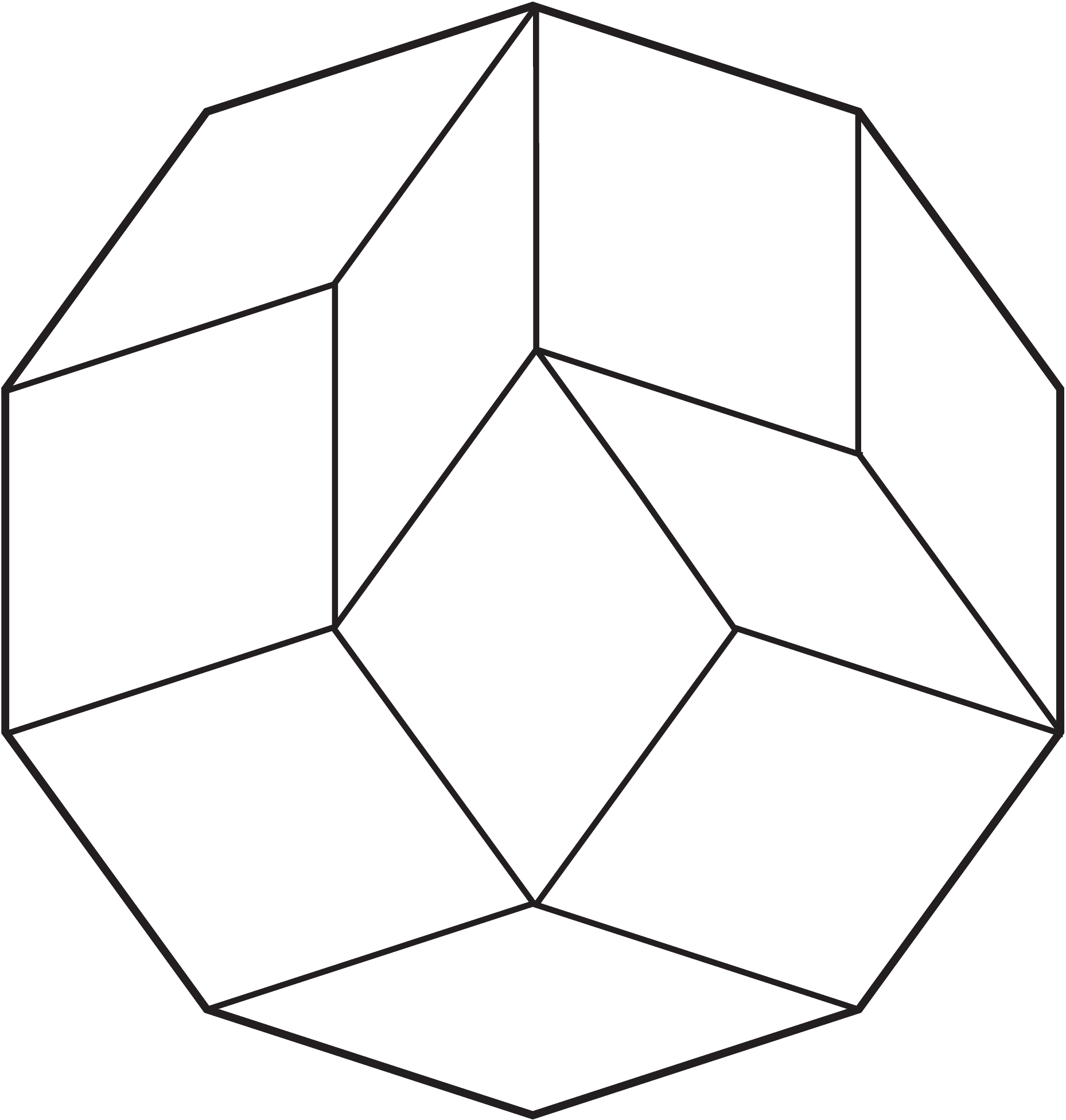}     
&     \includegraphics[width=1in]{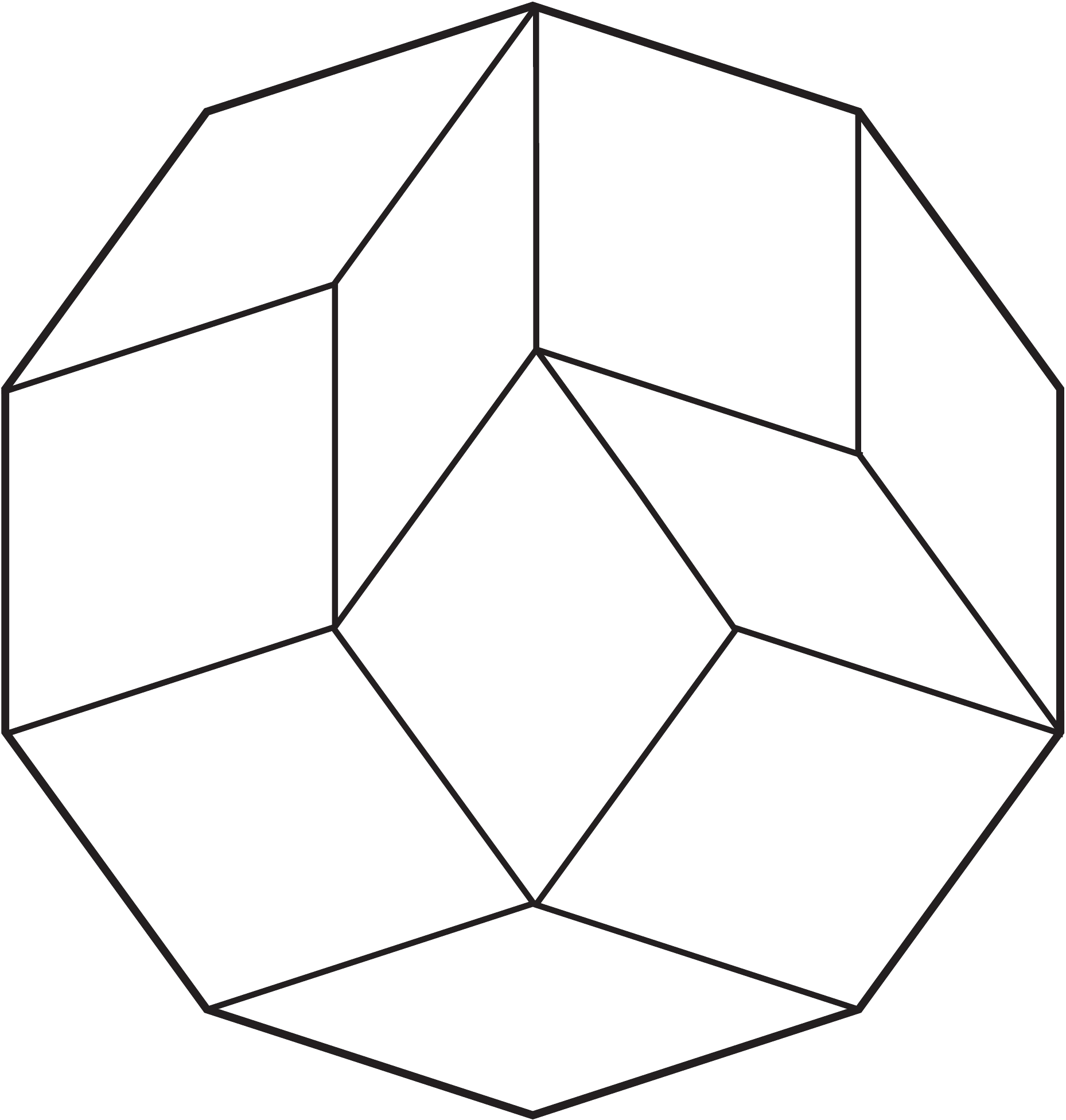} \\

  &    &   &     &    \\

      \includegraphics[width=1in]{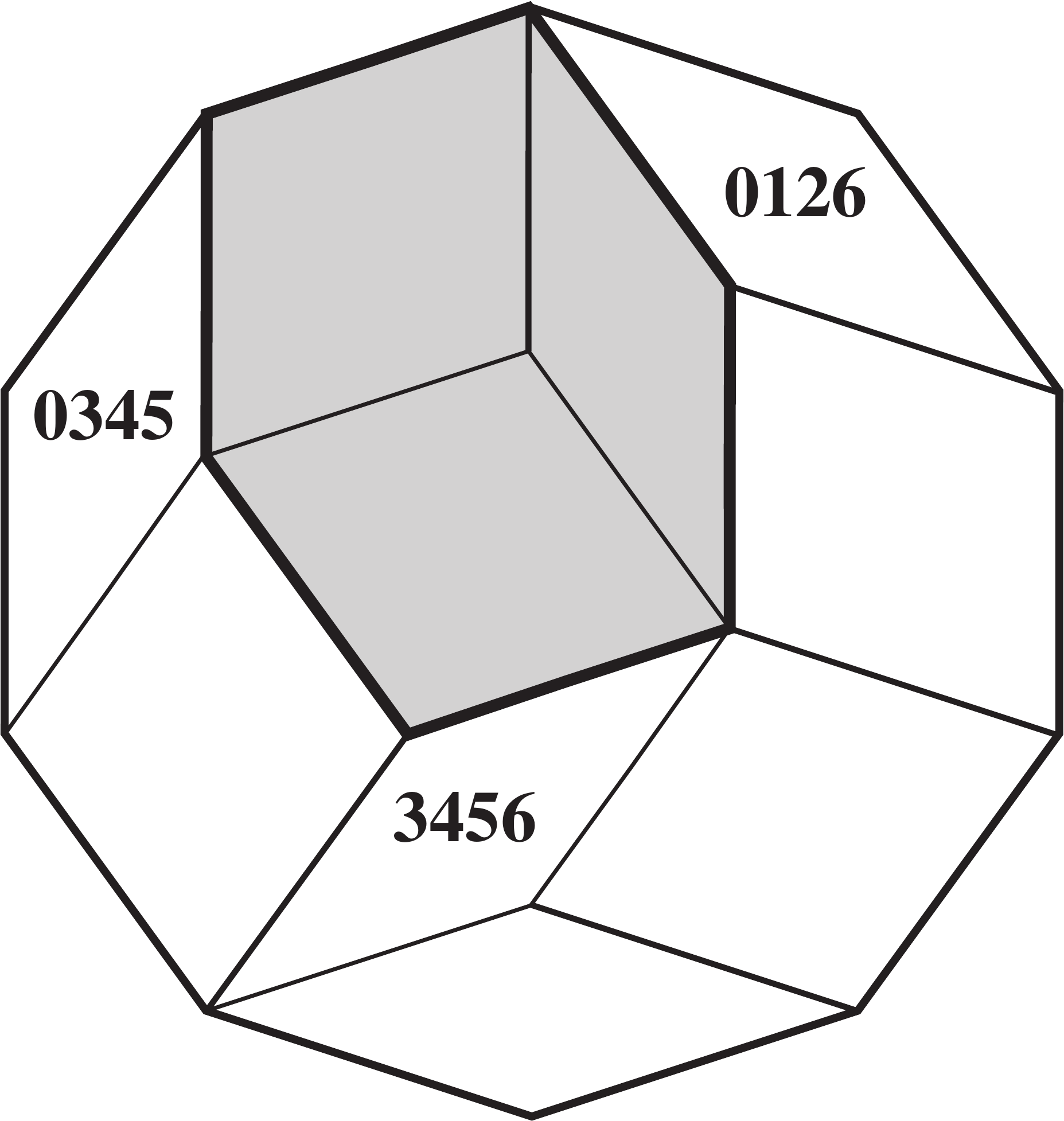}     
  &      \includegraphics[width=1in]{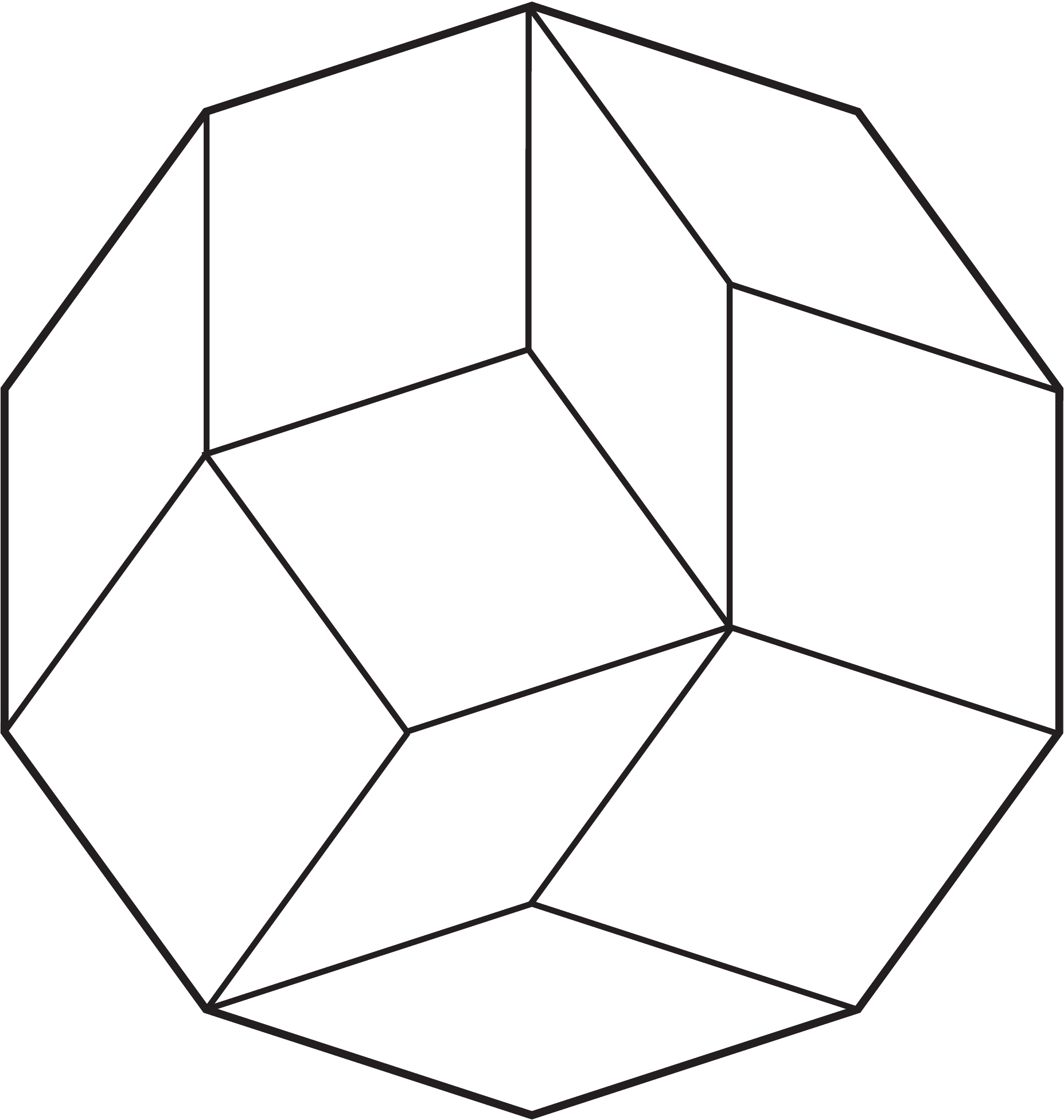}     
 &    \includegraphics[width=1in]{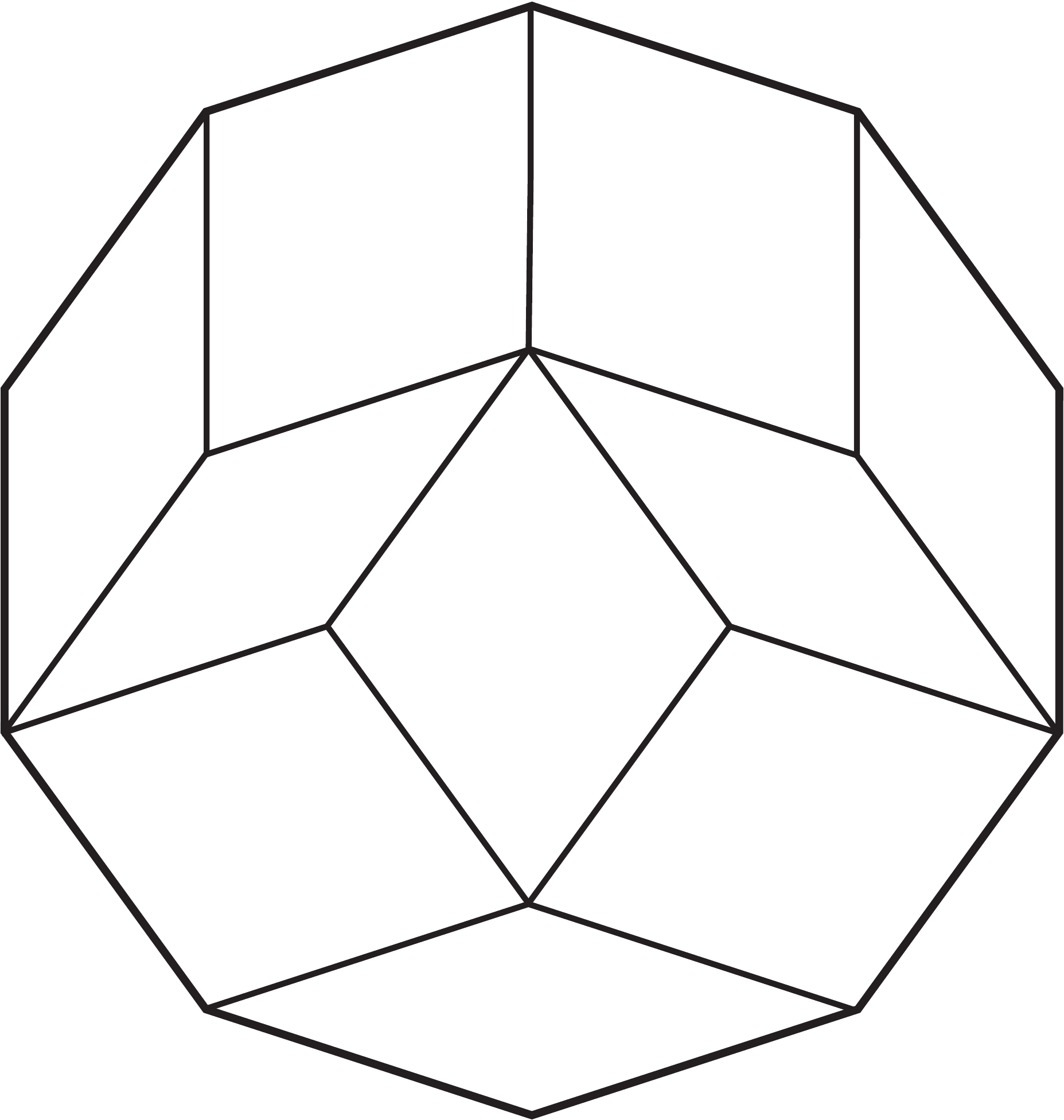}     
&     \includegraphics[width=1in]{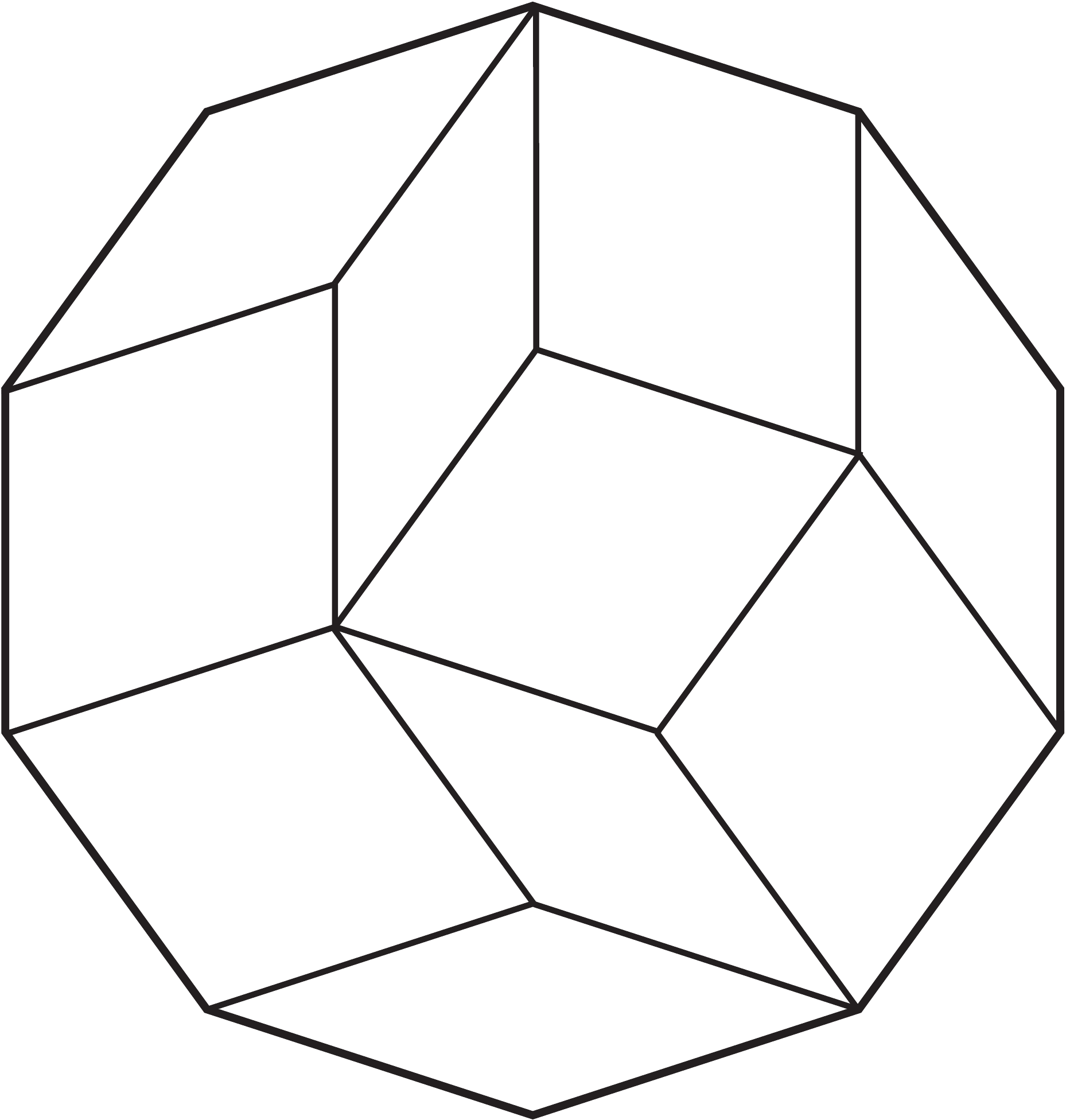}     
&     \includegraphics[width=1in]{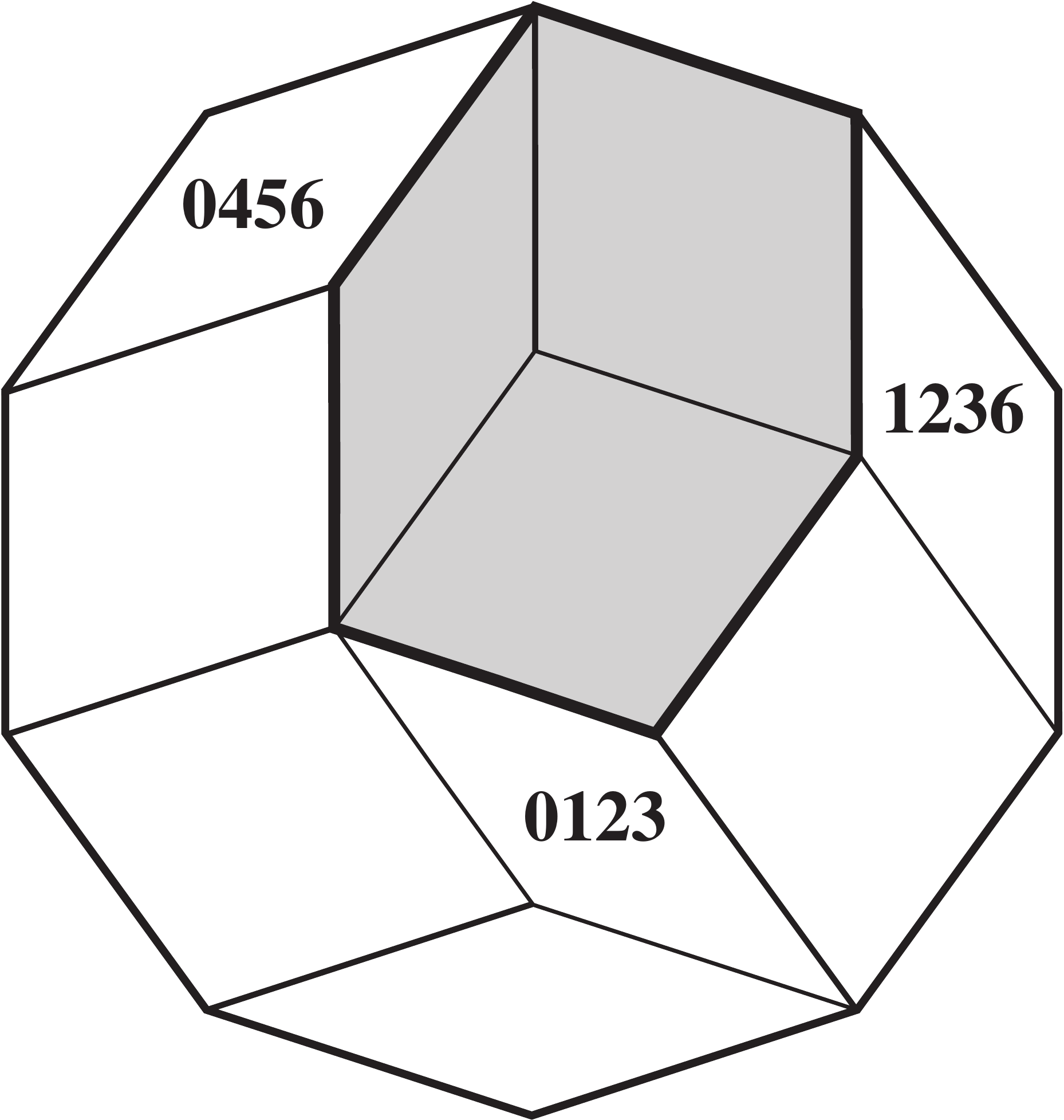} \\

02356   &    &   &     & 01346   \\

      \includegraphics[width=1in]{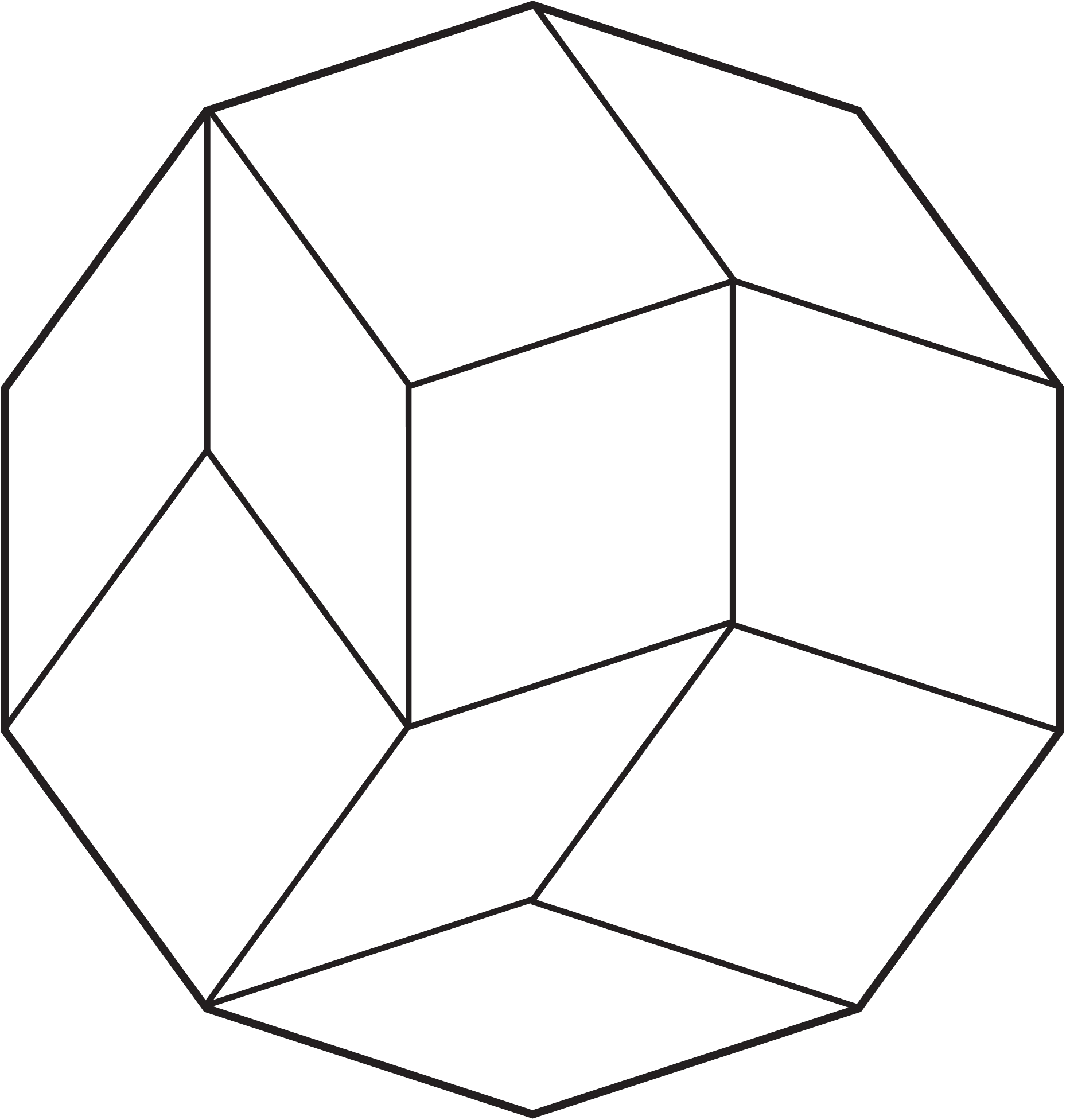}     
  &      \includegraphics[width=1in]{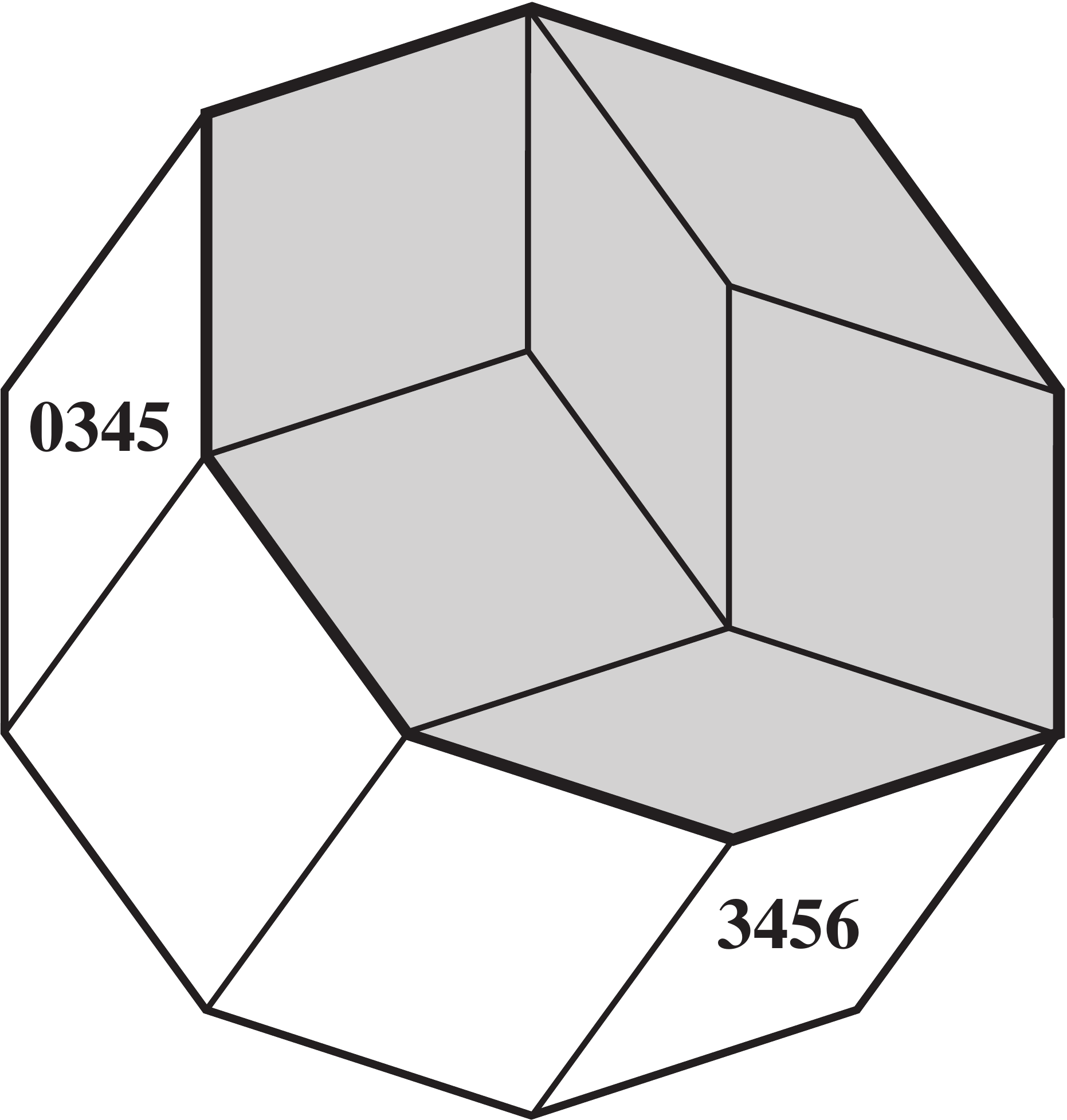}     
 &    \includegraphics[width=1in]{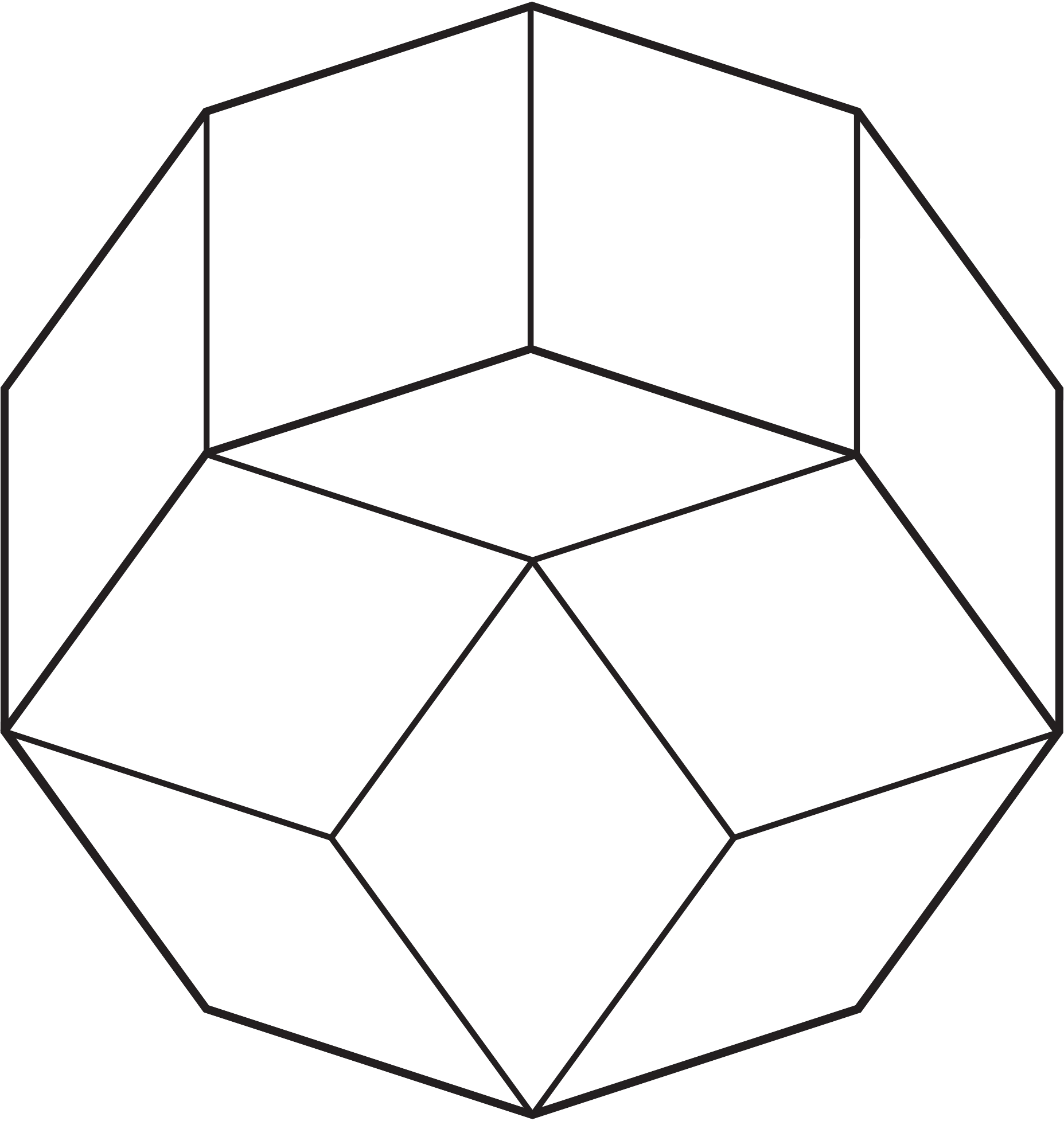}     
&     \includegraphics[width=1in]{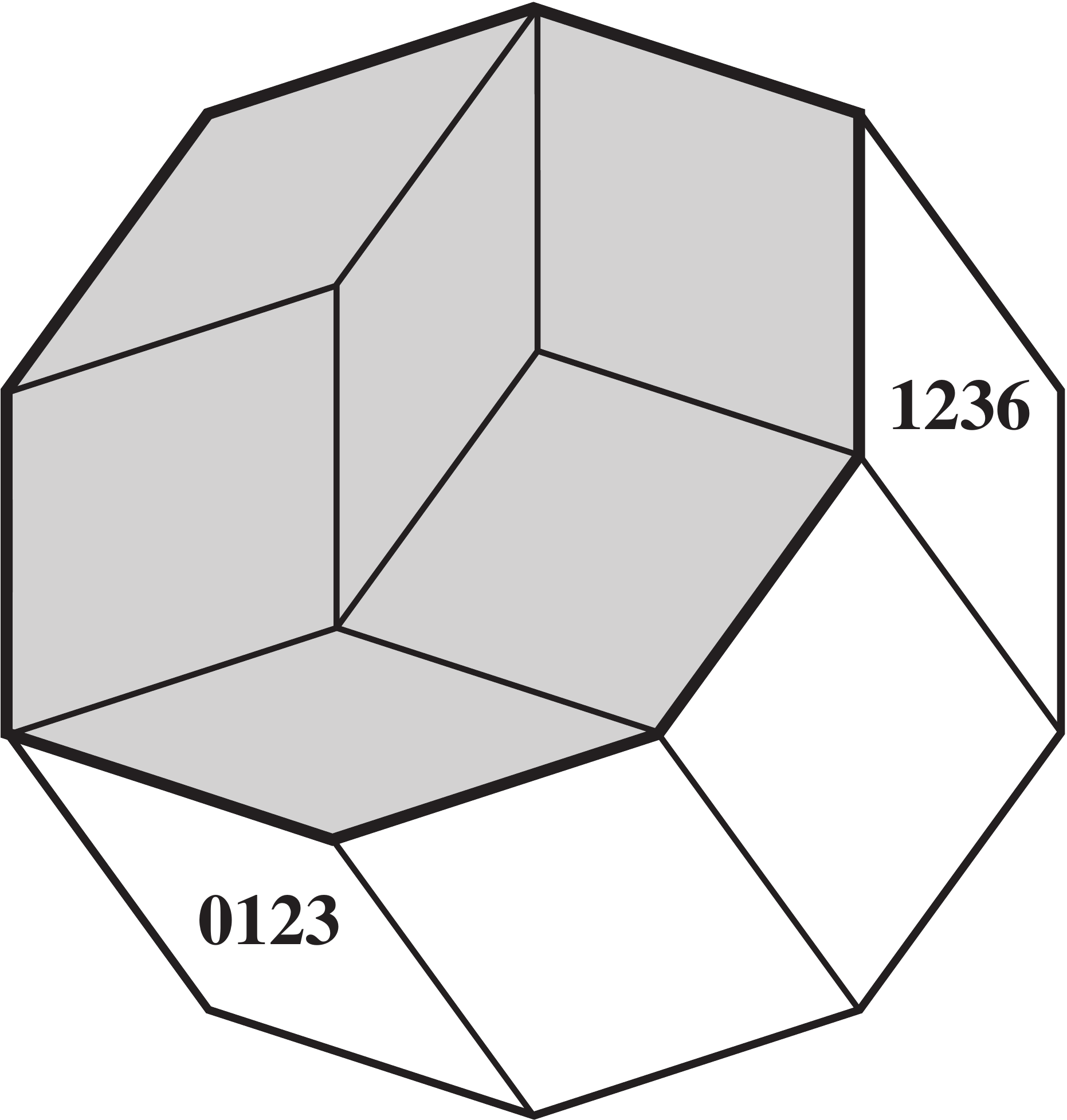}     
&     \includegraphics[width=1in]{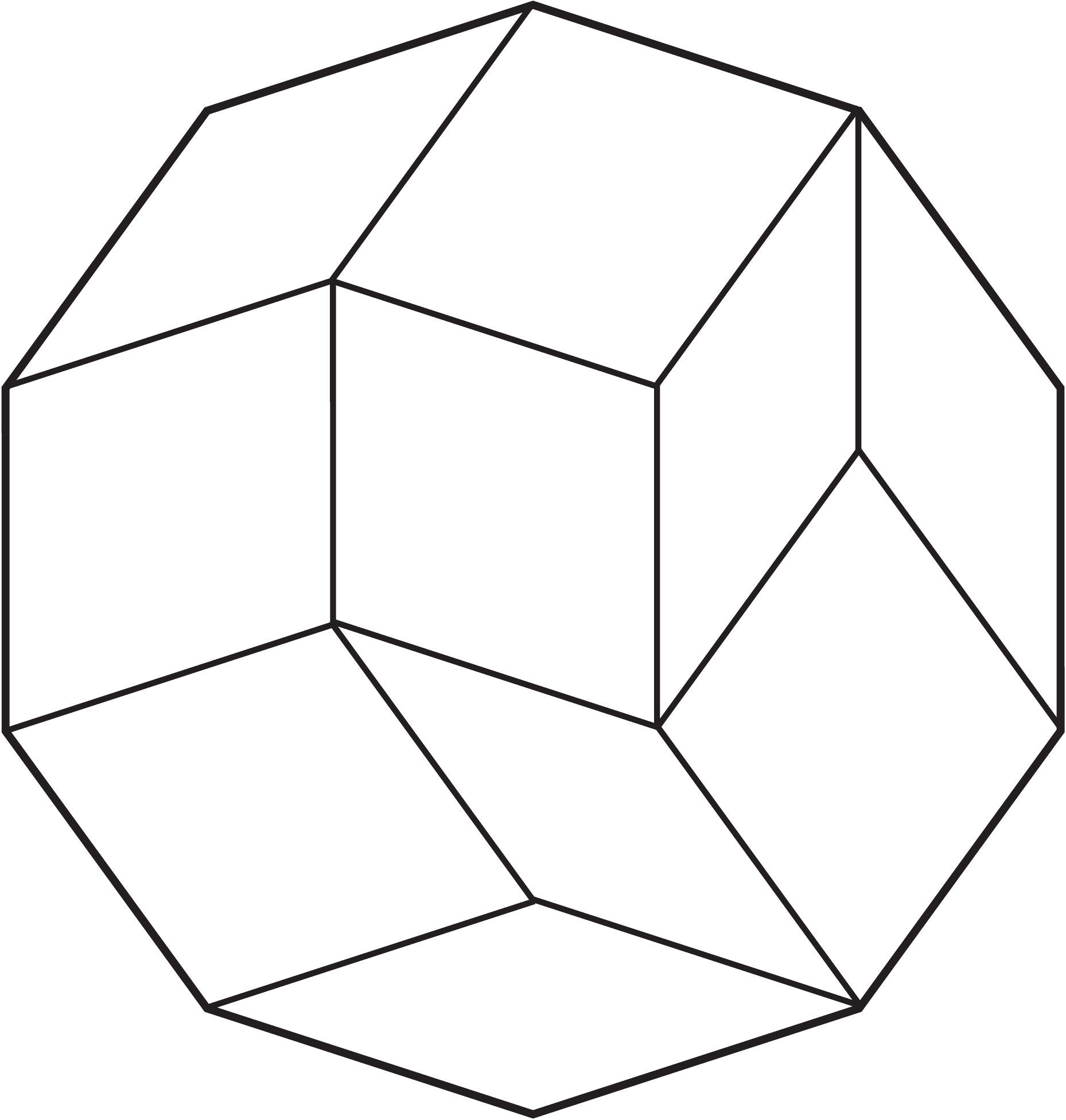} \\

  & 012356   &   &  013456   &    \\

     \includegraphics[width=1in]{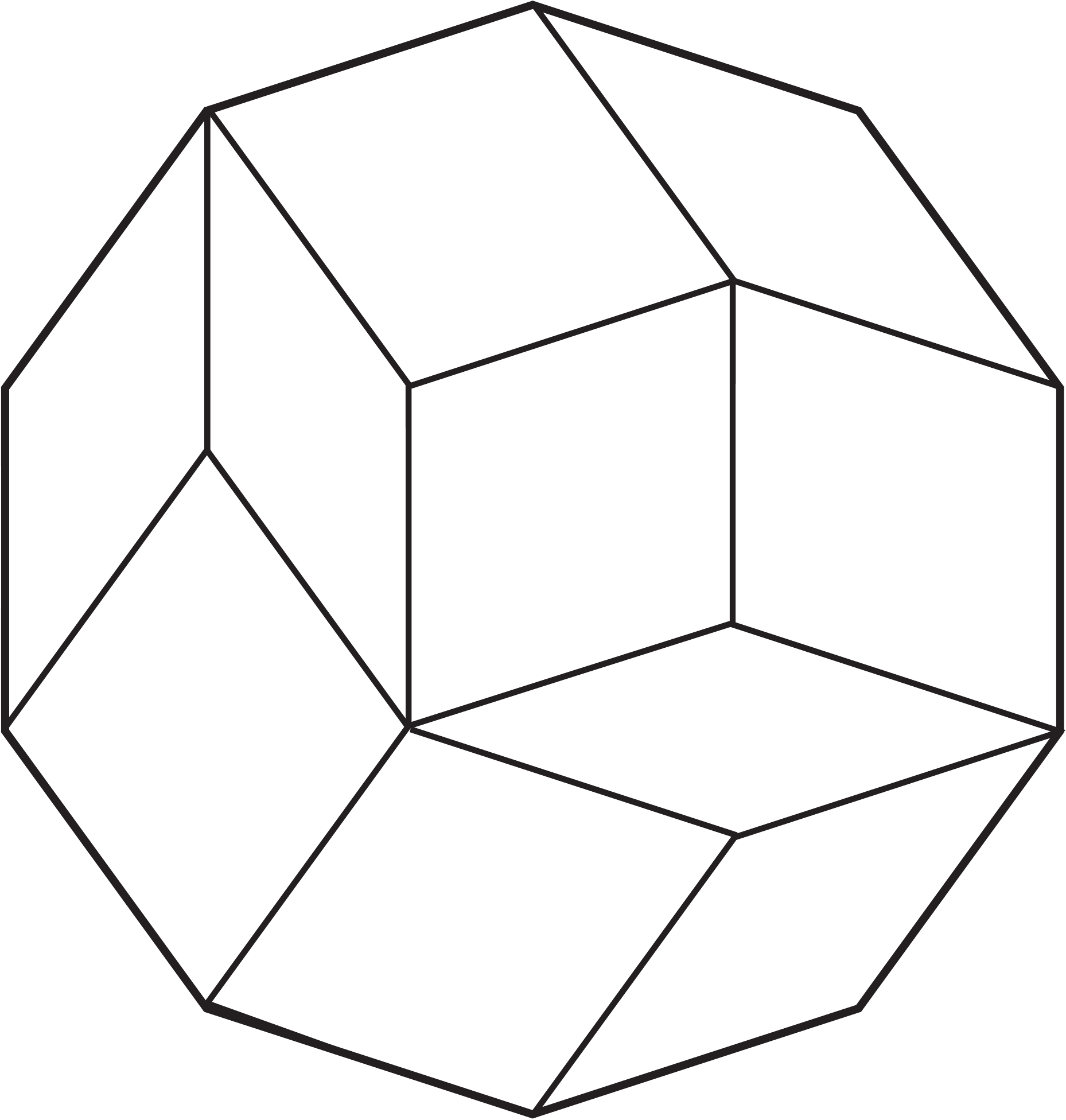}     
  &      \includegraphics[width=1in]{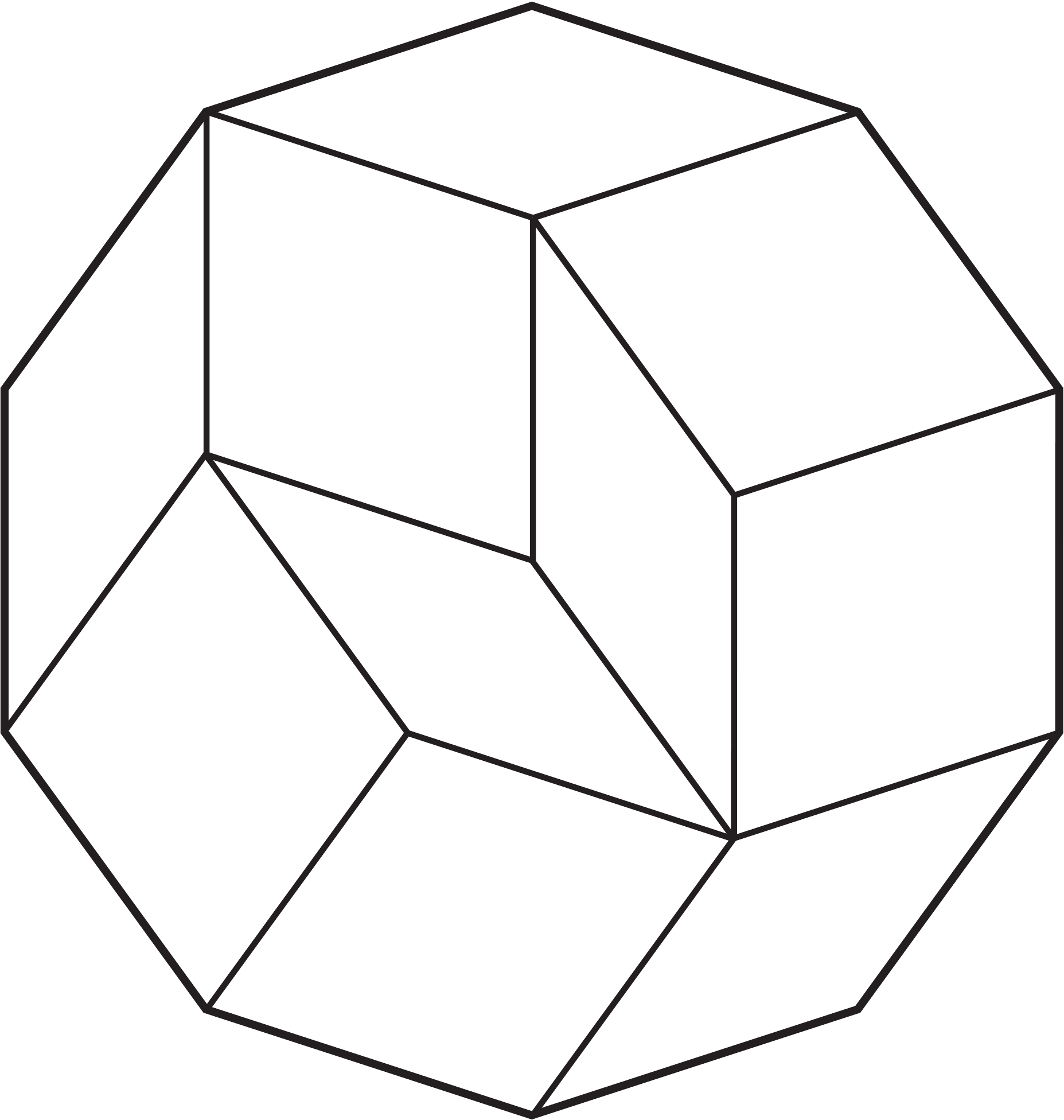}     
 &    \includegraphics[width=1in]{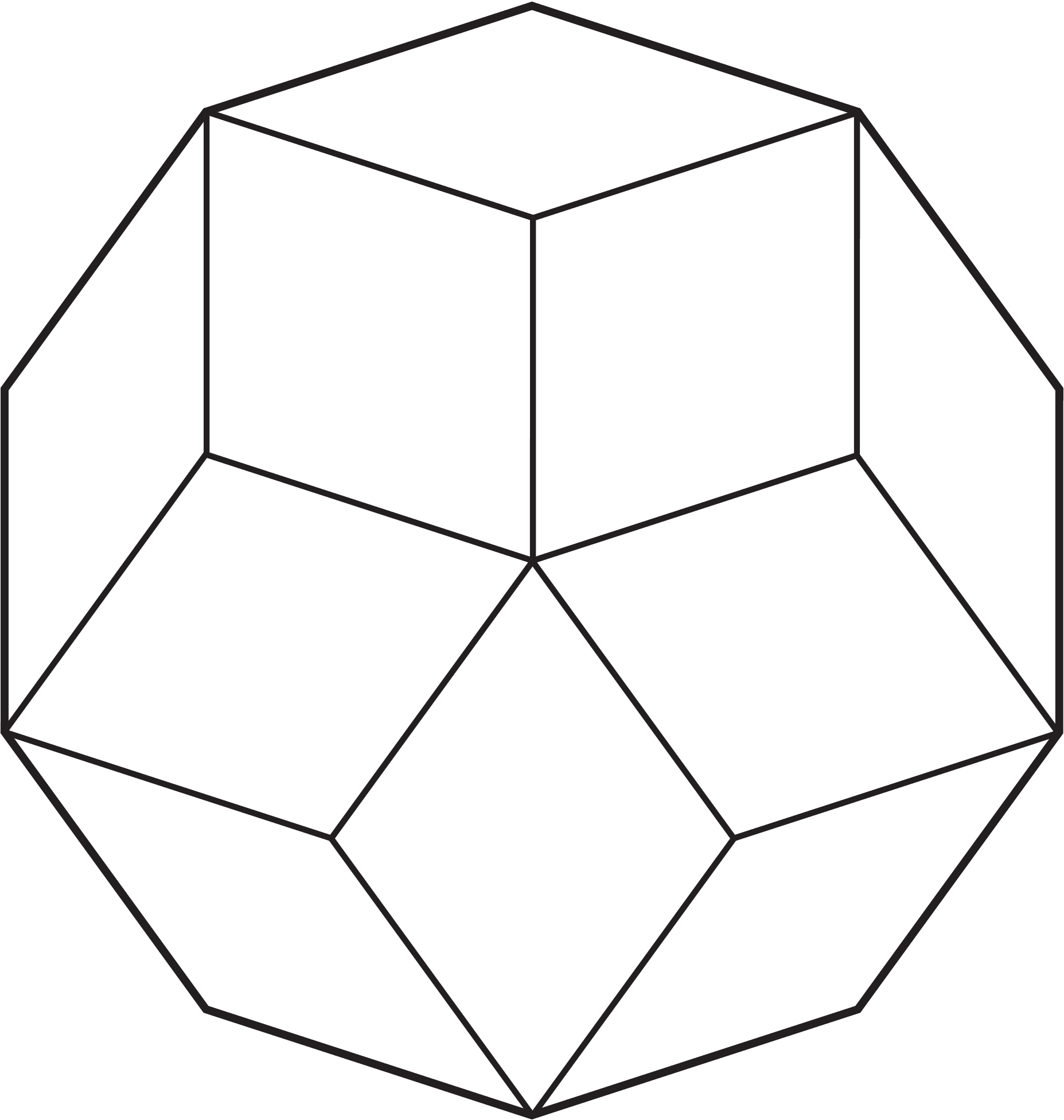}     
&     \includegraphics[width=1in]{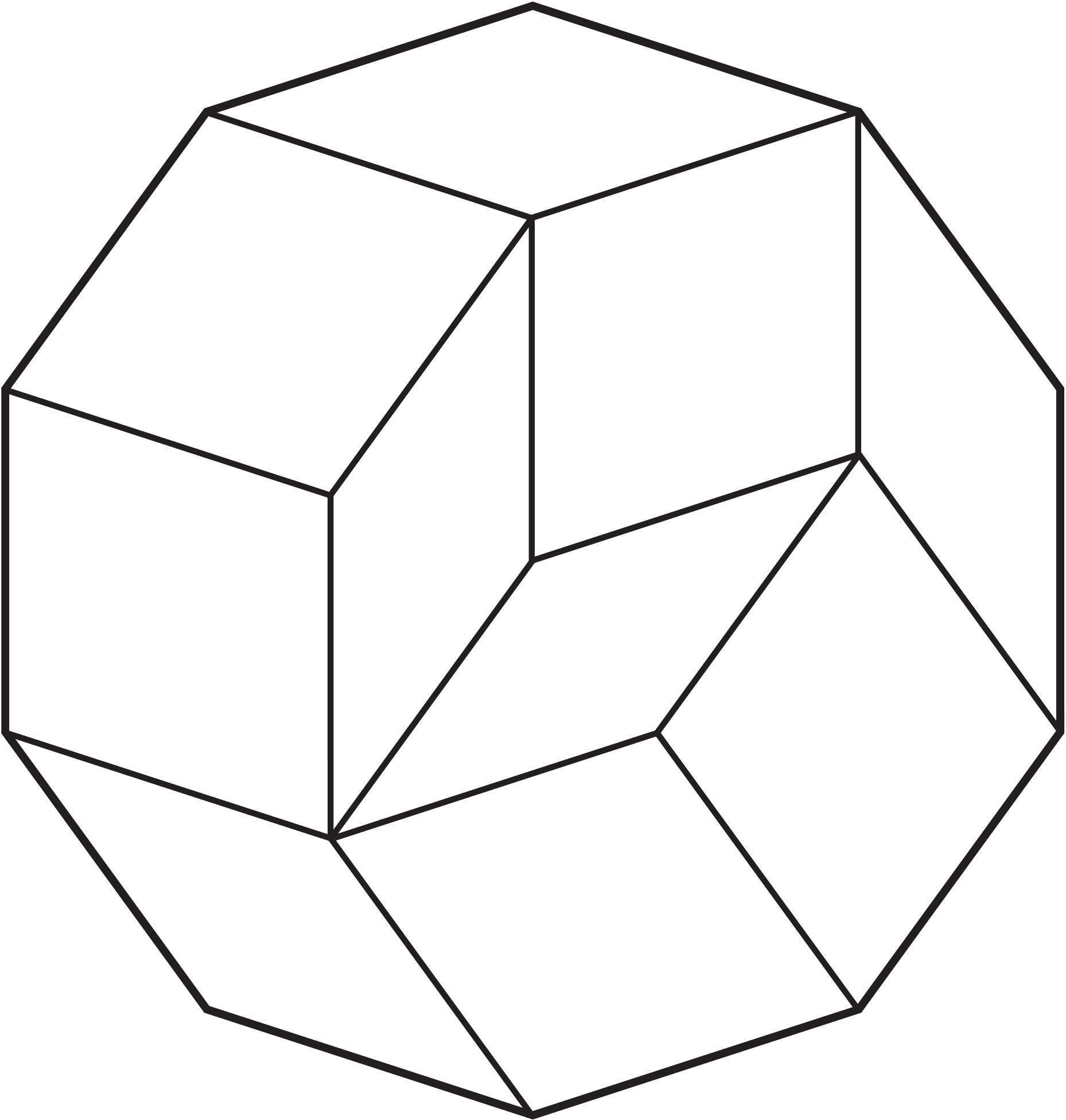}     
&     \includegraphics[width=1in]{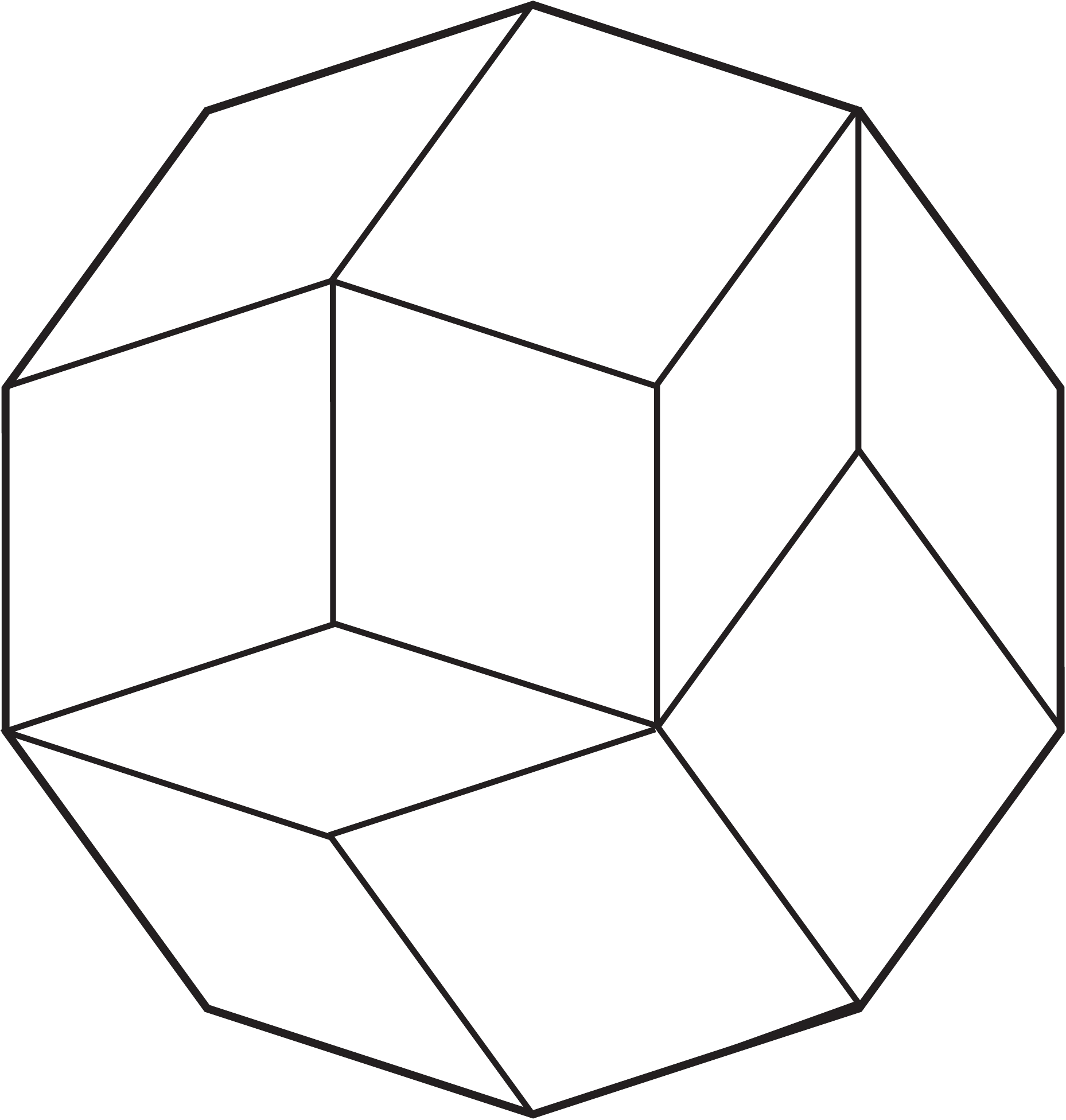} \\

  &    &   &     &    \\

     \includegraphics[width=1in]{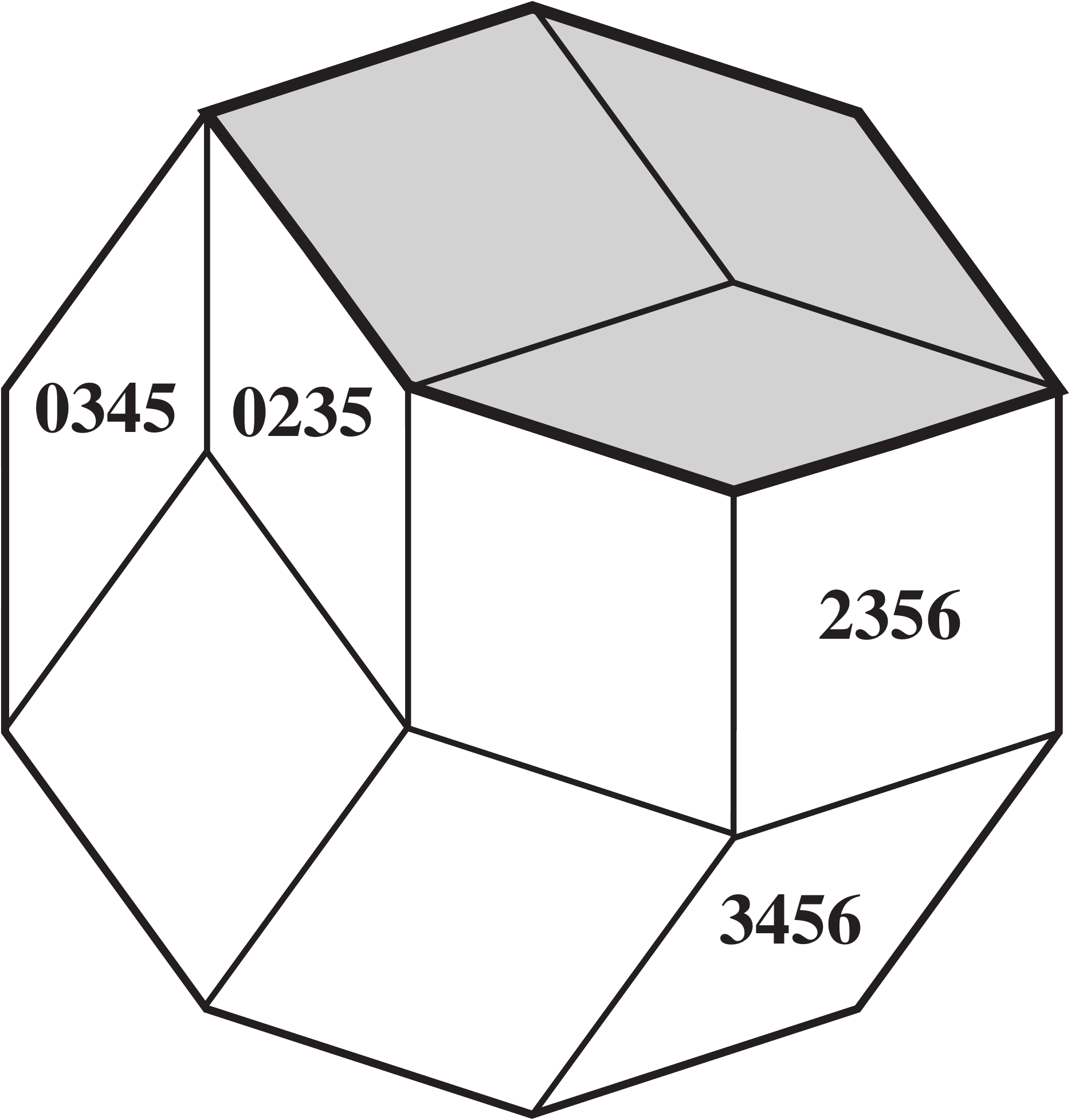}     
  &      \includegraphics[width=1in]{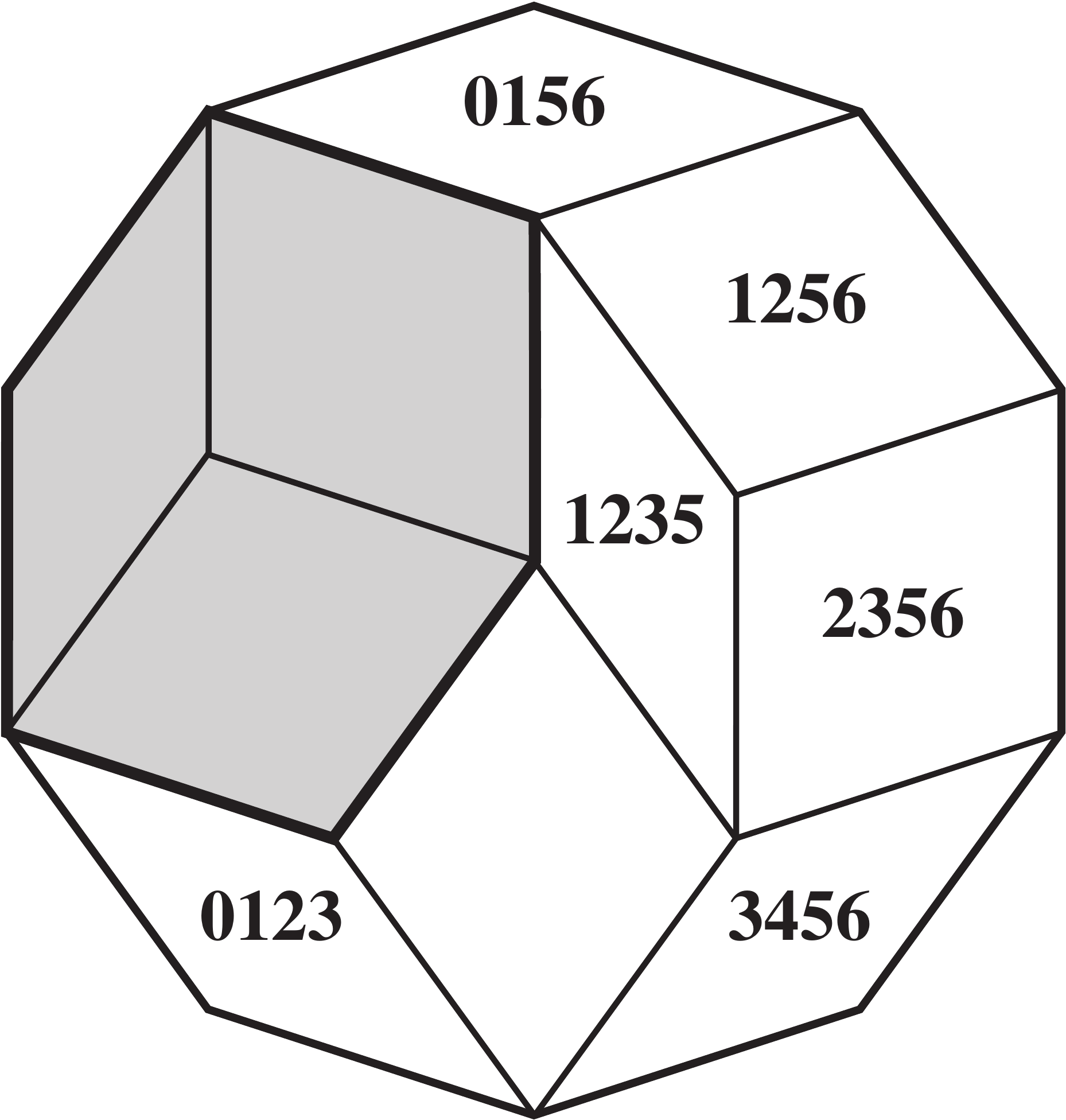}     
 &    \includegraphics[width=1in]{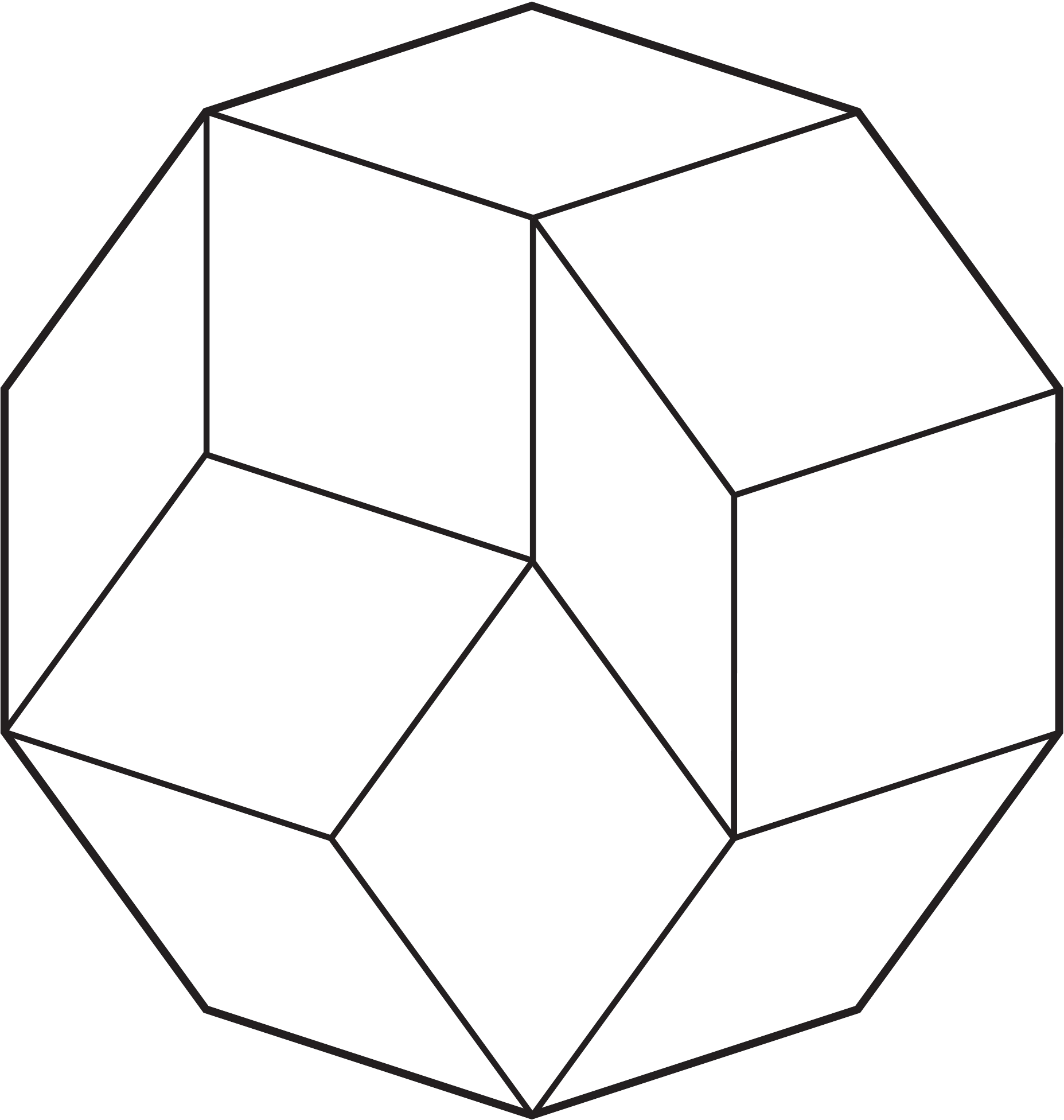}     
&     \includegraphics[width=1in]{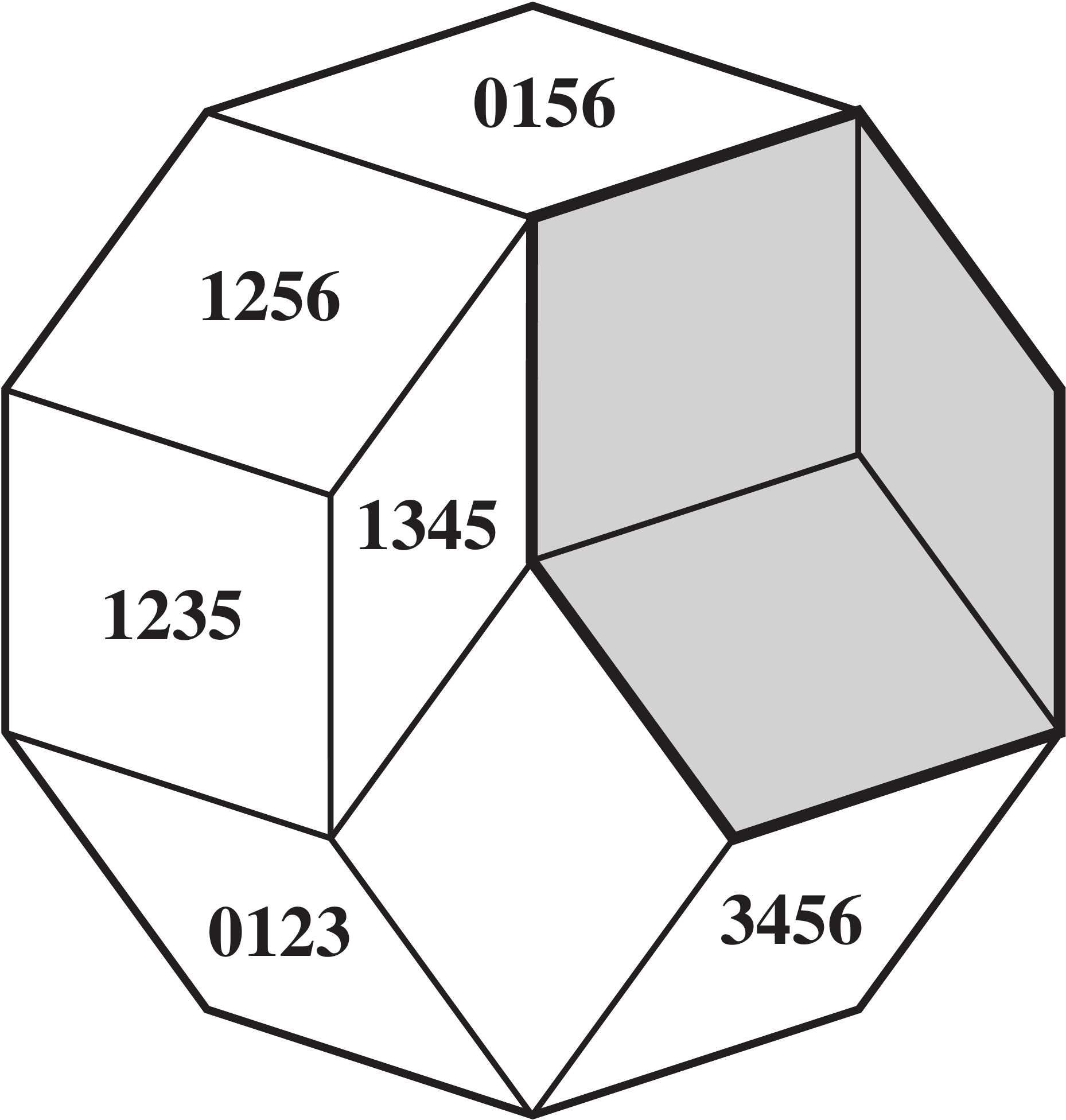}     
&     \includegraphics[width=1in]{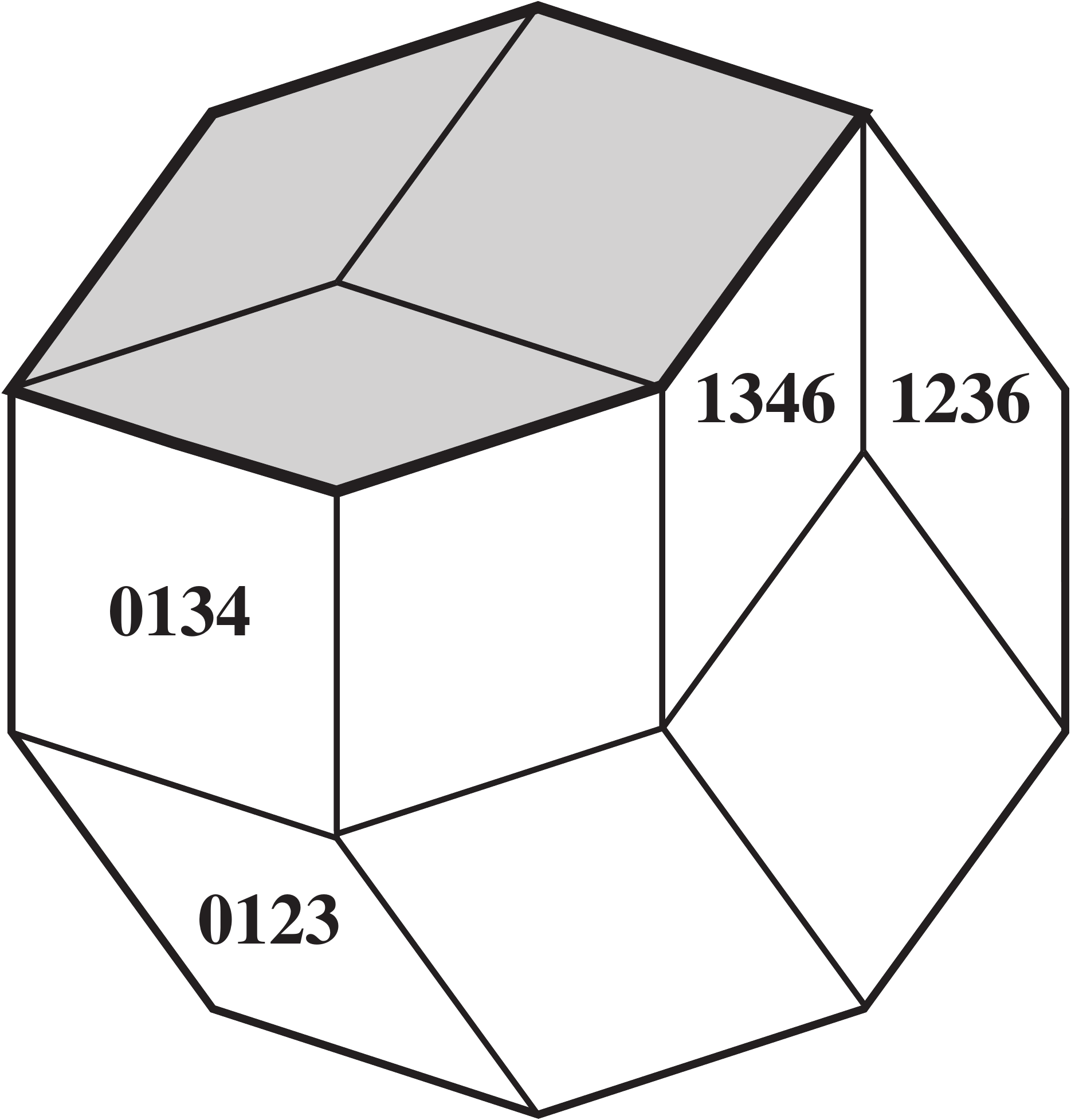} \\

 01256 & 01345   &   &   12356  &01456    \\

      \includegraphics[width=1in]{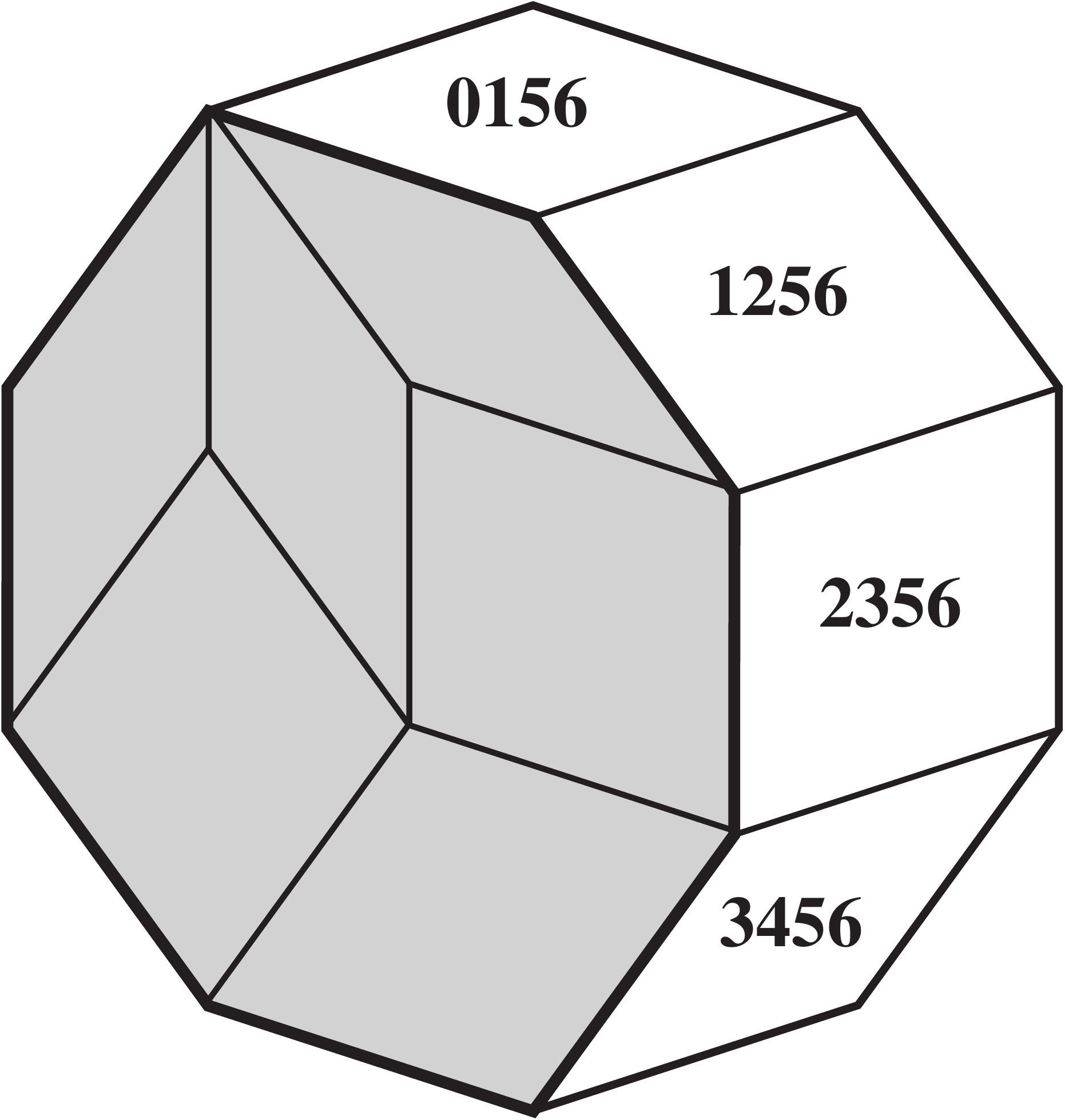}     
  &      \includegraphics[width=1in]{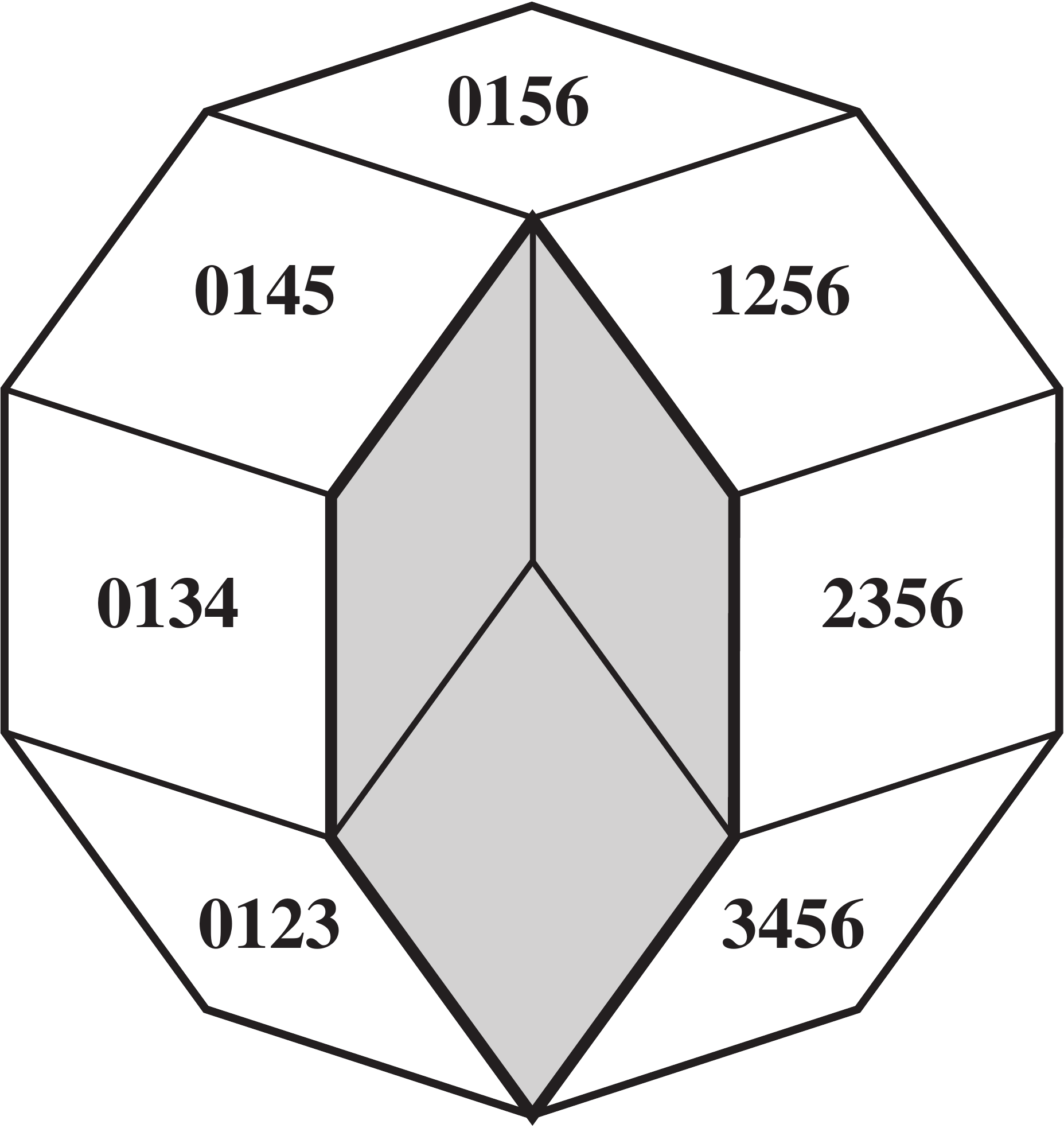}     
 &    \includegraphics[width=1in]{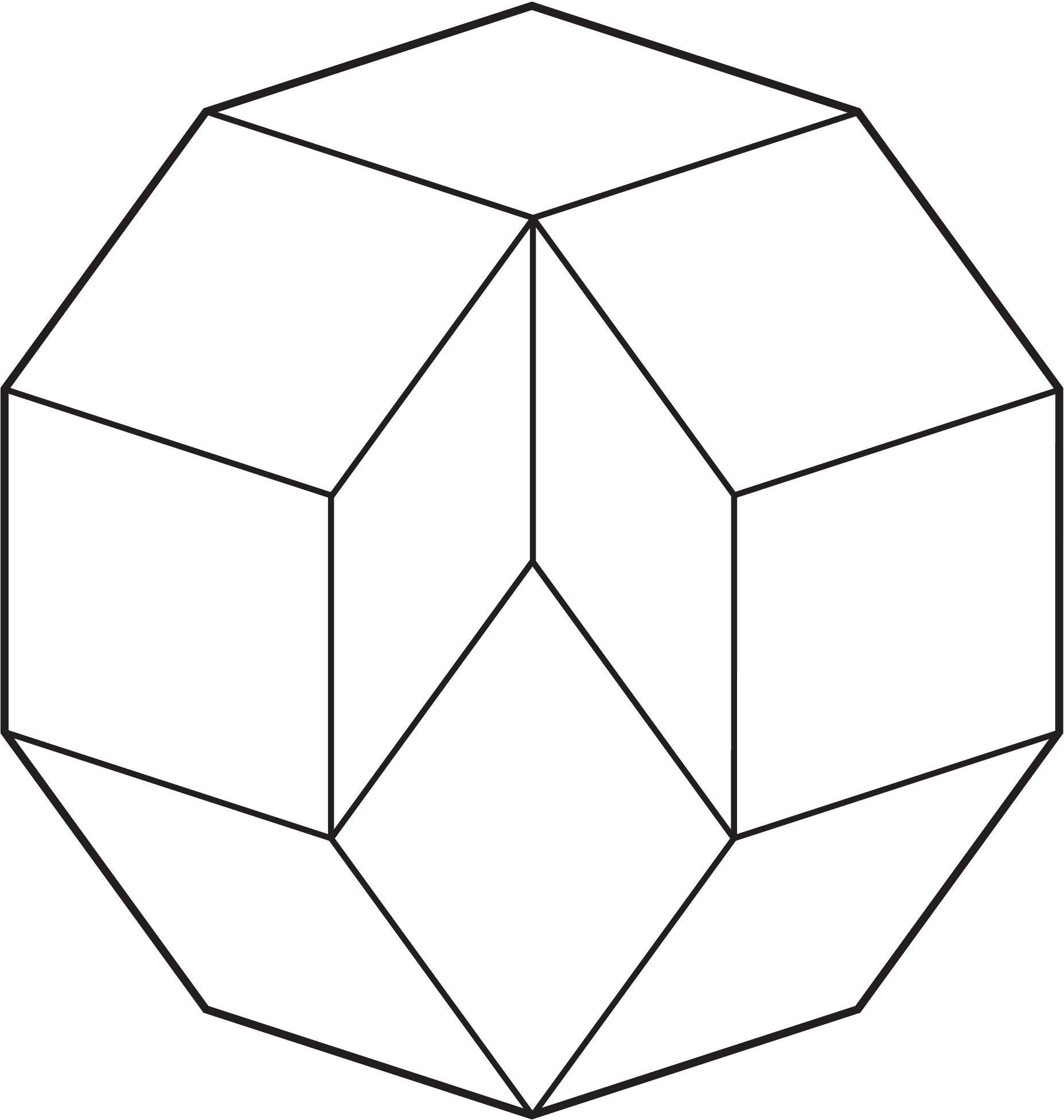}     
&     \includegraphics[width=1in]{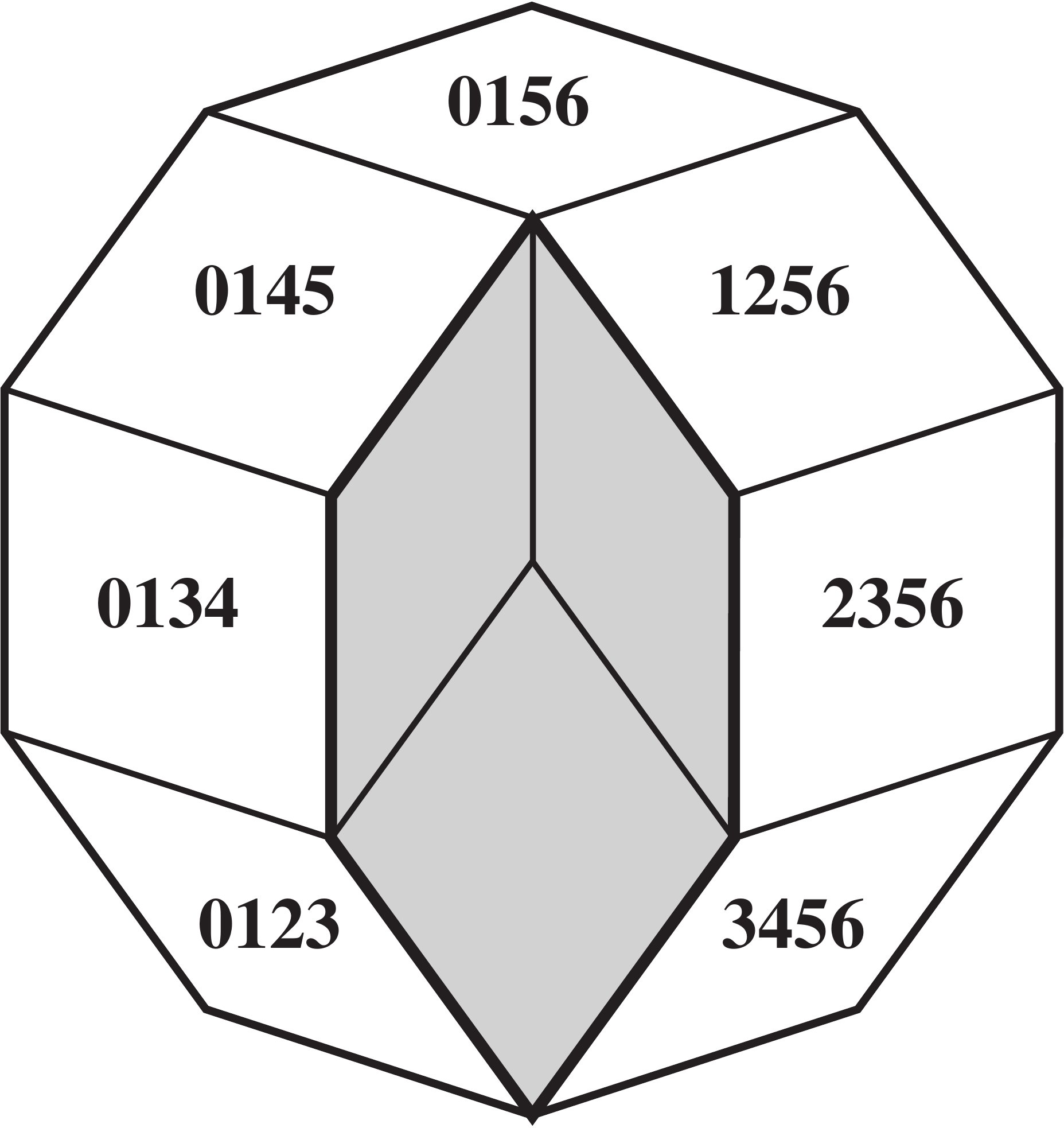}     
&     \includegraphics[width=1in]{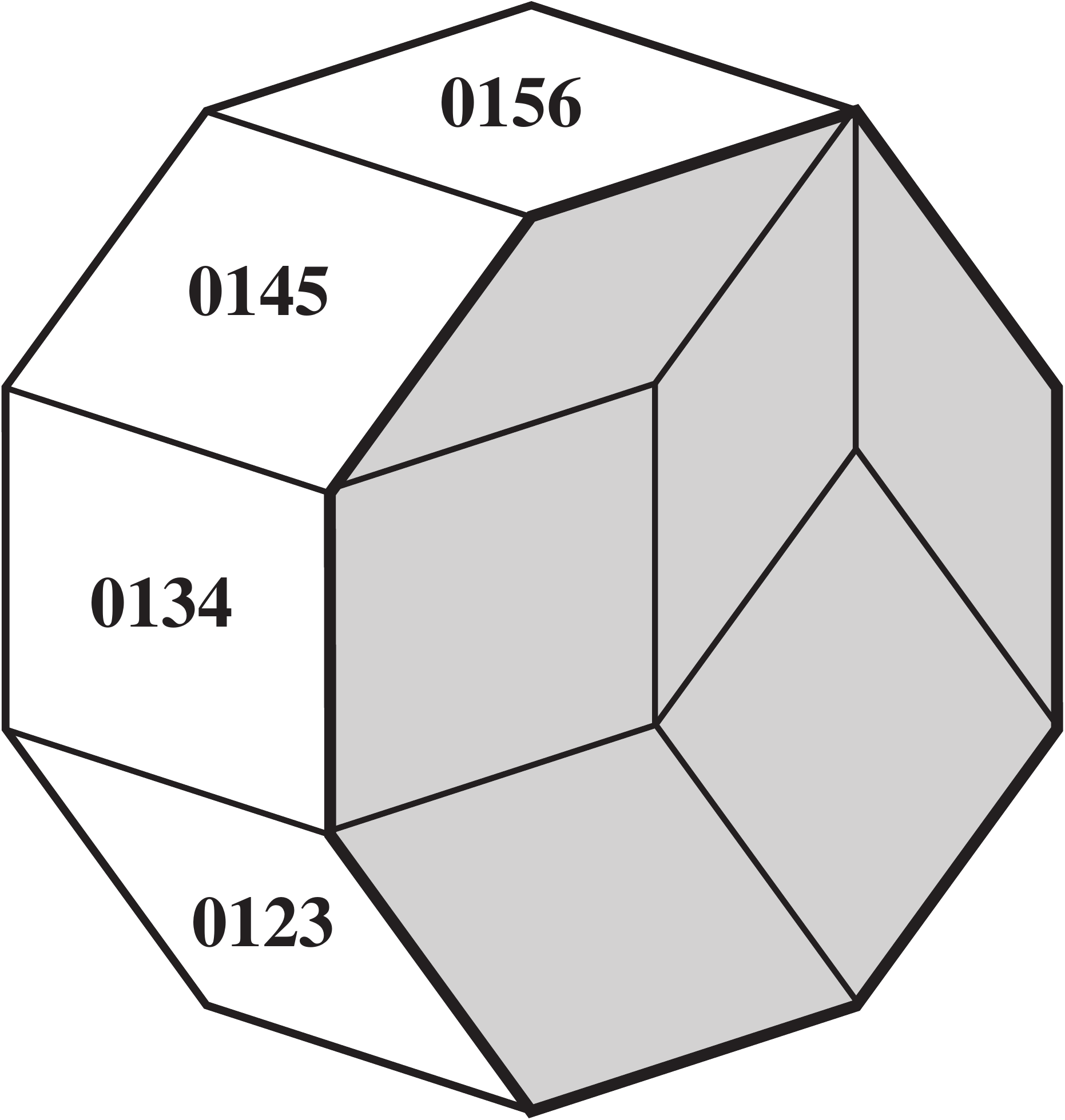} \\

 012345 & 12345   &   & 12345    & 123456   \\
 %6-0  &  6-1  &  6-2  &  6-3  &  6-4
  \label{fig:19}
%\end{figure}
 \end{array}
%\right)
\]

\begin{figure}[htbp] %  figure placement: here, top, bottom, or page
   \centering
   \includegraphics[width=6in]{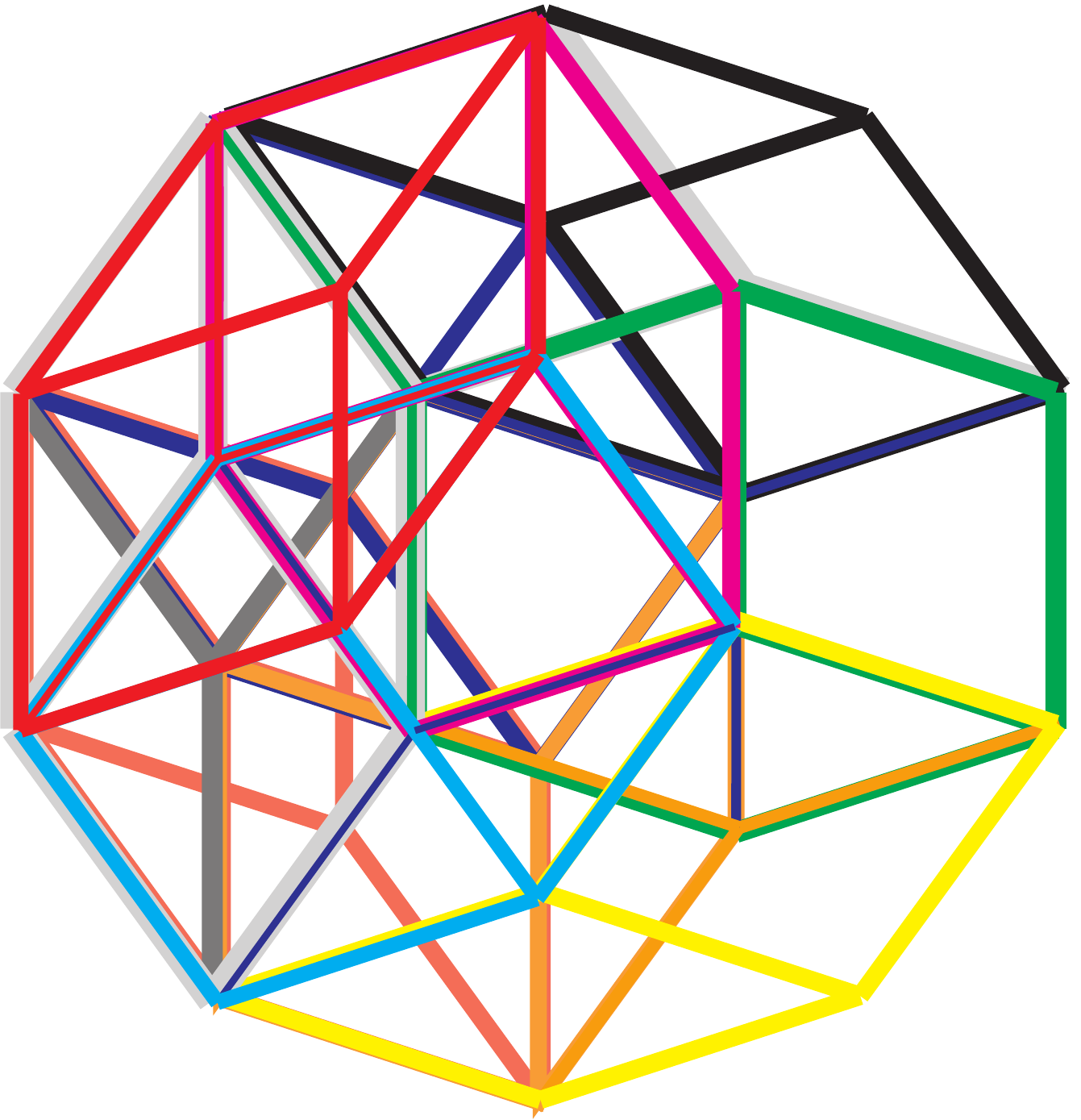} 
   \caption{Front 10 rhombi give $\sigma_2[5]$; back 10  rhombi give $\tau_2[5]$. The ball consists of 10 3-cubes.}
   \label{fig:21-a}
\end{figure}

\begin{figure}[htbp] %  figure placement: here, top, bottom, or page
   \centering
   \includegraphics[width=7in]{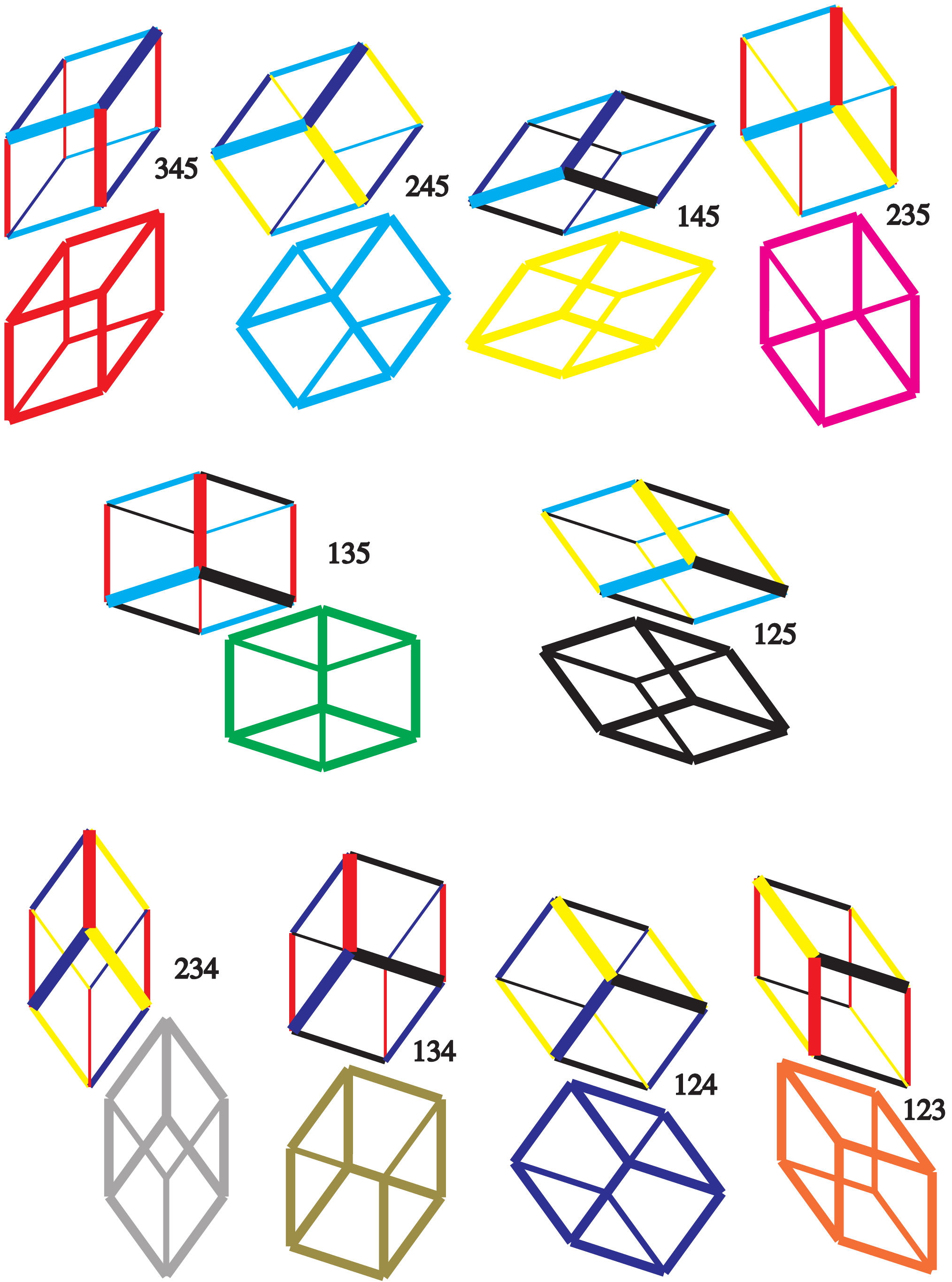} 
%   \caption{Fig21-new-1a}
%   \label{fig:21-a}
\end{figure}

\begin{figure}[htbp] %  figure placement: here, top, bottom, or page
   \centering
\includegraphics[width=7in]{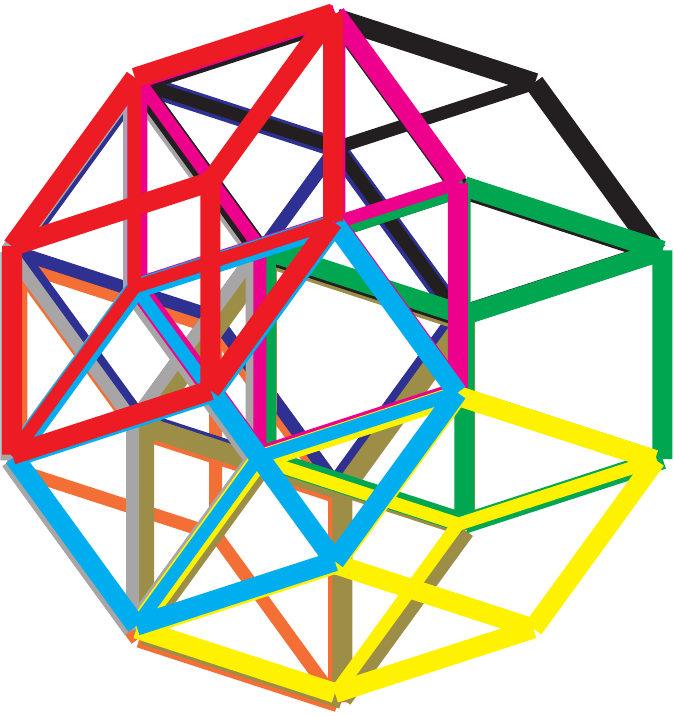}  
\end{figure} 
\begin{figure}[htbp] %  figure placement: here, top, bottom, or page
   \centering
\includegraphics[width=7in]{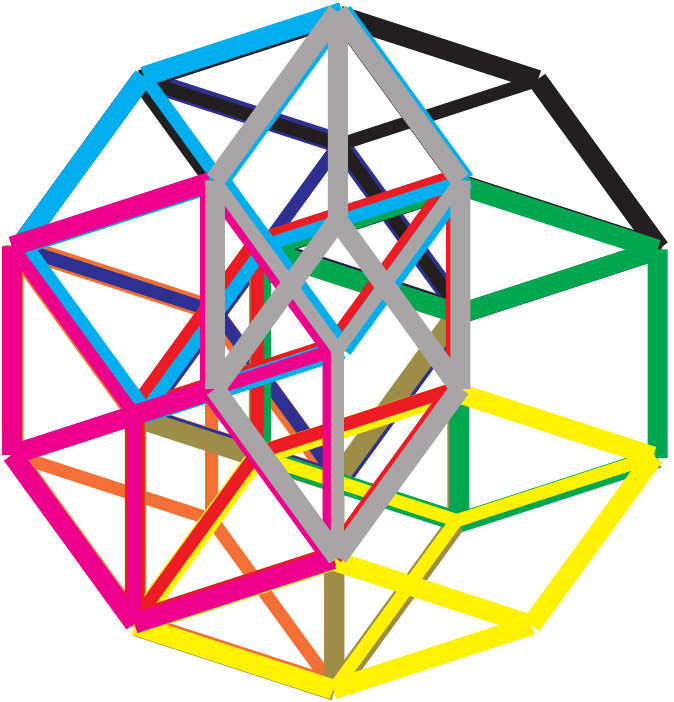}  
\end{figure} 
\begin{figure}[htbp] %  figure placement: here, top, bottom, or page
   \centering
 \includegraphics[width=7in]{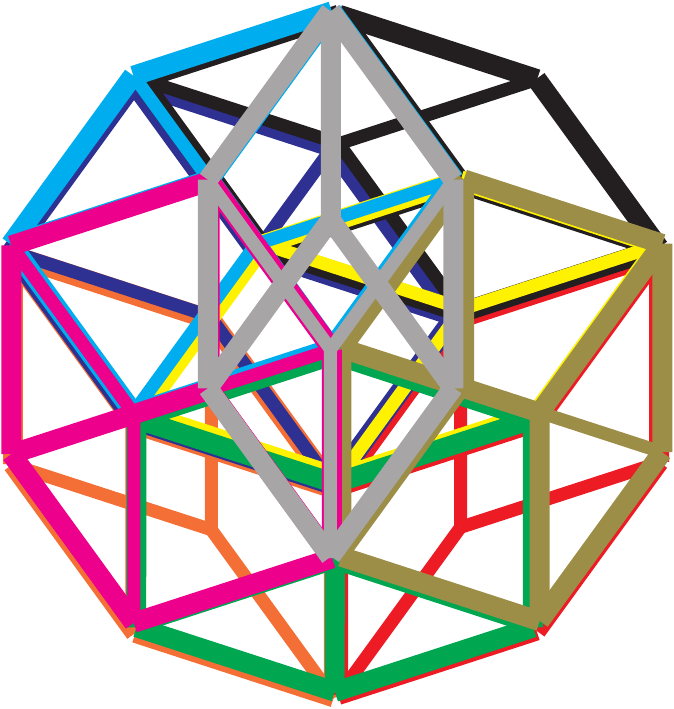} 
 \end{figure} 
\begin{figure}[htbp] %  figure placement: here, top, bottom, or page
   \centering
\includegraphics[width=7in]{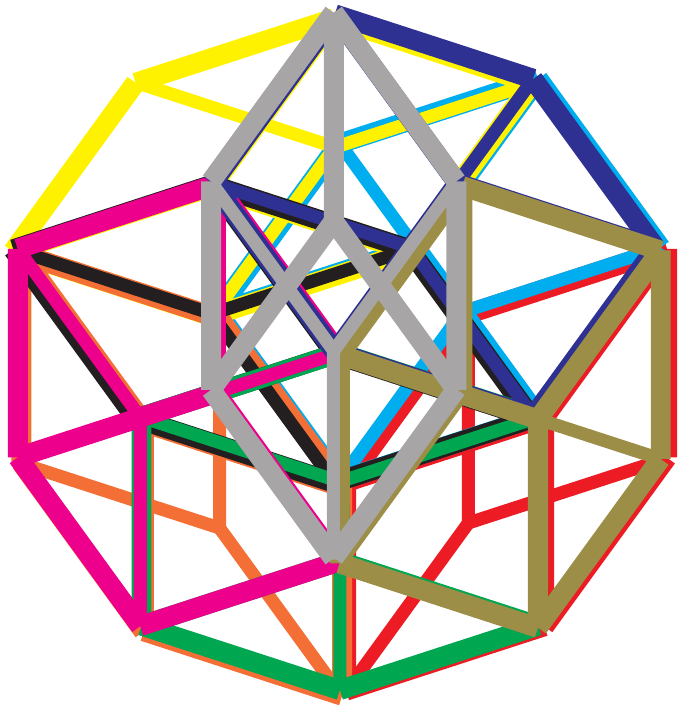}
\end{figure} 
\begin{figure}[htbp] %  figure placement: here, top, bottom, or page
   \centering
  % \caption{Cubes packing a 3-ball: edge coloured and FigC coloured}
 %  \label{fig:2222}
%\end{figure}
%\begin{figure}[htbp] %  figure placement: here, top, bottom, or page
 %  \centering%
%\includegraphics[width=7in]{edge3}\qquad\qquad    
 \includegraphics[width=7in]{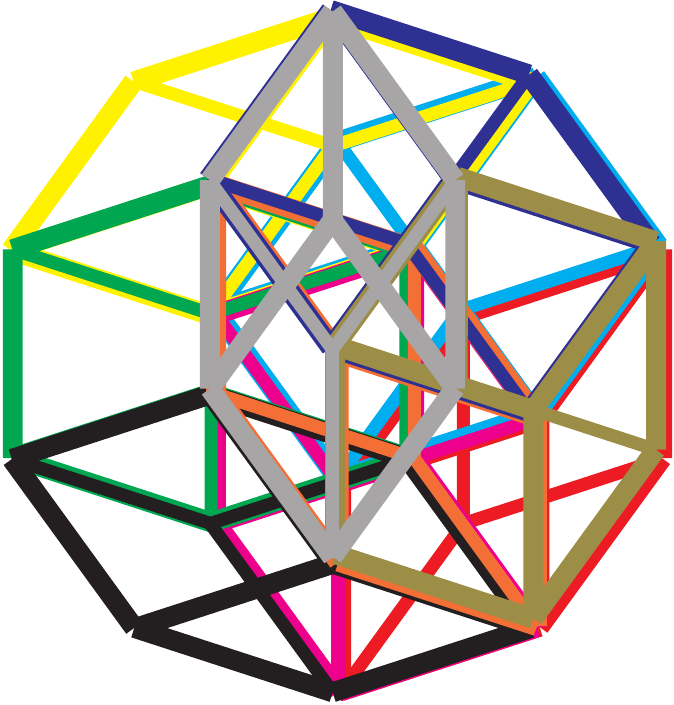} 
 \end{figure} 
\begin{figure}[htbp] %  figure placement: here, top, bottom, or page
   \centering
%includegraphics[width=7in]{edge5}\qquad\qquad     
\includegraphics[width=7in]{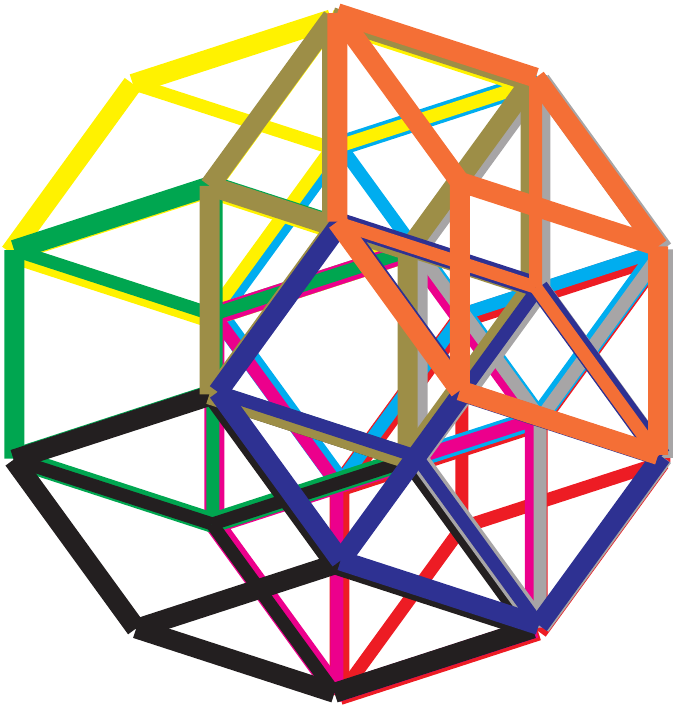}
 %  \caption{Cubes packing a 3-ball: edge coloured and FigC coloured}
 %  \label{fig:2222}
\end{figure}

\begin{figure}[htbp] %  figure placement: here, top, bottom, or page
   \centering
   \includegraphics[width=3.1in]{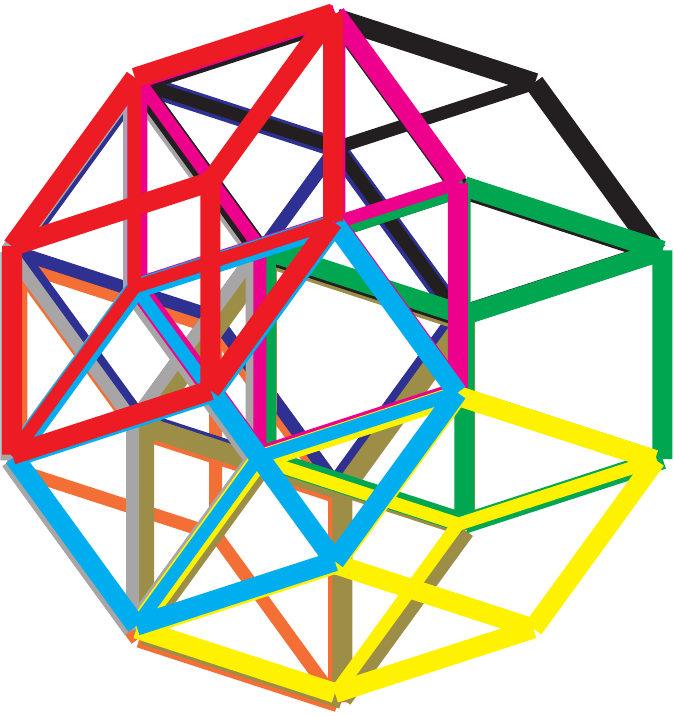} \qquad\qquad    \includegraphics[width=3.1in]{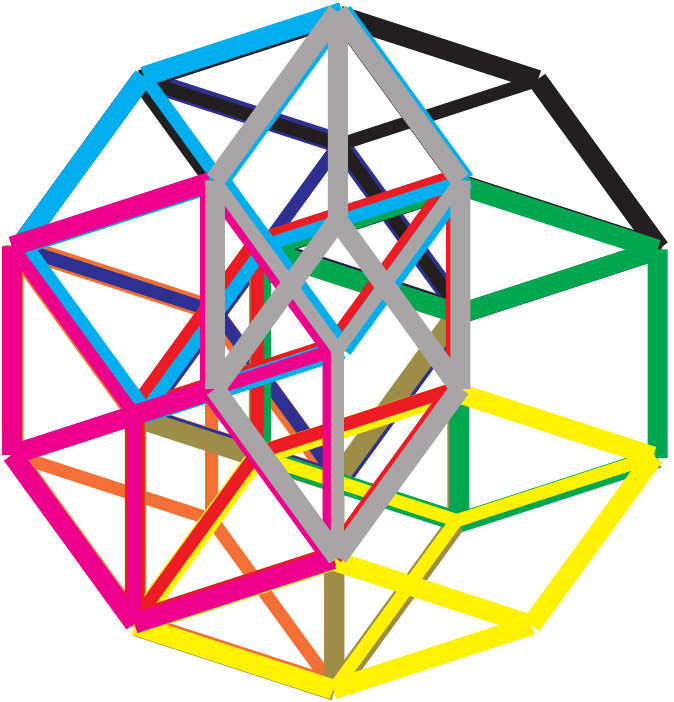}
\includegraphics[width=3.1in]{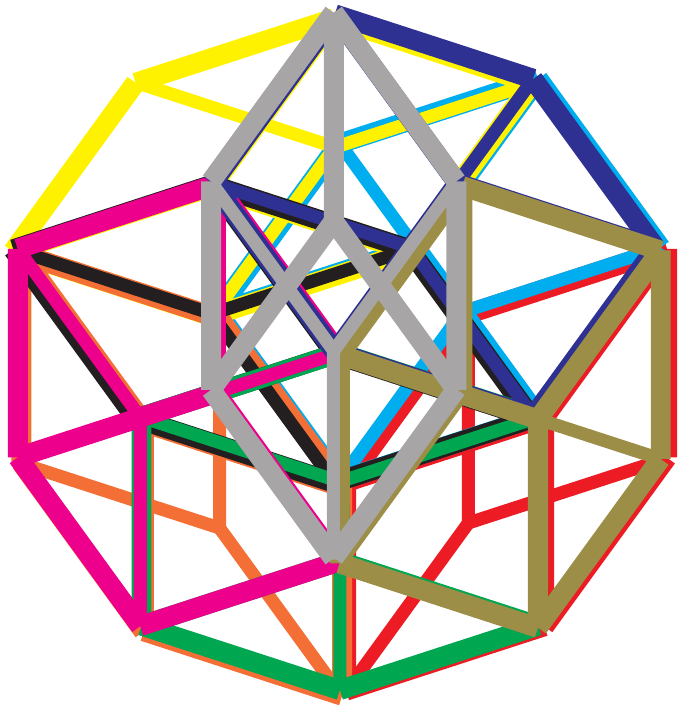}\qquad\qquad    \includegraphics[width=3.1in]{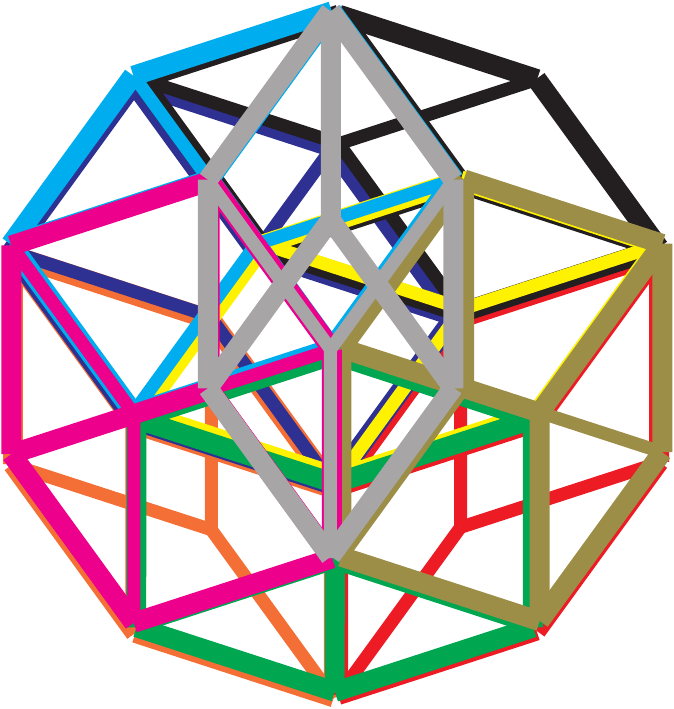}   
\includegraphics[width=3.1in]{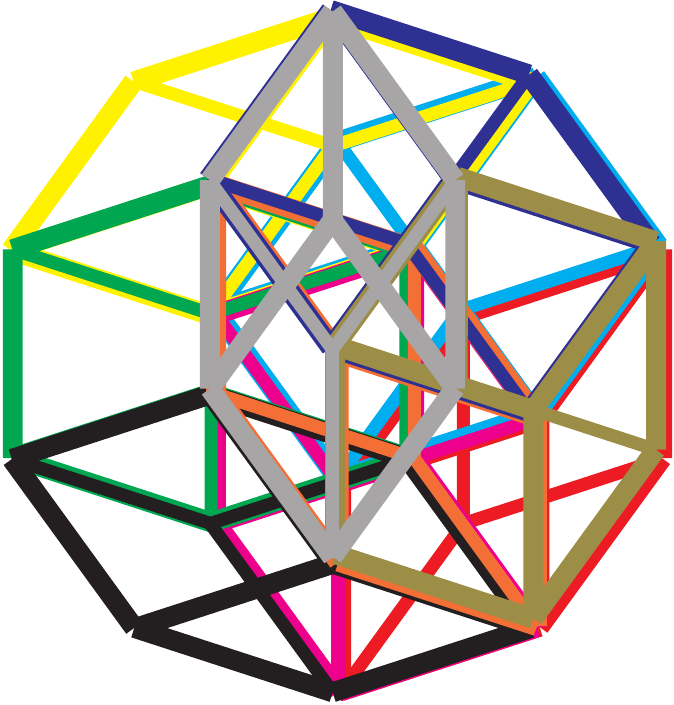}\qquad\qquad     \includegraphics[width=3.1in]{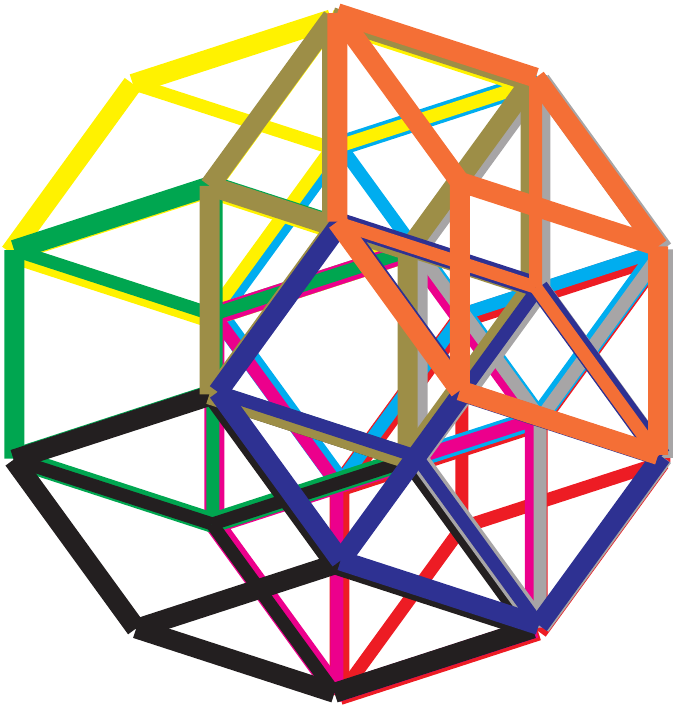} 
  % \caption{Cubes packing a 3-ball: FigC coloured and FigC coloured}
   \label{fig:2222}
\end{figure}
%\pagebreak

\begin{figure}[htbp] %  figure placement: here, top, bottom, or page
   \centering
   \includegraphics[width=3.1in]{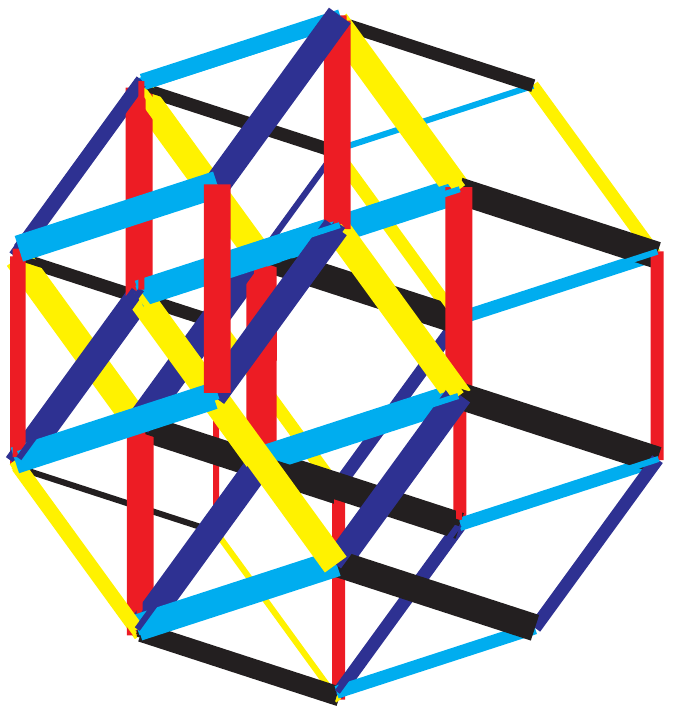} \qquad\qquad    \includegraphics[width=3.1in]{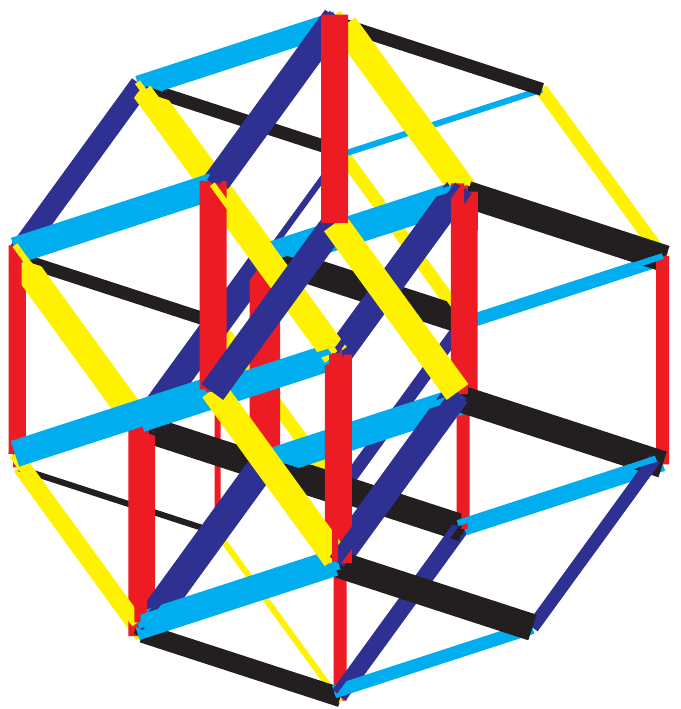}
\includegraphics[width=3.1in]{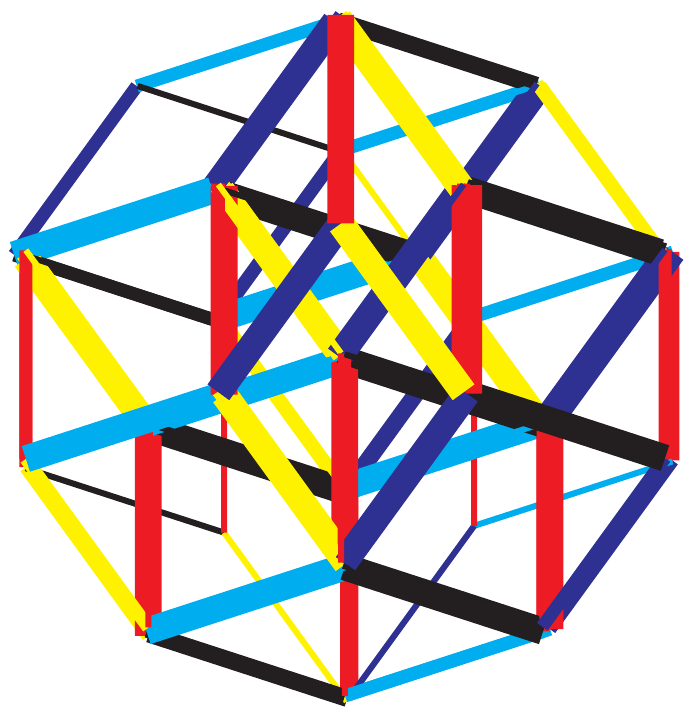}\qquad\qquad    \includegraphics[width=3.1in]{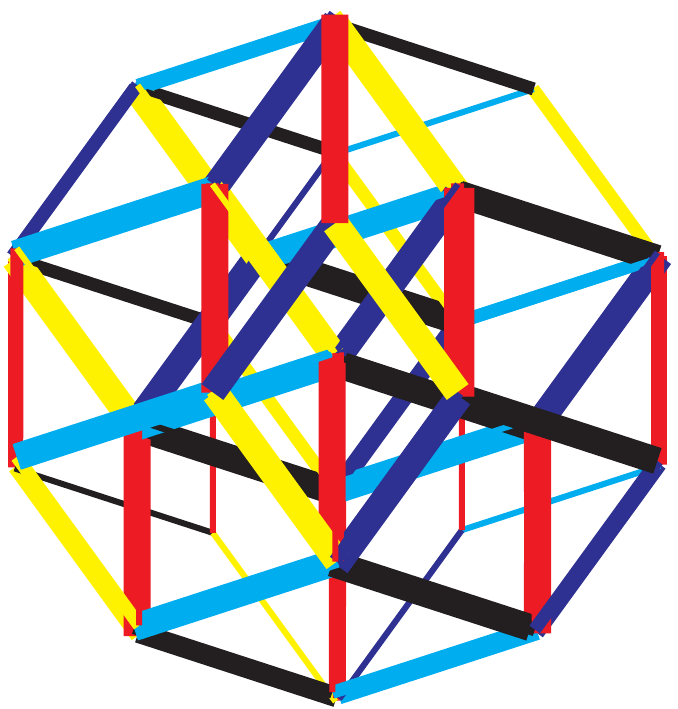}   
\includegraphics[width=3.1in]{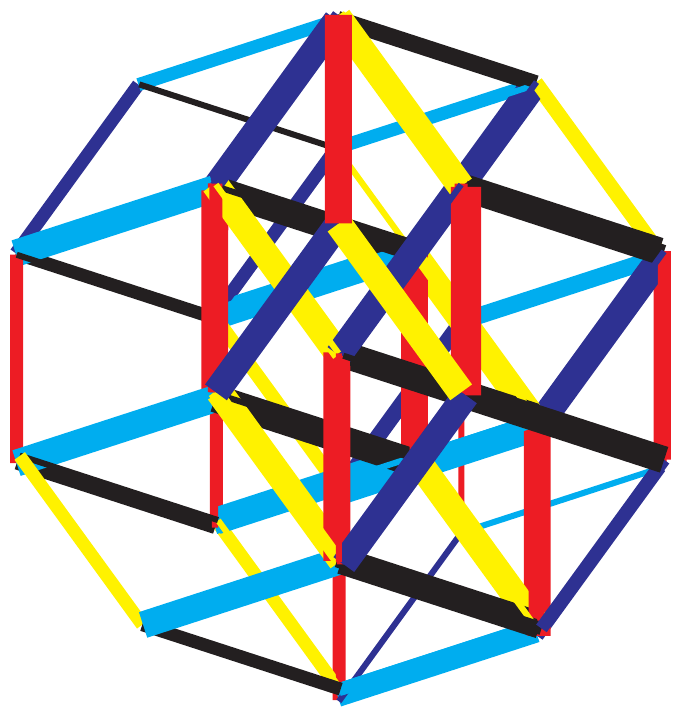}\qquad\qquad     \includegraphics[width=3.1in]{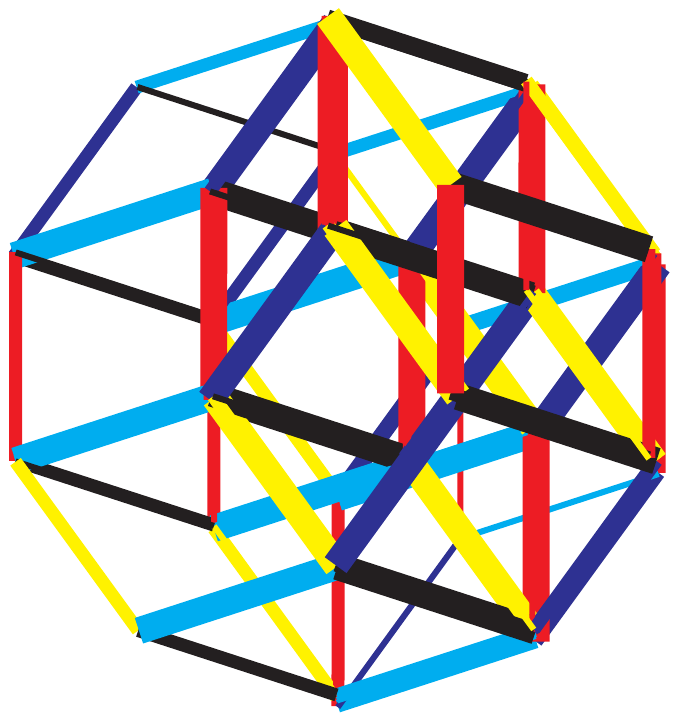} 
  % \caption{Cubes packing a 3-ball: FigE coloured and FigE coloured}
  % \label{fig:2222}
\end{figure}

  \end{document}